\documentclass[12pt,reqno]{amsart}
\usepackage[a4paper,margin=2.5cm]{geometry}
\usepackage{amssymb,amsmath,mathtools,bbm,tikz,tikz-cd,color,mathrsfs,longtable,array}
\usetikzlibrary{arrows.meta,decorations.markings}
\usepackage[bookmarks=true,%
    colorlinks=true,%
    linkcolor=blue,%
    citecolor=blue,%
    filecolor=blue,%
    menucolor=blue,%
    urlcolor=blue,%
    breaklinks=true]{hyperref}

\newtheorem{theorem}{Theorem}[section]
\theoremstyle{definition}
\newtheorem{proposition}[theorem]{Proposition}
\newtheorem{lemma}[theorem]{Lemma}
\newtheorem{definition}[theorem]{Definition}
\newtheorem{remark}[theorem]{Remark}
\newtheorem{corollary}[theorem]{Corollary}

\newtheorem{example}[theorem]{Example}

\def\a{\alpha}
\def\b{\beta}

\def\d{\delta}
\def\e{\epsilon}
\def\ve{\varepsilon}

\def\BZ{\mathbbm Z}
\def\BQ{\mathbbm Q}
\def\BR{\mathbbm R}
\def\BC{\mathbbm C}
\def\BK{\mathbbm K}

\def\BF{\mathbbm F}

\def\CN{\mathcal N}

\def\Hom{\mathrm{Hom}}

\def\p#1{\overline{#1}}

\def\LG{\mathrm{LG}}

\def\CJ{\mathcal{J}}

\def\N{\CN}

\def\1{\boldsymbol{1}}

\renewcommand{\phi}{\varphi}

\def\BE{\mathbbm E}

\def\CA{\mathcal A}
\def\CB{\mathcal B}
\def\CC{\mathcal C}

\def\CF{\mathcal F}
\def\CG{\mathcal G}
\def\CR{\mathcal R}
\def\CS{\mathcal S}
\def\CV{\mathcal V}
   \def\T{\mathsf{T}}

\def\tP{\widetilde{P}}

\def\bx{\bar x}
\def\by{\bar y}

\def\hate{\hat e}

\def\GL{\mathrm{GL}}
\def\diag{\mathrm{diag}}
\def\sgn{\mathrm{sgn}}
\def\spn{\mathrm{span}}
\def\Lam{\Lambda}               
\def\bD{\overline{D}}      
\def\DD{\mathsf D}          
\def\Der{\mathrm{Der}}
\def\rk{\mathrm{rank}}
\def\CT{\mathcal T}               
\def\im{\operatorname{im}}       
\def\corank{\mathrm{corank}}
\def\lk{\mathrm{lk}}
\def\sm{\smallsetminus}

\def\ind#1{\bigl[#1\bigr]}  
\def\tQ{\widetilde{Q}}
\def\WW{\mathsf W}
\def\FF{\mathfrak F}
\def\tpi{\tilde\pi}
\def\vs{\varsigma}
\def\Kash{\mathsf K}
\def\GLf{\mathsf G}
\def\Goe{\mathsf A}
\def\Stk{\mathsf S}
\def\CU{\mathcal U}
\def\CW{\mathcal W}
\def\sign{\operatorname{sign}}

\def\lb{[\![}\def\rb{]\!]}
\def\kk{\Bbbk}
\def\xs{\mathsf x}\def\ys{\mathsf y}\def\sv{\mathsf s}\def\us{\mathsf u}
\def\pth{\mathsf p}\def\mer{\mathsf r}
\def\tpth{\tilde{\mathsf p}}\def\tstar{\tilde\star}
\def\CZ{\mathcal Z}\def\CDp{\mathcal D}\def\CL{\mathcal L}
\def\vth{\Phi}\def\mhol{\mathsf m}
\def\out{\mathrm{out}}\def\inn{\mathrm{in}}
\begin{document}

\title[Deformed Alexander models of the Links--Gould invariant]{Deformed Alexander models of the Links--Gould invariant}

\author[M. Chrisman]{Micah Chrisman}
\address{Chrisman: Department of Mathematics, The Ohio State University at Marion, 1465 Mount Vernon Avenue, Marion, OH 43302}
\email{chrisman.76@osu.edu}

\author[M. Harper]{Matthew Harper}
\address{Harper: Department of Mathematics, Michigan State University, East Lansing, Michigan, USA}
\email{mrhmath@proton.me}

\author[B.-M. Kohli]{Ben-Michael Kohli}
\address{Kohli: Section de Math\'ematiques, Universit\'e de Gen\`eve, rue du Conseil-G\'en\'eral 7-9, 1205 Gen\`eve, Switzerland}
\email{bm.kohli@protonmail.ch}

\author[G. Tahar]{Guillaume Tahar}
\address{Tahar: Department of Mathematics with Computer Science, Guangdong Technion - Israel Institute of Technology (GTIIT), Shantou, Guangdong, China}
\email{guillaume.tahar@gtiit.edu.cn}


\date{\today}

\begin{abstract}
The Links--Gould invariant is a two-variable polynomial quantum link invariant arising from the Reshetikhin--Turaev construction applied to the quantum supergroup $U_q\mathfrak{gl}(2|1)$. It generalizes the Alexander polynomial and shares many of its topological properties.
\par
In a recent work by two of the authors, the Links--Gould polynomial is obtained by coupling two determinants which deform the unreduced Burau determinant of the Alexander polynomial. In this monograph, we show that each of them is the determinant of a family of linear forms, one for each crossing, on the lattice of flows of the diagram, in an arbitrary basis. Different bases deform four models of the Alexander polynomial into models of the Links--Gould polynomial: the Seifert matrix of the canonical Seifert surface, Alexander's face matrix, Kashaev's symmetric matrix, and Kauffman's state sum. The crossing variables are identified with the variable of the Alexander polynomial, localized at the bands of this surface; when they all equal $-1$, the lattice of flows gives the Goeritz matrices and the Gordon--Litherland forms of the state surfaces of the diagram, with their signatures. The coupling is carried by the same lattice: the Links--Gould polynomial is the mixed discriminant of bilinear forms, one for each crossing, which are read off a local system of rank one on the complement of the link and of the co-cores of the bands of the Seifert surface; this local system extends to the complement of the link at the free fermion point only.
\par
As a consequence, we reprove the first two Tait conjectures, on the crossing number and on the writhe of reduced alternating diagrams. These results illustrate the hybrid Alexander/Jones nature of the Links--Gould invariant. This study is a step toward a topological interpretation of this invariant and prepares the tools needed for a possible categorification.
\end{abstract}
\maketitle

{\footnotesize
\tableofcontents
}

\part*{Introduction}

The Links--Gould invariant $\LG^{2,1}(L;t_0,t_1)=\LG(L;t_0,t_1)\in\BZ[t_0^{\pm1},t_1^{\pm1}]$ of an oriented link $L$ is obtained from the Reshetikhin--Turaev construction applied to a four-dimensional representation of the Hopf superalgebra $U_q\mathfrak{gl}(2|1)$, see \cite{LinksGould, DeWitKauffmanLinks, DW01}. It is closely related to the Alexander polynomial $\Delta_L$, to which it specializes in two different ways, by \cite{DWIL} and by \cite{Kohli, KPM} respectively:
$$
\LG(L;t_0,-t_0^{-1}) = \Delta_L(t_0^2)\,, \qquad
\LG(L; t_0, t_0^{-1}) = \Delta_L(t_0)^2 .
$$
Several of the most distinctive three-dimensional properties of the Alexander polynomial have counterparts for $\LG$ which contain and improve upon the classical statements: the degree of $\LG$ bounds the genus of a knot from below \cite{Kohli-Tahar}, $\LG$ is monic on fibered links \cite{LNVdV25}, and $\LG$ is an alternating polynomial on alternating links \cite{HK}.

The last of these results rests on a determinant formula for $\LG$, obtained in \cite{HK} by extending the free-fermion model of Kauffman and Saleur \cite{KauffmanSaleur} for the Alexander--Conway polynomial to the $V_1$-polynomial of Garoufalidis and Kashaev \cite{GaroufalidisKashaev, GHKKST}, which is $\LG$: on every connected diagram $D$ of $L$ with crossings $1,\dots,c$ of signs $\ve_j$, and for every edge $\hate$ of $D$,
\begin{equation}\label{eq.intro.formula}
   \LG(L;t_0,t_1)=(-1)^{\mu-1}\ \BE_D\bigl[P_{D,\hate}(x)\,\tP_{D,\hate}(y)\bigr],
\end{equation}
where $\mu$ is the number of components of $L$, $P_{D,\hate}(x)$ and $\tP_{D,\hate}(y)$ are the determinants of two ``Burau matrices'' of the diagram in which an independent variable $x_j$, resp.\ $y_j$, is attached to the crossing $j$, and $\BE_D$ is the $\BK$-linear map determined by $x_j\mapsto t_0^{\ve_j}$, $y_j\mapsto t_1^{\ve_j}$ and $x_jy_j\mapsto1$ \cite[Theorem~7.8]{HK}.

In words, the formula does the following. We call $x_j$ and $y_j$ the \emph{crossing variables}, and $P_{D,\hate}$ and $\tP_{D,\hate}$ the \emph{edge polynomials} of the diagram, of the first and of the second \emph{flavor}. Each of the two edge polynomials has degree at most one in each of its variables. So every monomial of the product $P_{D,\hate}(x)\,\tP_{D,\hate}(y)$ contains, for each crossing $j$, exactly one of the four factors $1$, $x_j$, $y_j$ and $x_jy_j$, and $\BE_D$ replaces these factors by $1$, $t_0^{\ve_j}$, $t_1^{\ve_j}$ and $1$ respectively. The last rule is the point. Without it, $\BE_D$ would be the substitution $x_j=t_0^{\ve_j}$, $y_j=t_1^{\ve_j}$, and \eqref{eq.intro.formula} would be the product of the two Burau determinants, which is $\Delta_L(t_0)\Delta_L(t_1)$ up to a unit. This substitution sends $x_jy_j$ to $(t_0t_1)^{\ve_j}$ instead of $1$. So the two rules agree at $t_0t_1=1$, and in general \eqref{eq.intro.formula} is an expansion of $\LG$ around this \textit{free fermion point}. Since $\BE_D$ is not a substitution, the image of a product is not the product of the images: $\BE_D$ couples the two determinants, crossing by crossing, and we call it the \emph{coupling}. See Section~\ref{sec.formula} for a detailed exposition.

\emph{The trefoil.} The simplest example shows what happens. For the standard diagram of the right-handed trefoil, the closure of the braid $\sigma_1^3$, with its three positive crossings numbered along the braid, and for a suitable edge $\hate$, the edge polynomials are
\[
   P_{D,\hate}(x)=1-x_1+x_1x_2\qquad\text{and}\qquad\tP_{D,\hate}(y)=1-y_3+y_2y_3
\]
(Example~\ref{ex.trefoilflows}). Their product has nine monomials. Eight of them contain no product $x_jy_j$, and on these $\BE_D$ is the substitution $x_j=t_0$, $y_j=t_1$. The ninth one is $x_1x_2\cdot y_2y_3$, which contains the two variables of the crossing $2$; its image is $t_0t_1$, and not $t_0^2t_1^2$. Hence
\begin{align*}
   \LG(3_1;t_0,t_1)&=1-t_0-t_1+t_0^2+2t_0t_1+t_1^2-t_0^2t_1-t_0t_1^2\\
   &=(1-t_0+t_0^2)(1-t_1+t_1^2)+t_0t_1\,(1-t_0t_1):
\end{align*}
the product of the Alexander polynomials of the trefoil in the variables $t_0$ and $t_1$, up to units, and a correction which vanishes at the free fermion point. We shall come back to this example.

\subsection*{The question}

The one-variable Alexander polynomial has several classical expressions. The first is the determinant of the Alexander matrix of a presentation of the link group; this is what the Burau determinants are at the \emph{Alexander point} $x_j=t^{\ve_j}$, where $P_{D,\hate}(t^{\ve})$ is the determinant of the Alexander matrix of the Wirtinger presentation with one generator for each edge of $D$, cut along $\hate$. The second is $\det(tV-V^{\T})$, where $V$ is the Seifert matrix of a Seifert surface of $L$: a square matrix on the first homology of the surface, of size $2g+\mu-1$. Classical three-dimensional properties of $\Delta_L$, like the genus bound $\deg\Delta_L\le2g+\mu-1$ and the monicity on fibered links, are read off this second form, because the matrix is small and of degree one in $t$. The third is the oldest one: Alexander's matrix \cite{Alexander28}, whose rows are the crossings and whose columns are the faces of the diagram, whose expansion is Kauffman's state model \cite{KauffmanFKT} and of which Kashaev's symmetric matrix \cite{Kashaev21} is a quadratic version. At $t=-1$ these models meet the matrices of Goeritz \cite{Goeritz33} and the forms of Gordon and Litherland \cite{GordonLitherland}, which are defined for every spanning surface of the link, that is, for every compact surface whose boundary is the link, whether this surface is orientable, like a Seifert surface, or nonorientable, like most checkerboard surfaces.

The determinant formula \eqref{eq.intro.formula} is a deformation of the first of these expressions: at the Alexander point its two determinants are those of Alexander matrices of the link group, and in general every crossing has a variable of its own. This raises the question which we study. Can the other classical expressions be deformed in the same way? In other words: are the edge polynomials also determinants of matrices of Seifert type, and of matrices indexed by the faces of the diagram, in which the crossing variables are kept? Do they have a state model? And which properties of $\LG$ can be read off such expressions?

The purpose of this monograph is to give all these forms to the two determinants of \eqref{eq.intro.formula}, with the crossing variables included, and to draw the consequences for $\LG$. The point of view which organizes the exposition is that they are \emph{one} object, written in different bases. The cut diagram carries a lattice of rank $c$, the lattice of its integral \emph{flows}; each crossing defines a linear form on this lattice which contains the variable of the crossing and no other one; and the edge polynomials are the determinants of these $c$ forms in an arbitrary basis of the lattice (Theorem~\ref{thm.intro.flows}). The boundaries of the faces of the diagram form a basis, in which one finds Alexander's matrix and, after expansion, Kauffman's states. The Seifert circles, together with cycles of the canonical Seifert surface $F_D$ of the diagram, form another basis, in which one finds a deformation of the Seifert matrix. More generally there are bases for each of the $2^c$ ways of smoothing the crossings of the diagram, and behind each smoothing there is a spanning surface of the link, its state surface.

The following table lists the bases which are used, the classical models to which they lead, and what is obtained here.

\begin{center}
\small
\renewcommand{\arraystretch}{1.25}
\begin{tabular}{@{}>{\raggedright\arraybackslash}p{0.25\textwidth}>{\raggedright\arraybackslash}p{0.28\textwidth}>{\raggedright\arraybackslash}p{0.39\textwidth}@{}}
\hline
\emph{Basis of the lattice of flows}&\emph{Classical model}&\emph{In this monograph}\\
\hline
the boundaries of the faces&Alexander's matrix&the same matrix with the crossing variables, for each of the two determinants (Section~\ref{ssec.facematrices})\\
&Kashaev's matrix&its deformation by the crossing variables of the two flavors, whose coupled determinant is $\LG$ up to sign (Section~\ref{sec.faces})\\
&Kauffman's state sum&a state sum for $\LG$ over the pairs of states (Section~\ref{sec.kauffman})\\
the Seifert circles and cycles of the canonical surface&the matrix $tV-V^{\T}$ of the canonical surface&the deformed Seifert matrices; the degrees and the extreme coefficients of $\LG$ (Sections~\ref{sec.seifert} and~\ref{sec.degree})\\
the circles of a smoothing and cycles of its state surface&at $t=-1$: the Goeritz matrices and the forms of Gordon and Litherland&at $x=-1$: these forms for all the state surfaces, with their determinants and their signatures (Sections~\ref{sec.statesurfaces} and~\ref{sec.minusone})\\
\hline
\end{tabular}
\end{center}

The lattice of flows is also the group of the $1$-cycles which are carried by the diagram, and this is what makes the comparison with the surfaces exact. We show that the linking numbers of these cycles with the curves pushed off a state surface are sums of contributions of the crossings, each of which is the flow through one edge of the crossing times the number of passages through its band (Theorem~\ref{thm.intro.linking}). For the canonical Seifert surface of an \emph{arbitrary} diagram this identifies the linear form of a crossing: it is the contribution of the band of the crossing to $tV_F-V_F^{\T}$, in which $t$ has been replaced by the variable of the crossing. In this sense the crossing variable $x_j$ is the variable $t$ of the Alexander polynomial, localized at the band $j$ of the Seifert surface; the second determinant of \eqref{eq.intro.formula} is related in the same way to the Seifert surface on the other side of the plane of the diagram.

From this we obtain $\nu\times\nu$ matrices on $H_1(F_D)$, with $\nu=2g(F_D)+\mu-1$, whose determinants are the edge polynomials and which are $tV_F-V_F^{\T}$ at the Alexander point, and we show that the canonical genus degree bound, the extreme coefficients, the behavior under Murasugi sums along Seifert circles and, for alternating links, the fiberedness criterion and the Hopf invariant of the fibration come out of these matrices exactly as the classical statements come out of $tV_F-V_F^{\T}$.

In the basis of the faces we relate the two determinants to Kashaev's matrix; at the point where all the crossing variables are $-1$ we find the Goeritz matrices and the forms of Gordon and Litherland of all the state surfaces, and we obtain their determinants and their signatures from the lattice of flows. Expanding the face matrices gives the counterpart of Kauffman's state model: a state sum for $\LG$ over pairs of Kauffman states, which has no cancellation on alternating diagrams and in which the signature of an alternating link appears. On a reduced alternating diagram the two ends of this state sum, for the total degree in $t_0$ and $t_1$, are governed by the positive and by the negative crossings, which reproves the first two conjectures of Tait on alternating links.

The coupling $\BE_D$ itself has a place in this picture. It is not an evaluation of the crossing variables, because of the rule $x_jy_j\mapsto1$, but it can be carried inside the determinants: it turns the product of the two linear forms of a crossing into one bilinear form $\Psi_j$ on the lattice of flows, with coefficients in $\BZ[t_0^{\pm1},t_1^{\pm1}]$, and $\LG$ is the mixed discriminant of these $c$ forms, that is, the coefficient of $\lambda_1\cdots\lambda_c$ in $\det\bigl(\sum_j\lambda_j\Psi_j\bigr)$ (Theorem~\ref{thm.intro.flows}). These forms come from a local system of rank one. It lives neither on the Seifert surface nor on the complement of the link, but on the complement of the link and of the co-cores of the bands of $F_D$, the arcs which join the two strands at the crossings: its holonomy is $t_1$ or $t_0^{-1}$ around the boundaries of the disks of $F_D$, according to a coloring of the disks with two colors, and $(t_0t_1)^{\pm1}$ around the co-cores, so that it extends to the complement of the link at the free fermion point only. The two linear forms of a crossing are two lifts of the meridian of its co-core, and $\Psi_j$ is an interleaved product of these lifts (Theorem~\ref{thm.intro.cocores}).

What the deformed Seifert matrices prove is the genus bound for the \emph{canonical} surface of a diagram; the bound of \cite{Kohli-Tahar, LNVdV25} with the genus of the link, and the monicity theorem of \cite{LNVdV25}, are not reproved here, and we show on examples that for other Seifert surfaces the bound is a cancellation inside the coupling $\BE_D$ rather than a property of the two determinants (Section~\ref{ssec.beyond}). The determinant formula does lead to them by another route, modulo work of William Ballinger. Regrouped according to the superposition of the two states of a pair, the state sum of Theorem~\ref{thm.intro.kauffman} below is the case $N=2$ of a state model $Q_N$ which Ballinger has presented in several talks \cite{BallingerTalk}, and which was not known to be the Links--Gould invariant (Corollary~\ref{cor.ballinger}); and with his model Ballinger plans to prove -- in a unified and elegant way -- the genus bound for every Seifert surface, the multiplicativity of the top coefficient under arbitrary Murasugi sums and the monicity on fibered links. The Kinoshita--Terasaka knot shows that $\LG$ is not a function of the Seifert form, so no formula in terms of a Seifert matrix alone can exist (Remark~\ref{rem.KT}); the diagram enters through the crossing variables and the coupling, and this is where the quantum invariant differs from the classical ones.

\subsection*{Results}
The main results are stated here in the order of the text. The notation which a statement needs is introduced just before it, so that this part of the introduction can be read without the rest of the text, and each theorem collects statements which are proved, often in a more precise form, in the sections which are indicated before it.
\subsubsection*{The lattice of flows}
Let $D$ be a connected diagram with crossings $1,\dots,c$ and set of edges $E$, $|E|=2c$; at the crossing $j$ we write $o_j,u_j$ for the incoming and $o'_j,u'_j$ for the outgoing edges of the over-strand and of the under-strand (Figure~\ref{fig.crossing}). A \emph{flow} on $D$ is a vector $z\in\BZ^E$ such that, at every crossing, what flows in is what flows out: $z_{o_j}+z_{u_j}=z_{o'_j}+z_{u'_j}$. For example, every component of the link, every Seifert circle of the diagram and the boundary of every face define flows: the first two have the value $1$ on their own edges and $0$ on the other ones, and the boundary of a face has the value $1$ on the edges which have the face on their left, $-1$ on those which have it on their right and $0$ on the other ones. The flows which vanish on the cut edge $\hate$ form a lattice $\FF_{D,\hate}$ of rank $c$, the first homology group of the cut diagram. The \emph{crossing forms} of the two flavors are the linear forms
\[
   \pi_j(z)=x_j\,z_{u'_j}-z_{u_j}\qquad\text{and}\qquad\tpi_j(z)=y_j\,z_{o_j}-z_{o'_j}
\]
on the flows: the first one looks at the under-strand of the crossing $j$ and the second one at its over-strand. The first result is proved in Section~\ref{sec.flows} (Theorems~\ref{thm.flows}, \ref{thm.faceform} and~\ref{thm.mixed}, Proposition~\ref{prop.adapted}).

\begin{theorem}\label{thm.intro.flows}
Let $D$ be a connected diagram of a link $L$ with $\mu$ components and $\hate$ an edge of $D$. For every $\BZ$-basis $z_1,\dots,z_c$ of the lattice of flows $\FF_{D,\hate}$,
\[
   \det\bigl(\pi_j(z_i)\bigr)_{i,j}=\kappa\,P_{D,\hate}(x)\qquad\text{and}\qquad\det\bigl(\tpi_j(z_i)\bigr)_{i,j}=(-1)^{\mu-1}\kappa\,\tP_{D,\hate}(y),
\]
with the same sign $\kappa=\pm1$ in the two identities. Consequently, for every basis,
\[
   \LG(L;t_0,t_1)=\BE_D\Bigl[\det\bigl(\pi_j(z_i)\bigr)\det\bigl(\tpi_j(z_i)\bigr)\Bigr]=\ve_D\,\BE_D\Bigl[\det\Bigl(\sum_{j=1}^c\ve_j\,\pi_j(z_i)\,\tpi_j(z_k)\Bigr)_{i,k}\Bigr],
\]
where $\ve_D=\prod_j\ve_j$ is the product of the signs of the crossings. Moreover, let $\Psi_j$ be the bilinear form on the flows, with values in $\BZ[t_0^{\pm1},t_1^{\pm1}]$, which is the image of $\pi_j\otimes\tpi_j$ under $\BE_D$:
\[
   \Psi_j(z,z')=z_{u'_j}z'_{o_j}-t_0^{\ve_j}\,z_{u'_j}z'_{o'_j}-t_1^{\ve_j}\,z_{u_j}z'_{o_j}+z_{u_j}z'_{o'_j} .
\]
Then, for every basis, $\LG(L;t_0,t_1)$ is the coefficient of $\lambda_1\cdots\lambda_c$ in the polynomial
\[
   \det\Bigl(\sum_{j=1}^c\lambda_j\,\Psi_j(z_i,z_k)\Bigr)_{i,k=1,\dots,c}
\]
in the indeterminates $\lambda_1,\dots,\lambda_c$.
\end{theorem}

In the matrix $\bigl(\pi_j(z_i)\bigr)$ the variable $x_j$ occurs in the column $j$ only, where it multiplies a column of integers; a change of basis multiplies the two matrices on the left by the same matrix of $\GL_c(\BZ)$. The last statement carries the coupling inside the determinant: the coefficient in question is the \emph{mixed discriminant} of the forms $\Psi_1,\dots,\Psi_c$, the symmetric multilinear function of $c$ forms whose value on $c$ equal forms is $c!$ times their determinant (Section~\ref{ssec.mixed}). The polynomial $\det\bigl(\sum_j\lambda_j\Psi_j\bigr)$ itself is not an invariant of the link: it changes under a Reidemeister move of type I, and it depends on the edge $\hate$ (Proposition~\ref{prop.kink} and Example~\ref{ex.pencilunknot}). For the trefoil considered above it is $\lambda_1\lambda_2\lambda_3\,\LG(3_1)+(1-t_0t_1)\,\lambda_2^2\,\bigl((t_0-1)\lambda_1+(t_1-1)\lambda_3\bigr)$: the Links--Gould polynomial is the coefficient of $\lambda_1\lambda_2\lambda_3$, and the other terms vanish at the free fermion point (Example~\ref{ex.trefoilmixed}). Three kinds of bases are used in this monograph.

\emph{The faces.} The boundaries of the faces of $D$ which are not adjacent to $\hate$ form a basis of $\FF_{D,\hate}$. In this basis the two matrices are Alexander's matrix of the diagram with the crossing variables inserted, in two versions: the row of the crossing $j$ is the relation of its under-strand for the first flavor and that of its over-strand for the second one (Theorem~\ref{thm.faceform}).

\emph{The Seifert circles.} The Seifert circles other than the one which contains $\hate$, together with flows which represent a basis of the first homology of the canonical Seifert surface, that is, a basis of the cycle space of the \emph{Seifert graph} $\Gamma_D$, whose vertices are the $s$ Seifert circles and whose edges are the crossings, form a basis (Lemma~\ref{lem.lifts}). It is in this basis that we prove the theorem, by eliminating the generators of the edges along the Seifert circles, and it leads to the deformed Seifert matrices below.

\emph{The smoothings.} A \emph{smoothing} $\vs\in\{0,1\}^c$ smooths every crossing in one of the two possible ways, along the orientation ($\vs_j=0$) or not ($\vs_j=1$), and turns the diagram into a family of circles. For every smoothing there are bases formed by the circles and by flows which represent a basis of the quotient of the lattice by the circles (Proposition~\ref{prop.adapted}). The oriented smoothing gives the Seifert circles, and the two checkerboard smoothings give back the faces.

\emph{The trefoil.} The diagram of the trefoil considered above has five faces. The boundaries of the three faces which are not adjacent to the edge $\hate$ form a basis $z_1,z_2,z_3$ of the lattice of flows, and in this basis the two matrices of Theorem~\ref{thm.intro.flows} are
\[
   \bigl(\pi_j(z_i)\bigr)=\begin{pmatrix}-x_1&-1&0\\0&-x_2&-1\\1&1&1\end{pmatrix}\qquad\text{and}\qquad
   \bigl(\tpi_j(z_i)\bigr)=\begin{pmatrix}-1&-y_2&0\\0&-1&-y_3\\1&1&1\end{pmatrix}
\]
(Example~\ref{ex.trefoilflows}). Their determinants are the two edge polynomials $1-x_1+x_1x_2$ and $1-y_3+y_2y_3$. The variables $x_3$ and $y_1$ do not occur: the edge $\hate$ leaves the crossing $3$ as its under-strand and enters the crossing $1$ as its over-strand, and the flows of the lattice vanish on it.

\subsubsection*{State surfaces and the localization of linking numbers}
The quotient of the lattice of flows by the circles of a smoothing $\vs$ is the first homology group of a surface. The \emph{state surface} $F_\vs$ consists of a disk for each circle of $\vs$ and of a half-twisted band at each crossing; its boundary is the link, and it is the canonical Seifert surface $F_D$ for the oriented smoothing, a checkerboard surface for the checkerboard smoothings, and in general a nonorientable surface. We place all the disks below the plane of the diagram, the disk of an inner circle above the disk of every circle which encloses it (Figure~\ref{fig.statesurface}). The flows are then $1$-cycles which lie on all these surfaces at once: a flow follows the edges of the link and passes from one strand of a crossing to the other one along the co-core of its band. The second result is proved in Sections~\ref{sec.statesurfaces} and~\ref{sec.seifert} (Theorem~\ref{thm.localization}, Corollaries~\ref{cor.GLform} and~\ref{cor.seifertform}, Proposition~\ref{prop.otherside} and Theorem~\ref{thm.deformedseifert}).

\begin{theorem}\label{thm.intro.linking}
Let $D$ be a connected diagram, $\vs$ a smoothing of $D$ and $F_\vs$ its state surface. For a flow $z$ put $U_j(z)=z_{u'_j}+z_{u_j}$, and let $w^\vs_j(z)$ be the number of times the cycle $z$ runs through the band of the crossing $j$, which is $z_{o_j}-z_{u'_j}$ if $\vs_j=0$ and $-z_{o'_j}-z_{u'_j}$ if $\vs_j=1$.
\begin{enumerate}\renewcommand{\theenumi}{\roman{enumi}}
\item $H_1(F_\vs;\BZ)$ is the quotient of the lattice of all the flows of $D$ by the circles of $\vs$. For two flows $z,z'$, the form of Gordon and Litherland of $F_\vs$, that is, the linking number of the cycle $z$ with a curve on the surface in the class of $z'$, pushed off the surface to both of its sides, is
\[
   \CG_{F_\vs}\bigl([z],[z']\bigr)=\sum_{j=1}^c\theta_j\,U_j(z)\,w^\vs_j(z'),\qquad\theta_j=\ve_j(-1)^{\vs_j} .
\]
\item For the canonical Seifert surface $F_D$ of $D$, oriented by the link, let $e^+_j$ and $e^-_j$ be the edges of the under-strand of the crossing $j$ which lie on the left and on the right of its over-strand. The Seifert form $V_F$ of $F_D$ and its transpose are given, on the classes of two flows $z,z'$, by
\[
   V_F\bigl([z],[z']\bigr)=\sum_{j=1}^c\ve_j\,z_{e^+_j}\,\bigl(z'_{o_j}-z'_{u'_j}\bigr),\qquad V_F\bigl([z'],[z]\bigr)=\sum_{j=1}^c\ve_j\,z_{e^-_j}\,\bigl(z'_{o_j}-z'_{u'_j}\bigr).
\]
\item Consequently, with $V_j(z)=\ve_jz_{e^+_j}$ and $V'_j(z)=\ve_jz_{e^-_j}$, the contributions of the band $j$ to the Seifert form and to its transpose, the crossing form of the first flavor is
\[
   \pi_j=x_j^{\ind{\ve_j=1}}\,V_j-x_j^{\ind{\ve_j=-1}}\,V'_j ,
\]
which is $t^{-\ind{\ve_j=-1}}\bigl(t\,V_j-V'_j\bigr)$ at the Alexander point $x_j=t^{\ve_j}$. The crossing form $\tpi_j$ of the second flavor is given in the same way by the edges of the over-strand and by the Seifert form of the canonical surface $F^\uparrow_D$ placed above the plane of the diagram, in which the disk of an inner circle lies below the disks of the circles which enclose it.
\end{enumerate}
\end{theorem}

So the matrix of Theorem~\ref{thm.intro.flows} is obtained from $tV_F-V_F^{\T}$ by \emph{not} summing over the bands through which the second cycle passes, so that there is one column for each band instead of one column for each cycle, and by replacing $t$ by the variable of the band. The theorem holds for every diagram, whatever the orientations and the nesting of its Seifert circles; when a Seifert circle has crossings on both sides, the two surfaces $F_D$ and $F^\uparrow_D$ are the two Murasugi sums of the parts of the diagram inside and outside the circle, and their Seifert forms are different in general; we compare them in Corollary~\ref{cor.twoseifert}.

The proof rests on one observation: when the surface lies below the plane of the diagram, the faces of the diagram, which stay in this plane, are disks which span their boundaries and meet the surface nowhere else; so the linking number of a flow with a curve is read off the passages of the curve through the faces (Lemma~\ref{lem.spanning} and Figure~\ref{fig.unrolled}).

\subsubsection*{Deformed Seifert matrices}
In a basis of the second kind the rows of the Seifert circles can be eliminated, and what remains is a square matrix of size $\nu=c-s+1=2g(F_D)+\mu-1$ on the homology of the canonical surface. Fix a basis of the cycle space of $\Gamma_D$, encoded in a $c\times\nu$ integer matrix $Y$. The third result is proved in Sections~\ref{sec.flows} and~\ref{sec.seifert} (Propositions~\ref{prop.seifert} and~\ref{prop.seifertb}, Corollaries~\ref{cor.seifert} and~\ref{cor.seifertb}, Theorem~\ref{thm.alexpoint}, Corollary~\ref{cor.classical} and Proposition~\ref{prop.seifertalex}).

\begin{theorem}\label{thm.intro.seifert}
There are $c\times c$ matrices $K(x)$ and $K^b(y)$, with entries in $\BZ[x_1,\dots,x_c]$ and in $\BZ[y_1,\dots,y_c]$ respectively, which are read off from the order in which the crossings are met along the Seifert circles, and, for every spanning tree $\CT$ of $\Gamma_D$, $c\times\nu$ matrices $Y(x)$ and $Y(y)$ with entries in $\BZ[x_j^{\pm1}]$ and $\BZ[y_j^{\pm1}]$, with $Y(1)=Y$, such that
\[
   P_{D,\hate}(x)=x^{\CT}\det\bigl(Y^{\T}K(x)Y(x)\bigr),\qquad \tP_{D,\hate}(y)=y^{\CT}\det\bigl(Y(y)^{\T}K^b(y)Y\bigr),
\]
with no sign, where $x^{\CT}$ and $y^{\CT}$ are monomials in the variables of the edges of $\CT$. Consequently
\[
   \LG(L;t_0,t_1)=(-1)^{\mu-1}\,\BE_D\Bigl[x^{\CT}y^{\CT}\det\bigl(M_{D,\hate}(x)\bigr)\det\bigl(M^b_{D,\hate}(y)\bigr)\Bigr]
\]
with the $\nu\times\nu$ matrices $M_{D,\hate}(x)=Y^{\T}K(x)Y(x)$ and $M^b_{D,\hate}(y)=Y(y)^{\T}K^b(y)Y$. The rows of $Y^{\T}K(x)$ are the crossing forms of the flows which represent the basis, and the columns of $Y(x)$ form a basis of a deformation of the cycle space. At the Alexander point $x_j=t^{\ve_j}$, for the basis of the fundamental cycles of $\CT$ and for every diagram,
\[
   M_{D,\hate}(t^{\ve})=\bigl(tV_F-V_F^{\T}\bigr)\,\mathsf Z'_t\qquad\text{and}\qquad M^b_{D,\hate}(t^{\ve})=\mathsf Z'_t\,\bigl(tV_{F^\uparrow}-V_{F^\uparrow}^{\T}\bigr),
\]
where $V_F$ and $V_{F^\uparrow}$ are the Seifert matrices of the canonical surfaces below and above the plane and $\mathsf Z'_t$ is a diagonal matrix of units. The matrix $M_{D,\hate}(t^{\ve})$ presents the same $\BZ[t^{\pm1}]$-module as the cut Alexander matrix of $D$; its determinant is $\Delta_L(t)$ up to a unit, and its corank over $\BQ(t)$ is the rank of the Alexander module of the infinite cyclic cover of the link exterior classified by the total linking number.
\end{theorem}

The matrices $K(x)$ and $K^b(y)$ are built from \emph{jump relations}: walking along a Seifert circle, the generator attached to an edge jumps, at each crossing, by a multiple of the generator attached to the band of that crossing, the multiple being $x_j$ or $-1$ according to the side of the crossing on which the circle passes (Figure~\ref{fig.jump}). Since the theorem has no sign, it gives the classical formula with its sign, $\det(tV-V^{\T})=(-1)^{\mu-1}t^{n/2}\Delta_L(t)$ for every Seifert matrix $V$ of size $n$ and the Conway-normalized Alexander polynomial (Corollary~\ref{cor.classical}). For the closed $2$-braids $M_{D,\hate}(x)=\sum_jx_jV_j-V_F^{\T}$ with $V_F=\sum_jV_j$ the Seifert matrix and $V_j$ of rank one (Proposition~\ref{prop.twobraid}); the trefoil and the figure-eight knot are worked out for all cut edges in Examples~\ref{ex.trefoil} and~\ref{ex.figeight}.

\emph{The trefoil.} For the diagram of the trefoil considered above, $\nu=2$, and in a suitable basis the Seifert matrix of the canonical surface is $V_F=\bigl(\begin{smallmatrix}-1&-1\\0&-1\end{smallmatrix}\bigr)$. The deformed Seifert matrix of the first flavor is
\[
   M_{D,\hate}(x)=\begin{pmatrix}1-x_1&-x_1\\1+x_2-x_1&1-x_1\end{pmatrix}=x_1\begin{pmatrix}-1&-1\\-1&-1\end{pmatrix}+x_2\begin{pmatrix}0&0\\1&0\end{pmatrix}-V_F^{\T} .
\]
The two matrices which are multiplied by $x_1$ and by $x_2$ are the contributions of the bands of the crossings $1$ and $2$ to the Seifert form, and their sum is $V_F$; for this edge $\hate$ the band of the crossing $3$ does not contribute. So $M_{D,\hate}(x)$ is the matrix $tV_F-V_F^{\T}$ in which $t$ has been replaced, band by band, by the variable of the crossing. Its determinant is the edge polynomial $1-x_1+x_1x_2$, and for $x_1=x_2=t$ the matrix is $tV_F-V_F^{\T}$, of determinant $1-t+t^2$ (Examples~\ref{ex.trefoilseifert} and~\ref{ex.trefoil}).

\emph{Cochains, not classes.} The crossings are the bands of $F_D$, and the bands correspond to their co-cores, which are the relative $1$-cells of the surface. So the family $x=(x_1,\dots,x_c)$ of the crossing variables is a multiplicative relative $1$-\emph{cochain} of $(F_D,\partial F_D)$, and the ring of the polynomials in the crossing variables is the group ring of the group $\BZ^c$ of the relative cochains. The relative cohomology classes form the subgroup of the cocycles, which is the cycle space, of rank $\nu<c$; and one may ask whether the determinant formula can be written in terms of the classes alone. It cannot (Section~\ref{ssec.cochains}): the monomials of the edge polynomials do not lie in one coset of the cycle space, already for the trefoil; the Alexander point $x=t^{\ve}$ is not a cocycle, for any diagram; and as a local system on the surface, that is, modulo the gauge transformations of the disks, the Alexander point is trivial, while the edge polynomials are not gauge invariant.

The reason is the one given by Theorem~\ref{thm.intro.linking}: the variable $t$ of the Alexander polynomial does not live on the Seifert surface, it distinguishes its two sides, and the variable $x_j$ does this at the band $j$ only, independently of the other bands. The variable $t$ itself is not seen by the classes: the sign $\ve_j$ is the difference of the levels of the two Seifert circles joined by the band $j$ (Lemma~\ref{lem.level}), so that the Alexander point sends the monomial $x^n$ of every class $n$ to $1$. The values $t_0^{\ve_j}$ and $t_1^{\ve_j}$ which the coupling $\BE_D$ gives to the crossing variables are nevertheless holonomies of a local system of rank one; but it lives on the complement in $\BR^3$ of the link and of the co-cores of the bands, and not on the surface (Theorem~\ref{thm.intro.cocores}).

Two consequences are used later: the sums of the extremal coefficients of both edge polynomials are $\det V_F$ and $(-1)^\nu\det V_F$ for every diagram (Corollary~\ref{cor.extremesums}), and inserting $\sigma_i\sigma_i^{-1}$ into a braid word, which attaches a hollow handle to $F_D$ (Figure~\ref{fig.hollowhandle}), changes $P_{D,\hate}$ by a multiple $(x_ax_b-1)R(x)$ of $x_ax_b-1$, where $a$ and $b$ are the two new crossings and $R(x)$, the coefficient of $x_ax_b$ in the new edge polynomial, does not contain their variables; this is a trace of the handle which only the coupling $\BE_D$ removes (Proposition~\ref{prop.hollowhandle}).

\subsubsection*{Degrees, Murasugi sums and alternating links}
The fourth result is the degree bound and the extreme coefficients. For a set $S$ of crossings we write $x^S=\prod_{j\in S}x_j$, $y^S=\prod_{j\in S}y_j$ and $\ve(S)=\sum_{j\in S}\ve_j$. We grade the crossing variables by the signs of the crossings, $\deg_\ve x_j=\deg_\ve y_j=\ve_j$, so that $\BE_D[x^Sy^T]$ is a monomial of degree $\ve(S)-\ve(T)$ for the degree $\deg(t_0^at_1^b)=a-b$ of \cite{Kohli-Tahar}, and we write $g(L)$ for the genus of $L$ and $\mathrm{top}(L,q)$ for the coefficient of $p^{2g(L)+\mu-1}$ in $\LG(L)$, the top coefficient of \cite{LNVdV25}, in the variables $t_0=q/p$, $t_1=pq$ of that paper. The following is proved in Section~\ref{sec.degree} (Theorems~\ref{thm.degrees}, \ref{thm.canonical}, \ref{thm.tight}, \ref{thm.murasugi} and~\ref{thm.alternating}).

\begin{theorem}\label{thm.intro.degrees}
Let $D$ be a connected diagram of $L$, $\hate$ an edge of $D$ and $\nu=c-s+1=2g(F_D)+\mu-1$.
\begin{enumerate}\renewcommand{\theenumi}{\roman{enumi}}
\item Up to shifts of the rows and the columns, which are the same for the two flavors, every entry of $M_{D,\hate}(x)$ and of $M^b_{D,\hate}(y)$ is the sum of a homogeneous part of $\ve$-degree $0$ and of a homogeneous part of $\ve$-degree $1$. Consequently the monomials of $P_{D,\hate}$ and of $\tP_{D,\hate}$ have $\ve$-degrees in one and the same interval of length $\nu$, and the degrees of $\LG(L)$ lie in $[-\nu,\nu]$.
\item If $\Delta_L$ has span $\nu$ (for instance if $D$ is alternating or homogeneous), the degrees of $\LG(L)$ fill $[-\nu,\nu]$, so that $\mathrm{span}\,\LG(L)=4g(L)+2\mu-2$, and the coefficients of the two extreme degrees are sums over the monomials of $P_{D,\hate}$ and $\tP_{D,\hate}$ of extreme $\ve$-degree.
\item If a Seifert circle $C$ of $D$ has crossings on both sides, $F_D$ is the Murasugi sum along the disk of $C$ of the canonical surfaces of two diagrams $D_{\rm in}$ and $D_{\rm out}$ (Figure~\ref{fig.murasugi}), and the coefficient $\tau_D(q)$ of $p^{\nu}$ in $\LG(L)$, which is $\mathrm{top}(L,q)$ when $F_D$ has minimal genus, multiplies: $\tau_D=\tau_{D_{\rm in}}\tau_{D_{\rm out}}$. For the closure of a homogeneous braid in which $\sigma_i$ occurs $a_i$ times, with the exponent $\ve_i$, this gives $\mathrm{top}(L,q)=q^{\sum_i\ve_i(a_i-1)}$.
\item If $D$ is alternating, the sums of (ii) run over the extremal spanning trees of \cite{HK}, have nonnegative coefficients, and reduce to a single monomial $q^{\ve(S_0\triangle T_0)}$ exactly when $L$ is fibered, $S_0$ and $T_0$ being the crossing sets of the two extremal trees (Figure~\ref{fig.trees}). In that case the Hopf invariant $\lambda(\xi_F)$ of the plane field of the fiber $F$, which by \cite{LNVdV25} determines this monomial as $q^{\nu-2\lambda(\xi_F)}$, is read off the diagram:
\[
   \lambda(\xi_F)=\tfrac12\bigl(\nu-\ve(S_0\triangle T_0)\bigr).
\]
\end{enumerate}
\end{theorem}

The proof of (i) is the Seifert-matrix proof of $\deg\Delta_L\le c-s+1$ with the crossing variables carried along; the only geometric input is that every cycle of the Seifert graph lies in one region of the complement of the Seifert circles, whose boundary circles have two adjacent levels (Figure~\ref{fig.regions}). Part (i) minimized over the diagrams of $L$ gives $\deg_p\LG(L)\le2g_c(L)+\mu-1$ with $g_c(L)$ the canonical genus, which can be larger than $g(L)$; the hollow-handle law is the elementary case of the cancellation which, inside $\BE_D$, produces the bound with $g(L)$.

Part (iii) is a case of the theorem of Anghel, L\'opez Neumann and van der Veen on arbitrary Murasugi sums of minimal genus surfaces \cite{ALvdV26}, which they prove with the quantum supergroup; here it is the block structure of the extremal parts of the deformed Seifert matrices. What (iv) adds to the fibered case is the formula for the Hopf invariant: for a fibered alternating link, a quantity of the fibration is read off any alternating diagram as half the difference between $\nu$ and the signed number of crossings in which the two extremal spanning trees differ (Example~\ref{ex.alttrees} checks it against the values of \cite{LNVdV25} for $3_1$, $4_1$, $6_2$ and $6_3$). Theorem~\ref{thm.intro.kauffman} below identifies this signed number: it is minus the signature of the link.

\emph{The trefoil.} For the diagram of the trefoil considered above, $\nu=2$. In the variables $p$ and $q$, the polynomial which was computed at the beginning of this introduction is
\[
   \LG(3_1)=q^2\,p^{2}-(q+q^3)\,p+(1+2q^2)-(q+q^3)\,p^{-1}+q^2\,p^{-2} .
\]
Its degrees fill the interval $[-2,2]$, as (ii) says, and its top coefficient is $\mathrm{top}(3_1,q)=q^2$. This is the monomial which (iii) gives for the closure of the braid $\sigma_1^3$; and since it is $q^{\nu-2\lambda(\xi_F)}$, the Hopf invariant of the fiber of the trefoil is $\lambda(\xi_F)=0$.

\subsubsection*{The coupling: a local system on the complement of the co-cores}
In the space $\BR^3$ in which the state surfaces are placed, let $G$ be the union of the link and of the co-cores $\alpha_1,\dots,\alpha_c$ of the bands, the vertical segments which join the under-strand to the over-strand at the crossings. The lattice of all the flows of $D$ is the group of the $1$-cycles of the graph $G$, and its dual is carried by the complement $X_G=\BR^3\sm G$: for a local system of free modules of rank one on $X_G$, over a commutative ring $\kk$, the linear forms on the flows with values in $\kk$ are the lifted $1$-chains of $X_G$ which begin and end over a base point, the pairing being a twisted linking number. There are two base points, one above and one below the diagram, hence two such identifications. A local system of rank one on $X_G$ is given by a unit $\vth_f$ of $\kk$ for each face $f$ of $D$, up to a common factor, and its holonomy around an edge $e$ is $\mhol_e=\vth_{f_{\rm L}(e)}/\vth_{f_{\rm R}(e)}$, where $f_{\rm L}(e)$ and $f_{\rm R}(e)$ are the faces on the left and on the right of $e$. The fifth result is proved in Section~\ref{sec.cocores} (Propositions~\ref{prop.readings} and~\ref{prop.interleaved}, Corollary~\ref{cor.anypotential}, Theorem~\ref{thm.cocores} and Proposition~\ref{prop.paths}).

\begin{theorem}\label{thm.intro.cocores}
Let $D$ be a connected diagram of $L$ and $\hate$ an edge of $D$.
\begin{enumerate}\renewcommand{\theenumi}{\roman{enumi}}
\item For every local system of rank one on $X_G$ and every crossing $j$, the lifts of the meridian of the co-core $\alpha_j$ are the crossing forms, up to units. Lifted from a point below the over-strand, the meridian is the form $\pi_j$, in which $x_j$ is $\mhol_e^{-\ve_j}$ for an edge $e$ of the over-strand; lifted from a point above the under-strand, it is the form $\tpi_j$, in which $y_j$ is $\mhol_e^{\ve_j}$ for an edge $e$ of the under-strand.
\item Color the Seifert circles of $D$ with two colors, in such a way that every crossing joins two circles of different colors, and consider the local system with $\mhol_e=t_1$ on the edges of the first color and $\mhol_e=t_0^{-1}$ on the edges of the second color. For every crossing $j$ there is a combination $\CZ_j$ of products of two lifts of the meridian of $\alpha_j$, in which the starting point of the first lift depends on the position of the second one, whose bilinear form is the form $\Psi_j$ of Theorem~\ref{thm.intro.flows}, up to a unit; and $\LG(L;t_0,t_1)$ is the mixed discriminant of the chains $\CZ_1,\dots,\CZ_c$, up to a unit $t_0^at_1^b$.
\item The holonomy of this local system around every co-core is $(t_0t_1)^{\pm1}$, and its holonomy around the band of every crossing is trivial. In general, the local systems of rank one on $X_G$ correspond in this way to the linear maps $x_j\mapsto\xs_j$, $y_j\mapsto\ys_j$, $x_jy_j\mapsto\sv_j$ applied to $P_{D,\hate}\tP_{D,\hate}$; such a map is an evaluation at the crossing $j$, that is, $\sv_j=\xs_j\ys_j$, if and only if the holonomy around the co-core $\alpha_j$ is trivial, and $\sv_j=1$ if and only if the holonomy around the band of the crossing $j$ is trivial.
\end{enumerate}
\end{theorem}

So the two features of the coupling $\BE_D$ have a meaning in the complement of the co-cores: the two variables are, up to an inversion, the holonomies around the Seifert circles of the two colors, and the rule $x_jy_j\mapsto1$ says that the local system is trivial around the bands. It is not trivial around the co-cores, unless $t_0t_1=1$, and this is why $\BE_D$ is not an evaluation; the local systems which extend to the complement of the link give evaluations, hence products of two determinants of Alexander type. Theorem~\ref{thm.intro.cocores} is a description of the determinant formula, and not a definition of $\LG$ from a Seifert surface: the lifts are taken with respect to paths which join the meridians of the co-cores to the two base points, on the two sides of the plane of the diagram, and the result depends on these paths, through a term which is divisible by $1-t_0t_1$ (Section~\ref{ssec.vertical}). For the trefoil, the local system and the chains $\CZ_j$ are written out in Example~\ref{ex.trefoilcocores} (Figure~\ref{fig.trefoilcocores}).

\subsubsection*{The faces: Kashaev's matrix}
Kashaev \cite{Kashaev21} attached to a diagram a \emph{symmetric} matrix $\tau_D$ indexed by its faces, with entries in $\BZ[2\zeta]$, $2\zeta=t^{1/2}+t^{-1/2}$, and conjectured that its signature computes the Levine--Tristram signature; the conjecture was proved by Liu \cite{Liu25}, and Cimasoni--Ferretti \cite{CimasoniFerretti24} and Liu showed that the determinant of $\tau_D$ with two adjacent faces removed is $\Delta_L(t)^2$ up to a unit. In the lattice of flows the matrix of Kashaev is the coupled form of Theorem~\ref{thm.intro.flows}, written in the basis of the faces. Section~\ref{sec.faces} proves the following (Theorem~\ref{thm.kashaev}, Corollaries~\ref{cor.kashaev} and~\ref{cor.deformedkashaev}).

\begin{theorem}\label{thm.intro.faces}
Let $\CF$ be the set of the faces of $D$, and let $Q_D(x)$ and $\tQ_D(y)$ be the $c\times\CF$ matrices whose entries are the values $\pi_j(\partial f)$ and $\tpi_j(\partial f)$ of the crossing forms on the boundaries of the faces.
\begin{enumerate}\renewcommand{\theenumi}{\roman{enumi}}
\item Kashaev's matrix is the Gram matrix of the two flavors at the free fermion point: $\tau_D(t)=\WW_t^{-1}\,Q_D(t^{\ve})^{\T}\diag(\ve_1,\dots,\ve_c)\,\tQ_D(t^{-\ve})\,\WW_t$, where $\WW_t$ is the diagonal matrix with entries $(-t^{1/2})^{w(f)}$, $w(f)$ being the winding number of $D$ around the face $f$. Consequently $\det\tau_D(t)$ with two adjacent faces removed is $\ve_D\,\Delta_L(t)^2$, where $\Delta_L$ is the Conway-normalized Alexander polynomial, and the kernel of $\tau_D(t)$ is spanned by the two vectors $\bigl((-t^{\pm1/2})^{w(f)}\bigr)_f$ when $\Delta_L\ne0$.
\item The \emph{deformed Kashaev matrix} $\tau_D(x,y)=Q_D(x)^{\T}\diag(\ve_1,\dots,\ve_c)\,\tQ_D(y)$, a sum over the crossings of matrices of rank one which are bilinear in the two variables of the crossing, satisfies
\[
   \LG(L;t_0,t_1)=\ve_D\,\BE_D\bigl[\det\tau_{D,\hate}(x,y)\bigr]\qquad\text{and}\qquad\tau_D(t^{\ve},t^{-\ve})=\WW_t\,\tau_D(t)\,\WW_t^{-1}
\]
for every edge $\hate$, where $\tau_{D,\hate}$ is the matrix with the rows and the columns of the two faces adjacent to $\hate$ removed.
\end{enumerate}
\end{theorem}

The deformation breaks the symmetry of Kashaev's matrix: the rows of $\tau_D(x,y)$ belong to the first flavor and its columns to the second one, and $\BE_D$ couples them crossing by crossing, exactly as it couples the deformed Seifert matrices: it replaces the block of rank one of every crossing by the matrix of the form $\Psi_j$ of Theorem~\ref{thm.intro.flows}, and the determinant by the mixed discriminant. The Hermitian structure of $\tau_D$, which is the source of the signature, has no counterpart for $\LG$ (Remark~\ref{rem.deformedkashaev}). The deformed Kashaev matrix of the trefoil is written out in Example~\ref{ex.trefoilkashaev}.

\subsubsection*{The point $x=-1$: Goeritz matrices and the forms of Gordon and Litherland}
At the Alexander point the crossing variables depend on the orientation of the link, except for $t=-1$, where they are all equal to $-1$. There the crossing forms are $-U_j$ and $-O_j$, with $U_j(z)=z_{u'_j}+z_{u_j}$ and $O_j(z)=z_{o_j}+z_{o'_j}$, the oriented smoothing is no longer distinguished, and the lattice of flows contains the classical theory of the spanning surfaces of the link. For a smoothing $\vs$ let $\bar\vs$ be the complementary smoothing, let $\CG_F$ denote the form of Gordon and Litherland of a surface $F$, and put
\[
   \Goe_\vs=-\sum_{j=1}^c\theta_j\,w^\vs_j\otimes w^\vs_j,\qquad\Stk_\vs=-\sum_{j=1}^c\theta_j\,w^{\bar\vs}_j\otimes w^\vs_j,\qquad\Kash=\sum_{j=1}^c\ve_j\,U_j\otimes O_j ,
\]
three bilinear forms on the lattice of flows; $\Kash$ is the coupled form of Theorem~\ref{thm.intro.flows} at $x=y=-1$. We write $\sigma_L$ for the signature of $L$, with the sign convention for which the right-handed trefoil has signature $-2$. Section~\ref{sec.minusone} proves the following (Proposition~\ref{prop.formsminusone}, Theorem~\ref{thm.goeritz}, Proposition~\ref{prop.J}, Theorems~\ref{thm.determinant} and~\ref{thm.signature} and their corollaries).

\begin{theorem}\label{thm.intro.minusone}
Let $D$ be a connected diagram of $L$ and $\vs$ a smoothing of $D$, with state surfaces $F_\vs$ below and $F^\uparrow_\vs$ above the plane of the diagram.
\begin{enumerate}\renewcommand{\theenumi}{\roman{enumi}}
\item The three forms are symmetric, and $\CG_{F_\vs}=\Goe_\vs+\Stk_\vs$ and $\CG_{F^\uparrow_\vs}=\Goe_\vs-\Stk_\vs$ on the homology of the surfaces. The form $\Goe_\vs$ is a sum of Goeritz forms, one for each region of the complement of the circles of $\vs$, and $\Stk_\vs$ couples the regions which are separated by one circle; the two forms of Gordon and Litherland are exchanged by the involution which is $(-1)^d$ on the homology of the regions of depth $d$, and they coincide when no circle has crossings on both sides.
\item $\Kash=\Goe_\vs+\Goe_{\bar\vs}=\CG_{F_\vs}+\CG_{F_{\bar\vs}}$ for every smoothing. In the basis of the faces, colored in a checkerboard manner, $\Kash$ is the direct sum of the two Goeritz matrices of the diagram, and Kashaev's matrix at $t=-1$ is congruent to this direct sum by a diagonal matrix of signs, for every diagram.
\item The Gram matrices of the forms $\CG_{F_\vs}$ and $\CG_{F^\uparrow_\vs}$, for all the smoothings, present the same abelian group: the quotient of $\BZ^c$ by the image of the lattice of flows under $(U_1,\dots,U_c)$, which is the first homology group of the double covering of $S^3$ branched along $L$. In particular their determinants are $\pm P_{D,\hate}(-1,\dots,-1)$, of absolute value $|\Delta_L(-1)|$.
\item For every smoothing,
\[
   \sign\CG_{F_\vs}=\sign\CG_{F^\uparrow_\vs}=\sigma_L+\sum_{j:\,\vs_j=1}\ve_j .
\]
Consequently $\sign G_v(D)-\mu_v(D)=\sigma_L$ for the Goeritz matrix $G_v(D)$ of each of the two colors, with its correction term $\mu_v(D)$, and $\sign\Kash=2\sigma_L+\mathrm{wr}(D)$, where $\mathrm{wr}(D)$ is the writhe of the diagram.
\end{enumerate}
\end{theorem}

Part (iv) is the formula of Gordon and Litherland for the state surfaces, and its two consequences are their formula for the Goeritz matrices \cite{GordonLitherland} and the theorem of Cimasoni and Ferretti on the signature of Kashaev's matrix at $t=-1$ \cite{CimasoniFerretti24}, which is the case $\omega=-1$ of Kashaev's conjecture. The proof given here uses neither four-dimensional topology nor special diagrams: changing the smoothing of one crossing adds to the form the square $\ve_j\,U_j\otimes U_j$ of a linear form, and by (iii) and the Jordan curve theorem it changes the rank by exactly one, so that the signature changes by $\ve_j$; and for the oriented smoothing the form is $V_F+V_F^{\T}$. For an alternating diagram the forms of the two checkerboard surfaces are definite, and (iv) becomes Traczyk's formula for the signature (Corollary~\ref{cor.altdefinite}); the image of the lattice of flows in $\BZ^c$ is then the sum of the lattice of the integral flows and of the lattice of the integral cuts of the Tait graph (Proposition~\ref{prop.underlattice}).

\emph{The trefoil.} For the diagram of the trefoil considered above, all the crossings are positive and $\sigma_L=-2$. The state surface of the oriented smoothing is the canonical Seifert surface, and its form $V_F+V_F^{\T}$ has the signature $-2$. By (iv), the forms of the other state surfaces have the signature $-1$, $0$ or $1$, according to whether one, two or three crossings are smoothed against the orientation. In the last case the surface is a checkerboard surface, a M\"obius band, and its form is the matrix $(3)$. The determinants of all these forms are $\pm3$, as (iii) says, since $|\Delta_L(-1)|=3$; and the signature of $\Kash$ is $-1=2\sigma_L+\mathrm{wr}(D)$ (Example~\ref{ex.trefoilminusone}).

\subsubsection*{The Kauffman model}
The last group of results is the counterpart of Kauffman's state model \cite{KauffmanFKT}, which is the expansion of Alexander's determinant over the bijections between the crossings and the faces. A \emph{Kauffman state} $\mathfrak s$ of $(D,\hate)$ chooses one of the four corners of every crossing, in such a way that every face other than the two faces adjacent to $\hate$ contains exactly one of the chosen corners (Figures~\ref{fig.corners} and~\ref{fig.trefoilstates}). Let $X(\mathfrak s)$ be the set of crossings at which the chosen corner is adjacent to the outgoing under-edge $u'_j$, $Y(\mathfrak s)$ the set of those at which it is adjacent to the incoming over-edge $o_j$, and $V(\mathfrak s)=X(\mathfrak s)\triangle Y(\mathfrak s)$, which is the set of crossings at which the chosen corner lies between the two outgoing or between the two incoming edges. We write $x^S=\prod_{j\in S}x_j$ and $\ve(S)=\sum_{j\in S}\ve_j$ for a set $S$ of crossings. We write $b(\mathfrak s)$ for the number of crossings at which the chosen corner lies between the two incoming edges; $(-1)^{b(\mathfrak s)}$ is the sign of the bijection which the state defines between the crossings and the faces, up to a sign which does not depend on the state. Section~\ref{sec.kauffman} proves the following (Proposition~\ref{prop.statesum}, Theorem~\ref{thm.kauffmanalt}, Corollary~\ref{cor.statestrees}, Proposition~\ref{prop.signature} and Corollary~\ref{cor.extremalstates}, and Proposition~\ref{prop.signs} for the signs in (i)).

\begin{theorem}\label{thm.intro.kauffman}
Let $D$ be a connected diagram of $L$ and $\hate$ an edge of $D$. In the following sums $\mathfrak s$ and $\mathfrak s'$ run over the Kauffman states of $(D,\hate)$.
\begin{enumerate}\renewcommand{\theenumi}{\roman{enumi}}
\item The edge polynomials are
\[
   P_{D,\hate}(x)=(-1)^{\mu-1}\sum_{\mathfrak s}(-1)^{b(\mathfrak s)}\,x^{X(\mathfrak s)},\qquad \tP_{D,\hate}(y)=\sum_{\mathfrak s}(-1)^{b(\mathfrak s)}\,y^{Y(\mathfrak s)};
\]
the number $m_{D,\hate}=\ve(X(\mathfrak s))+\ve(Y(\mathfrak s))$ is the same for all states; and
\[
   \LG(L;t_0,t_1)=\sum_{(\mathfrak s,\mathfrak s')}(-1)^{b(\mathfrak s)+b(\mathfrak s')}\prod_{j\in X(\mathfrak s)\sm Y(\mathfrak s')}t_0^{\ve_j}\ \prod_{j\in Y(\mathfrak s')\sm X(\mathfrak s)}t_1^{\ve_j},
\]
the sum running over the ordered pairs of states.
\item If $D$ is alternating, then the sign $\sigma_{D,\hate}=(-1)^{|X(\mathfrak s)\cap Y(\mathfrak s)|}$ does not depend on the state, the edge polynomials are
\[
   P_{D,\hate}(x)=\sigma_{D,\hate}\sum_{\mathfrak s}(-1)^{|X(\mathfrak s)|}\,x^{X(\mathfrak s)},\qquad \tP_{D,\hate}(y)=\sigma_{D,\hate}\sum_{\mathfrak s}(-1)^{|Y(\mathfrak s)|}\,y^{Y(\mathfrak s)},
\]
and
\[
   (-1)^{\mu-1}\,\LG(L;-t_0,-t_1)=\sum_{(\mathfrak s,\mathfrak s')}\ \prod_{j\in X(\mathfrak s)\sm Y(\mathfrak s')}t_0^{\ve_j}\ \prod_{j\in Y(\mathfrak s')\sm X(\mathfrak s)}t_1^{\ve_j}
\]
is a sum of $\Delta_L(-1)^2$ monomials. If moreover $\hate$ leaves a crossing as its over-strand, then $\mathfrak s\mapsto X(\mathfrak s)$ and $\mathfrak s\mapsto Y(\mathfrak s)$ are bijections from the states onto the crossing sets of the under-trees and of the over-trees of \cite{HK}.
\item If $D$ is alternating, then $\ve(V(\mathfrak s))=-\sigma_L$ for every state $\mathfrak s$, and the sign of \textup{(ii)} is $\sigma_{D,\hate}=(-1)^{(m_{D,\hate}+\sigma_L)/2}$. If moreover $L$ is fibered, with fiber $F$, the two extremal trees of Theorem~\ref{thm.intro.degrees}~(iv) are the two flavors $X(\mathfrak s_0)$, $Y(\mathfrak s_0)$ of one state, so that
\[
   \mathrm{top}(L,q)=q^{-\sigma_L}\qquad\text{and}\qquad \lambda(\xi_F)=\tfrac12\bigl(\nu+\sigma_L\bigr)=g(L)+\tfrac12\bigl(\mu-1+\sigma_L\bigr).
\]
\end{enumerate}
\end{theorem}

Part (ii) is the theorem of \cite{HK} in Kauffman's model: the sign rule is Kasteleyn's rule for the planar bipartite graph formed by the crossings and the faces of $D$, and it holds with the crossing variables, for the two flavors and with the same global sign (Lemma~\ref{lem.kasteleyn}). In the state sum of (i) the weight of a pair is the product of the weights of two independent states, corrected at the crossings of $X(\mathfrak s)\cap Y(\mathfrak s')$; at $t_0t_1=1$ the correction disappears and the sum is the square of Kauffman's state sum; in general the state sum is the expansion of the mixed discriminant of Theorem~\ref{thm.intro.flows} in the basis of the faces (Remark~\ref{rem.statesum}). The relation $\ve(V(\mathfrak s))=-\sigma_L$ of (iii) rests on the signature formula of Theorem~\ref{thm.intro.minusone} for the black checkerboard surface.

If the pairs are grouped according to the two markers which they put at every crossing, the state sum becomes a sum over Kauffman states with multiplicity two, in which every cycle along which the two states of a pair differ contributes the factor $1+t_0t_1$ (Figure~\ref{fig.gammagraph}). This is the model $Q_2$ of Ballinger \cite{BallingerTalk}: his invariant $Q_2$ is the Links--Gould invariant, which was conjectured and not known (Corollary~\ref{cor.ballinger} and the discussion which follows it). At $t_0t_1=-1$, where $\LG(L;t_0,-t_0^{-1})=\Delta_L(t_0^2)$ for every link by the first of the two specializations recalled above, the pairs $(\mathfrak s,\mathfrak s)$ alone give this value: their sum is Kauffman's state sum for $\Delta_L(t_0^2)$, and the specialization is equivalent to the cancellation of all the other pairs (Corollary~\ref{cor.diagonal}); in the form of Corollary~\ref{cor.ballinger} this cancellation is visible, since $1+t_0t_1$ vanishes there (Remark~\ref{rem.ballinger}). So the two specializations of $\LG$ are the two extreme regimes of the model.

\emph{The trefoil.} For the diagram of the trefoil considered at the beginning of this introduction there are three states (Figure~\ref{fig.trefoilstates}), for which
\[
   \bigl(X(\mathfrak s),Y(\mathfrak s)\bigr)=\bigl(\{1,2\},\emptyset\bigr),\quad\bigl(\{1\},\{3\}\bigr),\quad\bigl(\emptyset,\{2,3\}\bigr)\qquad\text{and}\qquad(-1)^{b(\mathfrak s)}=1,\ -1,\ 1 .
\]
By (i) they give the three monomials $x_1x_2$, $-x_1$, $1$ of $P_{D,\hate}$ and the three monomials $1$, $-y_3$, $y_2y_3$ of $\tP_{D,\hate}$, and the nine pairs of states give the nine monomials of the product of the two polynomials. The only pair $(\mathfrak s,\mathfrak s')$ in which $X(\mathfrak s)$ meets $Y(\mathfrak s')$ consists of the first and of the third state, which have the crossing $2$ in common. Here $m_{D,\hate}=2$, and $\ve(V(\mathfrak s))=2=-\sigma_L$ for the three states. The sum of (ii) is
\[
   \LG(3_1;-t_0,-t_1)=1+t_0+t_1+t_0^2+2t_0t_1+t_1^2+t_0^2t_1+t_0t_1^2 ,
\]
a sum of $9=\Delta_L(-1)^2$ monomials. The three pairs $(\mathfrak s,\mathfrak s)$ contribute $t_0^2+t_0t_1+t_1^2$, both to this sum and to $\LG(3_1;t_0,t_1)$; at $t_1=-t_0^{-1}$ this contribution is $t_0^2-1+t_0^{-2}=\Delta_L(t_0^2)$, and the six other pairs cancel (Examples~\ref{ex.signs} and~\ref{ex.states}).

The signs in (i) are determined in Section~\ref{ssec.signs}: for knots by the values of the edge polynomials at the point where all the variables are $1$ and by the descending diagrams, and for links by induction on the number of components, a diagram with one more component being the smoothing of a crossing (Proposition~\ref{prop.signs}). By (iii), for a fibered alternating link the top coefficient of $\LG$ and the Hopf invariant of the fiber are determined by the genus, the number of components and the signature.

On an alternating diagram the states are the spanning trees of the Tait graph (Figure~\ref{fig.taittrefoil}), and Section~\ref{ssec.tait} writes the state sum on this graph, in the lattice $\BZ^c$ spanned by the crossings: the two flavors of a state are two orientations of the graph, the orientation of the link enters through the orientation defined by the over-strands, which is the sum of an integral flow on the positive crossings and of an integral cut on the negative ones, and $-\sigma_L$ is, for an alternating knot, four times the correction term of the spin structure of the double branched cover computed by Manolescu and Owens \cite{ManolescuOwens} (Theorem~\ref{thm.tait}, Proposition~\ref{prop.omega} and Remark~\ref{rem.branched}).

Section~\ref{ssec.gluing} gives a gluing formula: the Links--Gould polynomial of a link which is cut into two tangles by a Conway sphere is a bilinear function of the polynomials of three closures of each tangle. For an alternating diagram the two tangles are two pieces of the Tait graph with two common vertices, and the way in which such pieces are glued is what the lattice does not determine (Theorem~\ref{thm.gluing} and Corollary~\ref{cor.twosum}).

Section~\ref{ssec.lattice} then says to what extent the polynomial of an alternating link is computed from this lattice: for a knot it is determined by it, and the states are the $\mathrm{Spin}^c$ structures of the double branched cover, each of them represented by an orientation of the Tait graph which depends on the root (Corollary~\ref{cor.latticeknots} and Proposition~\ref{prop.spinc}). The dependence on the root is explicit: each $\mathrm{Spin}^c$ structure determines a smoothing of the diagram, its representatives change exactly when the root passes from one circle of this smoothing to another one, and $|\Delta_L(-1)|$ of the $\Delta_L(-1)^2$ terms of the state sum are the same for all the roots and are written in the lattice (Corollary~\ref{cor.classsmoothing}, Proposition~\ref{prop.nat} and Remark~\ref{rem.rootdependence}). These terms contain the two ends of the polynomial for the total degree in $t_0$ and $t_1$, which gives the last result (Theorem~\ref{thm.extreme} and Corollary~\ref{cor.crossings}).

\begin{theorem}\label{thm.intro.crossings}
Let $D$ be a connected diagram of $L$ with $c_+$ positive and $c_-$ negative crossings.
\begin{enumerate}\renewcommand{\theenumi}{\roman{enumi}}
\item Every monomial $t_0^at_1^b$ of $\LG(L)$ satisfies $-c_-\le a+b\le c_+$.
\item If $D$ is alternating and reduced, both bounds are attained. In the state sum of Theorem~\ref{thm.intro.kauffman}~(ii), the terms with $a+b=c_+$ are in bijection with the acyclic orientations of the Tait graph of $D$ which have a given vertex as their only source, and the terms with $a+b=-c_-$ with the acyclic orientations of the dual graph which have a given vertex as their only source.
\item Consequently, if $D$ is alternating and reduced, every diagram of $L$ has at least $c_+$ positive and at least $c_-$ negative crossings. In particular $D$ has the smallest number of crossings among the diagrams of $L$, and all the reduced alternating diagrams of $L$ have the same writhe, which is zero if $L$ is isotopic to its mirror image.
\end{enumerate}
\end{theorem}

Part (iii) contains the first two conjectures of Tait on alternating links, which were proved by Kauffman, Murasugi and Thistlethwaite with the Jones polynomial and the Kauffman polynomial \cite{Kauffman87, Murasugi87, Thistlethwaite87, Thistlethwaite88, Murasugi87b}. Here the Links--Gould polynomial plays the role of the Jones polynomial: (i) is a property of the determinant formula \eqref{eq.intro.formula} for every diagram, and (ii) comes from the positivity of the state sum; the total degree sees the positive and the negative crossings separately, so that the writhe comes with the number of crossings. Thus the two degrees of $\LG$ carry the two kinds of information: for an alternating link the degree of Theorem~\ref{thm.intro.degrees} gives the genus, as the degree of the Alexander polynomial does, and the total degree gives the crossings of a reduced alternating diagram, as the span of the Jones polynomial does. This illustrates the hybrid nature of the Links--Gould polynomial, between the Alexander polynomial and the Jones polynomial.

\subsection*{Conventions}\label{ssec.conventions}
We use the notation of \cite{HK}, which we recall in Section~\ref{sec.formula}. Matrices act on column vectors; for a matrix $M$ and index sets $T,S$ we write $M[T|S]$ for the submatrix with rows $T$ and columns $S$, and $M^{\T}$ for the transpose. For a square matrix $M$ whose rows and columns are indexed by a set $\CN$ and for $d\in\CN$, $M_{\hat d}$ denotes $M$ with the row and the column of $d$ removed. $e_f$ is the unit column vector of the index $f$, and $\1$ the column vector with all entries equal to $1$. Throughout, $\BK=\BQ(t_0,t_1)$, and the letter $t$ denotes an indeterminate. For a statement $P$, the bracket $\ind{P}$ is $1$ if $P$ holds and $0$ otherwise; thus $x_k^{\ind{\ve_k=1}}$ is $x_k$ at a positive crossing and $1$ at a negative one. The sign of a crossing of an oriented diagram is the usual one: the crossing is positive if, when it is drawn with both strands pointing upward, the over-strand runs from the lower left to the upper right. As in \cite{HK}, the letter $\sigma_i$ of a braid word is a positive crossing, and links are oriented. Flows are denoted by the letter $z$ and smoothings by the letter $\vs$; Kashaev's variable, which is written $x$ in the literature, is denoted by $\zeta$. An index of notation is given at the end of the text (page~\pageref{sec.notation}): it lists the symbols which are used in more than one place, with their meaning and the place where each of them is defined.

\subsection*{Structure of the monograph}
This monograph has three parts and is organized so that proofs appear directly after the statements.

Part~\ref{part.formula} sets up the determinant formula and the lattice of flows. Section~\ref{sec.formula} recalls the formula and its expansion around the free fermion point. Section~\ref{sec.alexander} contains the classical background: the Alexander matrix of a diagram, the level lemma, the identification of the colorings with the derivations of the link group, the Alexander polynomial and the nullity. Section~\ref{sec.flows} introduces the flows and the crossing forms, proves Theorem~\ref{thm.intro.flows} by elimination along the Seifert circles, treats the basis of the faces, and writes the coupling as a mixed discriminant on the lattice of flows.

Part~\ref{part.surfaces} is about surfaces. Section~\ref{sec.statesurfaces} constructs the state surfaces and proves the localization of the linking numbers, Theorem~\ref{thm.intro.linking}~(i) and~(ii). Section~\ref{sec.seifert} identifies the crossing forms with the Seifert form band by band, constructs the deformed Seifert matrices and proves Theorem~\ref{thm.intro.seifert}; it explains why the crossing variables are cochains, and contains the examples, the Kinoshita--Terasaka knot, and the laws for the hollow handle, the braid relation and Markov stabilization. Section~\ref{sec.degree} proves Theorem~\ref{thm.intro.degrees}, describes the monomials of the second edge polynomial as disk decompositions of the canonical surface, and discusses the other Seifert surfaces. Section~\ref{sec.cocores} reads the crossing forms and the coupling in the complement of the link and of the co-cores of the bands, proves Theorem~\ref{thm.intro.cocores}, and discusses the dependence of this description on the diagram.

Part~\ref{part.faces} is about the faces. Section~\ref{sec.faces} proves Theorem~\ref{thm.intro.faces}. Section~\ref{sec.minusone} proves Theorem~\ref{thm.intro.minusone}. Section~\ref{sec.kauffman} proves Theorem~\ref{thm.intro.kauffman}, determines the signs of the states, identifies the state sum with Ballinger's model for $N=2$, writes the state sum of an alternating diagram on its Tait graph, gives a gluing formula for $\LG$ along a Conway sphere, relates the states to the $\mathrm{Spin}^c$ structures of the double branched cover and proves Theorem~\ref{thm.intro.crossings}. The conclusion discusses the question which motivates this work, that of a topological construction of $\LG$ on the link exterior, and the companion question of a categorification of $\LG$.

The three parts need not be read in a row. Part~\ref{part.formula} is used everywhere; its last subsection, Section~\ref{ssec.mixed}, which is used in Section~\ref{sec.cocores}, can be skipped at first reading. The main dependencies between the other sections are the following. Sections~\ref{sec.statesurfaces}, \ref{sec.seifert} and~\ref{sec.degree} are to be read in this order. Section~\ref{sec.cocores} uses Sections~\ref{ssec.mixed} and~\ref{sec.statesurfaces}, and it is independent of Section~\ref{sec.degree}. In Part~\ref{part.faces}, Section~\ref{sec.faces} can be read directly after Part~\ref{part.formula}. Section~\ref{sec.minusone} uses in addition the state surfaces of Section~\ref{sec.statesurfaces} and, of Section~\ref{sec.faces}, the winding numbers and the definitions. Section~\ref{sec.kauffman} uses the winding numbers of Section~\ref{ssec.windings}, the $\ve$-degree of Section~\ref{sec.degree} and, for alternating diagrams, the results of Sections~\ref{ssec.alternating}, \ref{ssec.altminusone} and~\ref{ssec.underlattice}. The paragraph which opens each part says more.

The standard diagram of the right-handed trefoil is the running example of the text. It can be followed through the text: Example~\ref{ex.trefoilflows} (the flows and the crossing forms), Example~\ref{ex.trefoilmixed} (the mixed discriminant), Example~\ref{ex.trefoilseifert} (the Seifert matrix), Example~\ref{ex.trefoil} (the deformed Seifert matrices), Example~\ref{ex.trefoilcocores} (the local system of the coupling), Example~\ref{ex.trefoilkashaev} (the deformed Kashaev matrix), Example~\ref{ex.trefoilminusone} (the forms at $x=-1$), Examples~\ref{ex.signs} and~\ref{ex.states} (the states and the state sum) and Example~\ref{ex.tait} (the state sum on the Tait graph). The standard diagram of the figure-eight knot, which has crossings of both signs, plays the same role in Examples~\ref{ex.figeight}, \ref{ex.alttrees}, \ref{ex.figeightminusone}, \ref{ex.states}, \ref{ex.tait} and~\ref{ex.fig8classes}.

\subsection*{Tool Use}No AI was used in the preparation of this manuscript, and furthermore, we do not consent for AI to be trained on this manuscript.

\subsection*{Acknowledgments}
The authors thank William Ballinger, David Cimasoni, Luke Conners, Anthony Conway, Killian Davis, Livio Ferretti, Tamás Kálmán, Rinat Kashaev, Nicolai Reshetikhin, Davide Saccardo, Boris Shapiro, Jiebo Song, Arkady Vaintrob, Emmanuel Wagner and David Wahiche for helpful discussions. BK wishes to thank his brother Mathieu for his help in these troubled times. He is also grateful to GTIIT for its hospitality during his stay in Shantou in November 2026, when this long-term project was brought to completion. MH was partially supported by NSF-RTG grant \#DMS-2135960. The work of BMK was partially supported by the SNSF research program NCCR The Mathematics of Physics (SwissMAP), and the SNSF grant no. 200021-232258.

\clearpage
\part{The determinant formula and the lattice of flows}\label{part.formula}

This part contains the algebraic framework of the monograph. Section~\ref{sec.formula} recalls the determinant formula. Section~\ref{sec.alexander} recalls what the two determinants of the formula are at the Alexander point, and introduces the levels of the Seifert circles. Section~\ref{sec.flows} is the core of the part: it introduces the lattice of flows and the crossing forms, it proves that the two determinants are the determinants of these forms in every basis of the lattice (Theorem~\ref{thm.flows}), and it describes the bases which are used in Parts~\ref{part.surfaces} and~\ref{part.faces}. Its last subsection, Section~\ref{ssec.mixed}, carries the coupling of the two flavors inside the determinants; it is used in Section~\ref{sec.cocores}, and it can be skipped at first reading.

\section{The determinant formula}\label{sec.formula}

In this section we recall the determinant formula of \cite{HK} for an arbitrary diagram, together with the expansion around the free fermion point $t_0t_1=1$. The formula is the starting point of this monograph: everything that follows is a study of the two determinants which it contains, the edge polynomials $P_{D,\hate}$ and $\tP_{D,\hate}$, and of the linear map $\BE_D$ which couples them.

\subsection{Diagrams and their edge matrices}\label{ssec.diagrams}
\begin{figure}[ht]
\centering
\begin{tikzpicture}[>=Stealth,scale=0.9,every node/.style={font=\small}]
\begin{scope}
\draw (1,-1) -- (0.2,-0.2); \draw[->] (-0.2,0.2) -- (-1,1);
\draw[->] (-1,-1) -- (1,1);
\node[below left] at (-1,-1) {$o_j$}; \node[above right] at (1,1) {$o'_j$};
\node[below right] at (1,-1) {$u_j$}; \node[above left] at (-1,1) {$u'_j$};
\node at (0,-1.7) {$\ve_j=1$};
\end{scope}
\begin{scope}[xshift=3.6cm]
\draw[->] (-1,-1) to[out=50,in=-50] (-1,1);
\draw[->] (1,-1) to[out=130,in=-130] (1,1);
\node[below left] at (-1,-1) {$o_j$}; \node[above left] at (-1,1) {$u'_j$};
\node[below right] at (1,-1) {$u_j$}; \node[above right] at (1,1) {$o'_j$};
\node at (-1.1,0) {$A_j$}; \node at (1.1,0) {$B_j$};
\node at (0,-1.7) {smoothing};
\end{scope}
\begin{scope}[xshift=8cm]
\draw (-1,-1) -- (-0.2,-0.2); \draw[->] (0.2,0.2) -- (1,1);
\draw[->] (1,-1) -- (-1,1);
\node[below right] at (1,-1) {$o_j$}; \node[above left] at (-1,1) {$o'_j$};
\node[below left] at (-1,-1) {$u_j$}; \node[above right] at (1,1) {$u'_j$};
\node at (0,-1.7) {$\ve_j=-1$};
\end{scope}
\begin{scope}[xshift=11.6cm]
\draw[->] (-1,-1) to[out=50,in=-50] (-1,1);
\draw[->] (1,-1) to[out=130,in=-130] (1,1);
\node[below left] at (-1,-1) {$u_j$}; \node[above left] at (-1,1) {$o'_j$};
\node[below right] at (1,-1) {$o_j$}; \node[above right] at (1,1) {$u'_j$};
\node at (-1.1,0) {$B_j$}; \node at (1.1,0) {$A_j$};
\node at (0,-1.7) {smoothing};
\end{scope}
\end{tikzpicture}
\caption{The edges at a crossing $j$ of sign $\ve_j$, and the oriented smoothing: the edges $o_j,u'_j$ lie on the Seifert circle $A_j$ and the edges $u_j,o'_j$ on the Seifert circle $B_j$ (Sections~\ref{ssec.level} and~\ref{ssec.smoothings}).}
\label{fig.crossing}
\end{figure}
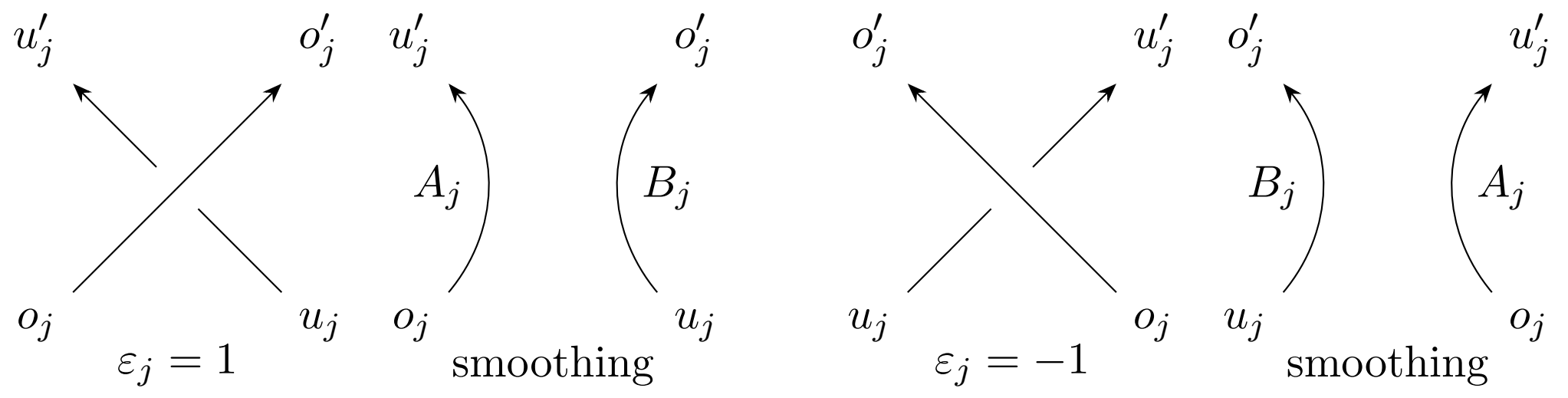

As in \cite[Section~7.1]{HK}, a \emph{diagram} is an oriented link diagram $D$ on the oriented sphere $S^2$ which is connected and has $c\ge1$ crossings, numbered $1,\dots,c$. The sign of the crossing $j$ is $\ve_j\in\{1,-1\}$. The crossings cut $D$ into $2c$ oriented arcs, the \emph{edges} of $D$; every edge leaves one crossing and enters one crossing. At the crossing $j$ we denote by $o_j,u_j$ the incoming edges of the over-strand and of the under-strand, and by $o'_j,u'_j$ the outgoing ones (Figure~\ref{fig.crossing}). Thus every edge is $o'_j$ or $u'_j$ for exactly one $j$, and $o_k$ or $u_k$ for exactly one $k$. The \emph{edge matrices} of $D$ are the matrices with rows and columns indexed by the edges
\begin{align}
   T^a_D(x)&=\sum_{j=1}^c\Bigl(E_{o'_jo_j}+(1-x_j)\,E_{u'_jo_j}+x_j\,E_{u'_ju_j}\Bigr),\label{eq.Ta}\\
   T^b_D(y)&=\sum_{j=1}^c\Bigl(y_j\,E_{o'_jo_j}+(1-y_j)\,E_{u'_jo_j}+E_{u'_ju_j}\Bigr),\label{eq.Tb}
\end{align}
where $x_1,\dots,x_c,y_1,\dots,y_c$ are indeterminates and $E_{e'e}$ is the matrix unit with entry $1$ in row $e'$ and column $e$ \cite[Definition~6.3 and Section~7.1]{HK}. This means that the row $o'_j$ of $T^a_D(x)$ is $e_{o_j}^{\T}$, and the row $u'_j$ is $(1-x_j)e_{o_j}^{\T}+x_je_{u_j}^{\T}$; the column $u_j$ of $T^b_D(y)$ is $e_{u'_j}$, and the column $o_j$ is $y_je_{o'_j}+(1-y_j)e_{u'_j}$. In particular
\begin{equation}\label{eq.sums}
   \text{the rows of $T^a_D(x)$ have sum $1$, and the columns of $T^b_D(y)$ have sum $1$,}
\end{equation}
so that $(1-T^a_D(x))\1=0$ and $\1^{\T}(1-T^b_D(y))=0$. We write
\begin{equation}\label{eq.PP}
   P_{D,\hate}(x)=\det\bigl((1-T^a_D(x))_{\hate}\bigr),\qquad \tP_{D,\hate}(y)=\det\bigl((1-T^b_D(y))_{\hate}\bigr)
\end{equation}
for the determinants of the edge matrices with the row and the column of an edge $\hate$ removed. The matrix $1-T^a_D(x)$ is the sum of a constant matrix and of the matrices $-x_j\,e_{u'_j}(e_{u_j}-e_{o_j})^{\T}$ of rank one: the variable $x_j$ occurs only in the row $u'_j$, in which it occurs linearly, and since the determinant is a linear function of each row, $P_{D,\hate}(x)$ has degree $\le1$ in each variable $x_j$. Likewise $1-T^b_D(y)$ is the sum of a constant matrix and of the matrices $-y_j\,(e_{o'_j}-e_{u'_j})e_{o_j}^{\T}$, the variable $y_j$ occurs only in the column $o_j$, and $\tP_{D,\hate}(y)$ has degree $\le1$ in each $y_j$ (compare \cite[(28) and Lemma~3.4]{HK}).

Let $\BE_D$ be the $\BK$-linear map, defined on the polynomials of degree $\le1$ in each of the variables $x_j,y_j$, with
\begin{equation}\label{eq.ED}
   \BE_D\Bigl[\prod_{j\in S}x_j\prod_{j\in S'}y_j\Bigr]=\prod_{j\in S\sm S'}t_0^{\ve_j}\ \prod_{j\in S'\sm S}t_1^{\ve_j}\qquad(S,S'\subseteq\{1,\dots,c\}).
\end{equation}
So $\BE_D$ is the linear substitution $x_j\mapsto t_0^{\ve_j}$ and $y_j\mapsto t_1^{\ve_j}$, except that the product $x_jy_j$ of the two variables of the same crossing is replaced by $1$.

\begin{theorem}[{\cite[Theorem~7.8]{HK}}]\label{thm.HK}
Let $L$ be an oriented link with $\mu$ components, let $D$ be a diagram of $L$ and let $\hate$ be an edge of $D$. Then
\[
   \LG(L;t_0,t_1)=V_1(L;t_0,t_1)=(-1)^{\mu-1}\ \BE_D\bigl[P_{D,\hate}(x)\,\tP_{D,\hate}(y)\bigr].
\]
\end{theorem}

Here $V_1(L;t_0,t_1)$ is the $V_1$-polynomial of Garoufalidis and Kashaev \cite{GaroufalidisKashaev} in the normalization of \cite{GHKST}, which coincides with the Links--Gould polynomial $\LG(L;t_0,t_1)$ in the variables of Ishii \cite{GHKKST}; see \cite[Section~2]{HK}. From now on we write $\LG(L)$ for this polynomial; $\LG$ of the unknot is $1$.

\subsection{The expansion around the free fermion point}
The map $\BE_D$ is linear but not multiplicative: it treats the two variables $x_j,y_j$ of a crossing together. The expansion of this subsection isolates this interaction. It writes $\LG$ as a sum over the sets $S$ of crossings at which the two determinants are coupled, and the term of $S$ is divisible by $(1-t_0t_1)^{|S|}$; at the free fermion point $t_0t_1=1$ only the term of the empty set remains, and it is a product of two Alexander polynomials.

Put
\begin{equation}\label{eq.xbar}
   \bx_j=t_0^{\ve_j},\qquad \by_j=t_1^{\ve_j},\qquad g_j=1-\bx_j\by_j=1-(t_0t_1)^{\ve_j}\qquad(j=1,\dots,c).
\end{equation}
Since $P_{D,\hate}$ and $\tP_{D,\hate}$ have degree $\le1$ in each variable, they can be written as
\begin{equation}\label{eq.ab}
   P_{D,\hate}(x)=\sum_{S\subseteq\{1,\dots,c\}}a_S\prod_{j\in S}(x_j-\bx_j),\qquad
   \tP_{D,\hate}(y)=\sum_{S\subseteq\{1,\dots,c\}}b_S\prod_{j\in S}(y_j-\by_j),
\end{equation}
with uniquely determined coefficients $a_S=a_S(t_0)\in\BZ[t_0^{\pm1}]$ and $b_S=b_S(t_1)\in\BZ[t_1^{\pm1}]$: the entries of the matrices are polynomials with integer coefficients in the variables, and $\bx_j,\by_j$ are units.

\begin{corollary}[{\cite[Corollary~5.8]{HK}}]\label{cor.expansion}
With the notation \eqref{eq.xbar} and \eqref{eq.ab},
\[
   \LG(L;t_0,t_1)=(-1)^{\mu-1}\sum_{S\subseteq\{1,\dots,c\}}\ \Bigl(\prod_{j\in S}g_j\Bigr)\,a_S(t_0)\,b_S(t_1).
\]
\end{corollary}
\begin{proof}
Put $\xi_j=x_j-\bx_j$ and $\xi'_j=y_j-\by_j$. By \eqref{eq.ED}, $\BE_D[\xi_j]=\BE_D[\xi'_j]=0$ and $\BE_D[\xi_j\xi'_j]=1-\bx_j\by_j=g_j$. The map $\BE_D$ is multiplicative over distinct crossings: a product $\prod_jf_j$ of polynomials $f_j$ of degree $\le1$ in $x_j,y_j$ and involving no other variables satisfies $\BE_D[\prod_jf_j]=\prod_j\BE_D[f_j]$, both sides being multilinear in the $f_j$ and equal on monomials. Hence
\[
   \BE_D\Bigl[\prod_{j\in S}\xi_j\prod_{j\in S'}\xi'_j\Bigr]=\prod_{j\in S}g_j\ \text{ if $S=S'$,\qquad and $0$ otherwise.}
\] Expanding the product $P_{D,\hate}(x)\tP_{D,\hate}(y)$ by \eqref{eq.ab} and applying Theorem~\ref{thm.HK} gives the claim. In \cite{HK} this is stated for closed braids; the proof is the same for every diagram.
\end{proof}

On the lattice of flows of Section~\ref{sec.flows} the coefficients $a_S$ and $b_S$ are themselves determinants of linear forms, and the expansion is that of a mixed discriminant (Remark~\ref{rem.mixed}~(ii)).

The term $S=\emptyset$ is $P_{D,\hate}(\bx)\tP_{D,\hate}(\by)$, the product of the two Burau determinants at $t_0$ and at $t_1$; it is the product $\Delta_L(t_0)\Delta_L(t_1)$ of Alexander polynomials up to a unit (Lemma~\ref{lem.alexpoly} below). All other terms are divisible by $1-t_0t_1$, because
\begin{equation}\label{eq.gdiv}
   g_j=1-t_0t_1\ \ \text{if $\ve_j=1$},\qquad g_j=1-(t_0t_1)^{-1}=-(t_0t_1)^{-1}(1-t_0t_1)\ \ \text{if $\ve_j=-1$}.
\end{equation}

\section{The Alexander matrix of a diagram and the nullity}\label{sec.alexander}

Throughout this section, $D$ is a diagram of a link $L$ with $\mu$ components in the sense of Section~\ref{ssec.diagrams}, with crossings $1,\dots,c$ and a set of edges $E$, $|E|=2c$. For a field $\BF$ and an element $t\in\BF$, $t\ne0$, we write $t^{\ve}$ for the point $x_j=t^{\ve_j}$ (resp.\ $y_j=t^{\ve_j}$), $j=1,\dots,c$, and we call $x=t^{\ve}$ the \emph{Alexander point} of $D$ at $t$. The main case is $\BF=\BQ(t)$ with $t$ an indeterminate. The material of this section is classical. We prove the level lemma, which describes the Seifert circles of $D$ and relates the two flavors of the determinant formula, and we record what the Alexander matrix of $D$ is: the Alexander matrix of the Wirtinger presentation of the link group, whose kernel is identified with the derivations of the group (Section~\ref{ssec.derivations}) and whose corank is the rank of the Alexander module of $L$ (Section~\ref{ssec.alexnullity}). The section has two roles in what follows. The level function governs the degrees of the edge polynomials and of the deformed Seifert matrices (Sections~\ref{ssec.regions} and~\ref{sec.degree}), and the identification of the Burau determinants with the Alexander polynomial (Lemma~\ref{lem.alexpoly}) is what relates each deformed model to the classical model which it deforms.

\subsection{The Alexander matrix of a diagram and its colorings}
At the Alexander point the edge matrices of Section~\ref{sec.formula} become matrices with entries in $\BF$. We first give names to these matrices and to the vectors of their kernels.

\begin{definition}\label{def.coloring}
The \emph{Alexander matrix} of $D$ at $t$ is the $|E|\times |E|$ matrix
\[
   A_D(t)=1-T^a_D(t^{\ve})\qquad\text{with entries in }\BF.
\]
An \emph{Alexander coloring} (or simply a coloring) of $D$ at $t$ is a vector $w\in\BF^E$ in the kernel of $A_D(t)$, that is, a function $e\mapsto w_e$ on the edges of $D$ such that at every crossing $j$
\begin{equation}\label{eq.coloring}
   w_{o'_j}=w_{o_j},\qquad w_{u'_j}=t^{\ve_j}\,w_{u_j}+(1-t^{\ve_j})\,w_{o_j}.
\end{equation}
For an edge $\hate$ we write $A_{D,\hate}(t)=A_D(t)_{\hate}$ for the Alexander matrix with the row and the column of $\hate$ removed (the \emph{cut matrix}), and similarly $B_D(t)=1-T^b_D(t^{\ve})$ and $B_{D,\hate}(t)=B_D(t)_{\hate}$.
\end{definition}

The equations \eqref{eq.coloring} are read off from \eqref{eq.Ta}: the row $o'_j$ of $A_D(t)$ is $e_{o'_j}^{\T}-e_{o_j}^{\T}$ and the row $u'_j$ is $e_{u'_j}^{\T}-t^{\ve_j}e_{u_j}^{\T}-(1-t^{\ve_j})e_{o_j}^{\T}$. The constant vector $\1$ is a coloring, by \eqref{eq.sums}. In the notation of Section~\ref{sec.formula}, $P_{D,\hate}(t^{\ve})=\det A_{D,\hate}(t)$ and $\tP_{D,\hate}(t^{\ve})=\det B_{D,\hate}(t)$.

\subsection{The level function}\label{ssec.level}
Smoothing every crossing of $D$ in the way which is compatible with the orientations produces the \emph{Seifert circles} of $D$: disjoint oriented simple closed curves on $S^2$, and every edge of $D$ lies on exactly one of them. At a crossing the two smoothed arcs run side by side in the same direction, with the region between them on the right of one arc and on the left of the other; since the two sides of a simple closed curve are its two complementary components, the two arcs lie on distinct circles: a crossing joins two distinct Seifert circles. We draw $D$ in the plane, by removing from $S^2$ a point which does not lie on $D$.

\begin{definition}\label{def.level}
For a Seifert circle $C$ of $D$ let $s(C)=1$ if $C$ is oriented clockwise and $s(C)=-1$ if it is oriented counterclockwise, and let
\[
   h(C)=\sum_{C'}s(C')\ +\ \frac{s(C)-1}{2}\ \in\BZ,
\]
the sum running over the Seifert circles $C'\ne C$ whose interior (the bounded component of the complement of $C'$ in the plane) contains $C$. The \emph{level} of an edge $e$ is $h(e)=h(C)$, where $C$ is the Seifert circle on which $e$ lies.
\end{definition}

\begin{figure}[ht]
\centering
\begin{tikzpicture}[>=Stealth,scale=0.75,every node/.style={font=\small},
  ccw/.style={postaction={decorate},decoration={markings,mark=at position 0.25 with {\arrow{>}},mark=at position 0.75 with {\arrow{>}}}},
  cw/.style={postaction={decorate},decoration={markings,mark=at position 0.25 with {\arrow{<}},mark=at position 0.75 with {\arrow{<}}}}]
\begin{scope}
\draw[ccw] (-1.35,0) circle (1.2); \draw[cw] (1.35,0) circle (1.2);
\node at (-1.35,0) {$C_{\rm L}$}; \node at (1.35,0) {$C_{\rm R}$};
\node at (0,0.55) {$\Omega$};
\node[font=\scriptsize] at (0,-2.6) {$s(C_{\rm L})=-1,\ s(C_{\rm R})=1$};
\node[font=\scriptsize] at (0,-3.15) {$h(C_{\rm L})=-1,\ h(C_{\rm R})=0$};
\end{scope}
\begin{scope}[xshift=7cm]
\draw[cw] (0,0) circle (2.1); \draw[cw] (-0.95,0) circle (0.9);
\node at (0.9,0.8) {$C_{\rm L}$}; \node at (-0.95,0) {$C_{\rm R}$};
\node at (-1.65,1.05) {$\Omega$};
\node[font=\scriptsize] at (0,-2.6) {$s(C_{\rm L})=s(C_{\rm R})=1$};
\node[font=\scriptsize] at (0,-3.15) {$h(C_{\rm L})=0,\ h(C_{\rm R})=1$};
\end{scope}
\begin{scope}[xshift=14cm]
\draw[ccw] (0,0) circle (2.1); \draw[ccw] (0.95,0) circle (0.9);
\node at (-0.9,0.8) {$C_{\rm R}$}; \node at (0.95,0) {$C_{\rm L}$};
\node at (1.65,1.05) {$\Omega$};
\node[font=\scriptsize] at (0,-2.6) {$s(C_{\rm L})=s(C_{\rm R})=-1$};
\node[font=\scriptsize] at (0,-3.15) {$h(C_{\rm L})=-2,\ h(C_{\rm R})=-1$};
\end{scope}
\end{tikzpicture}
\caption{The three cases in the proof of the level lemma. In each case the two arcs of the crossing are the arcs of $C_{\rm L}$ and $C_{\rm R}$ bounding the region $\Omega$, both pointing upward; the levels are those of the case where no other Seifert circle encloses $C_{\rm L}$ or $C_{\rm R}$.}
\label{fig.level}
\end{figure}
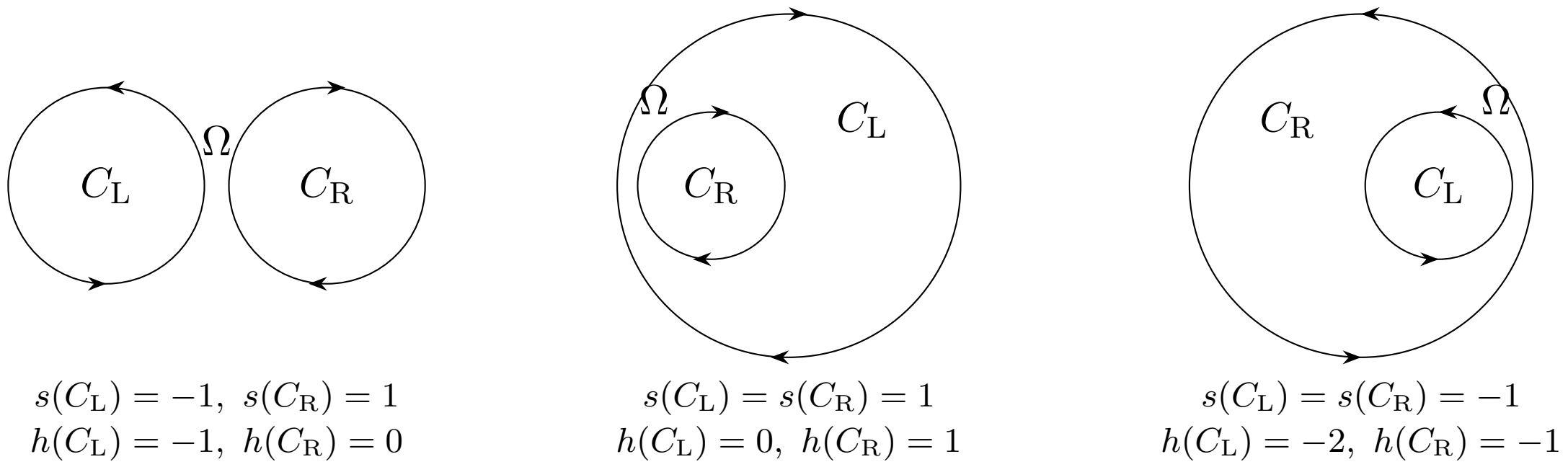

\begin{lemma}[the level lemma]\label{lem.level}
At every crossing $j$,
\begin{equation}\label{eq.levels}
   h(u'_j)=h(o_j),\qquad h(o'_j)=h(u_j),\qquad h(u_j)=h(o_j)+\ve_j .
\end{equation}
Consequently, with $\DD_t=\diag\bigl(t^{h(e)}\bigr)_{e\in E}$:
\begin{enumerate}\renewcommand{\theenumi}{\roman{enumi}}
\item $B_D(t)=\DD_t\,A_D(t)\,\DD_t^{-1}$, and $B_{D,\hate}(t)=(\DD_t)_{\hate}\,A_{D,\hate}(t)\,(\DD_t)_{\hate}^{-1}$ for every edge $\hate$;
\item the row vector $\ell_t=(t^{h(e)})_{e\in E}$ satisfies $\ell_t\,A_D(t)=0$.
\end{enumerate}
\end{lemma}
\begin{proof}
Draw the crossing $j$ with both strands pointing upward, and let $C_{\rm L}$ and $C_{\rm R}$ be the Seifert circles which contain the arcs on the left and on the right. The smoothing joins the incoming left edge to the outgoing left edge and the incoming right edge to the outgoing right edge. For $\ve_j=1$ the over-strand comes in from the left, so that $o_j,u_j,o'_j,u'_j$ are the incoming left, incoming right, outgoing right and outgoing left edges; for $\ve_j=-1$ it comes in from the right, and $o_j,u_j,o'_j,u'_j$ are the incoming right, incoming left, outgoing left and outgoing right edges. In both cases $o_j$ and $u'_j$ lie on one circle and $u_j$ and $o'_j$ on the other, which gives the first two identities of \eqref{eq.levels}, and the third one is equivalent to $h(C_{\rm R})=h(C_{\rm L})+1$ in both cases.

To prove $h(C_{\rm R})=h(C_{\rm L})+1$, let $\Omega$ be the component of the complement of the Seifert circles which contains the region between the two arcs at the crossing. It lies to the right of $C_{\rm L}$ and to the left of $C_{\rm R}$. A circle $C'$ different from $C_{\rm L},C_{\rm R}$ whose interior contains exactly one of them would separate $C_{\rm L}$ from $C_{\rm R}$, which is impossible since $\Omega$ is adjacent to both; so the circles enclosing $C_{\rm L}$ and those enclosing $C_{\rm R}$ differ at most by $C_{\rm L}$ or $C_{\rm R}$ themselves. There are three cases (Figure~\ref{fig.level}).
\begin{enumerate}
\item Neither circle encloses the other. Then $\Omega$ lies outside both. The exterior of $C_{\rm L}$ is on its right, so $C_{\rm L}$ is counterclockwise and $s(C_{\rm L})=-1$; the exterior of $C_{\rm R}$ is on its left, so $C_{\rm R}$ is clockwise and $s(C_{\rm R})=1$. The enclosing circles are the same, and $h(C_{\rm R})-h(C_{\rm L})=\frac{1-1}2-\frac{-1-1}2=1$.
\item $C_{\rm R}$ lies inside $C_{\rm L}$. Then $\Omega$ lies inside $C_{\rm L}$ and outside $C_{\rm R}$. The interior of $C_{\rm L}$ is on its right, so $s(C_{\rm L})=1$; the exterior of $C_{\rm R}$ is on its left, so its interior is on its right and $s(C_{\rm R})=1$. The circles enclosing $C_{\rm R}$ are those enclosing $C_{\rm L}$ together with $C_{\rm L}$, so $h(C_{\rm R})-h(C_{\rm L})=s(C_{\rm L})=1$.
\item $C_{\rm L}$ lies inside $C_{\rm R}$. Then $\Omega$ lies inside $C_{\rm R}$ and outside $C_{\rm L}$; as in the previous case $s(C_{\rm L})=s(C_{\rm R})=-1$, the circles enclosing $C_{\rm L}$ are those enclosing $C_{\rm R}$ together with $C_{\rm R}$, and $h(C_{\rm R})-h(C_{\rm L})=-s(C_{\rm R})=1$.
\end{enumerate}
This proves \eqref{eq.levels}.

(i) The entries of $\DD_tT^a_D(t^{\ve})\DD_t^{-1}$ at the crossing $j$ are $t^{h(o'_j)-h(o_j)}=t^{\ve_j}$ in position $(o'_j,o_j)$, $t^{h(u'_j)-h(o_j)}(1-t^{\ve_j})=1-t^{\ve_j}$ in position $(u'_j,o_j)$, and $t^{h(u'_j)-h(u_j)}t^{\ve_j}=1$ in position $(u'_j,u_j)$; by \eqref{eq.Tb} these are the entries of $T^b_D(t^{\ve})$. Conjugation by a diagonal matrix commutes with the deletion of the row and the column of the same index.

(ii) We check $\ell_tT^a_D(t^{\ve})=\ell_t$ column by column. The column $o_j$ of $T^a_D(t^{\ve})$ has the entries $1$ in row $o'_j$ and $1-t^{\ve_j}$ in row $u'_j$, and $t^{h(o'_j)}+(1-t^{\ve_j})t^{h(u'_j)}=t^{h(o_j)+\ve_j}+t^{h(o_j)}-t^{h(o_j)+\ve_j}=t^{h(o_j)}$; the column $u_j$ has the single entry $t^{\ve_j}$ in row $u'_j$, and $t^{\ve_j}t^{h(u'_j)}=t^{h(o_j)+\ve_j}=t^{h(u_j)}$.
\end{proof}

\begin{remark}\label{rem.level}
For the closure diagram of a braid, drawn with the strands running upward and closed on the right, the Seifert circles are the nested clockwise circles of the positions, the circle of the position $1$ being the outermost one, and $h(e)=p-1$ for an edge in the position $p$: $\DD_t$ is the matrix $\diag(1,t,\dots,t^{n-1})$ of \cite[Remark~5.9]{HK}. In general $h$ is determined by \eqref{eq.levels} up to an additive constant, because the graph whose vertices are the Seifert circles and whose edges are the crossings is connected. The function $h$ depends on the choice of the point at infinity, but the statements (i) and (ii) do not.
\end{remark}

\begin{lemma}[the cut matrix]\label{lem.cut}
Let $\hate$ be an edge of $D$.
\begin{enumerate}\renewcommand{\theenumi}{\roman{enumi}}
\item Restriction to $E\sm\{\hate\}$ is an isomorphism from $\{w\in\ker A_D(t):w_{\hate}=0\}$ onto $\ker A_{D,\hate}(t)$, and from $\{n\in\BF^E: n^{\T}A_D(t)=0,\ n_{\hate}=0\}$ onto the left kernel $\{n':n'^{\T}A_{D,\hate}(t)=0\}$ of the cut matrix.
\item $\corank A_{D,\hate}(t)=\corank A_D(t)-1$.
\item The same statements hold for $B_D(t)$ and $B_{D,\hate}(t)$; in particular $\corank B_{D,\hate}(t)=\corank A_{D,\hate}(t)$.
\end{enumerate}
\end{lemma}
\begin{proof}
(i) Let $w'\in\ker A_{D,\hate}(t)$ and let $w$ be its extension by $w_{\hate}=0$. For $e\ne\hate$, the entry $(A_D(t)w)_e$ is the row $e$ of $A_D(t)$ applied to $w$; since $w_{\hate}=0$, this is the row $e$ of the cut matrix applied to $w'$, which is zero. For the remaining entry we use the left null vector $\ell_t$ of Lemma~\ref{lem.level}: $0=\ell_tA_D(t)w=\sum_et^{h(e)}(A_D(t)w)_e=t^{h(\hate)}(A_D(t)w)_{\hate}$, so $(A_D(t)w)_{\hate}=0$. Thus $w\in\ker A_D(t)$, and restriction and extension by zero are inverse to each other. The statement for the left kernels is the same, with the right null vector $\1$ in place of $\ell_t$: if $n'^{\T}A_{D,\hate}(t)=0$ and $n$ is the extension of $n'$ by $n_{\hate}=0$, then $(n^{\T}A_D(t))_e=0$ for $e\ne\hate$, and $0=n^{\T}A_D(t)\1=\sum_e(n^{\T}A_D(t))_e=(n^{\T}A_D(t))_{\hate}$.

(ii) The map $w\mapsto w_{\hate}$ is onto $\BF$ on $\ker A_D(t)$, because $\1$ is in the kernel; its kernel is the subspace of (i). Hence $\dim\ker A_{D,\hate}(t)=\dim\ker A_D(t)-1$, and the coranks of the square matrices $A_D(t)$ and $A_{D,\hate}(t)$ are the dimensions of their kernels.

(iii) By Lemma~\ref{lem.level} (i), $B_D(t)$ is conjugate to $A_D(t)$ by a diagonal matrix, $\ell_t\DD_t^{-1}=\1^{\T}$ is a left null vector and $\DD_t\1$ a right null vector of $B_D(t)$, both with nonzero entries. The proof of (i) and (ii) applies directly.
\end{proof}

\subsection{Colorings are derivations}\label{ssec.derivations}
The matrix $A_D(t)$ has a classical meaning: it is the Alexander matrix of the Wirtinger presentation of the group of the link with one generator for each edge of the diagram, and its kernel, the space of the colorings of $D$ at $t$, is the space of the derivations of the group of the link with values in $\BF$, on which the group acts through the total linking number. In this subsection we prove these two statements; the Alexander polynomial and the nullity are treated in Section~\ref{ssec.alexnullity}.

Let $G=\pi_1(X)$ be the group of the link\label{not.group}, $X=S^3\sm L$. We use the Wirtinger presentation of $G$ in the form in which there is one generator for each edge of $D$; see \cite[Chapter~VI]{CrowellFox}, \cite[Section~3.D]{Rolfsen}, \cite[Section~3.A]{BurdeZieschang}. Let $F$ be the free group on generators $x_e$, $e\in E$, where $x_e$ stands for a meridian loop of the edge $e$, based at a point above the plane of the diagram. For each crossing $j$ put
\begin{equation}\label{eq.relators}
   R_j=x_{o'_j}\,x_{o_j}^{-1},\qquad r_j=x_{u'_j}\,\bigl(x_{o_j}^{\eta_j}\,x_{u_j}\,x_{o_j}^{-\eta_j}\bigr)^{-1},
\end{equation}
where $\eta_j\in\{1,-1\}$ is a sign. Then
\begin{equation}\label{eq.presentation}
   G\cong\langle x_e,\ e\in E\mid R_j,\ r_j,\ j=1,\dots,c\rangle .
\end{equation}
The relations $R_j$ identify the meridians of the edges which belong to the same arc of the diagram (the pieces of $D$ between two consecutive undercrossings), and $r_j$ is the Wirtinger relation at the crossing $j$, which says that the meridian of the outgoing under-edge is the conjugate of the meridian of the incoming under-edge by the meridian of the over-strand: adding the generators $x_{o'_j}$ and the relations $R_j$ to the usual Wirtinger presentation with one generator per arc is a Tietze transformation. Here $\eta_j=\ve_j$ for all $j$, or $\eta_j=-\ve_j$ for all $j$, depending on the conventions chosen for the orientation of $S^3$ and of the meridians; nothing below depends on this choice (Remark~\ref{rem.convention}). Note that every meridian $x_e$ is sent to $1$ by the total linking number homomorphism $\phi\colon G\to\BZ$.

Let $M$ be a left $G$-module. A \emph{derivation} of $G$ with values in $M$ is a map $d\colon G\to M$ with $d(gh)=d(g)+g\cdot d(h)$ for all $g,h\in G$; the set of derivations forms a vector space $\Der(G,M)$ if $M$ is a vector space over a field on which $G$ acts linearly. A derivation satisfies $d(1)=0$ and $d(g^{-1})=-g^{-1}\cdot d(g)$. For $m\in M$ the map $g\mapsto g\cdot m-m$ is a derivation, called \emph{principal}; the principal derivations form the subspace $\mathrm{PDer}(G,M)$, and $H^1(G;M)=\Der(G,M)/\mathrm{PDer}(G,M)$ is the first cohomology group of $G$ with coefficients in $M$ \cite[Chapter~III]{Brown}. We shall use two elementary facts.

\begin{lemma}\label{lem.der}
\begin{enumerate}\renewcommand{\theenumi}{\roman{enumi}}
\item Let $F$ be a free group on a set of generators $\{x_e\}$ and $M$ an $F$-module. For every family $(m_e)$ of elements of $M$ there is a unique derivation $d\colon F\to M$ with $d(x_e)=m_e$ for all $e$.
\item Let $G=F/N$, where $N$ is the normal subgroup of $F$ generated by elements $r_1,\dots,r_m$, and let $M$ be a $G$-module, viewed as an $F$-module. A derivation $d\colon F\to M$ vanishes on $N$ (and hence induces a derivation $G\to M$) if and only if $d(r_1)=\dots=d(r_m)=0$.
\end{enumerate}
\end{lemma}
\begin{proof}
(i) Consider the semidirect product $M\rtimes F$, that is, the set $M\times F$ with the product $(m,g)(m',g')=(m+g\cdot m',\,gg')$. A map $d\colon F\to M$ is a derivation if and only if $g\mapsto(d(g),g)$ is a group homomorphism from $F$ to $M\rtimes F$. Since $F$ is free, there is exactly one homomorphism with $x_e\mapsto(m_e,x_e)$; its first component is the required derivation, and every derivation with these values on the generators arises in this way.

(ii) The elements of $N$ act trivially on $M$. For $g\in F$ and $n\in N$, $d(gng^{-1})=d(g)+g\cdot d(n)+gn\cdot d(g^{-1})=d(g)+g\cdot d(n)-gng^{-1}\cdot d(g)=g\cdot d(n)$, and $d(nn')=d(n)+d(n')$, $d(n^{-1})=-d(n)$ for $n,n'\in N$. Hence $d$ vanishes on $N$, which consists of the products of conjugates of the $r_i^{\pm1}$, if and only if it vanishes on $r_1,\dots,r_m$.
\end{proof}

Now let $\BF$ be a field, $t\in\BF\sm\{0,1\}$, and let $\BF_t$ denote the vector space $\BF$ with the action of $G$ through the total linking number: $g\cdot m=t^{\phi(g)}m$. This action factors through the abelianization of $G$; every meridian acts by multiplication by $t$. The action is nontrivial since $t\ne1$, so a principal derivation $g\mapsto(t^{\phi(g)}-1)m$ vanishes only for $m=0$, and $\mathrm{PDer}(G,\BF_t)\cong\BF$ is one-dimensional.

\begin{proposition}\label{prop.colorings}
Assume that the relations \eqref{eq.relators} hold with $\eta_j=\ve_j$ for all $j$. For $w\in\BF^E$ let $d_w\colon F\to\BF_t$ be the derivation with $d_w(x_e)=w_e$. Then $w\mapsto d_w$ induces an isomorphism from the space $\ker A_D(t)$ of colorings of $D$ at $t$ onto $\Der(G,\BF_t)$. Consequently
\[
   \dim_{\BF}\ker A_D(t)=1+\dim_{\BF}H^1(G;\BF_t)\qquad\text{and}\qquad \corank A_{D,\hate}(t)=\dim_{\BF}H^1(G;\BF_t)
\]
for every edge $\hate$.
\end{proposition}
\begin{proof}
By Lemma~\ref{lem.der} (i), $w\mapsto d_w$ is an isomorphism from $\BF^E$ onto $\Der(F,\BF_t)$, and by Lemma~\ref{lem.der} (ii), $d_w$ induces a derivation of $G$ if and only if $d_w(R_j)=d_w(r_j)=0$ for all $j$. Conversely, every derivation of $G$ composed with $F\to G$ is a derivation of $F$ vanishing on the relators, so $\Der(G,\BF_t)$ is identified with the space of $w$ such that $d_w$ vanishes on all relators. We compute, using $d(g^{-1})=-g^{-1}\cdot d(g)$ and the fact that every generator acts by $t$:
\[
   d_w(R_j)=d_w(x_{o'_j})+x_{o'_j}\cdot d_w(x_{o_j}^{-1})=w_{o'_j}-t\cdot t^{-1}w_{o_j}=w_{o'_j}-w_{o_j}.
\]
For $r_j$ we write $r_j=x_{u'_j}\,q^{-1}$ with $q=x_{o_j}^{\ve_j}x_{u_j}x_{o_j}^{-\ve_j}$, so that $d_w(r_j)=w_{u'_j}+x_{u'_j}\cdot d_w(q^{-1})=w_{u'_j}-x_{u'_j}q^{-1}\cdot d_w(q)=w_{u'_j}-d_w(q)$, since $x_{u'_j}q^{-1}=r_j$ acts trivially. For $\ve_j=1$,
\[
   d_w(q)=d_w(x_{o_j})+x_{o_j}\cdot d_w(x_{u_j})+x_{o_j}x_{u_j}\cdot d_w(x_{o_j}^{-1})
         =w_{o_j}+t\,w_{u_j}-t^2\cdot t^{-1}w_{o_j}=t\,w_{u_j}+(1-t)\,w_{o_j},
\]
and for $\ve_j=-1$,
\begin{align*}
   d_w(q)&=d_w(x_{o_j}^{-1})+x_{o_j}^{-1}\cdot d_w(x_{u_j})+x_{o_j}^{-1}x_{u_j}\cdot d_w(x_{o_j})\\
         &=-t^{-1}w_{o_j}+t^{-1}w_{u_j}+w_{o_j}=t^{-1}\,w_{u_j}+(1-t^{-1})\,w_{o_j}.
\end{align*}
So $d_w(R_j)=0$ and $d_w(r_j)=0$ for all $j$ if and only if $w$ satisfies the equations \eqref{eq.coloring}, i.e.\ $w\in\ker A_D(t)$. This proves the first statement, and the formula for $\dim\ker A_D(t)$ follows because the principal derivations form a one-dimensional subspace. The last formula is Lemma~\ref{lem.cut} (ii).
\end{proof}

\begin{remark}[the convention]\label{rem.convention}
If the relations hold with $\eta_j=-\ve_j$, the same computation identifies $\ker A_D(t)$ with $\Der(G,\BF_{t^{-1}})$. For $\BF=\BQ(t)$ this makes no difference to the dimensions: the automorphism $t\mapsto t^{-1}$ of the field $\BQ(t)$ transforms $\BQ(t)_t$ into $\BQ(t)_{t^{-1}}$ and $\Der(G,\BQ(t)_t)$ bijectively and semilinearly onto $\Der(G,\BQ(t)_{t^{-1}})$, and likewise for the spaces of principal derivations. The same automorphism applied to the entries of $A_D(t)$ gives $A_D(t^{-1})$, so that $\dim\ker A_D(t)=\dim\ker A_D(t^{-1})$ in any case. We also note that the computation of the proof, read for the free group $F$ alone, identifies $A_D(t)$ with the Fox Jacobian matrix $\bigl(\partial R_j/\partial x_e,\ \partial r_j/\partial x_e\bigr)$ of the presentation \eqref{eq.presentation}, evaluated under the abelianization $x_e\mapsto t$ \cite{Fox, CrowellFox}: by the fundamental formula of Fox calculus, $d(r)=\sum_e\bigl(\partial r/\partial x_e\bigr)\cdot d(x_e)$ for every derivation $d$ of $F$ and every $r\in F$, where the group ring element $\partial r/\partial x_e$ acts through the abelianization; so the linear forms $w\mapsto d_w(R_j)$ and $w\mapsto d_w(r_j)$ computed in the proof, which are the rows $o'_j$ and $u'_j$ of $A_D(t)$ by \eqref{eq.coloring}, have the coefficients $\partial R_j/\partial x_e$ and $\partial r_j/\partial x_e$ evaluated at $x_e\mapsto t$. In other words, $A_D(t)$ is the Alexander matrix of the Wirtinger presentation of $G$ with one generator for each edge, the row of the edge $o'_j$ corresponding to the relator $R_j$ and the row of $u'_j$ to $r_j$.
\end{remark}

\subsection{The Alexander polynomial and the nullity}\label{ssec.alexnullity}
The identification of $A_D(t)$ with an Alexander matrix gives the classical interpretation of the Burau determinants of Section~\ref{sec.formula}, which is Lemma~\ref{lem.alexpoly}. Recall \cite[Chapter~VII]{CrowellFox} that the one-variable Alexander polynomial $\Delta_L(t)\in\BZ[t^{\pm1}]$ of $L$ is defined, up to multiplication by a unit $\pm t^m$, as a generator of the smallest principal ideal containing the first elementary ideal $E_1(G)$ of $G$ with respect to $\phi$: for any finite presentation $\langle x_1,\dots,x_n\mid r_1,\dots,r_m\rangle$ of $G$, $E_1(G)$ is the ideal of $\BZ[t^{\pm1}]$ generated by the $(n-1)\times(n-1)$ minors of the $m\times n$ Alexander matrix $\bigl(\partial r_i/\partial x_k\bigr)$ evaluated under $x_k\mapsto t^{\phi(x_k)}$, and it does not depend on the presentation. By Remark~\ref{rem.convention}, $A_D(t)$ is the Alexander matrix of the presentation \eqref{eq.presentation}, evaluated under $x_e\mapsto t$ if $\eta_j=\ve_j$ and under $x_e\mapsto t^{-1}$ if $\eta_j=-\ve_j$.

\begin{lemma}[the Alexander polynomial]\label{lem.alexpoly}
Let $D$ be a diagram of $L$ with $c$ crossings and let $\hate$ be an edge of $D$. Then
\[
   P_{D,\hate}(t^{\ve})=\det A_{D,\hate}(t)\qquad\text{and}\qquad\tP_{D,\hate}(t^{\ve})=\det B_{D,\hate}(t)
\]
are both equal to $\Delta_L(t)$ up to multiplication by a unit $\pm t^m$ of $\BZ[t^{\pm1}]$. Consequently, $\LG(L;t,t^{-1})$ equals $\Delta_L(t)^2$ up to multiplication by a unit. More precisely, if $\Delta_L(t)=\nabla_L(t^{1/2}-t^{-1/2})$ is the Conway-normalized Alexander polynomial, where $\nabla_L$ is the Conway polynomial, then
\[
   P_{D,\hate}(t^{\ve})=\tP_{D,\hate}(t^{\ve})=(-1)^{\mu-1}\,t^{m_{D,\hate}/2}\,\Delta_L(t)\qquad\text{and}\qquad \LG(L;t,t^{-1})=\Delta_L(t)^2,
\]
where $m_{D,\hate}$ is the integer of Proposition~\ref{prop.statesum}~\textup{(i)}.
\end{lemma}
\begin{proof}
\emph{The statements up to a unit.} Let $A=A_D(t)$, a $2c\times2c$ matrix over $\BZ[t^{\pm1}]$ with $A\1=0$ and $\ell_tA=0$ (Lemma~\ref{lem.level} (ii)); in particular $\det A=0$. Let $\operatorname{adj}(A)$ be the adjugate matrix, whose entry $(r,s)$ is $(-1)^{r+s}$ times the determinant of $A$ with the row $s$ and the column $r$ removed, so that $A\operatorname{adj}(A)=\operatorname{adj}(A)A=\det(A)\cdot1=0$. If $A$ has rank at most $2c-2$ over $\BQ(t)$, all its $(2c-1)\times(2c-1)$ minors vanish; so $E_1(G)=0$, $\Delta_L=0$, and $P_{D,\hate}(t^{\ve})=0$. Suppose that $A$ has rank $2c-1$. Then its kernel is $\BQ(t)\1$ and its left kernel is $\BQ(t)\ell_t$. Every column of $\operatorname{adj}(A)$ lies in the kernel of $A$, so the column $s$ is $v_s\1$ with $v_s\in\BQ(t)$, and in fact $v_s\in\BZ[t^{\pm1}]$, being an entry of $\operatorname{adj}(A)$. Every row of $\operatorname{adj}(A)$ lies in the left kernel of $A$; since all rows are equal to $(v_s)_s$, we get $(v_s)_s=c_0\,\ell_t$ for some $c_0\in\BQ(t)$, and $c_0=v_st^{-h(s)}\in\BZ[t^{\pm1}]$. Thus every $(2c-1)\times(2c-1)$ minor of $A$ obtained by removing one row and one column is $c_0$ times a unit, these minors generate $E_1(G)$, hence $E_1(G)=(c_0)$ and $\Delta_L(t)$ is $c_0$ up to a unit, as is $\det A_{D,\hate}(t)$, the minor with the row and the column of $\hate$ removed. This proves the claim for $P_{D,\hate}(t^{\ve})$ if $\eta_j=\ve_j$; if $\eta_j=-\ve_j$, the same argument gives that $\det A_{D,\hate}(t)$ is $\Delta_L(t^{-1})$ up to a unit, and $\Delta_L(t^{-1})=\pm t^m\Delta_L(t)$ by the symmetry of the Alexander polynomial \cite{Torres}. The claim for $\tP_{D,\hate}(t^{\ve})$ follows from Lemma~\ref{lem.level} (i), the conjugation by the diagonal matrix $(\DD_t)_{\hate}$ not changing the determinant. Finally, for $t_0=t$ and $t_1=t^{-1}$ all the $g_j$ of \eqref{eq.xbar} vanish, so by Corollary~\ref{cor.expansion} only the term $S=\emptyset$ survives: $\LG(L;t,t^{-1})=(-1)^{\mu-1}P_{D,\hate}(t^{\ve})\,\tP_{D,\hate}(t^{-\ve})$, which is $\Delta_L(t)\Delta_L(t^{-1})$, i.e.\ $\Delta_L(t)^2$, up to a unit.

\emph{The units.} They are determined in Section~\ref{ssec.signs}, where only the statements up to a unit are used: the first identity is Proposition~\ref{prop.signs}~(iv), and the second one follows from it and from Theorem~\ref{thm.HK} at the free fermion point $t_0=t$, $t_1=t^{-1}$, where $\BE_D$ is the evaluation at $x=t^{\ve}$, $y=t^{-\ve}$, because $\Delta_L(t^{-1})=(-1)^{\mu-1}\Delta_L(t)$.
\end{proof}

The last statement of the lemma is the identity $\LG(L;t_0,t_0^{-1})=\Delta_L(t_0)^2$ of \cite{Kohli, KPM} which was recalled in the introduction; compare also \cite[Remark~5.9]{HK}.

\emph{The nullity.} Let $X=S^3\sm L$ be the exterior of $L$ (or, equivalently for all that follows, the compact manifold obtained by removing an open tubular neighborhood of $L$), with fundamental group $G$, and let $\phi\colon G\to\BZ$ be the total linking number homomorphism, which sends every meridian to $1$. Let $p\colon\widetilde X\to X$ be the infinite cyclic covering classified by $\phi$. The group of deck transformations of $\widetilde X$ is infinite cyclic, and we let its generator $t$ be the deck transformation corresponding to $1\in\BZ$; the homology groups $H_q(\widetilde X;\BQ)$ are then modules over the ring $\BQ[t^{\pm1}]$. Since $X$ has the homotopy type of a finite CW complex, $C_*(\widetilde X;\BQ)$ is a complex of finitely generated free $\BQ[t^{\pm1}]$-modules, and the $H_q(\widetilde X;\BQ)$ are finitely generated modules over the principal ideal domain $\BQ[t^{\pm1}]$. The module $H_1(\widetilde X;\BQ)$ is the \emph{Alexander module} of $L$ (with rational coefficients).

\begin{definition}\label{def.nullity}
The \emph{nullity} of $L$ is the rank of its Alexander module,
\[
   \b(L)=\rk_{\BQ[t^{\pm1}]}H_1(\widetilde X;\BQ)=\dim_{\BQ(t)}\bigl(H_1(\widetilde X;\BQ)\otimes_{\BQ[t^{\pm1}]}\BQ(t)\bigr).
\]
\end{definition}

The nullity is the number of copies of $\BQ[t^{\pm1}]$ in the decomposition of the Alexander module into cyclic modules; $\b(L)=0$ means that the Alexander module is a torsion module. The following lemma relates the nullity to the cohomology of the group $G$ with coefficients in $\BQ(t)_t$, and hence, by Proposition~\ref{prop.colorings}, to the colorings of any diagram of $L$.

\begin{lemma}\label{lem.nullitycoh}
$\b(L)=\dim_{\BQ(t)}H^1(G;\BQ(t)_t)=\dim_{\BQ(t)}H_1(X;\BQ(t)_t)$, where $\BQ(t)_t$ is the local system on $X$ defined by the action of $G$ on $\BQ(t)$ through $\phi$.
\end{lemma}
\begin{proof}
For a connected CW complex $X$ with fundamental group $G$ and a $G$-module $M$, the cohomology $H^1(X;M)$ with local coefficients is isomorphic to $H^1(G;M)$ \cite[Chapters~II and~III]{Brown}: the cellular chain complex $C_*(\widehat X)$ of the universal cover $\widehat X$ of $X$ is a complex of free $\BZ G$-modules in which $C_2\to C_1\to C_0\to\BZ\to0$ is exact, because $\widehat X$ is connected and simply connected; it can be completed to a free resolution of $\BZ$ by changing only the modules of degree $\ge3$, and $H^1(\Hom_{\BZ G}(-,M))$ only depends on the degrees $\le2$ of the resolution. Next, let $M=\BQ(t)_t$. Since the action of $G$ on $M$ factors through $\phi$, the cochains with local coefficients are $C^q(X;M)=\Hom_{\BZ G}(C_q(\widehat X),M)=\Hom_{\BQ[t^{\pm1}]}(C_q(\widetilde X;\BQ),\BQ(t))$, where $\widehat X$ is the universal cover and $C_q(\widetilde X)=\BZ[t^{\pm1}]\otimes_{\BZ G}C_q(\widehat X)$; see \cite[Section~3.H]{Hatcher}. The $\BQ[t^{\pm1}]$-module $\BQ(t)$ is divisible, hence injective, $\BQ[t^{\pm1}]$ being a principal ideal domain. So the functor $\Hom_{\BQ[t^{\pm1}]}(-,\BQ(t))$ is exact, and
\[
   H^1(X;M)\cong\Hom_{\BQ[t^{\pm1}]}\bigl(H_1(\widetilde X;\BQ),\BQ(t)\bigr).
\]
A finitely generated $\BQ[t^{\pm1}]$-module $N$ is the direct sum of a free module of rank $\rk N$ and of a torsion module, and $\Hom_{\BQ[t^{\pm1}]}(N,\BQ(t))\cong\BQ(t)^{\rk N}$, because $\BQ(t)$ is torsion free. Hence $\dim H^1(X;M)=\b(L)$. Finally, $H_1(X;M)$ is the homology of $C_*(\widetilde X;\BQ)\otimes_{\BQ[t^{\pm1}]}\BQ(t)$, and $\BQ(t)$ is a flat $\BQ[t^{\pm1}]$-module, so $H_1(X;M)\cong H_1(\widetilde X;\BQ)\otimes_{\BQ[t^{\pm1}]}\BQ(t)$, of dimension $\b(L)$.
\end{proof}

\begin{corollary}\label{cor.corank}
Let $D$ be a diagram of $L$ and $\hate$ an edge of $D$. Over $\BF=\BQ(t)$,
\[
   \dim\ker A_D(t)=1+\b(L),\qquad \corank A_{D,\hate}(t)=\corank B_{D,\hate}(t)=\b(L).
\]
In particular $\b(L)=0$ if and only if $P_{D,\hate}(t^{\ve})=\det A_{D,\hate}(t)\ne0$ for one, or equivalently for every, edge $\hate$; and this is the case if and only if the one-variable Alexander polynomial of $L$ is nonzero.
\end{corollary}
\begin{proof}
The formulas follow from Proposition~\ref{prop.colorings}, Remark~\ref{rem.convention}, Lemma~\ref{lem.cut} and Lemma~\ref{lem.nullitycoh}. For the last statement, $\b(L)=0$ if and only if the Alexander module is a torsion module, if and only if its order, which is the one-variable Alexander polynomial of $L$ \cite{Milnor}, is nonzero; recall also that $P_{D,\hate}(t^{\ve})$ is this polynomial up to multiplication by a unit of $\BZ[t^{\pm1}]$ (Lemma~\ref{lem.alexpoly}).
\end{proof}

\section{Flows on a diagram and the crossing forms}\label{sec.flows}

The edge polynomials $P_{D,\hate}(x)$ and $\tP_{D,\hate}(y)$ are defined as determinants of matrices of size $2c-1$ which are indexed by the edges of the diagram. In this section we show that they are determinants of $c\times c$ matrices of a very simple kind, and that these matrices are the matrices of one and the same object in different bases. The cut diagram carries a lattice of rank $c$, the lattice $\FF_{D,\hate}$ of its integral flows. Each crossing $j$ defines a linear form $\pi_j$ on this lattice, which contains the variable $x_j$ and no other variable, and $P_{D,\hate}(x)$ is, up to sign, the determinant of the $c$ forms $\pi_1,\dots,\pi_c$ in an arbitrary basis of the lattice; the same holds for $\tP_{D,\hate}(y)$ with forms $\tpi_1,\dots,\tpi_c$, and the two signs are related (Theorem~\ref{thm.flows}). No basis is preferred. The models of the Alexander polynomial which are deformed in this monograph correspond to particular bases: the boundaries of the faces of the diagram give Alexander's matrix (Theorem~\ref{thm.faceform}), and the Seifert circles, together with lifts of the cycles of the Seifert graph, give the Seifert matrices of Part~\ref{part.surfaces}. More generally there is a family of bases for each way of smoothing the crossings of the diagram (Proposition~\ref{prop.adapted}).

Throughout this section $D$ is a diagram in the sense of Section~\ref{ssec.diagrams}, with crossings $1,\dots,c$ and set of edges $E$, $|E|=2c$, and $\hate$ is an edge of $D$. We write
\[
   \Lam_c=\BZ[x_1^{\pm1},\dots,x_c^{\pm1}],\qquad \Lam'_c=\BZ[y_1^{\pm1},\dots,y_c^{\pm1}].
\]

\subsection{The lattice of flows}\label{ssec.flowlattice}
We write $z=(z_e)_{e\in E}$ for the elements of $\BZ^E$, and we think of $z_e$ as the amount of a quantity which flows along the edge $e$, in the direction of $e$ if $z_e>0$ and against it if $z_e<0$.

\begin{definition}[flows]\label{def.flows}
A \emph{flow} on $D$ is a vector $z\in\BZ^E$ such that, at every crossing, what flows in is what flows out:
\begin{equation}\label{eq.flow}
   z_{o_j}+z_{u_j}=z_{o'_j}+z_{u'_j}\qquad(j=1,\dots,c).
\end{equation}
The flows form a subgroup $\FF_D$ of $\BZ^E$, the \emph{lattice of flows} of $D$. For an edge $\hate$ of $D$,
\[
   \FF_{D,\hate}=\{z\in\FF_D:\ z_{\hate}=0\}
\]
is the lattice of the flows of the \emph{cut diagram}: the flows which do not use the edge $\hate$.
\end{definition}

In the language of graph theory, $D$ is a connected graph with $c$ vertices of degree $4$ and $2c$ oriented edges; the lattice $\FF_D$ is its cycle space $H_1(D;\BZ)$, and $\FF_{D,\hate}$ is the cycle space $H_1(D-\hate;\BZ)$ of the graph obtained by removing the edge $\hate$. Here are two families of flows which exist for every diagram; a third one, the circles of the smoothings of $D$, is introduced in Section~\ref{ssec.smoothings}.

\emph{The components of the link.} For a component $K$ of $L$, the vector $\sum_{e\subseteq K}e_e$, where $e_e\in\BZ^E$ is the unit vector of the edge $e$, is a flow: at a crossing, each of the two strands which belongs to $K$ contributes $1$ to both sides of \eqref{eq.flow}.

\emph{The boundaries of the faces.}\label{not.faces} The \emph{faces} of $D$ are the connected components of $S^2\sm D$; we write $\CF$ for the set of faces. Since $D$ is connected, every face is an open disk, and since $D$ has $c$ vertices and $2c$ edges, Euler's formula $c-2c+|\CF|=2$ gives $|\CF|=c+2$. For an edge $e$ let $f_{\rm L}(e)$ and $f_{\rm R}(e)$ be the faces on the left and on the right of $e$. They are distinct: the edges of $D$ are traversed by the closed curves which make up $L$, so every edge lies on a closed trail of the graph $D$, hence is not a bridge, and the two sides of a nonbridge edge of a connected graph on the sphere belong to different faces. The \emph{boundary} of a face $f$ is the vector
\begin{equation}\label{eq.boundaryface}
   \partial f=\sum_{e\in E}\Bigl(\ind{f_{\rm L}(e)=f}-\ind{f_{\rm R}(e)=f}\Bigr)\,e_e\ \in\BZ^E :
\end{equation}
it is the sum of the edges which lie on the boundary of $f$, with the sign $+$ for those which have $f$ on their left and the sign $-$ for those which have $f$ on their right. It is convenient to record the same numbers edge by edge: for an edge $e$ we put
\begin{equation}\label{eq.m}
   m(e)=e_{f_{\rm L}(e)}-e_{f_{\rm R}(e)}\in\BZ^{\CF},
\end{equation}
a column vector indexed by the faces, so that the coordinate of $\partial f$ at the edge $e$ is the coordinate of $m(e)$ at the face $f$. In the Dehn presentation of the link group, in which there is one generator for each face, $m(e)$ is the meridian of the edge $e$ written additively in the face generators \cite[Appendix]{KauffmanFKT}.

\begin{lemma}[faces and flows]\label{lem.faces}
\begin{enumerate}\renewcommand{\theenumi}{\roman{enumi}}
\item The boundary $\partial f$ of every face $f$ is a flow.
\item Every flow $z\in\FF_D$ can be written as $z=\sum_{f\in\CF}k_f\,\partial f$ with $k\in\BZ^{\CF}$, and $k$ is unique up to adding the same integer to all its coordinates. We call $k$ a \emph{potential} of $z$; in terms of $k$,
\begin{equation}\label{eq.potential}
   z_e=k_{f_{\rm L}(e)}-k_{f_{\rm R}(e)}=m(e)^{\T}k\qquad\text{for every edge }e .
\end{equation}
In particular $\FF_D$ is a free abelian group of rank $c+1$.
\item For every edge $\hate$, the boundaries of the $c$ faces which are not adjacent to $\hate$ form a $\BZ$-basis of $\FF_{D,\hate}$. In particular $\FF_{D,\hate}$ is free of rank $c$.
\end{enumerate}
\end{lemma}
\begin{proof}
The diagram $D$ is the $1$-skeleton of a cell decomposition of the sphere: the $0$-cells are the crossings, the $1$-cells are the edges, with their orientations, and the $2$-cells are the faces, which we orient as $S^2$. Its cellular chain complex is
\[
   \BZ^{\CF}\ \xrightarrow{\ \partial_2\ }\ \BZ^{E}\ \xrightarrow{\ \partial_1\ }\ \BZ^{c},
\]
where $\partial_1e_e$ is the crossing which $e$ enters minus the crossing which $e$ leaves, and $\partial_2e_f=\partial f$: the boundary of an oriented region runs along its boundary edges with the region on its left, and the two sides of an edge belong to different faces. The coordinate of $\partial_1z$ at the crossing $j$ is $z_{o_j}+z_{u_j}-z_{o'_j}-z_{u'_j}$, so that $\FF_D=\ker\partial_1$ by \eqref{eq.flow}. Since $H_1(S^2;\BZ)=0$, $\ker\partial_1=\im\partial_2$: this is (i) and the existence of $k$ in (ii). Since $H_2(S^2;\BZ)=\BZ$ is generated by the sum of all the faces, $\ker\partial_2$ consists of the constant vectors: this is the uniqueness in (ii). Formula \eqref{eq.potential} is \eqref{eq.boundaryface} read at the edge $e$, and the rank of $\FF_D\cong\BZ^{\CF}/\BZ\1$ is $(c+2)-1$.

(iii) The boundary of a face which is not adjacent to $\hate$ has the coordinate $0$ at $\hate$ by \eqref{eq.boundaryface}, so it lies in $\FF_{D,\hate}$. Let $z\in\FF_{D,\hate}$ and let $k$ be the potential of $z$ with $k_{f_{\rm L}(\hate)}=0$. By \eqref{eq.potential}, $0=z_{\hate}=k_{f_{\rm L}(\hate)}-k_{f_{\rm R}(\hate)}$, so $k$ vanishes on the two faces adjacent to $\hate$, and $z$ is a combination of the boundaries of the other faces. If a combination $\sum_fk_f\,\partial f$ of the boundaries of the faces not adjacent to $\hate$ vanishes, then the vector $k$, extended by $0$ to the two faces adjacent to $\hate$, is constant by (ii), hence zero.
\end{proof}

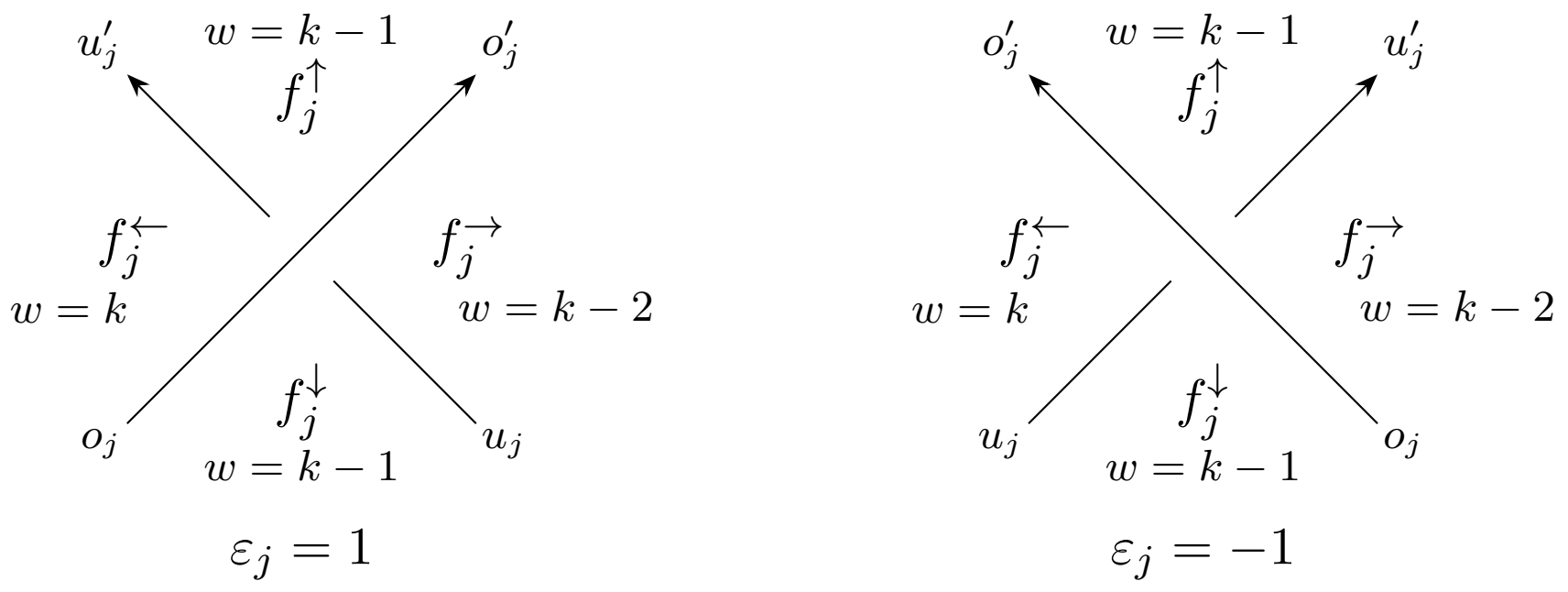
\begin{figure}[ht]
\centering
\begin{tikzpicture}[>=Stealth,scale=1.05,every node/.style={font=\small}]
\begin{scope}
\draw (1.2,-1.2) -- (0.22,-0.22); \draw[->] (-0.22,0.22) -- (-1.2,1.2);
\draw[->] (-1.2,-1.2) -- (1.2,1.2);
\node[below left,font=\scriptsize] at (-1.1,-1.1) {$o_j$}; \node[above right,font=\scriptsize] at (1.1,1.1) {$o'_j$};
\node[below right,font=\scriptsize] at (1.1,-1.1) {$u_j$}; \node[above left,font=\scriptsize] at (-1.1,1.1) {$u'_j$};
\node at (0,1.05) {$f^\uparrow_j$}; \node at (0,-1.05) {$f^\downarrow_j$};
\node at (-1.15,0) {$f^\leftarrow_j$}; \node at (1.15,0) {$f^\rightarrow_j$};
\node[font=\scriptsize] at (0,1.5) {$w=k-1$}; \node[font=\scriptsize] at (0,-1.5) {$w=k-1$};
\node[font=\scriptsize] at (-1.6,-0.4) {$w=k$}; \node[font=\scriptsize] at (1.75,-0.4) {$w=k-2$};
\node at (0,-2.1) {$\ve_j=1$};
\end{scope}
\begin{scope}[xshift=6.2cm]
\draw (-1.2,-1.2) -- (-0.22,-0.22); \draw[->] (0.22,0.22) -- (1.2,1.2);
\draw[->] (1.2,-1.2) -- (-1.2,1.2);
\node[below right,font=\scriptsize] at (1.1,-1.1) {$o_j$}; \node[above left,font=\scriptsize] at (-1.1,1.1) {$o'_j$};
\node[below left,font=\scriptsize] at (-1.1,-1.1) {$u_j$}; \node[above right,font=\scriptsize] at (1.1,1.1) {$u'_j$};
\node at (0,1.05) {$f^\uparrow_j$}; \node at (0,-1.05) {$f^\downarrow_j$};
\node at (-1.15,0) {$f^\leftarrow_j$}; \node at (1.15,0) {$f^\rightarrow_j$};
\node[font=\scriptsize] at (0,1.5) {$w=k-1$}; \node[font=\scriptsize] at (0,-1.5) {$w=k-1$};
\node[font=\scriptsize] at (-1.6,-0.4) {$w=k$}; \node[font=\scriptsize] at (1.75,-0.4) {$w=k-2$};
\node at (0,-2.1) {$\ve_j=-1$};
\end{scope}
\end{tikzpicture}
\caption{The four faces at a crossing, for the two signs, and their winding numbers $w$ (Lemma~\ref{lem.winding}): crossing an upward strand from right to left raises $w$ by one. A flow with potential $k$ has $z_e=k_{f_{\rm L}(e)}-k_{f_{\rm R}(e)}$; see \eqref{eq.mcrossing} and \eqref{eq.facetable}.}
\label{fig.faces}
\end{figure}

\emph{The four corners of a crossing.} At a crossing $j$, the four edges $o_j,u_j,o'_j,u'_j$ separate four corners. We call $f^\uparrow_j$ the face between the outgoing edges $o'_j,u'_j$, $f^\downarrow_j$ the face between the incoming edges $o_j,u_j$, and $f^\leftarrow_j$, $f^\rightarrow_j$ the faces on the left and on the right when the crossing is drawn with both strands pointing upward (Figure~\ref{fig.faces}); the four faces need not be distinct. For $\ve_j=1$ the over-strand comes in from the left, so $f^\leftarrow_j$ lies between $o_j$ and $u'_j$ and $f^\rightarrow_j$ between $u_j$ and $o'_j$; for $\ve_j=-1$ it comes in from the right, and $f^\leftarrow_j$ lies between $u_j$ and $o'_j$ and $f^\rightarrow_j$ between $o_j$ and $u'_j$. Writing $f$ for $e_f$, \eqref{eq.m} gives
\begin{equation}\label{eq.mcrossing}
\begin{aligned}
   \ve_j=1:&\quad m(o_j)=f^\leftarrow_j-f^\downarrow_j,\ \ m(u_j)=f^\downarrow_j-f^\rightarrow_j,\ \ m(o'_j)=f^\uparrow_j-f^\rightarrow_j,\ \ m(u'_j)=f^\leftarrow_j-f^\uparrow_j,\\
   \ve_j=-1:&\quad m(o_j)=f^\downarrow_j-f^\rightarrow_j,\ \ m(u_j)=f^\leftarrow_j-f^\downarrow_j,\ \ m(o'_j)=f^\leftarrow_j-f^\uparrow_j,\ \ m(u'_j)=f^\uparrow_j-f^\rightarrow_j.
\end{aligned}
\end{equation}
By \eqref{eq.potential} these formulas give the four numbers $z_{o_j},z_{u_j},z_{o'_j},z_{u'_j}$ of a flow $z$ in terms of the values of a potential at the four corners; for instance $z_{u'_j}=k(f^\leftarrow_j)-k(f^\uparrow_j)$ if $\ve_j=1$. In both cases
\begin{equation}\label{eq.vertexrelation}
   m(o_j)+m(u_j)=m(o'_j)+m(u'_j),
\end{equation}
which is the flow condition \eqref{eq.flow} for the boundaries of the faces, and also the relation of the Dehn presentation at the crossing $j$ (the product of the four face generators around the crossing, with alternating exponents, is trivial).

\subsection{The crossing forms and the main theorem}\label{ssec.crossingforms}
Each crossing of $D$ defines two linear forms on the lattice of flows, one for each flavor of the determinant formula. The main theorem of this section (Theorem~\ref{thm.flows}) says that the two edge polynomials are, up to a sign, the determinants of these forms.

\begin{definition}[crossing forms]\label{def.crossingforms}
For a crossing $j$ of $D$, the \emph{crossing forms} of the two flavors are the linear forms on $\FF_D$, with values in $\Lam_c$ and in $\Lam'_c$ respectively,
\begin{equation}\label{eq.crossingforms}
   \pi_j(z)=x_j\,z_{u'_j}-z_{u_j},\qquad \tpi_j(z)=y_j\,z_{o_j}-z_{o'_j} .
\end{equation}
We write $\Pi(x)=(\pi_1,\dots,\pi_c)\colon\FF_D\otimes\Lam_c\to\Lam_c^c$ and $\tilde\Pi(y)=(\tpi_1,\dots,\tpi_c)\colon\FF_D\otimes\Lam'_c\to\Lam_c'^c$ for the maps which they define.
\end{definition}

Thus $\pi_j$ looks at the under-strand of the crossing $j$: it is the flow through the outgoing under-edge, multiplied by the variable of the crossing, minus the flow through the incoming under-edge. The form $\tpi_j$ looks at the over-strand in the same way, with the roles of the incoming and of the outgoing edge exchanged. When the variables are set to $1$ the two forms coincide: by \eqref{eq.flow},
\begin{equation}\label{eq.piatone}
   \pi_j(z)\big|_{x_j=1}=z_{u'_j}-z_{u_j}=z_{o_j}-z_{o'_j}=\tpi_j(z)\big|_{y_j=1},
\end{equation}
which is the amount of flow which passes from the over-strand to the under-strand at the crossing $j$.

The forms $\pi_j$ arise as follows. Read the matrix $1-T^a_D(x)$ as a presentation matrix, with one generator $[e]$ for each edge: by \eqref{eq.Ta}, the row of the edge $o'_j$ is the relation $[o'_j]=[o_j]$ and the row of $u'_j$ is the relation $[u'_j]=(1-x_j)[o_j]+x_j[u_j]$. Put $\d_j=[u_j]-[o_j]$\label{not.delta}, the jump of the generator from the over-strand to the under-strand at the crossing $j$.

\begin{lemma}[the relation of a flow]\label{lem.walk}
Let $\CA_{\hate}$ be the $\Lam_c$-module with generators $[e]$, $e\in E$, and relations $[\hate]=0$ and the rows of $1-T^a_D(x)$ other than the row $\hate$. Then for every flow $z\in\FF_{D,\hate}$,
\[
   \sum_{j=1}^c\pi_j(z)\,\d_j=0\qquad\text{in }\CA_{\hate} .
\]
\end{lemma}
\begin{proof}
Put $p_j=[o_j]$, so that $[u_j]=p_j+\d_j$. Let $e\ne\hate$ be an edge; it leaves a crossing $k$ and enters a crossing $l$. The row $e$ of the presentation says that $[e]=p_k$ if $e=o'_k$, and that $[e]=(1-x_k)[o_k]+x_k[u_k]=p_k+x_k\d_k$ if $e=u'_k$. On the other hand $[e]=p_l$ if $e=o_l$ and $[e]=p_l+\d_l$ if $e=u_l$. Comparing the two expressions of $[e]$,
\[
   p_l-p_k=\ind{e=u'_k}\,x_k\,\d_k-\ind{e=u_l}\,\d_l .
\]
Multiply this relation by $z_e$ and sum over the edges $e\ne\hate$; since $z_{\hate}=0$ we may as well sum over all the edges. On the left-hand side the coefficient of $p_j$ is the flow of $z$ into the crossing $j$ minus the flow out of it, which is zero by \eqref{eq.flow}. On the right-hand side the coefficient of $\d_j$ is $x_jz_{u'_j}-z_{u_j}=\pi_j(z)$.
\end{proof}

So every flow of the cut diagram is a closed walk along which the jumps $\d_j$ add up to zero, the jump at the crossing $j$ being counted $\pi_j(z)$ times. The $c$ flows of a basis of $\FF_{D,\hate}$ give $c$ relations between the $c$ elements $\d_1,\dots,\d_c$, and the main result of this section is that the determinant of these relations is the edge polynomial.

\begin{theorem}[the edge polynomials on the lattice of flows]\label{thm.flows}
Let $D$ be a diagram of a link with $\mu$ components, $\hate$ an edge of $D$ and $z_1,\dots,z_c$ a $\BZ$-basis of the lattice $\FF_{D,\hate}$. Then there is a sign $\kappa=\pm1$ such that
\[
   \det\bigl(\pi_j(z_i)\bigr)_{i,j=1,\dots,c}=\kappa\,P_{D,\hate}(x)\qquad\text{and}\qquad \det\bigl(\tpi_j(z_i)\bigr)_{i,j=1,\dots,c}=(-1)^{\mu-1}\,\kappa\,\tP_{D,\hate}(y) .
\]
The sign $\kappa$ depends on the basis and on the numbering of the crossings, and it is the same in the two identities. Consequently
\begin{equation}\label{eq.LGflows}
   \LG(L;t_0,t_1)=\BE_D\Bigl[\det\bigl(\pi_j(z_i)\bigr)\,\det\bigl(\tpi_j(z_i)\bigr)\Bigr]
\end{equation}
for every basis of $\FF_{D,\hate}$, with no sign.
\end{theorem}

The theorem is proved in Section~\ref{ssec.elimination}. We first make some remarks on its statement.

\begin{remark}\label{rem.flows}
(i) \emph{Change of basis.} If $z'_1,\dots,z'_c$ is another basis of $\FF_{D,\hate}$, then $z'_i=\sum_lU_{il}z_l$ for a matrix $U\in\GL_c(\BZ)$, and the two matrices $\bigl(\pi_j(z'_i)\bigr)$ and $\bigl(\tpi_j(z'_i)\bigr)$ are $U$ times the matrices of the basis $z_1,\dots,z_c$. So the two determinants are multiplied by the same sign $\det U$: it suffices to prove the theorem for one basis, and the product of the two determinants in \eqref{eq.LGflows} does not depend on the basis. A change of the numbering of the crossings permutes the columns of the two matrices in the same way.

(ii) \emph{The variables sit on the columns.} The matrix of the theorem is
\[
   \bigl(\pi_j(z_i)\bigr)_{i,j}=Z_{u'}X-Z_u,\qquad X=\diag(x_1,\dots,x_c),
\]
where $Z_{u'}$ and $Z_u$ are the integer matrices with entries $(z_i)_{u'_j}$ and $(z_i)_{u_j}$: the variable $x_j$ occurs in the column $j$ only, where it multiplies a column of integers. This makes it visible that $P_{D,\hate}$ has degree at most one in each variable. Likewise $\bigl(\tpi_j(z_i)\bigr)=Z_o\,\diag(y_1,\dots,y_c)-Z_{o'}$.

(iii) \emph{The coupled form.} Let $\tau=\sum_{j=1}^c\ve_j\,\pi_j\otimes\tpi_j$, a bilinear form on $\FF_{D,\hate}$ with values in $\BZ[x,y]$. Its matrix in the basis $z_1,\dots,z_c$ is the product of the matrix $\bigl(\pi_j(z_i)\bigr)$, of $\diag(\ve_1,\dots,\ve_c)$ and of the transpose of $\bigl(\tpi_j(z_i)\bigr)$, so that \eqref{eq.LGflows} reads
\begin{equation}\label{eq.LGtau}
   \LG(L;t_0,t_1)=\ve_D\,\BE_D\Bigl[\det\bigl(\tau(z_i,z_k)\bigr)_{i,k=1,\dots,c}\Bigr],\qquad \ve_D=\prod_{j=1}^c\ve_j :
\end{equation}
up to the sign $\ve_D$, the Links--Gould polynomial is the image under $\BE_D$ of the discriminant of the form $\tau$ on the lattice of flows of the cut diagram. In the basis of the faces this is the deformed Kashaev matrix of Section~\ref{ssec.deformedkashaev}. The map $\BE_D$ can also be applied to $\tau$ term by term, before the determinant is taken; this turns the discriminant into a mixed discriminant and removes the sign $\ve_D$ (Section~\ref{ssec.mixed}).
\end{remark}

\subsection{Smoothings, circles and adapted bases}\label{ssec.smoothings}
The Seifert circles of Section~\ref{ssec.level} are the circles of one particular smoothing of the diagram, in which every crossing is smoothed along the orientations. In this subsection we consider all the $2^c$ smoothings at once. Each of them gives a family of bases of the lattice of flows of the cut diagram, which contain all the circles of the smoothing but one (Proposition~\ref{prop.adapted}). The bases of the oriented smoothing serve to prove the main theorem (Section~\ref{ssec.elimination}) and give the deformed Seifert matrices of Section~\ref{sec.seifert}; the bases of the other smoothings belong to the state surfaces of Section~\ref{sec.statesurfaces}.

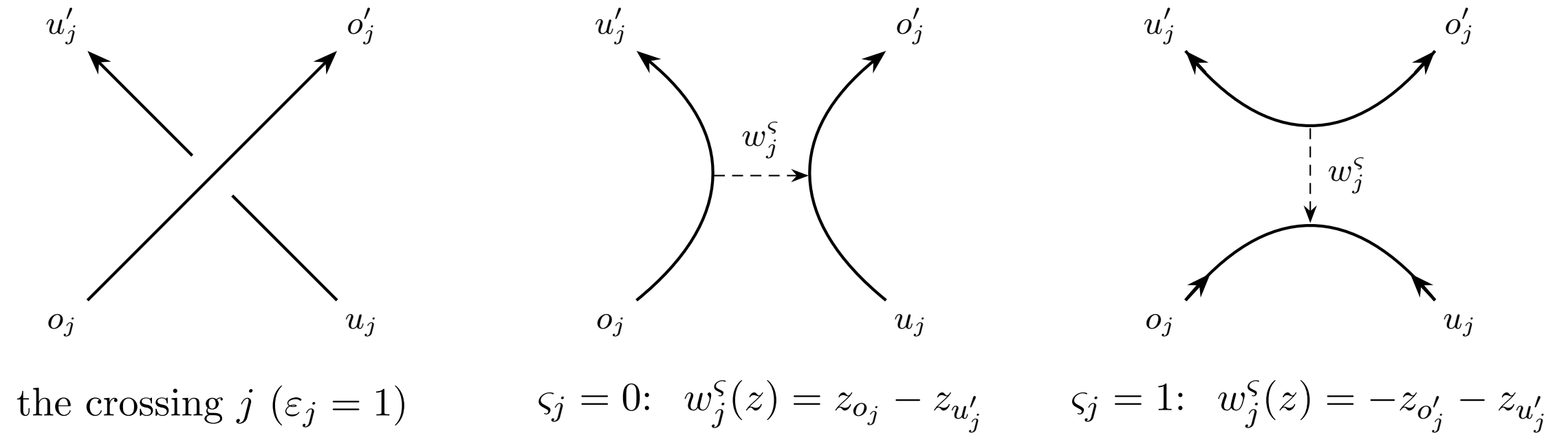
\begin{figure}[ht]
\centering
\begin{tikzpicture}[scale=1.2,>=Stealth,every node/.style={font=\small}]
\def\lab{%
  \node[font=\scriptsize] at (-1.2,-1.2) {$o_j$}; \node[font=\scriptsize] at (1.2,-1.2) {$u_j$}; \node[font=\scriptsize] at (1.2,1.2) {$o'_j$}; \node[font=\scriptsize] at (-1.2,1.2) {$u'_j$};}
\begin{scope}
  \draw[thick] (1,-1) -- (0.16,-0.16); \draw[thick,->] (-0.16,0.16) -- (-1,1);
  \draw[thick,->] (-1,-1) -- (1,1);
  \lab
  \node at (0,-1.85) {the crossing $j$ ($\ve_j=1$)};
\end{scope}
\begin{scope}[xshift=4.4cm]
  \draw[thick,->] (-1,-1) .. controls (-0.2,-0.35) and (-0.2,0.35) .. (-1,1);
  \draw[thick,->] (1,-1) .. controls (0.2,-0.35) and (0.2,0.35) .. (1,1);
  \draw[->,densely dashed] (-0.38,0) -- (0.38,0);
  \node[above,font=\scriptsize] at (0,0.02) {$w^\vs_j$};
  \lab
  \node at (0,-1.85) {$\vs_j=0$: \ $w^\vs_j(z)=z_{o_j}-z_{u'_j}$};
\end{scope}
\begin{scope}[xshift=8.8cm]
  \draw[thick] (-1,-1) .. controls (-0.35,-0.2) and (0.35,-0.2) .. (1,-1);
  \draw[thick] (-1,1) .. controls (-0.35,0.2) and (0.35,0.2) .. (1,1);
  \draw[->,thick] (-1,-1) -- (-0.80,-0.785); \draw[->,thick] (1,-1) -- (0.80,-0.785);
  \draw[->,thick] (0.84,0.825) -- (1,1); \draw[->,thick] (-0.84,0.825) -- (-1,1);
  \draw[->,densely dashed] (0,0.38) -- (0,-0.38);
  \node[right,font=\scriptsize] at (0.02,0) {$w^\vs_j$};
  \lab
  \node at (0,-1.85) {$\vs_j=1$: \ $w^\vs_j(z)=-z_{o'_j}-z_{u'_j}$};
\end{scope}
\end{tikzpicture}
\caption{A positive crossing and its two smoothings: the oriented smoothing ($\vs_j=0$), whose arcs join $o_j$ to $u'_j$ and $u_j$ to $o'_j$, and the other smoothing ($\vs_j=1$), whose arcs join $o_j$ to $u_j$ and $o'_j$ to $u'_j$. The dashed arrow is the direction in which the transfer $w^\vs_j$ of Definition~\ref{def.transfer} counts the flow which passes from one arc to the other: from the arc through $u'_j$ to the arc through $u_j$. For a negative crossing the picture is the same with the labels of Figure~\ref{fig.crossing}.}
\label{fig.smoothings}
\end{figure}

A crossing of an oriented diagram can be smoothed in two ways (Figure~\ref{fig.smoothings}). The \emph{oriented smoothing} replaces the crossing $j$ by two arcs which respect the orientations: one arc joins $o_j$ to $u'_j$ and the other one joins $u_j$ to $o'_j$. The other smoothing joins the two incoming edges $o_j,u_j$ by an arc and the two outgoing edges $o'_j,u'_j$ by another arc.

\begin{definition}[smoothings and their circles]\label{def.smoothing}
A \emph{smoothing} of $D$ is a vector $\vs\in\{0,1\}^c$: the crossing $j$ receives the oriented smoothing if $\vs_j=0$ and the other one if $\vs_j=1$. Smoothing all the crossings according to $\vs$ turns $D$ into a disjoint union of simple closed curves on $S^2$, the \emph{circles of $\vs$}; every edge of $D$ lies on exactly one of them. We write $s_\vs$ for the number of circles of $\vs$. The circles of the smoothing $\vs=0$ are the Seifert circles of $D$ (Section~\ref{ssec.level}), and $s_0=s$.
\end{definition}

\emph{Circles as flows.}\label{not.circleflows} Let $C$ be a circle of the smoothing $\vs$. Choose one of the two directions in which $C$ can be traveled, and let $z^C\in\BZ^E$ be the vector whose coordinate at an edge $e$ is $1$ if $e$ lies on $C$ and is traveled in its own direction, $-1$ if $e$ lies on $C$ and is traveled against its direction, and $0$ if $e$ does not lie on $C$. Then $z^C$ is a flow. Indeed, when $C$ passes through a crossing along an arc of the oriented smoothing, it travels an incoming and an outgoing edge of the crossing, both in their own direction or both against it, and this contributes the same number to the two sides of \eqref{eq.flow}; when it passes along an arc of the other smoothing, it travels two incoming edges, or two outgoing edges, one of them in its own direction and the other one against it, and this contributes $1-1=0$ to one side of \eqref{eq.flow} and nothing to the other side. The flow $z^C$ is determined by $C$ up to sign. For the Seifert circles we always travel in the direction of the link, so that $z^C=\sum_{e\subseteq C}e_e$; we then simply write $C$ for $z^C$.

\emph{Transfers.} A flow which arrives at a crossing along one of the two arcs of a smoothing may leave it along the other one. The following linear forms measure how much of it does.

\begin{definition}[transfers]\label{def.transfer}
Let $\vs$ be a smoothing of $D$. The \emph{transfer} of a flow $z\in\FF_D$ at the crossing $j$ is
\[
   w^\vs_j(z)=\begin{cases} z_{o_j}-z_{u'_j}&\text{if }\vs_j=0,\\ -z_{o'_j}-z_{u'_j}&\text{if }\vs_j=1,\end{cases}
\]
and $w^\vs=(w^\vs_1,\dots,w^\vs_c)\colon\FF_D\to\BZ^c$. For the smoothing $\vs=0$ we write $w_j$ and $w$ for $w^0_j$ and $w^0$.
\end{definition}

In both cases $w^\vs_j(z)$ is the amount of flow which passes, at the crossing $j$, from the arc through $u'_j$ to the arc through $u_j$ (Figure~\ref{fig.smoothings}). If $\vs_j=0$, the flow $z_{o_j}$ arrives along the arc which joins $o_j$ to $u'_j$ and the flow $z_{u'_j}$ leaves along it; the difference has passed to the other arc, and by \eqref{eq.flow} it is also $z_{o'_j}-z_{u_j}$, the flow which the other arc has gained. If $\vs_j=1$, the flow $z_{o'_j}+z_{u'_j}$ leaves the crossing along the two ends of the arc through $u'_j$; all of it has come from the arc through $u_j$, along whose two ends $z_{o_j}+z_{u_j}$ arrives; so the flow from the arc through $u'_j$ to the arc through $u_j$ is $-(z_{o'_j}+z_{u'_j})$. With the two linear forms
\begin{equation}\label{eq.UO}
   U_j(z)=z_{u'_j}+z_{u_j},\qquad O_j(z)=z_{o_j}+z_{o'_j},
\end{equation}
which are twice the flow along the under-strand and along the over-strand of the crossing $j$, averaged over the incoming and the outgoing edge, the transfers of both kinds are, by \eqref{eq.flow},
\begin{equation}\label{eq.transferUO}
   w^\vs_j=\tfrac12\bigl((-1)^{\vs_j}\,O_j-U_j\bigr).
\end{equation}

\begin{lemma}[circles and transfers]\label{lem.circles}
Let $\vs$ be a smoothing of $D$ and let $z^C$ be the flows of its circles, for some choice of the directions of travel.
\begin{enumerate}\renewcommand{\theenumi}{\roman{enumi}}
\item The flows $z^C$ are linearly independent, and a flow $z\in\FF_D$ is an integral combination of them if and only if $w^\vs(z)=0$.
\item The quotient $\FF_D/\langle z^C\rangle$ of the lattice of flows by the subgroup generated by the circles of $\vs$ is a free abelian group of rank $c+1-s_\vs$.
\end{enumerate}
\end{lemma}
\begin{proof}
(i) The flows $z^C$ have pairwise disjoint supports, since every edge lies on exactly one circle; so they are linearly independent. Let $C$ be a circle and $j$ a crossing with $\vs_j=0$. The form $w^\vs_j$ involves the two edges $o_j$ and $u'_j$, which are joined by an arc of the smoothing. If $C$ contains this arc, then $z^C_{o_j}=z^C_{u'_j}=\pm1$, and otherwise $z^C_{o_j}=z^C_{u'_j}=0$; in both cases $w^\vs_j(z^C)=0$. If $\vs_j=1$, the form $w^\vs_j$ involves the two outgoing edges $o'_j$ and $u'_j$, which are joined by an arc; if $C$ contains this arc it travels one of the two edges in its direction and the other one against it, so $z^C_{o'_j}+z^C_{u'_j}=0$, and otherwise both numbers vanish. Hence $w^\vs(z^C)=0$ for every circle. Conversely let $z$ be a flow with $w^\vs(z)=0$. At a crossing with $\vs_j=0$ this says that $z_{o_j}=z_{u'_j}$, and then $z_{u_j}=z_{o'_j}$ by \eqref{eq.flow}; at a crossing with $\vs_j=1$ it says that $z_{o'_j}=-z_{u'_j}$, and then $z_{o_j}=-z_{u_j}$ by \eqref{eq.flow}. In all cases, for the two edges $e,e'$ joined by an arc of the smoothing, $z_e=z_{e'}$ if they are traveled in the same sense by the circle through the arc, and $z_e=-z_{e'}$ otherwise; that is, $z^C_ez_e=z^C_{e'}z_{e'}$. So the number $a_C=z^C_ez_e$ is the same for all the edges $e$ of a circle $C$, and $z=\sum_Ca_Cz^C$.

(ii) By (i) the map $w^\vs$ induces an injection of $\FF_D/\langle z^C\rangle$ into $\BZ^c$, so the quotient is free; its rank is $(c+1)-s_\vs$ by Lemma~\ref{lem.faces}~(ii) and by the independence of the $s_\vs$ flows $z^C$.
\end{proof}

In Section~\ref{sec.statesurfaces} we shall see that $\FF_D/\langle z^C\rangle$ is the first homology group of a surface: the state surface $F_\vs$ of the smoothing, which consists of disks bounded by the circles of $\vs$ and of a band at each crossing, the transfer $w^\vs_j(z)$ being the number of times the cycle $z$ runs through the band of the crossing $j$.

\begin{proposition}[bases adapted to a smoothing]\label{prop.adapted}
Let $\vs$ be a smoothing of $D$, let $\hat C$ be the circle of $\vs$ which contains the edge $\hate$, and let $\nu_\vs=c+1-s_\vs$. Let $z_1,\dots,z_{\nu_\vs}\in\FF_{D,\hate}$ be flows whose classes form a basis of $\FF_D/\langle z^C\rangle$. Then the flows $z_1,\dots,z_{\nu_\vs}$ together with the $s_\vs-1$ flows $z^C$ of the circles $C\ne\hat C$ form a $\BZ$-basis of $\FF_{D,\hate}$. Such flows $z_i$ exist: every class of $\FF_D/\langle z^C\rangle$ contains a flow which vanishes at $\hate$.
\end{proposition}
\begin{proof}
The coordinate of $z^{\hat C}$ at $\hate$ is $\pm1$, and the other circles do not contain $\hate$. So a flow $z$ and the flow $z\mp z_{\hate}z^{\hat C}$, which vanishes at $\hate$, have the same class; this proves the last statement. Let $z\in\FF_{D,\hate}$. Its class is a combination $\sum_in_i[z_i]$, so $z-\sum_in_iz_i=\sum_Ca_Cz^C$ for some integers $a_C$. The coordinate at $\hate$ of the left-hand side is $0$ and that of the right-hand side is $\pm a_{\hat C}$; hence $a_{\hat C}=0$, and $z$ is a combination of the $z_i$ and of the $z^C$ with $C\ne\hat C$. If $\sum_in_iz_i+\sum_{C\ne\hat C}a_Cz^C=0$, then $\sum_in_i[z_i]=0$, so all the $n_i$ vanish, and then all the $a_C$ vanish by Lemma~\ref{lem.circles}~(i).
\end{proof}

By Theorem~\ref{thm.flows} there is thus one determinant formula for the edge polynomials for each smoothing of the diagram: the rows of the matrix $\bigl(\pi_j(z_i)\bigr)$ are of two kinds, the rows of the circles of the smoothing and the rows of the lifted classes. We now make this explicit for the oriented smoothing, which is the case that we use to prove the theorem; the basis of the faces is recovered in Remark~\ref{rem.checkerboardbasis}.

\emph{The Seifert graph.}\label{not.seifertgraph} Let $C_1,\dots,C_s$ be the Seifert circles of $D$. At every crossing $j$ the edges $o_j$ and $u'_j$, which are joined by an arc of the oriented smoothing, lie on one Seifert circle, which we call $A_j$, and the edges $u_j$ and $o'_j$ lie on another one, which we call $B_j$ (Figure~\ref{fig.crossing}; the two circles are distinct by Section~\ref{ssec.level}). The \emph{Seifert graph} $\Gamma_D$ is the graph with vertices $C_1,\dots,C_s$ and one edge for each crossing $j$, joining $A_j$ to $B_j$; we orient this edge from $A_j$ to $B_j$. Since $D$ is connected, $\Gamma_D$ is connected. Let
\[
   \partial\colon\BZ^c\to\BZ^s,\qquad \partial e_j=e_{B_j}-e_{A_j},
\]
be its incidence map, and let $H=\ker\partial\subseteq\BZ^c$ be the \emph{cycle space} of $\Gamma_D$, a free abelian group of rank
\[
   \nu=c-s+1 .
\]
For $\gamma\in H$ we write $w_j(\gamma)$ for its coordinates: $\gamma=\sum_jw_j(\gamma)e_j$ runs $w_j(\gamma)$ times through the edge $j$ of $\Gamma_D$, from $A_j$ to $B_j$. For a Seifert circle $C$ and a crossing $j$ on $C$ we put
\[
   s_C(j)=1\ \text{ if }A_j=C,\qquad s_C(j)=-1\ \text{ if }B_j=C,
\]
so that $\partial\gamma=0$ means that
\begin{equation}\label{eq.sw}
   \sum_{j\text{ on }C}s_C(j)\,w_j(\gamma)=0\qquad\text{for every Seifert circle }C:
\end{equation}
the cycle $\gamma$ leaves the vertex $C$ as often as it enters it.

We now walk along the Seifert circles. Every edge $e$ of $D$ lies on exactly one Seifert circle $C$ and enters exactly one crossing $j$, \emph{the crossing at the end of $e$}: either $e=o_j$, in which case $C=A_j$ and the next edge on $C$ is $u'_j$, or $e=u_j$, in which case $C=B_j$ and the next edge on $C$ is $o'_j$. Conversely a crossing $j$ lies on exactly two Seifert circles, $A_j$ (at the end of $o_j$) and $B_j$ (at the end of $u_j$), so a crossing is met exactly once when going around a circle which contains it. We choose a \emph{base edge} $b_C$ on every Seifert circle $C$, with the only requirement that
\begin{equation}\label{eq.baseedge}
   b_{\hat C}=\hate,\qquad\text{where $\hat C$ is the Seifert circle containing $\hate$.}
\end{equation}
For an edge $e$ on the circle $C$, let $\mathrm{Pre}_C(e)$ be the set of crossings at the ends of the edges of $C$ which are met when going along $C$, in the direction of its orientation, from $b_C$ to $e$, not counting the crossing at the end of $e$ itself; thus $\mathrm{Pre}_C(b_C)=\emptyset$. For two crossings $j,k$ on $C$ we write
\[
   k\prec_Cj\qquad\text{if $k\in\mathrm{Pre}_C(e)$, where $e$ is the edge of $C$ which ends at $j$,}
\]
that is, if $k$ is met before $j$ when going along $C$ from its base edge; $\prec_C$ is a total order on the crossings of $C$.

\begin{lemma}[the cycles of the Seifert graph as flows]\label{lem.lifts}
\begin{enumerate}\renewcommand{\theenumi}{\roman{enumi}}
\item For every flow $z\in\FF_D$ the vector $w(z)=\bigl(z_{o_j}-z_{u'_j}\bigr)_j$ of its transfers for the oriented smoothing is a cycle of the Seifert graph, and $w\colon\FF_D\to H$ is onto, with kernel the subgroup generated by the Seifert circles. Thus $w$ induces an isomorphism of $\FF_D/\langle C_1,\dots,C_s\rangle$ onto $H$.
\item For every cycle $\gamma\in H$ there is exactly one flow $z^\gamma$ with $w(z^\gamma)=\gamma$ which vanishes on all the base edges, the \emph{lift} of $\gamma$. On an edge $e$ of the Seifert circle $C$,
\begin{equation}\label{eq.lift}
   z^\gamma_e=-\sum_{k\in\mathrm{Pre}_C(e)}s_C(k)\,w_k(\gamma) .
\end{equation}
In particular $z^\gamma\in\FF_{D,\hate}$, and for every crossing $k$
\begin{equation}\label{eq.gh}
   z^\gamma_{u'_k}=\sum_{j:\,k\prec_{A_k}j}s_{A_k}(j)\,w_j(\gamma),\qquad
   z^\gamma_{u_k}=-\sum_{j:\,j\prec_{B_k}k}s_{B_k}(j)\,w_j(\gamma) .
\end{equation}
\item If $\gamma_1,\dots,\gamma_\nu$ is a basis of $H$, then the lifts $z^{\gamma_1},\dots,z^{\gamma_\nu}$ and the Seifert circles $C\ne\hat C$ form a $\BZ$-basis of $\FF_{D,\hate}$, which we call a \emph{Seifert basis}.
\end{enumerate}
\end{lemma}
\begin{proof}
Let $z$ be a flow and put $w_j=w_j(z)$. Going along a Seifert circle $C$ and passing a crossing $j$ with $A_j=C$, one goes from the edge $o_j$ to the edge $u'_j$, and $z_{u'_j}=z_{o_j}-w_j$ by the definition of the transfer; passing a crossing $j$ with $B_j=C$, one goes from $u_j$ to $o'_j$, and $z_{o'_j}=z_{u_j}+w_j$ by \eqref{eq.flow}. So
\begin{equation}\label{eq.walkcircle}
   \text{along $C$, the flow $z_e$ changes by $-s_C(j)\,w_j$ at each crossing $j$ of $C$.}
\end{equation}
After a full turn the flow is back to its value, hence $\sum_{j\text{ on }C}s_C(j)w_j=0$ for every $C$, which says that $w(z)\in H$ by \eqref{eq.sw}. The kernel of $w$ is generated by the Seifert circles by Lemma~\ref{lem.circles}~(i). Conversely let $\gamma\in H$, and define $z^\gamma$ by \eqref{eq.lift}. By the definition of $\mathrm{Pre}_C$, the vector $z^\gamma$ vanishes on the base edges and satisfies \eqref{eq.walkcircle} with $w_j=w_j(\gamma)$ at every crossing of $C$ except possibly at the last one, the crossing at which one comes back to $b_C$; and there it satisfies it as well, because the sum of all the changes is zero by \eqref{eq.sw}. At a crossing $j$, \eqref{eq.walkcircle} for the two circles $A_j$ and $B_j$ gives $z^\gamma_{u'_j}=z^\gamma_{o_j}-w_j(\gamma)$ and $z^\gamma_{o'_j}=z^\gamma_{u_j}+w_j(\gamma)$; adding the two, $z^\gamma$ satisfies \eqref{eq.flow}, so it is a flow, and the first identity says that $w_j(z^\gamma)=w_j(\gamma)$. This proves that $w$ is onto, which completes (i), and the existence in (ii). A flow with given transfers which vanishes on the base edges is determined on every circle by \eqref{eq.walkcircle}, which is the uniqueness; and $z^\gamma\in\FF_{D,\hate}$ because $\hate$ is a base edge. For \eqref{eq.gh}: the edge $u'_k$ follows the crossing $k$ on the circle $A_k$, so if $u'_k$ is not the base edge, $\mathrm{Pre}_{A_k}(u'_k)$ consists of $k$ and of the crossings $j\prec_{A_k}k$, and \eqref{eq.lift} and \eqref{eq.sw} give $z^\gamma_{u'_k}=-\sum_{j\preceq k}s_{A_k}(j)w_j(\gamma)=\sum_{j\succ k}s_{A_k}(j)w_j(\gamma)$; if $u'_k$ is the base edge of $A_k$, then $k$ is the last crossing of $A_k$ and both sides vanish. The edge $u_k$ ends at the crossing $k$ on the circle $B_k$, so $\mathrm{Pre}_{B_k}(u_k)=\{j:j\prec_{B_k}k\}$.

(iii) By (i) the classes of the lifts form a basis of $\FF_D/\langle C_1,\dots,C_s\rangle$, and the lifts lie in $\FF_{D,\hate}$; so this is Proposition~\ref{prop.adapted} for the smoothing $\vs=0$.
\end{proof}

By \eqref{eq.gh}, $z^\gamma_{u'_k}$ is the flow of $\gamma$ through the bands of the circle $A_k$ which are attached after the band $k$, between $k$ and the base edge, and $-z^\gamma_{u_k}$ is the flow through the bands of $B_k$ which are attached before $k$. We fix a $c\times\nu$ integer matrix $Y$\label{not.Y} whose columns $\gamma_1,\dots,\gamma_\nu$ form a basis of $H$. A convenient choice is the basis of \emph{fundamental cycles} of a spanning tree $\CT$ of $\Gamma_D$: for every edge $j\notin\CT$, the cycle consisting of $j$ and of the path in $\CT$ from $B_j$ back to $A_j$. Let $E_{\CT}$ be the $c\times(s-1)$ matrix with columns $e_j$, $j\in\CT$. Then the $c\times c$ matrix $(Y\,|\,E_{\CT})$ is unimodular: for the fundamental cycles, after ordering the edges of $\Gamma_D$ with the edges not in $\CT$ first, it has the form $\bigl(\begin{smallmatrix}1&0\\ *&1\end{smallmatrix}\bigr)$; and for another basis $Y'=YQ$ with $Q\in\GL_\nu(\BZ)$, so that $(Y'\,|\,E_{\CT})=(Y\,|\,E_{\CT})\bigl(\begin{smallmatrix}Q&0\\0&1\end{smallmatrix}\bigr)$.

\subsection{Elimination along the Seifert circles}\label{ssec.elimination}
We now prove Theorem~\ref{thm.flows}. It suffices to prove it for one basis, by Remark~\ref{rem.flows}~(i), and we take a Seifert basis (Lemma~\ref{lem.lifts}~(iii)). We compute $P_{D,\hate}(x)$ and $\tP_{D,\hate}(y)$ by eliminating the generators of the edges along the Seifert circles (Propositions~\ref{prop.seifert} and~\ref{prop.seifertb}), and we then recognize the crossing forms in the result (Lemma~\ref{lem.rowsareflows}). We keep the notation of Section~\ref{ssec.smoothings}: the Seifert circles $C$ with their base edges $b_C$, the sets $\mathrm{Pre}_C(e)$, the orders $\prec_C$ and the signs $s_C(j)$.

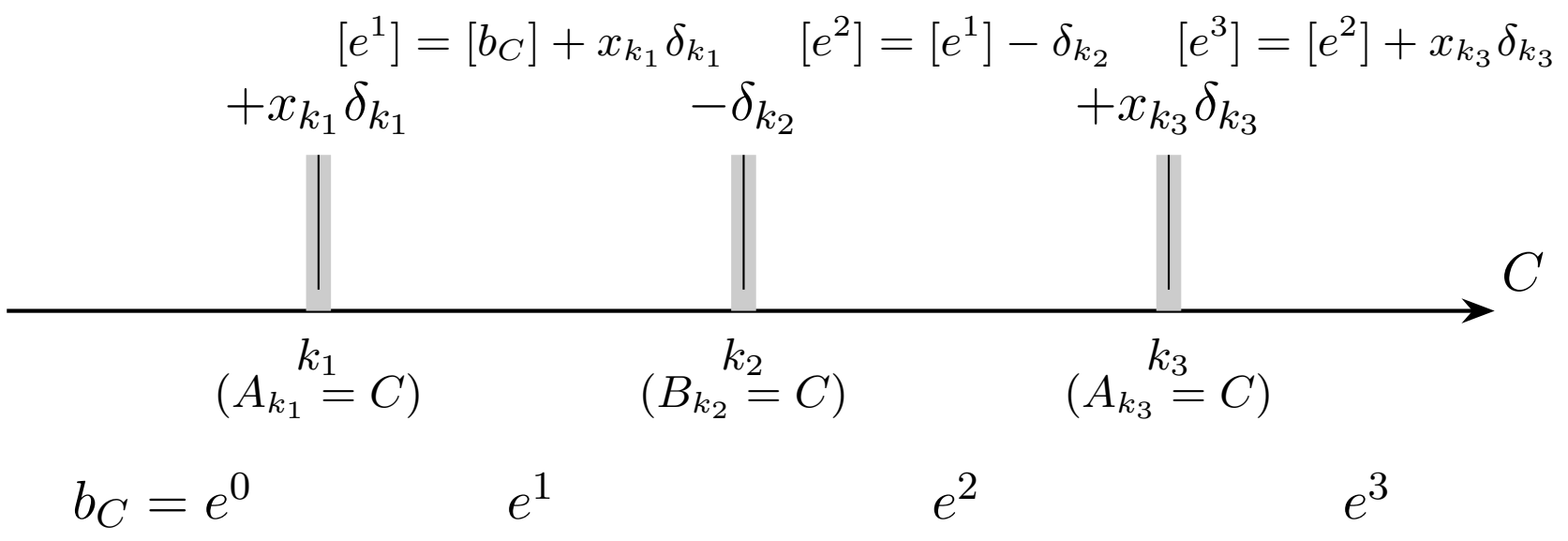
\begin{figure}[ht]
\centering
\begin{tikzpicture}[>=Stealth,scale=1,every node/.style={font=\small}]
\draw[->,thick] (0,0) -- (10.5,0);
\node[above] at (10.7,0) {$C$};
\foreach \x/\lab/\jump/\typ in {2.2/k_1/+x_{k_1}\d_{k_1}/A_{k_1}=C, 5.2/k_2/-\d_{k_2}/B_{k_2}=C, 8.2/k_3/+x_{k_3}\d_{k_3}/A_{k_3}=C}{
  \draw[line width=5pt,gray!40] (\x,0) -- (\x,1.1);
  \draw (\x,0.15) -- (\x,1.1);
  \node[above] at (\x,1.1) {$\jump$};
  \node[below,font=\scriptsize] at (\x,-0.05) {$\lab$};
  \node[font=\scriptsize] at (\x,-0.6) {$(\typ)$};
}
\foreach \x/\e in {1.1/b_C=e^0, 3.7/e^1, 6.7/e^2, 9.6/e^3}{ \node[below] at (\x,-1.0) {$\e$}; }
\node[font=\scriptsize] at (3.7,1.9) {$[e^1]=[b_C]+x_{k_1}\d_{k_1}$};
\node[font=\scriptsize] at (6.7,1.9) {$[e^2]=[e^1]-\d_{k_2}$};
\node[font=\scriptsize] at (9.6,1.9) {$[e^3]=[e^2]+x_{k_3}\d_{k_3}$};
\end{tikzpicture}
\caption{The jump relations along a Seifert circle $C$ with base edge $b_C$ (Step~2 of the proof of Proposition~\ref{prop.seifert}): at a crossing $k$ with $A_k=C$ the generator jumps by $x_k\d_k$, at a crossing with $B_k=C$ it jumps by $-\d_k$, and $[e]=[b_C]+J_C(e)$ is the sum of the jumps met on the way from $b_C$ to $e$; the full turn $\CC_C$ is the sum of all the jumps of $C$.}
\label{fig.jump}
\end{figure}

Let $\Lam_c^c$ be the free $\Lam_c$-module of rank $c$, with one basis vector $\d_j$ for each crossing $j$; we think of $\d_j$ as the generator attached to the band of the crossing $j$ (in the proof of Proposition~\ref{prop.seifert} it will be the jump $[u_j]-[o_j]$ across this band), and we identify $\Lam_c^c$ at $x=1$ with $\BZ^c\otimes\Lam_c$, where $\BZ^c$ is the group in which the cycle space $H$ lives, by $\d_j\mapsto e_j$. For an edge $e$ on the circle $C$ we define the \emph{jump sum}
\begin{equation}\label{eq.jumpsum}
   J_C(e)=\sum_{k\in\mathrm{Pre}_C(e),\ A_k=C}x_k\d_k\ -\sum_{k\in\mathrm{Pre}_C(e),\ B_k=C}\d_k\ \in\Lam_c^c ,
\end{equation}
and the \emph{full turn} of the circle $C$,
\begin{equation}\label{eq.fullturn}
   \CC_C=\sum_{k:\,A_k=C}x_k\d_k\ -\sum_{k:\,B_k=C}\d_k\ \in\Lam_c^c .
\end{equation}
Finally let $K(x)$ be the $c\times c$ matrix over $\Lam_c$ whose row $j$ is the vector
\begin{equation}\label{eq.Krow}
   \d_j+J_{A_j}(o_j)-J_{B_j}(u_j)\ \in\Lam_c^c
\end{equation}
written in the basis $\d_1,\dots,\d_c$: the entry $K(x)_{jk}$ is the coefficient of $\d_k$ in \eqref{eq.Krow}. Explicitly, $K(x)_{jj}=1$, because $j$ is the crossing at the end of $o_j$ and of $u_j$ and is therefore excluded from $\mathrm{Pre}_{A_j}(o_j)$ and from $\mathrm{Pre}_{B_j}(u_j)$, and for $k\ne j$
\[
   K(x)_{jk}=\theta_{A_j,o_j}(k)-\theta_{B_j,u_j}(k),\qquad
   \theta_{C,e}(k)=\begin{cases}x_k&\text{if $k\in\mathrm{Pre}_C(e)$ and $A_k=C$,}\\ -1&\text{if $k\in\mathrm{Pre}_C(e)$ and $B_k=C$,}\\ 0&\text{if $k\notin\mathrm{Pre}_C(e)$.}\end{cases}
\]
Thus $K(x)_{jk}\ne0$ only if the crossing $k$ lies on $A_j$ or on $B_j$ and is met before the crossing $j$ when going around that circle from its base edge; the matrix $K(x)$ is read off from the cyclic order of the crossings along the Seifert circles. The jump sums are illustrated in Figure~\ref{fig.jump}.

Let $\CC(x)$ be the $(s-1)\times c$ matrix whose rows, indexed by the Seifert circles $C\ne\hat C$, are the full turns $\CC_C$ written in the basis $\d_1,\dots,\d_c$: by \eqref{eq.fullturn},
\begin{equation}\label{eq.CC}
   \CC(x)_{C,k}=\begin{cases}x_k&\text{if $A_k=C$,}\\-1&\text{if $B_k=C$,}\\0&\text{if the crossing $k$ does not lie on $C$.}\end{cases}
\end{equation}
Thus $\CC(x)$ is a deformation of the incidence matrix of the Seifert graph. Indeed, let $\partial_A$ and $\partial_B$ be the $s\times c$ matrices with entries $(\partial_A)_{C,k}=\ind{A_k=C}$ and $(\partial_B)_{C,k}=\ind{B_k=C}$, so that the incidence map $\partial$, written as an $s\times c$ matrix with $\partial_{C,k}=\ind{B_k=C}-\ind{A_k=C}$, is $\partial=\partial_B-\partial_A$, and let $X=\diag(x_1,\dots,x_c)$. Then $\CC(x)$ is the matrix $\partial_AX-\partial_B$ with the row $\hat C$ removed, and in particular
\[
   \CC(1)=-\hat\partial ,
\]
where $\hat\partial$ denotes $\partial$ with the row $\hat C$ removed. Since the columns of $Y$ are cycles, $\partial Y=0$, and hence $\hat\partial Y=0$.

\begin{proposition}\label{prop.seifert}
With the notation above,
\[
   P_{D,\hate}(x)=\det\begin{pmatrix}K(x)&-\hat\partial^{\T}\\-\CC(x)&0\end{pmatrix},
\]
the rows and columns of the matrix being indexed by the crossings followed by the Seifert circles $C\ne\hat C$.
\end{proposition}
\begin{proof}
Let $A=1-T^a_D(x)$ and let $\hat A$ be the matrix obtained from $A$ by replacing its row $\hate$ by the unit row vector $e_{\hate}^{\T}$; expanding $\det\hat A$ along this row gives $\det\hat A=\det A_{\hate}=P_{D,\hate}(x)$. We transform $\hat A$ by operations which do not change the determinant. The matrix $\hat A$ is a presentation matrix: its rows are relations between the generators $[e]$, $e\in E$, and we read the row $u'_j$ of $A$, which is $e_{u'_j}^{\T}-(1-x_j)e_{o_j}^{\T}-x_je_{u_j}^{\T}$, as the relation $[u'_j]=(1-x_j)[o_j]+x_j[u_j]$, the row $o'_j$ as $[o'_j]=[o_j]$, and the row $\hate$ of $\hat A$ as $[\hate]=0$.

\emph{Step 1: new generators.} Let $\hat A'$ be the $(2c+c)\times(2c+c)$ matrix $\bigl(\begin{smallmatrix}\hat A&0\\N&1\end{smallmatrix}\bigr)$, whose last $c$ rows, indexed by $\d_1,\dots,\d_c$, are $N_j=e_{o_j}^{\T}-e_{u_j}^{\T}$ in the first $2c$ columns, followed by the identity matrix. The new row $\d_j$ is the relation $\d_j=[u_j]-[o_j]$, and $\det\hat A'=\det\hat A$.

\emph{Step 2: the jump relations.} For every crossing $j$ with $u'_j\ne\hate$ we subtract $x_j$ times the row $\d_j$ from the row $u'_j$, which becomes $e_{u'_j}^{\T}-e_{o_j}^{\T}-x_je_{\d_j}^{\T}$; for every crossing $j$ with $o'_j\ne\hate$ we add the row $\d_j$ to the row $o'_j$, which becomes $e_{o'_j}^{\T}-e_{u_j}^{\T}+e_{\d_j}^{\T}$. After these row operations, for every edge $e\ne\hate$ the row $e$ expresses $[e]$ as the generator of the edge preceding $e$ on its Seifert circle plus a \emph{jump}: $[u'_j]=[o_j]+x_j\d_j$ on the circle $A_j$ and $[o'_j]=[u_j]-\d_j$ on the circle $B_j$. In other words, going around a Seifert circle $C$, the generator jumps by $x_k\d_k$ at a crossing $k$ with $A_k=C$ and by $-\d_k$ at a crossing $k$ with $B_k=C$, which are exactly the terms of \eqref{eq.jumpsum} and \eqref{eq.fullturn}.

\emph{Step 3: elimination of the edges.} Let $C$ be a Seifert circle and let $b_C=e^0,e^1,\dots,e^{m-1}$ be its edges in the order in which they are met along $C$. For $i=1,\dots,m-1$ in this order, we clear the column $e^i$ of the matrix by subtracting multiples of the row $e^i$ from the other rows, and then expand the determinant along the column $e^i$, which now contains only the entry $1$ in the row $e^i$: the determinant is unchanged, and the row and the column of $e^i$ disappear. By induction on $i$, when the column $e^i$ is cleared the row $e^i$ is $e_{e^i}^{\T}-e_{b_C}^{\T}-J_C(e^i)$ (we identify the vector $J_C(e)\in\Lam_c^c$ with the row vector of its coordinates in the columns $\d_1,\dots,\d_c$): indeed, the crossing at the end of $e^{i-1}$ is the crossing $k$ which $e^i$ leaves, and $J_C(e^i)$ is $J_C(e^{i-1})$ plus the jump of $e^i$. Thus, in all remaining rows, the generator $[e^i]$ is replaced by $[b_C]+J_C(e^i)$. If $C\ne\hat C$, the row of the base edge $b_C$, which is $e_{b_C}^{\T}-e_{e^{m-1}}^{\T}-(\text{jump of }b_C)$, becomes $-\CC_C$ after the substitution, the crossing at the end of $e^{m-1}$ being the one which $b_C$ leaves: the entry of this row in the column $b_C$ vanishes, and the row is the relation ``the full turn of $C$ is zero''. If $C=\hat C$, then $b_C=\hate$ and the row $\hate$ is $e_{\hate}^{\T}$; we clear the column $\hate$ using this row, which replaces $[\hate]$ by $0$ in the other rows, and expand along the column $\hate$. We do this for all Seifert circles. The remaining rows are the rows $\d_j$, $j=1,\dots,c$, and the rows $b_C$, $C\ne\hat C$; the remaining columns are the columns $\d_j$ and $b_C$, $C\ne\hat C$. The row $\d_j$ was $e_{\d_j}^{\T}+e_{o_j}^{\T}-e_{u_j}^{\T}$ and has become
\[
   e_{\d_j}^{\T}+J_{A_j}(o_j)-J_{B_j}(u_j)+e_{b_{A_j}}^{\T}-e_{b_{B_j}}^{\T},\qquad\text{with }e_{b_{\hat C}}=e_{\hate}\mapsto0 ,
\]
whose $\d$-part is the row $j$ of $K(x)$ by \eqref{eq.Krow}, and whose entry in the column $b_C$ is $\ind{A_j=C}-\ind{B_j=C}=-\hat\partial_{C,j}$. The row $b_C$ has become $-\CC_C$, with no entry in the columns $b_{C'}$. Ordering the rows and the columns in the same way (the $\d_j$ first, then the $b_C$) gives the matrix of the statement.
\end{proof}

The second flavor is treated in the same way, with the roles of the two sides exchanged. Let $K^b(y)$ be the $c\times c$ matrix over $\Lam'_c$ whose row $j$ is the coordinate vector of
\begin{equation}\label{eq.Kbrow}
   \d_j-J^b_{B_j}(u_j)+y_j\,J^b_{A_j}(o_j),\qquad J^b_C(e)=\sum_{k\in\mathrm{Pre}_C(e),\ A_k=C}\d_k\ -\sum_{k\in\mathrm{Pre}_C(e),\ B_k=C}\d_k ,
\end{equation}
and let $\CC(y)$ be the matrix \eqref{eq.CC} with $x$ replaced by $y$.

\begin{proposition}\label{prop.seifertb}
With the notation above,
\[
   \tP_{D,\hate}(y)=\det\begin{pmatrix}K^b(y)&\CC(y)^{\T}\\ \hat\partial&0\end{pmatrix},
\]
the rows and columns of the matrix being indexed by the crossings followed by the Seifert circles $C\ne\hat C$.
\end{proposition}
\begin{proof}
By \eqref{eq.Tb}, the row $o'_j$ of $1-T^b_D(y)$ is $e_{o'_j}^{\T}-y_je_{o_j}^{\T}$ and the row $u'_j$ is $e_{u'_j}^{\T}-(1-y_j)e_{o_j}^{\T}-e_{u_j}^{\T}$. We repeat the proof of Proposition~\ref{prop.seifert} with the new generators $\d_j=[u_j]-y_j[o_j]$, i.e.\ with the rows $N_j=y_je_{o_j}^{\T}-e_{u_j}^{\T}$ in Step~1. In Step~2, subtracting the row $\d_j$ from the row $u'_j$ gives $e_{u'_j}^{\T}-e_{o_j}^{\T}-e_{\d_j}^{\T}$, and adding it to the row $o'_j$ gives $e_{o'_j}^{\T}-e_{u_j}^{\T}+e_{\d_j}^{\T}$: the jumps are now $+\d_k$ at a crossing $k$ with $A_k=C$ and $-\d_k$ at a crossing with $B_k=C$, which gives the jump sums $J^b_C(e)$ of \eqref{eq.Kbrow} and the full turns $\sum_{A_k=C}\d_k-\sum_{B_k=C}\d_k$, whose coordinate vector is $-\hat\partial_C$, the row $C$ of $-\hat\partial$, for $C\ne\hat C$, after Step~3. So the row $b_C$ becomes $\hat\partial_C$ (in Step~3 of the proof of Proposition~\ref{prop.seifert} the row $b_C$ became minus the full turn). The row $\d_j$, which was $e_{\d_j}^{\T}-e_{u_j}^{\T}+y_je_{o_j}^{\T}$, becomes $e_{\d_j}^{\T}-J^b_{B_j}(u_j)+y_jJ^b_{A_j}(o_j)-e_{b_{B_j}}^{\T}+y_je_{b_{A_j}}^{\T}$ (with $e_{\hate}\mapsto0$), whose $\d$-part is the row $j$ of $K^b(y)$ and whose entry in the column $b_C$ is $y_j\ind{A_j=C}-\ind{B_j=C}=\CC(y)_{C,j}$.
\end{proof}

It remains to recognize the crossing forms in the matrices $K(x)$, $\CC(x)$, $K^b(y)$ and $\CC(y)$.

\begin{lemma}[the rows are flows]\label{lem.rowsareflows}
Let $\gamma\in H$ be a cycle of the Seifert graph, written as a column vector with the coordinates $w_j(\gamma)$, let $z^\gamma$ be its lift, and let $C$ be a Seifert circle. Then, for every crossing $k$,
\begin{align*}
   \bigl(\gamma^{\T}K(x)\bigr)_k&=\pi_k(z^\gamma),& \CC(x)_{C,k}&=\pi_k(C)\quad(C\ne\hat C),\\
   \bigl(K^b(y)\,\gamma\bigr)_k&=-\tpi_k(z^\gamma),& \CC(y)_{C,k}&=\tpi_k(C)\quad(C\ne\hat C).
\end{align*}
\end{lemma}
\begin{proof}
The flow $C$ of a Seifert circle is $1$ on the edges of $C$ and $0$ on the others, so $\pi_k(C)=x_k\ind{u'_k\subseteq C}-\ind{u_k\subseteq C}=x_k\ind{A_k=C}-\ind{B_k=C}$, which is $\CC(x)_{C,k}$ by \eqref{eq.CC}; and $\tpi_k(C)=y_k\ind{o_k\subseteq C}-\ind{o'_k\subseteq C}=y_k\ind{A_k=C}-\ind{B_k=C}=\CC(y)_{C,k}$.

Write $w_j$ for $w_j(\gamma)$. By \eqref{eq.Krow} and \eqref{eq.jumpsum}, the coefficient of $\d_k$ in the row $j\ne k$ of $K(x)$ is $x_k$ if $A_k=A_j$ and $k\prec_{A_j}j$, it is $-1$ if $B_k=A_j$ and $k\prec_{A_j}j$, it is $-x_k$ if $A_k=B_j$ and $k\prec_{B_j}j$, and it is $+1$ if $B_k=B_j$ and $k\prec_{B_j}j$; in the first and third cases the circle is $A_k$ and the coefficient is $s_{A_k}(j)\,x_k$, and in the other two cases the circle is $B_k$ and the coefficient is $-s_{B_k}(j)$. Together with $K(x)_{kk}=1$ this gives
\[
   \bigl(\gamma^{\T}K(x)\bigr)_k=\sum_jw_jK(x)_{jk}=x_k\sum_{j:\,k\prec_{A_k}j}s_{A_k}(j)\,w_j\ +\ w_k-\sum_{j:\,k\prec_{B_k}j}s_{B_k}(j)\,w_j .
\]
The first sum is $z^\gamma_{u'_k}$ by \eqref{eq.gh}. By \eqref{eq.sw} for the circle $B_k$, on which $s_{B_k}(k)=-1$, the remaining terms are $w_k-\sum_{j\succ k}s_{B_k}(j)w_j=\sum_{j\prec k}s_{B_k}(j)w_j=-z^\gamma_{u_k}$, again by \eqref{eq.gh}. Hence $\bigl(\gamma^{\T}K(x)\bigr)_k=x_kz^\gamma_{u'_k}-z^\gamma_{u_k}=\pi_k(z^\gamma)$.

For the second flavor, by \eqref{eq.Kbrow} the row $k$ of $K^b(y)$ is
\[
   \d_k-\sum_{j\prec_{B_k}k}s_{B_k}(j)\,\d_j+y_k\sum_{j\prec_{A_k}k}s_{A_k}(j)\,\d_j ,
\]
so that
\[
   \bigl(K^b(y)\gamma\bigr)_k=w_k-\sum_{j\prec_{B_k}k}s_{B_k}(j)\,w_j+y_k\sum_{j\prec_{A_k}k}s_{A_k}(j)\,w_j=w_k+z^\gamma_{u_k}-y_k\,z^\gamma_{o_k},
\]
by \eqref{eq.gh} and by \eqref{eq.lift} for the edge $o_k$ of the circle $A_k$, for which $\mathrm{Pre}_{A_k}(o_k)=\{j:j\prec_{A_k}k\}$. Since $z^\gamma_{o'_k}=z^\gamma_{u_k}+w_k$ by \eqref{eq.walkcircle}, this is $z^\gamma_{o'_k}-y_kz^\gamma_{o_k}=-\tpi_k(z^\gamma)$.
\end{proof}

\begin{proof}[Proof of Theorem~\ref{thm.flows}]
Let $\CT$ be a spanning tree of the Seifert graph and $Y$ a basis of its cycle space, with columns $\gamma_1,\dots,\gamma_\nu$. We order the Seifert basis as $z^{\gamma_1},\dots,z^{\gamma_\nu}$, followed by the circles $C\ne\hat C$. By Lemma~\ref{lem.rowsareflows}, the matrices of the crossing forms in this basis are the $c\times c$ matrices
\begin{equation}\label{eq.N}
   N_{D,\hate}(x)=\bigl(\pi_k(z_i)\bigr)_{i,k}=\begin{pmatrix}Y^{\T}K(x)\\ \CC(x)\end{pmatrix},\qquad
   \tilde N_{D,\hate}(y)=\bigl(\tpi_k(z_i)\bigr)_{i,k}=\begin{pmatrix}-\bigl(K^b(y)Y\bigr)^{\T}\\ \CC(y)\end{pmatrix}.
\end{equation}
Let $U$ be the $c\times c$ matrix with the rows of $Y^{\T}$ followed by those of $E_{\CT}^{\T}$; it is unimodular, being the transpose of $(Y\,|\,E_{\CT})$. Multiplying the first $c$ rows of the matrix of Proposition~\ref{prop.seifert} by $U$ multiplies the determinant by $\det U=\pm1$ and gives
\[
   \begin{pmatrix}Y^{\T}K(x)&-(\hat\partial Y)^{\T}\\ E_{\CT}^{\T}K(x)&-(\hat\partial E_{\CT})^{\T}\\-\CC(x)&0\end{pmatrix}
   =\begin{pmatrix}Y^{\T}K(x)&0\\ E_{\CT}^{\T}K(x)&-(\hat\partial E_{\CT})^{\T}\\-\CC(x)&0\end{pmatrix},
\]
because the columns of $Y$ are cycles, $\hat\partial Y=0$. The last $s-1$ columns have their nonzero entries in the middle block $-(\hat\partial E_{\CT})^{\T}$, which is, up to sign, the transpose of the incidence matrix of the tree $\CT$ with the row of the vertex $\hat C$ removed, a unimodular matrix. Exchange the last two blocks of rows, which have $s-1$ rows each: the matrix becomes block triangular, and its determinant is
\begin{equation}\label{eq.signN}
\begin{aligned}
   \det U\cdot P_{D,\hate}(x)&=(-1)^{s-1}\det\begin{pmatrix}Y^{\T}K(x)\\-\CC(x)\end{pmatrix}\det\bigl(-(\hat\partial E_{\CT})^{\T}\bigr)\\
   &=(-1)^{s-1}\det(\hat\partial E_{\CT})\,\det N_{D,\hate}(x) .
\end{aligned}
\end{equation}
Hence
\begin{equation}\label{eq.kappaSeifert}
   \det N_{D,\hate}(x)=\kappa\,P_{D,\hate}(x)\qquad\text{with}\qquad\kappa=(-1)^{s-1}\det(Y\,|\,E_{\CT})\,\det(\hat\partial E_{\CT})=\pm1 .
\end{equation}
For the second flavor, multiply the matrix of Proposition~\ref{prop.seifertb} on the right by the matrix $\diag\bigl((Y\,|\,E_{\CT}),1\bigr)$, of determinant $\det U$: it becomes
\[
   \begin{pmatrix}K^b(y)Y&K^b(y)E_{\CT}&\CC(y)^{\T}\\ 0&\hat\partial E_{\CT}&0\end{pmatrix},
\]
since $\hat\partial Y=0$. Exchanging the last two blocks of columns, which have $s-1$ columns each, makes this matrix block triangular, so that
\[
   \det U\cdot\tP_{D,\hate}(y)=(-1)^{s-1}\det(\hat\partial E_{\CT})\,\det\bigl(K^b(y)Y\ \big|\ \CC(y)^{\T}\bigr).
\]
The matrix $\bigl(K^b(y)Y\,|\,\CC(y)^{\T}\bigr)$ is the transpose of $\tilde N_{D,\hate}(y)$ with the signs of its first $\nu$ rows changed. Therefore
\[
   \det\tilde N_{D,\hate}(y)=(-1)^{\nu}\,\kappa\,\tP_{D,\hate}(y),
\]
with the same sign $\kappa$ as in \eqref{eq.kappaSeifert}. Finally, the oriented smoothing of one crossing changes the number of components of a link diagram by one: it merges two components if the two strands of the crossing belong to different components, and it splits a component into two otherwise. Smoothing the $c$ crossings one after the other turns the $\mu$ components of $L$ into the $s$ Seifert circles, so that $s\equiv\mu+c$ modulo $2$, and $\nu=c-s+1\equiv\mu-1$ modulo $2$: $(-1)^{\nu}=(-1)^{\mu-1}$. This proves the two identities of the theorem for the Seifert basis, hence for every basis by Remark~\ref{rem.flows}~(i). The formula for $\LG$ follows from Theorem~\ref{thm.HK}: $\BE_D\bigl[\det(\pi_j(z_i))\det(\tpi_j(z_i))\bigr]=\kappa^2(-1)^{\mu-1}\,\BE_D\bigl[P_{D,\hate}\,\tP_{D,\hate}\bigr]=\LG(L;t_0,t_1)$.
\end{proof}

In the matrix $N_{D,\hate}(x)$ of \eqref{eq.N} the rows are of two kinds. The last $s-1$ rows, those of the Seifert circles, form the matrix $\CC(x)$, a deformation of the incidence matrix of the Seifert graph. The first $\nu$ rows are those of the cycles of the Seifert graph; by \eqref{eq.gh}, their entries are
\begin{equation}\label{eq.Ncolumns}
   \pi_k(z^{\gamma})=x_k\,g_k(\gamma)+h_k(\gamma),\qquad g_k(\gamma)=z^\gamma_{u'_k},\quad h_k(\gamma)=-z^\gamma_{u_k} .
\end{equation}
The passage from this $c\times c$ matrix to a $\nu\times\nu$ matrix on the cycle space, and its interpretation as a deformed Seifert form of the canonical surface, are the subject of Part~\ref{part.surfaces}.

\subsection{The basis of the faces: Alexander's matrix}\label{ssec.facematrices}
The Alexander polynomial has a third classical expression besides the Wirtinger presentation and the Seifert matrix, which is in fact the oldest one: Alexander's original matrix \cite{Alexander28}, whose rows are the crossings and whose columns are the faces of the diagram, each crossing contributing the entries $t,-t,1,-1$ at the four faces around it; the determinant of the matrix with two adjacent columns removed is $\Delta_L(t)$ up to a unit. In the lattice of flows this is the basis of Lemma~\ref{lem.faces}~(iii).

\begin{definition}[the face matrices]\label{def.facematrices}
The \emph{face matrices} of $D$ are the $c\times\CF$ matrices $Q_D(x)$ and $\tQ_D(y)$, with entries in $\BZ[x]$ and $\BZ[y]$, whose entries are the values of the crossing forms on the boundaries of the faces,
\begin{equation}\label{eq.facematrixentries}
   Q_D(x)_{j,f}=\pi_j(\partial f),\qquad \tQ_D(y)_{j,f}=\tpi_j(\partial f).
\end{equation}
In terms of the vectors $m(e)$ of \eqref{eq.m}, their rows are, for $j=1,\dots,c$,
\begin{equation}\label{eq.facerows}
   a_j(x_j)^{\T}=\bigl(x_j\,m(u'_j)-m(u_j)\bigr)^{\T}\qquad\text{and}\qquad b_j(y_j)^{\T}=\bigl(y_j\,m(o_j)-m(o'_j)\bigr)^{\T}.
\end{equation}
For an edge $\hate$ of $D$, $Q_{D,\hate}(x)$ and $\tQ_{D,\hate}(y)$ denote the $c\times c$ matrices obtained by removing the two columns of the faces $f_{\rm L}(\hate)$ and $f_{\rm R}(\hate)$ adjacent to $\hate$.
\end{definition}

By \eqref{eq.mcrossing}, the rows are, in the basis $(f^\uparrow_j,f^\leftarrow_j,f^\downarrow_j,f^\rightarrow_j)$ of the four faces at the crossing (coinciding faces receive the sum of their entries),
\begin{equation}\label{eq.facetable}
\begin{array}{c|cccc|cccc}
   & \multicolumn{4}{c|}{a_j(x_j)} & \multicolumn{4}{c}{b_j(y_j)}\\
   & f^\uparrow_j & f^\leftarrow_j & f^\downarrow_j & f^\rightarrow_j & f^\uparrow_j & f^\leftarrow_j & f^\downarrow_j & f^\rightarrow_j\\ \hline
   \ve_j=1 & -x_j & x_j & -1 & 1 & -1 & y_j & -y_j & 1\\
   \ve_j=-1 & x_j & -1 & 1 & -x_j & 1 & -1 & y_j & -y_j
\end{array}
\end{equation}
The row $a_j$ vanishes exactly on the vectors $k\in\BF^{\CF}$ with $m(u_j)^{\T}k=x_j\,m(u'_j)^{\T}k$: it is the relation ``the meridian of $u_j$ is $x_j$ times the meridian of $u'_j$'' of the under-strand, which is the relation $P_{D,u_j}=x_j\,P_{D,u'_j}$ satisfied by the cut polynomials (Lemma~\ref{lem.understrand}), written in the face generators; likewise $b_j$ is the relation ``the meridian of $o'_j$ is $y_j$ times the meridian of $o_j$'' of the over-strand. The first flavor deforms the under-strand of the crossing and the second one its over-strand. At the Alexander point $x_j=t^{\ve_j}$, $Q_{D,\hate}(t^{\ve})$ is Alexander's matrix \cite{Alexander28}: each crossing contributes $\pm t^{\ve_j}$ at the two corners adjacent to the outgoing under-edge $u'_j$ and $\pm1$ at the two corners adjacent to the incoming under-edge $u_j$, the two faces adjacent to $\hate$ being omitted.

\begin{theorem}[face form]\label{thm.faceform}
For every diagram $D$ of a link with $\mu$ components and every edge $\hate$ of $D$,
\[
   P_{D,\hate}(x)=\kappa\,\det Q_{D,\hate}(x)\qquad\text{and}\qquad \tP_{D,\hate}(y)=(-1)^{\mu-1}\kappa\,\det\tQ_{D,\hate}(y)
\]
as polynomials in the crossing variables, where the sign $\kappa=\pm1$ depends only on the orders of the crossings and of the faces, which are the same for the rows and for the columns of the two face matrices.
\end{theorem}
\begin{proof}
By Lemma~\ref{lem.faces}~(iii) the boundaries of the faces which are not adjacent to $\hate$ form a basis of $\FF_{D,\hate}$, and by \eqref{eq.facematrixentries} the matrices of the crossing forms in this basis are the transposes of $Q_{D,\hate}(x)$ and of $\tQ_{D,\hate}(y)$. So this is Theorem~\ref{thm.flows} for the basis of the faces.
\end{proof}

Thus the edge polynomials are Alexander's face matrices with the crossing variables inserted; expanding their determinants over the bijections between the crossings and the faces gives the Kauffman states of Section~\ref{sec.kauffman}, where the sign $\kappa$ is determined (Proposition~\ref{prop.statesum}~(ii)).

\begin{remark}[the faces as a basis adapted to a smoothing]\label{rem.checkerboardbasis}
Color the faces of $D$ black and white in such a way that the two faces adjacent to an edge have different colors. At every crossing two opposite corners are black, and there is exactly one smoothing $\vs_{\rm b}$ of $D$ whose arcs cut off the black corners. Its circles are the boundaries of the black faces, pushed slightly into them: the flow $z^C$ of the circle in the black face $f$ is $\pm\partial f$. The boundaries of the white faces are flows whose classes generate the quotient of $\FF_D$ by the circles, by Lemma~\ref{lem.faces}~(ii). So the basis of the faces is, up to signs, a basis adapted to the smoothing $\vs_{\rm b}$ in the sense of Proposition~\ref{prop.adapted}, and also to the complementary smoothing $\vs_{\rm w}$, with the roles of the two colors exchanged. Alexander's matrix and the matrix $N_{D,\hate}(x)$ of the Seifert circles are two instances of the same construction; we come back to the checkerboard smoothings in Section~\ref{sec.minusone}.
\end{remark}

The second flavor is the first flavor of the \emph{inverse diagram} $\bD$, obtained from $D$ by reversing the orientation and changing all crossings. The crossings of $\bD$ are those of $D$, with $\bar o_j=u'_j$, $\bar u_j=o'_j$, $\bar o'_j=u_j$ and $\bar u'_j=o_j$ (the outgoing under-edge of $D$ becomes the incoming over-edge, and so on); the edges of $\bD$ are those of $D$ with the opposite orientation, so that $m_{\bD}(e)=-m(e)$; and the faces of $\bD$ are those of $D$.

\begin{lemma}[the inverse diagram]\label{lem.inverse}
$T^a_{\bD}(y)=T^b_D(y)^{\T}$, hence $P_{\bD,\hate}(y)=\tP_{D,\hate}(y)$ for every edge $\hate$; and $Q_{\bD}(y)=-\tQ_D(y)$.
\end{lemma}
\begin{proof}
By \eqref{eq.Ta} applied to $\bD$, $T^a_{\bD}(y)=\sum_j\bigl(E_{u_ju'_j}+(1-y_j)E_{o_ju'_j}+y_jE_{o_jo'_j}\bigr)$, which is the transpose of $T^b_D(y)$ in \eqref{eq.Tb}; removing the row and the column of $\hate$ commutes with transposition. The row $j$ of $Q_{\bD}(y)$ is $y_jm_{\bD}(\bar u'_j)-m_{\bD}(\bar u_j)=-y_jm(o_j)+m(o'_j)=-b_j(y_j)$.
\end{proof}

We end this section with the relations between the edge polynomials of the different cut edges. They say that the edge polynomials themselves form a flow, on which all the crossing forms vanish.

\begin{lemma}\label{lem.understrand}
For every crossing $j$ of $D$,
\[
   P_{D,u_j}(x)=x_j\,P_{D,u'_j}(x)\qquad\text{and}\qquad P_{D,o_j}(x)=P_{D,o'_j}(x)+(1-x_j)\,P_{D,u'_j}(x).
\]
\end{lemma}
\begin{proof}
Let $A=1-T^a_D(x)$, a $2c\times2c$ matrix over the field $\BQ(x_1,\dots,x_c)$ with $A\1=0$, and let $\operatorname{adj}(A)$ be its adjugate, whose diagonal entry $(\hate,\hate)$ is $\det A_{\hate}=P_{D,\hate}(x)$. If $A$ has rank at most $2c-2$, then $\operatorname{adj}(A)=0$ and all the $P_{D,\hate}$ vanish. Otherwise $\ker A$ is spanned by $\1$ and the left kernel of $A$ by a row vector $n$, and the identities $A\operatorname{adj}(A)=\operatorname{adj}(A)A=0$ show, as in the proof of Lemma~\ref{lem.alexpoly}, that $\operatorname{adj}(A)=\lambda\,\1n$ for some $\lambda\ne0$; hence $P_{D,\hate}=\lambda\,n_{\hate}$ for every edge $\hate$. By \eqref{eq.Ta}, the column $u_j$ of $A$ has the entries $1$ in the row $u_j$ and $-x_j$ in the row $u'_j$, and the column $o_j$ has the entries $1$ in the row $o_j$, $-1$ in the row $o'_j$ and $-(1-x_j)$ in the row $u'_j$. The equations $(nA)_{u_j}=0$ and $(nA)_{o_j}=0$ read $n_{u_j}=x_jn_{u'_j}$ and $n_{o_j}=n_{o'_j}+(1-x_j)n_{u'_j}$.
\end{proof}

\begin{proposition}[the flow of the cut polynomials]\label{prop.potential}
The vector $v(x)=\bigl(P_{D,e}(x)\bigr)_{e\in E}$ is a flow on $D$ with coefficients in $\BZ[x]$, and $\pi_j\bigl(v(x)\bigr)=0$ for every crossing $j$. Consequently there is a vector $k(x)\in\BQ(x)^{\CF}$, unique up to adding a multiple of $\1$, such that
\[
   P_{D,e}(x)=k_{f_{\rm L}(e)}(x)-k_{f_{\rm R}(e)}(x)\quad\text{for every edge $e$},\qquad\text{and}\qquad Q_D(x)\,k(x)=0 .
\]
Thus $\ker Q_D(x)\supseteq\spn(\1,k(x))$, with equality if $P_{D,\hate}\ne0$ for some $\hate$.
\end{proposition}
\begin{proof}
By Lemma~\ref{lem.understrand}, $v_{u_j}=x_jv_{u'_j}$ and $v_{o_j}=v_{o'_j}+(1-x_j)v_{u'_j}$ at every crossing. The first identity is $\pi_j(v)=0$, and the sum of the two is $v_{o_j}+v_{u_j}=v_{o'_j}+v_{u'_j}$: the vector $v$ satisfies \eqref{eq.flow}. Lemma~\ref{lem.faces}~(ii) and its proof hold with coefficients in the field $\BQ(x)$ in place of $\BZ$; so $v=\sum_fk_f\,\partial f$ for a potential $k\in\BQ(x)^{\CF}$, unique up to adding multiples of $\1$, and $v_e=k_{f_{\rm L}(e)}-k_{f_{\rm R}(e)}$ by \eqref{eq.potential}. By \eqref{eq.facematrixentries}, the entry $j$ of $Q_D(x)k$ is $\sum_fk_f\,\pi_j(\partial f)=\pi_j(v)=0$; and $Q_D(x)\1=\Pi(x)\bigl(\sum_f\partial f\bigr)=0$. If some $P_{D,\hate}\ne0$, then $v\ne0$, $k\notin\BQ(x)\1$, and $Q_{D,\hate}(x)$ is invertible by Theorem~\ref{thm.faceform}, so that $\ker Q_D(x)$ has dimension $2$.
\end{proof}

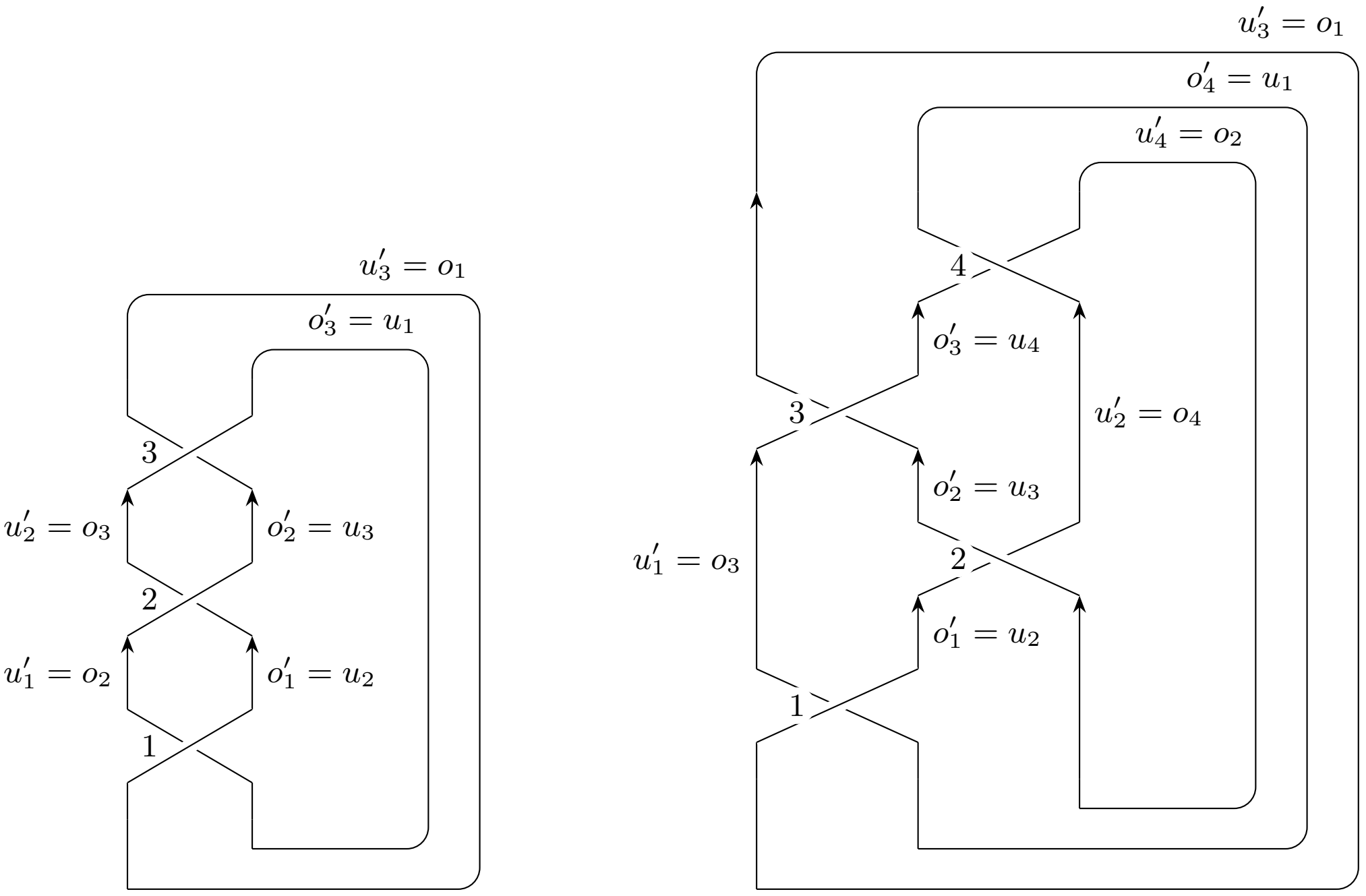
\begin{figure}[ht]
\centering
\begin{tikzpicture}[>=Stealth,scale=0.72]
\draw (0.00,-0.00) -- (0.00,0.50);
\draw[->] (0.00,1.50) -- (0.00,2.50);
\node[left,font=\scriptsize] at (0.00,2.00) {$u'_1=o_2$};
\draw[->] (0.00,3.50) -- (0.00,4.50);
\node[left,font=\scriptsize] at (0.00,4.00) {$u'_2=o_3$};
\draw (0.00,5.50) -- (0.00,6.00);
\draw[rounded corners=6pt] (0.00,6.00) -- (0.00,7.15) -- (4.80,7.15) -- (4.80,-0.95) -- (0.00,-0.95);
\draw (0.00,-0.95) -- (0.00,0.00);
\node[above,font=\scriptsize] at (3.90,7.15) {$u'_3=o_1$};
\draw (1.70,-0.00) -- (1.70,0.50);
\draw[->] (1.70,1.50) -- (1.70,2.50);
\node[right,font=\scriptsize] at (1.70,2.00) {$o'_1=u_2$};
\draw[->] (1.70,3.50) -- (1.70,4.50);
\node[right,font=\scriptsize] at (1.70,4.00) {$o'_2=u_3$};
\draw (1.70,5.50) -- (1.70,6.00);
\draw[rounded corners=6pt] (1.70,6.00) -- (1.70,6.40) -- (4.10,6.40) -- (4.10,-0.40) -- (1.70,-0.40);
\draw (1.70,-0.40) -- (1.70,0.00);
\node[above,font=\scriptsize] at (3.20,6.40) {$o'_3=u_1$};
\draw (1.70,0.50) -- (0.00,1.50);
\draw[preaction={draw=white,line width=4pt}] (0.00,0.50) -- (1.70,1.50);
\node[font=\scriptsize,circle,inner sep=1pt,fill=white] at (0.30,1.00) {$1$};
\draw (1.70,2.50) -- (0.00,3.50);
\draw[preaction={draw=white,line width=4pt}] (0.00,2.50) -- (1.70,3.50);
\node[font=\scriptsize,circle,inner sep=1pt,fill=white] at (0.30,3.00) {$2$};
\draw (1.70,4.50) -- (0.00,5.50);
\draw[preaction={draw=white,line width=4pt}] (0.00,4.50) -- (1.70,5.50);
\node[font=\scriptsize,circle,inner sep=1pt,fill=white] at (0.30,5.00) {$3$};
\end{tikzpicture}\hspace{1.2cm}
\begin{tikzpicture}[>=Stealth,scale=0.72]
\draw (0.00,-0.00) -- (0.00,0.50);
\draw[->] (0.00,1.50) -- (0.00,4.50);
\node[left,font=\scriptsize] at (0.00,3.00) {$u'_1=o_3$};
\draw[->] (0.00,5.50) -- (0.00,8.00);
\draw[rounded corners=6pt] (0.00,8.00) -- (0.00,9.90) -- (8.20,9.90) -- (8.20,-1.50) -- (0.00,-1.50);
\draw (0.00,-1.50) -- (0.00,0.00);
\node[above,font=\scriptsize] at (7.30,9.90) {$u'_3=o_1$};
\draw (2.20,-0.00) -- (2.20,0.50);
\draw[->] (2.20,1.50) -- (2.20,2.50);
\node[right,font=\scriptsize] at (2.20,2.00) {$o'_1=u_2$};
\draw[->] (2.20,3.50) -- (2.20,4.50);
\node[right,font=\scriptsize] at (2.20,4.00) {$o'_2=u_3$};
\draw[->] (2.20,5.50) -- (2.20,6.50);
\node[right,font=\scriptsize] at (2.20,6.00) {$o'_3=u_4$};
\draw (2.20,7.50) -- (2.20,8.00);
\draw[rounded corners=6pt] (2.20,8.00) -- (2.20,9.15) -- (7.50,9.15) -- (7.50,-0.95) -- (2.20,-0.95);
\draw (2.20,-0.95) -- (2.20,0.00);
\node[above,font=\scriptsize] at (6.60,9.15) {$o'_4=u_1$};
\draw[->] (4.40,-0.00) -- (4.40,2.50);
\draw[->] (4.40,3.50) -- (4.40,6.50);
\node[right,font=\scriptsize] at (4.40,5.00) {$u'_2=o_4$};
\draw (4.40,7.50) -- (4.40,8.00);
\draw[rounded corners=6pt] (4.40,8.00) -- (4.40,8.40) -- (6.80,8.40) -- (6.80,-0.40) -- (4.40,-0.40);
\draw (4.40,-0.40) -- (4.40,0.00);
\node[above,font=\scriptsize] at (5.90,8.40) {$u'_4=o_2$};
\draw (2.20,0.50) -- (0.00,1.50);
\draw[preaction={draw=white,line width=4pt}] (0.00,0.50) -- (2.20,1.50);
\node[font=\scriptsize,circle,inner sep=1pt,fill=white] at (0.55,1.00) {$1$};
\draw (2.20,2.50) -- (4.40,3.50);
\draw[preaction={draw=white,line width=4pt}] (4.40,2.50) -- (2.20,3.50);
\node[font=\scriptsize,circle,inner sep=1pt,fill=white] at (2.75,3.00) {$2$};
\draw (2.20,4.50) -- (0.00,5.50);
\draw[preaction={draw=white,line width=4pt}] (0.00,4.50) -- (2.20,5.50);
\node[font=\scriptsize,circle,inner sep=1pt,fill=white] at (0.55,5.00) {$3$};
\draw (2.20,6.50) -- (4.40,7.50);
\draw[preaction={draw=white,line width=4pt}] (4.40,6.50) -- (2.20,7.50);
\node[font=\scriptsize,circle,inner sep=1pt,fill=white] at (2.75,7.00) {$4$};
\end{tikzpicture}
\caption{The closed braids $\sigma_1^3$ (left) and $(\sigma_1\sigma_2^{-1})^2$ (right), with the crossings numbered in the order of the braid words and the edges labeled; the strands run upward and the closing arcs are on the right.}
\label{fig.diagrams}
\end{figure}

\begin{figure}[ht]
\centering
\begin{tikzpicture}[>=Stealth,scale=0.6]
\draw[rounded corners=6pt,postaction={decorate},decoration={markings,mark=at position 0.12 with {\arrow{>}}}] (0.00,-0.95) -- (0.00,7.15) -- (4.80,7.15) -- (4.80,-0.95) -- (0.00,-0.95) -- cycle;
\node[above,font=\scriptsize] at (4.20,7.15) {$C_1$};
\draw[rounded corners=6pt,postaction={decorate},decoration={markings,mark=at position 0.12 with {\arrow{>}}}] (1.70,-0.40) -- (1.70,6.40) -- (4.10,6.40) -- (4.10,-0.40) -- (1.70,-0.40) -- cycle;
\node[above,font=\scriptsize] at (3.50,6.40) {$C_2$};
\draw[line width=5pt,gray!40] (0.00,1.00) -- (1.70,1.00);
\draw[->] (0.15,1.00) -- (1.55,1.00);
\node[above,font=\scriptsize] at (0.85,1.10) {$1$};
\draw[line width=5pt,gray!40] (0.00,3.00) -- (1.70,3.00);
\draw[->] (0.15,3.00) -- (1.55,3.00);
\node[above,font=\scriptsize] at (0.85,3.10) {$2$};
\draw[line width=5pt,gray!40] (0.00,5.00) -- (1.70,5.00);
\draw[->] (0.15,5.00) -- (1.55,5.00);
\node[above,font=\scriptsize] at (0.85,5.10) {$3$};
\end{tikzpicture}\hspace{1.2cm}
\begin{tikzpicture}[>=Stealth,scale=0.6]
\draw[rounded corners=6pt,postaction={decorate},decoration={markings,mark=at position 0.12 with {\arrow{>}}}] (0.00,-1.50) -- (0.00,9.90) -- (8.20,9.90) -- (8.20,-1.50) -- (0.00,-1.50) -- cycle;
\node[above,font=\scriptsize] at (7.60,9.90) {$C_1$};
\draw[rounded corners=6pt,postaction={decorate},decoration={markings,mark=at position 0.12 with {\arrow{>}}}] (2.20,-0.95) -- (2.20,9.15) -- (7.50,9.15) -- (7.50,-0.95) -- (2.20,-0.95) -- cycle;
\node[above,font=\scriptsize] at (6.90,9.15) {$C_2$};
\draw[rounded corners=6pt,postaction={decorate},decoration={markings,mark=at position 0.12 with {\arrow{>}}}] (4.40,-0.40) -- (4.40,8.40) -- (6.80,8.40) -- (6.80,-0.40) -- (4.40,-0.40) -- cycle;
\node[above,font=\scriptsize] at (6.20,8.40) {$C_3$};
\draw[line width=5pt,gray!40] (0.00,1.00) -- (2.20,1.00);
\draw[->] (0.15,1.00) -- (2.05,1.00);
\node[above,font=\scriptsize] at (1.10,1.10) {$1$};
\draw[line width=5pt,gray!40] (2.20,3.00) -- (4.40,3.00);
\draw[->] (4.25,3.00) -- (2.35,3.00);
\node[above,font=\scriptsize] at (3.30,3.10) {$2$};
\draw[line width=5pt,gray!40] (0.00,5.00) -- (2.20,5.00);
\draw[->] (0.15,5.00) -- (2.05,5.00);
\node[above,font=\scriptsize] at (1.10,5.10) {$3$};
\draw[line width=5pt,gray!40] (2.20,7.00) -- (4.40,7.00);
\draw[->] (4.25,7.00) -- (2.35,7.00);
\node[above,font=\scriptsize] at (3.30,7.10) {$4$};
\end{tikzpicture}
\caption{The Seifert circles of the two diagrams of Figure~\ref{fig.diagrams} (the circle $C_p$ is the strand in position $p$, closed up) and the bands of the canonical Seifert surface, one for each crossing $j$, drawn as an arrow from $A_j$ to $B_j$. The Seifert graph $\Gamma_D$ has the circles as vertices and the bands as edges: two vertices joined by three parallel edges for the trefoil, and a path $C_1-C_2-C_3$ with doubled edges for the figure-eight knot.}
\label{fig.circles}
\end{figure}
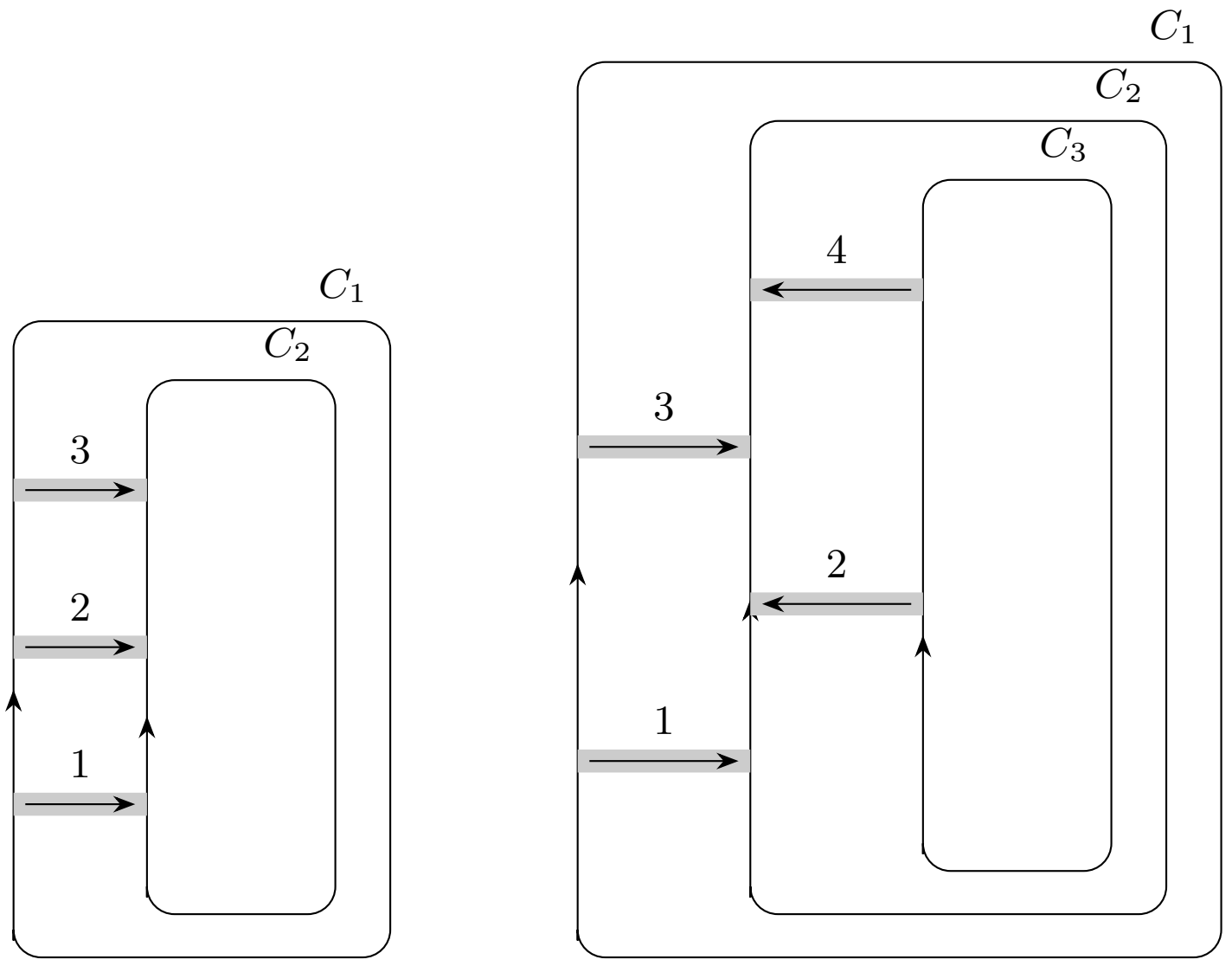

\begin{figure}[ht]
\centering
\begin{tikzpicture}[>=Stealth,scale=0.7]
\def\trefoilbase{%
   \draw (0,0) -- (0,0.5); \draw[->] (0,1.5) -- (0,2.5); \draw[->] (0,3.5) -- (0,4.5); \draw (0,5.5) -- (0,6);
   \draw[rounded corners=5pt] (0,6) -- (0,7.1) -- (4.0,7.1) -- (4.0,-1.1) -- (0,-1.1) -- (0,0);
   \draw (1.7,0) -- (1.7,0.5); \draw[->] (1.7,1.5) -- (1.7,2.5); \draw[->] (1.7,3.5) -- (1.7,4.5); \draw (1.7,5.5) -- (1.7,6);
   \draw[rounded corners=5pt] (1.7,6) -- (1.7,6.5) -- (2.9,6.5) -- (2.9,-0.5) -- (1.7,-0.5) -- (1.7,0);
   \foreach \y in {1,3,5} { \draw (1.7,{\y-0.5}) -- (0,{\y+0.5}); \draw[preaction={draw=white,line width=3pt}] (0,{\y-0.5}) -- (1.7,{\y+0.5}); }}
\trefoilbase
\foreach \y/\j in {1/1,3/2,5/3} \node[font=\scriptsize,circle,inner sep=0.6pt,fill=white] at (0.3,\y) {$\j$};
\node[font=\small] at (-1.9,5.2) {$f_0$};
\node[font=\small] at (0.85,2) {$g_1$};
\node[font=\small] at (0.85,4) {$g_2$};
\node[font=\small] at (3.45,3) {$g_3$};
\node[font=\small] at (2.3,3) {$f_\infty$};
\node[font=\scriptsize,left] at (0,2.0) {$u'_1=o_2$};
\node[font=\scriptsize,left] at (0,4.0) {$u'_2=o_3$};
\node[font=\scriptsize,right] at (1.7,1.75) {$u_2$};
\node[font=\scriptsize,right] at (1.7,4.25) {$u_3$};
\node[font=\scriptsize,above] at (2.0,7.1) {$\hate=u'_3=o_1$};
\node[font=\scriptsize,above] at (2.35,6.5) {$u_1$};
\begin{scope}[xshift=6.2cm,yshift=3cm]
\node[font=\small,anchor=west] at (0,2.7) {$w(f_0)=0,\qquad w(g_j)=-1,\qquad w(f_\infty)=-2$};
\node[font=\small,anchor=west] at (0,1.2) {$f^\uparrow_j=g_j,\qquad f^\downarrow_j=g_{j-1},\qquad f^\leftarrow_j=f_0,\qquad f^\rightarrow_j=f_\infty$};
\node[font=\small,anchor=west] at (0,-0.3) {$a_j=x_j(f_0-g_j)+(f_\infty-g_{j-1})$};
\node[font=\small,anchor=west] at (0,-1.8) {$b_j=y_j(f_0-g_{j-1})+(f_\infty-g_j)$};
\end{scope}
\end{tikzpicture}
\caption{The faces of the standard diagram of the trefoil (Examples~\ref{ex.trefoilflows} and~\ref{ex.trefoilkashaev}): the unbounded face $f_0$, the bigons $g_1,g_2,g_3$ between the two strands and the face $f_\infty$, with their winding numbers, the four faces at the crossing $j$ (with $g_0=g_3$) and the rows $a_j$, $b_j$ of the two face matrices. The two faces adjacent to the cut edge $\hate=u'_3=o_1$ are $f_0$ and $g_3$.}
\label{fig.trefoilfaces}
\end{figure}

\begin{example}[the trefoil]\label{ex.trefoilflows}
Let $D$ be the closure of the braid $\sigma_1^3$, the standard diagram of the right-handed trefoil, with the crossings $1,2,3$ in the order of the braid word (Figure~\ref{fig.diagrams}, left, and Figure~\ref{fig.trefoilfaces}). All crossings are positive, and with the indices taken modulo $3$ the six edges are $u'_j=o_{j+1}$ and $o'_j=u_{j+1}$. There are two Seifert circles (Figure~\ref{fig.circles}, left): $C_1$, which contains the edges $o_j$, and $C_2$, which contains the edges $u_j$; so $A_j=C_1$ and $B_j=C_2$ for all $j$, and $\nu=2$. The five faces are the unbounded face $f_0$, outside $C_1$, the face $f_\infty$ inside $C_2$, and the three bigons $g_1,g_2,g_3$ between the two strands, $g_j$ lying between the crossings $j$ and $j+1$. Their boundaries are
\[
   \partial f_0=e_{o_1}+e_{o_2}+e_{o_3},\qquad \partial f_\infty=-(e_{u_1}+e_{u_2}+e_{u_3}),\qquad \partial g_j=e_{o'_j}-e_{u'_j} ,
\]
the two circles being clockwise. Take the cut edge $\hate=u'_3=o_1$. The faces adjacent to $\hate$ are $f_0$ and $g_3$, and the basis of the faces of $\FF_{D,\hate}$ is $\partial g_1,\partial g_2,\partial f_\infty$. Since $\partial g_1=e_{u_2}-e_{u'_1}$, the crossing forms take the values $\pi_1(\partial g_1)=-x_1$, $\pi_2(\partial g_1)=-1$ and $\pi_3(\partial g_1)=0$ on it; in the same way one finds the matrices of the two flavors in this basis,
\[
   \bigl(\pi_j(z_i)\bigr)=\begin{pmatrix}-x_1&-1&0\\0&-x_2&-1\\1&1&1\end{pmatrix},\qquad
   \bigl(\tpi_j(z_i)\bigr)=\begin{pmatrix}-1&-y_2&0\\0&-1&-y_3\\1&1&1\end{pmatrix},
\]
which are the transposes of $Q_{D,\hate}(x)$ and of $\tQ_{D,\hate}(y)$. Their determinants are
\[
   P_{D,\hate}(x)=1-x_1+x_1x_2\qquad\text{and}\qquad\tP_{D,\hate}(y)=1-y_3+y_2y_3 ,
\]
with $\kappa=1$. For the Seifert basis we take the spanning tree $\CT=\{1\}$ of the Seifert graph, which has two vertices and three parallel edges, the basis $\gamma_1=e_2-e_1$, $\gamma_2=e_3-e_1$ of its cycle space, and the base edge $u_1$ on $C_2$. The lifts are $z^{\gamma_1}=e_{u'_1}-e_{u_2}=-\partial g_1$ and $z^{\gamma_2}=e_{u'_1}+e_{u'_2}-e_{u_2}-e_{u_3}=-\partial g_1-\partial g_2$: the cycle $\gamma_2$ enters $C_1$ through the band $1$, runs along $C_1$ through the edges $u'_1=o_2$ and $u'_2=o_3$, passes to $C_2$ through the band $3$ and returns to the band $1$ along $C_2$, against the edges $u_3$ and $u_2$. The remaining vector of the Seifert basis is $C_2=-\partial f_\infty$. In this basis
\[
   N_{D,\hate}(x)=\begin{pmatrix}x_1&1&0\\x_1&1+x_2&1\\-1&-1&-1\end{pmatrix},\qquad \det N_{D,\hate}(x)=-(1-x_1+x_1x_2),
\]
in accordance with \eqref{eq.kappaSeifert}: $\kappa=(-1)^{s-1}\det(Y\,|\,E_{\CT})\det(\hat\partial E_{\CT})=(-1)\cdot1\cdot1$, the matrix $(Y\,|\,E_{\CT})$ having the determinant $1$ and $\hat\partial E_{\CT}$ being the $1\times1$ matrix $(\partial e_1)_{C_2}=1$.
\end{example}

\subsection{The coupling on the lattice of flows: a mixed discriminant}\label{ssec.mixed}
In \eqref{eq.LGflows} and \eqref{eq.LGtau} the Links--Gould polynomial is the image under $\BE_D$ of a determinant. Since $\BE_D$ is not multiplicative, this image is not the determinant of a matrix with entries in $\BZ[t_0^{\pm1},t_1^{\pm1}]$. In this subsection we carry $\BE_D$ inside the determinant. It turns the product $\pi_j\otimes\tpi_j$ of the two crossing forms of a crossing into \emph{one} bilinear form $\Psi_j$ on the lattice of flows, whose four coefficients are $1$, $-t_0^{\ve_j}$, $-t_1^{\ve_j}$ and $1$, and the Links--Gould polynomial is obtained from the $c$ forms $\Psi_1,\dots,\Psi_c$ by the multilinear operation which is the polarization of the determinant (Theorem~\ref{thm.mixed}). So the coupling itself can be written on the lattice of flows, in an arbitrary basis and with no sign. The price is that $\LG$ is one coefficient of a determinant and not a determinant: we show at the end of the subsection that the determinant in question is not an invariant of the link (Proposition~\ref{prop.kink} and Example~\ref{ex.pencilunknot}).

\begin{definition}[mixed discriminants]\label{def.mixed}
Let $\kk$ be a commutative ring and let $\mathsf M_1,\dots,\mathsf M_c$ be $c\times c$ matrices with entries in $\kk$. We denote by
\[
   \lb\mathsf M_1,\dots,\mathsf M_c\rb\ \in\kk
\]
the coefficient of $\lambda_1\cdots\lambda_c$ in the polynomial $\det(\lambda_1\mathsf M_1+\dots+\lambda_c\mathsf M_c)\in\kk[\lambda_1,\dots,\lambda_c]$, where $\lambda_1,\dots,\lambda_c$ are indeterminates. If $\mathsf B_1,\dots,\mathsf B_c$ are bilinear forms on $\FF_D$ with values in $\kk$ and $\hate$ is an edge of $D$, we put
\[
   \lb\mathsf B_1,\dots,\mathsf B_c\rb_{\hate}=\lb\mathsf M_1,\dots,\mathsf M_c\rb,\qquad \mathsf M_j=\bigl(\mathsf B_j(z_i,z_k)\bigr)_{i,k=1,\dots,c},
\]
where $z_1,\dots,z_c$ is a $\BZ$-basis of the lattice $\FF_{D,\hate}$.
\end{definition}

Expanding the determinant over the permutations and collecting the terms which contain each of the variables $\lambda_1,\dots,\lambda_c$ exactly once, one finds
\begin{equation}\label{eq.mixedexplicit}
   \lb\mathsf M_1,\dots,\mathsf M_c\rb=\sum_{\sigma,\sigma'}\sgn(\sigma)\,\sgn(\sigma')\prod_{j=1}^c\bigl(\mathsf M_j\bigr)_{\sigma(j),\sigma'(j)},
\end{equation}
the sum running over all the pairs of permutations of $\{1,\dots,c\}$: in the term of $\det\bigl(\sum_j\lambda_j\mathsf M_j\bigr)$ which belongs to a permutation $\varrho$, the monomial $\lambda_1\cdots\lambda_c$ is obtained by taking $\lambda_j\mathsf M_j$ in one row $\sigma(j)$ for each $j$, and its coefficient is $\sgn(\varrho)\prod_j(\mathsf M_j)_{\sigma(j),\varrho\sigma(j)}$; we put $\sigma'=\varrho\sigma$. Up to the factor $c!$, the element $\lb\mathsf M_1,\dots,\mathsf M_c\rb$ is the \emph{mixed discriminant} of the matrices $\mathsf M_j$ \cite{Bapat89}: it is symmetric and multilinear in $\mathsf M_1,\dots,\mathsf M_c$, and $\lb\mathsf M,\dots,\mathsf M\rb=c!\det\mathsf M$. By \eqref{eq.mixedexplicit},
\begin{equation}\label{eq.mixedproperties}
   \lb\mathsf M_1^{\T},\dots,\mathsf M_c^{\T}\rb=\lb\mathsf M_1,\dots,\mathsf M_c\rb\quad\text{and}\quad\lb U\mathsf M_1V,\dots,U\mathsf M_cV\rb=\det U\,\det V\,\lb\mathsf M_1,\dots,\mathsf M_c\rb
\end{equation}
for all $c\times c$ matrices $U,V$; the second identity holds because $\det\bigl(\sum_j\lambda_jU\mathsf M_jV\bigr)$ is the product of $\det U\det V$ and of $\det\bigl(\sum_j\lambda_j\mathsf M_j\bigr)$. In particular $\lb\mathsf B_1,\dots,\mathsf B_c\rb_{\hate}$ does not depend on the basis of $\FF_{D,\hate}$: a change of basis replaces every matrix $\mathsf M_j$ by $U\mathsf M_jU^{\T}$ with $U\in\GL_c(\BZ)$, and $(\det U)^2=1$.

\begin{lemma}\label{lem.cauchybinet}
Suppose that $\mathsf M_j=\sum_{r=1}^{n_j}v_{j,r}\,w_{j,r}^{\T}$ with column vectors $v_{j,r},w_{j,r}\in\kk^c$, for $j=1,\dots,c$. Then
\[
   \lb\mathsf M_1,\dots,\mathsf M_c\rb=\sum_{r_1,\dots,r_c}\det\bigl(v_{1,r_1},\dots,v_{c,r_c}\bigr)\,\det\bigl(w_{1,r_1},\dots,w_{c,r_c}\bigr),
\]
the sum running over all the choices of $r_j\in\{1,\dots,n_j\}$. In particular, for matrices of rank one,
\begin{align*}
   \lb v_1w_1^{\T},\dots,v_cw_c^{\T}\rb&=\det(v_1,\dots,v_c)\,\det(w_1,\dots,w_c),\\
   \det\Bigl(\sum_{j=1}^c\lambda_j\,v_jw_j^{\T}\Bigr)&=\lambda_1\cdots\lambda_c\,\det(v_1,\dots,v_c)\,\det(w_1,\dots,w_c).
\end{align*}
\end{lemma}
\begin{proof}
Both sides of the first formula are multilinear in the matrices $\mathsf M_j$, so it suffices to prove it when $n_j=1$ for all $j$, that is, to prove the last two formulas. Let $\mathsf V$ and $\mathsf W$ be the $c\times c$ matrices with the columns $v_1,\dots,v_c$ and $w_1,\dots,w_c$. Then $\sum_j\lambda_j\,v_jw_j^{\T}=\mathsf V\,\diag(\lambda_1,\dots,\lambda_c)\,\mathsf W^{\T}$, whose determinant is $\lambda_1\cdots\lambda_c\det\mathsf V\det\mathsf W$.
\end{proof}

\begin{definition}[coupled crossing forms]\label{def.coupledforms}
Let $\kk$ be a commutative ring, $j$ a crossing of $D$ and $\sv,\xs,\ys\in\kk$. We denote by $\Psi_j[\sv;\xs,\ys]$ the bilinear form on $\FF_D$ with values in $\kk$
\begin{equation}\label{eq.Psigeneral}
   \Psi_j[\sv;\xs,\ys](z,z')=\sv\,z_{u'_j}z'_{o_j}-\xs\,z_{u'_j}z'_{o'_j}-\ys\,z_{u_j}z'_{o_j}+z_{u_j}z'_{o'_j} .
\end{equation}
The \emph{coupled crossing form} of the crossing $j$ is the form $\Psi_j=\Psi_j[1;t_0^{\ve_j},t_1^{\ve_j}]$ with values in $\BZ[t_0^{\pm1},t_1^{\pm1}]$:
\begin{equation}\label{eq.Psi}
   \Psi_j(z,z')=z_{u'_j}z'_{o_j}-t_0^{\ve_j}\,z_{u'_j}z'_{o'_j}-t_1^{\ve_j}\,z_{u_j}z'_{o_j}+z_{u_j}z'_{o'_j} .
\end{equation}
\end{definition}

For an element $\xs$ of a commutative ring we write $\pi_j^{(\xs)}$ for the crossing form $\pi_j$ in which the variable $x_j$ has been replaced by $\xs$, and likewise $\tpi_j^{(\ys)}$; thus $\pi_j^{(\xs)}(z)=\xs\,z_{u'_j}-z_{u_j}$ and $\tpi_j^{(\ys)}(z)=\ys\,z_{o_j}-z_{o'_j}$. By \eqref{eq.crossingforms}, the product of the two crossing forms of the crossing $j$ is
\[
   \pi_j(z)\,\tpi_j(z')=x_jy_j\,z_{u'_j}z'_{o_j}-x_j\,z_{u'_j}z'_{o'_j}-y_j\,z_{u_j}z'_{o_j}+z_{u_j}z'_{o'_j},
\]
so that $\Psi_j[\sv;\xs,\ys](z,z')$ is the image of $\pi_j(z)\tpi_j(z')$ under the linear map which sends $x_j$ to $\xs$, $y_j$ to $\ys$ and $x_jy_j$ to $\sv$. In particular
\begin{equation}\label{eq.rankone}
   \Psi_j[\xs\ys;\xs,\ys]=\pi_j^{(\xs)}\otimes\tpi_j^{(\ys)}\qquad\text{and}\qquad\Psi_j(z,z')=\BE_D\bigl[\pi_j(z)\,\tpi_j(z')\bigr].
\end{equation}
In the coordinates $(z_{u'_j},z_{u_j})$ of $z$ and $(z'_{o_j},z'_{o'_j})$ of $z'$, the form $\Psi_j[\sv;\xs,\ys]$ has the matrix $\bigl(\begin{smallmatrix}\sv&-\xs\\-\ys&1\end{smallmatrix}\bigr)$, of determinant $\sv-\xs\ys$: the first variable of the form only sees the under-strand of the crossing, and the second one its over-strand.

\begin{theorem}[the Links--Gould polynomial as a mixed discriminant]\label{thm.mixed}
Let $D$ be a diagram of a link $L$ with $\mu$ components and $\hate$ an edge of $D$.
\begin{enumerate}\renewcommand{\theenumi}{\roman{enumi}}
\item Let $\kk$ be a commutative ring and $\sv_j,\xs_j,\ys_j\in\kk$ for $j=1,\dots,c$. Let $\BE'$ be the $\BZ$-linear map, from the polynomials of degree $\le1$ in each of the variables $x_j,y_j$ with integer coefficients to $\kk$, such that
\[
   \BE'\Bigl[\,\prod_{j\in S}x_j\prod_{j\in S'}y_j\Bigr]=\prod_{j\in S\sm S'}\xs_j\ \prod_{j\in S'\sm S}\ys_j\ \prod_{j\in S\cap S'}\sv_j\qquad(S,S'\subseteq\{1,\dots,c\}),
\]
that is, $x_j\mapsto\xs_j$, $y_j\mapsto\ys_j$ and $x_jy_j\mapsto\sv_j$. Then
\[
   \bigl\lb\,\Psi_1[\sv_1;\xs_1,\ys_1],\dots,\Psi_c[\sv_c;\xs_c,\ys_c]\,\bigr\rb_{\hate}=(-1)^{\mu-1}\,\BE'\bigl[P_{D,\hate}(x)\,\tP_{D,\hate}(y)\bigr].
\]
\item The Links--Gould polynomial is the mixed discriminant of the coupled crossing forms:
\[
   \LG(L;t_0,t_1)=\lb\Psi_1,\dots,\Psi_c\rb_{\hate} .
\]
In other words, for every basis $z_1,\dots,z_c$ of $\FF_{D,\hate}$, the polynomial $\LG(L;t_0,t_1)$ is the coefficient of $\lambda_1\cdots\lambda_c$ in
\[
   \det\Bigl(\sum_{j=1}^c\lambda_j\,\Psi_j(z_i,z_k)\Bigr)_{i,k=1,\dots,c} .
\]
\end{enumerate}
\end{theorem}
\begin{proof}
(i) Let $z_1,\dots,z_c$ be a basis of $\FF_{D,\hate}$. By Theorem~\ref{thm.flows}, in which the sign $\kappa$ occurs twice,
\begin{align*}
   (-1)^{\mu-1}P_{D,\hate}(x)\,\tP_{D,\hate}(y)&=\det\bigl(\pi_j(z_i)\bigr)\det\bigl(\tpi_j(z_i)\bigr)\\
   &=\sum_{\sigma,\sigma'}\sgn(\sigma)\,\sgn(\sigma')\prod_{j=1}^c\pi_j(z_{\sigma(j)})\,\tpi_j(z_{\sigma'(j)}),
\end{align*}
the sum running over the pairs of permutations of $\{1,\dots,c\}$. The factor $f_j=\pi_j(z_{\sigma(j)})\tpi_j(z_{\sigma'(j)})$ is a polynomial of degree $\le1$ in $x_j$ and in $y_j$ which contains no other variable. The map $\BE'$ is multiplicative over distinct crossings, $\BE'[\prod_jf_j]=\prod_j\BE'[f_j]$ for such polynomials, both sides being multilinear in the $f_j$ and equal on the monomials; and $\BE'[f_j]=\Psi_j[\sv_j;\xs_j,\ys_j](z_{\sigma(j)},z_{\sigma'(j)})$, as we have just seen. So the image of the right-hand side under $\BE'$ is the right-hand side of \eqref{eq.mixedexplicit} for the matrices of the forms $\Psi_j[\sv_j;\xs_j,\ys_j]$ in the basis $z_1,\dots,z_c$.

(ii) For $\kk=\BZ[t_0^{\pm1},t_1^{\pm1}]$ and $\sv_j=1$, $\xs_j=t_0^{\ve_j}$, $\ys_j=t_1^{\ve_j}$, the map $\BE'$ is the restriction of $\BE_D$ to the polynomials with integer coefficients, by \eqref{eq.ED}, and the forms are the coupled crossing forms. So (ii) follows from (i) and Theorem~\ref{thm.HK}.
\end{proof}

Thus $\BE_D$ acts on the coupled form $\tau=\sum_j\ve_j\,\pi_j\otimes\tpi_j$ of Remark~\ref{rem.flows}~(iii) term by term: it replaces the form $\pi_j\otimes\tpi_j$, of rank one, by the form $\Psi_j$, and the discriminant by the mixed discriminant. There is no sign in the formula: by Lemma~\ref{lem.cauchybinet}, the discriminant of $\tau$ is $\ve_D$ times the mixed discriminant of the forms $\pi_j\otimes\tpi_j$, which removes the sign $\ve_D$ of \eqref{eq.LGtau}. The three models of Parts~\ref{part.surfaces} and~\ref{part.faces} use Theorem~\ref{thm.flows} for each flavor separately, with $\BE_D$ applied at the end; Theorem~\ref{thm.mixed} will be used in Section~\ref{sec.cocores}, where the forms $\Psi_j$ are read off the complement of the link and of the co-cores of the bands of its Seifert surface.

\begin{example}[the trefoil]\label{ex.trefoilmixed}
Let $D$ be the diagram of the trefoil of Example~\ref{ex.trefoilflows}, with the cut edge $\hate=u'_3=o_1$ and the basis $\partial g_1,\partial g_2,\partial f_\infty$ of $\FF_{D,\hate}$. The column $j$ of the matrices $\bigl(\pi_j(z_i)\bigr)$ and $\bigl(\tpi_j(z_i)\bigr)$ of that example gives the values of $\pi_j$ and of $\tpi_j$ on the basis, and by \eqref{eq.rankone} the matrices of the coupled crossing forms are obtained from the products of these columns by $x_j\mapsto t_0$, $y_j\mapsto t_1$, $x_jy_j\mapsto1$:
\[
   \Psi_1=\begin{pmatrix}t_0&0&-t_0\\0&0&0\\-1&0&1\end{pmatrix},\qquad
   \Psi_2=\begin{pmatrix}t_1&1&-1\\1&t_0&-t_0\\-t_1&-1&1\end{pmatrix},\qquad
   \Psi_3=\begin{pmatrix}0&0&0\\0&t_1&-1\\0&-t_1&1\end{pmatrix}.
\]
The first and the third matrix have rank one: since $\hate=o_1=u'_3$, the forms $\tpi_1$ and $\pi_3$ contain no variable on $\FF_{D,\hate}$, so that $\BE_D$ acts on $\pi_1\otimes\tpi_1$ and on $\pi_3\otimes\tpi_3$ as an evaluation. The second matrix has rank two. One finds
\begin{align*}
   \det\bigl(\lambda_1\Psi_1+\lambda_2\Psi_2+\lambda_3\Psi_3\bigr)
   &=\det\begin{pmatrix}\lambda_1t_0+\lambda_2t_1&\lambda_2&-\lambda_1t_0-\lambda_2\\ \lambda_2&\lambda_2t_0+\lambda_3t_1&-\lambda_2t_0-\lambda_3\\ -\lambda_1-\lambda_2t_1&-\lambda_2-\lambda_3t_1&\lambda_1+\lambda_2+\lambda_3\end{pmatrix}\\
   &=\lambda_1\lambda_2\lambda_3\,\LG(3_1)+(1-t_0t_1)\,\lambda_2^2\,\bigl((t_0-1)\lambda_1+(t_1-1)\lambda_3\bigr),
\end{align*}
where $\LG(3_1)=1-t_0-t_1+t_0^2+2t_0t_1+t_1^2-t_0^2t_1-t_0t_1^2$ is the Links--Gould polynomial of the right-handed trefoil, which is also $\BE_D[P_{D,\hate}\tP_{D,\hate}]=\BE_D[(1-x_1+x_1x_2)(1-y_3+y_2y_3)]$.
\end{example}

\begin{remark}\label{rem.mixed}
(i) \emph{The rank, and the free fermion point.} With the notation \eqref{eq.xbar}, the form $\Psi_j$ comes from the matrix $\bigl(\begin{smallmatrix}1&-\bx_j\\-\by_j&1\end{smallmatrix}\bigr)$, whose determinant is $g_j=1-(t_0t_1)^{\ve_j}$. So $\Psi_j$ has rank at most two, and at $t_0t_1=1$ it is the form $\pi_j^{(\bx_j)}\otimes\tpi_j^{(\by_j)}$ of rank one, by \eqref{eq.rankone}: there the mixed discriminant is the product of the two determinants $\det\bigl(\pi_j^{(\bx_j)}(z_i)\bigr)$ and $\det\bigl(\tpi_j^{(\by_j)}(z_i)\bigr)$, by Lemma~\ref{lem.cauchybinet}, which is the free fermion point of Section~\ref{sec.formula}.

(ii) \emph{The expansion around the free fermion point.} In general
\begin{equation}\label{eq.Psisplit}
   \Psi_j(z,z')=\pi_j^{(\bx_j)}(z)\,\tpi_j^{(\by_j)}(z')+g_j\,z_{u'_j}\,z'_{o_j} ,
\end{equation}
a sum of two forms of rank one. By Lemma~\ref{lem.cauchybinet}, the mixed discriminant of Theorem~\ref{thm.mixed} is therefore the sum, over the subsets $S$ of $\{1,\dots,c\}$, of $\prod_{j\in S}g_j$ times the product of two determinants: that of the matrix whose column $j$ is $\bigl((z_i)_{u'_j}\bigr)_i$ if $j\in S$ and $\bigl(\pi_j^{(\bx_j)}(z_i)\bigr)_i$ if $j\notin S$, and that of the matrix whose column $j$ is $\bigl((z_i)_{o_j}\bigr)_i$ if $j\in S$ and $\bigl(\tpi_j^{(\by_j)}(z_i)\bigr)_i$ if $j\notin S$. This is the expansion of Corollary~\ref{cor.expansion}, and it identifies its coefficients: since $\pi_j=\pi_j^{(\bx_j)}+(x_j-\bx_j)\,z_{u'_j}$ and $\tpi_j=\tpi_j^{(\by_j)}+(y_j-\by_j)\,z_{o_j}$, the two determinants are $\kappa\,a_S(t_0)$ and $(-1)^{\mu-1}\kappa\,b_S(t_1)$, by Theorem~\ref{thm.flows} and by the definition \eqref{eq.ab} of $a_S$ and $b_S$. So $a_S$ and $b_S$ are, up to sign, the determinants of the crossing forms at the free fermion point in which, at the crossings of $S$, the crossing form has been replaced by the flow through the outgoing under-edge, resp.\ through the incoming over-edge.

(iii) \emph{An evaluation and a projection.} The map $\BE_D$ is not an evaluation of the crossing variables. Theorem~\ref{thm.mixed} writes it as an evaluation followed by a projection: the forms $\Psi_j$ are defined over $\BZ[t_0^{\pm1},t_1^{\pm1}]$, and the Links--Gould polynomial is the part of the determinant of their sum in which every crossing has been used exactly once.
\end{remark}

The mirror image and the inversion of a link act on the coupled crossing forms in a simple way, which gives a proof of two known symmetries of $\LG$ \cite{DeWitKauffmanLinks, Ishii04}.

\begin{corollary}[mirror image and inverse]\label{cor.symmetries}
Let $L^*$ be the mirror image of $L$ and let $-L$ be the link $L$ with the orientations of all its components reversed. Then
\[
   \LG(L^*;t_0,t_1)=\LG(L;t_0^{-1},t_1^{-1})\qquad\text{and}\qquad\LG(-L;t_0,t_1)=\LG(L;t_1,t_0).
\]
\end{corollary}
\begin{proof}
Let $D^*$ be the diagram obtained from $D$ by changing all the crossings. It is a diagram of $L^*$ with the same edges and the same flows as $D$ and with the opposite signs, and at every crossing the over-strand and the under-strand are exchanged: the edges $o_j,u_j,o'_j,u'_j$ of $D^*$ are the edges $u_j,o_j,u'_j,o'_j$ of $D$. By \eqref{eq.Psi} the coupled crossing form of the crossing $j$ of $D^*$ is
\[
   z_{o'_j}z'_{u_j}-t_0^{-\ve_j}\,z_{o'_j}z'_{u'_j}-t_1^{-\ve_j}\,z_{o_j}z'_{u_j}+z_{o_j}z'_{u'_j},
\]
in the notation of $D$. This is the value at $(z',z)$ of the form $\Psi_j$ of $D$ in which $t_0,t_1$ have been replaced by $t_0^{-1},t_1^{-1}$. So the matrices of the forms of $D^*$ are the transposes of those of $D$ with the inverted variables, and the first identity follows from \eqref{eq.mixedproperties} and Theorem~\ref{thm.mixed}. Let $-D$ be the diagram $D$ with all the orientations reversed, a diagram of $-L$. It has the same signs as $D$ and, when its edges are identified with those of $D$, the same flows; and at every crossing the incoming and the outgoing edges are exchanged. So its coupled crossing form at the crossing $j$ is $z_{u_j}z'_{o'_j}-t_0^{\ve_j}z_{u_j}z'_{o_j}-t_1^{\ve_j}z_{u'_j}z'_{o'_j}+z_{u'_j}z'_{o_j}$, in the notation of $D$, which is the form $\Psi_j$ of $D$ with $t_0$ and $t_1$ exchanged. The second identity follows from Theorem~\ref{thm.mixed}.
\end{proof}

Since $\LG$ does not distinguish a link from its inverse \cite{DeWitKauffmanLinks}, the second identity is the symmetry of $\LG$ in its two variables, and the two identities together give
\begin{equation}\label{eq.mirror}
   \LG(L;t_1^{-1},t_0^{-1})=\LG(-L^*;t_0,t_1)=\LG(L^*;t_0,t_1) .
\end{equation}

\emph{The determinant of the pencil.} It is natural to ask whether the whole polynomial of Theorem~\ref{thm.mixed}~(ii), and not only one of its coefficients, is an invariant of the link. It is not.

\begin{definition}\label{def.pencil}
For a diagram $D$, an edge $\hate$ of $D$ and a basis $z_1,\dots,z_c$ of $\FF_{D,\hate}$, the \emph{determinant of the pencil} of the coupled crossing forms is
\[
   \CDp_{D,\hate}(\lambda)=\det\Bigl(\sum_{j=1}^c\lambda_j\,\Psi_j(z_i,z_k)\Bigr)_{i,k=1,\dots,c}\ \in\ \BZ[t_0^{\pm1},t_1^{\pm1}][\lambda_1,\dots,\lambda_c] .
\]
More generally, if $\FF'$ is a subgroup of $\FF_D$ which is a direct summand of rank $n$, we write $\CDp_{\FF'}(\lambda)$ for the determinant of the $n\times n$ matrix $\bigl(\sum_j\lambda_j\Psi_j(z_i,z_k)\bigr)$, where $z_1,\dots,z_n$ is a basis of $\FF'$, with $\CDp_{\FF'}=1$ if $n=0$.
\end{definition}

These determinants do not depend on the bases, a change of basis multiplying them by the square of a determinant $\pm1$. The polynomial $\CDp_{D,\hate}(\lambda)$ is homogeneous of degree $c$ in $\lambda$, and its degree in each $\lambda_j$ is at most $2$, because the matrix of $\Psi_j$ has rank at most two. By Theorem~\ref{thm.mixed}, the coefficient of $\lambda_1\cdots\lambda_c$ is $\LG(L)$. All its other coefficients are divisible by $1-t_0t_1$: at $t_1=t_0^{-1}$ the forms $\Psi_j$ have rank one (Remark~\ref{rem.mixed}~(i)), and by Lemma~\ref{lem.cauchybinet} the determinant of the pencil is then $\lambda_1\cdots\lambda_c$ times an element of $\BZ[t_0^{\pm1}]$. In Example~\ref{ex.trefoilmixed}, $\CDp_{D,\hate}(\lambda)$ has three monomials.

\begin{proposition}[a move of type I]\label{prop.kink}
Let $e\ne\hate$ be an edge of $D$, and let $D^+$ be a diagram obtained from $D$ by a Reidemeister move of type I on the edge $e$, which creates a new crossing, numbered $0$, of sign $\ve_0$; we regard $\hate$ as an edge of $D^+$. Let $\FF_{D,\hate,e}$ be the lattice of the flows of $D$ which vanish at $\hate$ and at $e$. Then
\[
   \CDp_{D^+,\hate}(\lambda,\lambda_0)=\lambda_0\,\CDp_{D,\hate}(\lambda)+\lambda_0^2\,\bigl(1-(t_0t_1)^{\ve_0}\bigr)\,\CDp_{\FF_{D,\hate,e}}(\lambda)
\]
if the two faces adjacent to $e$ are not the two faces adjacent to $\hate$, and $\CDp_{D^+,\hate}(\lambda,\lambda_0)=\lambda_0\,\CDp_{D,\hate}(\lambda)$ otherwise.
\end{proposition}

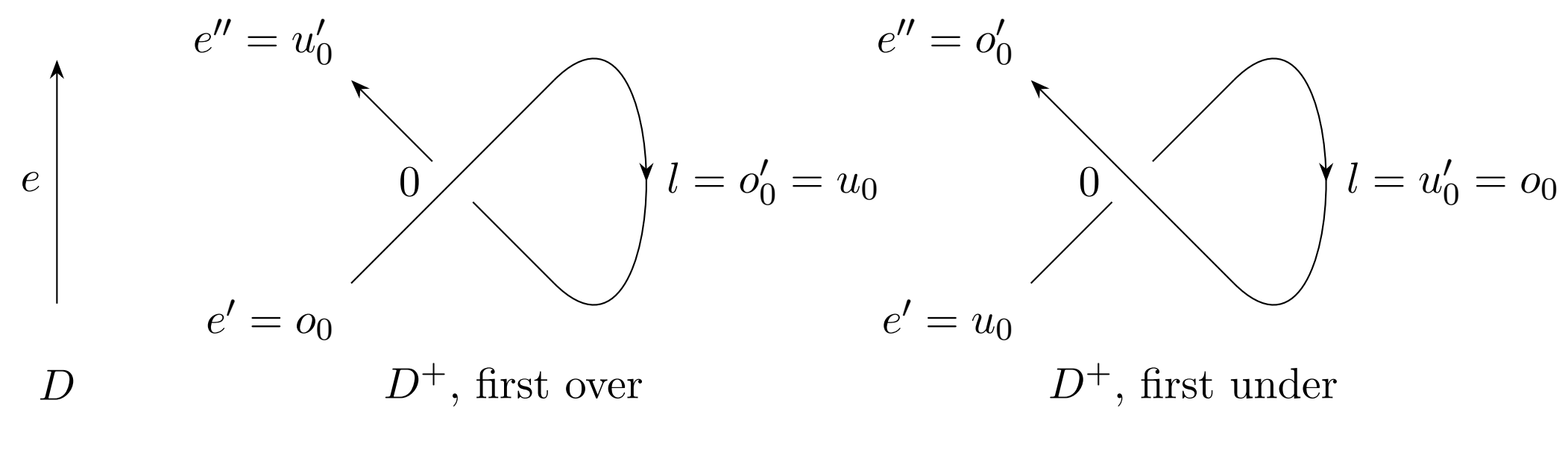
\begin{figure}[ht]
\centering
\begin{tikzpicture}[>=Stealth,scale=0.9,every node/.style={font=\small},
  mid/.style={postaction={decorate},decoration={markings,mark=at position 0.5 with {\arrow{>}}}}]
\begin{scope}
\draw[->] (0,-1.2) -- (0,1.2);
\node[left] at (0,0) {$e$};
\node at (0,-2.0) {$D$};
\end{scope}
\begin{scope}[xshift=3.9cm]
\draw (1,-1) -- (0.2,-0.2); \draw[->] (-0.2,0.2) -- (-1,1);
\draw (-1,-1) -- (1,1);
\draw[mid] (1,1) to[out=45,in=-45,looseness=2.2] (1,-1);
\node[below left] at (-1,-1) {$e'=o_0$}; \node[above left] at (-1,1) {$e''=u'_0$};
\node[right] at (1.95,0) {$l=o'_0=u_0$};
\node[left=4pt] at (0,0) {$0$};
\node at (0.6,-2.0) {$D^+$, first over};
\end{scope}
\begin{scope}[xshift=10.6cm]
\draw (-1,-1) -- (-0.2,-0.2); \draw (0.2,0.2) -- (1,1);
\draw[->] (1,-1) -- (-1,1);
\draw[mid] (1,1) to[out=45,in=-45,looseness=2.2] (1,-1);
\node[below left] at (-1,-1) {$e'=u_0$}; \node[above left] at (-1,1) {$e''=o'_0$};
\node[right] at (1.95,0) {$l=u'_0=o_0$};
\node[left=4pt] at (0,0) {$0$};
\node at (0.6,-2.0) {$D^+$, first under};
\end{scope}
\end{tikzpicture}
\caption{A Reidemeister move of type I on the edge $e$ of $D$ (Proposition~\ref{prop.kink}): it replaces $e$ by the edges $e'$ and $e''$ and by the loop $l$, and it creates the crossing $0$. In the middle the strand passes first over itself, on the right first under itself. Here the loop lies on the right of the strand, and the sign of the new crossing is $\ve_0=1$ in the middle and $\ve_0=-1$ on the right; when the loop lies on the left of the strand, the two signs are exchanged.}
\label{fig.kink}
\end{figure}

\begin{proof}
The move replaces the edge $e$ by three edges (Figure~\ref{fig.kink}): an edge $e'$ which enters the new crossing, an edge $e''$ which leaves it, and a loop $l$ which leaves the new crossing and comes back to it. If the strand passes first over itself, then $(o_0,u_0,o'_0,u'_0)=(e',l,l,e'')$, and if it passes first under itself, then $(o_0,u_0,o'_0,u'_0)=(l,e',e'',l)$. In both cases the condition \eqref{eq.flow} at the new crossing is $z_{e'}=z_{e''}$. So the flows of $D^+$ are the vectors obtained from a flow $z$ of $D$ by putting $z_{e'}=z_{e''}=z_e$ and by giving an arbitrary value to $z_l$, and
\[
   \FF_{D^+,\hate}=\FF_{D,\hate}\oplus\BZ\,z_0,
\]
where a flow of $D$ is extended by $z_l=0$ and $z_0$ is the flow of $D^+$ which is $1$ at $l$ and $0$ at the other edges. For a flow $z$ of $D^+$ we write $z_e$ for the common value of $z_{e'}$ and $z_{e''}$. By \eqref{eq.Psi}, in the first case
\[
   \Psi_0(z,z')=z_ez'_e-t_0^{\ve_0}\,z_ez'_l-t_1^{\ve_0}\,z_lz'_e+z_lz'_l=\bigl(z_l-t_0^{\ve_0}z_e\bigr)\bigl(z'_l-t_1^{\ve_0}z'_e\bigr)+g_0\,z_ez'_e
\]
with $g_0=1-(t_0t_1)^{\ve_0}$, and in the second case $\Psi_0(z,z')=z_lz'_l-t_0^{\ve_0}z_lz'_e-t_1^{\ve_0}z_ez'_l+z_ez'_e=\bigl(z_l-t_1^{\ve_0}z_e\bigr)\bigl(z'_l-t_0^{\ve_0}z'_e\bigr)+g_0\,z_ez'_e$. So in both cases $\Psi_0(z,z')=h(z)\,h'(z')+g_0\,z_ez'_e$ with two linear forms $h,h'$ on the flows of $D^+$ such that $h(z_0)=h'(z_0)=1$; the forms $\Psi_j$ with $j\ne0$ do not involve the edge $l$, and on the flows of $D$ they are the coupled crossing forms of $D$. Let $z_1,\dots,z_c$ be a basis of $\FF_{D,\hate}$, let $\mathsf M=\bigl(\sum_{j\ne0}\lambda_j\Psi_j(z_i,z_k)\bigr)$, of determinant $\CDp_{D,\hate}(\lambda)$, and let $v$ be the column vector with the coordinates $v_i=(z_i)_e$. In the basis $z_0,z_1,\dots,z_c$ of $\FF_{D^+,\hate}$ the matrix of the pencil is
\[
   \begin{pmatrix}\lambda_0&\lambda_0\,h'(z_k)\\[2pt] \lambda_0\,h(z_i)&\mathsf M_{ik}+\lambda_0\,h(z_i)\,h'(z_k)+\lambda_0\,g_0\,v_iv_k\end{pmatrix}.
\]
Subtracting $h(z_i)$ times the first row from the row of $z_i$, for every $i$, gives
\[
   \CDp_{D^+,\hate}(\lambda,\lambda_0)=\lambda_0\,\det\bigl(\mathsf M+\lambda_0\,g_0\,vv^{\T}\bigr).
\]
If the two faces adjacent to $e$ are the two faces adjacent to $\hate$, then the boundaries of the other faces vanish at $e$, so $v=0$ by Lemma~\ref{lem.faces}~(iii), and we are done. Otherwise one of the faces adjacent to $e$ is not adjacent to $\hate$, its boundary is an element of $\FF_{D,\hate}$ whose coordinate at $e$ is $\pm1$, and the linear form $z\mapsto z_e$ maps $\FF_{D,\hate}$ onto $\BZ$. Its kernel $\FF_{D,\hate,e}$ is then a direct summand of rank $c-1$, and we can choose the basis so that $(z_1)_e=1$ and $z_2,\dots,z_c$ is a basis of $\FF_{D,\hate,e}$. Then $vv^{\T}$ is the matrix unit $E_{11}$, and $\det(\mathsf M+\lambda_0g_0E_{11})=\det\mathsf M+\lambda_0g_0\det\mathsf M_{\hat1}$, where $\mathsf M_{\hat1}$, the matrix $\mathsf M$ with its first row and its first column removed, is the matrix of the pencil of $D$ on $\FF_{D,\hate,e}$.
\end{proof}

\begin{example}[two diagrams of the unknot]\label{ex.pencilunknot}
Let $D_0$ be the closure of the braid $\sigma_1$ on two strands, a diagram of the unknot with one crossing, and let $e_1,e_2$ be its two edges, which close the first and the second strand. For each of them the lattice $\FF_{D_0,e_i}$ is generated by the other one, on which $\Psi_1$ takes the value $1$ by \eqref{eq.Psi}; so $\CDp_{D_0,e_1}=\CDp_{D_0,e_2}=\lambda_1$. Let $D$ be the closure of the braid $\sigma_1\sigma_2$ on three strands, with its four edges: the edges $e_1,e_2,e_3$ which close the three strands, and the edge $e_4$ between the two crossings. It is obtained from $D_0$ by a move of type I on $e_2$, which creates the crossing $2$, and also from the closure of $\sigma_2$ by a move of type I which creates the crossing $1$. The lattice of the flows of $D_0$ which vanish at $e_1$ and at $e_2$ is zero, and the faces adjacent to $e_1$ and to $e_2$ are not the same. By Proposition~\ref{prop.kink} and by the symmetry of the diagram,
\[
   \CDp_{D,e_1}(\lambda)=\lambda_1\lambda_2+(1-t_0t_1)\,\lambda_2^2,\qquad \CDp_{D,e_3}(\lambda)=\lambda_1\lambda_2+(1-t_0t_1)\,\lambda_1^2,
\]
whereas $\CDp_{D,e_2}(\lambda)=\CDp_{D,e_4}(\lambda)=\lambda_1\lambda_2$. Indeed, the edges of the two crossings are $(o_1,u_1,o'_1,u'_1)=(e_1,e_2,e_4,e_1)$ and $(o_2,u_2,o'_2,u'_2)=(e_4,e_3,e_3,e_2)$, the flows of $D$ are the vectors with $z_{e_2}=z_{e_4}$, and on the flows which vanish at $e_2$ and $e_4$ the two forms are $\Psi_1(z,z')=z_{e_1}z'_{e_1}$ and $\Psi_2(z,z')=z_{e_3}z'_{e_3}$. So the determinant of the pencil depends on the cut edge, and so does its value $\det\bigl(\sum_j\Psi_j\bigr)$ at $\lambda=(1,\dots,1)$, which is $1$ or $2-t_0t_1$ for this diagram of the unknot. More generally, let $\mathsf w_0,\mathsf w_1,\mathsf w_2$ be elements of a field which contains $\BQ(t_0,t_1)$, and replace $\lambda_j^m$ by $\mathsf w_m$ in $\CDp_{D,\hate}(\lambda)$, for every $j$. For this diagram the result is $\mathsf w_1^2+(1-t_0t_1)\,\mathsf w_0\mathsf w_2$ for $\hate=e_1$ and $\mathsf w_1^2$ for $\hate=e_2$. If the result is to be the same for all the cut edges, then $\mathsf w_0\mathsf w_2=0$; and then, for every diagram, only the monomial $\lambda_1\cdots\lambda_c$ contributes, because a monomial of degree $c$ in $c$ variables which is not $\lambda_1\cdots\lambda_c$ contains a variable with the exponent $0$ and a variable with an exponent $\ge2$. So no specialization of this kind extracts from the pencil an invariant other than $\mathsf w_1^c\,\LG(L)$.
\end{example}

In the terms of Remark~\ref{rem.mixed}~(iii): in the determinant of the pencil every crossing is used zero, one or two times, the number of times being the exponent of its variable $\lambda_j$; the terms in which some crossing is not used exactly once vanish at the free fermion point $t_0t_1=1$; and the Links--Gould polynomial is the part in which every crossing is used once. A move of type I adds a term in which the new crossing is used twice, with the factor $g_0=1-(t_0t_1)^{\ve_0}$, and this term is not zero in general.

\part{State surfaces and deformed Seifert matrices}\label{part.surfaces}

In this part the lattice of flows is related to surfaces in space. Section~\ref{sec.statesurfaces} places the link and the state surfaces of all the smoothings of a diagram in space, and proves the localization theorem, on which the whole part rests. Section~\ref{sec.seifert} applies it to the canonical Seifert surface: the crossing forms are the Seifert form, band by band, and the edge polynomials are the determinants of deformed Seifert matrices. Section~\ref{sec.degree} draws the consequences for the degrees and for the extreme coefficients of the Links--Gould polynomial. Section~\ref{sec.cocores} turns to the coupling of the two flavors, which it reads off a local system on the complement of the link and of the co-cores of the bands. Sections~\ref{sec.degree} and~\ref{sec.cocores} are independent of each other.

\section{State surfaces and the localization of linking numbers}\label{sec.statesurfaces}

Every smoothing $\vs$ of a diagram $D$ defines a surface $F_\vs$ in $\BR^3$ whose boundary is the link: disks bounded by the circles of $\vs$, and a half-twisted band at each crossing. For the oriented smoothing this is the canonical Seifert surface of $D$, produced by Seifert's algorithm; the other smoothings give the \emph{state surfaces} of the diagram, which need not be orientable, and among them the two checkerboard surfaces. In this section we place all these surfaces on one side of the plane of the diagram, and we show that the lattice of flows computes their linking forms: the first homology group of $F_\vs$ is the quotient of $\FF_D$ by the circles of $\vs$ (Proposition~\ref{prop.statehomology}), and the linking number of a flow with a curve pushed off the surface is a sum of contributions of the crossings through whose bands the curve passes, each contribution being the flow through one of the two edges of the under-strand of the crossing (Theorem~\ref{thm.localization}). This gives closed formulas for the Gordon--Litherland form of every state surface (Corollary~\ref{cor.GLform}) and for the Seifert form of the canonical surface of every diagram (Corollary~\ref{cor.seifertform}). The proof rests on one observation: when the surface lies below the plane of the diagram, the faces of the diagram, which stay in this plane, are disks which span their boundaries and do not meet the surface elsewhere (Lemma~\ref{lem.spanning}).

In the next section these formulas identify the crossing forms $\pi_j$ of Section~\ref{sec.flows} with the Seifert form of the canonical surface, cut into the contributions of the individual bands. For the other smoothings they give the forms of Gordon and Litherland which are studied in Section~\ref{sec.minusone}, at the point where all the crossing variables are equal to $-1$; and the position in space of the link, of the co-cores and of the surfaces is used again in Section~\ref{sec.cocores}.

\subsection{The link, the co-cores and the state surfaces}\label{ssec.statesurface}
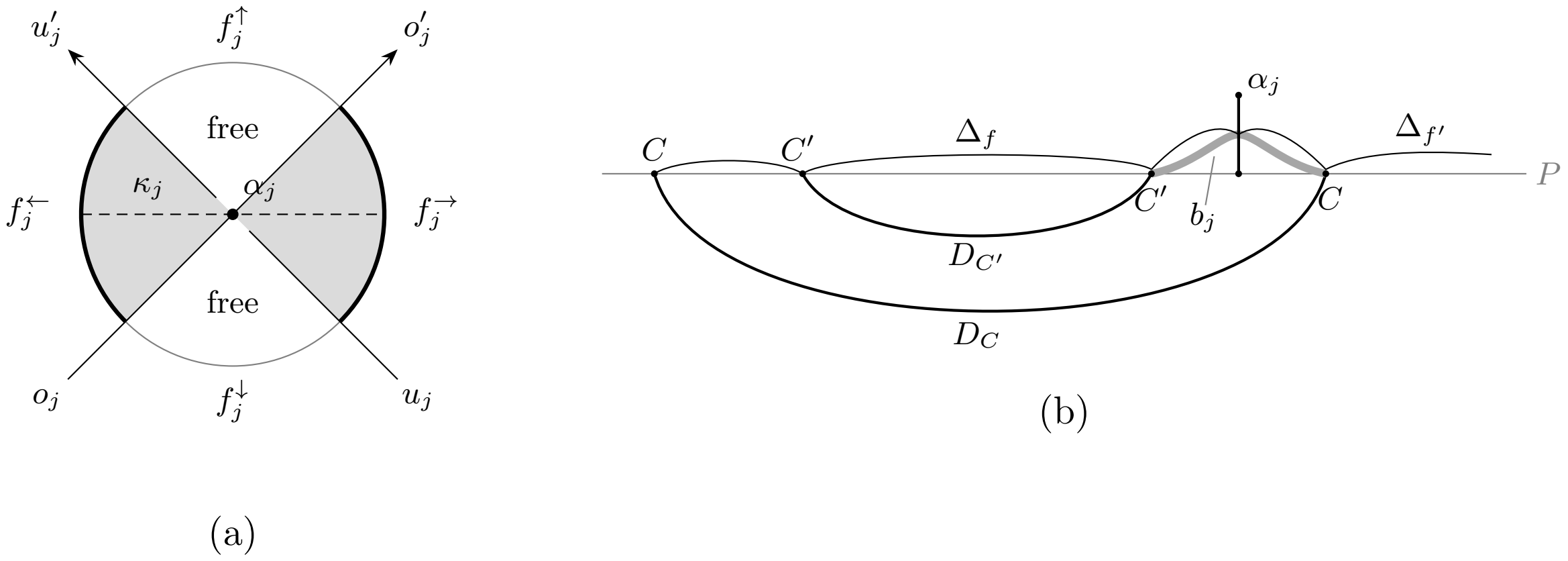
\begin{figure}[ht]
\centering
\begin{tikzpicture}[>=Stealth,every node/.style={font=\small}]
\begin{scope}[scale=1.15]
\fill[gray!28] (0,0) -- (135:1.3) arc (135:225:1.3) -- cycle;
\fill[gray!28] (0,0) -- (-45:1.3) arc (-45:45:1.3) -- cycle;
\draw[gray] (0,0) circle (1.3);
\draw[very thick] (135:1.3) arc (135:225:1.3);
\draw[very thick] (-45:1.3) arc (-45:45:1.3);
\draw (315:2.0) -- (315:0.2); \draw[->] (135:0.2) -- (135:2.0);
\draw[->] (225:2.0) -- (45:2.0);
\draw[densely dashed] (180:1.3) -- (0:1.3);
\fill (0,0) circle (1.4pt);
\node[font=\scriptsize] at (225:2.25) {$o_j$}; \node[font=\scriptsize] at (45:2.25) {$o'_j$};
\node[font=\scriptsize] at (315:2.25) {$u_j$}; \node[font=\scriptsize] at (135:2.25) {$u'_j$};
\node[font=\scriptsize] at (0,0.75) {free}; \node[font=\scriptsize] at (0,-0.75) {free};
\node[font=\scriptsize] at (-0.72,0.22) {$\kappa_j$};
\node[font=\scriptsize] at (0.24,0.2) {$\alpha_j$};
\node[font=\scriptsize] at (-1.75,0) {$f^\leftarrow_j$}; \node[font=\scriptsize] at (1.75,0) {$f^\rightarrow_j$};
\node[font=\scriptsize] at (0,1.6) {$f^\uparrow_j$}; \node[font=\scriptsize] at (0,-1.6) {$f^\downarrow_j$};
\node at (0,-2.75) {(a)};
\end{scope}
\begin{scope}[xshift=3.9cm,yshift=0.4cm,scale=0.86]
\draw[gray] (-0.3,0) -- (10.3,0); \node[gray,font=\scriptsize] at (10.55,0) {$P$};
\draw[thick] (0.3,0) .. controls (0.9,-2.1) and (7.4,-2.1) .. (8,0);
\draw[thick] (2,0) .. controls (2.5,-0.95) and (5.5,-0.95) .. (6,0);
\node[font=\scriptsize] at (4,-0.95) {$D_{C'}$}; \node[font=\scriptsize] at (4,-1.85) {$D_C$};
\draw[line width=2.2pt,gray!70] (6,0) .. controls (6.5,0.1) and (6.8,0.45) .. (7,0.45) .. controls (7.2,0.45) and (7.5,0.1) .. (8,0);
\draw[thick] (7,0) -- (7,0.9); \node[font=\scriptsize] at (7.3,1.02) {$\alpha_j$};
\fill (7,0) circle (1.1pt); \fill (7,0.9) circle (1.1pt);
\draw (2,0) .. controls (2.4,0.28) and (5.6,0.28) .. (6,0.05);
\draw (6,0.05) .. controls (6.3,0.36) and (6.7,0.64) .. (7,0.45);
\draw (0.3,0) .. controls (0.6,0.2) and (1.7,0.2) .. (2,0);
\draw (7,0.45) .. controls (7.3,0.64) and (7.7,0.36) .. (8,0.05) .. controls (8.5,0.3) and (9.4,0.25) .. (9.9,0.22);
\node[font=\scriptsize] at (4,0.45) {$\Delta_f$}; \node[font=\scriptsize] at (9.1,0.5) {$\Delta_{f'}$};
\fill (0.3,0) circle (1.1pt); \fill (2,0) circle (1.1pt); \fill (6,0) circle (1.1pt); \fill (8,0) circle (1.1pt);
\node[font=\scriptsize] at (0.3,0.28) {$C$}; \node[font=\scriptsize] at (1.95,0.3) {$C'$};
\node[font=\scriptsize] at (6.0,-0.28) {$C'$}; \node[font=\scriptsize] at (8.05,-0.28) {$C$};
\node[font=\scriptsize] at (6.6,-0.5) {$b_j$}; \draw[gray,thin] (6.62,-0.36) -- (6.72,0.2);
\node at (5,-2.75) {(b)};
\end{scope}
\end{tikzpicture}
\caption{The state surface of a smoothing. (a) A positive crossing seen from above, with the disk $\CB_j$, the two arcs of the oriented smoothing on $\partial \CB_j$ (thick), the two covered quadrants (gray), over which the band $b_j$ lies, the two free quadrants, the core $\kappa_j$ of the band (dashed) and the co-core $\alpha_j$, which is vertical and is seen as a point. (b) A vertical section along the core of the band of a crossing which joins a circle $C$ to a circle $C'$ inside it: the disks $D_C$ and $D_{C'}$ hang below the plane $P$, the inner one above the outer one; the band (gray), seen along its core, rises from the two arcs of the smoothing to the middle of the co-core; the face disks $\Delta_f$ and $\Delta_{f'}$ of Lemma~\ref{lem.spanning} (thin) lie above everything, and they contain the co-core in their boundary. At $C'$ the band continues the disk across the arc; at $C$ it is folded back over the disk.}
\label{fig.statesurface}
\end{figure}

In this subsection we fix, once and for all, the position in space of the link, of the co-cores of its crossings and of the state surfaces of all the smoothings of the diagram. The construction is explicit, with coordinates, because the proofs of this section compute linking numbers by counting how many times a curve passes through a disk, and for this one has to know which of two surfaces lies above the other one. The picture to keep in mind is that of Figure~\ref{fig.statesurface}: the link lies in the plane of the diagram, except for a small bridge at each crossing; the disks of the circles of a smoothing hang below this plane, like bowls which are stacked into one another; and the band of a crossing rises above the plane, from the under-strand to the over-strand.

Let $D$ be a diagram. We remove from $S^2$ a point which does not lie on $D$ and regard $D$ as a diagram in the horizontal plane $P=\BR^2\times\{0\}$ of $\BR^3$; the vertical direction is that of $e_z=(0,0,1)$, and \emph{above} and \emph{below} refer to it. The faces of $D$ are now the bounded components of $P\sm D$ and one unbounded component.

\emph{The link.} Fix $r_0>0$ so small that the closed disks $\CB_j\subset P$ of radius $r_0$ centered at the crossings are pairwise disjoint; after an isotopy of the plane we may assume that $D\cap \CB_j$ consists of two perpendicular diameters of $\CB_j$, the two strands of the crossing $j$. The four radii of $\CB_j$ on these diameters are the \emph{half-strands} of the crossing: they are the ends of the edges $o_j,u_j,o'_j,u'_j$, and we denote them by the same letters. Going around the crossing, half-strands of the over-strand alternate with half-strands of the under-strand. The four closed quarters of $\CB_j$ between consecutive half-strands are the \emph{quadrants} of the crossing; each of them is bounded by a half-strand of the over-strand and a half-strand of the under-strand, and it is a corner of one of the four faces $f^\uparrow_j,f^\leftarrow_j,f^\downarrow_j,f^\rightarrow_j$ of Section~\ref{ssec.flowlattice}. Fix a height $h_0>0$ and a continuous decreasing function $\chi\colon[0,r_0]\to[0,1]$ with $\chi(0)=1$ and $\chi(r_0)=0$. We realize the link $L$ as follows: outside the disks $\CB_j$ it is the diagram $D$ itself, in the plane $P$; inside $\CB_j$ the under-strand stays in $P$, and the point of the over-strand at the distance $r$ from the crossing is lifted to the height $h_0\chi(r)$. So $L$ lies in the plane $P$ except for $c$ small bridges, and its vertical projection is $D$.

\emph{The co-cores and the graph $G$.}\label{not.G} The \emph{co-core} of the crossing $j$ is the vertical segment $\alpha_j$ which joins the under-strand to the over-strand at the crossing, from the height $0$ to the height $h_0$; we orient it upward. Let
\[
   G=L\cup\alpha_1\cup\dots\cup\alpha_c\subset\BR^3 .
\]
It is a graph with $2c$ vertices, the end points of the co-cores, and $3c$ edges: the $2c$ edges of $D$, which are now arcs of $L$ from one end of a co-core to one end of a co-core, and the $c$ co-cores. Contracting the co-cores gives back the graph $D$. A flow $z\in\FF_D$ determines a $1$-cycle of $G$, which we denote by the same letter:
\begin{equation}\label{eq.flowcycle}
   z=\sum_{e\in E}z_e\,e+\sum_{j=1}^c\bigl(z_{u_j}-z_{u'_j}\bigr)\,\alpha_j .
\end{equation}
Indeed, at the lower end of $\alpha_j$ the edge $u_j$ arrives and the edges $u'_j$ and $\alpha_j$ leave, and at its upper end $o_j$ and $\alpha_j$ arrive and $o'_j$ leaves; so the chain \eqref{eq.flowcycle} is a cycle if and only if $z_{u_j}-z_{u'_j}=z_{o'_j}-z_{o_j}$ for every $j$, which is \eqref{eq.flow}. Conversely the coefficients of a $1$-cycle of $G$ on the edges of $L$ form a flow, which determines its coefficients on the co-cores. Thus
\begin{equation}\label{eq.flowsG}
   \FF_D=Z_1(G;\BZ)=H_1(G;\BZ):
\end{equation}
the lattice of flows is the group of the $1$-cycles of the graph $G$ in $\BR^3$. For instance the boundary $\partial f$ of a face is the closed curve which follows the edges of $L$ around $f$ and, at each corner of $f$, passes from one strand to the other along the co-core.

\emph{The circles of a smoothing.} Let $\vs$ be a smoothing of $D$ (Definition~\ref{def.smoothing}). We draw its arcs on the circles $\partial \CB_j$: the four ends of the half-strands cut $\partial \CB_j$ into four arcs, one in each quadrant, and the two arcs of the smoothing of the crossing $j$ are two opposite ones. If $\vs_j=0$ they are the arcs in the quadrants of $f^\leftarrow_j$ and $f^\rightarrow_j$, each of which lies between an incoming and an outgoing edge; if $\vs_j=1$ they are the arcs in the quadrants of $f^\uparrow_j$ and $f^\downarrow_j$. We call these two quadrants the \emph{covered quadrants} of the crossing $j$ and the other two its \emph{free quadrants}; covered and free quadrants alternate around the crossing (Figure~\ref{fig.statesurface}~(a)). A circle $C$ of $\vs$ is thus a simple closed curve in $P$ made of pieces of edges of $D$, outside the disks $\CB_j$, and of arcs of the circles $\partial \CB_j$. We write $\Omega_C\subset P$ for the closed disk which it bounds. Two circles of $\vs$ are disjoint; either the two disks which they bound are disjoint, or one of the circles lies in the interior of the disk bounded by the other one.

\emph{The disks.} For every circle $C$ of $\vs$ choose a continuous function $d_C\colon\Omega_C\to[0,\infty)$ which vanishes exactly on $C$, in such a way that
\begin{equation}\label{eq.stacking}
   d_C>d_{C'}\ \text{ on }\Omega_{C'}\qquad\text{whenever the circle $C'$ lies in the interior of $\Omega_C$.}
\end{equation}
(For instance $d_C(p)=2^{-n(C)}\min\bigl(\rho,\mathrm{dist}(p,C)\bigr)$, where $n(C)$ is the number of circles which enclose $C$ and $\rho$ is smaller than the distance between any two circles of $\vs$.) The disk of the circle $C$ is the graph of $-d_C$,
\[
   D_C=\bigl\{(p,-d_C(p)):\ p\in\Omega_C\bigr\} .
\]
It lies below the plane $P$ and meets it exactly along its boundary $C$. By \eqref{eq.stacking} the disks of two different circles are disjoint: \emph{the disk of an inner circle lies above the disk of every circle which encloses it} (Figure~\ref{fig.statesurface}~(b)). This is the convention which we use for all state surfaces, whatever the orientations of the circles; see Section~\ref{ssec.otherside} for the other choice.

\begin{figure}[ht]
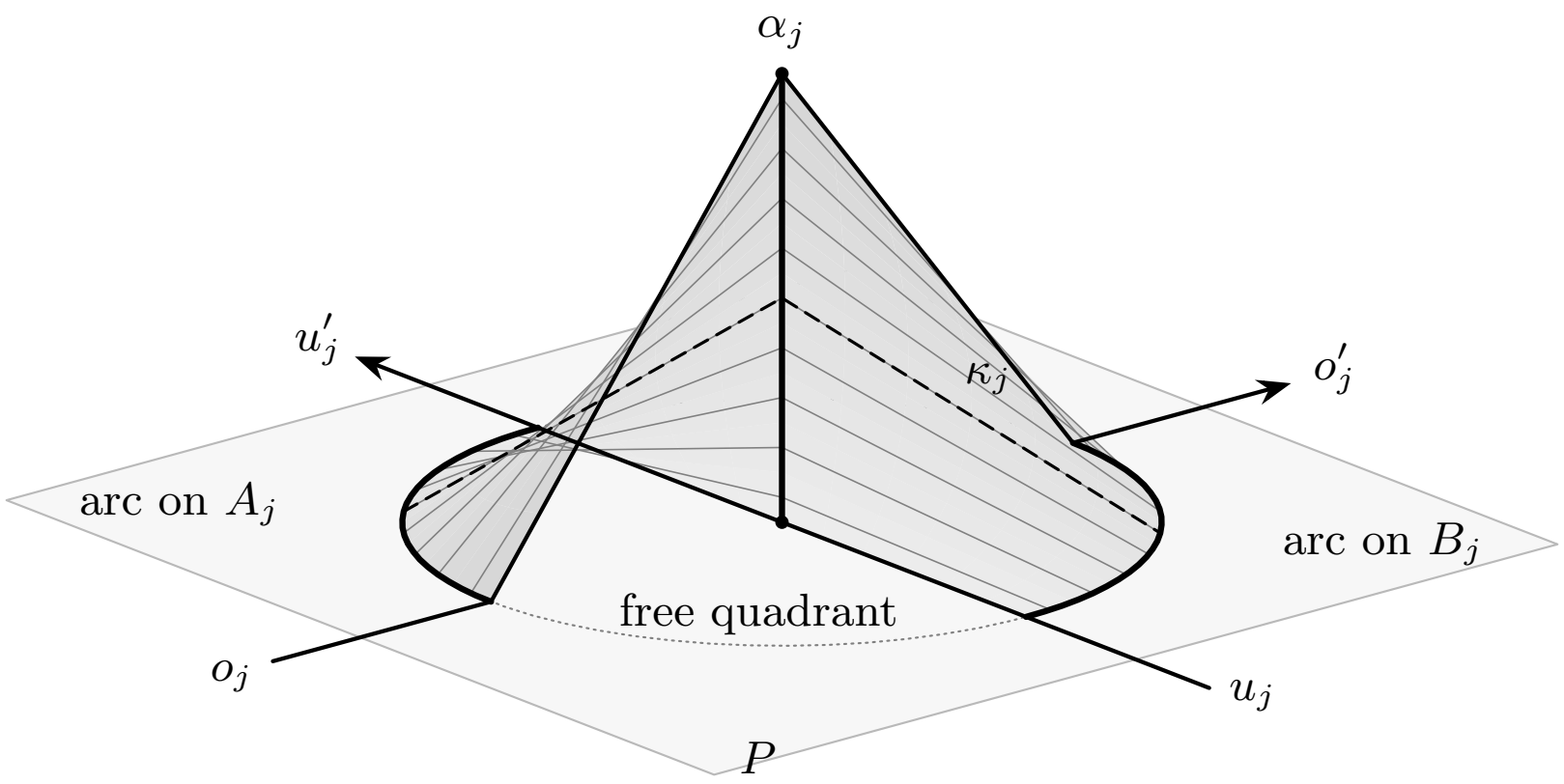

\centering

\caption{The band $b_j$ of a positive crossing $j$ for the oriented smoothing, seen from the side of the free quadrant of $f^\downarrow_j$ and drawn for $\chi(r)=1-r/r_0$. The under-strand lies in the plane $P$, the over-strand rises to the height $h_0$ above the crossing, and the co-core $\alpha_j$ joins them. Each half-band is swept by the segments which join a point of the co-core to a point of an arc of the smoothing (thick): as the point climbs the co-core, the segment turns by a quarter of a turn, from the direction of the under-strand to that of the over-strand. The two segments which start at the middle point of the co-core form the core $\kappa_j$ (dashed). One sees the same side of the band everywhere: the upper side of the half-band on the right, and the lower side of the half-band on the left.}
\label{fig.band3d}
\end{figure}

\emph{The bands.} Let $p$ be a point of a quadrant of $\CB_j$, different from the crossing, and let $r(p)$ be its distance from the crossing. We let $\tau(p)\in[-1,1]$ be the angular coordinate of $p$ in its quadrant, normalized so that $\tau=1$ on the half-strand of the over-strand, $\tau=-1$ on the half-strand of the under-strand, and $\tau$ is an affine function of the angle in between; and we put
\begin{equation}\label{eq.ruled}
   \beta(p)=\tfrac12\,h_0\,\chi\bigl(r(p)\bigr)\bigl(1+\tau(p)\bigr)\ \in[0,h_0] .
\end{equation}
The function $\beta$ interpolates, in each quadrant, between the heights of the two half-strands which bound it: $\beta=h_0\chi(r)$ on the half-strand of the over-strand and $\beta=0$ on that of the under-strand; and $\beta=0$ on the arc of $\partial \CB_j$. When $p$ tends to the crossing along the ray of angular coordinate $\tau$, $\beta(p)$ tends to $\frac12h_0(1+\tau)$, which runs through $[0,h_0]$: the closure of the graph of $\beta$ over a quadrant contains the co-core $\alpha_j$. The \emph{band} of the crossing $j$ is
\[
   b_j=\alpha_j\ \cup\ \bigl\{(p,\beta(p)):\ p\text{ in a covered quadrant of }\CB_j,\ p\ne\text{the crossing}\bigr\} .
\]
It is a disk. Its boundary consists of the two strands of $L$ over $\CB_j$ and of the two arcs of the smoothing of the crossing $j$; it lies above the plane $P$ and meets it exactly in the under-strand and in these two arcs. The co-core cuts $b_j$ into two \emph{half-bands}, one over each covered quadrant. Near the co-core the band is a piece of a helicoid: at the height $\frac12h_0(1+\tau)$ it leaves $\alpha_j$ horizontally in the two opposite directions of angular coordinate $\tau$ in the two covered quadrants, which turn by a quarter of a turn, from the direction of the under-strand to that of the over-strand, as the height increases from $0$ to $h_0$ (Figure~\ref{fig.band3d}). The \emph{core} of the band is the curve $\kappa_j\subset b_j$ over the bisectors of the two covered quadrants, at the height $\frac12h_0\chi(r)$; it joins the middle points of the two arcs of the smoothing and passes through the middle point of the co-core.

\begin{definition}[state surfaces]\label{def.statesurface}
The \emph{state surface} of the smoothing $\vs$ of $D$ is
\[
   F_\vs=\bigcup_{C}D_C\ \cup\ \bigcup_{j=1}^cb_j\ \subset\BR^3 ,
\]
the union of the disks of the circles of $\vs$ and of the bands of the crossings. For the oriented smoothing $\vs=0$ we write $F_D=F_0$ and call it the \emph{canonical Seifert surface} of the diagram $D$.
\end{definition}

\begin{lemma}\label{lem.statesurface}
The set $F_\vs$ is a compact connected surface embedded in $\BR^3$, with boundary $L$. The band $b_j$ meets the disks exactly along the two arcs of the smoothing of the crossing $j$, and the graph $G$ is contained in $F_\vs$. The surface $F_D=F_0$ is orientable.
\end{lemma}
\begin{proof}
The disks $D_C$ are pairwise disjoint and meet the plane $P$ along the circles of $\vs$. The band $b_j$ lies over $\CB_j$, above $P$, and meets $P$ in the under-strand and in the two arcs of the smoothing; the under-strand meets the circles of $\vs$ only at its two end points, which are end points of these arcs. So $b_j$ meets the union of the disks exactly along the two arcs, and two different bands are disjoint. Along an arc, a disk and a half-band are glued along a common boundary arc, so that $F_\vs$ is a surface near the points of the arc. The boundary of $F_\vs$ consists of the parts of the circles of $\vs$ which are not arcs of the smoothing, that is, of the diagram outside the disks $\CB_j$, and of the parts of the boundaries of the bands which are not such arcs, that is, of the strands of $L$ over the disks $\CB_j$: together this is $L$. The surface contains the co-cores, hence $G$, and it is connected because $D$ is. For $\vs=0$, orient every disk $D_C$ so that its boundary is the Seifert circle $C$ with the orientation of $L$. A band joins two disks along two arcs of Seifert circles which are traveled by $L$ in the same direction on both sides of the crossing; so the orientations of the two disks extend over the band, and $F_D$ is orientable, with oriented boundary $L$.
\end{proof}

So $F_D$ is the surface which Seifert's algorithm produces from $D$ \cite[Chapter~2]{Lickorish}: disks bounded by the Seifert circles, the disks of nested circles being stacked, joined by a half-twisted band at each crossing. The description of the algorithm leaves a choice in the stacking of nested disks, which matters as soon as a Seifert circle has crossings on both sides (Section~\ref{ssec.otherside}); here all the disks are pushed to the same side of the plane of the diagram. If one prefers to think of the disks as flat disks at different heights, the disk of each circle lying above the disks of the circles which enclose it, and of the bands as rising from the outer circle to the inner one, one obtains the same surface up to an isotopy which moves the points vertically.

\subsection{The homology of a state surface}\label{ssec.statehomology}
We now compute the first homology group of a state surface in terms of flows. The result is what makes the flows a common language for all the state surfaces: every flow is a $1$-cycle of the graph $G$, which lies on all of them; every homology class of $F_\vs$ is the class of a flow; and the flows which are homologous to zero on $F_\vs$ are the combinations of the circles of $\vs$.

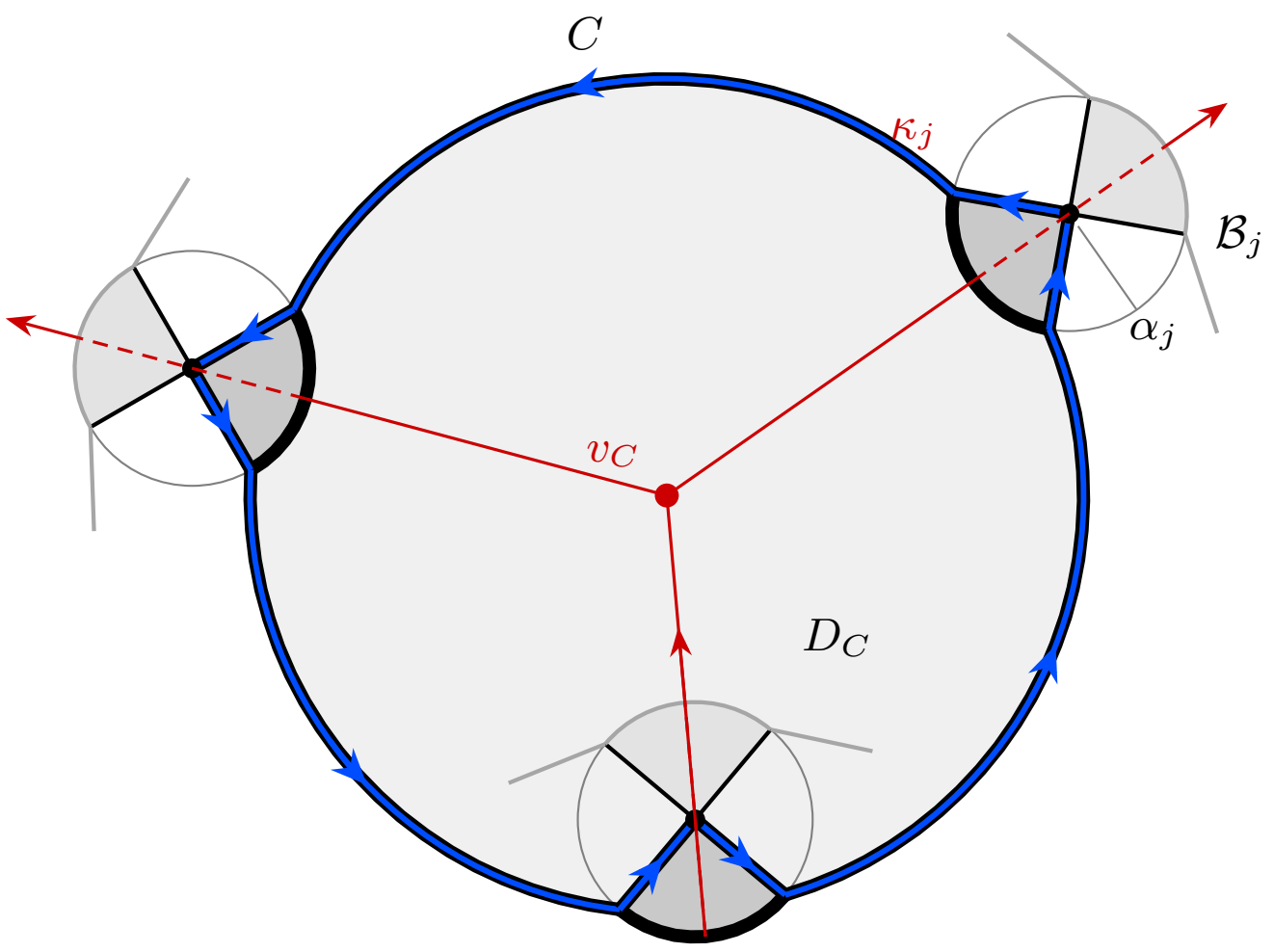
\begin{figure}[ht]
\centering
\begin{tikzpicture}[>=Stealth,scale=1.3,every node/.style={font=\small}]
\fill[gray!12] (2.017,0.877) -- (1.939,0.897) -- (1.863,0.926) -- (1.792,0.965) -- (1.727,1.013) -- (1.668,1.069) -- (1.617,1.132) -- (1.575,1.202) -- (1.543,1.276) -- (1.520,1.354) -- (1.507,1.434) -- (1.506,1.515) -- (1.515,1.596) -- (1.396,1.700) -- (1.270,1.796) -- (1.138,1.883) -- (1.000,1.959) -- (0.857,2.026) -- (0.709,2.082) -- (0.558,2.128) -- (0.404,2.163) -- (0.248,2.186) -- (0.090,2.198) -- (-0.068,2.199) -- (-0.226,2.188) -- (-0.382,2.167) -- (-0.537,2.134) -- (-0.688,2.090) -- (-0.836,2.035) -- (-0.980,1.969) -- (-1.119,1.894) -- (-1.252,1.809) -- (-1.379,1.714) -- (-1.498,1.611) -- (-1.610,1.499) -- (-1.714,1.380) -- (-1.808,1.253) -- (-1.893,1.120) -- (-1.969,0.981) -- (-1.933,0.909) -- (-1.907,0.832) -- (-1.891,0.752) -- (-1.886,0.671) -- (-1.891,0.591) -- (-1.907,0.511) -- (-1.933,0.434) -- (-1.969,0.361) -- (-2.014,0.294) -- (-2.067,0.233) -- (-2.128,0.180) -- (-2.196,0.135) -- (-2.200,-0.024) -- (-2.192,-0.183) -- (-2.173,-0.341) -- (-2.143,-0.497) -- (-2.101,-0.651) -- (-2.049,-0.801) -- (-1.986,-0.947) -- (-1.912,-1.088) -- (-1.828,-1.223) -- (-1.735,-1.352) -- (-1.633,-1.474) -- (-1.522,-1.588) -- (-1.404,-1.694) -- (-1.277,-1.791) -- (-1.145,-1.879) -- (-1.006,-1.957) -- (-0.862,-2.024) -- (-0.714,-2.081) -- (-0.561,-2.127) -- (-0.406,-2.162) -- (-0.249,-2.186) -- (-0.183,-2.234) -- (-0.112,-2.273) -- (-0.037,-2.302) -- (0.042,-2.322) -- (0.123,-2.330) -- (0.204,-2.329) -- (0.284,-2.316) -- (0.362,-2.294) -- (0.436,-2.261) -- (0.505,-2.219) -- (0.569,-2.168) -- (0.625,-2.109) -- (0.772,-2.060) -- (0.915,-2.001) -- (1.054,-1.931) -- (1.187,-1.852) -- (1.315,-1.764) -- (1.436,-1.667) -- (1.550,-1.561) -- (1.656,-1.448) -- (1.754,-1.328) -- (1.843,-1.201) -- (1.923,-1.068) -- (1.994,-0.930) -- (2.054,-0.787) -- (2.105,-0.640) -- (2.145,-0.490) -- (2.174,-0.338) -- (2.192,-0.184) -- (2.200,-0.029) -- (2.196,0.126) -- (2.182,0.281) -- (2.157,0.434) -- (2.121,0.585) -- (2.074,0.733) -- cycle;
\fill[gray!42] (2.125,1.488) -- (1.515,1.596) arc (170.00:260.00:0.620) -- cycle;
\fill[gray!22] (2.125,1.488) -- (2.736,1.380) arc (350.00:440.00:0.620) -- cycle;
\draw[gray,thin] (2.125,1.488) circle (0.620);
\fill[gray!42] (-2.506,0.671) -- (-2.196,0.135) arc (300.00:390.00:0.620) -- cycle;
\fill[gray!22] (-2.506,0.671) -- (-2.816,1.208) arc (480.00:570.00:0.620) -- cycle;
\draw[gray,thin] (-2.506,0.671) circle (0.620);
\fill[gray!42] (0.150,-1.711) -- (-0.249,-2.186) arc (230.00:320.00:0.620) -- cycle;
\fill[gray!22] (0.150,-1.711) -- (0.548,-1.236) arc (410.00:500.00:0.620) -- cycle;
\draw[gray,thin] (0.150,-1.711) circle (0.620);
\draw[line width=2.6pt] (2.017,0.877) -- (1.939,0.897) -- (1.863,0.926) -- (1.792,0.965) -- (1.727,1.013) -- (1.668,1.069) -- (1.617,1.132) -- (1.575,1.202) -- (1.543,1.276) -- (1.520,1.354) -- (1.507,1.434) -- (1.506,1.515) -- (1.515,1.596) -- (1.396,1.700) -- (1.270,1.796) -- (1.138,1.883) -- (1.000,1.959) -- (0.857,2.026) -- (0.709,2.082) -- (0.558,2.128) -- (0.404,2.163) -- (0.248,2.186) -- (0.090,2.198) -- (-0.068,2.199) -- (-0.226,2.188) -- (-0.382,2.167) -- (-0.537,2.134) -- (-0.688,2.090) -- (-0.836,2.035) -- (-0.980,1.969) -- (-1.119,1.894) -- (-1.252,1.809) -- (-1.379,1.714) -- (-1.498,1.611) -- (-1.610,1.499) -- (-1.714,1.380) -- (-1.808,1.253) -- (-1.893,1.120) -- (-1.969,0.981) -- (-1.933,0.909) -- (-1.907,0.832) -- (-1.891,0.752) -- (-1.886,0.671) -- (-1.891,0.591) -- (-1.907,0.511) -- (-1.933,0.434) -- (-1.969,0.361) -- (-2.014,0.294) -- (-2.067,0.233) -- (-2.128,0.180) -- (-2.196,0.135) -- (-2.200,-0.024) -- (-2.192,-0.183) -- (-2.173,-0.341) -- (-2.143,-0.497) -- (-2.101,-0.651) -- (-2.049,-0.801) -- (-1.986,-0.947) -- (-1.912,-1.088) -- (-1.828,-1.223) -- (-1.735,-1.352) -- (-1.633,-1.474) -- (-1.522,-1.588) -- (-1.404,-1.694) -- (-1.277,-1.791) -- (-1.145,-1.879) -- (-1.006,-1.957) -- (-0.862,-2.024) -- (-0.714,-2.081) -- (-0.561,-2.127) -- (-0.406,-2.162) -- (-0.249,-2.186) -- (-0.183,-2.234) -- (-0.112,-2.273) -- (-0.037,-2.302) -- (0.042,-2.322) -- (0.123,-2.330) -- (0.204,-2.329) -- (0.284,-2.316) -- (0.362,-2.294) -- (0.436,-2.261) -- (0.505,-2.219) -- (0.569,-2.168) -- (0.625,-2.109) -- (0.772,-2.060) -- (0.915,-2.001) -- (1.054,-1.931) -- (1.187,-1.852) -- (1.315,-1.764) -- (1.436,-1.667) -- (1.550,-1.561) -- (1.656,-1.448) -- (1.754,-1.328) -- (1.843,-1.201) -- (1.923,-1.068) -- (1.994,-0.930) -- (2.054,-0.787) -- (2.105,-0.640) -- (2.145,-0.490) -- (2.174,-0.338) -- (2.192,-0.184) -- (2.200,-0.029) -- (2.196,0.126) -- (2.182,0.281) -- (2.157,0.434) -- (2.121,0.585) -- (2.074,0.733) -- cycle;
\draw[thick] (1.515,1.596) -- (2.736,1.380);
\draw[thick] (2.017,0.877) -- (2.233,2.099);
\draw[thick,gray!70] (2.736,1.380) arc (350.00:440.00:0.620);
\draw[thick,gray!70] (2.736,1.380) -- (2.906,0.857);
\draw[thick,gray!70] (2.233,2.099) -- (1.799,2.437);
\draw[thick] (-2.196,0.135) -- (-2.816,1.208);
\draw[thick] (-1.969,0.981) -- (-3.043,0.361);
\draw[thick,gray!70] (-2.816,1.208) arc (480.00:570.00:0.620);
\draw[thick,gray!70] (-2.816,1.208) -- (-2.524,1.675);
\draw[thick,gray!70] (-3.043,0.361) -- (-3.024,-0.188);
\draw[thick] (-0.249,-2.186) -- (0.548,-1.236);
\draw[thick] (0.625,-2.109) -- (-0.325,-1.312);
\draw[thick,gray!70] (0.548,-1.236) arc (410.00:500.00:0.620);
\draw[thick,gray!70] (0.548,-1.236) -- (1.086,-1.350);
\draw[thick,gray!70] (-0.325,-1.312) -- (-0.835,-1.518);
\draw[line width=2.6pt] (2.017,0.877) -- (2.125,1.488) -- (1.515,1.596);
\draw[blue!70!cyan,line width=1.1pt] (2.017,0.877) -- (2.125,1.488) -- (1.515,1.596);
\draw[blue!70!cyan,line width=1.1pt,-{Stealth[length=7pt,width=6pt]}] (2.017,0.877) -- (2.080,1.232);
\draw[blue!70!cyan,line width=1.1pt,-{Stealth[length=7pt,width=6pt]}] (2.125,1.488) -- (1.710,1.561);
\draw[blue!70!cyan,line width=1.1pt] (1.515,1.596) arc (46.49:153.51:2.200);
\draw[blue!70!cyan,line width=1.1pt,-{Stealth[length=7pt,width=6pt]}] (-0.230,2.188) arc (96.00:104.00:2.200);
\draw[line width=2.6pt] (-1.969,0.981) -- (-2.506,0.671) -- (-2.196,0.135);
\draw[blue!70!cyan,line width=1.1pt] (-1.969,0.981) -- (-2.506,0.671) -- (-2.196,0.135);
\draw[blue!70!cyan,line width=1.1pt,-{Stealth[length=7pt,width=6pt]}] (-1.969,0.981) -- (-2.280,0.802);
\draw[blue!70!cyan,line width=1.1pt,-{Stealth[length=7pt,width=6pt]}] (-2.506,0.671) -- (-2.295,0.306);
\draw[blue!70!cyan,line width=1.1pt] (-2.196,0.135) arc (176.49:263.51:2.200);
\draw[blue!70!cyan,line width=1.1pt,-{Stealth[length=7pt,width=6pt]}] (-1.780,-1.293) arc (216.00:224.00:2.200);
\draw[line width=2.6pt] (-0.249,-2.186) -- (0.150,-1.711) -- (0.625,-2.109);
\draw[blue!70!cyan,line width=1.1pt] (-0.249,-2.186) -- (0.150,-1.711) -- (0.625,-2.109);
\draw[blue!70!cyan,line width=1.1pt,-{Stealth[length=7pt,width=6pt]}] (-0.249,-2.186) -- (-0.018,-1.910);
\draw[blue!70!cyan,line width=1.1pt,-{Stealth[length=7pt,width=6pt]}] (0.150,-1.711) -- (0.473,-1.982);
\draw[blue!70!cyan,line width=1.1pt] (0.625,-2.109) arc (286.49:383.51:2.200);
\draw[blue!70!cyan,line width=1.1pt,-{Stealth[length=7pt,width=6pt]}] (1.924,-1.067) arc (331.00:339.00:2.200);
\fill (2.125,1.488) circle (1.5pt);
\fill (-2.506,0.671) circle (1.5pt);
\fill (0.150,-1.711) circle (1.5pt);
\fill[red!80!black] (0,0) circle (1.8pt);
\draw[red!80!black,semithick] (0,0) -- (1.617,1.132);
\draw[red!80!black,semithick,densely dashed] (1.617,1.132) -- (2.633,1.844);
\draw[red!80!black,semithick,->] (2.633,1.844) -- (2.961,2.073);
\draw[red!80!black,semithick] (0,0) -- (-1.907,0.511);
\draw[red!80!black,semithick,densely dashed] (-1.907,0.511) -- (-3.105,0.832);
\draw[red!80!black,semithick,->] (-3.105,0.832) -- (-3.491,0.935);
\draw[red!80!black,semithick] (0,0) -- (0.204,-2.329);
\draw[red!80!black,semithick,densely dashed] (0.204,-2.329) -- (0.096,-1.093);
\draw[red!80!black,semithick,->] (0.096,-1.093) -- (0.061,-0.695);
\node[font=\scriptsize,red!80!black] at (-0.28,0.22) {$v_C$};
\node[font=\scriptsize] at (0.9,-0.75) {$D_C$};
\node[font=\scriptsize] at (-0.431,2.442) {$C$};
\node[font=\scriptsize] at (2.573,0.849) {$\alpha_j$};
\draw[gray,thin] (2.481,0.980) -- (2.171,1.422);
\node[font=\scriptsize,red!80!black] at (1.305,1.906) {$\kappa_j$};
\node[font=\scriptsize] at (3.030,1.364) {$\CB_j$};
\end{tikzpicture}
\caption{The $2$-cell of a circle $C$ of a smoothing, seen from above. The circle (thick) passes through three crossings; two of them lie outside $C$ and one lies inside. The cell consists of the disk $D_C$ (light gray) and of the three half-bands which are attached to it along the arcs of the smoothing on $C$ (dark gray); for the crossing inside $C$ the half-band is folded back over the disk. The boundary of the cell (blue) is the cycle $z^C$ of the graph $G$: it follows $C$, except that at each arc of the smoothing it runs along a half-strand to the crossing, along the co-core, which is seen as a point, and back along the other half-strand. In red: the vertex $v_C$ of the state graph $\Gamma_\vs$, the paths which join it to the middle points of the three arcs, and the cores of the bands (dashed), which lead to the other half-bands (lighter) and to the disks of the circles to which these are attached.}
\label{fig.cell}
\end{figure}

Removing the graph $G$ from $F_\vs$ leaves, for every circle $C$ of $\vs$, the interior of the disk $D_C$ together with the open half-bands which are attached to $D_C$ along the arcs of $C$: an open disk (Figure~\ref{fig.cell}). So $F_\vs$ is a cell complex with $1$-skeleton $G$ and one $2$-cell for each circle of $\vs$. The boundary of the $2$-cell of $C$ follows $C$, except that at every arc of the smoothing on $C$ it runs around the half-band which is attached there: along one half-strand to the crossing, along the co-core, and back along the other half-strand. It is the $1$-cycle of $G$ which is traversed by the circle $C$, that is, the flow $z^C$ of Section~\ref{ssec.smoothings}, up to sign.

\begin{proposition}[the homology of a state surface]\label{prop.statehomology}
For every smoothing $\vs$ of $D$, the inclusion of $G$ into $F_\vs$ induces an isomorphism
\[
   H_1(F_\vs;\BZ)\cong\FF_D\big/\langle z^C:\ C\text{ a circle of }\vs\rangle ,
\]
a free abelian group of rank $c+1-s_\vs$. Under this isomorphism the transfer $w^\vs_j(z)$ of a flow $z$ (Definition~\ref{def.transfer}) is the algebraic number of times the cycle $z$ runs through the band $b_j$, from the half-band which contains the half-strand $u'_j$ to the half-band which contains the half-strand $u_j$.
\end{proposition}
\begin{proof}
The first homology group of a cell complex is the group of the $1$-cycles of its $1$-skeleton modulo the boundaries of its $2$-cells; so the isomorphism follows from \eqref{eq.flowsG} and from the description of the $2$-cells, and the rank from Lemma~\ref{lem.circles}. For the last statement let $\alpha'_j\subset b_j$ be a parallel copy of the co-core inside the half-band of $u'_j$: an arc which joins a point $q$ of the half-strand $u'_j$ to a point $q'$ of the other half-strand which bounds this half-band, which is $o_j$ if $\vs_j=0$ and $o'_j$ if $\vs_j=1$. A closed curve on $F_\vs$ which runs through the band from the half-band of $u'_j$ to that of $u_j$ crosses $\alpha'_j$ once, toward the co-core, and the algebraic number of such crossings depends only on the homology class of the curve. A cycle $z$ of $G$ meets $\alpha'_j$ only at $q$ and at $q'$. It passes $z_{u'_j}$ times through $q$ in the direction of the edge $u'_j$, which leads away from the crossing; and it passes through $q'$ in the direction of the edge, $z_{o_j}$ times toward the crossing if $\vs_j=0$ and $z_{o'_j}$ times away from it if $\vs_j=1$. Pushing $z$ slightly into the interior of the surface, we see that it crosses $\alpha'_j$ toward the co-core $z_{o_j}-z_{u'_j}$ times if $\vs_j=0$ and $-z_{o'_j}-z_{u'_j}$ times if $\vs_j=1$, which is $w^\vs_j(z)$ in both cases.
\end{proof}

For the canonical surface this is the familiar description of $H_1(F_D)$: by Lemma~\ref{lem.lifts}~(i) the transfers identify $\FF_D/\langle C_1,\dots,C_s\rangle$ with the cycle space $H$ of the Seifert graph, the surface $F_D$ retracts onto a copy of $\Gamma_D$ whose vertices are the disks and whose edges are the cores of the bands, and $\nu=c-s+1$ is the rank of $H_1(F_D)$, which is $2g(F_D)+\mu-1$.

More generally, let $\Gamma_\vs$ be the \emph{state graph}\label{not.stategraph} of the smoothing: it has one vertex for each circle of $\vs$ and one edge for each crossing $j$, which joins the two circles (possibly the same one) on which the two arcs of the smoothing of $j$ lie. We embed it in $F_\vs$: the vertex of $C$ is a point in the interior of $D_C$, and the edge $j$ is the core $\kappa_j$ of the band, extended at each of its two ends by a path in the disk on which it arrives, from the middle point of the arc to the vertex of that disk; the interiors of these paths lie in the interiors of the disks (Figure~\ref{fig.cell}). The surface $F_\vs$ retracts by deformation onto $\Gamma_\vs$, each disk with its half-bands retracting onto the star formed by the paths which it contains. We say that a closed curve on $F_\vs$ is in \emph{standard position} if it is a closed walk in $\Gamma_\vs$. Every class in $H_1(F_\vs;\BZ)$ is a sum of classes of such curves, and a curve $b$ in standard position runs through the cores of the bands; we write $w^\vs_j(b)$ for the number of times it runs through $\kappa_j$ from the half-band of $u'_j$ to the half-band of $u_j$, minus the number of times it runs through it in the other direction. By Proposition~\ref{prop.statehomology} these are the transfers of any flow in the class of $b$.

\subsection{The faces span their boundaries}\label{ssec.spanning}
This subsection contains the observation on which the whole section rests. To compute the linking number of a flow $z$ with a closed curve one needs a surface, or a $2$-chain, whose boundary is $z$. The faces of the diagram provide one: every flow is a combination of the boundaries $\partial f$ of the faces (Lemma~\ref{lem.faces}), and the face $f$ itself, lifted slightly above the plane of the diagram, is a disk $\Delta_f$ with boundary $\partial f$. All the state surfaces lie below these disks. So the linking number of a flow with a curve which stays close to a state surface is obtained by counting the passages of the curve through the disks $\Delta_f$, and this can be done on the diagram.

We use the following conventions for linking numbers. For two disjoint oriented closed curves, or more generally two disjoint $1$-cycles, $a$ and $c$ in $\BR^3$, the linking number $\lk(a,c)$ is symmetric and bilinear, and it can be computed from an oriented surface, or a $2$-chain, $\Delta$ with boundary $a$: it is the intersection number $\Delta\cdot c$, the number of times $c$ passes through $\Delta$ in the direction of its positive normal minus the number of times it passes in the opposite direction. The positive normal $n$ of an oriented surface is such that $(v_1,v_2,n)$ is a positive basis of $\BR^3$ for every positive basis $(v_1,v_2)$ of a tangent plane, and the boundary is oriented by the rule ``outward normal first''. For instance, the boundary of a horizontal disk is traveled counterclockwise, as seen from above, if the positive normal of the disk is $e_z$, and it has the linking number $+1$ with a vertical line which is traveled upward through the disk. With this convention, $\lk(a,c)$ is the sum of the signs of the crossings at which $c$ passes over $a$ in a projection, and two curves which cross positively twice have linking number $+1$.

Let $f$ be a bounded face of $D$. Choose a continuous function $\psi_f\ge0$ on the closure of $f$ in $P$ which vanishes exactly on the boundary of $f$, and define a function $u_f$ on $f$ by
\begin{equation}\label{eq.uf}
   u_f(p)=\beta(p)+\psi_f(p)\ \text{ if $p$ lies in a disk $\CB_j$},\qquad u_f(p)=\psi_f(p)\ \text{ otherwise,}
\end{equation}
with $\beta$ as in \eqref{eq.ruled}; it is continuous because $\beta=0$ on the circles $\partial \CB_j$. Let $\Delta_f$ be the closure of the graph of $u_f$ in $\BR^3$, oriented by the upward normal. So $\Delta_f$ is the face $f$, lifted slightly above the plane, and at each corner of $f$ it rises along the quadrant from the half-strand of the under-strand to that of the over-strand, as the bands do (Figure~\ref{fig.statesurface}~(b)).

\begin{lemma}[the faces span their boundaries]\label{lem.spanning}
Let $f$ be a bounded face of $D$ and $\vs$ a smoothing of $D$.
\begin{enumerate}\renewcommand{\theenumi}{\roman{enumi}}
\item The boundary of the $2$-chain $\Delta_f$ is the cycle $\partial f$ of $G$.
\item Over the face $f$, the surface $F_\vs$ lies strictly below $\Delta_f$: if $p\in f$ and $(p,h)\in F_\vs$, then $h<u_f(p)$. In particular $\Delta_f$ meets $F_\vs$ only along the graph $G$.
\item Let $c$ be an oriented closed curve in $\BR^3$ which is disjoint from $G$ and which crosses the vertical walls over the boundary of $f$ a finite number of times. Then $\lk(\partial f,c)$ is the number of times $c$, while it travels over $f$, passes from below $\Delta_f$ to above it, minus the number of times it passes from above to below. Consequently, if a maximal arc of $c$ over $f$ is below $\Delta_f$ at both of its ends, or above it at both ends, it contributes $0$ to $\lk(\partial f,c)$; it contributes $+1$ if it enters below and leaves above, and $-1$ in the opposite case.
\end{enumerate}
\end{lemma}
\begin{proof}
(i) The graph of $u_f$ over the open disk $f$ is an open disk, and we determine what its closure adds over the boundary of $f$. Over a point $p$ of an edge of $D$ outside the disks $\CB_j$, $u_f$ tends to $\psi_f(p)=0$: the closure contains $p$, which is a point of $L$. Over a point $p\ne$ the crossing of a half-strand in $\CB_j$, $\psi_f$ tends to $0$ and $\tau$ tends to $1$ or to $-1$, so $u_f$ tends to $h_0\chi(r)$ on a half-strand of the over-strand and to $0$ on a half-strand of the under-strand: the closure contains the point of $L$ above $p$. Over the crossing $j$, at a corner of $f$, $u_f$ tends to $\frac12h_0(1+\tau)$ along the ray of angular coordinate $\tau$, since $\psi_f$ tends to $0$: the closure contains the co-core $\alpha_j$. So the boundary of $\Delta_f$ is a $1$-cycle of $G$. On the edges of $L$ it is the boundary of the face $f$, traveled counterclockwise because $\Delta_f$ is oriented by the upward normal: each boundary edge is traveled in its direction if $f$ lies on its left and against it otherwise, which is $\partial f$ by \eqref{eq.boundaryface}. A cycle of $G$ is determined by its coefficients on the edges of $L$, by \eqref{eq.flowsG}.

(ii) Let $p\in f$. The points of the disks $D_C$ over $p$ have heights $-d_C(p)\le0$, and $u_f(p)\ge\psi_f(p)>0$. The band $b_j$ lies over the two covered quadrants of $\CB_j$ and over the crossing itself; if $p$ lies in a covered quadrant, the point of $b_j$ over $p$ has the height $\beta(p)<\beta(p)+\psi_f(p)=u_f(p)$. There are no other points of $F_\vs$ over $p$.

(iii) By (i), $\lk(\partial f,c)=\Delta_f\cdot c$. The interior of $\Delta_f$ is a graph over $f$, oriented by the upward normal, so $\Delta_f\cdot c$ is the number of times $c$ crosses this graph upward minus the number of times it crosses it downward. A maximal arc of $c$ over $f$ begins and ends at points which lie over the boundary of $f$, or it is all of $c$. Since $c$ is disjoint from $G$, each of its ends lies above or below the point of $L$, or the co-core, which $\Delta_f$ contains over the same point of the boundary of $f$; so the arc is above or below the graph near each of its ends, and its algebraic number of crossings with the graph is $1$, $-1$ or $0$ as stated. A closed curve which lies entirely over $f$ crosses the graph as often upward as downward.
\end{proof}

Since every flow is a combination of boundaries of faces with a potential which vanishes on the unbounded face (Lemma~\ref{lem.faces}~(ii)), part (iii) computes the linking number of every flow with every closed curve $c$ in the complement of $G$:
\begin{equation}\label{eq.lkpotential}
   \lk(z,c)=\sum_{f\text{ bounded}}k_f\ \Delta_f\cdot c\qquad\text{if }z=\sum_fk_f\,\partial f\text{ with }k=0\text{ on the unbounded face.}
\end{equation}

\subsection{Pushing a band off the surface}\label{ssec.pushoff}
The forms of Gordon and Litherland and of Seifert are defined by pushing curves off the surface. In this subsection we describe the push-offs of a curve which runs through the band of a crossing, and we determine how they pass through the face disks of the four faces around the crossing (Lemma~\ref{lem.passage}). This local computation is the heart of the proof of the localization theorem.

Let $b$ be a closed curve on a surface $F\subset\BR^3$, in the interior of $F$. A small tubular neighborhood of $F$ along $b$ is cut by $F$ into two sides, locally; pushing $b$ slightly off $F$ to both sides gives a $1$-cycle $\tau b$ in the complement of $F$, which runs twice along $b$ and consists of two closed curves parallel to $b$ if the sides do not get exchanged when $b$ is traveled once, and of one closed curve otherwise \cite{GordonLitherland}. If $F$ is oriented, the two curves are the push-offs $b^+$ and $b^-$ of $b$ to the positive and to the negative side of $F$. All that matters for the linking numbers with cycles lying on $F$ is the class of $\tau b$, or of $b^{\pm}$, up to deformations in the complement of $F$, and we now describe convenient representatives for a curve $b$ in standard position on a state surface $F_\vs$. Fix a small number $\varrho>0$ and a radius $r_1<r_0$.

\emph{Along a path in a disk.} Where $b$ follows a path in a disk $D_C$, the two push-offs are the copies of this path which are displaced vertically by $\varrho$, above and below the disk.

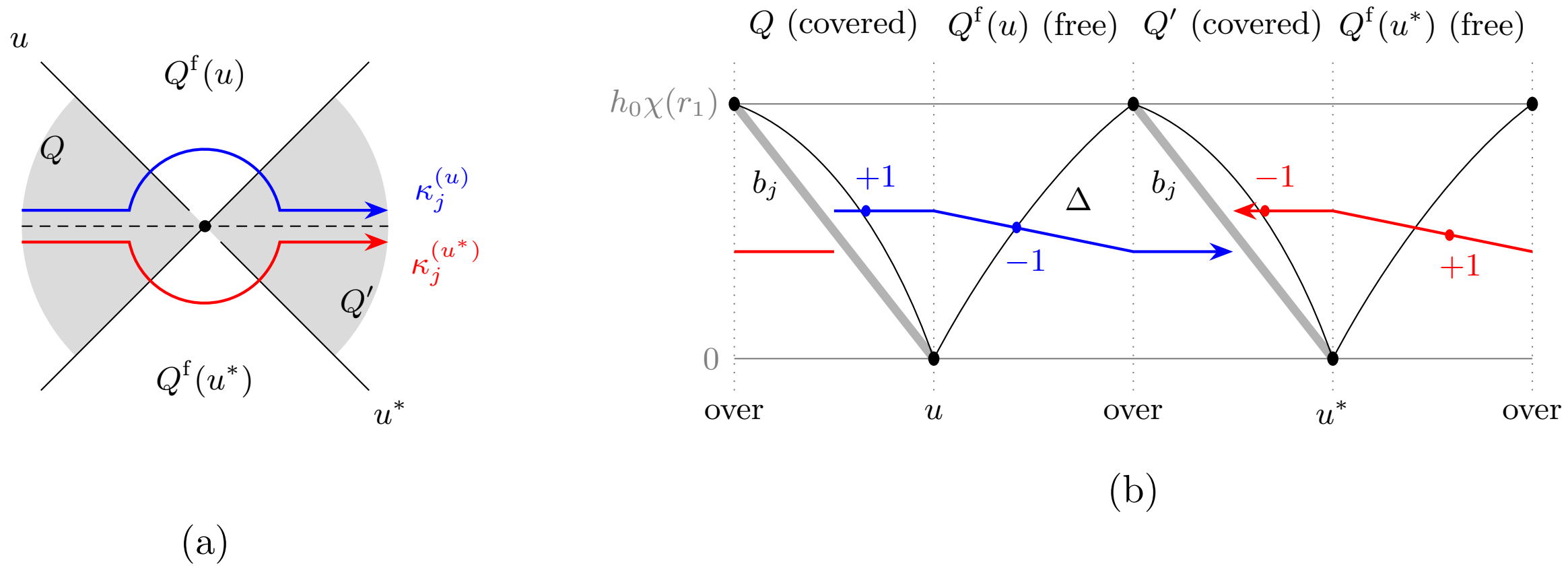
\begin{figure}[ht]
\centering
\begin{tikzpicture}[>=Stealth,every node/.style={font=\small}]
\begin{scope}[scale=1.2]
\fill[gray!28] (0,0) -- (135:1.5) arc (135:225:1.5) -- cycle;
\fill[gray!28] (0,0) -- (-45:1.5) arc (-45:45:1.5) -- cycle;
\draw (315:1.9) -- (315:0.18); \draw (135:0.18) -- (135:1.9);
\draw (225:1.9) -- (45:1.9);
\draw[densely dashed] (180:1.5) -- (0:1.5);
\fill (0,0) circle (1.4pt);
\draw[blue,thick,->] (-1.5,0.13) -- (-0.62,0.13) arc (168:12:0.63) -- (1.5,0.13);
\draw[red,thick,->] (-1.5,-0.13) -- (-0.62,-0.13) arc (192:348:0.63) -- (1.5,-0.13);
\node[font=\scriptsize] at (-1.25,0.62) {$Q$}; \node[font=\scriptsize] at (1.25,-0.62) {$Q'$};
\node[font=\scriptsize] at (135:2.15) {$u$}; \node[font=\scriptsize] at (315:2.15) {$u^*$};
\node[font=\scriptsize] at (0,1.25) {$Q^{\rm f}(u)$}; \node[font=\scriptsize] at (0,-1.25) {$Q^{\rm f}(u^*)$};
\node[blue,font=\scriptsize] at (1.95,0.3) {$\kappa_j^{(u)}$}; \node[red,font=\scriptsize] at (1.98,-0.3) {$\kappa_j^{(u^*)}$};
\node at (0,-2.6) {(a)};
\end{scope}
\begin{scope}[xshift=5.2cm,yshift=-1.3cm,xscale=0.98,yscale=1.25]
\foreach \x in {0,2,4,6,8}{ \draw[gray,thin,dotted] (\x,-0.25) -- (\x,2.35); }
\draw[gray,thin] (0,0) -- (8,0); \draw[gray,thin] (0,2) -- (8,2);
\node[gray,font=\scriptsize,left] at (0,0) {$0$}; \node[gray,font=\scriptsize,left] at (0,2) {$h_0\chi(r_1)$};
\draw[line width=2.4pt,gray!60] (0,2) -- (2,0); \draw[line width=2.4pt,gray!60] (4,2) -- (6,0);
\draw (0,2) .. controls (0.9,1.75) and (1.6,0.9) .. (2,0);
\draw (2,0) .. controls (2.5,0.75) and (3.3,1.6) .. (4,2);
\draw (4,2) .. controls (4.9,1.75) and (5.6,0.9) .. (6,0);
\draw (6,0) .. controls (6.5,0.75) and (7.3,1.6) .. (8,2);
\foreach \x in {0,4,8}{ \fill (\x,2) circle (1.6pt); }
\foreach \x in {2,6}{ \fill (\x,0) circle (1.6pt); }
\node[font=\scriptsize] at (2,-0.42) {$u$}; \node[font=\scriptsize] at (6,-0.42) {$u^*$};
\node[font=\scriptsize] at (0,-0.42) {over}; \node[font=\scriptsize] at (4,-0.42) {over}; \node[font=\scriptsize] at (8,-0.42) {over};
\draw[blue,thick,->] (1,1.16) -- (2,1.16) -- (4,0.84) -- (5,0.84);
\draw[red,thick] (1,0.84) -- (0,0.84); \draw[red,thick,->] (8,0.84) -- (6,1.16) -- (5,1.16);
\fill[blue] (1.32,1.16) circle (1.3pt); \fill[blue] (2.83,1.03) circle (1.3pt);
\fill[red] (7.17,0.97) circle (1.3pt); \fill[red] (5.32,1.16) circle (1.3pt);
\node[blue,font=\scriptsize] at (1.42,1.42) {$+1$}; \node[blue,font=\scriptsize] at (2.92,0.76) {$-1$};
\node[red,font=\scriptsize] at (7.28,0.72) {$+1$}; \node[red,font=\scriptsize] at (5.42,1.42) {$-1$};
\node[font=\scriptsize] at (1,2.62) {$Q$ (covered)}; \node[font=\scriptsize] at (3,2.62) {$Q^{\rm f}(u)$ (free)};
\node[font=\scriptsize] at (5,2.62) {$Q'$ (covered)}; \node[font=\scriptsize] at (7,2.62) {$Q^{\rm f}(u^*)$ (free)};
\node[font=\scriptsize] at (0.32,1.36) {$b_j$}; \node[font=\scriptsize] at (4.32,1.36) {$b_j$};
\node[font=\scriptsize] at (3.45,1.25) {$\Delta$};
\node at (4,-1.05) {(b)};
\end{scope}
\end{tikzpicture}
\caption{The two push-offs of a passage along the core of the band $b_j$, from the covered quadrant $Q$ to the covered quadrant $Q'$. (a) Seen from above: $\kappa_j^{(u)}$ (blue) turns around the co-core through the free quadrant next to the half-strand $u$ of the under-strand, and $\kappa_j^{(u^*)}$ (red) through the other free quadrant. (b) The vertical cylinder of radius $r_1$ around the co-core, unrolled: the horizontal coordinate is the angle around the crossing and the vertical coordinate is the height. The dots are the four half-strands, those of the under-strand at the height $0$ and those of the over-strand at the height $h_0\chi(r_1)$; the thick gray lines are the two half-bands and the thin curves are the face disks $\Delta_f$ of the four corners, which coincide with the half-bands on the half-strands and lie above them elsewhere. The curve $\kappa_j^{(u)}$ starts above the half-band over $Q$ and below the face disk, passes above the face disk of $Q$ before it crosses over $u$, passes below the face disk of $Q^{\rm f}(u)$ before it crosses under the over-strand, and ends below the half-band over $Q'$; $\kappa_j^{(u^*)}$ goes the other way around (the left and right edges of the picture are identified). The signs are those of Lemma~\ref{lem.passage}.}
\label{fig.unrolled}
\end{figure}

\emph{Along the core of a band.} Let $Q$ and $Q'$ be the two covered quadrants of the crossing $j$, let $u$ be the half-strand of the under-strand which bounds $Q$ and $u^*$ the one which bounds $Q'$; thus $\{u,u^*\}=\{u_j,u'_j\}$. Each of them separates a covered quadrant from a free one: we write $Q^{\rm f}(u)$ for the free quadrant on the other side of $u$, and $Q^{\rm f}(u^*)$ for the free quadrant on the other side of $u^*$. Going around the crossing one meets $Q$, $u$, $Q^{\rm f}(u)$, a half-strand of the over-strand, $Q'$, $u^*$, $Q^{\rm f}(u^*)$ and the other half-strand of the over-strand. Consider a passage of $b$ along the core $\kappa_j$ from the arc in $Q$ to the arc in $Q'$. The band is vertical along the co-core, so its two sides at the middle point of the co-core face the two free quadrants, and we name the two push-offs of the passage after the half-strands $u$ and $u^*$ (Figure~\ref{fig.unrolled}):
\begin{itemize}
\item[$\kappa_j^{(u)}$:] over the bisector of $Q$, from the arc to the distance $r_1$ from the crossing, it is the core displaced upward by $\varrho$, above the half-band over $Q$. At the distance $r_1$ from the crossing it turns around the co-core: it crosses over the half-strand $u$, which lies in the plane $P$, into the free quadrant $Q^{\rm f}(u)$, where it descends by $2\varrho$, then crosses under the half-strand of the over-strand, which is at the height $h_0\chi(r_1)$, into $Q'$, and reaches the bisector of $Q'$ below the half-band over $Q'$. Over the bisector of $Q'$, from the distance $r_1$ to the arc, it is the core displaced downward by $\varrho$.
\item[$\kappa_j^{(u^*)}$:] over the bisector of $Q$ it is the core displaced downward by $\varrho$, below the half-band over $Q$. At the distance $r_1$ from the crossing it turns around the co-core the other way: it crosses under the half-strand of the over-strand which bounds $Q$ into the free quadrant $Q^{\rm f}(u^*)$, where it ascends by $2\varrho$, then crosses over the half-strand $u^*$ into $Q'$, above the half-band over $Q'$, and over the bisector of $Q'$ it is the core displaced upward by $\varrho$.
\end{itemize}
While they turn around the co-core, at the distance $r_1$ from it, both curves stay at heights between $\frac12h_0\chi(r_1)-\varrho$ and $\frac12h_0\chi(r_1)+\varrho$. We take $\varrho<\frac14h_0\chi(r_1)$, so that these heights lie strictly between the height $0$ of the half-strands of the under-strand and the height $h_0\chi(r_1)$ of those of the over-strand. Then the two curves do not meet the band: over $Q$ the half-band is the graph of $\beta$, which at the distance $r_1$ decreases from $h_0\chi(r_1)$ on the half-strand of the over-strand to $\frac12h_0\chi(r_1)$ on the bisector and to $0$ on $u$, so that $\kappa_j^{(u)}$, which moves from the bisector to $u$ at the height $\frac12h_0\chi(r_1)+\varrho$, stays above it, and $\kappa_j^{(u^*)}$, which moves from the bisector to the half-strand of the over-strand at the height $\frac12h_0\chi(r_1)-\varrho$, stays below it; the same holds over $Q'$; and there is no band over the free quadrants. Near the middle point of an arc of the smoothing, where a passage along a core is preceded or followed by a path in a disk, the push-offs of the two pieces are modified and joined, each to the one which lies on the same side of the surface, by short paths in the complement of the surface, at heights between $-\varrho$ and $\varrho$. When the band continues the disk across the arc, the push-off above the half-band is joined to the push-off above the disk, and the push-off below the half-band to the one below the disk. When the band is folded back over the disk (Figure~\ref{fig.statesurface}~(b), at the circle $C$), the push-off above the half-band passes over the arc and is joined to the push-off below the disk, and the push-off below the half-band turns back inside the fold and is joined to the push-off above the disk.

\begin{lemma}[one passage through a band]\label{lem.passage}
Let $\vs$ be a smoothing of $D$ and $j$ a crossing, with the notation above, and let $f(Q)$, $f(Q')$, $f(Q^{\rm f}(u))$ and $f(Q^{\rm f}(u^*))$ be the faces of $D$ whose corners are the four quadrants. If $\varrho$ is small enough, then for a passage along the core $\kappa_j$ from $Q$ to $Q'$:
\begin{enumerate}\renewcommand{\theenumi}{\roman{enumi}}
\item the push-off $\kappa_j^{(u)}$ passes once from below to above the face disk of $f(Q)$, once from above to below the face disk of $f(Q^{\rm f}(u))$, and it is below the face disk of $f(Q')$ all along;
\item the push-off $\kappa_j^{(u^*)}$ is below the face disk of $f(Q)$ all along, it passes once from below to above the face disk of $f(Q^{\rm f}(u^*))$, and once from above to below the face disk of $f(Q')$.
\end{enumerate}
Here ``once'' is meant algebraically, in the sense of Lemma~\ref{lem.spanning}~(iii), and a face disk is omitted if the face is the unbounded one. At its two ends, over the middle points of the two arcs of the smoothing, each of the two push-offs is below the face disk of the face in which it lies.
\end{lemma}
\begin{proof}
Let $f=f(Q)$. Over the bisector of $Q$, between the distances $r_1$ and $r_0$ from the crossing, the function $\psi_f$ is positive, hence bounded from below by a positive number; we take $\varrho$ smaller than this number, and likewise for $Q'$. Then over the bisector of $Q$ the curve $\kappa_j^{(u)}$, at the height $\beta+\varrho$, is below the face disk, which is at the height $\beta+\psi_f$, and so is $\kappa_j^{(u^*)}$, at the height $\beta-\varrho$. This holds in particular at the two ends of the passage.

(i) At the distance $r_1$, the curve $\kappa_j^{(u)}$ moves in $Q$ from the bisector, where it is below the face disk of $f$, to the half-strand $u$, at the constant height $\frac12h_0\chi(r_1)+\varrho>0$. On $u$ the face disk is at the height $0$, the height of the under-strand. So the curve leaves the cylinder over $f$ above the face disk: it has passed from below to above. It enters the cylinder over the face $f(Q^{\rm f}(u))$ over the same half-strand $u$, above the face disk of that face, which also contains $u$; it leaves it under the half-strand of the over-strand, at a height smaller than $h_0\chi(r_1)$, the height of this half-strand and of the face disk over it: it has passed from above to below. It then enters the cylinder over $f(Q')$ below the face disk, and it stays below the half-band over $Q'$, which is below the face disk by Lemma~\ref{lem.spanning}~(ii).

(ii) The curve $\kappa_j^{(u^*)}$ stays below the half-band over $Q$, hence below the face disk of $f(Q)$, until it leaves the cylinder over $f(Q)$ under the half-strand of the over-strand. It enters the cylinder over $f(Q^{\rm f}(u^*))$ below the face disk, and leaves it over the half-strand $u^*$, above the face disk: it has passed from below to above. It enters the cylinder over $f(Q')$ over $u^*$, above the face disk, and it reaches the bisector of $Q'$ at the height $\beta+\varrho$, below the face disk: it has passed from above to below.
\end{proof}

Figure~\ref{fig.unrolled} shows the whole argument. Its content is that the six sheets which contain the co-core, the two half-bands and the four face disks, are arranged around it in such a way that the face disk of a covered quadrant lies on the same side of the band as the face disk of the free quadrant which is separated from it by a half-strand of the \emph{under}-strand. This is where the stacking of the surface below the plane of the diagram enters, and it is the reason why the edges of the under-strand, and not those of the over-strand, appear in the formulas which follow.

\subsection{The localization theorem}\label{ssec.localization}
We can now state and prove the main result of this section. It says that the linking number of a flow with the push-off of a curve on a state surface is a sum over the crossings: each passage of the curve through the band of a crossing contributes the flow through the two edges of the under-strand of this crossing, with a sign which only depends on the band.

Recall the linear forms $U_j(z)=z_{u'_j}+z_{u_j}$ of \eqref{eq.UO} and the transfers $w^\vs_j$ of Definition~\ref{def.transfer}.

\begin{definition}\label{def.bandsign}
The \emph{sign} of the band of the crossing $j$ in the state surface $F_\vs$ is
\[
   \theta_j=\theta_j(\vs)=\ve_j\,(-1)^{\vs_j}\ \in\{1,-1\}.
\]
\end{definition}

The sign $\theta_j$ does not depend on the orientation of the link: reversing the orientation of one of the two strands of the crossing $j$ changes the sign $\ve_j$ of the crossing, and it exchanges the oriented smoothing with the other one. It records the sense of the half twist of the band, and in terms of the diagram it says on which side of the under-strand the covered quadrants lie:
\begin{equation}\label{eq.thetaside}
   \theta_j=1\iff\left\{\begin{array}{l}\text{the covered quadrant next to $u'_j$ lies on the left of $u'_j$, and}\\[1pt] \text{the covered quadrant next to $u_j$ lies on the right of $u_j$,}\end{array}\right.
\end{equation}
and for $\theta_j=-1$ the two sides are exchanged; the four cases are drawn in Figure~\ref{fig.cornersmoothing} of Section~\ref{sec.minusone}. Indeed, by \eqref{eq.mcrossing}, the faces on the left and on the right of $u'_j$ are $f^\leftarrow_j$ and $f^\uparrow_j$ if $\ve_j=1$, and $f^\uparrow_j$ and $f^\rightarrow_j$ if $\ve_j=-1$; those of $u_j$ are $f^\downarrow_j$ and $f^\rightarrow_j$ if $\ve_j=1$, and $f^\leftarrow_j$ and $f^\downarrow_j$ if $\ve_j=-1$; and the covered quadrants are those of $f^\leftarrow_j,f^\rightarrow_j$ if $\vs_j=0$ and those of $f^\uparrow_j,f^\downarrow_j$ if $\vs_j=1$.

\begin{theorem}[localization of linking numbers]\label{thm.localization}
Let $\vs$ be a smoothing of a diagram $D$, let $F_\vs$ be its state surface, let $b$ be a closed curve on $F_\vs$ in standard position and let $z\in\FF_D$ be a flow, regarded as a $1$-cycle of the graph $G\subset F_\vs$. Then
\[
   \lk(z,\tau b)=\sum_{j=1}^c\theta_j\,U_j(z)\,w^\vs_j(b) .
\]
More precisely, a passage of $b$ along the core of the band $b_j$, from the half-band of $u'_j$ to the half-band of $u_j$, contributes $\theta_j\,z_{u'_j}$ to $\lk(z,\tau b)$ through the push-off $\kappa_j^{(u'_j)}$ on the side of $u'_j$, and $\theta_j\,z_{u_j}$ through the push-off $\kappa_j^{(u_j)}$ on the side of $u_j$; a passage in the opposite direction contributes the opposite numbers; and the parts of $\tau b$ which are pushed off the disks do not contribute.
\end{theorem}
\begin{proof}
Let $k$ be the potential of $z$ which vanishes on the unbounded face. By \eqref{eq.lkpotential}, $\lk(z,\tau b)=\sum_fk_f\,\Delta_f\cdot\tau b$, the sum running over the bounded faces, and by Lemma~\ref{lem.spanning}~(iii) the number $\Delta_f\cdot\tau b$ is obtained by following $\tau b$ and recording, each time it travels over the face $f$, whether it enters and leaves the cylinder over $f$ below or above the face disk $\Delta_f$. We examine first the parts of $\tau b$ which are pushed off the disks, and then the push-offs of its passages along the cores of the bands.

We first show that, if $\varrho$ is small enough, the parts of $\tau b$ which are pushed off the disks are below the face disk of every bounded face over which they travel. Such a part is a copy of a path in a disk $D_C$, displaced vertically by $\varrho$, and it is joined at its ends to the push-offs of the passages which precede and follow it. A point of it which is displaced downward, or which lies over a point $p$ with $d_C(p)>\varrho$, has a negative height, and the face disks are at positive heights. The remaining points lie at heights between $0$ and $\varrho$ over points $p$ which are close to the two ends of the path, the middle points of two arcs of the smoothing; each of these two middle points lies in the interior of a face $f$, where $u_f\ge\psi_f>0$, and if $\varrho$ is small enough then $u_f>\varrho$ near them. The same applies to the short pieces which join the push-offs near these middle points. So these parts of $\tau b$ never pass from one side of a face disk to the other, and by Lemma~\ref{lem.spanning}~(iii) they contribute nothing; moreover, by the last statement of Lemma~\ref{lem.passage}, the push-offs of the passages begin and end below the face disks, so that the contributions of the passages can be computed one passage at a time.

Consider now a passage along $\kappa_j$ from the half-band of $u'_j$ to the half-band of $u_j$, and apply Lemma~\ref{lem.passage} with $u=u'_j$ and $u^*=u_j$. Let $f^{\rm c}$ and $f^{\rm f}$ be the faces of the covered and of the free quadrant next to $u'_j$. By Lemma~\ref{lem.passage}~(i) the push-off $\kappa_j^{(u'_j)}$ contributes $k_{f^{\rm c}}-k_{f^{\rm f}}$. The edge $u'_j$ separates these two faces, and $z_{u'_j}$ is the value of $k$ on its left minus the value on its right by \eqref{eq.potential}; by \eqref{eq.thetaside} the covered quadrant is on the left of $u'_j$ exactly when $\theta_j=1$. So the contribution is $\theta_j\,z_{u'_j}$. Let now $f^{\rm c}$ and $f^{\rm f}$ be the faces of the covered and of the free quadrant next to $u_j$. By Lemma~\ref{lem.passage}~(ii) the push-off $\kappa_j^{(u_j)}$ contributes $k_{f^{\rm f}}-k_{f^{\rm c}}$, and by \eqref{eq.thetaside} the covered quadrant is on the right of $u_j$ exactly when $\theta_j=1$; so this contribution is $\theta_j\,z_{u_j}$. (If one of the faces is the unbounded one, its face disk is missing and its potential is zero, and the formulas remain true.) For a passage in the other direction the two push-offs are the same curves traveled backward, and the contributions change sign. Summing over the passages of $b$ gives the formula: by definition $w^\vs_j(b)$ is the number of passages from the half-band of $u'_j$ to that of $u_j$ minus the number of passages in the other direction, and each of them contributes $\theta_j\bigl(z_{u'_j}+z_{u_j}\bigr)=\theta_j\,U_j(z)$ in total.
\end{proof}

The first consequence concerns the form of Gordon and Litherland\label{not.GL} \cite{GordonLitherland}. For a compact surface $F\subset\BR^3$, orientable or not, it is the symmetric bilinear form
\[
   \CG_F\colon H_1(F;\BZ)\times H_1(F;\BZ)\to\BZ,\qquad \CG_F([a],[b])=\lk(a,\tau b),
\]
where $a$ and $b$ are closed curves on $F$; if $F$ is orientable, $\CG_F([a],[b])=\lk(a,b^+)+\lk(a,b^-)$ is the symmetrized Seifert form.

\begin{corollary}[the Gordon--Litherland form of a state surface]\label{cor.GLform}
For every smoothing $\vs$ of $D$ and all flows $z,z'\in\FF_D$, with classes $[z],[z']$ in $H_1(F_\vs;\BZ)=\FF_D/\langle z^C\rangle$,
\begin{equation}\label{eq.GLform}
   \CG_{F_\vs}\bigl([z],[z']\bigr)=\sum_{j=1}^c\theta_j\,U_j(z)\,w^\vs_j(z') .
\end{equation}
In particular the right-hand side is a symmetric bilinear form on $\FF_D$ which vanishes when $z$ or $z'$ is a circle of $\vs$.
\end{corollary}
\begin{proof}
The class $[z']$ is a sum of classes of curves $b$ in standard position, whose transfers add up to $w^\vs(z')$ by Proposition~\ref{prop.statehomology}; and the cycle $z$ of $G$ lies on $F_\vs$ and represents the class $[z]$. Both sides of the formula of Theorem~\ref{thm.localization} are additive in $b$.
\end{proof}

The symmetry of the right-hand side of \eqref{eq.GLform} and its vanishing on the circles are proved algebraically in Section~\ref{sec.minusone}, where this form is studied together with the Goeritz matrices and Kashaev's matrix.

\begin{figure}[ht]
\centering
\begin{tikzpicture}[>=Stealth,scale=0.95,every node/.style={font=\small}]
\begin{scope}
\draw (1.2,-1.2) -- (0.2,-0.2); \draw[->] (-0.2,0.2) -- (-1.2,1.2);
\draw[->] (-1.2,-1.2) -- (1.2,1.2);
\node[font=\scriptsize] at (-2.0,-1.3) {$o_j=\hat e^+_j$}; \node[font=\scriptsize] at (2.0,1.3) {$o'_j=\hat e^-_j$};
\node[font=\scriptsize] at (2.0,-1.3) {$u_j=e^-_j$}; \node[font=\scriptsize] at (-2.0,1.3) {$u'_j=e^+_j$};
\node[font=\scriptsize] at (-1.25,0) {$C_{\rm L}$}; \node[font=\scriptsize] at (1.25,0) {$C_{\rm R}$};
\node at (0,-2.0) {$\ve_j=1$};
\end{scope}
\begin{scope}[xshift=7.6cm]
\draw (-1.2,-1.2) -- (-0.2,-0.2); \draw[->] (0.2,0.2) -- (1.2,1.2);
\draw[->] (1.2,-1.2) -- (-1.2,1.2);
\node[font=\scriptsize] at (2.0,-1.3) {$o_j=\hat e^-_j$}; \node[font=\scriptsize] at (-2.0,1.3) {$o'_j=\hat e^+_j$};
\node[font=\scriptsize] at (-2.0,-1.3) {$u_j=e^+_j$}; \node[font=\scriptsize] at (2.0,1.3) {$u'_j=e^-_j$};
\node[font=\scriptsize] at (-1.25,0) {$C_{\rm L}$}; \node[font=\scriptsize] at (1.25,0) {$C_{\rm R}$};
\node at (0,-2.0) {$\ve_j=-1$};
\end{scope}
\end{tikzpicture}
\caption{The edges $e^{\pm}_j$ of the under-strand and $\hat e^{\pm}_j$ of the over-strand of a crossing. The edges with the sign $+$ are the two edges on the left, on the Seifert circle $C_{\rm L}$; those with the sign $-$ are on the right, on $C_{\rm R}$. The Seifert form of the canonical surface below the plane reads the flow through $e^+_j$, and its transpose the flow through $e^-_j$; for the surface above the plane they read the flows through $\hat e^-_j$ and $\hat e^+_j$ (Corollary~\ref{cor.seifertform}, Proposition~\ref{prop.otherside} and Theorem~\ref{thm.deformedseifert}).}
\label{fig.sides}
\end{figure}
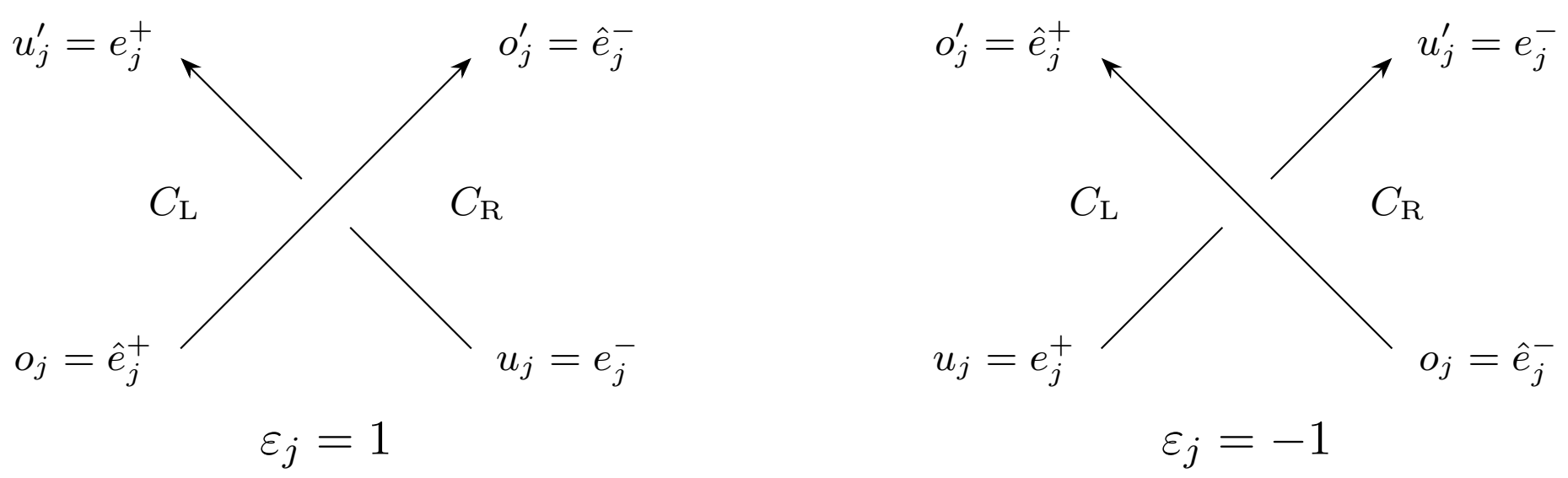

For the canonical surface the theorem gives more, because the two push-offs can be told apart globally. We orient $F_D$ so that its oriented boundary is $L$, and we write $b^+$ and $b^-$ for the push-offs of a curve $b$ on $F_D$ to the positive and to the negative side. The \emph{Seifert form} of $F_D$ is the bilinear form\label{not.VF} $V_F([a],[b])=\lk(a,b^+)$ on $H_1(F_D;\BZ)$; it satisfies $\lk(a,b^-)=V_F([b],[a])$, since $\lk(a,b^-)=\lk(a^+,b)$. For a crossing $j$ let $e^+_j$ be the edge of the under-strand which lies on the left of the over-strand, and $e^-_j$ the one which lies on its right (Figure~\ref{fig.sides}):
\begin{equation}\label{eq.eplus}
   e^+_j=\begin{cases}u'_j&\text{if }\ve_j=1,\\ u_j&\text{if }\ve_j=-1,\end{cases}\qquad\qquad e^-_j=\begin{cases}u_j&\text{if }\ve_j=1,\\ u'_j&\text{if }\ve_j=-1.\end{cases}
\end{equation}
Recall that $w_j(z)=z_{o_j}-z_{u'_j}$ is the transfer of the oriented smoothing, the number of times the cycle $z$ runs through the band of the crossing $j$ from the Seifert circle $A_j$ to the Seifert circle $B_j$.

\begin{corollary}[the Seifert form of the canonical surface]\label{cor.seifertform}
Let $D$ be a diagram and $F_D$ its canonical Seifert surface. For all flows $z,z'\in\FF_D$, and for a cycle $b$ on $F_D$ in standard position which represents the class of $z'$,
\[
   \lk(z,b^+)=\sum_{j=1}^c\ve_j\,z_{e^+_j}\,w_j(z'),\qquad \lk(z,b^-)=\sum_{j=1}^c\ve_j\,z_{e^-_j}\,w_j(z') .
\]
Hence the Seifert form of $F_D$ is given, on the classes of the flows, by
\begin{equation}\label{eq.seifertform}
   V_F\bigl([z],[z']\bigr)=\sum_{j=1}^c\ve_j\,z_{e^+_j}\,\bigl(z'_{o_j}-z'_{u'_j}\bigr),
\end{equation}
and the second formula computes $V_F([z'],[z])$.
\end{corollary}
\begin{proof}
For $\vs=0$ we have $\theta_j=\ve_j$ and $w^\vs_j=w_j$. By Theorem~\ref{thm.localization}, a passage of $b$ through the band of the crossing $j$ from $A_j$ to $B_j$, that is, from the half-band of $u'_j$ to that of $u_j$, contributes $\ve_jz_{u'_j}$ through the push-off on the side of $u'_j$ and $\ve_jz_{u_j}$ through the push-off on the side of $u_j$. So it remains to show that \emph{at the middle point of the co-core $\alpha_j$ the positive normal of $F_D$ points to the side of $e^+_j$}, that is, into the free quadrant next to $e^+_j$. Along the co-core the band is vertical: at the height $\frac12h_0(1+\tau)$ it leaves $\alpha_j$ in the two directions of angular coordinate $\tau$ of the covered quadrants (Section~\ref{ssec.statesurface}). So its normal along $\alpha_j$ is horizontal and perpendicular to directions which lie in the covered quadrants: it points into one of the two free quadrants, and by continuity into the same one at all heights. At the upper end of $\alpha_j$, which is a point of the over-strand, the band lies below the over-strand; so the outward normal of the surface is $e_z$ there, and the positive normal is $e_z\times t_o$, where $t_o$ is the direction of the over-strand: the horizontal vector which points to the left of the over-strand. Hence, all along the co-core, the positive normal points into the free quadrant which lies on the left of the over-strand. The line of the over-strand cuts the disk $\CB_j$ into two halves, each of which consists of a covered and a free quadrant, and the free quadrant in the left half is bounded by a half-strand of the over-strand and by the half-strand of the under-strand which lies on the left of the over-strand, that is, by the end of the edge $e^+_j$.
\end{proof}

\begin{example}[the Seifert matrix of the trefoil]\label{ex.trefoilseifert}
For the diagram of the trefoil of Example~\ref{ex.trefoilflows} all crossings are positive and $e^+_j=u'_j$, so that $V_F([z],[z'])=\sum_jz_{u'_j}w_j(z')$. The lifts of the cycles $\gamma_1=e_2-e_1$ and $\gamma_2=e_3-e_1$ are $z^{\gamma_1}=e_{u'_1}-e_{u_2}$ and $z^{\gamma_2}=e_{u'_1}+e_{u'_2}-e_{u_2}-e_{u_3}$; their coordinates at the edges $(u'_1,u'_2,u'_3)$ are $(1,0,0)$ and $(1,1,0)$, and their transfers are the coordinates $(-1,1,0)$ and $(-1,0,1)$ of the cycles. Hence
\[
   V_F=\bigl(V_F(\gamma_a,\gamma_b)\bigr)_{a,b}=\begin{pmatrix}1\cdot(-1)&1\cdot(-1)\\ 1\cdot(-1)+1\cdot1&1\cdot(-1)+1\cdot0\end{pmatrix}=\begin{pmatrix}-1&-1\\0&-1\end{pmatrix},
\]
the familiar Seifert matrix of the right-handed trefoil, with $\det(tV_F-V_F^{\T})=1-t+t^2$.
\end{example}

\subsection{The state surfaces above the plane}\label{ssec.otherside}
The construction of Section~\ref{ssec.statesurface} pushes all the disks below the plane of the diagram. Pushing them above it gives another surface with the same boundary, in which the disk of an inner circle lies \emph{below} the disks of the circles which enclose it. We need this second surface because the second flavor of the determinant formula belongs to it, as the first flavor belongs to $F_D$ (Theorem~\ref{thm.deformedseifert}). We define it by reflection. Let $D^*$ be the diagram obtained from $D$ by changing all the crossings; it has the same edges as $D$, with the same orientations, and the same smoothings, but the over-strand of $D$ is the under-strand of $D^*$ and vice versa, and the signs of the crossings are opposite. Let $\varrho_P$ be the reflection of $\BR^3$ in the plane $P$. The image under $\varrho_P$ of the link of $D^*$, realized as in Section~\ref{ssec.statesurface}, coincides with $D$ outside the disks $\CB_j$; over $\CB_j$ the over-strand of $D$ lies in $P$ and the under-strand dips below it. Lifting the point over $p\in\CB_j$ by $h_0\chi(r(p))$, for every crossing $j$, moves this link onto $L$; let $\phi$ be this vertical isotopy of $\BR^3$, which is the identity outside the cylinders over the disks $\CB_j$.

\begin{definition}\label{def.uppersurface}
For a smoothing $\vs$ of $D$, the \emph{state surface above the plane} is $F^\uparrow_\vs=\phi\,\varrho_P\bigl(F_\vs(D^*)\bigr)$, where $F_\vs(D^*)$ is the state surface of the smoothing $\vs$ of the diagram $D^*$. For $\vs=0$ we write $F^\uparrow_D=F^\uparrow_0$ and call it the canonical Seifert surface of $D$ above the plane; we orient it so that its oriented boundary is $L$.
\end{definition}

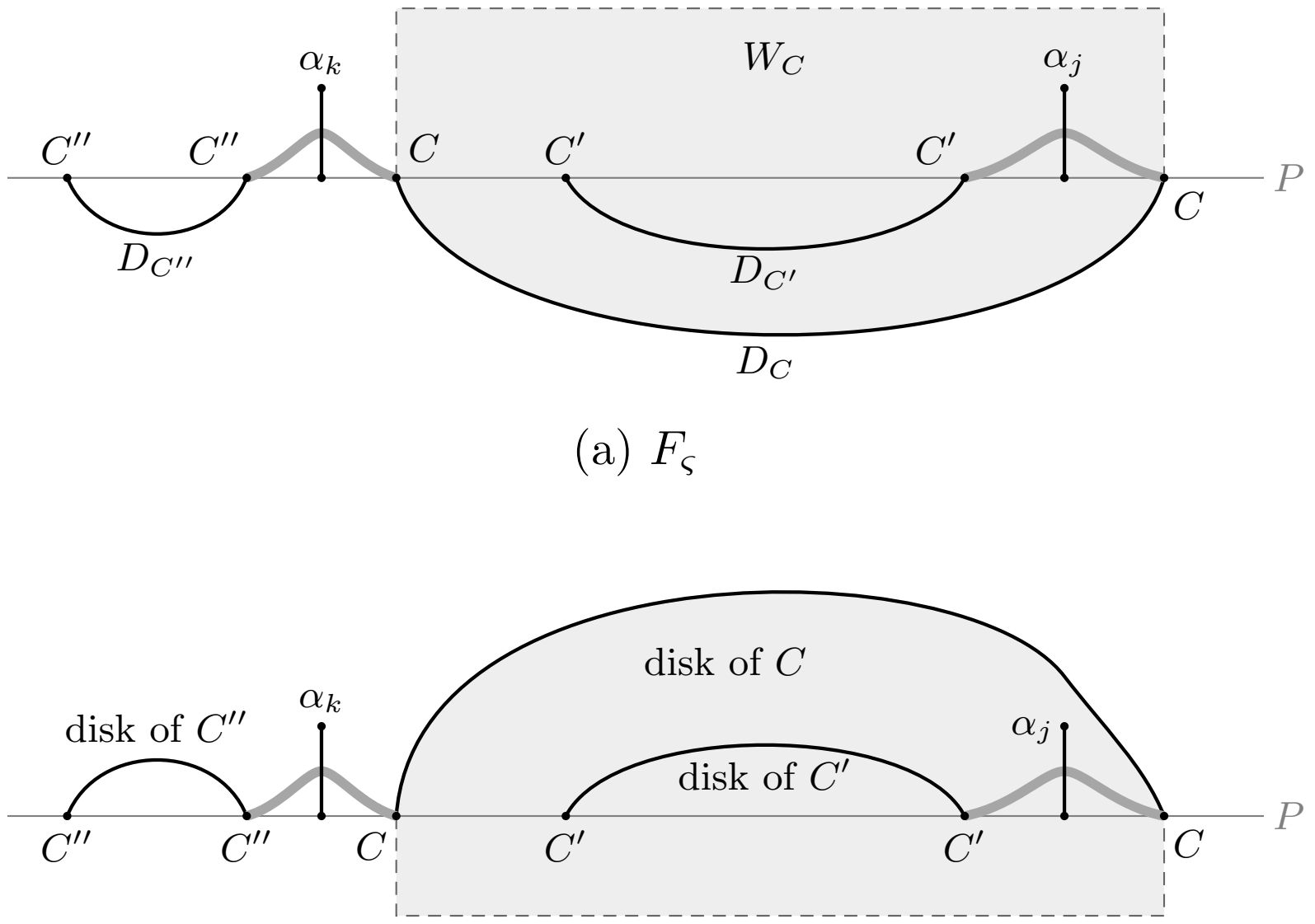
\begin{figure}[ht]
\centering
\begin{tikzpicture}[>=Stealth,every node/.style={font=\small},scale=0.8]
\begin{scope}
\fill[gray!13] (0.3,0) .. controls (0.9,-2.1) and (7.4,-2.1) .. (8,0) -- (8,1.7) -- (0.3,1.7) -- cycle;
\draw[densely dashed,gray!80!black] (0.3,0) -- (0.3,1.7) -- (8,1.7) -- (8,0);
\node[font=\scriptsize] at (4.1,1.2) {$W_C$};
\draw[gray] (-3.6,0) -- (9.0,0); \node[gray,font=\scriptsize] at (9.25,0) {$P$};
\draw[thick] (0.3,0) .. controls (0.9,-2.1) and (7.4,-2.1) .. (8,0);
\draw[thick] (2,0) .. controls (2.5,-0.95) and (5.5,-0.95) .. (6,0);
\draw[thick] (-3,0) .. controls (-2.7,-0.75) and (-1.5,-0.75) .. (-1.2,0);
\node[font=\scriptsize] at (4,-0.95) {$D_{C'}$}; \node[font=\scriptsize] at (4,-1.85) {$D_C$}; \node[font=\scriptsize] at (-2.1,-0.82) {$D_{C''}$};
\draw[line width=2.2pt,gray!70] (6,0) .. controls (6.5,0.1) and (6.8,0.45) .. (7,0.45) .. controls (7.2,0.45) and (7.5,0.1) .. (8,0);
\draw[line width=2.2pt,gray!70] (-1.2,0) .. controls (-0.85,0.1) and (-0.6,0.45) .. (-0.45,0.45) .. controls (-0.3,0.45) and (-0.05,0.1) .. (0.3,0);
\draw[thick] (7,0) -- (7,0.9); \draw[thick] (-0.45,0) -- (-0.45,0.9);
\fill (7,0) circle (1.2pt); \fill (7,0.9) circle (1.2pt); \fill (-0.45,0) circle (1.2pt); \fill (-0.45,0.9) circle (1.2pt);
\foreach \x in {-3,-1.2,0.3,2,6,8} \fill (\x,0) circle (1.2pt);
\node[font=\scriptsize] at (-3,0.3) {$C''$}; \node[font=\scriptsize] at (-1.52,0.3) {$C''$};
\node[font=\scriptsize] at (0.58,0.3) {$C$}; \node[font=\scriptsize] at (8.25,-0.28) {$C$};
\node[font=\scriptsize] at (2,0.3) {$C'$}; \node[font=\scriptsize] at (5.72,0.3) {$C'$};
\node[font=\scriptsize] at (7,1.15) {$\alpha_j$}; \node[font=\scriptsize] at (-0.45,1.15) {$\alpha_k$};
\node at (2.7,-2.75) {(a) $F_\vs$};
\end{scope}
\begin{scope}[yshift=-6.4cm]
\fill[gray!13] (0.3,0) .. controls (0.5,2.7) and (6.0,2.7) .. (7.0,1.4) .. controls (7.4,0.88) and (7.8,0.5) .. (8,0) -- (8,-1.0) -- (0.3,-1.0) -- cycle;
\draw[densely dashed,gray!80!black] (0.3,0) -- (0.3,-1.0) -- (8,-1.0) -- (8,0);
\draw[gray] (-3.6,0) -- (9.0,0); \node[gray,font=\scriptsize] at (9.25,0) {$P$};
\draw[thick] (0.3,0) .. controls (0.5,2.7) and (6.0,2.7) .. (7.0,1.4) .. controls (7.4,0.88) and (7.8,0.5) .. (8,0);
\draw[thick] (2,0) .. controls (2.5,0.95) and (5.5,0.95) .. (6,0);
\draw[thick] (-3,0) .. controls (-2.7,0.75) and (-1.5,0.75) .. (-1.2,0);
\node[font=\scriptsize] at (4,0.42) {disk of $C'$}; \node[font=\scriptsize] at (3.6,1.55) {disk of $C$}; \node[font=\scriptsize] at (-2.1,0.9) {disk of $C''$};
\draw[line width=2.2pt,gray!70] (6,0) .. controls (6.5,0.1) and (6.8,0.45) .. (7,0.45) .. controls (7.2,0.45) and (7.5,0.1) .. (8,0);
\draw[line width=2.2pt,gray!70] (-1.2,0) .. controls (-0.85,0.1) and (-0.6,0.45) .. (-0.45,0.45) .. controls (-0.3,0.45) and (-0.05,0.1) .. (0.3,0);
\draw[thick] (7,0) -- (7,0.9); \draw[thick] (-0.45,0) -- (-0.45,0.9);
\fill (7,0) circle (1.2pt); \fill (7,0.9) circle (1.2pt); \fill (-0.45,0) circle (1.2pt); \fill (-0.45,0.9) circle (1.2pt);
\foreach \x in {-3,-1.2,0.3,2,6,8} \fill (\x,0) circle (1.2pt);
\node[font=\scriptsize] at (-3,-0.3) {$C''$}; \node[font=\scriptsize] at (-1.2,-0.3) {$C''$};
\node[font=\scriptsize] at (0.05,-0.3) {$C$}; \node[font=\scriptsize] at (8.25,-0.28) {$C$};
\node[font=\scriptsize] at (2,-0.3) {$C'$}; \node[font=\scriptsize] at (6,-0.3) {$C'$};
\node[font=\scriptsize] at (6.68,0.86) {$\alpha_j$}; \node[font=\scriptsize] at (-0.45,1.15) {$\alpha_k$};
\node at (2.7,-1.75) {(b) $F^\uparrow_\vs$};
\end{scope}
\end{tikzpicture}
\caption{Vertical sections of the two state surfaces of a smoothing $\vs$, along the cores of the bands of two crossings $j$ and $k$. The circle $C$ encloses the circle $C'$; the crossing $j$ lies inside $C$, between $C'$ and $C$, and the crossing $k$ lies outside $C$, between $C$ and a third circle $C''$. (a) The surface $F_\vs$: the disks hang below the plane $P$, the disk of the inner circle $C'$ above the disk of $C$, and the bands (gray) lie above $P$. The part of the surface which lies inside $C$ is contained in the region $W_C$ above the disk $D_C$ (shaded; see Section~\ref{ssec.murasugi}). (b) The surface $F^\uparrow_\vs$: the bands are the same, the disks lie above $P$, the disk of $C'$ below the disk of $C$, and the disk of $C$ passes above the band of the crossing $j$. The part of the surface which lies inside $C$ is now below the disk of $C$ (shaded).}
\label{fig.twosides}
\end{figure}

Explicitly, $F^\uparrow_\vs$ has the same bands as $F_\vs$. Indeed the angular coordinate of the diagram $D^*$ is $-\tau$, so the band of the crossing $j$ of $D^*$ is the graph of $\frac12h_0\chi(r)(1-\tau)$ over the covered quadrants; its reflection is the graph of $-\frac12h_0\chi(r)(1-\tau)$, and $\phi$ lifts it by $h_0\chi(r)$ onto the graph of $\frac12h_0\chi(r)(1+\tau)=\beta$, which is $b_j$. The disks are different: the disk of a circle $C$ in $F^\uparrow_\vs$ is the graph of $d_C$ over $\Omega_C$, lifted by $h_0\chi(r)$ over the disks $\CB_j$ which are contained in $\Omega_C$, so that it lies above the plane $P$ and passes above the bands of the crossings inside $C$ (Figure~\ref{fig.twosides}).

The map $\phi\varrho_P$ sends the graph $G$ of $D^*$ onto the graph $G$ of $D$, each edge onto the same edge and each co-core onto the co-core of the same crossing with the opposite orientation; so it sends the cycle of a flow $z$ onto the cycle of the same flow $z$, by \eqref{eq.flowcycle} and \eqref{eq.flow}. Hence $F^\uparrow_\vs$ contains $G$, its first homology group is again $\FF_D/\langle z^C\rangle$, and it has the same $2$-cells as $F_\vs$. The two surfaces are nevertheless not isotopic in general, and their linking forms differ.

\begin{proposition}[the surfaces above the plane]\label{prop.otherside}
Let $\vs$ be a smoothing of $D$ and $z,z'\in\FF_D$. Recall that $O_j(z)=z_{o_j}+z_{o'_j}$.
\begin{enumerate}\renewcommand{\theenumi}{\roman{enumi}}
\item The Gordon--Litherland form of the state surface above the plane is
\begin{equation}\label{eq.GLup}
   \CG_{F^\uparrow_\vs}\bigl([z],[z']\bigr)=-\sum_{j=1}^c\ve_j\,O_j(z)\,w^\vs_j(z') .
\end{equation}
\item For a crossing $j$ let $\hat e^+_j$ be the edge of the over-strand which lies on the left of the under-strand and $\hat e^-_j$ the one which lies on its right: $\hat e^+_j=o_j$, $\hat e^-_j=o'_j$ if $\ve_j=1$, and $\hat e^+_j=o'_j$, $\hat e^-_j=o_j$ if $\ve_j=-1$. The Seifert form of $F^\uparrow_D$ is given by
\begin{equation}\label{eq.seifertup}
   V_{F^\uparrow}\bigl([z],[z']\bigr)=-\sum_{j=1}^c\ve_j\,z_{\hat e^-_j}\,w_j(z'),\qquad V_{F^\uparrow}\bigl([z'],[z]\bigr)=-\sum_{j=1}^c\ve_j\,z_{\hat e^+_j}\,w_j(z') .
\end{equation}
\end{enumerate}
\end{proposition}
\begin{proof}
A reflection changes the signs of the linking numbers and an isotopy preserves them. We apply Theorem~\ref{thm.localization} and Corollary~\ref{cor.seifertform} to $D^*$ and translate. The crossing $j$ of $D^*$ has the sign $-\ve_j$, its edges are $o^*_j=u_j$, $u^*_j=o_j$, $o'^*_j=u'_j$ and $u'^*_j=o'_j$, and the smoothing $\vs$ has the same arcs for $D^*$ as for $D$. Hence the band sign of $D^*$ is $\theta^*_j=-\theta_j$, the form $U_j$ of $D^*$ is $z_{o'_j}+z_{o_j}=O_j(z)$, and the transfer of $D^*$ is $z_{u_j}-z_{o'_j}=-w^\vs_j(z)$ if $\vs_j=0$, by \eqref{eq.flow}, and $-z_{u'_j}-z_{o'_j}=w^\vs_j(z)$ if $\vs_j=1$; in both cases $\theta^*_j$ times the transfer of $D^*$ is $\ve_j\,w^\vs_j(z)$.

(i) For a curve $b$ in standard position on $F_\vs(D^*)$ in the class of $z'$, the map $\phi\varrho_P$ sends $\tau b$ onto the double push-off of the image of $b$. So $\CG_{F^\uparrow_\vs}([z],[z'])=-\lk(z,\tau b)$, the linking number being computed for $D^*$, and Theorem~\ref{thm.localization} gives $-\sum_j\theta^*_j\,O_j(z)\,w^{\vs,*}_j(z')=-\sum_j\ve_jO_j(z)w^\vs_j(z')$.

(ii) The map $\phi\varrho_P$ sends the link of $D^*$ onto $L$ with its orientation, hence the oriented surface $F_{D^*}$ onto the oriented surface $F^\uparrow_D$; since it reverses the orientation of $\BR^3$, it exchanges the positive and the negative sides. So $V_{F^\uparrow}([z],[z'])$ is minus the linking number of $z$ with the push-off of $b$ to the negative side of $F_{D^*}$, which by Corollary~\ref{cor.seifertform} is $-\sum_j(-\ve_j)\,z_{e^{-,*}_j}\,w^*_j(z')$, where $w^*_j=-w_j$ and $e^{-,*}_j$ is the edge of the under-strand of $D^*$ on the right of the over-strand of $D^*$: this is the edge of the over-strand of $D$ on the right of the under-strand of $D$, namely $\hat e^-_j$. The second formula is obtained in the same way from the push-off to the positive side of $F_{D^*}$.
\end{proof}

So the edges of the under-strands govern the surfaces below the plane and those of the over-strands the surfaces above it. If a Seifert circle $C$ has crossings on both sides, $F_D$ and $F^\uparrow_D$ are the two Murasugi sums, along the disk of $C$, of the canonical surfaces of the parts of the diagram inside and outside $C$ (Section~\ref{ssec.murasugi}): in both of them the two summands lie on opposite sides of the disk, and the two surfaces differ by the side on which each summand lies. We shall see that the first flavor of the determinant formula belongs to $F_D$ and the second one to $F^\uparrow_D$ (Theorem~\ref{thm.deformedseifert}), and that the Gordon--Litherland forms of the two families of state surfaces are related by an involution (Proposition~\ref{prop.J}). When no circle of $\vs$ has crossings on both sides, for instance for the two checkerboard surfaces of Section~\ref{sec.minusone}, the two surfaces $F_\vs$ and $F^\uparrow_\vs$ are isotopic: each disk can then be moved from one side of the plane to the other without meeting the rest of the surface, a circle which encloses other circles being the outermost one, whose disk is moved through the point at infinity. Accordingly the two formulas \eqref{eq.GLform} and \eqref{eq.GLup} agree in that case (Proposition~\ref{prop.J}).

\section{The deformed Seifert form}\label{sec.seifert}

The one-variable Alexander polynomial is $\det(tV-V^{\T})$ up to a unit, where $V$ is a Seifert matrix of the link: the matrix, in a basis of the first homology group of a Seifert surface, of the Seifert form $(a,b)\mapsto\lk(a,b^+)$. In this section we show that the determinant formula has this form for the canonical Seifert surface $F_D$ of every diagram, with the crossing variables included. The key is that the crossing forms of Section~\ref{sec.flows} \emph{are} the Seifert form: by the localization theorem, the Seifert form of $F_D$ is a sum over the bands of the surface, the contribution of the band of a crossing being the product of a flow through an edge of its under-strand and of the number of passages through the band, and the crossing form $\pi_j$ is the contribution of the band $j$ to $tV_F-V_F^{\T}$ in which $t$ has been replaced by the variable of the crossing (Theorem~\ref{thm.deformedseifert}). The second flavor is related in the same way to the canonical surface on the other side of the plane of the diagram.

We then pass from the $c\times c$ matrix of the crossing forms in a Seifert basis to a $\nu\times\nu$ matrix on $H_1(F_D)$, the \emph{deformed Seifert matrix} $M_{D,\hate}(x)$, whose determinant is the edge polynomial up to an explicit monomial and with no sign (Corollary~\ref{cor.seifert}); at the Alexander point it is $tV_F-V_F^{\T}$ up to units on the columns, for every diagram (Theorem~\ref{thm.alexpoint}). Section~\ref{ssec.cochains} explains in which sense the crossing variables are cochains on the surface and why they cannot be replaced by cohomology classes. The section ends with examples, with the Kinoshita--Terasaka knot, which shows that no formula in terms of the Seifert form alone can exist, and with the behavior of the edge polynomials under the moves which relate the closed braid diagrams of a link.

Throughout, $D$ is a diagram, $\hate$ an edge of $D$, and we use the notation of Sections~\ref{sec.flows} and~\ref{sec.statesurfaces}: the Seifert circles $C$ of $D$ and the circle $\hat C$ of $\hate$, the Seifert graph $\Gamma_D$, its cycle space $H$ of rank $\nu=c-s+1$, the transfers $w_j(z)=z_{o_j}-z_{u'_j}$ of the oriented smoothing, and the canonical Seifert surfaces $F_D$ and $F^\uparrow_D$ of $D$ below and above the plane, with $H_1(F_D)=H_1(F^\uparrow_D)=\FF_D/\langle C_1,\dots,C_s\rangle\cong H$.

\subsection{The crossing forms are the Seifert form, band by band}\label{ssec.bandbyband}
This subsection contains the observation on which the whole section rests: each of the two crossing forms of a crossing is a combination of two linear forms, which are the contributions of the band of this crossing to the Seifert form of a canonical surface and to its transpose. We begin by giving names to these contributions.

Recall from \eqref{eq.eplus} and Proposition~\ref{prop.otherside} the four edges $e^+_j,e^-_j$ (the edges of the under-strand on the left and on the right of the over-strand) and $\hat e^+_j,\hat e^-_j$ (the edges of the over-strand on the left and on the right of the under-strand) of a crossing $j$. When the crossing is drawn with both strands pointing upward, $e^+_j$ and $\hat e^+_j$ are the two edges on the left, which lie on the Seifert circle $C_{\rm L}$ on the left of the crossing, and $e^-_j$ and $\hat e^-_j$ are the two edges on the right, which lie on the Seifert circle $C_{\rm R}$ (Figure~\ref{fig.sides}); by the proof of Lemma~\ref{lem.level}, $h(C_{\rm R})=h(C_{\rm L})+1$. For a flow $z\in\FF_D$ we put
\begin{equation}\label{eq.Vj}
   V_j(z)=\ve_j\,z_{e^+_j},\qquad V'_j(z)=\ve_j\,z_{e^-_j},\qquad V^\uparrow_j(z)=-\ve_j\,z_{\hat e^-_j},\qquad V'^{\uparrow}_j(z)=-\ve_j\,z_{\hat e^+_j} .
\end{equation}
So $V_j(z)$ and $V'_j(z)$ are the flows of $z$ through the two edges of the under-strand of the crossing $j$, with the sign of the crossing, and $V^\uparrow_j(z)$ and $V'^{\uparrow}_j(z)$ are the flows through the two edges of its over-strand, with the opposite sign. By the localization theorem, the products of these numbers and of the number $w_j(z')$ of passages of a second flow $z'$ through the band of the crossing $j$ are the contributions of this band to the Seifert forms of the two canonical surfaces. This is part (i) of the following theorem; parts (ii) and (iii) express the crossing forms in terms of the four forms \eqref{eq.Vj}.

\begin{theorem}[the crossing forms and the Seifert forms]\label{thm.deformedseifert}
Let $D$ be a diagram, and let $V_F$ and $V_{F^\uparrow}$ be the Seifert forms of its canonical Seifert surfaces $F_D$ and $F^\uparrow_D$ below and above the plane.
\begin{enumerate}\renewcommand{\theenumi}{\roman{enumi}}
\item \emph{The Seifert forms, band by band.} For all flows $z,z'\in\FF_D$,
\begin{align*}
   V_F\bigl([z],[z']\bigr)&=\sum_{j=1}^cV_j(z)\,w_j(z'),& V_F\bigl([z'],[z]\bigr)&=\sum_{j=1}^cV'_j(z)\,w_j(z'),\\
   V_{F^\uparrow}\bigl([z],[z']\bigr)&=\sum_{j=1}^cV^\uparrow_j(z)\,w_j(z'),& V_{F^\uparrow}\bigl([z'],[z]\bigr)&=\sum_{j=1}^cV'^{\uparrow}_j(z)\,w_j(z') .
\end{align*}
\item \emph{The crossing forms.} For every crossing $j$,
\[
   \pi_j=x_j^{\ind{\ve_j=1}}\,V_j-x_j^{\ind{\ve_j=-1}}\,V'_j,\qquad \tpi_j=y_j^{\ind{\ve_j=-1}}\,V^\uparrow_j-y_j^{\ind{\ve_j=1}}\,V'^{\uparrow}_j .
\]
\item \emph{Special points.} At $x_j=y_j=1$ both crossing forms are $V_j-V'_j=V^\uparrow_j-V'^{\uparrow}_j$, and at the Alexander point $x_j=y_j=t^{\ve_j}$
\[
   \pi_j=t^{-\ind{\ve_j=-1}}\bigl(t\,V_j-V'_j\bigr),\qquad \tpi_j=t^{-\ind{\ve_j=-1}}\bigl(V^\uparrow_j-t\,V'^{\uparrow}_j\bigr).
\]
\end{enumerate}
\end{theorem}
\begin{proof}
(i) is Corollary~\ref{cor.seifertform} and Proposition~\ref{prop.otherside}~(ii). (ii) If $\ve_j=1$, then $e^+_j=u'_j$ and $e^-_j=u_j$, so $x_jV_j-V'_j=x_jz_{u'_j}-z_{u_j}$; if $\ve_j=-1$, then $e^+_j=u_j$ and $e^-_j=u'_j$, so $V_j-x_jV'_j=-z_{u_j}+x_jz_{u'_j}$. In both cases this is $\pi_j(z)$. If $\ve_j=1$, then $\hat e^-_j=o'_j$ and $\hat e^+_j=o_j$, so $V^\uparrow_j-y_jV'^{\uparrow}_j=-z_{o'_j}+y_jz_{o_j}$; if $\ve_j=-1$, then $\hat e^-_j=o_j$ and $\hat e^+_j=o'_j$, so $y_jV^\uparrow_j-V'^{\uparrow}_j=y_jz_{o_j}-z_{o'_j}$. In both cases this is $\tpi_j(z)$. (iii) follows from (ii): for $\ve_j=1$, $\pi_j=tV_j-V'_j$ at $x_j=t$, and for $\ve_j=-1$, $\pi_j=V_j-t^{-1}V'_j$ at $x_j=t^{-1}$; likewise for $\tpi_j$. At $x_j=y_j=1$ the two forms coincide by \eqref{eq.piatone}.
\end{proof}

Together with Theorem~\ref{thm.flows} this gives the determinant formula the following form. Let $z_1,\dots,z_c$ be a basis of $\FF_{D,\hate}$. Then, up to sign, $P_{D,\hate}(x)$ is the determinant of the $c\times c$ matrix with the entries
\begin{equation}\label{eq.PVV}
   \pi_k(z_i)=x_k^{\ind{\ve_k=1}}\,V_k(z_i)-x_k^{\ind{\ve_k=-1}}\,V'_k(z_i),
\end{equation}
and $\tP_{D,\hate}(y)$ is, up to sign, the determinant of the matrix with the entries
\[
   \tpi_k(z_i)=y_k^{\ind{\ve_k=-1}}\,V^\uparrow_k(z_i)-y_k^{\ind{\ve_k=1}}\,V'^{\uparrow}_k(z_i) .
\]
Compare with the classical matrix: by (i), the entry of $tV_F-V_F^{\T}$ at two cycles $z,z'$ is
\[
   t\,V_F\bigl([z],[z']\bigr)-V_F\bigl([z'],[z]\bigr)=\sum_{k=1}^c\bigl(t\,V_k(z)-V'_k(z)\bigr)\,w_k(z') .
\]
The deformation consists in \emph{not} summing over the bands $k$ through which the second cycle passes, so that there is one column for each band instead of one column for each cycle, and in replacing $t$ by the variable of the band: by $x_k$ if the crossing is positive, and by $x_k^{-1}$, after multiplication of the column by $x_k$, if it is negative. In this sense the crossing variable $x_k$ is the variable $t$ of the Alexander polynomial, localized at the band of the crossing $k$. In the matrix \eqref{eq.PVV} the rows of the Seifert circles $C\ne\hat C$ of a Seifert basis are the rows of the deformed incidence matrix $\CC(x)$; they will now be eliminated.

\subsection{Regions and levels}\label{ssec.regions}
\begin{figure}[ht]
\centering
\begin{tikzpicture}[>=Stealth,scale=1.15,every node/.style={font=\small},
  cw/.style={postaction={decorate},decoration={markings,mark=at position 0.12 with {\arrow{<}},mark=at position 0.62 with {\arrow{<}}}},
  ccw/.style={postaction={decorate},decoration={markings,mark=at position 0.12 with {\arrow{>}},mark=at position 0.62 with {\arrow{>}}}}]
\fill[gray!12] (0,0) circle (3.2);
\fill[white] (-1.3,0) circle (1.25);
\fill[white] (1.8,0.1) circle (0.95);
\fill[gray!30] (-1.3,0) circle (1.25);
\fill[white] (-1.3,0.1) circle (0.5);
\draw[line width=6pt,gray!60] (-1.3,1.25) -- (-1.3,3.2);
\draw[line width=6pt,gray!60] (-2.45,-0.45) -- (-3.1,-0.75);
\draw[line width=6pt,gray!60] (-0.05,0.05) -- (0.85,0.1);
\draw[line width=6pt,gray!60] (-1.8,0.1) -- (-2.55,0.1);
\draw[line width=6pt,gray!60] (-0.8,0.1) -- (-0.05,0.1);
\draw[cw,thick] (0,0) circle (3.2);
\draw[cw,thick] (-1.3,0) circle (1.25);
\draw[ccw,thick] (1.8,0.1) circle (0.95);
\draw[cw,thick] (-1.3,0.1) circle (0.5);
\node at (0,-2.7) {$C_0$, $h=0$};
\node at (-1.3,-0.85) {$C_1$, $h=1$};
\node at (1.8,0.1) {$C_2$, $h=0$};
\node[font=\scriptsize] at (-1.3,0.1) {$C_3$};
\node[font=\scriptsize] at (-1.3,0.85) {$\Omega_2$: $1,2$};
\node at (0.6,2.2) {$\Omega_1$: levels $0,1$};
\node at (0,3.6) {unbounded region: levels $-1,0$};
\end{tikzpicture}
\caption{Four Seifert circles, oriented as indicated ($s=1$ for the clockwise circles $C_0,C_1,C_3$ and $s=-1$ for $C_2$), with levels $h(C_0)=0$, $h(C_1)=1$, $h(C_2)=0$, $h(C_3)=2$, and the bands of the canonical surface (gray). The circles bounding the region $\Omega_1$ (light gray) have levels $0$ and $1$, those bounding $\Omega_2$ (darker) have levels $1$ and $2$: every band joins two circles of adjacent levels bounding the same region, and every cycle of the Seifert graph stays in one region (Lemma~\ref{lem.regions}).}
\label{fig.regions}
\end{figure}
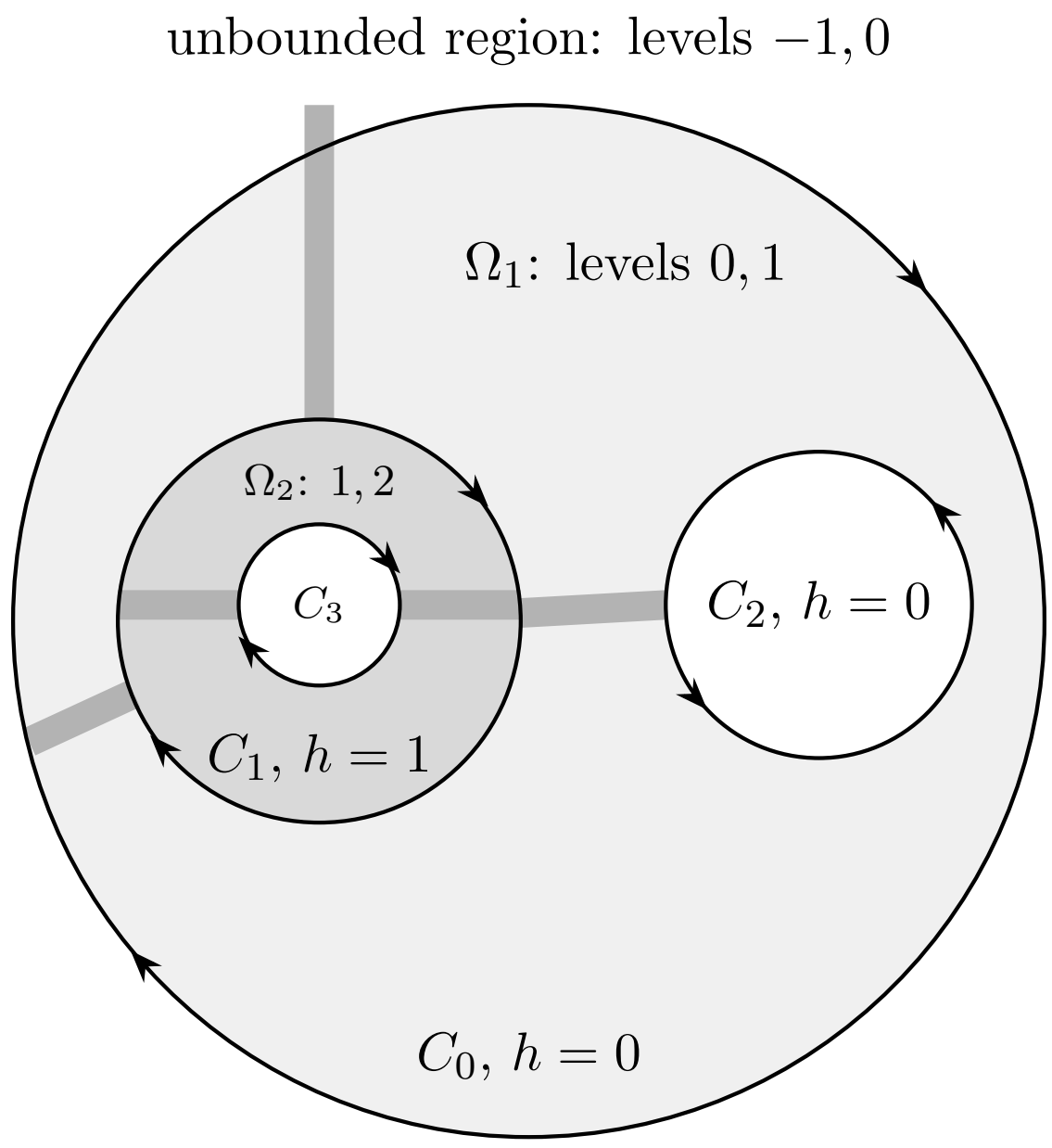

We shall need to know how the cycles of the Seifert graph are placed with respect to the levels of Section~\ref{ssec.level}. The answer is that the Seifert circles of a cycle have only two levels, and that these two levels are adjacent (Lemma~\ref{lem.regions}). We use it twice: in Section~\ref{ssec.alexpoint}, to compare the deformed Seifert matrices at the Alexander point with the classical matrix, and in Section~\ref{sec.degree}, to compute the degrees of their entries. We draw $D$ in the plane and use the levels $h(C)\in\BZ$ of the Seifert circles of Definition~\ref{def.level}, which satisfy $h(B_j)=h(A_j)+\ve_j$ at every crossing (Lemma~\ref{lem.level}). The complement in the plane of the union of the Seifert circles is a disjoint union of open \emph{regions}. A bounded region $\Omega$ has an outer boundary circle $C_{\rm out}$, the smallest Seifert circle whose interior contains $\Omega$, and inner boundary circles, the maximal Seifert circles contained in the interior of $C_{\rm out}$ and disjoint from $\Omega$; the unbounded region is bounded by the outermost Seifert circles. At a crossing $j$, the two Seifert circles $A_j$ and $B_j$ bound the region $\Omega$ which contains the crossing itself (the region of the proof of Lemma~\ref{lem.level}): we say that the crossing, or the band of $F_D$ at this crossing, \emph{lies in $\Omega$} (Figure~\ref{fig.regions}).

\begin{lemma}\label{lem.regions}
\begin{enumerate}\renewcommand{\theenumi}{\roman{enumi}}
\item For every region $\Omega$ there is an integer $H_\Omega$ such that every Seifert circle bounding $\Omega$ has level $H_\Omega$ or $H_\Omega+1$.
\item Every cycle of the Seifert graph $\Gamma_D$ consists of crossings lying in one and the same region $\Omega$. Hence the Seifert circles of a cycle $\gamma$ of $\Gamma_D$ have two adjacent levels $H_\gamma$ and $H_\gamma+1$, where $H_\gamma=H_\Omega$.
\end{enumerate}
\end{lemma}
\begin{proof}
(i) Let $C_{\rm in}$ be an inner boundary circle of a bounded region with outer circle $C_{\rm out}$. The Seifert circles whose interior contains $C_{\rm in}$ are $C_{\rm out}$ and those whose interior contains $C_{\rm out}$, so by Definition~\ref{def.level}
\[
   h(C_{\rm in})=h(C_{\rm out})-\tfrac{s(C_{\rm out})-1}2+s(C_{\rm out})+\tfrac{s(C_{\rm in})-1}2=h(C_{\rm out})+\tfrac{s(C_{\rm out})+s(C_{\rm in})}2 ,
\]
which is $h(C_{\rm out})$ or $h(C_{\rm out})+s(C_{\rm out})$. With $H_\Omega=h(C_{\rm out})+\frac{s(C_{\rm out})-1}2$, both $h(C_{\rm out})=H_\Omega+\frac{1-s(C_{\rm out})}2$ and $h(C_{\rm in})$ lie in $\{H_\Omega,H_\Omega+1\}$. For the unbounded region, the outermost circles have levels $\frac{s(C)-1}2\in\{-1,0\}$.

(ii) Every Seifert circle $C$ bounds exactly two regions, the one directly inside it and the one directly outside it, and separates the plane: a crossing lying in a region inside $C$ joins two circles which are $C$ or lie inside $C$, a crossing lying in a region outside $C$ joins two circles which are $C$ or lie outside $C$. Thus the regions, ordered by the nesting of the circles, form a tree, $\Gamma_D$ is the union of the subgraphs $\Gamma_\Omega$ (the circles bounding $\Omega$ and the crossings lying in $\Omega$), and two of these subgraphs meet in at most one vertex. Let $\gamma$ be a cycle of $\Gamma_D$ (a closed path without repeated vertices) and suppose that it contains crossings lying in two different regions $\Omega\ne\Omega'$. Some circle $C$ separates $\Omega$ from $\Omega'$; the two crossings lie on different sides of $C$, so the cycle, which is connected, passes through the vertex $C$ at least twice---once to go from one side to the other and once to come back---which is impossible. So all crossings of $\gamma$ lie in one region $\Omega$, all its circles bound $\Omega$, and (i) applies.
\end{proof}

\subsection{The deformed cycle space and the deformed Seifert matrices}\label{ssec.deformedmatrices}
The matrix of the crossing forms in a Seifert basis has size $c$, whereas a Seifert matrix of the canonical surface has size $\nu=c-s+1$, the rank of $H_1(F_D)$. In this subsection we eliminate the $s-1$ rows of the Seifert circles. What remains is a matrix of size $\nu$, the deformed Seifert matrix, and the determinant only changes by an explicit monomial. The price of the elimination is that the cycles of the Seifert graph, along which the classical Seifert matrix adds up the contributions of the bands, are replaced by \emph{deformed cycles}: the vectors of the kernel of the rows which are eliminated.

Let $Y$ be a basis of the cycle space $H$ of the Seifert graph, with columns $\gamma_1,\dots,\gamma_\nu$, and let $N_{D,\hate}(x)$ be the matrix \eqref{eq.N} of the crossing forms in the Seifert basis $z^{\gamma_1},\dots,z^{\gamma_\nu}$, $C\ne\hat C$. Its last $s-1$ rows form the matrix $\CC(x)$ of \eqref{eq.CC}, whose row $C$ is $\bigl(\pi_k(C)\bigr)_k$. The kernel of $\CC(x)$ is a deformation of the cycle space, and the Schur complement of $N_{D,\hate}(x)$ with respect to the rows of the circles is a square matrix of size $\nu$.

\begin{corollary}\label{cor.seifert}
Let $\CT$ be a spanning tree of $\Gamma_D$, rooted at $\hat C$, and for every circle $C\ne\hat C$ let $j(C)$ be the edge of $\CT$ joining $C$ to its parent in the tree.
\begin{enumerate}\renewcommand{\theenumi}{\roman{enumi}}
\item $P_{D,\hate}(x)=\kappa\,\det N_{D,\hate}(x)$, where $N_{D,\hate}(x)=\bigl(\begin{smallmatrix}Y^{\T}K(x)\\ \CC(x)\end{smallmatrix}\bigr)$ is the matrix of the crossing forms in the Seifert basis and
\[
   \kappa=(-1)^{s-1}\det(Y\,|\,E_{\CT})\,\det(\hat\partial E_{\CT})=\pm1 .
\]
\item The $(s-1)\times(s-1)$ matrix $\CC(x)E_{\CT}$ is invertible over $\Lam_c$, with determinant $\pm x^{\CT}$, where
\[
   x^{\CT}=\prod_{C\ne\hat C,\ A_{j(C)}=C}x_{j(C)} .
\]
The $c\times\nu$ matrix
\begin{equation}\label{eq.Yx}
   Y(x)=Y-E_{\CT}\bigl(\CC(x)E_{\CT}\bigr)^{-1}\CC(x)Y
\end{equation}
has entries in $\Lam_c$, satisfies $Y(1)=Y$, and its columns form a basis of $\ker\CC(x)\subseteq\Lam_c^c$, the \emph{deformed cycle space}.
\item $P_{D,\hate}(x)=x^{\CT}\det\bigl(Y^{\T}K(x)Y(x)\bigr)$, a $\nu\times\nu$ determinant with no sign in front.
\end{enumerate}
\end{corollary}
\begin{proof}
(i) This is \eqref{eq.kappaSeifert}.

(ii) The spanning tree serves to put the matrix $\CC(x)E_{\CT}$ in triangular form. Order the circles $C\ne\hat C$ so that every circle comes before its parent. The column $j(C)$ of $\CC(x)E_{\CT}$ has, by \eqref{eq.CC}, the entry $x_{j(C)}$ in the row $A_{j(C)}$ and the entry $-1$ in the row $B_{j(C)}$, when these rows exist; these two circles are $C$ and its parent, and the row of the parent is missing if the parent is $\hat C$. With this ordering of the rows and the columns, $\CC(x)E_{\CT}$ is therefore lower triangular, with diagonal entry $x_{j(C)}$ if $A_{j(C)}=C$ and $-1$ if $B_{j(C)}=C$. Hence $\det(\CC(x)E_{\CT})=\pm x^{\CT}$, a unit of $\Lam_c$, and $Y(x)$ has entries in $\Lam_c$. With the same ordering the matrix $\hat\partial E_{\CT}$ is lower triangular as well, with diagonal entry $-1$ if $A_{j(C)}=C$ and $1$ if $B_{j(C)}=C$, so that the product of the two diagonal entries of a circle $C$ is $-x_{j(C)}$ or $-1$, and
\begin{equation}\label{eq.signtree}
   \det(\hat\partial E_{\CT})\,\det\bigl(\CC(x)E_{\CT}\bigr)=(-1)^{s-1}\,x^{\CT}
\end{equation}
for every ordering of the circles and of the edges of $\CT$, the same in the two matrices. We have $\CC(x)Y(x)=\CC(x)Y-\CC(x)Y=0$, and $Y(1)=Y$ because $\CC(1)Y=-\hat\partial Y=0$. If $w\in\Lam_c^c$ satisfies $\CC(x)w=0$, write $w=Y\eta+E_{\CT}\zeta$ with $\eta\in\Lam_c^{\nu}$ and $\zeta\in\Lam_c^{s-1}$, which is possible since $(Y\,|\,E_{\CT})$ is unimodular; then $\CC(x)Y\eta+\CC(x)E_{\CT}\zeta=0$ gives $\zeta=-(\CC(x)E_{\CT})^{-1}\CC(x)Y\eta$ and $w=Y(x)\eta$. The rows of $Y(x)$ indexed by the edges $j\notin\CT$ are those of $Y$, which form the identity matrix for the basis of fundamental cycles, and in general a unimodular matrix (since $(Y\,|\,E_{\CT})$ is unimodular); so $Y(x)\eta=0$ implies $\eta=0$, and the columns of $Y(x)$ form a basis of $\ker\CC(x)$.

(iii) We now eliminate the rows of the circles, by passing to the Schur complement of the block $\CC(x)E_{\CT}$. Multiplying the matrix $N_{D,\hate}(x)$ on the right by the unimodular matrix $(Y\,|\,E_{\CT})$ gives
\[
   \begin{pmatrix}Y^{\T}K(x)Y&Y^{\T}K(x)E_{\CT}\\ \CC(x)Y&\CC(x)E_{\CT}\end{pmatrix},
\]
whose lower right block is invertible by (ii). Its determinant is $\det(\CC(x)E_{\CT})$ times the determinant of the Schur complement
\[
   Y^{\T}K(x)Y-Y^{\T}K(x)E_{\CT}\bigl(\CC(x)E_{\CT}\bigr)^{-1}\CC(x)Y=Y^{\T}K(x)Y(x) .
\]
Thus $\det N_{D,\hate}(x)\det(Y\,|\,E_{\CT})=\det(\CC(x)E_{\CT})\det\bigl(Y^{\T}K(x)Y(x)\bigr)$, and by (i) and \eqref{eq.signtree}
\[
   P_{D,\hate}(x)=(-1)^{s-1}\det(\hat\partial E_{\CT})\det\bigl(\CC(x)E_{\CT}\bigr)\det\bigl(Y^{\T}K(x)Y(x)\bigr)=x^{\CT}\det\bigl(Y^{\T}K(x)Y(x)\bigr).
\]
\end{proof}

The second flavor has the same structure, with the deformation on the other side. Let $Y(y)$ be the matrix \eqref{eq.Yx} with $x$ replaced by $y$, and recall the matrices $K^b(y)$ of \eqref{eq.Kbrow} and $\tilde N_{D,\hate}(y)$ of \eqref{eq.N}.

\begin{corollary}\label{cor.seifertb}
With $\CT$ and $\kappa$ as in Corollary~\ref{cor.seifert},
\[
   \tP_{D,\hate}(y)=(-1)^{\mu-1}\kappa\,\det\tilde N_{D,\hate}(y)\qquad\text{and}\qquad\tP_{D,\hate}(y)=y^{\CT}\det\bigl(Y(y)^{\T}K^b(y)Y\bigr),
\]
the second formula without any sign, where $y^{\CT}$ is the monomial $x^{\CT}$ in the variables $y$.
\end{corollary}
\begin{proof}
The first formula was proved in the proof of Theorem~\ref{thm.flows}, in the form
\begin{equation}\label{eq.signNb}
   \det(Y\,|\,E_{\CT})\cdot\tP_{D,\hate}(y)=(-1)^{s-1}\det(\hat\partial E_{\CT})\,\det\bigl(K^b(y)Y\ \big|\ \CC(y)^{\T}\bigr).
\end{equation}
For the second one, multiply the $c\times c$ matrix $\bigl(K^b(y)Y\,|\,\CC(y)^{\T}\bigr)$ on the left by the matrix $W$ with the rows of $Y(y)^{\T}$ followed by those of $E_{\CT}^{\T}$. By \eqref{eq.Yx}, $W$ is the product of a unipotent matrix and of the transpose of $(Y\,|\,E_{\CT})$, so $\det W=\det(Y\,|\,E_{\CT})$. Since $Y(y)^{\T}\CC(y)^{\T}=(\CC(y)Y(y))^{\T}=0$, the result is
\[
   \begin{pmatrix}Y(y)^{\T}K^b(y)Y&0\\ E_{\CT}^{\T}K^b(y)Y&(\CC(y)E_{\CT})^{\T}\end{pmatrix},
\]
so that $\det(Y\,|\,E_{\CT})\det\bigl(K^b(y)Y\,|\,\CC(y)^{\T}\bigr)=\det\bigl(Y(y)^{\T}K^b(y)Y\bigr)\det\bigl(\CC(y)E_{\CT}\bigr)$. With \eqref{eq.signNb} and \eqref{eq.signtree}, and since $\det(Y\,|\,E_{\CT})^2=1$,
\[
   \tP_{D,\hate}(y)=(-1)^{s-1}\det(\hat\partial E_{\CT})\det\bigl(\CC(y)E_{\CT}\bigr)\det\bigl(Y(y)^{\T}K^b(y)Y\bigr)=y^{\CT}\det\bigl(Y(y)^{\T}K^b(y)Y\bigr).
\]
\end{proof}

\begin{definition}[deformed Seifert matrices]\label{def.deformedseifertmatrices}
The \emph{deformed Seifert matrices} of the diagram $D$, for the cut edge $\hate$, the basis $Y$ of the cycle space and the spanning tree $\CT$, are the $\nu\times\nu$ matrices
\begin{equation}\label{eq.MD}
   M_{D,\hate}(x)=Y^{\T}K(x)Y(x)\quad\text{and}\quad M^b_{D,\hate}(y)=Y(y)^{\T}K^b(y)Y ,
\end{equation}
with entries in $\Lam_c$ and in $\Lam'_c$. By Lemma~\ref{lem.rowsareflows}, their entries are
\begin{equation}\label{eq.Mentries}
   M_{D,\hate}(x)_{a,b}=\sum_{k=1}^c\pi_k\bigl(z^{\gamma_a}\bigr)\,Y(x)_{k,b},\qquad M^b_{D,\hate}(y)_{a,b}=-\sum_{k=1}^cY(y)_{k,a}\,\tpi_k\bigl(z^{\gamma_b}\bigr).
\end{equation}
\end{definition}

In words: the entry $(a,b)$ of $M_{D,\hate}(x)$ pairs the cycle $\gamma_a$, through the crossing forms of its lift, with the $b$-th deformed cycle. In the entry $\sum_k\bigl(t\,V_k-V'_k\bigr)(z^{\gamma_a})\,w_k(\gamma_b)$ of the classical matrix $tV_F-V_F^{\T}$ (Section~\ref{ssec.bandbyband}), the second factor is the cycle $\gamma_b$ itself. The two matrices are compared in Theorem~\ref{thm.alexpoint}.

By Corollaries~\ref{cor.seifert}~(iii) and~\ref{cor.seifertb},
\begin{equation}\label{eq.PM}
   P_{D,\hate}(x)=x^{\CT}\det M_{D,\hate}(x)\qquad\text{and}\qquad\tP_{D,\hate}(y)=y^{\CT}\det M^b_{D,\hate}(y),
\end{equation}
without any sign; in particular the right-hand sides are polynomials of degree $\le1$ in each variable. So Theorem~\ref{thm.HK} takes the following form.

\begin{corollary}[Seifert-matrix form of the determinant formula]\label{cor.seifertformula}
For every diagram $D$ of a link $L$ with $\mu$ components, every edge $\hate$ of $D$, every basis $Y$ of the cycle space of $\Gamma_D$ and every spanning tree $\CT$ of $\Gamma_D$,
\[
   \LG(L;t_0,t_1)=(-1)^{\mu-1}\ \BE_D\Bigl[x^{\CT}y^{\CT}\det\bigl(M_{D,\hate}(x)\bigr)\det\bigl(M^b_{D,\hate}(y)\bigr)\Bigr].
\]
\end{corollary}

\subsection{The deformed Seifert form, and why the variables are cochains}\label{ssec.cochains}
This subsection is about the meaning of the deformed Seifert matrices. It answers two questions. Of which object, defined without any choice of bases, is $M_{D,\hate}(x)$ the matrix? And what kind of object on the surface $F_D$ is the family of the crossing variables?

The matrix $M_{D,\hate}(x)$ is the matrix of a bilinear pairing which does not depend on any basis. By Lemma~\ref{lem.rowsareflows}, $\CC(x)w$ is the vector with the coordinates $\sum_k\pi_k(C)\,w_k$, $C\ne\hat C$; so the deformed cycle space is
\[
   H_x=\ker\CC(x)=\Bigl\{w\in\Lam_c^c:\ \sum_{k=1}^c\pi_k(C)\,w_k=0\ \text{for every Seifert circle }C\ne\hat C\Bigr\},
\]
a free $\Lam_c$-module of rank $\nu$ by Corollary~\ref{cor.seifert}~(ii). Since $\FF_D=\FF_{D,\hate}\oplus\BZ\hat C$, the group $H_1(F_D)=\FF_D/\langle C_1,\dots,C_s\rangle$ is also the quotient of $\FF_{D,\hate}$ by the Seifert circles $C\ne\hat C$. Hence the formula
\begin{equation}\label{eq.Theta}
   \Theta_x\bigl([z],w\bigr)=\sum_{k=1}^c\pi_k(z)\,w_k\qquad(z\in\FF_{D,\hate},\ w\in H_x)
\end{equation}
defines a pairing $\Theta_x\colon H_1(F_D;\BZ)\times H_x\to\Lam_c$, additive in the first variable and $\Lam_c$-linear in the second one, which we call the \emph{deformed Seifert form} of the diagram cut at $\hate$. By \eqref{eq.Mentries}, $M_{D,\hate}(x)$ is its matrix in the basis $\gamma_1,\dots,\gamma_\nu$ of $H_1(F_D)$ and the basis of $H_x$ formed by the columns of $Y(x)$. At $x=1$ the deformed cycle space is $H\otimes\Lam_c$, the cycle space itself, because $\CC(1)=-\hat\partial$ and $\ker\hat\partial=\ker\partial$ (the sum of the rows of $\partial$ is zero); and by Theorem~\ref{thm.deformedseifert}
\[
   \Theta_1\bigl([z],[z']\bigr)=\sum_k\bigl(V_k(z)-V'_k(z)\bigr)w_k(z')=V_F\bigl([z],[z']\bigr)-V_F\bigl([z'],[z]\bigr),
\]
the intersection form of $F_D$, up to the sign convention. For general $x$,
\[
   \Theta_x\bigl([z],w\bigr)=\sum_k\Bigl(x_k^{\ind{\ve_k=1}}\,V_k(z)-x_k^{\ind{\ve_k=-1}}\,V'_k(z)\Bigr)w_k :
\]
the Seifert form and its transpose, with the variable of each band inserted. In the same way the second flavor defines a pairing $\tilde\Theta_y(w,[z])=-\sum_kw_k\,\tpi_k(z)$ on $H_y\times H_1(F^\uparrow_D;\BZ)$, whose matrix is $M^b_{D,\hate}(y)$.

\emph{Cochains.} What kind of object is the family $x=(x_1,\dots,x_c)$ of the crossing variables, from the point of view of the surface $F_D$? The crossings are the bands, and the bands correspond to their co-cores $\alpha_1,\dots,\alpha_c$. In the cell structure of $F_D$ of Section~\ref{ssec.statehomology}, the $1$-cells which do not lie on the boundary $L$ are the co-cores, and all the $0$-cells lie on $L$. So the cellular chain complex of the pair $(F_D,\partial F_D)$ is
\[
   C_2(F_D,\partial F_D)=\BZ^s\ \xrightarrow{\ \partial^{\T}\ }\ C_1(F_D,\partial F_D)=\bigoplus_{j=1}^c\BZ\,\alpha_j=\BZ^c\ \longrightarrow\ C_0(F_D,\partial F_D)=0 ,
\]
where the boundary of the $2$-cell of a Seifert circle $C$ is $\sum_j\bigl(\ind{B_j=C}-\ind{A_j=C}\bigr)\alpha_j$ by \eqref{eq.flowcycle}: this is the transpose of the incidence map $\partial$ of the Seifert graph. The group $\BZ^c$ in which the cycle space $H$ lives is therefore the group of the relative $1$-chains of the surface, and also, by means of the basis of the co-cores, the group $C^1(F_D,\partial F_D;\BZ)$ of the relative $1$-\emph{cochains}. A monomial $x^n=\prod_jx_j^{n_j}$ is the cochain which takes the value $n_j$ on $\alpha_j$, and
\[
   \Lam_c=\BZ\bigl[C^1(F_D,\partial F_D;\BZ)\bigr]
\]
is the group ring of the group of relative cochains. A point $x\in(\BC^\times)^c$ at which the polynomials are evaluated is a multiplicative cochain: it assigns the number $x_j$ to the co-core $\alpha_j$.

Since there are no relative $0$-cells, there are no coboundaries, and the cohomology group $H^1(F_D,\partial F_D;\BZ)$ is the subgroup of the cocycles: the cochains $n$ which vanish on the boundaries of the $2$-cells, that is, with $\partial n=0$. This is the cycle space $H$. (Lefschetz duality identifies $H^1(F_D,\partial F_D)$ with $H_1(F_D)$, and in these terms it is the identity of $H$: a cycle $\gamma$ is the cocycle $\alpha_j\mapsto w_j(\gamma)$, which counts the passages of $\gamma$ through the co-cores.) So the classes form a subgroup of rank $\nu=c-s+1$ of the group of cochains, of rank $c$. One might hope that the determinant formula only involves the classes, that is, that the edge polynomials lie, up to a monomial, in the group ring $\BZ[H^1(F_D,\partial F_D;\BZ)]=\BZ[H]$ of a group which depends on the surface alone, or that their values depend only on the class of $x$ in a suitable sense. This is not the case, for three reasons.

\emph{(a) The polynomials are not combinations of cycle monomials.} For the trefoil of Example~\ref{ex.trefoilflows}, $H=\{n\in\BZ^3:n_1+n_2+n_3=0\}$ and $P_{D,\hate}(x)=1-x_1+x_1x_2$. The exponents $(0,0,0)$, $(1,0,0)$ and $(1,1,0)$ of its monomials lie in three different cosets of $H$ in $\BZ^3$, so that $P_{D,\hate}$ is not the product of a monomial and of an element of $\BZ[H]$. The $c$ variables are independent, and the coupling $\BE_D$ of the two flavors uses each of them separately: it replaces $x_jy_j$ by $1$ crossing by crossing.

\emph{(b) The Alexander point is not a cocycle.} A multiplicative cochain $x$ is a cocycle if it takes the value $1$ on the boundary of every $2$-cell, that is, if $\prod_{A_j=C}x_j=\prod_{B_j=C}x_j$ for every Seifert circle $C$. Let $x=t^{\ve}$ be the Alexander point, with $t$ not a root of unity, and let $C$ be a Seifert circle of maximal level. By Lemma~\ref{lem.level}, $h(B_j)=h(A_j)+\ve_j$; so $\ve_j=-1$ for the crossings with $A_j=C$ and $\ve_j=1$ for those with $B_j=C$, and the condition for $C$ reads $t^{-n}=1$, where $n\ge1$ is the number of crossings on $C$. So the point at which the determinant formula gives the Alexander polynomial is a cocycle for no diagram.

\emph{(c) The Alexander point is trivial as a local system.} There is a dual way of forgetting part of the cochain. A function on the co-cores is the same thing as a local system of rank one on $F_D$ which is trivialized on each of the disks $D_C$, with the holonomy $x_j$ through the band $j$ from $A_j$ to $B_j$. Forgetting the trivializations means dividing by the gauge group $(\BC^\times)^s$ of the disks, which acts by $x_j\mapsto\lambda_{B_j}x_j\lambda_{A_j}^{-1}$; the quotient is the torus $H^1(F_D;\BC^\times)$ of the local systems on $F_D$, again of dimension $\nu$. But $\ve_j=h(B_j)-h(A_j)$, so that the Alexander point $x=t^{\ve}$ is the gauge transform of the trivial cochain $x=1$ by $\lambda_C=t^{h(C)}$: as a local system on the surface it is trivial. This is the classical fact that the matrix $tV-V^{\T}$ involves no local system on the Seifert surface. And the edge polynomials are not invariant under the gauge group, even up to units: for the trefoil it acts by $x\mapsto\lambda x$, and $1-\lambda x_1+\lambda^2x_1x_2$ is not a multiple of $1-x_1+x_1x_2$.

So neither the subgroup of the cocycles nor the quotient by the gauge group carries the crossing variables: the determinant formula is a function of the cochain itself. Theorem~\ref{thm.deformedseifert} shows why. The variable $t$ of the Alexander polynomial does not live on the Seifert surface; it is the generator of the group of deck transformations of the infinite cyclic covering of the link exterior, which is dual to the surface, and in the matrix $tV-V^{\T}$ it distinguishes the two sides of the surface, the push-offs $b^+$ and $b^-$. The crossing variable $x_j$ plays this role at the band $j$ only: it multiplies the linking numbers with the push-off of the band $j$ to one of its two sides. So the $c$ bands carry independent variables, whereas for a class the values of the cochain on the bands of a Seifert circle are tied together: they satisfy the cocycle condition of (b), or they are only defined up to the gauge transformations of (c).

The variable $t$ itself is not seen by the classes. Indeed, $\ve_j=h(B_j)-h(A_j)$ by Lemma~\ref{lem.level}: the vector $\ve=\partial^{\T}h$ is the coboundary of the level function $h$ on the Seifert circles, which is the fact behind (b) and (c). Hence $\sum_j\ve_jn_j=\sum_Ch(C)\,(\partial n)_C$ for every $n\in\BZ^c$. So the Alexander point sends the monomial $x^n$ to a power of $t$ which depends only on $\partial n$, that is, on the coset of $n$ modulo the cycle space $H=\ker\partial$, and it sends every element of the group ring $\BZ[H]$ of the classes to an integer: a polynomial in the classes alone would not contain the variable $t$.

With both flavors there are two variables at each crossing: $x_j$ belongs to the relation of the under-strand, at the lower end of the co-core $\alpha_j$, and $y_j$ to that of the over-strand, at its upper end (Definition~\ref{def.crossingforms}). In Section~\ref{sec.cocores} we shall see that the values $t_0^{\ve_j}$ and $t_1^{\ve_j}$ which the coupling $\BE_D$ gives to the crossing variables are nevertheless holonomies of a local system of rank one, which lives on the complement in $\BR^3$ of the link and of the co-cores, and not on the surface.

\begin{remark}[twisted homology]\label{rem.twisted}
The rows of $\CC(x)$ are, up to sign, the boundaries of the $2$-cells of $F_D$ in the relative chain complex with coefficients twisted by the cochain $x$, the cell of $\hat C$ being left out; the quotient $\Lam_c^c/\langle\CC_C:C\ne\hat C\rangle$ of the relative $1$-chains by these boundaries is a twisted relative homology group of the surface cut at $\hate$, free of rank $\nu$, and $H_x$ is its dual. In these terms the deformed Seifert form is the homomorphism $H_1(F_D;\BZ)\otimes\Lam_c\to\Lam_c^c/\langle\CC_C\rangle$ induced by $\Pi(x)$, which at $x=1$ is the natural map $H_1(F_D)\to H_1(F_D,\partial F_D)$, the adjoint of the intersection form. Its determinant, with respect to a basis of $H_1(F_D)$ and to the volume form of the target defined by the co-cores and the full turns, is $\pm P_{D,\hate}(x)$; the classes of the co-cores of the bands outside a spanning tree $\CT$ form a basis of the target of volume $\det(\CC(x)E_{\CT})=\pm x^{\CT}$, which is where the monomial $x^{\CT}$ of Corollary~\ref{cor.seifert}~(iii) comes from.
\end{remark}

\subsection{The Alexander point}\label{ssec.alexpoint}
In this subsection we evaluate the deformed Seifert matrices at the Alexander point $x=t^{\ve}$, and we find the classical matrix $tV_F-V_F^{\T}$ of the canonical surface, with units on its columns (Theorem~\ref{thm.alexpoint}). There are two steps: we first determine the deformed cycles at this point (Lemma~\ref{lem.Yt}), and we then compute the entries of the matrix with Theorem~\ref{thm.deformedseifert}.

At the Alexander point the deformed cycles are the cycles themselves, rescaled.

\begin{lemma}\label{lem.Yt}
Let $D$ be a connected diagram, $\CT$ a spanning tree of $\Gamma_D$ and $Y$ the basis of fundamental cycles of $\CT$. Then
\[
   Y(t^{\ve})=\mathsf Z_t^{-1}\,Y\,\mathsf Z'_t,\qquad \mathsf Z_t=\diag\bigl(t^{-\ind{\ve_k=-1}}\bigr)_{k=1,\dots,c},\quad \mathsf Z'_t=\diag\bigl(t^{-\ind{\ve_i=-1}}\bigr)_{i\notin\CT}.
\]
\end{lemma}
\begin{proof}
The matrix $\mathsf Z_t^{-1}Y\mathsf Z'_t$ has entries in $\BZ[t^{\pm1}]$ and its rows at the crossings $i\notin\CT$ form the identity matrix, like those of $Y$ and of $Y(t^{\ve})$; by Corollary~\ref{cor.seifert} (ii) and its proof, which applies verbatim at $x=t^{\ve}$ since $\det(\CC(t^{\ve})E_{\CT})$ is a unit, a vector of $\ker\CC(t^{\ve})$ is determined by its coordinates at these crossings, so it suffices to show that $\CC(t^{\ve})\mathsf Z_t^{-1}Y=0$. The entry $(C,k)$ of $\CC(t^{\ve})\mathsf Z_t^{-1}$ is $t^{\ve_k+\ind{\ve_k=-1}}\ind{A_k=C}-t^{\ind{\ve_k=-1}}\ind{B_k=C}=t^{\ind{\ve_k=1}}\ind{A_k=C}-t^{\ind{\ve_k=-1}}\ind{B_k=C}$. Let $i\notin\CT$ and let $k$ be a band of the fundamental cycle $\gamma_i$ lying on the circle $C$. By Lemma~\ref{lem.regions}~(ii), the circles of $\gamma_i$ have two adjacent levels $H_i$ and $H_i+1$, and $h(B_k)=h(A_k)+\ve_k$ (Lemma~\ref{lem.level}). If $h(C)=H_i$, then $\ve_k=1$ when $A_k=C$ and $\ve_k=-1$ when $B_k=C$, so the exponent of $t$ in the entry $(C,k)$ is $1$ in both cases; if $h(C)=H_i+1$, it is $0$ in both cases. Hence $(\CC(t^{\ve})\mathsf Z_t^{-1}Y)_{C,i}=t^{\ind{h(C)=H_i}}\sum_k\bigl(\ind{A_k=C}-\ind{B_k=C}\bigr)Y_{ki}=-t^{\ind{h(C)=H_i}}(\hat\partial Y)_{C,i}=0$.
\end{proof}

Lemma~\ref{lem.Yx} below is the statement about the $\ve$-degrees of the entries of $Y(x)$ which underlies this: for two bands $i,k$ of the same fundamental cycle, $h(B_i)-h(B_k)=\ind{\ve_k=-1}-\ind{\ve_i=-1}$.

\begin{theorem}[the Alexander point]\label{thm.alexpoint}
Let $D$ be a diagram, $\hate$ an edge of $D$, $\CT$ a spanning tree of $\Gamma_D$ and $Y$ the basis of fundamental cycles of $\CT$. Let $V_F$ and $V_{F^\uparrow}$ be the Seifert matrices of the canonical Seifert surfaces $F_D$ and $F^\uparrow_D$ in this basis. Then
\[
   M_{D,\hate}(t^{\ve})=\bigl(tV_F-V_F^{\T}\bigr)\,\mathsf Z'_t\qquad\text{and}\qquad M^b_{D,\hate}(t^{\ve})=\mathsf Z'_t\,\bigl(tV_{F^\uparrow}-V_{F^\uparrow}^{\T}\bigr),
\]
with the diagonal matrix $\mathsf Z'_t$ of Lemma~\ref{lem.Yt}. Consequently
\[
   P_{D,\hate}(t^{\ve})=t^{\ve(\CT')-\nu_-}\det\bigl(tV_F-V_F^{\T}\bigr),\qquad \tP_{D,\hate}(t^{\ve})=t^{\ve(\CT')-\nu_-}\det\bigl(tV_{F^\uparrow}-V_{F^\uparrow}^{\T}\bigr),
\]
where $\nu_-$ is the number of negative crossings not in $\CT$ and $\ve(\CT')$ is the exponent of $t$ in the value $t^{\ve(\CT')}$ of the monomial $x^{\CT}$ of Corollary~\ref{cor.seifert} at the Alexander point, as in Theorem~\ref{thm.degrees}. In particular $M_{D,\hate}(1)=V_F-V_F^{\T}$.
\end{theorem}
\begin{proof}
Let $i_b\notin\CT$ be the crossing of the fundamental cycle $\gamma_b$. By Lemma~\ref{lem.Yt}, $Y(t^{\ve})_{k,b}=t^{\ind{\ve_k=-1}}\,Y_{k,b}\,t^{-\ind{\ve_{i_b}=-1}}$, and $Y_{k,b}=w_k(\gamma_b)$. By \eqref{eq.Mentries} and Theorem~\ref{thm.deformedseifert}~(iii) and~(i),
\begin{align*}
   M_{D,\hate}(t^{\ve})_{a,b}&=\sum_kt^{-\ind{\ve_k=-1}}\bigl(tV_k-V'_k\bigr)(z^{\gamma_a})\ t^{\ind{\ve_k=-1}}\,w_k(\gamma_b)\,t^{-\ind{\ve_{i_b}=-1}}\\
   &=\bigl(t\,V_F(\gamma_a,\gamma_b)-V_F(\gamma_b,\gamma_a)\bigr)\,t^{-\ind{\ve_{i_b}=-1}} ,
\end{align*}
which is the entry $(a,b)$ of $(tV_F-V_F^{\T})\mathsf Z'_t$. In the same way
\begin{align*}
   M^b_{D,\hate}(t^{\ve})_{a,b}&=-\sum_kt^{-\ind{\ve_{i_a}=-1}}\,w_k(\gamma_a)\,\bigl(V^\uparrow_k-tV'^{\uparrow}_k\bigr)(z^{\gamma_b})\\
   &=t^{-\ind{\ve_{i_a}=-1}}\bigl(t\,V_{F^\uparrow}(\gamma_a,\gamma_b)-V_{F^\uparrow}(\gamma_b,\gamma_a)\bigr),
\end{align*}
because $\sum_kV^\uparrow_k(z^{\gamma_b})w_k(\gamma_a)=V_{F^\uparrow}(\gamma_b,\gamma_a)$ and $\sum_kV'^{\uparrow}_k(z^{\gamma_b})w_k(\gamma_a)=V_{F^\uparrow}(\gamma_a,\gamma_b)$. The formulas for the edge polynomials follow from Corollaries~\ref{cor.seifert}~(iii) and~\ref{cor.seifertb}, the monomial $x^{\CT}$ having the value $t^{\ve(\CT')}$ at the Alexander point and $\det\mathsf Z'_t$ being $t^{-\nu_-}$. At $t=1$ the matrix $\mathsf Z'_t$ is the identity.
\end{proof}

So for every diagram the deformed Seifert matrix at the Alexander point is the classical matrix $tV_F-V_F^{\T}$ of the canonical surface, up to units on the columns, and the second flavor gives the classical matrix of the canonical surface on the other side of the plane; since $P_{D,\hate}(t^{\ve})=\tP_{D,\hate}(t^{\ve})$ by Lemma~\ref{lem.level}~(i), the two classical determinants coincide. In particular the degree bound of Theorem~\ref{thm.degrees} below is literally the classical bound $\deg\det(tV_F-V_F^{\T})\le\nu$ read at the Alexander point, and the extreme coefficients of $P_{D,\hate}(t^{\ve})$ are $\det V_F$ and $(-1)^\nu\det V_F$; we draw the consequences for the extremal monomials of the edge polynomials in Corollary~\ref{cor.extremesums}.

Since there is no sign in Theorem~\ref{thm.alexpoint}, and since the edge polynomials at the Alexander point are known exactly (Lemma~\ref{lem.alexpoly}), we obtain the classical formula for the Alexander polynomial in terms of a Seifert matrix with its sign, in the conventions of Section~\ref{ssec.localization}.

\begin{corollary}[the classical formula]\label{cor.classical}
Let $V$ be a Seifert matrix of a link $L$ with $\mu$ components, that is, the matrix $\bigl(\lk(\gamma_a,\gamma_b^+)\bigr)_{a,b}$ of a basis of the first homology of a connected Seifert surface of $L$, with the convention of Section~\ref{ssec.localization} for the positive side, and let $n$ be its size. Then
\[
   \det\bigl(tV-V^{\T}\bigr)=(-1)^{\mu-1}\,t^{n/2}\,\Delta_L(t),
\]
where $\Delta_L(t)=\nabla_L(t^{1/2}-t^{-1/2})$ is the Conway-normalized Alexander polynomial.
\end{corollary}
\begin{proof}
For every square matrix $V$ of size $n$, the polynomial $f(t)=\det(tV-V^{\T})$ satisfies $f(t)=\det(tV^{\T}-V)=(-t)^nf(t^{-1})$. Let first $V=V_F$ for the canonical surface of a diagram $D$ of $L$, so that $n=\nu$. By Theorem~\ref{thm.alexpoint} and Lemma~\ref{lem.alexpoly}, $f(t)=(-1)^{\mu-1}t^a\Delta_L(t)$ with $a=\frac12m_{D,\hate}-\ve(\CT')+\nu_-$. If $\Delta_L=0$, both sides of the formula vanish. Otherwise, since $\Delta_L(t^{-1})=(-1)^{\mu-1}\Delta_L(t)$ and $\nu\equiv\mu-1$ modulo $2$, the symmetry of $f$ reads $t^a\Delta_L(t)=t^{\nu-a}\Delta_L(t)$, so that $a=\nu/2$. For another Seifert matrix $V$ of $L$, of size $n$, we use the classical fact that $t^{-n/2}\det(tV-V^{\T})=\det(t^{1/2}V-t^{-1/2}V^{\T})$ is the same for all the Seifert matrices of $L$: any two of them are $S$-equivalent, and this determinant is not changed by a congruence or by an elementary enlargement \cite[Chapter~8]{Lickorish}. Every link has a connected diagram, which gives the formula in general.
\end{proof}

So $\det(t^{1/2}V-t^{-1/2}V^{\T})=\nabla_L(t^{-1/2}-t^{1/2})$: with our conventions for the crossings and for the push-offs, the determinant satisfies the skein relation of the Conway polynomial in the variable $t^{-1/2}-t^{1/2}$, and the sign $(-1)^{\mu-1}$ is the price of writing $\Delta_L(t)=\nabla_L(t^{1/2}-t^{-1/2})$ as in Section~\ref{ssec.signs}. For the positive Hopf link, for instance, $V=(-1)$ (Proposition~\ref{prop.twobraid} for $n=2$) and $\det(tV-V^{\T})=1-t=-t^{1/2}(t^{1/2}-t^{-1/2})$.

\subsection{The Alexander module, and the dependence on the choices}\label{ssec.choices}
At the Alexander point the deformed Seifert matrices are Alexander matrices in the usual sense. We prove this first (Proposition~\ref{prop.seifertalex}), and we then determine how the matrix $M_{D,\hate}(x)$ depends on the choices which enter its definition (Lemma~\ref{lem.choices}).

\begin{proposition}\label{prop.seifertalex}
Let $t$ be an indeterminate and let $M_{D,\hate}(t^{\ve})$ be the matrix $M_{D,\hate}(x)$ at $x_j=t^{\ve_j}$, $j=1,\dots,c$, a $\nu\times\nu$ matrix over $\BZ[t^{\pm1}]$.
\begin{enumerate}\renewcommand{\theenumi}{\roman{enumi}}
\item The $\BZ[t^{\pm1}]$-modules presented by $M_{D,\hate}(t^{\ve})$ and by the cut Alexander matrix $A_{D,\hate}(t)$ of Definition~\ref{def.coloring} are isomorphic.
\item $\det M_{D,\hate}(t^{\ve})$ is $\Delta_L(t)$ up to multiplication by a unit, and the corank of $M_{D,\hate}(t^{\ve})$ over $\BQ(t)$ is the nullity $\b(L)$ of Definition~\ref{def.nullity}, the rank of the Alexander module of the total linking number cover.
\end{enumerate}
The same holds for $M^b_{D,\hate}(t^{\ve})$ and $B_{D,\hate}(t)$.
\end{proposition}
\begin{proof}
(i) Consider the operations which lead from the matrix $\hat A$ of the proof of Proposition~\ref{prop.seifert} to the matrix $M_{D,\hate}(x)$: those of that proof, those of the proof of Theorem~\ref{thm.flows} which lead to \eqref{eq.kappaSeifert}, and those of the proof of Corollary~\ref{cor.seifert}~(iii). They are row operations (which do not change the module presented by a matrix), the adjunction of generators together with relations expressing them in terms of the old ones (Step~1 of the proof of Proposition~\ref{prop.seifert}), the elimination of generators which are expressed by one of the relations in terms of the others (Step~3, the elimination of the generators $b_C$ in the proof of Theorem~\ref{thm.flows} and of the $E_{\CT}$-coordinates in Corollary~\ref{cor.seifert}~(iii)), and changes of basis of the generators (the multiplication by $(Y\,|\,E_{\CT})$); none of them changes the isomorphism class of the presented module. The eliminations use the invertibility of $(\hat\partial E_{\CT})^{\T}$ over $\BZ$ and of $\CC(x)E_{\CT}$ over $\Lam_c$, and the latter remains invertible after the specialization $x_j\mapsto t^{\ve_j}$, its determinant being a monomial. Thus $M_{D,\hate}(t^{\ve})$ presents the same module as the matrix $\hat A$ at $x=t^{\ve}$, which is $A_D(t)$ with the row $\hate$ replaced by $e_{\hate}^{\T}$; this matrix presents $\BZ[t^{\pm1}]^E/(\text{rows of $A_D(t)$ other than $\hate$},\ e_{\hate})$, which is the module presented by $A_{D,\hate}(t)$.

(ii) The determinant follows from Corollary~\ref{cor.seifert} (iii) and Lemma~\ref{lem.alexpoly}, and the corank from (i), since the corank of a square matrix over a field is the dimension of its cokernel, and the modules presented by the two matrices, tensored with $\BQ(t)$, are isomorphic (Corollary~\ref{cor.corank}).
\end{proof}

Together with Theorem~\ref{thm.alexpoint} this proves, for canonical surfaces, the classical fact that $tV-V^{\T}$ presents the Alexander module of the total linking number cover (see \cite[Chapter~6]{Lickorish} and \cite[Chapter~8]{BurdeZieschang} for arbitrary Seifert surfaces). The deformed matrix $M_{D,\hate}(x)$ inserts the crossing variables into this picture, and Corollary~\ref{cor.seifertformula} couples the two flavors crossing by crossing through $\BE_D$: by \eqref{eq.ED}, $\BE_D$ evaluates the first factor at $x=\bx$ and the second at $y=\by$ and replaces each product $(x_j-\bx_j)(y_j-\by_j)$ by $g_j=1-(t_0t_1)^{\ve_j}$ (compare the proof of Corollary~\ref{cor.expansion}).

The matrix $M_{D,\hate}(x)$ depends on the choices which have been made: the base edges $b_C$ of the Seifert circles $C\ne\hat C$ (the base edge of $\hat C$ is $\hate$), the basis $Y$ of the cycle space and the spanning tree $\CT$. The dependence is mild.

\begin{lemma}[dependence on the choices]\label{lem.choices}
\begin{enumerate}\renewcommand{\theenumi}{\roman{enumi}}
\item The matrix $M_{D,\hate}(x)$ does not depend on the base edges $b_C$, $C\ne\hat C$.
\item If $Y$ is replaced by $YQ$ with $Q\in\GL_\nu(\BZ)$, then $M_{D,\hate}(x)$ is replaced by $Q^{\T}M_{D,\hate}(x)Q$.
\item If $\CT$ is replaced by another spanning tree, then $M_{D,\hate}(x)$ is replaced by $M_{D,\hate}(x)R(x)$ for a matrix $R(x)\in\GL_\nu(\Lam_c)$ with $R(1)=1$.
\end{enumerate}
The same holds for $M^b_{D,\hate}(y)$: it does not depend on the base edges, it is replaced by $Q^{\T}M^b_{D,\hate}(y)Q$ under $Y\mapsto YQ$, and by $R(y)^{\T}M^b_{D,\hate}(y)$ under a change of spanning tree, where $R(y)$ is the matrix of (iii) with $x$ replaced by $y$.
\end{lemma}
\begin{proof}
(i) By \eqref{eq.Mentries}, $M_{D,\hate}(x)_{a,b}=\sum_k\pi_k(z^{\gamma_a})Y(x)_{k,b}$, and the matrix $Y(x)$ does not involve the base edges. Let $C\ne\hat C$. If the base edge of $C$ is changed, the lift $z^{\gamma_a}$ is replaced by another flow with the same transfers, which still vanishes on the base edges of the other circles; by Lemma~\ref{lem.circles}~(i) the difference is a combination of Seifert circles, and since it vanishes on the base edges of the circles other than $C$ it is an integral multiple of the flow $C$. Now $\sum_k\pi_k(C)Y(x)_{k,b}$ is the entry $(C,b)$ of $\CC(x)Y(x)=0$. So the matrix $M_{D,\hate}(x)$ does not change.

(ii) By \eqref{eq.Yx}, replacing $Y$ by $YQ$ replaces $Y(x)$ by $Y(x)Q$.

(iii) Let $Y'(x)$ be the matrix \eqref{eq.Yx} of the other tree. The columns of $Y(x)$ and those of $Y'(x)$ are bases of $\ker\CC(x)$ (Corollary~\ref{cor.seifert} (ii)), so $Y'(x)=Y(x)R(x)$ with $R(x)\in\GL_\nu(\Lam_c)$, and $R(1)=1$ because $Y'(1)=Y=Y(1)$ and the columns of $Y$ are linearly independent. The statements for $M^b_{D,\hate}(y)$ are proved in the same way from the second formula of \eqref{eq.Mentries}: a change of the base edge of $C$ changes the lifts by multiples of $C$, and $\sum_kY(y)_{k,a}\tpi_k(C)$ is the entry $(C,a)$ of $\CC(y)Y(y)=0$; and $Y\mapsto YQ$, resp.\ a change of tree, replaces $Y(y)$ by $Y(y)Q$, resp.\ by $Y(y)R(y)$.
\end{proof}

In the terms of Section~\ref{ssec.cochains}, (ii) and (iii) are changes of basis of the two modules $H_1(F_D)$ and $H_x$ on which the deformed Seifert form is defined, and (i) says that this form is well defined on the quotient of $\FF_{D,\hate}$ by the Seifert circles.

\subsection{Examples}\label{ssec.examples}
We work out the deformed Seifert matrices of the standard diagrams of the trefoil and of the figure-eight knot, for all their cut edges, and of the closures of the $2$-braids $\sigma_1^n$. We use the following notation. Let $\hate$ be an edge of $D$ which lies on the under-strand of a crossing; let $j_0$ be the crossing which $\hate$ leaves as the under-strand, $\hate=u'_{j_0}$, if there is one, and otherwise the crossing which $\hate$ enters as the under-strand, $\hate=u_{j_0}$. (An edge of the form $o'_j=o_k$, which leaves a crossing and enters the next one as the over-strand, lies on no under-strand.) For a crossing $j$ let
\[
   f_j=e_j^{\T}Y\in\BZ^{1\times\nu}
\]
be the row $j$ of $Y$, the row vector of the transfers $w_j(\gamma_1),\dots,w_j(\gamma_\nu)$ of the basis cycles, and let $g_j\in\BZ^\nu$ be the column vector with the entries $V_j(z^{\gamma_a})$ of \eqref{eq.Vj}; by Theorem~\ref{thm.deformedseifert}~(i), the Seifert matrix of $F_D$ in the basis $Y$ is $V_F=\sum_jg_jf_j$, the sum over the bands of matrices of rank one. Recall that the two edges of an under-strand have cut polynomials which differ by a unit, $P_{D,u_j}=x_jP_{D,u'_j}$ (Lemma~\ref{lem.understrand}). The matrices of the examples were computed from Corollary~\ref{cor.seifert}~(iii) by computer.

\begin{example}[the trefoil]\label{ex.trefoil}
Let $D$ be the closure of the braid $\sigma_1^3$, with the notation of Examples~\ref{ex.trefoilflows} and~\ref{ex.trefoilseifert}: the spanning tree $\CT=\{1\}$ and the basis $\gamma_1=e_2-e_1$, $\gamma_2=e_3-e_1$ of the cycle space, so that
\[
   f_1=(-1,-1),\qquad f_2=(1,0),\qquad f_3=(0,1),\qquad V_F=\begin{pmatrix}-1&-1\\0&-1\end{pmatrix}.
\]
Consider first $\hate=u'_3=o_1$, so that $j_0=3$ and $\hat C=C_1$. Then $\mathrm{Pre}_{C_1}(o_j)=\{1,\dots,j-1\}$, and with the base edge $u_1$ on $C_2$ also $\mathrm{Pre}_{C_2}(u_j)=\{1,\dots,j-1\}$; the row $j$ of $K(x)$ is $\d_j+\sum_{k<j}(1+x_k)\d_k$, the full turn of $C_2$ is $-\d_1-\d_2-\d_3$, so that $\CC(x)E_{\CT}=-1$ and $Y(x)=Y$. One finds $P_{D,\hate}(x)=1-x_1+x_1x_2$ and
\begin{align*}
   M_{D,\hate}(x)&=\begin{pmatrix}1-x_1&-x_1\\1+x_2-x_1&1-x_1\end{pmatrix}\\
   &=x_1\begin{pmatrix}-1&-1\\-1&-1\end{pmatrix}+x_2\begin{pmatrix}0&0\\1&0\end{pmatrix}+\begin{pmatrix}1&0\\1&1\end{pmatrix}
   =x_1g_1f_1+x_2g_2f_2-V_F^{\T}
\end{align*}
with $g_1=(1,1)^{\T}$ and $g_2=(0,1)^{\T}$, the vectors of the flows of the two lifts through $u'_1$ and through $u'_2$ (and $g_3=0$). So the deformed Seifert matrix is $x_1V_1+x_2V_2-V_F^{\T}$ with $V_1=g_1f_1$ and $V_2=g_2f_2$ the contributions of the bands $1$ and $2$ to the Seifert form, and $V_1+V_2=V_F$: it is obtained from $tV_F-V_F^{\T}$ by letting the variable of each crossing multiply the share of its band; at the Alexander point it is $tV_F-V_F^{\T}$.

The other five cut edges are treated in the same way, with the same tree and basis. We group them according to the crossing $j_0$.
\begin{itemize}
\item $j_0=1$, $\hate=u'_1=o_2$ or $\hate=u_1=o'_3$: $P_{D,u'_1}=1-x_2+x_2x_3$, $P_{D,u_1}=x_1P_{D,u'_1}$, and
\[
   M_{D,u'_1}(x)=M_{D,u_1}(x)=\begin{pmatrix}1-x_2&-x_3\\1&1-x_3\end{pmatrix}=x_2g_2f_2+x_3g_3f_3-V_F^{\T},\quad g_2=\tbinom{-1}{0},\ g_3=\tbinom{-1}{-1}.
\]
\item $j_0=2$, $\hate=u'_2=o_3$ or $\hate=u_2=o'_1$: $P_{D,u'_2}=1-x_3+x_1x_3$, $P_{D,u_2}=x_2P_{D,u'_2}$,
\[
   M_{D,u'_2}(x)=\begin{pmatrix}1-x_1&-x_1\\1&1-x_3\end{pmatrix}=x_1g_1f_1+x_3g_3f_3-V_F^{\T},\quad g_1=\tbinom{1}{0},\ g_3=\tbinom{0}{-1},
\]
and
\[
   M_{D,u_2}(x)=\begin{pmatrix}x_1^{-1}x_2-x_2&x_1^{-1}x_3-1-x_3\\x_1^{-1}x_2&x_1^{-1}x_3-x_3\end{pmatrix},\qquad \det M_{D,u_2}(x)=x_1^{-1}P_{D,u_2}(x).
\]
\item $j_0=3$, $\hate=u'_3=o_1$ (treated above) or $\hate=u_3=o'_2$: $P_{D,u_3}=x_3P_{D,u'_3}$, and
\[
   M_{D,u_3}(x)=\begin{pmatrix}1-x_2&-x_3\\x_1^{-1}x_2&x_1^{-1}x_3-x_3\end{pmatrix},\qquad \det M_{D,u_3}(x)=x_1^{-1}P_{D,u_3}(x).
\]
\end{itemize}
For the three edges $u'_{j_0}$ the deformed Seifert matrix is thus of the form $\sum_{j\ne j_0}x_jg_jf_j-V_F^{\T}$ with $\sum_{j\ne j_0}g_jf_j=V_F$, exactly as for $\hate=u'_3$. For $u_2$ and $u_3$ the inverse $x_1^{-1}$ of the variable of the tree edge enters through $Y(x)$, and the determinant is $P_{D,\hate}(x)$ times the unit $x_1^{-1}$ (recall that $P_{D,u_{j_0}}=x_{j_0}P_{D,u'_{j_0}}$).
\end{example}

\begin{example}[the figure-eight knot]\label{ex.figeight}
Let $D$ be the closure of $(\sigma_1\sigma_2^{-1})^2$, the standard diagram of the figure-eight knot, with the crossings $1,2,3,4$ in the order of the braid word (Figure~\ref{fig.diagrams}, right); $\ve_1=\ve_3=1$ and $\ve_2=\ve_4=-1$. The Seifert circles are the three strands of the braid closed up (Figure~\ref{fig.circles}, right): $C_1$ contains $o_1=u'_3$ and $u'_1=o_3$; $C_2$ contains $u_1=o'_4$, $o'_1=u_2$, $o'_2=u_3$ and $o'_3=u_4$; $C_3$ contains $o_2=u'_4$ and $u'_2=o_4$. Thus $A_1=A_3=C_1$, $A_2=A_4=C_3$, $B_j=C_2$ for all $j$, and $\nu=2$. We take $\CT=\{1,2\}$ and the basis $\gamma_1=e_3-e_1$, $\gamma_2=e_4-e_2$ of the cycle space, so that
\[
   f_1=(-1,0),\qquad f_2=(0,-1),\qquad f_3=(1,0),\qquad f_4=(0,1),\qquad V_F=\begin{pmatrix}-1&-1\\0&1\end{pmatrix},
\]
the Seifert matrix of $F_D$ in this basis, given by Theorem~\ref{thm.deformedseifert}~(i); it satisfies $\det(tV_F-V_F^{\T})=-t^2+3t-1$. Here the Seifert circle $C_2$ has crossings on both sides, and the canonical surface above the plane has a different Seifert matrix in the same basis, $V_{F^\uparrow}=\bigl(\begin{smallmatrix}-1&0\\1&1\end{smallmatrix}\bigr)$, with the same determinant $\det(tV_{F^\uparrow}-V_{F^\uparrow}^{\T})=-t^2+3t-1$. Every edge of $D$ lies on an under-strand, and we group the eight cut edges according to the crossing $j_0$, i.e.\ into the four under-strands $\{u'_{j_0},u_{j_0}\}$. In each case $P_{D,u_{j_0}}=x_{j_0}P_{D,u'_{j_0}}$, and the deformed Seifert matrices (computed from Corollary~\ref{cor.seifert} (iii)) are the following.
\begin{itemize}
\item $j_0=1$ ($\ve_{j_0}=1$): $P_{D,u'_1}=x_2+x_4-x_2x_3-x_2x_4+x_2x_3x_4$, and
\[
   M_{D,u'_1}(x)=M_{D,u_1}(x)=\begin{pmatrix}1-x_3&-x_2^{-1}x_4\\1&1-x_4\end{pmatrix},\qquad \det M_{D,u'_1}(x)=x_2^{-1}P_{D,u'_1}(x).
\]
\item $j_0=3$ ($\ve_{j_0}=1$): $P_{D,u'_3}=x_2+x_4-x_1x_4-x_2x_4+x_1x_2x_4$, and
\[
   M_{D,u'_3}(x)=\begin{pmatrix}1-x_1&-x_1-(1-x_1)x_2^{-1}x_4\\1&1-x_4\end{pmatrix},\quad
   M_{D,u_3}(x)=\begin{pmatrix}x_1^{-1}x_3-x_3&-1\\x_1^{-1}x_3&x_2^{-1}x_4-x_4\end{pmatrix},
\]
with determinants $x_2^{-1}P_{D,u'_3}(x)$ and $x_1^{-1}x_2^{-1}P_{D,u_3}(x)$.
\item $j_0=2$ ($\ve_{j_0}=-1$): $P_{D,u'_2}=x_1+x_3-x_1x_3-x_3x_4+x_1x_3x_4$, and
\[
   M_{D,u'_2}(x)=M_{D,u_2}(x)=\begin{pmatrix}x_1^{-1}x_3-x_3&-1\\1&1-x_4\end{pmatrix},\qquad \det M_{D,u'_2}(x)=x_1^{-1}P_{D,u'_2}(x).
\]
\item $j_0=4$ ($\ve_{j_0}=-1$): $P_{D,u'_4}=x_1+x_3-x_1x_2-x_1x_3+x_1x_2x_3$, and
\[
   M_{D,u'_4}(x)=\begin{pmatrix}x_1^{-1}x_3-x_3&-1\\1-x_2+x_1^{-1}x_2x_3&1-x_2\end{pmatrix},\qquad
   M_{D,u_4}(x)=\begin{pmatrix}1-x_3&-x_2^{-1}x_4\\x_1^{-1}x_3&x_2^{-1}x_4-x_4\end{pmatrix},
\]
with determinants $x_1^{-1}P_{D,u'_4}(x)$ and $x_1^{-1}x_2^{-1}P_{D,u_4}(x)$.
\end{itemize}
Here the variables $x_1,x_2$ of the tree edges enter through $Y(x)$, with their inverses, and the entries are no longer of degree one in the crossing variables: products such as $x_2^{-1}x_4$ and $x_1^{-1}x_2x_3$ occur. At the Alexander point each of the eight matrices presents the Alexander module of the figure-eight knot (Proposition~\ref{prop.seifertalex}), and its determinant is $-t^2+3t-1$ up to a unit; for instance $M_{D,u'_1}(t^{\ve})=\bigl(\begin{smallmatrix}1-t&-1\\1&1-t^{-1}\end{smallmatrix}\bigr)=(tV_F-V_F^{\T})\diag(1,t^{-1})$, in accordance with Theorem~\ref{thm.alexpoint}.
\end{example}

For the closures of the $2$-braids the deformed Seifert matrix can be written down in general, and it has the structure observed for the trefoil.

\begin{proposition}[closed $2$-braids]\label{prop.twobraid}
Let $n\ge2$, let $D$ be the closure of the braid $\sigma_1^n$, with the crossings $1,\dots,n$ in the order of the braid word, and let $\hate=u'_n=o_1$, so that $j_0=n$. With the spanning tree $\CT=\{1\}$ and the basis $\gamma_a=e_a-e_1$, $a=2,\dots,n$, of the cycle space, so that $f_1=-(1,\dots,1)$ and $f_b=e_b^{\T}$ for $2\le b\le n$,
\begin{equation}\label{eq.twobraid}
   M_{D,\hate}(x)=\sum_{j=1}^{n-1}x_jg_jf_j-V^{\T},\qquad g_1=(1,\dots,1)^{\T},\quad g_b=\sum_{a>b}e_a\ \ (2\le b\le n-1),
\end{equation}
and $V=\sum_{j=1}^{n-1}g_jf_j=-U$, where $U$ is the $(n-1)\times(n-1)$ matrix with entries $U_{ab}=\ind{a\le b}$; explicitly, $M_{D,\hate}(x)_{ab}=\ind{a=b}-x_1+\ind{a>b}(1+x_b)$ for $2\le a,b\le n$. The matrix $-U$ is the Seifert matrix $V_F$ of $F_D$, the standard Seifert surface of the $(2,n)$ torus link, with respect to the basis $\gamma_2,\dots,\gamma_n$. Thus $M_{D,\hate}(x)$ is obtained from $tV_F-V_F^{\T}$ by distributing the Seifert matrix $V_F=\sum_{j\ne j_0}g_jf_j$ over the crossings $j\ne j_0$, the variable $x_j$ multiplying the rank-one summand $g_jf_j$; moreover $g_j=V_F(e_j-e_n)$ for $2\le j\le n-1$ and $g_1=V_F(0-e_n)$, where $e_j\in\BZ^{n-1}$ is the coordinate vector of $\gamma_j$ and $0$ stands for the tree edge $1$: so $g_jf_j$ is the Seifert form evaluated on the cycle which runs through the band $j$ and back through the band $j_0$ of the cut edge.
\end{proposition}
\begin{proof}[Proof of Proposition~\ref{prop.twobraid}]
The Seifert circles are $C_1\ni o_j=u'_{j-1}$ and $C_2\ni u_j=o'_{j-1}$ (indices modulo $n$), with $A_j=C_1$ and $B_j=C_2$ for all $j$ and $\hat C=C_1$; with the base edge $u_1$ of $C_2$, $\mathrm{Pre}_{C_1}(o_j)=\mathrm{Pre}_{C_2}(u_j)=\{1,\dots,j-1\}$. By \eqref{eq.jumpsum} and \eqref{eq.Krow}, $K(x)=1+\Lambda(1+X)$, where $\Lambda$ is the strictly lower triangular $n\times n$ matrix with all entries below the diagonal equal to $1$ and $X=\diag(x_1,\dots,x_n)$. The only full turn is $\CC_{C_2}=-\sum_k\d_k$, so $\CC(x)=-(1,\dots,1)$, $\CC(x)E_{\CT}=-1$, and $Y(x)=Y+E_{\CT}\CC(x)Y=Y$ because every column of $Y$ has entries $1,-1$ and $0$ otherwise. With $Y_{jb}=\ind{j=b}-\ind{j=1}$ one finds, for $2\le a,b\le n$,
\[
   (Y^{\T}Y)_{ab}=1+\ind{a=b},\qquad (Y^{\T}\Lambda Y)_{ab}=\ind{b<a}-1,\qquad (Y^{\T}\Lambda XY)_{ab}=\ind{b<a}\,x_b-x_1,
\]
because $(Y^{\T}NY)_{ab}=N_{ab}-N_{a1}-N_{1b}+N_{11}$ for every $n\times n$ matrix $N$, and $\Lambda_{ab}=\ind{b<a}$. The sum of the three is the stated formula for $M_{D,\hate}(x)=Y^{\T}K(x)Y$. In this formula, the coefficient of $x_1$ is the matrix with all entries $-1$, which is $g_1f_1$; the coefficient of $x_b$, $2\le b\le n-1$, is the matrix with the entries $\ind{a>b}$ in the column $b$ and $0$ elsewhere, which is $g_bf_b$; $x_n$ does not occur; and the constant term $\ind{a\ge b}$ is $U^{\T}=-V^{\T}$ with $V=\sum_jg_jf_j=-U$. Finally, $-U$ is the Seifert matrix of $F_D$ in the basis $\gamma_2,\dots,\gamma_n$, by Theorem~\ref{thm.deformedseifert}~(i): the lift of $\gamma_a$ is $\sum_{j<a}(e_{u'_j}-e_{u_{j+1}})$, as in Example~\ref{ex.trefoilflows}, so that $V_j(z^{\gamma_a})=z^{\gamma_a}_{u'_j}=\ind{j<a}$ is the entry $a$ of $g_j$ for $j\le n-1$ and vanishes for $j=n$, and $V_F=\sum_jg_jf_j$ (for $n=3$ it is the matrix of Example~\ref{ex.trefoilseifert}). The last claim follows from $-U(e_j-e_n)=\1_{a>j}$ for $2\le j\le n-1$ and $-U(0-e_n)=\1$, where $\1_{a>j}$ denotes the vector with entries $\ind{a>j}$.
\end{proof}

\begin{remark}\label{rem.notgeneral}
The structure of \eqref{eq.twobraid}, in which each crossing variable multiplies a matrix of rank one and the constant term is the negative transpose of the Seifert matrix $V_F$, is special to these diagrams. We have checked by computer that for the standard diagrams of $6_2$ and $6_3$ (closures of $3$-braids, with $\nu=4$) the deformed Seifert matrix $M_{D,\hate}(x)$ is not equivalent over $\Lam_c$, for any cut edge, to a matrix $\sum_{j\ne j_0}x_j^{\pm1}g_jf_j-V^{\T}$ with $g_j\in\BZ^\nu$ and $V=\sum_{j\ne j_0}g_jf_j\in\{V_F,V_F^{\T}\}$, for any choice of the exponents: no such matrix has determinant $P_{D,\hate}(x)$ up to a unit of $\Lam_c$. The entries of the Seifert form between bands which are attached in interleaved order to a Seifert circle are not attached to any single crossing, and the way the crossing variables deform the Seifert form is less simple than \eqref{eq.twobraid} suggests. The rank-one structure does hold at the level of the $c\times c$ matrix $N_{D,\hate}(x)$ of the crossing forms, for every diagram (Remark~\ref{rem.flows}~(ii)), and the rank-one summand of $x_k$ is the contribution of the band $k$ to the Seifert form or to its transpose (Theorem~\ref{thm.deformedseifert}); it is lost in passing to the Schur complement $M_{D,\hate}(x)$.
\end{remark}

\subsection{A counterexample: the Kinoshita--Terasaka knot}
The deformed Seifert matrices depend on the diagram, through the crossing variables, and not only on the Seifert form of the canonical surface. The following remark shows that this cannot be avoided: there is no formula for the Links--Gould polynomial in terms of a Seifert matrix alone.

\begin{remark}[$\LG$ is not a function of the Seifert form]\label{rem.KT}
The Alexander polynomial, the signatures and the Blanchfield pairing of a knot are functions of the Seifert matrix of any of its Seifert surfaces, and the Seifert matrices of all Seifert surfaces of a knot are $S$-equivalent; conversely, the $S$-equivalence class of a Seifert matrix is determined by the Blanchfield pairing on the Alexander module \cite{Trotter, Kearton}. Nothing of the kind holds for $\LG$. Let $K$ be the Kinoshita--Terasaka knot $11n_{42}$ \cite{KinoshitaTerasaka}, whose Alexander polynomial is $1$. Its Alexander module is trivial, so its Seifert matrices are $S$-equivalent to the empty matrix, like those of the unknot. But $\LG(K)\ne1=\LG(\text{unknot})$: from the planar diagram code of $11n_{42}$ in the Knot Atlas \cite{KnotAtlas}, Vogel's algorithm gives the closed $4$-braid $\sigma_1^{-1}\sigma_2^2\sigma_3^{-1}\sigma_2^{-1}\sigma_1\sigma_2^{-1}\sigma_1\sigma_2^{-2}\sigma_3\sigma_2\sigma_1^{-1}$, which represents $K$ or its mirror image, and the state sum of \cite[Definition~2.2]{HK} gives
\begin{align*}
   \LG(K)={}&-23+18(t_0^{-1}+t_1^{-1})+9(t_0+t_1)-3(t_0^{-2}+t_1^{-2})-8t_0^{-1}t_1^{-1}-9(t_0^{-1}t_1+t_0t_1^{-1})\\
   &+2(t_0^{2}+t_1^{2})+6t_0t_1-(t_0^{-3}+t_1^{-3})-7(t_0^{-2}t_1^{-1}+t_0^{-1}t_1^{-2})+2(t_0^{-2}t_1+t_0t_1^{-2})\\
   &+(t_0^{-1}t_1^{2}+t_0^{2}t_1^{-1})-(t_0^{3}+t_1^{3})-11(t_0t_1^{2}+t_0^{2}t_1)+6(t_0^{-3}t_1^{-1}+t_0^{-1}t_1^{-3})+14t_0^{-2}t_1^{-2}\\
   &+3(t_0t_1^{3}+t_0^{3}t_1)+10t_0^{2}t_1^{2}-(t_0^{-4}t_1^{-1}+t_0^{-1}t_1^{-4})-7(t_0^{-3}t_1^{-2}+t_0^{-2}t_1^{-3})\\
   &-2(t_0^{2}t_1^{3}+t_0^{3}t_1^{2})+(t_0^{-4}t_1^{-2}+t_0^{-2}t_1^{-4})+2t_0^{-3}t_1^{-3},
\end{align*}
which is symmetric in $t_0,t_1$ and restricts to $1=\Delta_K(t)^2$ on the curve $t_1=t_0^{-1}$, as it must by Lemma~\ref{lem.alexpoly}. After the appropriate change of variables, this confirms the computation of \cite[Section~4.9]{DeWitKauffmanLinks}. So $\LG$ is not determined by the $S$-equivalence class of the Seifert form, and there is no formula expressing $\LG$ in terms of a Seifert matrix alone: any such formula, valid for all Seifert surfaces of all knots, would be invariant under the enlargements of Seifert matrices, which are realized by attaching hollow handles to a surface, and hence under $S$-equivalence. The deformed Seifert matrices of Corollary~\ref{cor.seifertformula} depend on the diagram through the crossing variables, and the coupling of the two flavors by $\BE_D$ is where the quantum invariant differs from the classical ones.
\end{remark}

\subsection{Moves}\label{ssec.moves}
Two closed braid diagrams of the same link are related by conjugation, Markov stabilization and the braid relations (Markov's theorem), and the invariance of $\LG$ under these moves is built into the $R$-matrix proof of \eqref{eq.intro.formula}. The moves have a transparent effect on the edge polynomials, which we record because they govern the dependence of the deformed Seifert matrices on the diagram; the Reidemeister move of type II, which attaches a hollow handle to the surface, is included because it is the elementary enlargement of a Seifert surface (Remark~\ref{rem.KT}).

The proofs use the following elimination. In the matrix $1-T^a_D(x)$, the row of an edge which leaves a crossing expresses the generator of that edge in terms of the generators of the incoming edges of the crossing ($[o'_j]=[o_j]$ and $[u'_j]=(1-x_j)[o_j]+x_j[u_j]$ by \eqref{eq.Ta}). If $w$ is a subword of the braid word none of whose internal edges is the cut edge, the generators of the internal edges of $w$ can be eliminated from the presentation as in Step~3 of the proof of Proposition~\ref{prop.seifert}, without changing the determinant of the cut matrix; the rows of the outgoing edges of $w$ then express them in terms of the incoming edges of $w$ by the product of the \emph{blocks} of its crossings, where the block of a positive crossing of two adjacent strands is $\rho(x_j)=\bigl(\begin{smallmatrix}1-x_j&x_j\\1&0\end{smallmatrix}\bigr)$ and that of a negative one is $\rho^-(x_j)=\bigl(\begin{smallmatrix}0&1\\x_j&1-x_j\end{smallmatrix}\bigr)$, acting on the column of the generators of the left and right incoming edges (the over-strand of a positive crossing comes in from the left), and the blocks are composed in the order of the word, the block of a later crossing being applied after that of an earlier one. The cut polynomial $P_{D,\hate}(x)$ depends on $w$ only through this product. For the second flavor the same holds with the transposed blocks $\rho(y_j)^{\T}$ and $\rho^-(y_j)^{\T}$, by \eqref{eq.Tb}.

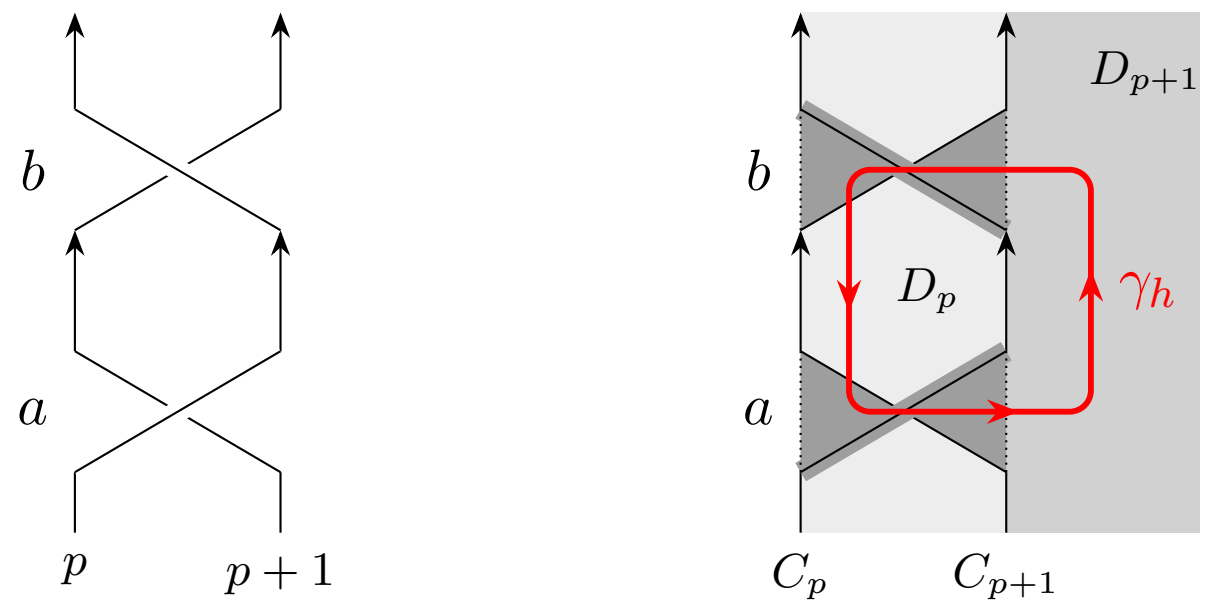
\begin{figure}[ht]
\centering
\begin{tikzpicture}[>=Stealth,scale=0.8]
\def\bigon{%
   \draw (0,0) -- (0,0.5); \draw[->] (0,1.5) -- (0,2.5); \draw[->] (0,3.5) -- (0,4.3);
   \draw (1.7,0) -- (1.7,0.5); \draw[->] (1.7,1.5) -- (1.7,2.5); \draw[->] (1.7,3.5) -- (1.7,4.3);
   \draw (1.7,0.5) -- (0,1.5); \draw[preaction={draw=white,line width=4pt}] (0,0.5) -- (1.7,1.5);
   \draw (0,2.5) -- (1.7,3.5); \draw[preaction={draw=white,line width=4pt}] (1.7,2.5) -- (0,3.5);}
\begin{scope}
\bigon
\node[font=\small] at (-0.35,1.0) {$a$}; \node[font=\small] at (-0.35,3.0) {$b$};
\node[font=\scriptsize,below] at (0,0) {$p$}; \node[font=\scriptsize,below] at (1.7,0) {$p+1$};
\end{scope}
\begin{scope}[xshift=6cm]
\fill[black!7] (0,0) rectangle (1.7,4.3);
\fill[black!18] (1.7,0) rectangle (3.3,4.3);
\fill[black!38] (0,0.5) -- (0.85,1) -- (0,1.5) -- cycle; \fill[black!38] (1.7,0.5) -- (0.85,1) -- (1.7,1.5) -- cycle;
\fill[black!38] (0,2.5) -- (0.85,3) -- (0,3.5) -- cycle; \fill[black!38] (1.7,2.5) -- (0.85,3) -- (1.7,3.5) -- cycle;
\draw (0,0) -- (0,0.5); \draw[->] (0,1.5) -- (0,2.5); \draw[->] (0,3.5) -- (0,4.3);
\draw (1.7,0) -- (1.7,0.5); \draw[->] (1.7,1.5) -- (1.7,2.5); \draw[->] (1.7,3.5) -- (1.7,4.3);
\draw (1.7,0.5) -- (0,1.5); \draw[preaction={draw=black!38,line width=4pt}] (0,0.5) -- (1.7,1.5);
\draw (0,2.5) -- (1.7,3.5); \draw[preaction={draw=black!38,line width=4pt}] (1.7,2.5) -- (0,3.5);
\draw[densely dotted] (0,0.5) -- (0,1.5) (1.7,0.5) -- (1.7,1.5) (0,2.5) -- (0,3.5) (1.7,2.5) -- (1.7,3.5);
\draw[red,line width=1.1pt,rounded corners=4pt,postaction={decorate,decoration={markings,mark=at position 0.16 with {\arrow{Stealth[length=6pt]}},mark=at position 0.38 with {\arrow{Stealth[length=6pt]}},mark=at position 0.88 with {\arrow{Stealth[length=6pt]}}}}] (0.4,1) -- (2.4,1) -- (2.4,3) -- (0.4,3) -- cycle;
\node[font=\small] at (-0.35,1.0) {$a$}; \node[font=\small] at (-0.35,3.0) {$b$};
\node[font=\scriptsize,below] at (0,0) {$C_p$}; \node[font=\scriptsize,below] at (1.7,0) {$C_{p+1}$};
\node[font=\scriptsize] at (1.05,2.0) {$D_p$}; \node[font=\scriptsize] at (2.85,3.8) {$D_{p+1}$};
\node[font=\small,red,right] at (2.45,2.0) {$\gamma_h$};
\end{scope}
\end{tikzpicture}
\caption{The hollow handle of Proposition~\ref{prop.hollowhandle}. Left: the bigon $\sigma_p\sigma_p^{-1}$ inserted in the braid, with the new crossings $a$ (positive) and $b$ (negative) between the strands in the positions $p$ and $p+1$. Right: the canonical surface near the bigon, seen from above: the disk $D_p$, on the right of the strand $C_p$, the disk $D_{p+1}$, which lies above it, on the right of $C_{p+1}$, and the bands of the two new crossings (dark), which have opposite half twists. The loop $\gamma_h=e_a+e_b$ runs through the band $a$ from $D_p$ to $D_{p+1}$ and back through the band $b$. The two half twists cancel, and $\gamma_h$ bounds a disk in the complement of the surface, the window between the two bands; compressing the surface along this disk gives back $F_D$. In this sense the surface is $F_D$ with a hollow handle attached.}
\label{fig.hollowhandle}
\end{figure}

\begin{proposition}[hollow handles: the move of type II]\label{prop.hollowhandle}
Let $D'$ be obtained from the closed braid diagram $D$ by inserting $\sigma_i\sigma_i^{-1}$ (or $\sigma_i^{-1}\sigma_i$) at some place of the braid word, with the cut edge outside the bigon, and let $a,b$ be the two new crossings, of signs $\ve_a=-\ve_b$. Then
\[
   P_{D',\hate}(x)=P_{D,\hate}(x)+(x_ax_b-1)\,R(x),\qquad \tP_{D',\hate}(y)=\tP_{D,\hate}(y)+(y_ay_b-1)\,\tilde R(y),
\]
with $R$ and $\tilde R$ polynomials in the variables of the crossings of $D$, namely the coefficients of $x_ax_b$ in $P_{D',\hate}$ and of $y_ay_b$ in $\tP_{D',\hate}$. Consequently
\[
   \BE_{D'}\bigl[P_{D',\hate}\,\tP_{D',\hate}\bigr]=\BE_D\bigl[P_{D,\hate}\,\tP_{D,\hate}\bigr]:
\]
the invariance of the determinant formula under the move of type II is the vanishing of $\BE_{D'}$ on the multiples of $x_ax_b-1$ and of $y_ay_b-1$ by polynomials which do not contain the variables of the other flavor at $a,b$. Moreover $F_{D'}$ is $F_D$ with a hollow handle attached (Figure~\ref{fig.hollowhandle}), $\nu(D')=\nu(D)+2$, and in the basis of $H_1(F_{D'})$ consisting of the basis cycles of $F_D$, the fundamental cycle $\gamma_a$ of the band $a$ and the loop $\gamma_h=e_a+e_b$ of the handle (the spanning tree being unchanged), the row of $\gamma_h$ in $M_{D',\hate}(x)$ is $(0,\dots,0,-x_a,1-x_a)$, and the block of the other rows and columns is the deformed Seifert matrix of the diagram $D_a$ obtained from $D$ by inserting the single crossing $a$.
\end{proposition}
\begin{proof}
The product of the blocks of the bigon is $\rho^-(x_b)\rho(x_a)=\bigl(\begin{smallmatrix}1&0\\1-x_ax_b&x_ax_b\end{smallmatrix}\bigr)=1+(x_ax_b-1)\bigl(\begin{smallmatrix}0&0\\-1&1\end{smallmatrix}\bigr)$ for $\sigma_i\sigma_i^{-1}$, and $\rho(x_a)\rho^-(x_b)=1+(x_ax_b-1)\bigl(\begin{smallmatrix}1&-1\\0&0\end{smallmatrix}\bigr)$ for $\sigma_i^{-1}\sigma_i$. When this product is the identity, the rows of the outgoing edges of the bigon identify them with its incoming edges, and eliminating them gives the cut matrix of $D$; so $P_{D',\hate}=P_{D,\hate}$ at $x_ax_b=1$. Since $P_{D',\hate}$ has degree $\le1$ in $x_a$ and in $x_b$ and $P_{D,\hate}$ does not contain them, the difference, written as $\alpha+\beta x_a+\gamma x_b+\d x_ax_b$ with $\alpha,\beta,\gamma,\d$ free of $x_a,x_b$, vanishes identically at $x_b=x_a^{-1}$, which forces $\beta=\gamma=0$ and $\alpha=-\d$. The same argument applies to $\tP$: the product of the transposed blocks of the bigon is $\rho^-(y_b)^{\T}\rho(y_a)^{\T}=\bigl(\rho(y_a)\rho^-(y_b)\bigr)^{\T}$ for $\sigma_i\sigma_i^{-1}$ and $\rho(y_a)^{\T}\rho^-(y_b)^{\T}=\bigl(\rho^-(y_b)\rho(y_a)\bigr)^{\T}$ for $\sigma_i^{-1}\sigma_i$, and by the two formulas above both are the identity when $y_ay_b=1$. For the coupling, $\BE_{D'}$ is multiplicative over distinct crossings (proof of Corollary~\ref{cor.expansion}), and $\BE_{D'}[x_ax_b]=t_0^{\ve_a+\ve_b}=1$, $\BE_{D'}[y_ay_b]=1$, $\BE_{D'}[x_ax_by_ay_b]=1$; so $\BE_{D'}[(x_ax_b-1)f]=0$ when $f$ does not contain $y_a,y_b$, $\BE_{D'}[(y_ay_b-1)f]=0$ when $f$ does not contain $x_a,x_b$, and $\BE_{D'}[(x_ax_b-1)(y_ay_b-1)f]=(1-1-1+1)\BE_{D'}[f]=0$ when $f$ contains none of the four; expanding the product $P_{D'}\tP_{D'}$ gives $\BE_{D'}[P_D\tP_D]=\BE_D[P_D\tP_D]$.

For the matrix statement, say the bigon is on the positions $p,p+1$ with $a$ before $b$, so that $A_a=B_b=C_p$ and $B_a=A_b=C_{p+1}$ and, choosing the base edges of $C_p,C_{p+1}$ outside the bigon (Lemma~\ref{lem.choices} (i)), $\mathrm{Pre}_{C_{p+1}}(o_b)=\mathrm{Pre}_{C_{p+1}}(u_a)\cup\{a\}$ and $\mathrm{Pre}_{C_p}(u_b)=\mathrm{Pre}_{C_p}(o_a)\cup\{a\}$. By \eqref{eq.Krow} the sum of the rows $a$ and $b$ of $K'(x)$ is $\d_a+\d_b+J_{C_p}(o_a)-J_{C_{p+1}}(u_a)+J_{C_{p+1}}(o_b)-J_{C_p}(u_b)=\d_a+\d_b+(\text{jump of $a$ on $C_{p+1}$})-(\text{jump of $a$ on $C_p$})=\d_a+\d_b-\d_a-x_a\d_a=\d_b-x_a\d_a$. The rows $a$ and $b$ of $Y'(x)$ are those of $Y'$ (Corollary~\ref{cor.seifert} (ii)), namely $(0,\dots,0,1,1)$ and $(0,\dots,0,0,1)$ in the stated basis, so the row of $\gamma_h$ in $M_{D',\hate}(x)=Y'^{\T}K'(x)Y'(x)$ is $(\d_b-x_a\d_a)Y'(x)=(0,\dots,0,-x_a,1-x_a)$. The other basis cycles do not run through $b$, and the row $b$ of $Y'(x)$ vanishes in their columns, so for the block of these cycles only the entries $K'(x)_{jk}$ with $j,k\ne b$ matter; these are the jump sums of $D_a$, since removing $b$ from the orders along $C_p$ and $C_{p+1}$ gives the orders of $D_a$, and $Y'(x)$ restricted to these rows and columns is $Y_{D_a}(x)$, because $\CC'(x)E_{\CT}=\CC_{D_a}(x)E_{\CT}$ and the column $b$ of $\CC'(x)$ is paired with the vanishing row $b$ of $Y'$.
\end{proof}

In the classical picture a hollow handle enlarges the Seifert matrix to $\bigl(\begin{smallmatrix}V&\xi&0\\0&*&1\\0&0&0\end{smallmatrix}\bigr)$, up to conventions, and $\det(tV'-V'^{\T})$ is $\det(tV-V^{\T})$ times a unit: the handle is invisible to the Alexander polynomial. In the deformed picture the handle leaves a trace, $(x_ax_b-1)R(x)$, which is killed only by the coupling with the second flavor. This is the mechanism observed on examples in Section~\ref{ssec.beyond}, by which the determinant formula fails to see the genus of noncanonical Seifert surfaces term by term.

\begin{proposition}[the braid relation]\label{prop.braidrelation}
Let $D'$ be obtained from $D$ by replacing a subword $\sigma_i\sigma_{i+1}\sigma_i$ by $\sigma_{i+1}\sigma_i\sigma_{i+1}$, the cut edge being outside the subword; let $a,b,c$ be the three crossings of the subword of $D$ in their order and $a',b',c'$ those of $D'$, matched by the pairs of strands which cross, so that $a\leftrightarrow c'$, $b\leftrightarrow b'$, $c\leftrightarrow a'$, and give matched crossings the same variables. Then
\[
   P_{D',\hate}=P_{D,\hate}+(x_a-x_b)\,R,\qquad \tP_{D',\hate}=\tP_{D,\hate}+(y_b-y_c)\,\tilde R,
\]
and for the inverse relation $\sigma_i^{-1}\sigma_{i+1}^{-1}\sigma_i^{-1}\to\sigma_{i+1}^{-1}\sigma_i^{-1}\sigma_{i+1}^{-1}$ the same holds with $(x_b-x_c)$ and $(y_a-y_b)$.
\end{proposition}
\begin{proof}
With $\rho_1(x)=\rho(x)\oplus1$ and $\rho_2(x)=1\oplus\rho(x)$ acting on the three strands, a direct computation gives
\[
   \rho_1(x_c)\rho_2(x_b)\rho_1(x_a)=\rho_2(x_a)\rho_1(x_b)\rho_2(x_c)\qquad\text{whenever }x_a=x_b ,
\]
a deformed Yang--Baxter equation (for $x_a=x_b=x_c$ it is the braid relation for the Burau representation; the two sides differ for generic $x_a,x_b$). Indeed, both products have the second row $(1-x_a,\,x_a,\,0)$ and the third row $(1,0,0)$; the first row of the left-hand side is $\bigl(1-x_a+x_ax_c-x_bx_c,\ x_a(1-x_c),\ x_bx_c\bigr)$ and that of the right-hand side is $\bigl(1-x_b,\ x_b(1-x_c),\ x_bx_c\bigr)$. So the difference of the two sides is $(x_a-x_b)(1-x_c)$ times the matrix whose first row is $(-1,1,0)$ and whose other rows vanish. The left-hand side is the product of the blocks of $\sigma_i\sigma_{i+1}\sigma_i$ with the variables $x_a,x_b,x_c$ in the order of the word, and the right-hand side that of $\sigma_{i+1}\sigma_i\sigma_{i+1}$ with the matched variables. Hence $P_{D',\hate}=P_{D,\hate}$ at $x_a=x_b$, and the difference is divisible by $x_a-x_b$. For the second flavor the blocks are transposed: the product for $D$ is the transpose of $\rho_1(y_a)\rho_2(y_b)\rho_1(y_c)$ and that for $D'$ the transpose of $\rho_2(y_c)\rho_1(y_b)\rho_2(y_a)$, which agree when $y_b=y_c$ by the same identity. For the inverse relation the same computation with $\rho^-$ gives the identities under $x_b=x_c$, resp.\ $y_a=y_b$.
\end{proof}

\begin{remark}\label{rem.markovetc}
(i) \emph{Markov stabilization} adds a strand and a crossing $k_0=\sigma_n^{\pm1}$, hence a disk $D_{n+1}$ and a band, and does not change $\nu$. The new band is a tree edge, the row of $C_{n+1}$ in $N_{D',\hate}(x)$ has the single entry $-1$ (if $\ve_{k_0}=1$, when $B_{k_0}=C_{n+1}$) or $x_{k_0}$ (if $\ve_{k_0}=-1$) in the column $k_0$, and expanding along it gives $P_{D',\hate}=\pm x_{k_0}^{\ind{\ve_{k_0}=-1}}P_{D,\hate}$. The sign is $+$, by Corollary~\ref{cor.seifert}~(iii): the new band lies on no cycle, the deformed Seifert matrix is unchanged, and the monomial $x^{\CT}$ is multiplied by $x_{k_0}$ if $A_{k_0}=C_{n+1}$, that is, if $\ve_{k_0}=-1$. In the same way $\tP_{D',\hate}=y_{k_0}^{\ind{\ve_{k_0}=-1}}\tP_{D,\hate}$. Conjugation is a change of the cut edge, Lemma~\ref{lem.understrand}.

(ii) \emph{Band generators.} By a theorem of Rudolph \cite{Rudolph83}, every Seifert surface is isotopic to a braided surface in a wider sense, in which the bands may join nonadjacent disks, passing in front of the disks in between; the boundary is the closure of a product of \emph{band generators} $\sigma_i\sigma_{i+1}\cdots\sigma_{j-2}\sigma_{j-1}^{\pm1}\sigma_{j-2}^{-1}\cdots\sigma_i^{-1}$, one for each band, joining the disks at the positions $i<j$ (with our convention for the crossings, the strand of the band passes over the strands of the disks in between, going and returning), and the canonical surface of this closed braid is the braided surface with a hollow handle attached for each disk passed by a band (the two crossings of a band with the strand of an intermediate disk are the two bands of a handle). The Seifert form of such a surface can still be cut into contributions of its bands, by closing each band up through its two disks and through the outside of the surface and pushing it off to either side, and one can form the matrices of the linking numbers of a basis of cycles with these curves, with one pair of variables for each band, exactly as the matrices \eqref{eq.PVV} are formed for a canonical surface. One may ask whether the image under $\BE$ of the product of the two determinants so obtained is, up to sign, $\LG$ of the boundary. It is not: already for a surface with three disks and four bands, two of which join the outer disk to the inner one, the two sides differ, and the right-hand side has degrees outside $[-\nu,\nu]$; for a braided surface of genus one with four disks whose boundary has Links--Gould polynomial $1$, the right-hand side is a Laurent polynomial with fourteen terms. The crossings which an expanded band generator makes with the intermediate strands carry variables of their own in \eqref{eq.intro.formula}, and these cannot be dispensed with: by Proposition~\ref{prop.hollowhandle} their contribution is a cancellation inside $\BE_D$, not a feature of the surface. We come back to this in Section~\ref{ssec.beyond}.
\end{remark}
\section{The \texorpdfstring{$\ve$}{epsilon}-degree: a diagrammatic genus bound and the top coefficient}\label{sec.degree}

Two of the three-dimensional properties of $\LG$ recalled in the introduction are statements about degrees. Let $g(L)$ be the genus of $L$ and $\mu$ its number of components, and give the monomials of $\BZ[t_0^{\pm1},t_1^{\pm1}]$ the degree
\begin{equation}\label{eq.KTdegree}
   \deg\bigl(t_0^at_1^b\bigr)=a-b
\end{equation}
of \cite{Kohli-Tahar}, so that $\mathrm{span}\,\LG(L)$ is the difference between the largest and the smallest degree of a monomial of $\LG(L)$. Kohli and Tahar proved that $\mathrm{span}\,\LG(K)\le4g(K)$ for every knot $K$, and L\'opez Neumann and van der Veen that $\LG(L)$ has degrees in $[-(2g+\mu-1),\,2g+\mu-1]$ for every link and that the coefficient of the extreme degree is a single monomial, with coefficient one, when $L$ is fibered \cite{Kohli-Tahar, LNVdV25}. (In the variables $p,q$ of \cite{LNVdV25}, $t_0=p^{-1}q$ and $t_1=pq$, so that $\deg(t_0^at_1^b)=-\deg_p(t_0^at_1^b)$, and the extreme degree $-(2g+\mu-1)$ is the top $p$-degree $2g+\mu-1=1-\chi(F)$ of a minimal genus Seifert surface $F$.) These are the counterparts of $\deg\Delta_L\le2g+\mu-1$ and of the monicity of the Alexander polynomial of a fibered link, and the proof of the bound in \cite{Kohli-Tahar}, on which \cite{LNVdV25} builds, is representation-theoretic: a Seifert surface of $L$ is presented as a tangle, and the degree is counted from the specific form of that tangle.

In this section we show what the deformed Seifert matrices say about these degrees. The classical argument for $\deg\Delta_L\le\nu$, where $\nu=\rk H_1(F_D)=c-s+1$ for the canonical Seifert surface $F_D$ of a diagram, is that $\Delta_L$ is the determinant of the $\nu\times\nu$ matrix $tV-V^{\T}$, each entry of which has degree at most one in $t$. We prove that the deformed Seifert matrices $M_{D,\hate}(x)$ and $M^b_{D,\hate}(y)$ have the same property with respect to the \emph{$\ve$-degree} $\deg_\ve x^S=\sum_{j\in S}\ve_j$ of the crossing variables, up to an explicit shift of each row and column by the levels $h$ of Section~\ref{ssec.level}, and that these shifts are the same for the two flavors.

Consequently (Theorem~\ref{thm.canonical}) $\LG(L)$ has degrees in $[-\nu,\nu]$ for every connected diagram $D$ of $L$: this is the bound of \cite{Kohli-Tahar, LNVdV25} with the genus of the canonical surface of a diagram in place of the genus of the link, exactly as $\deg\Delta_L\le c-s+1$ is the classical bound with the genus of $F_D$. When $\Delta_L$ has span $\nu$---for instance when $D$ is alternating or homogeneous---the bound is sharp and the extreme coefficients of $\LG$ are read off from the extremal monomials of the two determinants (Theorem~\ref{thm.tight}); for alternating links these monomials are the spanning trees of \cite{HK}, the extreme coefficient is a polynomial in $q$ with nonnegative coefficients, and it is a single monomial exactly when the link is fibered, which gives a combinatorial formula for the Hopf invariant appearing in \cite{LNVdV25} (Theorem~\ref{thm.alternating}).

Section~\ref{ssec.sectors} describes the monomials of $\tP_{D,\hate}$ as the ways of cutting $F_D$ into a disk along arcs of its Seifert disks, the $\ve$-degree being, up to an additive constant, minus the sum of the levels of the pieces (Proposition~\ref{prop.sectors}), and at the end of the section we explain why the determinant formula does not see the genus of other Seifert surfaces term by term. Section~\ref{sec.seifert} has shown that the deformed Seifert matrix at the Alexander point is exactly $tV_F-V_F^{\T}$ up to units on the columns, for every diagram (Theorem~\ref{thm.alexpoint}), and we shall see that the sums of the extremal coefficients of the edge polynomials are $\det V_F$ and $(-1)^\nu\det V_F$ (Corollary~\ref{cor.extremesums}).

\subsection{Degrees}\label{ssec.epsdegree}
We introduce the grading of the crossing variables on which this section is based, and we show how the coupling $\BE_D$ turns the $\ve$-degrees of the monomials of the two edge polynomials into the degrees \eqref{eq.KTdegree} of the monomials of $\LG$ (Lemma~\ref{lem.pq}). This reduces the degree bounds for $\LG$ to the problem of locating the $\ve$-degrees of the monomials of $P_{D,\hate}$ and of $\tP_{D,\hate}$.

Let $D$ be a connected diagram with crossings $1,\dots,c$ of signs $\ve_j$. For $S\subseteq\{1,\dots,c\}$ put $\ve(S)=\sum_{j\in S}\ve_j$, $x^S=\prod_{j\in S}x_j$ and $y^S=\prod_{j\in S}y_j$. We grade the rings $\Lam_c=\BZ[x_j^{\pm1}]$ and $\Lam'_c=\BZ[y_j^{\pm1}]$ by
\begin{equation}\label{eq.epsdegree}
   \deg_\ve x_j=\deg_\ve y_j=\ve_j,\qquad\text{so that }\ \deg_\ve x^S=\deg_\ve y^S=\ve(S),
\end{equation}
and we call a Laurent polynomial \emph{homogeneous of $\ve$-degree $\ell$} if all its monomials have $\ve$-degree $\ell$ (the zero polynomial is homogeneous of every degree). For a nonzero $f$ we write $\deg_\ve^-f\le\deg_\ve^+f$ for the smallest and the largest $\ve$-degree of a monomial of $f$. We write the edge polynomials in the monomial basis,
\begin{equation}\label{eq.monomials}
   P_{D,\hate}(x)=\sum_SP_S\,x^S,\qquad \tP_{D,\hate}(y)=\sum_T\tP_T\,y^T,\qquad P_S,\tP_T\in\BZ ,
\end{equation}
not to be confused with the coefficients $a_S,b_S$ of the expansion \eqref{eq.ab} around the free fermion point. The substitution $x_j\mapsto t^{\ve_j}$ sends $x^S$ to $t^{\ve(S)}$, so the $\ve$-degree is the degree in $t$ at the Alexander point: $P_{D,\hate}(t^{\ve})=\sum_\ell\bigl(\sum_{\ve(S)=\ell}P_S\bigr)t^\ell$, and the $t$-degrees occurring in $P_{D,\hate}(t^{\ve})$ are $\ve$-degrees of monomials of $P_{D,\hate}(x)$, the converse being false in general because of cancellations.

\begin{lemma}\label{lem.pq}
For all $S,T\subseteq\{1,\dots,c\}$, $\BE_D[x^Sy^T]$ is a single monomial, of degree $\ve(S)-\ve(T)$ in the sense of \eqref{eq.KTdegree}; in the variables $t_0=p^{-1}q$, $t_1=pq$,
\[
   \BE_D[x^Sy^T]=p^{\,\ve(T)-\ve(S)}\,q^{\,\ve(S\triangle T)},\qquad S\triangle T=(S\sm T)\cup(T\sm S).
\]
Consequently, if $P_{D,\hate}$ and $\tP_{D,\hate}$ are nonzero, the degrees of the monomials of $\LG(L)$ lie in $\bigl[\deg_\ve^-P_{D,\hate}-\deg^+_\ve\tP_{D,\hate},\ \deg^+_\ve P_{D,\hate}-\deg^-_\ve\tP_{D,\hate}\bigr]$, and the part of $\LG(L)$ of degree $\deg^-_\ve P_{D,\hate}-\deg^+_\ve\tP_{D,\hate}$, which is the coefficient of $p^{\deg^+_\ve\tP-\deg^-_\ve P}$ in $\LG(L;p,q)$, is
\begin{equation}\label{eq.topformula}
   (-1)^{\mu-1}\sum_{\ve(S)=\deg^-_\ve P,\ \ \ve(T)=\deg^+_\ve\tP}P_S\,\tP_T\;q^{\,\ve(S\triangle T)}\in\BZ[q^{\pm1}],
\end{equation}
and similarly at the other end, with $\deg^+_\ve P$ and $\deg^-_\ve\tP$.
\end{lemma}
\begin{proof}
By \eqref{eq.ED}, $\BE_D[x^Sy^T]=t_0^{\ve(S\sm T)}t_1^{\ve(T\sm S)}$, of degree $\ve(S\sm T)-\ve(T\sm S)=\ve(S)-\ve(T)$, and $t_0^at_1^b=p^{b-a}q^{a+b}$ with $a+b=\ve(S\sm T)+\ve(T\sm S)=\ve(S\triangle T)$. The rest follows from Theorem~\ref{thm.HK}: $\LG(L)=(-1)^{\mu-1}\sum_{S,T}P_S\tP_T\,\BE_D[x^Sy^T]$, and the monomials of degree $\deg^-_\ve P-\deg^+_\ve\tP$ can only come from the pairs with $\ve(S)=\deg^-_\ve P$ and $\ve(T)=\deg^+_\ve\tP$.
\end{proof}

The sum \eqref{eq.topformula} may vanish, in which case the degree of $\LG$ is smaller than the bound of the lemma; we come back to this in Section~\ref{ssec.beyond}. The first task is to locate the $\ve$-degrees of the two edge polynomials.

\subsection{The \texorpdfstring{$\ve$}{epsilon}-degrees of the deformed Seifert matrices}\label{ssec.degreesM}
The classical bound $\deg\Delta_L\le\nu$ holds because every entry of the matrix $tV-V^{\T}$ has degree at most one in $t$. In this subsection we prove the corresponding property of the deformed Seifert matrices: up to a shift which depends only on the row and on the column, the monomials of every entry have two possible $\ve$-degrees, which are adjacent (Proposition~\ref{prop.degreesM}). There are two ingredients. The entries of the deformed cycles are homogeneous (Lemma~\ref{lem.Yx}); and the lift of a cycle of the Seifert graph only meets the Seifert circles of this cycle, which have two adjacent levels (Lemma~\ref{lem.regions}). The theorem on the $\ve$-degrees of the edge polynomials (Theorem~\ref{thm.degrees}) and the genus bound (Theorem~\ref{thm.canonical}) follow.

We fix a spanning tree $\CT$ of $\Gamma_D$ and take for $Y$ the basis of fundamental cycles of $\CT$: for every crossing $i\notin\CT$, the column $\gamma_i$ of $Y$ is the cycle through the band $i$ and the path of $\CT$ from $B_i$ back to $A_i$, which has entry $1$ at $i$ and entries $0$ at the other crossings not in $\CT$. Thus the rows of $Y$ indexed by the crossings $i\notin\CT$ form the identity matrix, and by Corollary~\ref{cor.seifert} (ii) the same holds for $Y(x)$ and $Y(y)$: a vector in $\ker\CC(x)$ is determined by its coordinates at the crossings not in $\CT$. By Lemma~\ref{lem.regions}, the circles of $\gamma_i$ have the levels $H_i$ and $H_i+1$, where we write $H_i=H_{\gamma_i}$.

\begin{lemma}\label{lem.Yx}
For all crossings $k$ and $i\notin\CT$, the entry $Y(x)_{k,i}$ is homogeneous of $\ve$-degree $h(B_i)-h(B_k)$, and so is $Y(y)_{k,i}$.
\end{lemma}
\begin{proof}
The idea is to keep, of each entry of the deformed cycle, only the homogeneous component of the expected $\ve$-degree. We show that the vector so obtained still lies in the kernel of $\CC(x)$ and that it has the same coordinates as the deformed cycle at the crossings which are not in $\CT$; so it is the deformed cycle itself. Let $w=Y(x)e_i\in\ker\CC(x)$ and let $w'\in\Lam_c^c$ be the vector whose entry $w'_k$ is the homogeneous component of $w_k$ of $\ve$-degree $h(B_i)-h(B_k)$. For a circle $C\ne\hat C$, the equation $(\CC(x)w)_C=\sum_{A_k=C}x_kw_k-\sum_{B_k=C}w_k=0$ holds in each $\ve$-degree separately. Its homogeneous component of $\ve$-degree $h(B_i)-h(C)$ is $\sum_{A_k=C}x_kw'_k-\sum_{B_k=C}w'_k$: for $A_k=C$ we have $\deg_\ve x_k=\ve_k=h(B_k)-h(C)$, so the component of $x_kw_k$ of $\ve$-degree $h(B_i)-h(C)$ is $x_k$ times the component of $w_k$ of $\ve$-degree $h(B_i)-h(B_k)$, and for $B_k=C$ the component of $w_k$ of $\ve$-degree $h(B_i)-h(C)=h(B_i)-h(B_k)$ is $w'_k$. Hence $\CC(x)w'=0$. At the crossings not in $\CT$, $w$ has the entries $w_i=1$, which is homogeneous of $\ve$-degree $0=h(B_i)-h(B_i)$, and $w_k=0$ for $k\ne i$; so $w'$ has the same entries there, and $w'=w$ because a vector of $\ker\CC(x)$ is determined by these entries. The proof for $Y(y)$ is the same.
\end{proof}

\begin{proposition}\label{prop.degreesM}
Index the rows of $M_{D,\hate}(x)=Y^{\T}K(x)Y(x)$ by the fundamental cycles $\gamma$ and its columns by the crossings $i\notin\CT$, and the other way round for $M^b_{D,\hate}(y)=Y(y)^{\T}K^b(y)Y$. Then every monomial of the entry $M_{D,\hate}(x)_{\gamma,i}$, and every monomial of $M^b_{D,\hate}(y)_{i,\gamma}$, has $\ve$-degree $h(B_i)-h(C)$ for some Seifert circle $C$ of the cycle $\gamma$; in particular its $\ve$-degree is $h(B_i)-H_\gamma-1$ or $h(B_i)-H_\gamma$.
\end{proposition}
\begin{proof}
By \eqref{eq.Mentries}, $M_{D,\hate}(x)_{\gamma,i}=\sum_k\bigl(x_k\,z^\gamma_{u'_k}-z^\gamma_{u_k}\bigr)\,Y(x)_{k,i}$, where $z^\gamma$ is the lift of $\gamma$. By \eqref{eq.lift}, the lift of a cycle vanishes on the edges of the Seifert circles through which the cycle does not pass. So the term $x_k\,z^\gamma_{u'_k}\,Y(x)_{k,i}$ vanishes unless the circle $A_k$, on which the edge $u'_k$ lies, is a circle of $\gamma$, and by Lemma~\ref{lem.Yx} it is homogeneous of $\ve$-degree $\ve_k+h(B_i)-h(B_k)=h(B_i)-h(A_k)$; the term $z^\gamma_{u_k}\,Y(x)_{k,i}$ vanishes unless $B_k$, on which $u_k$ lies, is a circle of $\gamma$, and it is homogeneous of $\ve$-degree $h(B_i)-h(B_k)$. For the second flavor, $M^b_{D,\hate}(y)_{i,\gamma}=-\sum_kY(y)_{k,i}\bigl(y_k\,z^\gamma_{o_k}-z^\gamma_{o'_k}\bigr)$ by \eqref{eq.Mentries}. The term with $z^\gamma_{o_k}$ vanishes unless $A_k$, on which $o_k$ lies, is a circle of $\gamma$, and it is homogeneous of $\ve$-degree $h(B_i)-h(B_k)+\ve_k=h(B_i)-h(A_k)$; the term with $z^\gamma_{o'_k}$ vanishes unless $B_k$, on which $o'_k$ lies, is a circle of $\gamma$, and it is homogeneous of $\ve$-degree $h(B_i)-h(B_k)$. The last claim is Lemma~\ref{lem.regions} (ii).
\end{proof}

Thus, after multiplying the row $\gamma$ by a monomial of $\ve$-degree $H_\gamma+1$ and the column $i$ by a monomial of $\ve$-degree $-h(B_i)$, every entry of $M_{D,\hate}(x)$ is a sum of a homogeneous part of $\ve$-degree $0$ and of a homogeneous part of $\ve$-degree $1$, like the entries of $tV-V^{\T}$; and the same monomials work for $M^b_{D,\hate}(y)$, with rows and columns exchanged. For the closed $2$-braids all the fundamental cycles have the same two circles, all the shifts are equal, and this is the form \eqref{eq.twobraid}.

\begin{theorem}[the $\ve$-degrees of the edge polynomials]\label{thm.degrees}
Let $D$ be a connected diagram, $\hate$ an edge of $D$, $\nu=c-s+1$, $\CT$ a spanning tree of $\Gamma_D$ rooted at $\hat C$, and let
\[
   \delta_{D,\hate}=\ve(\CT')+\sum_{i\notin\CT}\bigl(h(B_i)-H_i\bigr),\qquad \CT'=\{j(C):\,C\ne\hat C,\ A_{j(C)}=C\},
\]
so that $\ve(\CT')$ is the $\ve$-degree of the monomial $x^{\CT}$ of Corollary~\ref{cor.seifert} and each term $h(B_i)-H_i$ is $0$ or $1$. Then every monomial of $P_{D,\hate}(x)$ and every monomial of $\tP_{D,\hate}(y)$ has $\ve$-degree in the interval
\[
   I_{D,\hate}=[\,\delta_{D,\hate}-\nu,\ \delta_{D,\hate}\,]
\]
of length $\nu$. In particular $\deg^+_\ve P_{D,\hate}-\deg^-_\ve\tP_{D,\hate}\le\nu$ and $\deg^+_\ve\tP_{D,\hate}-\deg^-_\ve P_{D,\hate}\le\nu$.
\end{theorem}
\begin{proof}
By Corollary~\ref{cor.seifert} (iii), $P_{D,\hate}(x)=x^{\CT}\det M_{D,\hate}(x)$, and the determinant is the sum of the products $\pm\prod_\gamma M_{D,\hate}(x)_{\gamma,\sigma(\gamma)}$ over the bijections $\sigma$ from the fundamental cycles to the crossings not in $\CT$. By Proposition~\ref{prop.degreesM}, every monomial of the product has $\ve$-degree between $\sum_\gamma(h(B_{\sigma(\gamma)})-H_\gamma-1)$ and $\sum_\gamma(h(B_{\sigma(\gamma)})-H_\gamma)$, and since $\sigma$ is a bijection, $\sum_\gamma h(B_{\sigma(\gamma)})=\sum_{i\notin\CT}h(B_i)$ does not depend on $\sigma$. So all monomials of $\det M_{D,\hate}(x)$ have $\ve$-degree in $[\,d-\nu,d\,]$ with $d=\sum_{i\notin\CT}(h(B_i)-H_i)$, and those of $P_{D,\hate}(x)$ in $[\,d+\ve(\CT')-\nu,\ d+\ve(\CT')\,]=I_{D,\hate}$. For $\tP_{D,\hate}(y)=y^{\CT}\det M^b_{D,\hate}(y)$ (Corollary~\ref{cor.seifertb}, with the same monomial $y^{\CT}$) the same computation applies, the entry $(i,\gamma)$ of $M^b_{D,\hate}(y)$ having the same $\ve$-degrees as the entry $(\gamma,i)$ of $M_{D,\hate}(x)$.
\end{proof}

The interval $I_{D,\hate}$ does not depend on the choice of the point at infinity used to define the levels, which changes all $h$ by the same constant (Remark~\ref{rem.level}); it may depend on the tree $\CT$ and it depends on $\hate$, but it contains the $\ve$-degrees of both edge polynomials for every choice. Lemma~\ref{lem.pq} and Theorem~\ref{thm.degrees} give:

\begin{theorem}[diagrammatic genus bound]\label{thm.canonical}
Let $D$ be a connected diagram of $L$ with $c$ crossings and $s$ Seifert circles, and let $\nu=c-s+1=2g(F_D)+\mu-1$ be the first Betti number of its canonical Seifert surface. Then all monomials of $\LG(L;t_0,t_1)$ have degree in $[-\nu,\nu]$ for the degree \eqref{eq.KTdegree}, i.e.\ $\deg_p\LG(L;p,q)\in[-\nu,\nu]$, and
\[
   \mathrm{span}\,\LG(L)\le2\nu=4g(F_D)+2\mu-2 .
\]
\end{theorem}
\begin{proof}
If one of the two edge polynomials vanishes, then $\LG(L)=0$ by Theorem~\ref{thm.HK}, and there is nothing to prove. Otherwise, by Lemma~\ref{lem.pq}, the degrees of the monomials of $\LG(L)$ lie between $\deg^-_\ve P_{D,\hate}-\deg^+_\ve\tP_{D,\hate}$ and $\deg^+_\ve P_{D,\hate}-\deg^-_\ve\tP_{D,\hate}$; by Theorem~\ref{thm.degrees} the first number is at least $-\nu$ and the second one is at most $\nu$. The degree in $p$ of a monomial is the opposite of its degree \eqref{eq.KTdegree}, since $t_0^at_1^b=p^{b-a}q^{a+b}$. Finally $\nu=2g(F_D)+\mu-1$, because $F_D$ is a connected orientable surface with $\mu$ boundary components which is made of $s$ disks and $c$ bands, so that $2-2g(F_D)-\mu=\chi(F_D)=s-c$.
\end{proof}

For instance, let $D$ be the standard diagram of the trefoil, the closure of $\sigma_1^3$, with the cut edge $\hate=u'_3=o_1$ and the spanning tree $\CT=\{1\}$ (Example~\ref{ex.trefoilflows}). Here $\nu=2$, the set $\CT'$ is empty and $h(B_i)-H_i=1$ for $i=2,3$, so that $\delta_{D,\hate}=2$ and $I_{D,\hate}=[0,2]$. Indeed $P_{D,\hate}=1-x_1+x_1x_2$ and $\tP_{D,\hate}=1-y_3+y_2y_3$, and the monomials of both polynomials have the $\ve$-degrees $0$, $1$ and $2$. By Lemma~\ref{lem.pq}, the coefficient of $p^{2}$ in $\LG$ comes from the single pair $S=\emptyset$, $T=\{2,3\}$, and it is $q^{\ve(\{2,3\})}=q^2$.

The bound of Theorem~\ref{thm.canonical} is the bound of Kohli--Tahar and L\'opez Neumann--van der Veen, $\deg_p\LG(L)\in[-(2g(L)+\mu-1),2g(L)+\mu-1]$, with $g(L)$ replaced by the genus of the canonical surface of an arbitrary diagram of $L$, in the same way as $\deg\Delta_L\le c-s+1$ is the classical bound $\deg\Delta_L\le2g(L)+\mu-1$ with the same replacement; the proof is the Seifert-matrix proof of the classical bound, with the crossing variables carried along. The bound of \cite{Kohli-Tahar, LNVdV25} is stronger, since $g(L)\le g(F_D)$ with strict inequality for many diagrams; we discuss in Section~\ref{ssec.beyond} why the stronger bound is not visible term by term in the determinant formula. At the Alexander point the deformed matrix is $tV_F-V_F^{\T}$ times a diagonal matrix of units (Theorem~\ref{thm.alexpoint}), so that the bound is the classical one in the strict sense, and whether it is attained is decided by $\det V_F$ (Corollary~\ref{cor.extremesums} below).

\subsection{Tight diagrams: sharpness and the top coefficient}\label{ssec.tight}
Call the connected diagram $D$ \emph{tight} if the Alexander polynomial of $L$ has span $\nu=c-s+1$; since $\mathrm{span}\,\Delta_L\le2g(L)+\mu-1\le\nu$, this means that the canonical Seifert surface of $D$ has minimal genus and that $\Delta_L$ detects the genus. Reduced alternating diagrams are tight (Crowell \cite{Crowell59} and Murasugi \cite{Murasugi58}), and so are the homogeneous diagrams of Cromwell \cite{Cromwell89}, which include the positive diagrams and the closures of homogeneous braids. For a tight diagram the two sides of Theorem~\ref{thm.degrees} meet.

\begin{theorem}[tight diagrams]\label{thm.tight}
Let $D$ be a tight connected diagram of $L$, let $\hate$ be an edge of $D$ and $I_{D,\hate}=[\delta-\nu,\delta]$ the interval of Theorem~\ref{thm.degrees}. Then:
\begin{enumerate}\renewcommand{\theenumi}{\roman{enumi}}
\item The $\ve$-degrees of the monomials of $P_{D,\hate}$ and of $\tP_{D,\hate}$ both fill $I_{D,\hate}$ exactly: $\deg^-_\ve P=\deg^-_\ve\tP=\delta-\nu$ and $\deg^+_\ve P=\deg^+_\ve\tP=\delta$. Moreover $\sum_{\ve(S)=\ell}P_S=\sum_{\ve(T)=\ell}\tP_T$ for every $\ell$, these sums being the coefficients of the polynomial
\[
   P_{D,\hate}(t^{\ve})=\tP_{D,\hate}(t^{\ve})=(-1)^{\mu-1}t^{m_{D,\hate}/2}\Delta_L(t)
\]
of Lemma~\ref{lem.alexpoly}, and they are nonzero at both ends $\ell=\delta-\nu$ and $\ell=\delta$ of $I_{D,\hate}$.
\item $\LG(L)$ has degrees exactly in $[-\nu,\nu]$: the bound of Theorem~\ref{thm.canonical} is attained, and $\mathrm{span}\,\LG(L)=4g(L)+2\mu-2$, which is the bound of \cite{Kohli-Tahar, LNVdV25}.
\item The coefficient of $p^{\nu}$ in $\LG(L;p,q)$, the top coefficient $\mathrm{top}(L,q)$ of \cite{LNVdV25}, is
\[
   \mathrm{top}(L,q)=(-1)^{\mu-1}\sum_{\ve(S)=\delta-\nu}\ \sum_{\ve(T)=\delta}P_S\,\tP_T\;q^{\,\ve(S\triangle T)},
\]
the sum running over the monomials of $P_{D,\hate}$ of smallest and of $\tP_{D,\hate}$ of largest $\ve$-degree, and $\mathrm{top}(L,1)$ is the square of the leading coefficient of $\Delta_L$ (the Alexander polynomial normalized so that $\LG(L;p,1)=\Delta_L(p)^2$). The coefficient of $p^{-\nu}$ is given by the same formula with the roles of smallest and largest exchanged.
\end{enumerate}
\end{theorem}
\begin{proof}
(i) The two determinants $P_{D,\hate}(t^{\ve})=\det A_{D,\hate}(t)$ and $\tP_{D,\hate}(t^{\ve})=\det B_{D,\hate}(t)$ are equal (Lemma~\ref{lem.level} (i)), and by Lemma~\ref{lem.alexpoly} they are equal to $\Delta_L(t)$ up to a unit. The $t$-degrees of the monomials of $P_{D,\hate}(t^{\ve})$ form a set of span $\mathrm{span}\,\Delta_L=\nu$ which is contained in the set of $\ve$-degrees of the monomials of $P_{D,\hate}(x)$, itself contained in $I_{D,\hate}$ by Theorem~\ref{thm.degrees}. An interval of length $\nu$ containing a set of span $\nu$ is the convex hull of that set; so the extreme $t$-degrees of $P_{D,\hate}(t^{\ve})$ are the ends of $I_{D,\hate}$, and so are the extreme $\ve$-degrees of $P_{D,\hate}(x)$. The same holds for $\tP_{D,\hate}$. The coefficient of $t^\ell$ in $P_{D,\hate}(t^{\ve})$ is $\sum_{\ve(S)=\ell}P_S$.

(ii) and (iii) By Lemma~\ref{lem.pq} and (i), the degrees of $\LG(L)$ lie in $[-\nu,\nu]$, and the part of degree $-\nu$, i.e.\ the coefficient of $p^{\nu}$, is the sum \eqref{eq.topformula} with $\deg^-_\ve P=\delta-\nu$ and $\deg^+_\ve\tP=\delta$. At $q=1$ this sum is $(-1)^{\mu-1}\bigl(\sum_{\ve(S)=\delta-\nu}P_S\bigr)\bigl(\sum_{\ve(T)=\delta}\tP_T\bigr)\ne0$ by (i), so the coefficient of $p^\nu$ is a nonzero polynomial in $q$ and $\LG(L)$ has a monomial of degree $-\nu$; by symmetry (or by the same argument at the other end) it has one of degree $\nu$. Finally $\LG(L;p,1)=\LG(L;p^{-1},p)=\Delta_L(p)^2$ by the identity recalled in the introduction and $\Delta_L(p^{-1})^2=\Delta_L(p)^2$ (the Alexander polynomial being symmetric up to sign), and the coefficient of $p^{\nu}$ in $\Delta_L(p)^2$ is the square of the leading coefficient of $\Delta_L$, the span of $\Delta_L$ being $\nu$. Finally $\nu=2g(L)+\mu-1$, because $\nu=\mathrm{span}\,\Delta_L\le2g(L)+\mu-1\le2g(F_D)+\mu-1=\nu$.
\end{proof}

The exact identity at the Alexander point (Theorem~\ref{thm.alexpoint}) determines the sums of the extremal coefficients of every diagram, whether it is tight or not.

\begin{corollary}[extremal coefficients and fiberedness]\label{cor.extremesums}
Let $D$ be a connected diagram of $L$, $\hate$ an edge of $D$, $\CT$ a spanning tree of $\Gamma_D$ and $I_{D,\hate}=[\delta-\nu,\delta]$ the interval of Theorem~\ref{thm.degrees}. Then
\[
   \sum_{\ve(S)=\delta}P_S=\det V_F,\qquad \sum_{\ve(S)=\delta-\nu}P_S=(-1)^\nu\det V_F,
\]
where $V_F$ is the Seifert matrix of $F_D$ in the basis of fundamental cycles of $\CT$, and the same holds for the coefficients $\tP_T$ of $\tP_{D,\hate}$, with the same signs. Consequently the coefficient of $p^{\nu}$ in $\LG(L;p,q)$ takes the value $(\det V_F)^2$ at $q=1$, and so does the coefficient of $p^{-\nu}$. In particular: if $\det V_F\ne0$, the bound of Theorem~\ref{thm.canonical} is attained, $\deg_p\LG(L)=\nu$, and $F_D$ has minimal genus; if $F_D$ is a fiber of $L$, for instance if $D$ is the closure of a homogeneous braid \cite{Stallings78}, then $\det V_F=\pm1$ and the sums of the extremal coefficients of both edge polynomials are $\pm1$; and if $D$ is loose, i.e.\ $\mathrm{span}\,\Delta_L<\nu$, then $\det V_F=0$ and all four sums vanish.
\end{corollary}
\begin{proof}
The circles of a fundamental cycle $\gamma_i$ have the levels $H_i$ and $H_i+1$ (Lemma~\ref{lem.regions} (ii)) and $h(B_i)=h(A_i)+\ve_i$, so $h(B_i)-H_i=\ind{\ve_i=1}$, and $\delta-\nu=\ve(\CT')+\nu_+-\nu=\ve(\CT')-\nu_-$ with $\nu_\pm$ the numbers of positive and negative crossings not in $\CT$. By Theorem~\ref{thm.alexpoint}, $P_{D,\hate}(t^{\ve})=t^{\delta-\nu}\det(tV_F-V_F^{\T})$, and $\det(tV_F-V_F^{\T})=\det(V_F)\,t^\nu+\dots+(-1)^\nu\det(V_F)$. The coefficient of $t^\ell$ in $P_{D,\hate}(t^{\ve})$ is $\sum_{\ve(S)=\ell}P_S$ (Section~\ref{ssec.epsdegree}), which gives the two sums for $P_{D,\hate}$; for $\tP_{D,\hate}$ they are the same, because $\tP_{D,\hate}(t^{\ve})=P_{D,\hate}(t^{\ve})$ (Lemmas~\ref{lem.level} (i) and~\ref{lem.alexpoly}) and the $\ve$-degrees of both polynomials lie in $I_{D,\hate}$ (Theorem~\ref{thm.degrees}).

We turn to the coefficient of $p^{\nu}$. By Lemma~\ref{lem.pq}, the coefficient of $p^\nu$ in $\LG(L;p,q)$ is $(-1)^{\mu-1}\sum P_S\tP_T\,q^{\ve(S\triangle T)}$ over $\ve(S)=\delta-\nu$ and $\ve(T)=\delta$, whose value at $q=1$ is the product of the two sums. This product is $(\det V_F)^2$: the coefficients of $t^{\delta-\nu}$ and of $t^{\delta}$ in $P_{D,\hate}(t^{\ve})$ are $(-1)^{\nu}\det V_F$ and $\det V_F$, the second being also the coefficient of $t^{\delta}$ in $\tP_{D,\hate}(t^{\ve})=P_{D,\hate}(t^{\ve})$, so the value is $(-1)^{\mu-1+\nu}(\det V_F)^2=(\det V_F)^2$, since $\nu=2g(F_D)+\mu-1$; the same computation applies to the coefficient of $p^{-\nu}$.

It remains to prove the three consequences. If $\det V_F\ne0$ this is nonzero, so $\deg_p\LG(L)=\nu$, and $\nu\le2g(L)+\mu-1\le2g(F_D)+\mu-1=\nu$ by the bound of \cite{Kohli-Tahar, LNVdV25}. If $F_D$ is a fiber, then $V_F$ is unimodular, $V_F^{-1}V_F^{\T}$ being, up to conventions, the matrix of the monodromy of the fibration on $H_1(F_D)$ \cite{BurdeZieschang}; and the closure of a homogeneous braid is fibered with fiber $F_D$ by Stallings \cite{Stallings78}. Finally $\mathrm{span}\,\Delta_L=\mathrm{span}\det(tV_F-V_F^{\T})$ is $\nu$ if $\det V_F\ne0$ and smaller otherwise.
\end{proof}

\begin{remark}[loose diagrams]\label{rem.loose}
If $D$ is not tight, the sum \eqref{eq.topformula} at the ends of $I_{D,\hate}$ may vanish, and then the degree of $\LG(L)$ is smaller than the bound of Theorem~\ref{thm.canonical}. The closure of the $3$-braid $\sigma_1^3\sigma_2\sigma_1^{-1}\sigma_2$ is a diagram of the knot $5_2$ with $c=6$, $s=3$, $\nu=4$, while $\Delta_{5_2}=2t^2-3t+2$ has span $2$: for a suitable cut edge the $\ve$-degrees of both edge polynomials span the interval $[1,4]$, so that $\deg^+_\ve\tP-\deg^-_\ve P=3$, but the sum \eqref{eq.topformula} vanishes, and $\LG(5_2;p,q)$ has $p$-degree $2=2g(5_2)$. Even for fibered knots the sum can vanish on a loose diagram: for the closed $3$-braid diagram $\sigma_1^3\sigma_2^{-1}\sigma_1^{-3}\sigma_2^{-1}$ of $8_{20}$, $\nu=6$, some cut edges give $\deg^+_\ve\tP-\deg^-_\ve P=5$, and $\deg_p^+\LG=4=2g(8_{20})$. On loose diagrams the extremal monomials need not be sign-coherent either: for a disk with four bands bounding a fibered knot of genus $2$ with a loose $16$-crossing diagram ($\nu=8$, $\mathrm{span}\,\Delta=4$), the monomials of $P_{D,\hate}$ of smallest $\ve$-degree are nine, with coefficients $\pm1$ adding up to $-1$, while the top coefficient of $\LG$ is a single monomial, as it must be by \cite{LNVdV25}. And the genus bound itself need not be sharp: the Kinoshita--Terasaka knot of Remark~\ref{rem.KT} has genus $2$ and $\deg_p\LG\in[-3,3]$. In these examples the vanishing of \eqref{eq.topformula} is a cancellation between monomials of $P_{D,\hate}$ of the same $\ve$-degree and opposite signs, and between monomials of $\tP_{D,\hate}$; on tight diagrams such cancellations cannot occur at the extreme $\ve$-degrees, by (i) of the theorem. For a loose diagram the cancellation is complete at both ends: the sums $\sum_{\ve(S)=\delta}P_S$ and $\sum_{\ve(S)=\delta-\nu}P_S$, and the corresponding sums for $\tP_{D,\hate}$, are $\pm\det V_F=0$ (Corollary~\ref{cor.extremesums}), so that \eqref{eq.topformula} vanishes at $q=1$.
\end{remark}

\subsection{Murasugi sums along Seifert circles}\label{ssec.murasugi}
The leading coefficient of the Alexander polynomial is multiplicative under Murasugi sums of minimal genus Seifert surfaces, because the Seifert matrix of a Murasugi sum is block triangular, while the polynomial itself is not multiplicative. Anghel, L\'opez Neumann and van der Veen have recently proved the analogous statement for $\LG$: the top coefficient of $\LG$, a Laurent polynomial in $q$, is multiplicative under Murasugi sums of minimal genus Seifert surfaces \cite{ALvdV26}, extending the plumbing case of \cite{LNVdV25}; their proof uses the quantum supergroup. The deformed Seifert matrices give a short proof of this statement for the Murasugi sums which are visible on a diagram, by the same mechanism as for the Alexander polynomial: the extremal homogeneous components of the edge polynomials are determinants of block triangular matrices.

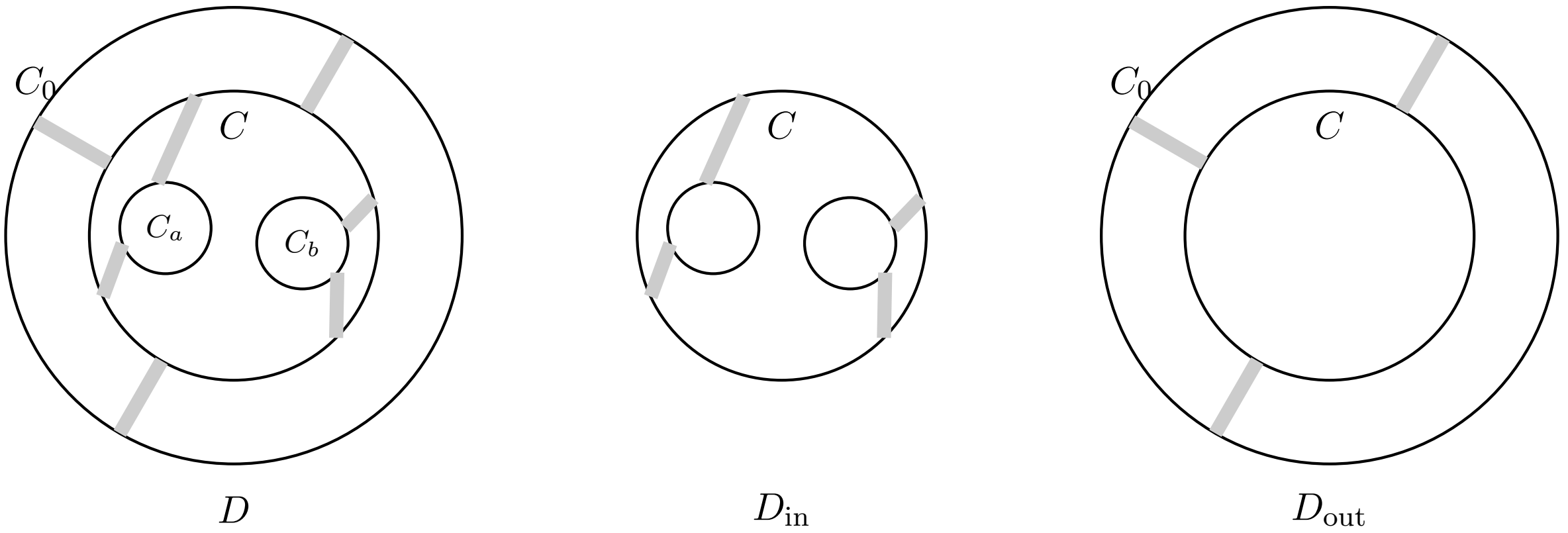
\begin{figure}[ht]
\centering
\begin{tikzpicture}[>=Stealth,scale=0.76,every node/.style={font=\small}]
\begin{scope}
\draw[thick] (0,0) circle (3); \draw[thick] (0,0) circle (1.9); \draw[thick] (-0.9,0.1) circle (0.6); \draw[thick] (0.9,-0.1) circle (0.6);
\node at (-2.6,2.0) {$C_0$}; \node at (0,1.45) {$C$}; \node[font=\scriptsize] at (-0.9,0.1) {$C_a$}; \node[font=\scriptsize] at (0.9,-0.1) {$C_b$};
\foreach \ang in {60,150,240}{ \draw[line width=4pt,gray!40] ({1.9*cos(\ang)},{1.9*sin(\ang)}) -- ({3*cos(\ang)},{3*sin(\ang)}); }
\foreach \ang in {100,200}{ \draw[line width=4pt,gray!40] ({-0.9+0.6*cos(\ang)},{0.1+0.6*sin(\ang)}) -- ({1.9*cos(\ang+5)},{1.9*sin(\ang+5)}); }
\foreach \ang in {-40,20}{ \draw[line width=4pt,gray!40] ({0.9+0.6*cos(\ang)},{-0.1+0.6*sin(\ang)}) -- ({1.9*cos(\ang-5)},{1.9*sin(\ang-5)}); }
\node at (0,-3.6) {$D$};
\end{scope}
\begin{scope}[xshift=7.2cm]
\draw[thick] (0,0) circle (1.9); \draw[thick] (-0.9,0.1) circle (0.6); \draw[thick] (0.9,-0.1) circle (0.6);
\node at (0,1.45) {$C$};
\foreach \ang in {100,200}{ \draw[line width=4pt,gray!40] ({-0.9+0.6*cos(\ang)},{0.1+0.6*sin(\ang)}) -- ({1.9*cos(\ang+5)},{1.9*sin(\ang+5)}); }
\foreach \ang in {-40,20}{ \draw[line width=4pt,gray!40] ({0.9+0.6*cos(\ang)},{-0.1+0.6*sin(\ang)}) -- ({1.9*cos(\ang-5)},{1.9*sin(\ang-5)}); }
\node at (0,-3.6) {$D_{\rm in}$};
\end{scope}
\begin{scope}[xshift=14.4cm]
\draw[thick] (0,0) circle (3); \draw[thick] (0,0) circle (1.9);
\node at (-2.6,2.0) {$C_0$}; \node at (0,1.45) {$C$};
\foreach \ang in {60,150,240}{ \draw[line width=4pt,gray!40] ({1.9*cos(\ang)},{1.9*sin(\ang)}) -- ({3*cos(\ang)},{3*sin(\ang)}); }
\node at (0,-3.6) {$D_{\rm out}$};
\end{scope}
\end{tikzpicture}
\caption{The decomposition of a diagram along a Seifert circle $C$ (Section~\ref{ssec.murasugi}): the Seifert circles of $D$ with the bands of $F_D$ (gray), and the diagrams $D_{\rm in}$ and $D_{\rm out}$, whose canonical surfaces are Murasugi summed along $D_C$ to give $F_D$.}
\label{fig.murasugi}
\end{figure}

Let $D$ be a connected diagram and $C$ a Seifert circle of $D$ with crossings on both sides. The Seifert circles other than $C$ lie inside or outside $C$, and no crossing joins a circle inside $C$ to a circle outside $C$. Let $D_{\rm in}$ be the diagram formed by $C$, the circles inside $C$ and the crossings between these circles, and $D_{\rm out}$ the diagram formed by $C$, the circles outside $C$ and the crossings between them (Figure~\ref{fig.murasugi}); both are connected, since $\Gamma_D$ is the union of $\Gamma_{D_{\rm in}}$ and $\Gamma_{D_{\rm out}}$, which meet in the vertex $C$ only, and $\nu(D)=\nu(D_{\rm in})+\nu(D_{\rm out})$. Let $L_{\rm in}$ and $L_{\rm out}$ be the links they represent.

The canonical surface $F_D$ is the Murasugi sum of $F_{D_{\rm in}}$ and $F_{D_{\rm out}}$ along the disk $D_C$ \cite{Murasugi63, Gabai83}. Indeed, in the construction of Section~\ref{ssec.statesurface}, let $W_C$ be the set of the points which lie over the disk $\Omega_C$ bounded by $C$ and above the disk $D_C$ (Figure~\ref{fig.twosides}~(a)). The part of $F_D$ which comes from $D_{\rm in}$, that is, the disks of the circles inside $C$ and the bands of the crossings inside $C$, lies in $W_C$, apart from the arcs along which these bands are attached to $D_C$: the disks of the inner circles lie between $D_C$ and the plane $P$, and the bands lie over $\Omega_C$, above $P$. The rest of $F_D$ lies outside $W_C$: the disks of the circles outside $C$ and the bands of the crossings outside $C$ lie over the complement of $\Omega_C$, and the disks of the circles which enclose $C$ lie below $D_C$. So the boundary of $W_C$, cut off at a large height, is a sphere which meets $F_D$ in the disk $D_C$ and separates the two parts: $F_D$ is the union of the canonical surface $F_{D_{\rm in}}$, inside this sphere, and of $F_{D_{\rm out}}$, outside it, glued along $D_C$. For the surface $F^\uparrow_D$ above the plane the same holds with the two sides of the disk exchanged: the part which comes from $D_{\rm in}$ lies below the disk of $C$ (Section~\ref{ssec.otherside} and Figure~\ref{fig.twosides}~(b)). By Gabai's theorem \cite{Gabai83}, $F_D$ has minimal genus if and only if $F_{D_{\rm in}}$ and $F_{D_{\rm out}}$ have.

For a connected diagram $D$ of $L$ let $\tau_D(q)\in\BZ[q^{\pm1}]$\label{not.tauD} be the coefficient of $p^{\nu(D)}$ in $\LG(L;p,q)$, the polynomial \eqref{eq.topformula} at the upper end of $I_{D,\hate}$, which does not depend on $\hate$. By the genus bound of \cite{Kohli-Tahar, LNVdV25}, $\tau_D$ is the top coefficient of $\LG(L)$ in the sense of \cite{ALvdV26} if $F_D$ has minimal genus, and $\tau_D=0$ otherwise.

\begin{theorem}[Murasugi sums along a Seifert circle]\label{thm.murasugi}
Let $D$ be a connected diagram of $L$ and $C$ a Seifert circle of $D$ with crossings on both sides. Then
\[
   \tau_D(q)=\tau_{D_{\rm in}}(q)\,\tau_{D_{\rm out}}(q).
\]
In particular, if $F_D$ has minimal genus, the top coefficient of $\LG(L)$ is the product of the top coefficients of $\LG(L_{\rm in})$ and $\LG(L_{\rm out})$.
\end{theorem}
\begin{proof}
The proof is modeled on that of the multiplicativity of the leading coefficient of the Alexander polynomial, in which the Seifert matrix of a Murasugi sum is block triangular. We choose bases adapted to the decomposition, and we show in turn: that the deformed Seifert matrix of $D$ has the deformed Seifert matrices of $D_{\rm in}$ and of $D_{\rm out}$ as its diagonal blocks; that the entries of each of the two other blocks are homogeneous, those of one block having the upper and those of the other block the lower of the two possible $\ve$-degrees; that the two extremal homogeneous components of the matrix are therefore block triangular, so that the extremal components of the edge polynomials of $D$ are the products of those of $D_{\rm in}$ and of $D_{\rm out}$; and that the top coefficient factors accordingly.

Choose the cut edge $\hate$ on $C$, spanning trees $\CT_{\rm in}$ and $\CT_{\rm out}$ of $\Gamma_{D_{\rm in}}$ and $\Gamma_{D_{\rm out}}$ rooted at $C$, and $\CT=\CT_{\rm in}\cup\CT_{\rm out}$, a spanning tree of $\Gamma_D$ rooted at $\hat C=C$ whose fundamental cycles are those of $\CT_{\rm in}$ and of $\CT_{\rm out}$. We use the basis $Y$ of fundamental cycles throughout; $\nu_{\rm in}$, $\nu_{\rm out}$, $\CT'_{\rm in}$, $\CT'_{\rm out}$, $\delta_{\rm in}$, $\delta_{\rm out}$ denote the quantities of Theorem~\ref{thm.degrees} for $D_{\rm in}$ and $D_{\rm out}$ with these choices, and we call a band or a cycle \emph{inner} or \emph{outer} according to the subdiagram it belongs to.

\emph{Block structure.} The order $\prec_C$ of the crossings of $C$ in $D$, from $\hate$, restricts to the orders $\prec_C$ of $D_{\rm in}$ and of $D_{\rm out}$, whose base edges are the edges of $C$ containing $\hate$; the orders along the other circles are the same in $D$ as in the subdiagram containing them. Hence $K(x)_{jk}$ is the entry of $K_{D_{\rm in}}(x)$ when $j,k$ are inner, that of $K_{D_{\rm out}}(x)$ when they are outer, and for $j$ inner and $k$ outer, or the other way round, it vanishes unless $j$ and $k$ both lie on $C$, in which case it is $s_C(j)\bigl(x_k\ind{A_k=C}-\ind{B_k=C}\bigr)\ind{k\prec_Cj}$ by \eqref{eq.Krow}. The matrix $\CC(x)E_{\CT}$, whose rows are the circles other than $C$ and whose columns are the tree bands, is the direct sum of $\CC_{D_{\rm in}}(x)E_{\CT_{\rm in}}$ and $\CC_{D_{\rm out}}(x)E_{\CT_{\rm out}}$, the tree band of a circle inside $C$ having its entries in the rows of that circle and of its parent, which is inside $C$ or is $C$ itself; so by \eqref{eq.Yx} the deformed cycle $Y(x)e_i$ of an inner crossing $i\notin\CT$ is the deformed cycle of $i$ in $D_{\rm in}$, supported on the inner bands, and likewise for the outer crossings. Writing $M(x)=Y^{\T}K(x)Y(x)$ in blocks according to inner and outer cycles and crossings,
\begin{gather*}
   M(x)=\begin{pmatrix}M_{D_{\rm in},\hate}(x)&M_{\rm io}(x)\\ M_{\rm oi}(x)&M_{D_{\rm out},\hate}(x)\end{pmatrix},\\
   M_{\rm io}(x)_{\gamma,i}=\sum_{k\text{ outer, on }C}\bigl(x_k\ind{A_k=C}-\ind{B_k=C}\bigr)\,Y(x)_{k,i}\sum_{j\text{ inner, on }C,\ k\prec_Cj}s_C(j)\,w_j(\gamma)
\end{gather*}
for an inner cycle $\gamma$ and an outer crossing $i\notin\CT$, and symmetrically for $M_{\rm oi}(x)$.

\emph{Degrees.} By Lemma~\ref{lem.Yx}, $Y(x)_{k,i}$ is homogeneous of $\ve$-degree $h(B_i)-h(B_k)$, and $x_k\ind{A_k=C}-\ind{B_k=C}$ is homogeneous of $\ve$-degree $h(B_k)-h(C)$ in both cases, since $\ve_k=h(B_k)-h(A_k)$; so $M_{\rm io}(x)_{\gamma,i}$ is homogeneous of $\ve$-degree $h(B_i)-h(C)$, and vanishes unless $\gamma$ passes through $C$. An inner cycle through $C$ lies in the region $\Omega_{\rm in}$ directly inside $C$, of which $C$ is the outer circle, and an outer cycle through $C$ in the region $\Omega_{\rm out}$ directly outside $C$, of which $C$ is an inner circle; by the proof of Lemma~\ref{lem.regions}, $h(C)=H_{\Omega_{\rm in}}+\ind{s(C)=-1}=H_{\Omega_{\rm out}}+\ind{s(C)=1}$. Hence, in the terms of Proposition~\ref{prop.degreesM}, if $C$ is clockwise, the entries of $M_{\rm io}(x)$ have the upper of the two possible $\ve$-degrees, $h(B_i)-H_\gamma$, and those of $M_{\rm oi}(x)$ the lower one, $h(B_i)-H_\gamma-1$; if $C$ is counterclockwise it is the other way round.

\emph{Extremal components.} Write every entry of $M(x)$ as the sum $M(x)_{\gamma,i}=M^{[0]}_{\gamma,i}+M^{[1]}_{\gamma,i}$ of its homogeneous components of $\ve$-degrees $h(B_i)-H_\gamma-1$ and $h(B_i)-H_\gamma$. As in the proof of Theorem~\ref{thm.degrees}, every monomial of $\det M(x)$ has $\ve$-degree in $[d-\nu,d]$ with $d=\sum_{i\notin\CT}h(B_i)-\sum_\gamma H_\gamma$, and the components of $\det M(x)$ of $\ve$-degrees $d-\nu$ and $d$ are $\det M^{[0]}$ and $\det M^{[1]}$. If $C$ is clockwise, $M^{[0]}$ has vanishing block $(\mathrm{in},\mathrm{out})$ and $M^{[1]}$ has vanishing block $(\mathrm{out},\mathrm{in})$; if $C$ is counterclockwise, the other blocks vanish. In both cases $M^{[0]}$ and $M^{[1]}$ are block triangular, and their diagonal blocks are the matrices $M^{[0]}$ and $M^{[1]}$ of $D_{\rm in}$ and of $D_{\rm out}$: the levels of the circles of $D_{\rm in}$, resp.\ $D_{\rm out}$, computed in $D$ differ from those computed in the subdiagram by a common constant (the circles inside $C$ are contained in the same circles in both diagrams, except for the circles containing $C$, which contain all of them; the circles outside $C$ are contained in the same circles in both), so that the components $M^{[0]},M^{[1]}$ of the diagonal blocks are the same in $D$ as in the subdiagrams. Hence $\det M^{[0]}=\det M^{[0]}_{\rm in}\det M^{[0]}_{\rm out}$ and $\det M^{[1]}=\det M^{[1]}_{\rm in}\det M^{[1]}_{\rm out}$.

\emph{The edge polynomials.} Since $P_{D,\hate}=x^{\CT}\det M(x)$ with $x^{\CT}=x^{\CT_{\rm in}}x^{\CT_{\rm out}}$ (Corollary~\ref{cor.seifert}~(iii), applied to the three diagrams), and since $\delta=\delta_{\rm in}+\delta_{\rm out}$ and $\nu=\nu_{\rm in}+\nu_{\rm out}$, the homogeneous components of $P_{D,\hate}$ of $\ve$-degrees $\delta-\nu$ and $\delta$ are the products of the corresponding components of $P_{D_{\rm in},\hate}$ and $P_{D_{\rm out},\hate}$.

\emph{The second flavor.} The same holds for the second flavor: $\tP_{D,\hate}=y^{\CT}\det M^b_{D,\hate}(y)$ with $M^b(y)=Y(y)^{\T}K^b(y)Y$ (Corollary~\ref{cor.seifertb}), the entry $(k,j)$ of $K^b(y)$ for $k$ outer and $j$ inner vanishes unless both lie on $C$, where it is $s_C(j)\bigl(y_k\ind{A_k=C}-\ind{B_k=C}\bigr)\ind{j\prec_Ck}$ by \eqref{eq.Kbrow}, of $\ve$-degree $h(B_k)-h(C)$, so that the entries of the block $(\mathrm{out},\mathrm{in})$ of $M^b(y)$ are homogeneous of $\ve$-degree $h(B_i)-h(C)$ for the row $i$, and the argument above applies to the two components of the entries of $M^b(y)$, which have the same $\ve$-degrees as those of $M(x)$ (Proposition~\ref{prop.degreesM}).

\emph{The top coefficient.} By Lemma~\ref{lem.pq}, $\tau_D(q)=(-1)^{\mu-1}\sum P_S\tP_T\,q^{\ve(S\triangle T)}$, the sum over the monomials $x^S$ of the component of $P_{D,\hate}$ of $\ve$-degree $\delta-\nu$ and the monomials $y^T$ of the component of $\tP_{D,\hate}$ of $\ve$-degree $\delta$. Writing $S=S_{\rm in}\sqcup S_{\rm out}$ and $T=T_{\rm in}\sqcup T_{\rm out}$, we have $P_S=P^{\rm in}_{S_{\rm in}}P^{\rm out}_{S_{\rm out}}$ and $\tP_T=\tP^{\rm in}_{T_{\rm in}}\tP^{\rm out}_{T_{\rm out}}$, where $P^{\rm in}_{S_{\rm in}}$ etc.\ are the coefficients of the extremal components of the edge polynomials of the subdiagrams, and $\ve(S\triangle T)=\ve(S_{\rm in}\triangle T_{\rm in})+\ve(S_{\rm out}\triangle T_{\rm out})$. The sum therefore factors: with $\mu_{\rm in}$ and $\mu_{\rm out}$ the numbers of components of $L_{\rm in}$ and $L_{\rm out}$,
\[
   \tau_D=(-1)^{\mu-1}(-1)^{\mu_{\rm in}-1}(-1)^{\mu_{\rm out}-1}\,\tau_{D_{\rm in}}\tau_{D_{\rm out}}=\tau_{D_{\rm in}}\tau_{D_{\rm out}},
\]
because $\nu\equiv\mu-1$ modulo $2$ for each of the three connected diagrams and $\nu=\nu_{\rm in}+\nu_{\rm out}$.

\emph{The value at $q=1$.} The identity $\tau_D(1)=\tau_{D_{\rm in}}(1)\,\tau_{D_{\rm out}}(1)$ is the classical multiplicativity of the determinant of the Seifert matrix. By Corollary~\ref{cor.extremesums}, $\tau_D(1)=(\det V_F)^2$ for every connected diagram. The Seifert matrix $V_F$ of $F_D$ is block triangular with respect to the inner and outer cycles. Indeed, by Theorem~\ref{thm.deformedseifert}~(i), $V_F(\gamma,\gamma')=\sum_j\ve_j\,z^\gamma_{e^+_j}\,w_j(\gamma')$, where the edge $e^+_j$ lies on the Seifert circle on the left of the crossing $j$, the one of lower level (Section~\ref{ssec.bandbyband}), and the lift $z^\gamma$ of a cycle vanishes on the circles through which the cycle does not pass. Suppose that $C$ is clockwise, so that $h(C)=H_{\Omega_{\rm out}}+1$: at every outer crossing $j$ on $C$ the circle $C$ is the one of higher level, and $e^+_j$ lies on the other circle of the crossing, which is outside $C$. For an inner cycle $\gamma$ and an outer cycle $\gamma'$, the sum for $V_F(\gamma,\gamma')$ runs over the outer bands $j$ of $\gamma'$, and $z^\gamma$ vanishes on the circles outside $C$; hence $V_F(\gamma,\gamma')=0$. If $C$ is counterclockwise, then $h(C)=H_{\Omega_{\rm in}}+1$, at every inner crossing on $C$ the edge $e^+_j$ lies on a circle inside $C$, and $V_F(\gamma',\gamma)=0$ for the same reason. The diagonal blocks of $V_F$ are the Seifert matrices of $F_{D_{\rm in}}$ and of $F_{D_{\rm out}}$: the lift of an inner cycle in $D$ is constant along $C$ between two consecutive inner crossings, so that it is the lift of the same cycle in $D_{\rm in}$, and likewise for the outer cycles. Hence $\det V_{F_D}=\det V_{F_{D_{\rm in}}}\det V_{F_{D_{\rm out}}}$. The last statement of the theorem follows from the genus bound of \cite{Kohli-Tahar, LNVdV25} and Gabai's theorem, as explained above.
\end{proof}

\begin{remark}\label{rem.murasugi}
(i) Iterating the theorem over all Seifert circles decomposes $F_D$ into the Murasugi sum of the surfaces $F_{D_\Omega}$ of the regions $\Omega$ of $D$, where $D_\Omega$ is the diagram formed by the circles bounding $\Omega$ and the crossings lying in $\Omega$, and gives $\tau_D=\prod_\Omega\tau_{D_\Omega}$. For a closed braid diagram the $D_\Omega$ are the closed $2$-braid diagrams formed by the letters $\sigma_i^{\pm1}$ of the braid word, one for each $i$, which represent torus links $T(2,a_i)$ and unlinks of two components; this is the counterpart of the decomposition of the braided surface into the surfaces of these $2$-braids, and of $\det V_F=\prod_i\det V_{F_{D_{\Omega_i}}}$.

(ii) The Murasugi sums along a Seifert circle are the planar star products of Murasugi and Przytycki, under which the top coefficient of the HOMFLY polynomial is multiplicative \cite{MurasugiPrzytycki89}; Theorem~\ref{thm.murasugi} is the analogue for $\LG$, and the two statements have the same mechanism. The theorem covers the Murasugi sums which are visible on a diagram, in which the two summands and the sum are canonical surfaces, and the restriction bears on the summands, not on the summing disk. A canonical surface is free, i.e.\ its complement is a handlebody, whereas a minimal genus Seifert surface need not be free, nor a free one canonical \cite{Moriah87,Livingston88,KobayashiKobayashi96}. For instance the annulus $A$ whose core is a nontrivial knot $K$ is not free, its complement being homotopy equivalent to that of $K$, and the plumbing of $A$ with a Hopf band is a genus one Seifert surface of a Whitehead double of $K$: a Murasugi sum of minimal genus surfaces which is not a sum of canonical surfaces. The free genus of Whitehead doubles can moreover be arbitrarily large \cite{Moriah87,Livingston88}, and a double of free genus larger than one has no diagram whose canonical surface has minimal genus, so that the theorem says nothing about its top coefficient, which \cite{ALvdV26} expresses as the product of those of $\partial A$ and of the Hopf link. Thus the theorem does not reprove the result of \cite{ALvdV26} for arbitrary Murasugi sums of minimal genus surfaces. For the sums which it covers it has no restriction on $\det V_{F_D}$: when $\det V_{F_D}=0$, where $\tau_D(1)=0$, the identity $\tau_D=\tau_{D_{\rm in}}\tau_{D_{\rm out}}$ still holds with the sign $+$, in agreement with the theorem of \cite{ALvdV26} when $F_D$ has minimal genus, and both sides vanish when it has not. A determinant formula for $\LG$ attached to an arbitrary Seifert surface, of the kind asked for in Section~\ref{ssec.beyond}, would give the general statement by the argument above, which uses only the block structure of the extremal parts of the matrices.

(iii) The proof shows a little more than the statement: with the cut edge on $C$, the extremal homogeneous components of the edge polynomials of $D$ are the products of those of $D_{\rm in}$ and $D_{\rm out}$. The edge polynomials themselves, and $\LG$, are not multiplicative, exactly as $\det(tV-V^{\T})$ is not for a block triangular $V$: the blocks $M_{\rm io}$ and $M_{\rm oi}$ both contribute to the middle degrees.
\end{remark}

\subsection{Alternating links: positivity, fiberedness and the Hopf invariant}\label{ssec.alternating}
For alternating diagrams the extremal monomials of the two edge polynomials are spanning trees, by the theorem of \cite{HK} which proves that $\LG$ of an alternating link is alternating, and the sum \eqref{eq.topformula} has no cancellations. We recall the statement \cite[Definition~8.2 and Proposition~8.5]{HK}. Let $D$ be a connected alternating diagram. Then every edge of $D$ is the over-strand at exactly one of its ends; let $\hate$ be an edge which leaves a crossing as its over-strand and enters the next crossing as its under-strand, $\hate=o'_{j_1}=u_{k_0}$. An \emph{under-tree} of $(D,\hate)$ is a map $\tau$ which assigns to every crossing $j$ one of the two edges $u_j,u'_j$ of its under-strand, such that, following from any crossing the edge assigned to it to the crossing at its other end (and from $k_0$ along $\hate$ to a new vertex $*$ if $\tau(k_0)=\hate$), one eventually arrives at $*$; put $S(\tau)=\{j:\tau(j)=u_j\}$. An \emph{over-tree} $\tau'$ is defined in the same way with the edges $o_j,o'_j$ of the over-strands, and $S'(\tau')=\{j:\tau'(j)=o'_j\}$. By \cite[Proposition~8.5]{HK}, with a sign $\sigma=\pm1$ independent of the trees (it is the sign $\sigma_{D,\hate}$ of Theorem~\ref{thm.kauffmanalt}, see Corollary~\ref{cor.statestrees}),
\begin{equation}\label{eq.HKtrees}
   P_{D,\hate}(x)=\sigma\sum_{\tau}(-1)^{|S(\tau)|}x^{S(\tau)},\qquad \tP_{D,\hate}(y)=\sigma\sum_{\tau'}(-1)^{|S'(\tau')|}y^{S'(\tau')},
\end{equation}
the sums running over the under-trees and the over-trees of $(D,\hate)$; the map $\tau\mapsto S(\tau)$ is injective, so that $P_S=\pm1$ if $S=S(\tau)$ for an under-tree and $P_S=0$ otherwise, and likewise for $\tP_T$.

\begin{theorem}[alternating links]\label{thm.alternating}
Let $D$ be a connected alternating diagram of $L$, with $\nu=c-s+1=2g(L)+\mu-1$, let $\hate=o'_{j_1}=u_{k_0}$ be as above, and let $I_{D,\hate}=[\delta-\nu,\delta]$ be the interval of Theorem~\ref{thm.degrees}. Let $N$ be the absolute value of the leading coefficient of $\Delta_L$. Then:
\begin{enumerate}\renewcommand{\theenumi}{\roman{enumi}}
\item There are exactly $N$ under-trees $\tau$ of smallest $\ve$-degree $\ve(S(\tau))=\delta-\nu$ and exactly $N$ over-trees $\tau'$ of largest $\ve$-degree $\ve(S'(\tau'))=\delta$, and
\[
   \mathrm{top}(L,q)=\sum_{(\tau,\tau')}q^{\,\ve(S(\tau)\triangle S'(\tau'))},
\]
the sum running over the $N^2$ pairs of such trees: the top coefficient of $\LG(L)$ is a polynomial in $q$ with nonnegative integer coefficients and $\mathrm{top}(L,1)=N^2$. The same holds for the coefficient of $p^{-\nu}$, with the under-trees of largest and the over-trees of smallest $\ve$-degree.
\item $L$ is fibered if and only if $N=1$, i.e.\ if and only if $(D,\hate)$ has a unique under-tree $\tau_0$ of smallest and a unique over-tree $\tau'_0$ of largest $\ve$-degree, and then
\[
   \mathrm{top}(L,q)=q^{\,\ve(S(\tau_0)\triangle S'(\tau'_0))} .
\]
\item If $L$ is fibered, with fiber $F$, the Hopf invariant $\lambda(\xi_F)$ of the plane field of \cite{LNVdV25} is
\[
   \lambda(\xi_F)=\tfrac12\bigl(\nu-\ve(S(\tau_0)\triangle S'(\tau'_0))\bigr)=\tfrac12\bigl(c-s+1-\ve(S(\tau_0)\triangle S'(\tau'_0))\bigr).
\]
\end{enumerate}
\end{theorem}

Part (iii) expresses a quantity of the fibration, the Hopf invariant of its plane field, in terms of the diagram alone: the two extremal trees are determined by the signs of the crossings and the Seifert circles (Example~\ref{ex.alttrees}). In Section~\ref{sec.kauffman} we show that the two extremal trees are the two flavors of one Kauffman state and that $\ve(S(\tau_0)\triangle S'(\tau'_0))$ is minus the signature $\sigma_L$ of $L$, so that $\mathrm{top}(L,q)=q^{-\sigma_L}$ and $\lambda(\xi_F)=\frac12(\nu+\sigma_L)$ (Proposition~\ref{prop.signature} and Corollary~\ref{cor.extremalstates}).

\begin{proof}
The theorem follows from Theorem~\ref{thm.tight} once we know two things: that alternating diagrams are tight, and that in \eqref{eq.HKtrees} all the trees of a given $\ve$-degree have the same sign, so that there are no cancellations. Removing a nugatory crossing from an alternating diagram removes one crossing and one Seifert circle and does not change $\nu$. After removing all of them we obtain a reduced alternating diagram of $L$ with the same $\nu$, and for such a diagram $\mathrm{span}\,\Delta_L=\nu$ by \cite{Crowell59, Murasugi58}. So the diagram $D$ is tight and Theorem~\ref{thm.tight} applies.

By \eqref{eq.HKtrees}, $P_{D,\hate}(t^{\ve})=\sigma\sum_\tau(-1)^{|S(\tau)|}t^{\ve(S(\tau))}$, and $|S|\equiv\ve(S)\pmod2$ for every $S$, since $\ve(S)=|S|-2\,|\{j\in S:\ve_j=-1\}|$; so all the under-trees of a given $\ve$-degree $\ell$ contribute with the same sign $\sigma(-1)^\ell$ to the coefficient of $t^\ell$, which is therefore $\pm$ the number of under-trees of $\ve$-degree $\ell$. At $\ell=\delta-\nu$ this coefficient is, by Theorem~\ref{thm.tight} (i), an extreme coefficient of $\Delta_L$, of absolute value $N$ (the extreme coefficients of $\Delta_L$ are equal up to sign, by the symmetry of the Alexander polynomial). The same holds for the over-trees at $\ell=\delta$. In the sum of Theorem~\ref{thm.tight} (iii), $P_S\tP_T=\sigma^2(-1)^{|S|+|T|}=(-1)^{\ve(S)+\ve(T)}=(-1)^{\delta-\nu+\delta}$ is the same for all $N^2$ pairs, so $\mathrm{top}(L,q)=\pm\sum_{(\tau,\tau')}q^{\ve(S(\tau)\triangle S'(\tau'))}$, and the sign is $+$ because $\mathrm{top}(L,1)=N^2>0$ by Theorem~\ref{thm.tight} (iii). This proves (i). By Murasugi's theorem \cite{Murasugi63}, an alternating link is fibered if and only if the leading coefficient of its Alexander polynomial is $\pm1$, which gives (ii). For (iii), the fiber $F$ of a fibered link has minimal genus, so $1-\chi(F)=2g(L)+\mu-1=\nu$ and \cite[Corollary~1]{LNVdV25} gives $\mathrm{top}(L,q)=q^{\nu-2\lambda(\xi_F)}$; comparing the exponents with (ii) gives the formula.
\end{proof}

\begin{figure}[ht]
\centering
\begin{tikzpicture}[>=Stealth,scale=0.5]
\begin{scope}
\draw (0.00,-0.00) -- (0.00,0.50);
\draw[->] (0.00,1.50) -- (0.00,2.50);
\draw[->,red,line width=2pt] (0.00,3.50) -- (0.00,4.50);
\node[left,font=\scriptsize] at (0.00,4.00) {$u'_2$};
\draw (0.00,5.50) -- (0.00,6.00);
\draw[rounded corners=6pt,red,line width=2pt] (0.00,6.00) -- (0.00,7.15) -- (4.80,7.15) -- (4.80,-0.95) -- (0.00,-0.95);
\draw[red,line width=2pt] (0.00,-0.95) -- (0.00,0.00);
\node[above,font=\scriptsize] at (3.90,7.15) {$u'_3$};
\draw (1.70,-0.00) -- (1.70,0.50);
\draw[->] (1.70,1.50) -- (1.70,2.50);
\draw[->] (1.70,3.50) -- (1.70,4.50);
\draw (1.70,5.50) -- (1.70,6.00);
\draw[rounded corners=6pt,red,line width=2pt] (1.70,6.00) -- (1.70,6.40) -- (4.10,6.40) -- (4.10,-0.40) -- (1.70,-0.40);
\draw[red,line width=2pt] (1.70,-0.40) -- (1.70,0.00);
\node[above,font=\scriptsize] at (3.20,6.40) {$u_1=\hate$};
\draw (1.70,0.50) -- (0.00,1.50);
\draw[preaction={draw=white,line width=4pt}] (0.00,0.50) -- (1.70,1.50);
\node[font=\scriptsize,circle,inner sep=1pt,fill=white] at (0.30,1.00) {$1$};
\draw (1.70,2.50) -- (0.00,3.50);
\draw[preaction={draw=white,line width=4pt}] (0.00,2.50) -- (1.70,3.50);
\node[font=\scriptsize,circle,inner sep=1pt,fill=white] at (0.30,3.00) {$2$};
\draw (1.70,4.50) -- (0.00,5.50);
\draw[preaction={draw=white,line width=4pt}] (0.00,4.50) -- (1.70,5.50);
\node[font=\scriptsize,circle,inner sep=1pt,fill=white] at (0.30,5.00) {$3$};
\node[font=\small] at (2.4,-1.8) {$\tau_0$: $S(\tau_0)=\{1\}$};
\end{scope}
\begin{scope}[xshift=7cm]
\draw (0.00,-0.00) -- (0.00,0.50);
\draw[->] (0.00,1.50) -- (0.00,2.50);
\draw[->] (0.00,3.50) -- (0.00,4.50);
\draw (0.00,5.50) -- (0.00,6.00);
\draw[rounded corners=6pt] (0.00,6.00) -- (0.00,7.15) -- (4.80,7.15) -- (4.80,-0.95) -- (0.00,-0.95);
\draw (0.00,-0.95) -- (0.00,0.00);
\draw (1.70,-0.00) -- (1.70,0.50);
\draw[->,blue,line width=2pt] (1.70,1.50) -- (1.70,2.50);
\node[right,font=\scriptsize] at (1.70,2.00) {$o'_1$};
\draw[->,blue,line width=2pt] (1.70,3.50) -- (1.70,4.50);
\node[right,font=\scriptsize] at (1.70,4.00) {$o'_2$};
\draw (1.70,5.50) -- (1.70,6.00);
\draw[rounded corners=6pt,blue,line width=2pt] (1.70,6.00) -- (1.70,6.40) -- (4.10,6.40) -- (4.10,-0.40) -- (1.70,-0.40);
\draw[blue,line width=2pt] (1.70,-0.40) -- (1.70,0.00);
\node[above,font=\scriptsize] at (3.20,6.40) {$o'_3=\hate$};
\draw (1.70,0.50) -- (0.00,1.50);
\draw[preaction={draw=white,line width=4pt}] (0.00,0.50) -- (1.70,1.50);
\node[font=\scriptsize,circle,inner sep=1pt,fill=white] at (0.30,1.00) {$1$};
\draw (1.70,2.50) -- (0.00,3.50);
\draw[preaction={draw=white,line width=4pt}] (0.00,2.50) -- (1.70,3.50);
\node[font=\scriptsize,circle,inner sep=1pt,fill=white] at (0.30,3.00) {$2$};
\draw (1.70,4.50) -- (0.00,5.50);
\draw[preaction={draw=white,line width=4pt}] (0.00,4.50) -- (1.70,5.50);
\node[font=\scriptsize,circle,inner sep=1pt,fill=white] at (0.30,5.00) {$3$};
\node[font=\small] at (2.4,-1.8) {$\tau'_0$: $S'(\tau'_0)=\{1,2,3\}$};
\end{scope}
\begin{scope}[xshift=14.5cm,yshift=-1.2cm,scale=0.8]
\draw (0.00,-0.00) -- (0.00,0.50);
\draw[->] (0.00,1.50) -- (0.00,4.50);
\draw[->] (0.00,5.50) -- (0.00,8.00);
\draw[rounded corners=6pt,red,line width=2pt] (0.00,8.00) -- (0.00,9.90) -- (8.20,9.90) -- (8.20,-1.50) -- (0.00,-1.50);
\draw[red,line width=2pt] (0.00,-1.50) -- (0.00,0.00);
\node[above,font=\scriptsize] at (7.30,9.90) {$u'_3$};
\draw (2.20,-0.00) -- (2.20,0.50);
\draw[->,red,line width=2pt] (2.20,1.50) -- (2.20,2.50);
\node[right,font=\scriptsize] at (2.20,2.00) {$u_2$};
\draw[->] (2.20,3.50) -- (2.20,4.50);
\draw[->,red,line width=2pt] (2.20,5.50) -- (2.20,6.50);
\node[right,font=\scriptsize] at (2.20,6.00) {$u_4$};
\draw (2.20,7.50) -- (2.20,8.00);
\draw[rounded corners=6pt,red,line width=2pt] (2.20,8.00) -- (2.20,9.15) -- (7.50,9.15) -- (7.50,-0.95) -- (2.20,-0.95);
\draw[red,line width=2pt] (2.20,-0.95) -- (2.20,0.00);
\node[above,font=\scriptsize] at (6.60,9.15) {$u_1=\hate$};
\draw[->] (4.40,-0.00) -- (4.40,2.50);
\draw[->] (4.40,3.50) -- (4.40,6.50);
\draw (4.40,7.50) -- (4.40,8.00);
\draw[rounded corners=6pt] (4.40,8.00) -- (4.40,8.40) -- (6.80,8.40) -- (6.80,-0.40) -- (4.40,-0.40);
\draw (4.40,-0.40) -- (4.40,0.00);
\draw (2.20,0.50) -- (0.00,1.50);
\draw[preaction={draw=white,line width=4pt}] (0.00,0.50) -- (2.20,1.50);
\node[font=\scriptsize,circle,inner sep=1pt,fill=white] at (0.55,1.00) {$1$};
\draw (2.20,2.50) -- (4.40,3.50);
\draw[preaction={draw=white,line width=4pt}] (4.40,2.50) -- (2.20,3.50);
\node[font=\scriptsize,circle,inner sep=1pt,fill=white] at (2.75,3.00) {$2$};
\draw (2.20,4.50) -- (0.00,5.50);
\draw[preaction={draw=white,line width=4pt}] (0.00,4.50) -- (2.20,5.50);
\node[font=\scriptsize,circle,inner sep=1pt,fill=white] at (0.55,5.00) {$3$};
\draw (2.20,6.50) -- (4.40,7.50);
\draw[preaction={draw=white,line width=4pt}] (4.40,6.50) -- (2.20,7.50);
\node[font=\scriptsize,circle,inner sep=1pt,fill=white] at (2.75,7.00) {$4$};
\node[font=\small] at (4.1,-2.4) {$\tau_0$: $S(\tau_0)=\{1,2,4\}$};
\end{scope}
\begin{scope}[xshift=22.5cm,yshift=-1.2cm,scale=0.8]
\draw (0.00,-0.00) -- (0.00,0.50);
\draw[->] (0.00,1.50) -- (0.00,4.50);
\draw[->] (0.00,5.50) -- (0.00,8.00);
\draw[rounded corners=6pt] (0.00,8.00) -- (0.00,9.90) -- (8.20,9.90) -- (8.20,-1.50) -- (0.00,-1.50);
\draw (0.00,-1.50) -- (0.00,0.00);
\draw (2.20,-0.00) -- (2.20,0.50);
\draw[->,blue,line width=2pt] (2.20,1.50) -- (2.20,2.50);
\node[right,font=\scriptsize] at (2.20,2.00) {$o'_1$};
\draw[->] (2.20,3.50) -- (2.20,4.50);
\draw[->,blue,line width=2pt] (2.20,5.50) -- (2.20,6.50);
\node[right,font=\scriptsize] at (2.20,6.00) {$o'_3$};
\draw (2.20,7.50) -- (2.20,8.00);
\draw[rounded corners=6pt,blue,line width=2pt] (2.20,8.00) -- (2.20,9.15) -- (7.50,9.15) -- (7.50,-0.95) -- (2.20,-0.95);
\draw[blue,line width=2pt] (2.20,-0.95) -- (2.20,0.00);
\node[above,font=\scriptsize] at (6.60,9.15) {$o'_4=\hate$};
\draw[->,blue,line width=2pt] (4.40,-0.00) -- (4.40,2.50);
\node[right,font=\scriptsize] at (4.40,1.2) {$o_2$};
\draw[->] (4.40,3.50) -- (4.40,6.50);
\draw (4.40,7.50) -- (4.40,8.00);
\draw[rounded corners=6pt,blue,line width=2pt] (4.40,8.00) -- (4.40,8.40) -- (6.80,8.40) -- (6.80,-0.40) -- (4.40,-0.40);
\draw[blue,line width=2pt] (4.40,-0.40) -- (4.40,0.00);
\draw (2.20,0.50) -- (0.00,1.50);
\draw[preaction={draw=white,line width=4pt}] (0.00,0.50) -- (2.20,1.50);
\node[font=\scriptsize,circle,inner sep=1pt,fill=white] at (0.55,1.00) {$1$};
\draw (2.20,2.50) -- (4.40,3.50);
\draw[preaction={draw=white,line width=4pt}] (4.40,2.50) -- (2.20,3.50);
\node[font=\scriptsize,circle,inner sep=1pt,fill=white] at (2.75,3.00) {$2$};
\draw (2.20,4.50) -- (0.00,5.50);
\draw[preaction={draw=white,line width=4pt}] (0.00,4.50) -- (2.20,5.50);
\node[font=\scriptsize,circle,inner sep=1pt,fill=white] at (0.55,5.00) {$3$};
\draw (2.20,6.50) -- (4.40,7.50);
\draw[preaction={draw=white,line width=4pt}] (4.40,6.50) -- (2.20,7.50);
\node[font=\scriptsize,circle,inner sep=1pt,fill=white] at (2.75,7.00) {$4$};
\node[font=\small] at (4.1,-2.4) {$\tau'_0$: $S'(\tau'_0)=\{1,3,4\}$};
\end{scope}
\end{tikzpicture}
\caption{The extremal trees of the trefoil (left) and of the figure-eight knot (right), for the cut edge $\hate=u_1$. The under-tree $\tau_0$ of smallest $\ve$-degree assigns to every crossing the red edge of its under-strand, the over-tree $\tau'_0$ of largest $\ve$-degree the blue edge of its over-strand; from every crossing, following the assigned edges leads to $\hate$. The crossings of $S(\tau_0)\triangle S'(\tau'_0)$ are those at which $\tau_0$ chooses the incoming under-edge $u_j$ or $\tau'_0$ chooses the outgoing over-edge $o'_j$, but not both: $\{2,3\}$ in both cases, of $\ve$-degree $2$ for the trefoil and $0$ for the figure-eight knot, so that $\mathrm{top}(3_1,q)=q^2$ and $\mathrm{top}(4_1,q)=1$.}
\label{fig.trees}
\end{figure}
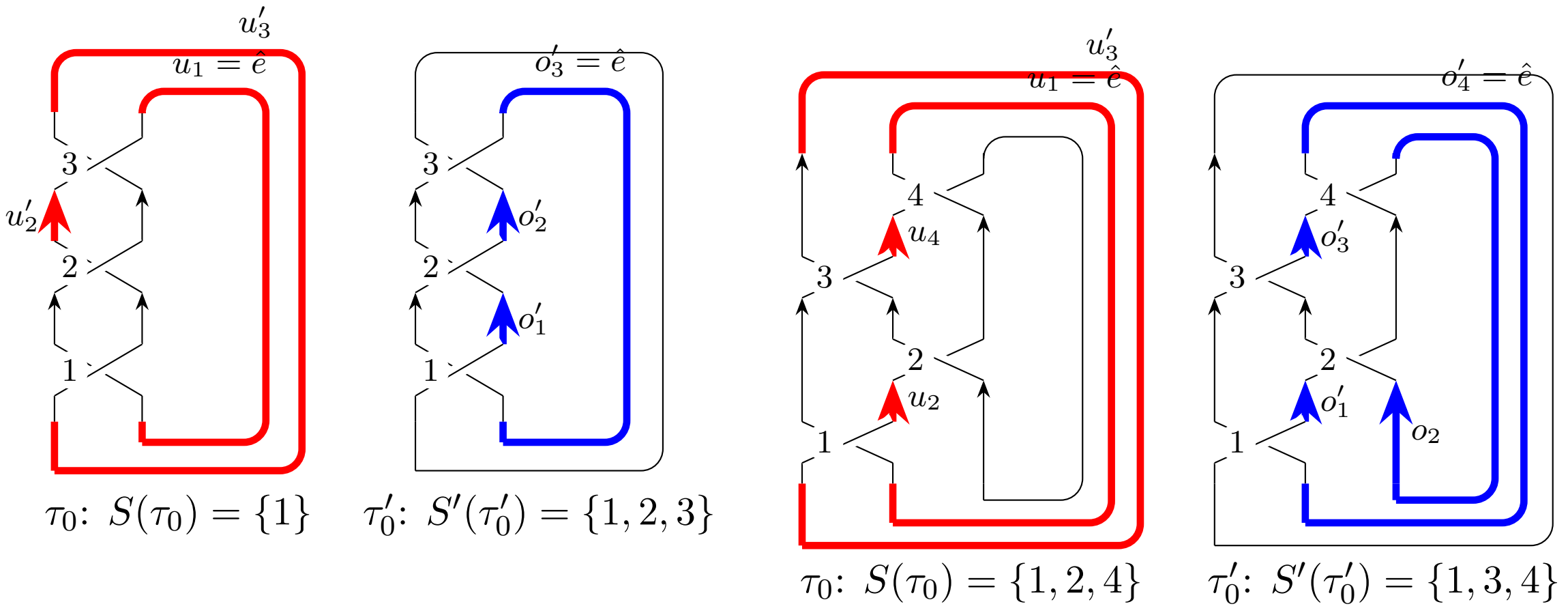

\begin{example}\label{ex.alttrees}
For the standard diagrams of the four alternating fibered knots $3_1$, $4_1$, $6_2$ and $6_3$, drawn as the closures of the braids $\sigma_1^3$, $\sigma_1\sigma_2^{-1}\sigma_1\sigma_2^{-1}$, $\sigma_1^{-1}\sigma_2\sigma_1^{-1}\sigma_2^3$ and $\sigma_1^{-1}\sigma_2^2\sigma_1^{-2}\sigma_2$, with the crossings numbered in the order of the braid words and the cut edge $\hate=u_1$ (the edge entering the first crossing as its under-strand, which leaves the last crossing as its over-strand), the extremal trees, computed from \eqref{eq.HKtrees}, are the following. Figure~\ref{fig.trees} shows the first two.
\[
\begin{array}{c|c|c|c|c|c|c}
   K & \ve & \nu & S(\tau_0) & S'(\tau'_0) & \ve(S(\tau_0)\triangle S'(\tau'_0)) & \lambda(\xi_F)\\ \hline
   3_1 & +++ & 2 & \{1\} & \{1,2,3\} & 2 & 0\\
   4_1 & +-+- & 2 & \{1,2,4\} & \{1,3,4\} & 0 & 1\\
   6_2 & -+-+++ & 4 & \{1,3\} & \{2,3,5,6\} & 2 & 1\\
   6_3 & -++--+ & 4 & \{1,4,5\} & \{2,3,5\} & 0 & 2
\end{array}
\]
The top coefficients $q^2,1,q^2,1$ agree with the table of \cite{LNVdV25}. For the trefoil, $P_{D,u_1}=x_1-x_1x_2+x_1x_2x_3$ and $\tP_{D,u_1}=y_3-y_2y_3+y_1y_2y_3$ (Example~\ref{ex.trefoil} and Lemma~\ref{lem.understrand}), the three under-trees have $\ve$-degrees $1,2,3$, the unique one of $\ve$-degree $1$ is $S(\tau_0)=\{1\}$, the unique over-tree of $\ve$-degree $3$ is $S'(\tau'_0)=\{1,2,3\}$, and $\mathrm{top}(3_1,q)=q^{\ve(\{2,3\})}=q^2$. For the figure-eight knot, $P_{D,u_1}=x_1x_2+x_1x_4+x_1x_2x_3x_4-x_1x_2x_4-x_1x_2x_3$ has the $\ve$-degrees $0,0,0,-1,1$, the extremal under-tree is $S(\tau_0)=\{1,2,4\}$ of $\ve$-degree $-1$, the extremal over-tree is $S'(\tau'_0)=\{1,3,4\}$ of $\ve$-degree $1$, and $S(\tau_0)\triangle S'(\tau'_0)=\{2,3\}$ has $\ve$-degree $\ve_2+\ve_3=0$: $\mathrm{top}(4_1,q)=1$ and $\lambda(\xi_F)=1$.
\end{example}

\begin{remark}\label{rem.homogeneous}
Theorem~\ref{thm.tight} applies to every tight diagram, and the only input of Theorem~\ref{thm.alternating} beyond it is the sign-coherence \eqref{eq.HKtrees} of the monomials of the two edge polynomials. On all the tight nonalternating diagrams we have computed (closures of positive and of homogeneous $3$- and $4$-braids with up to ten crossings, among them the torus knots $T(2,5)$, $T(2,7)$ and $T(3,4)$ and the closures of the positive $3$-braids $\sigma_1^2\sigma_2^2\sigma_1\sigma_2$ and $\sigma_1^3\sigma_2^3$, and a few hundred random one-disk band surfaces with bands of framing $\pm1$), the monomials of $P_{D,\hate}$ of smallest $\ve$-degree all have the same sign, those of $\tP_{D,\hate}$ of largest $\ve$-degree as well, and the number of each is the leading coefficient of $\Delta_L$ in absolute value; all these knots are fibered, the extremal monomials are unique and $\mathrm{top}(L,q)$ is the monomial $q^{\ve(S_0\triangle T_0)}$ of \cite{LNVdV25}. We do not know whether the sign-coherence of the extremal monomials holds on every tight diagram; Cromwell's homogeneous diagrams, which are Murasugi sums of alternating diagrams, are the natural test case. For every diagram the sum of the extremal monomials of $P_{D,\hate}$ is $\pm\det V_F$ (Corollary~\ref{cor.extremesums}), so that their sign-coherence is equivalent to their number being $|\det V_F|$, the absolute value of the leading coefficient of $\Delta_L$.
\end{remark}

Theorems~\ref{thm.murasugi} and~\ref{thm.alternating} together make the top coefficient explicit for the closures of homogeneous braids, which are Murasugi sums of alternating $2$-braid diagrams along their Seifert circles.

\begin{corollary}[homogeneous braids]\label{cor.homogeneous}
Let $D$ be the closure of a homogeneous braid on $n$ strands in which the generator $\sigma_i$ occurs $a_i\ge1$ times, always with the exponent $\ve_i\in\{\pm1\}$ $(i=1,\dots,n-1)$, and let $L$ be the link it represents, which is fibered with fiber $F_D$ \cite{Stallings78}, so that $\nu=\sum_{i=1}^{n-1}(a_i-1)=2g(L)+\mu-1$. Then
\[
   \mathrm{top}(L,q)=q^{\sum_i\ve_i(a_i-1)},\qquad\text{and}\qquad \lambda(\xi_{F_D})=\sum_{i:\,\ve_i=-1}(a_i-1)
\]
for the Hopf invariant of the plane field of the fibration \cite{LNVdV25}.
\end{corollary}
\begin{proof}
Iterating Theorem~\ref{thm.murasugi} over the Seifert circles $C_2,\dots,C_{n-1}$ of $D$ (Remark~\ref{rem.murasugi} (i)) gives $\tau_D=\prod_i\tau_{D_i}$, where $D_i$ is the closed $2$-braid diagram formed by the $a_i$ letters $\sigma_i^{\ve_i}$; and $\tau_D=\mathrm{top}(L,q)$ since $F_D$ has minimal genus. If $a_i\ge2$, $D_i$ is a reduced alternating diagram of the torus link $T(2,\ve_ia_i)$ with $\nu_i=a_i-1$, and by Theorem~\ref{thm.alternating} (ii) $\tau_{D_i}=q^{\ve(S(\tau_0)\triangle S'(\tau'_0))}$ with $S(\tau_0)\triangle S'(\tau'_0)=\{2,\dots,a_i\}$ for the cut edge $u_1$, exactly as for the trefoil in Example~\ref{ex.alttrees}; so $\tau_{D_i}=q^{\ve_i(a_i-1)}$. If $a_i=1$, $D_i$ is a diagram of the unknot with $\nu_i=0$ and $\tau_{D_i}=1$. The formula for $\lambda$ follows from $\mathrm{top}(L,q)=q^{\nu-2\lambda(\xi_F)}$ \cite[Corollary~1]{LNVdV25}.
\end{proof}

For the four knots of Example~\ref{ex.alttrees} this gives $q^2,1,q^2,1$ and $\lambda=0,1,1,2$ again; and for the positive braids $\mathrm{top}(L,q)=q^{\nu}$ and $\lambda=0$. In the language of \cite{LNVdV25, ALvdV26}, the corollary says that the top coefficient of a homogeneous braid closure is the product of the top coefficients $q^{\pm1}$ of the $\sum_i(a_i-1)$ Hopf bands of which its fiber is the plumbing (Stallings), the exponent being the number of positive Hopf bands minus the number of negative ones.

\subsection{The monomials of \texorpdfstring{$\tP_{D,\hate}$}{P~} as disk decompositions of the canonical surface}\label{ssec.sectors}
The proof of Theorem~\ref{thm.degrees} locates the $\ve$-degrees but does not say which monomials occur. For the second flavor the matrix of Proposition~\ref{prop.seifertb} gives a complete and geometric answer, which we include because it explains the role of the levels: the monomials of $\tP_{D,\hate}$ are the ways of cutting the canonical Seifert surface into a disk along arcs of its Seifert disks, and the $\ve$-degree of a monomial is, up to an additive constant, minus the sum of the levels of the pieces.

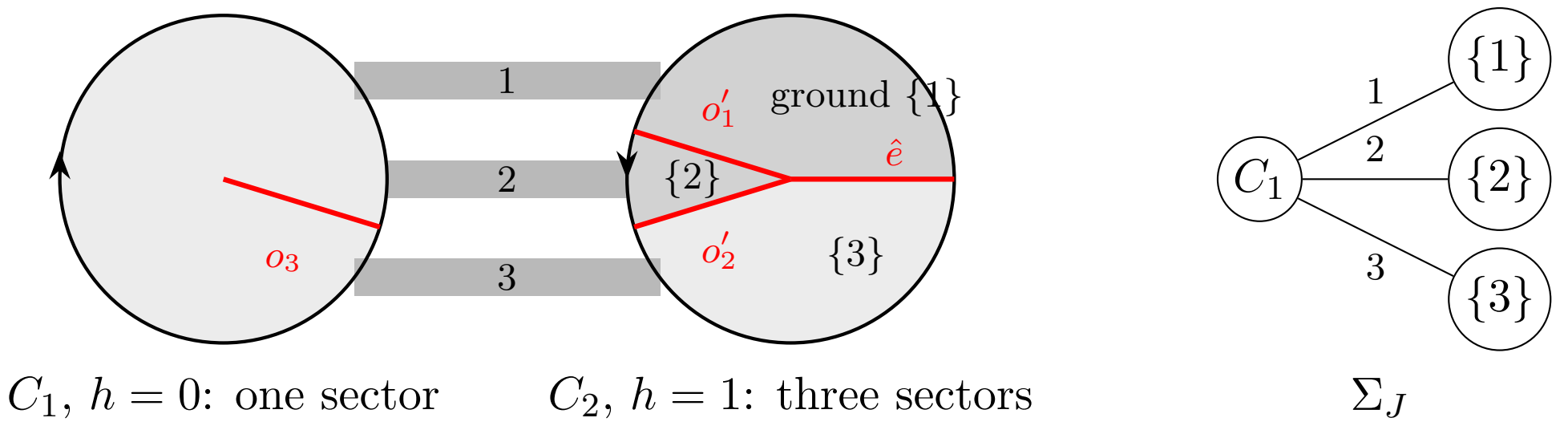
\begin{figure}[ht]
\centering
\begin{tikzpicture}[>=Stealth,scale=0.9,every node/.style={font=\small}]
\fill[gray!15] (0,0) circle (1.5);
\fill[gray!15] (5.2,0) circle (1.5);
\fill[gray!35] (5.2,0) -- ({5.2+1.5*cos(163)},{1.5*sin(163)}) arc (163:197:1.5) -- cycle;
\fill[gray!35] (5.2,0) -- (6.7,0) arc (0:163:1.5) -- cycle;
\foreach \y/\lab in {0.9/1, 0/2, -0.9/3} {
  \fill[gray!55] ({sqrt(2.25-\y*\y)},\y+0.17) -- ({5.2-sqrt(2.25-\y*\y)},\y+0.17) -- ({5.2-sqrt(2.25-\y*\y)},\y-0.17) -- ({sqrt(2.25-\y*\y)},\y-0.17) -- cycle;
  \node[font=\scriptsize] at (2.6,\y) {$\lab$};
}
\draw[thick,postaction={decorate},decoration={markings,mark=at position 0.5 with {\arrow{<}}}] (0,0) circle (1.5);
\draw[thick,postaction={decorate},decoration={markings,mark=at position 0.5 with {\arrow{>}}}] (5.2,0) circle (1.5);
\draw[red,very thick] (0,0) -- ({1.5*cos(-17)},{1.5*sin(-17)});
\draw[red,very thick] (5.2,0) -- ({5.2+1.5*cos(163)},{1.5*sin(163)});
\draw[red,very thick] (5.2,0) -- ({5.2+1.5*cos(197)},{1.5*sin(197)});
\draw[red,very thick] (5.2,0) -- (6.7,0);
\node[font=\scriptsize,red] at (0.55,-0.75) {$o_3$};
\node[font=\scriptsize,red] at (4.55,0.65) {$o'_1$};
\node[font=\scriptsize,red] at (4.55,-0.65) {$o'_2$};
\node[font=\scriptsize,red] at (6.15,0.25) {$\hate$};
\node at (0,-2.0) {$C_1$, $h=0$: one sector};
\node at (5.2,-2.0) {$C_2$, $h=1$: three sectors};
\node[font=\scriptsize] at (5.8,-0.7) {$\{3\}$};
\node[font=\scriptsize] at (5.9,0.75) {ground $\{1\}$};
\node[font=\scriptsize] at (4.3,0.0) {$\{2\}$};
\begin{scope}[xshift=9.5cm]
\node[draw,circle,inner sep=1.5pt] (a) at (0,0) {$C_1$};
\node[draw,circle,inner sep=1.5pt] (g) at (2.2,1.1) {$\{1\}$};
\node[draw,circle,inner sep=1.5pt] (b) at (2.2,0) {$\{2\}$};
\node[draw,circle,inner sep=1.5pt] (c) at (2.2,-1.1) {$\{3\}$};
\draw (a) -- node[above,font=\scriptsize] {$1$} (g); \draw (a) -- node[above,font=\scriptsize] {$2$} (b); \draw (a) -- node[below,font=\scriptsize] {$3$} (c);
\node at (1.1,-2.0) {$\Sigma_J$};
\end{scope}
\end{tikzpicture}
\caption{The monomial $y_3$ of $\tP_{D,u_1}=y_3-y_2y_3+y_1y_2y_3$ for the standard diagram of the trefoil (Example~\ref{ex.alttrees}), as a disk decomposition of its canonical surface: the two Seifert disks, of levels $0$ and $1$, joined by the three bands. The chosen edges of $J=\{3\}$ are $o_3=u'_2$ on $C_1$ and $o'_1=u_2$, $o'_2=u_3$ on $C_2$, together with $\hate=u_1=o'_3$; the radial arcs (red) cut $C_1$ into one sector containing all three crossings and $C_2$ into three sectors containing one crossing each. The sector graph is a tree, and the $\ve$-degree of $y_3$ is $\sum_jh(B_j)+h(\hat C)-\sum_\varpi h(\varpi)=3+1-3=1$.}
\label{fig.sectors}
\end{figure}

Let $J\subseteq\{1,\dots,c\}$ and choose at every crossing $j$ the over-edge $o_j$ if $j\in J$ and the over-edge $o'_j$ if $j\notin J$: the \emph{chosen edges} of $J$. On the circle $\hat C$ we also regard $\hate$ as chosen. Say that $J$ is \emph{admissible} if the chosen edges are pairwise distinct (an edge $o'_j=o_k$ can be chosen twice), $\hate$ is chosen only in its capacity as the cut edge (i.e.\ $j\in J$ if $o'_j=\hate$ and $j\notin J$ if $o_j=\hate$), and every Seifert circle contains a chosen edge. For an admissible $J$, the chosen edges cut every Seifert circle $C$ into arcs; we call \emph{sector} the arc which starts with a chosen edge $e$ and contains the edges of $C$ from $e$ up to, but excluding, the next chosen edge, together with the crossings at the ends of these edges: the sector \emph{contains} these crossings. Every crossing of $C$ lies in exactly one sector of $C$, and a sector is a sector of the disk of $F_D$ bounded by $C$, between two radial arcs drawn from the center of the disk to the chosen edges (Figure~\ref{fig.sectors}). The sector of $\hat C$ which starts with $\hate$ is the \emph{ground sector}. The \emph{sector graph} $\Sigma_J$ has the sectors as vertices and one edge for every crossing $k$, joining the sector of $A_k$ containing $k$ to the sector of $B_k$ containing $k$; it is obtained from the Seifert graph $\Gamma_D$ by splitting the vertex $C$ into as many vertices as there are chosen edges on $C$ (plus one for $\hat C$), and it has $c+1$ vertices and $c$ edges. Equivalently, $\Sigma_J$ is the graph of the surface obtained from $F_D$ by cutting each Seifert disk along its radial arcs: a sector of $C$ is a piece of the disk of $C$, and a crossing is the band joining two pieces.

\begin{proposition}\label{prop.sectors}
The monomial $y^J$ occurs in $\tP_{D,\hate}(y)$ if and only if $J$ is admissible and the sector graph $\Sigma_J$ is a tree, i.e.\ if and only if the radial arcs of $J$ cut $F_D$ into a disk; its coefficient is then $\pm1$. The $\ve$-degree of such a monomial is
\[
   \ve(J)=\sum_{j=1}^ch(B_j)+h(\hat C)-\sum_{\text{sectors }\varpi}h(\varpi),
\]
where $h(\varpi)$ is the level of the Seifert circle of the sector $\varpi$ and the sum runs over all sectors, the ground sector included.
\end{proposition}
\begin{proof}
We expand the determinant of the matrix of Proposition~\ref{prop.seifertb}. Each of its first $c$ rows is the sum of a row without variables and of $y_j$ times another row without variables, so the coefficient of the monomial $y^J$ is the determinant of an integer matrix; by row operations we identify this matrix with the incidence matrix of the sector graph $\Sigma_J$, with one row removed, and the matrix-tree theorem concludes.

\emph{The expansion.} For an edge $e$ on the circle $C$ put $[e]^b=e_{b_C}+J^b_C(e)\in\BZ^c\oplus\BZ^{s-1}$, the jump sum $J^b_C(e)$ of \eqref{eq.Kbrow} being written in the basis $\d_1,\dots,\d_c$, with the conventions $e_{b_{\hat C}}=0$ and $[\hate]^b=0$. By the jump relations of the second flavor (proof of Proposition~\ref{prop.seifertb}), $[o'_j]^b=[u_j]^b-\d_j$, so the row $j$ of the matrix $W^b=\bigl(\begin{smallmatrix}K^b(y)&\CC(y)^{\T}\\\hat\partial&0\end{smallmatrix}\bigr)$ of Proposition~\ref{prop.seifertb} is $y_j[o_j]^b-[o'_j]^b$: it is $u_j+y_jv_j$ with the integer vectors $u_j=-[o'_j]^b$ and $v_j=[o_j]^b$, which involve no variable. Hence $\tP_{D,\hate}(y)=\det W^b=\sum_Jy^J\det W^b_J$, where $W^b_J$ has the rows $v_j=[o_j]^b$ for $j\in J$, $u_j=-[o'_j]^b$ for $j\notin J$, and the rows $\hat\partial_C$, $C\ne\hat C$.

\emph{The sets which are not admissible.} If a chosen edge is chosen twice, or if $\hate$ is chosen by a crossing, two rows of $W^b_J$ are equal up to sign or one row vanishes, and $\det W^b_J=0$. Otherwise the rows of $W^b_J$ are the vectors $\pm[e]^b$, $e$ running over the chosen edges other than $\hate$, together with the $\hat\partial_C$. The column $b_C$ of $W^b_J$ ($C\ne\hat C$) has the entry $\pm1$ in the rows of the chosen edges on $C$ and $0$ elsewhere; so $\det W^b_J=0$ if $C$ has no chosen edge. Hence $\det W^b_J=0$ unless $J$ is admissible.

\emph{The admissible sets.} Let $J$ be admissible. If $e_1,\dots,e_r$ are the chosen edges of $C$ in the order in which they are met along $C$ from $b_C$, replace the rows of $e_2,\dots,e_r$ by $[e_a]^b-[e_{a-1}]^b$, which is the sum of the jumps $\pm\d_k$ at the crossings $k$ of the sector starting at $e_{a-1}$, and the row $\hat\partial_C=-\sum_{k\text{ on }C}(\text{jump at }k)$ by $\hat\partial_C+[e_r]^b-[e_1]^b$, which is minus the sum of the jumps at the crossings of the sector starting at $e_r$; then expand along the column $b_C$, whose only remaining entry is in the row of $e_1$. For $\hat C$, whose base edge is $\hate$ with $[\hate]^b=0$, the rows of the chosen edges are the jump sums $J^b_{\hat C}(e_a)$, and the differences are the sums of the jumps over the sectors starting at $\hate,e_1,\dots,e_{r-1}$: every sector of $\hat C$ except the ground sector gets a row. We obtain, up to sign, a $c\times c$ matrix with rows indexed by the sectors other than the ground sector and columns by the crossings, whose entry is $\pm1$ if the crossing lies in the sector (the sign being $+$ if the sector is on $A_k$ and $-$ if it is on $B_k$) and $0$ otherwise: it is the incidence matrix of the graph $\Sigma_J$ with the row of the ground sector removed. By the matrix-tree theorem (the reduced incidence matrix of a graph with $c+1$ vertices and $c$ edges is invertible if and only if the graph is a tree, and its determinant is then $\pm1$), the first claim follows.

\emph{The $\ve$-degree.} The condition $j\in J$ means that the chosen edge of $j$ lies on $A_j$ and $j\notin J$ that it lies on $B_j$, so $\ve(J)=\sum_{j\in J}(h(B_j)-h(A_j))=\sum_jh(B_j)-\sum_jh(C_j)$, where $C_j$ is the circle of the chosen edge of $j$; and $\sum_jh(C_j)$ is the sum of the levels of the sectors other than the ground sector, each sector being counted once for the chosen edge which starts it.
\end{proof}

The monomials of $P_{D,\hate}(x)$ have the same description in terms of the under-edges $u_j,u'_j$: if $\bD$ is the diagram obtained from $D$ by reversing all orientations and changing all crossings, the edges $(o_j,u_j,o'_j,u'_j)$ of $D$ are the edges $(u'_j,o'_j,u_j,o_j)$ of $\bD$, so that $T^a_{\bD}(x)=T^b_D(x)^{\T}$ by \eqref{eq.Ta} and \eqref{eq.Tb} and $P_{D,\hate}(x)=\tP_{\bD,\hate}(x)$; the diagram $\bD$ has the same Seifert circles as $D$, with the opposite orientations, the same $A_j$ and $B_j$, the signs $-\ve_j$ and the levels $-h-1$. By Theorem~\ref{thm.degrees}, Proposition~\ref{prop.sectors} contains the following combinatorial statement, for which we have no direct proof: over all the ways of cutting $F_D$ into a disk along radial arcs of its Seifert disks, one arc for each over-edge of each crossing, the sum of the levels of the sectors varies in an interval of length at most $\nu$.

\subsection{Beyond the canonical surface}\label{ssec.beyond}
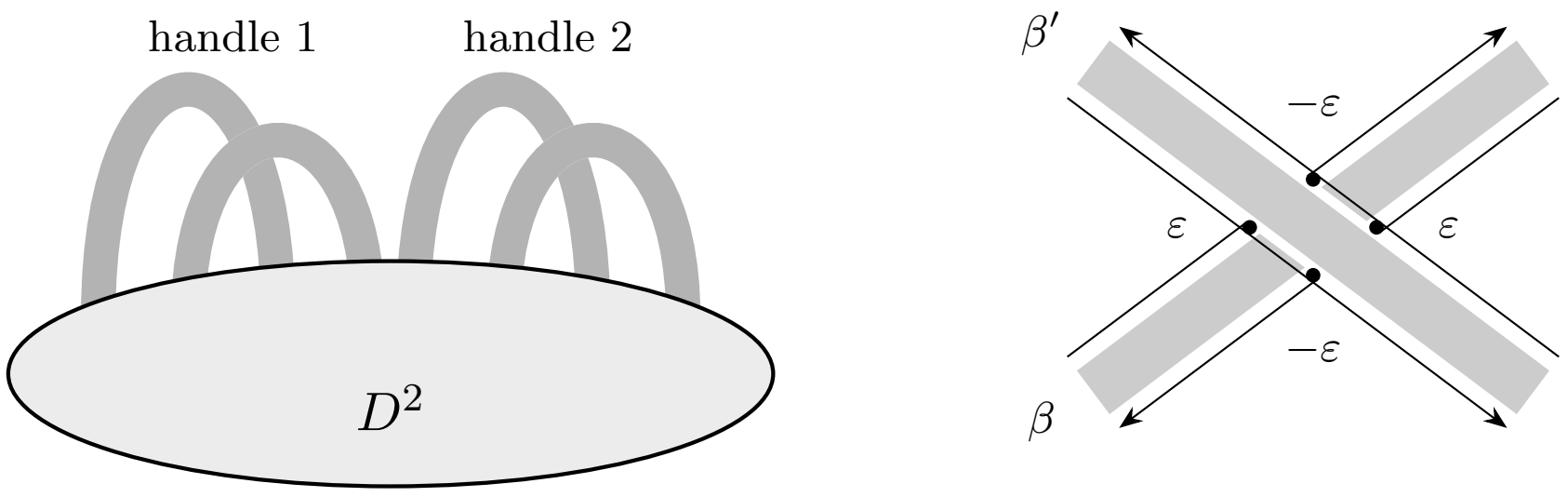
\begin{figure}[ht]
\centering
\begin{tikzpicture}[>=Stealth,scale=0.8,every node/.style={font=\small}]
\draw[line width=7pt,gray!60] (-2.6,0.5) .. controls (-2.6,3.2) and (-1.0,3.2) .. (-1.0,0.5);
\draw[line width=7pt,white] (-1.8,0.5) .. controls (-1.8,2.6) and (-0.2,2.6) .. (-0.2,0.5);
\draw[line width=7pt,gray!60] (-1.8,0.5) .. controls (-1.8,2.6) and (-0.2,2.6) .. (-0.2,0.5);
\draw[line width=7pt,gray!60] (0.2,0.5) .. controls (0.2,3.2) and (1.8,3.2) .. (1.8,0.5);
\draw[line width=7pt,white] (1.0,0.5) .. controls (1.0,2.6) and (2.6,2.6) .. (2.6,0.5);
\draw[line width=7pt,gray!60] (1.0,0.5) .. controls (1.0,2.6) and (2.6,2.6) .. (2.6,0.5);
\fill[gray!15] (0,0) ellipse (3.4 and 1.0);
\draw[thick] (0,0) ellipse (3.4 and 1.0);
\node at (0,-0.3) {$D^2$};
\node[font=\scriptsize] at (-1.4,3.0) {handle 1};
\node[font=\scriptsize] at (1.4,3.0) {handle 2};
\begin{scope}[xshift=8.2cm,yshift=1.3cm,scale=1.15]
\draw[line width=11pt,gray!40] (-1.7,-1.275) -- (1.7,1.275);
\draw[->] (-1.9,-1.0) -- (1.5,1.55); \draw[<-] (-1.5,-1.55) -- (1.9,1.0);
\draw[line width=15pt,white] (-1.7,1.275) -- (1.7,-1.275);
\draw[line width=11pt,gray!40] (-1.7,1.275) -- (1.7,-1.275);
\draw[->] (-1.9,1.0) -- (1.5,-1.55); \draw[<-] (-1.5,1.55) -- (1.9,-1.0);
\foreach \x/\y in {0/0.37, 0/-0.37, 0.49/0, -0.49/0} {\fill (\x,\y) circle (1.6pt);}
\node[font=\scriptsize,fill=white,inner sep=1pt] at (0,0.95) {$-\ve$};
\node[font=\scriptsize,fill=white,inner sep=1pt] at (0,-0.95) {$-\ve$};
\node[font=\scriptsize,fill=white,inner sep=1pt] at (-1.05,0) {$\ve$};
\node[font=\scriptsize,fill=white,inner sep=1pt] at (1.05,0) {$\ve$};
\node[font=\scriptsize] at (-2.1,1.5) {$\beta'$}; \node[font=\scriptsize] at (-2.1,-1.5) {$\beta$};
\end{scope}
\end{tikzpicture}
\caption{Left: a disk with two handles, each made of two bands attached in the interleaved pattern of \cite{Kohli-Tahar}; the bands may cross, twist and link. Right: where a band passes over another, the two antiparallel edges of each band produce four crossings of the link diagram, of signs $\ve,-\ve,-\ve,\ve$; the variables of two crossings of opposite signs have $\ve$-degree $0$ as a pair, and $1-x_ax_b$ is in the kernel of $\BE_D$ against most monomials of the second flavor.}
\label{fig.discband}
\end{figure}

Theorem~\ref{thm.canonical} is the genus bound for the canonical surface of the diagram. The bound of \cite{Kohli-Tahar, LNVdV25} holds for every Seifert surface, and its proof presents an arbitrary Seifert surface $F$ as the thickening of a tangle with $2g$ components and uses a basis of the tensor product $V\otimes V^*$ carried by the two edges of a band in which the braiding is independent of the parameter $\alpha$ \cite[Theorem~2.13]{Kohli-Tahar}, so that the degree is carried by the $4g$ places where the bands are attached to the disk. One may ask whether the determinant formula sees this stronger bound term by term, i.e.\ whether for the diagram $D=\partial F$ of a disk with bands and a suitable cut edge the $\ve$-degrees of $P_{D,\hate}$ and $\tP_{D,\hate}$ satisfy $\deg^+_\ve\tP-\deg^-_\ve P\le2g+\mu-1$. They do not. For the diagrams of Figure~\ref{fig.discband}, a disk with four bands forming two handles of genus one, with one crossing between the two bands of each handle and two crossings (a clasp or a pass) between bands of different handles, we computed the $\ve$-degrees of the two edge polynomials for all $2c$ cut edges: the smallest value of $\max(\deg^+_\ve\tP-\deg^-_\ve P,\ \deg^+_\ve P-\deg^-_\ve\tP)$ over the cut edges is $5=2g+1$, whereas $\deg_p\LG\in[-4,4]$ by \cite{Kohli-Tahar}. The same happens for the expansion of Corollary~\ref{cor.expansion} around the free fermion point, whose terms $g^Sa_S(t_0)b_S(t_1)$ have degrees exceeding $2g$. So the bound of \cite{Kohli-Tahar} is a cancellation inside $\BE_D$: the extremal monomials of $P_{D,\hate}\tP_{D,\hate}$ are killed by the coupling.

The mechanism of the cancellation is visible on a single crossing of two bands. The two edges of a band are antiparallel, so where a band $\beta$ passes over a band $\beta'$ the diagram has four crossings, of signs $\ve,-\ve,-\ve,\ve$ (Figure~\ref{fig.discband}, right). For two crossings $a,b$ of opposite signs, the $\ve$-degree of $x_ax_b$ is $0$, and $\BE_D[(1-x_ax_b)\,x^Sy^T]=0$ whenever $a,b\notin S$ and $T$ contains both or none of $a,b$, because $\BE_D[x_ax_b\,y^T]=t_0^{\ve_a+\ve_b}\BE_D[y^T]=\BE_D[y^T]$ in the first case and $\BE_D[x_ax_by_ay_b\,y^{T'}]=\BE_D[y_ay_by^{T'}]$ in the second. In the examples the extremal monomials of $P_{D,\hate}$ come in pairs $x^S$, $x^{S\cup\{a,b\}}$ with opposite coefficients, i.e.\ the extremal part of $P_{D,\hate}$ is divisible by factors $1-x_ax_b$ attached to pairs of crossings of a band crossing, and these factors are killed by $\BE_D$ against all but a few monomials of $\tP_{D,\hate}$, which are themselves not extremal. A proof of the bound of \cite{Kohli-Tahar} along these lines would require controlling this divisibility for the two flavors simultaneously.

Proposition~\ref{prop.hollowhandle} makes the divisibility a theorem in the elementary case: when a hollow handle is attached to the canonical surface of a closed braid by inserting $\sigma_i\sigma_i^{-1}$ into the braid word, the whole of $P_{D,\hate}$ changes by $(x_ax_b-1)R(x)$ and the whole of $\tP_{D,\hate}$ by $(y_ay_b-1)\tilde R(y)$, with $R,\tilde R$ free of the variables of the two new crossings; the monomials of the enlarged diagram whose $\ve$-degrees leave the interval of the old one come from these terms, and the coupling kills them. By Rudolph's theorem \cite{Rudolph83}, every Seifert surface $F$ is isotopic to a braided surface whose bands may join nonadjacent disks, and the canonical surface of the closed braid obtained by expanding the band generators is $F$ with one hollow handle for each disk passed by a band (Remark~\ref{rem.markovetc} (ii)). So the genus bound of \cite{Kohli-Tahar, LNVdV25} for $F$ would follow from the determinant formula if one could control, in the expanded braid word, the divisibility by the factors $x_ax_b-1$ and $y_ay_b-1$ attached to the passages of the bands over the intermediate strands. These canceling pairs are not adjacent in the word, so Proposition~\ref{prop.hollowhandle} does not apply to them as it stands, and the naive extension of the matrices \eqref{eq.PVV} to the braided surface itself, with one pair of variables per band and none for the passages, does not compute $\LG$ (Remark~\ref{rem.markovetc} (ii)). We have not found such a proof, and we leave it as a question whether the $\alpha$-free basis of \cite{Kohli-Tahar} has a counterpart in the free-fermion model which makes the cancellation visible. Another route from the determinant formula to the bound for arbitrary Seifert surfaces, through the state sum of Section~\ref{sec.kauffman} and the state model of Ballinger \cite{BallingerTalk}, is described after Corollary~\ref{cor.ballinger}.

The passages can be followed one step further, and the computation shows what distinguishes them from hollow handles. Let the band $k$ of $F$ pass in front of the disk $D_m$, and let $g$ and $r$ be the two crossings of its edges with $\partial D_m$, the first positive and the second negative, consecutive along $\partial D_m$; the edge $e$ of $\partial D_m$ between them is the outgoing under-edge of $g$ and the incoming under-edge of $r$. In the cut matrix of the first flavor the row of $e$ is the row of the block $\rho(x_g)$ which contains $x_g$, and the column of $e$ has exactly two nonzero entries, its diagonal entry $1$ and the entry $-x_r$ in the row of the outgoing under-edge of $r$, because the column of $\rho^-(x_r)$ at the incoming under-edge is $(0,x_r)^{\T}$; a term of the determinant which takes $x_g$ from the row of $e$ must cover the column of $e$ by that entry. Hence every monomial of $P_{D,\hate}$ containing $x_g$ contains $x_r$,
\[
P_{D,\hate}=P_{00}+x_rP_{01}+x_gx_rP_{11},
\]
whereas the transposed blocks have no such column and $\tP_{D,\hate}=\sum_{\beta}y_g^{\beta_1}y_r^{\beta_2}\tP_{\beta}$ has all four coefficients in general. Writing $P(a,b)$ for $P_{D,\hate}$ at $x_g=a$, $x_r=b$, the coupling on the pair $g,r$ is
\[
\BE_{g,r}[P_{D,\hate}\tP_{D,\hate}]=P(t_0,t_0^{-1})\,\tP_{00}+t_1P(t_1^{-1},t_0^{-1})\,\tP_{10}+t_1^{-1}P(t_0,t_1)\,\tP_{01}+P(t_1^{-1},t_1)\,\tP_{11},
\]
and since $P(a,b)$ depends on $a$ only through $ab$,
\begin{multline*}
\BE_{g,r}[P_{D,\hate}\tP_{D,\hate}]=P(t_0,t_0^{-1})\bigl(\tP_{00}+t_1\tP_{10}\bigr)+P(t_1^{-1},t_1)\bigl(\tP_{11}+t_1^{-1}\tP_{01}\bigr)\\
+(q^2-1)\,P_{11}\bigl(t_1^{-1}\tP_{01}-t_0^{-1}\tP_{10}\bigr).
\end{multline*}
For an adjacent pair $\sigma_i\sigma_i^{-1}$ one has $P_{01}=\tP_{10}=\tP_{01}=0$ as well and the right-hand side is $(P_{00}+P_{11})(\tP_{00}+\tP_{11})$, which is Proposition~\ref{prop.hollowhandle}. For a passage the mixed coefficients $\tP_{10}$, $\tP_{01}$ of the second flavor are present: they carry the last term, the trace of the handle, and they evaluate the first flavor at two different values, $t_0$ and $t_1^{-1}$, of the passage variable $x_g$ on the locus $x_gx_r=1$, on which the product of the blocks of the band generator is a conjugate of the block of its twist, so that $P$ is the determinant of a matrix with one row per band of $F$ whose entries carry the passage weight. We computed the bilinear form $\BE_{\rm passages}[P_{D,\hate}\tP_{D,\hate}]$ in the band variables for some two hundred random braided surfaces with three to five disks: its rank is $1$ for a majority of them, but it is as large as $3$ for a single passage and $5$ for three, so the passages cannot be integrated out one by one into a product $A(x)B(y)$ of two polynomials in the band variables; and the $p$-degrees of the three terms above can exceed $\nu(F)$ individually, like those of $P_{D,\hate}\tP_{D,\hate}$, and respect the bound only in the sum. A determinant formula attached to $F$ would therefore have to reorganize the coupling globally rather than passage by passage.

\emph{The determinant of the pencil.} By Theorem~\ref{thm.mixed}, $\LG(L)$ is the coefficient of $\lambda_1\cdots\lambda_c$ in the determinant $\CDp_{D,\hate}(\lambda)$ of the pencil of the coupled crossing forms (Definition~\ref{def.pencil}), a determinant of a matrix whose entries are polynomials in $t_0^{\pm1}$, $t_1^{\pm1}$ and $\lambda$. One may ask whether the bound of \cite{Kohli-Tahar} is a property of this determinant. It is not: the bound distinguishes the coefficient of $\lambda_1\cdots\lambda_c$ from the other ones. For the closure of $\sigma_1^3\sigma_2\sigma_1^{-1}\sigma_2$, the diagram of $5_2$ of Remark~\ref{rem.loose}, for which $\nu(F_D)=4$, the degrees in $p$ of the monomials of $\LG$ lie in $[-2,2]$, as they must for a knot of genus one, whereas those of $\CDp_{D,\hate}(\lambda)$ range from $-4$ to $3$ when $\hate$ is the edge which closes the first strand of the braid, and they leave $[-2,2]$ for eleven of the twelve cut edges. In all the examples which we computed ($124$ pairs of a diagram and a cut edge), the span in $p$ of $\CDp_{D,\hate}(\lambda)$ is at most $2\nu(F_D)$, the bound of the canonical surface; we have not tried to prove this.

\section{The complement of the co-cores: a local system for the coupling}\label{sec.cocores}

In the last three sections each flavor of the determinant formula was read on the canonical surface $F_D$: the crossing forms are the Seifert form of $F_D$, cut into the contributions of its bands (Theorems~\ref{thm.localization} and~\ref{thm.deformedseifert}), and the coupling $\BE_D$ was applied at the end, to the product of two determinants. We saw in Section~\ref{ssec.cochains} that the crossing variables are cochains on this surface, and that the formula is not a function of the local system on the surface which they define. In this section we show that the coupling does come from a local system of rank one, on another space: the complement
\[
   X_G=\BR^3\sm G,\qquad G=L\cup\alpha_1\cup\dots\cup\alpha_c,
\]
of the link together with the co-cores of the bands of $F_D$ (Section~\ref{ssec.statesurface}). The lattice of flows is the group $H_1(G;\BZ)$, by \eqref{eq.flowsG}, and the linear forms on it are given by the lifted $1$-chains of $X_G$, through a twisted version of the linking number (Lemma~\ref{lem.iota}). Under this dictionary:
\begin{itemize}
\item the two crossing forms $\pi_j$ and $\tpi_j$ are lifts of one and the same closed curve, the meridian of the co-core $\alpha_j$: the lift starts below the over-strand for $\pi_j$ and above the under-strand for $\tpi_j$, and the crossing variable is the holonomy of the local system around an edge of the other strand (Proposition~\ref{prop.readings});
\item the coupled crossing form $\Psi_j$ of Section~\ref{ssec.mixed} is an \emph{interleaved} product of these two lifts, in which the starting point of the first lift depends on the position of the second one (Proposition~\ref{prop.interleaved});
\item for the local system whose holonomy around an edge of $L$ is $t_1$ or $t_0^{-1}$, according to a coloring of the Seifert circles of $D$ with two colors, the Links--Gould polynomial is the mixed discriminant of the interleaved chains of the $c$ co-cores (Theorem~\ref{thm.cocores}).
\end{itemize}
The holonomy of this local system around every co-core is $(t_0t_1)^{\pm1}$. So it extends to the complement of the link exactly at the free fermion point $t_0t_1=1$, where it is the local system of the Alexander polynomial; and the rule $x_jy_j\mapsto1$, which distinguishes $\BE_D$ from an evaluation of the crossing variables, says that its holonomy around every band of $F_D$ is trivial (Remark~\ref{rem.dictionary}). This is a description of the coupling in the complement of the link and of the co-cores. It is not a definition of $\LG$ from the Seifert surface: the construction also uses the plane of the diagram, through the paths which join the meridians of the co-cores to two base points, and we show in Section~\ref{ssec.vertical} that this dependence cannot be removed in an obvious way.

Throughout this section $D$ is a diagram, placed in the plane $P\subset\BR^3$ as in Section~\ref{ssec.statesurface}, with its set of faces $\CF$, and $\kk$ is a commutative ring; in the applications $\kk=\BZ[t_0^{\pm1},t_1^{\pm1}]$. We use the mixed discriminants $\lb\,\cdot\,\rb$ and the forms $\Psi_j[\sv;\xs,\ys]$ of Section~\ref{ssec.mixed}.

\subsection{The complement of the graph \texorpdfstring{$G$}{G}}\label{ssec.XG}
The aim of this subsection is a combinatorial description of the paths in $X_G$: up to homotopy they are words in $c+2$ letters, one for each face of the diagram (Lemma~\ref{lem.words}). The reason is that the faces, lifted as in Section~\ref{ssec.spanning}, form a surface $\Sigma$ which contains $G$ and which cuts $X_G$ into two contractible pieces, the part above it and the part below it. So a path is determined, up to homotopy, by the sequence of the faces through which it passes from one piece to the other one.

We first describe $X_G$ by means of the disks $\Delta_f$ of Section~\ref{ssec.spanning}. These were defined for the bounded faces. For the unbounded face $f$ we choose in the same way a continuous function $\psi_f\ge0$ on its closure which vanishes exactly on its boundary, and which is bounded; we define $u_f$ by \eqref{eq.uf} and let $\Delta_f$ be the closure of its graph. As in the proof of Lemma~\ref{lem.spanning}~(i), for every face $f$ the closure $\Delta_f$ of the graph of $u_f$ consists of this graph, which we denote by $\Delta^\circ_f$, of the points of $L$ over the edges on the boundary of $f$, and of the co-cores at the corners of $f$. Let
\[
   \Sigma=\bigcup_{f\in\CF}\Delta_f .
\]
Let $u$ be the function on the plane $P$ without the crossings which is equal to $u_f$ on every face $f$ and whose value at a point $p$ of $D$ is the height of the point of $L$ over $p$. It is continuous, by the proof of Lemma~\ref{lem.spanning}~(i), it is bounded, and when $p$ tends to a crossing all the limits of $u(p)$ lie in $[0,h_0]$. So $\Sigma$ is the union of the graph of $u$ and of the co-cores: over a point of $P$ which is not a crossing $\Sigma$ has exactly one point, over a crossing $j$ it contains exactly the co-core $\alpha_j$, and
\[
   G\subset\Sigma,\qquad \Sigma\sm G=\bigsqcup_{f\in\CF}\Delta^\circ_f .
\]
We think of $\Sigma$ as the plane of the diagram, in which every crossing has been pulled apart into its co-core. We say that a point $(p,h)$ of $\BR^3$ lies \emph{above} $\Sigma$ if $h>u(p)$, or if $p$ is a crossing and $h>h_0$, and \emph{below} $\Sigma$ if $h<u(p)$, or if $p$ is a crossing and $h<0$. We denote by $X_+$ and $X_-$ the sets of the points above and below $\Sigma$.

\begin{lemma}[the two sides of $\Sigma$]\label{lem.XG}
\begin{enumerate}\renewcommand{\theenumi}{\roman{enumi}}
\item The space $X_G$ is the disjoint union of $X_+$, of $X_-$ and of the $c+2$ surfaces $\Delta^\circ_f$, $f\in\CF$, each of which is connected.
\item The sets $X_+$ and $X_-$ are open and contractible.
\item Let $W\subset P$ be an open disk contained in a face $f$. Then the set $W\times\BR$ of the points of $\BR^3$ over $W$ is a contractible open subset of $X_G$, so are its intersections with $X_+$ and with $X_-$, and its other points lie on $\Delta^\circ_f$. Every point of $\Delta^\circ_f$ lies in such a set.
\end{enumerate}
\end{lemma}
\begin{proof}
(i) Let $(p,h)\in\BR^3$. If $p$ is a crossing $j$, the point lies above $\Sigma$, below $\Sigma$ or on $\alpha_j\subset G$. Otherwise it lies above $\Sigma$, below $\Sigma$, or it is the point $(p,u(p))$ of $\Sigma$, which belongs to $L$ if $p\in D$ and to $\Delta^\circ_f$ if $p$ lies in the face $f$. The surface $\Delta^\circ_f$ is homeomorphic to $f$, which is connected.

(ii) Let $u^+$ be the function on $P$ which is equal to $u$ outside the crossings and to $h_0$ at the crossings. It is upper semicontinuous, because the limits of $u$ at a crossing are at most $h_0$; so $X_+=\{(p,h):h>u^+(p)\}$ is open. Let $h_1$ be an upper bound of $u^+$. The maps $(p,h)\mapsto\bigl(p,\,h+s\max(0,h_1+1-h)\bigr)$, $0\le s\le1$, send $X_+$ into itself, since they do not lower any point; they deform $X_+$ into the half-space $\{h\ge h_1+1\}$, which they fix pointwise. So $X_+$ is contractible. For $X_-$ the argument is the same, with the function which is equal to $u$ outside the crossings and to $0$ at the crossings, which is lower semicontinuous and nonnegative, and with the half-space $\{h\le-1\}$.

(iii) The set $W\times\BR$ is open and contractible, and it does not meet $G$, whose vertical projection is $D$. Its points above $\Sigma$ are the points $(p,h)$ with $p\in W$ and $h>u_f(p)$; they form a set homeomorphic to $W\times(0,\infty)$, and likewise for its points below $\Sigma$; the other points of $W\times\BR$ lie on $\Delta^\circ_f$. A point of $\Delta^\circ_f$ lies over a point of the open set $f$, hence over a small disk contained in $f$.
\end{proof}

\emph{The base points and the paths $\pth_f$.}\label{not.basepoints} We fix two base points $\star_+\in X_+$ and $\star_-\in X_-$. For a face $f$ we denote by $\pth_f$ a path in $X_G$ from $\star_+$ to $\star_-$ which meets $\Sigma$ exactly once, at a point of $\Delta^\circ_f$: it goes \emph{down through the face $f$} (Figure~\ref{fig.XG}, left). Its homotopy class, relative to its ends, does not depend on the choices. Indeed, such a path consists of a path in $X_+$, of a path which crosses $\Delta^\circ_f$ at a point $q$ inside a set $W\times\BR$ of Lemma~\ref{lem.XG}~(iii), and of a path in $X_-$; the first and the last one are determined up to homotopy by their ends, because $X_+$ and $X_-$ are simply connected, and the middle one is homotopic in $W\times\BR$ to a vertical segment through $q$. When $q$ moves along a path in $\Delta^\circ_f$, which is connected, the vertical segments through its points form a homotopy, in which the upper ends stay in $X_+$ and the lower ends in $X_-$. We denote by $\pth_f^{-1}$ the path $\pth_f$ traveled backward, from $\star_-$ to $\star_+$: it goes up through $f$. A \emph{word} is a sequence of letters $\pth_f^{\pm1}$ in which the exponents $1$ and $-1$ alternate; it defines a path in $X_G$, the product of its letters, with its ends in $\{\star_+,\star_-\}$.

\begin{figure}[ht]
\centering
\begin{tikzpicture}[>=Stealth,scale=1.0,every node/.style={font=\small}]
\begin{scope}
\draw[thick] (-2.9,0) -- (2.9,0);
\foreach \x in {-2.2,0,2.2} {\fill[white] (\x,0) circle (3.2pt); \draw[thick] (\x,0) circle (3.2pt);}
\draw[thick] (-0.08,-0.08) -- (0.08,0.08); \draw[thick] (-0.08,0.08) -- (0.08,-0.08);
\node[above] at (0,0.12) {$e$};
\node[below,font=\scriptsize] at (-0.68,-0.05) {$\Delta_{f_{\rm L}(e)}$}; \node[below,font=\scriptsize] at (0.72,-0.05) {$\Delta_{f_{\rm R}(e)}$};
\fill (0,1.7) circle (1.6pt) node[above] {$\star_+$}; \fill (0,-1.9) circle (1.6pt) node[below] {$\star_-$};
\draw[blue,semithick,->] (0,1.7) .. controls (-1.3,1.3) and (-1.3,0.6) .. (-1.3,0.12);
\draw[blue,semithick] (-1.3,0.12) .. controls (-1.3,-0.9) and (-1.3,-1.5) .. (0,-1.9);
\draw[blue,semithick,->] (0,1.7) .. controls (1.3,1.3) and (1.3,0.6) .. (1.3,0.12);
\draw[blue,semithick] (1.3,0.12) .. controls (1.3,-0.9) and (1.3,-1.5) .. (0,-1.9);
\node[blue,left] at (-1.25,0.95) {$\pth_{f_{\rm L}(e)}$}; \node[blue,right] at (1.25,0.95) {$\pth_{f_{\rm R}(e)}$};
\node at (-2.6,1.3) {$X_+$}; \node at (-2.6,-1.3) {$X_-$}; \node[right] at (2.9,0) {$\Sigma$};
\end{scope}
\begin{scope}[xshift=9.0cm,yshift=-1.45cm,scale=1.22]
\draw[thick,->] (-2.75,0) -- (2.75,0);
\node[below] at (-2.2,0) {$u_j$}; \node[below] at (2.2,0) {$u'_j$};
\draw[blue,semithick,densely dashed] plot[domain=0:180,samples=50,variable=\t] ({1.15*cos(\t)+0.713*sin(\t)},{1.1+0.414*sin(\t)});
\draw[line width=2.2pt] (0,0) -- (0,2.2);
\draw[line width=4pt,white] (-1.364,1.408) -- (1.364,2.992);
\draw[thick] (-1.364,1.408) -- (1.364,2.992);
\node[right] at (1.364,3.02) {$\hat e^+_j$}; \node[left] at (-1.45,1.5) {$\hat e^-_j$};
\node[right,font=\scriptsize] at (0.05,0.3) {$\alpha_j$};
\draw[blue,semithick,postaction={decorate},decoration={markings,mark=at position 0.3 with {\arrow{Stealth}}}] plot[domain=0:-180,samples=50,variable=\t] ({1.15*cos(\t)+0.713*sin(\t)},{1.1+0.414*sin(\t)});
\fill[blue] (-1.15,1.1) circle (1.9pt); \fill[blue] (1.15,1.1) circle (1.9pt);
\draw[blue,fill=white] (0.713,1.514) circle (1.9pt); \draw[blue,fill=white] (-0.713,0.686) circle (1.9pt);
\node[blue,font=\scriptsize] at (-1.2,0.84) {$0$}; \node[blue,font=\scriptsize] at (1.22,0.86) {$2$};
\node[blue,font=\scriptsize] at (0.62,1.72) {$1$}; \node[blue,font=\scriptsize] at (-0.72,0.47) {$3$};
\fill (0,3.8) circle (1.6pt) node[above] {$\star_+$}; \fill (0,-1.3) circle (1.6pt) node[below] {$\star_-$};
\draw[densely dotted] (-1.15,1.1) .. controls (-2.7,1.6) and (-1.9,3.5) .. (0,3.8);
\draw[densely dotted] (1.15,1.1) .. controls (2.6,1.3) and (2.5,3.6) .. (0,3.8);
\draw[densely dotted] (-0.713,0.686) .. controls (-1.8,0.3) and (-1.5,-1.0) .. (0,-1.3);
\draw[densely dotted] (0.713,1.514) .. controls (2.4,1.0) and (1.9,-1.0) .. (0,-1.3);
\end{scope}
\end{tikzpicture}
\caption{Left: a section of $X_G$ by a vertical plane which crosses an edge $e$ of $L$; the edge points away from the reader, the face on its left is on the left. The path $\pth_f$ goes from $\star_+$ down through the disk $\Delta_f$ to $\star_-$, and the loop $\pth_{f_{\rm L}(e)}\pth_{f_{\rm R}(e)}^{-1}$ goes once around $e$. Right: the co-core $\alpha_j$, between the under-strand, which lies in the plane of the figure, and the over-strand, which crosses this plane; the edge $\hat e^+_j$ of the over-strand lies behind the plane and $\hat e^-_j$ in front of it. The meridian $\mer_j$ of $\alpha_j$ is represented by a small horizontal circle around $\alpha_j$ (Lemma~\ref{lem.meridian}). It passes above the under-strand at the points $0$ and $2$, which are joined to $\star_+$ inside $X_+$, and below the over-strand at the points $1$ and $3$, which are joined to $\star_-$ inside $X_-$ (dotted).}
\label{fig.XG}
\end{figure}

\begin{lemma}[paths are words in the faces]\label{lem.words}
\begin{enumerate}\renewcommand{\theenumi}{\roman{enumi}}
\item Every path in $X_G$ with its ends in $\{\star_+,\star_-\}$ is homotopic, relative to its ends, to the path of a word.
\item Let $W\subset P$ be an open disk contained in a face $f$, and let $\gamma$ be a path in $W\times\BR$ from a point of $X_+$ to a point of $X_-$. Then the path from $\star_+$ to $\star_-$ which consists of a path in $X_+$, of $\gamma$ and of a path in $X_-$ is homotopic to $\pth_f$.
\item Let $\gamma$ be a loop in $X_G$ which is the path of the word $\pth_{f_1}^{d_1}\cdots\pth_{f_n}^{d_n}$, and let $z$ be a flow with potential $k$, viewed as a $1$-cycle of $G$ by \eqref{eq.flowcycle}. Then
\[
   \lk(z,\gamma)=-\sum_{i=1}^nd_i\,k_{f_i} .
\]
In particular $\lk(\partial f,\gamma)$ is minus the sum of the exponents $d_i$ of the letters with $f_i=f$.
\end{enumerate}
\end{lemma}
\begin{proof}
(i) By Lemma~\ref{lem.XG}, the sets $X_+$ and $X_-$ and the sets $W\times\BR$ of its part (iii) form an open cover of $X_G$. So the path can be cut into finitely many arcs, each of which lies in one of these sets. If an end of an arc lies on $\Sigma$, at a point of some $\Delta^\circ_f$, we move it slightly along the path or, if the path stays on $\Sigma$ near this point, we first push the path off $\Sigma$ there; so we may assume that the ends of the arcs lie in $X_+\cup X_-$. An arc in a set $W\times\BR$ whose two ends lie in $X_+$ is homotopic, relative to its ends, to a path in $(W\times\BR)\cap X_+$, because $W\times\BR$ is simply connected and $(W\times\BR)\cap X_+$ is connected; the same holds for $X_-$. An arc in $W\times\BR$ from $X_+$ to $X_-$ is homotopic in $W\times\BR$ to a path which crosses $\Sigma$ once, at a point of the disk $\Delta^\circ_f$ of the face $f$ which contains $W$. So the path is homotopic to a product of paths in $X_+$, of paths in $X_-$ and of paths which cross one of the surfaces $\Delta^\circ_f$ once; joining the points where it passes from one of these pieces to the next to $\star_+$ or to $\star_-$, inside $X_+$ or $X_-$, we obtain the path of a word.

(ii) In the contractible set $W\times\BR$ the path $\gamma$ is homotopic, relative to its ends, to a path which crosses $\Delta^\circ_f$ once.

(iii) In a word which defines a loop there are as many exponents $1$ as exponents $-1$, so the right-hand side does not change when a constant is added to $k$, and we may assume that $k$ vanishes on the unbounded face. Then $\lk(z,\gamma)=\sum_fk_f\,\Delta_f\cdot\gamma$, the sum running over the bounded faces, by \eqref{eq.lkpotential}. The path of the letter $\pth_f$ crosses $\Delta_f$ once, downward, and no other disk, and $\Delta_f$ is oriented by the upward normal; so $\Delta_f\cdot\gamma$ is minus the sum of the exponents of the letters with $f_i=f$. The last statement is the case $k=e_f$.
\end{proof}

For an edge $e$ of $D$, the word $\pth_{f_{\rm L}(e)}\,\pth_{f_{\rm R}(e)}^{-1}$ is a loop based at $\star_+$ which goes down through the face on the left of $e$ and comes back up through the face on its right (Figure~\ref{fig.XG}, left). By Lemma~\ref{lem.words}~(iii) and \eqref{eq.potential}, its linking number with a flow $z$ is $-(k_{f_{\rm L}(e)}-k_{f_{\rm R}(e)})=-z_e$: it is a meridian of the edge $e$ of $L$, the one which has the linking number $-1$ with $e$.

\begin{remark}[the face generators]\label{rem.dehn}
Let us add to $\BR^3$ a point at infinity; this changes neither the fundamental group of $X_G$ nor its local systems. When every co-core is shrunk to a point, $G$ becomes the diagram $D$, on a sphere $S^2\subset S^3$, and the complement of $G$ in $S^3$ has the homotopy type of $S^3\sm D$, since $G$ and $D$ have the same regular neighborhoods. Now $S^3\sm D$ is the union of two open balls and of the $c+2$ faces of $D$, which are open disks, and it retracts by deformation onto a graph with two vertices, the centers of the balls, and $c+2$ edges, one through each face. So the group $\pi_1(X_G,\star_+)$ is free of rank $c+1$, and by Lemma~\ref{lem.words}~(i) it is generated by the $c+1$ loops $\pth_f\,\pth_{f_0}^{-1}$, $f\ne f_0$, where $f_0$ is any face; hence these loops are a basis. These are the generators of the Dehn presentation of the group of $L$ which was mentioned in Section~\ref{ssec.flowlattice}; its relators are the meridians of the co-cores, which we now describe. We shall not use this remark in the proofs.
\end{remark}

\subsection{Local systems of rank one and their chains}\label{ssec.localsystems}
A local system of rank one on $X_G$ is given by one unit for each face of the diagram: the factor by which a vector of the fiber is multiplied when it is carried down through this face. We call such a family of units a multiplicative potential. In this subsection we define the local system of a multiplicative potential and the chains with coefficients in it which are carried by the paths $\pth_f$, and we show that these chains are the same thing as the linear forms on the lattice of flows (Lemma~\ref{lem.iota}). This is the dictionary which is used in the rest of the section.

\begin{definition}[multiplicative potentials]\label{def.potential}
A \emph{multiplicative potential} on $D$ with values in $\kk$ is a family $\vth=(\vth_f)_{f\in\CF}$ of units of $\kk$, one for each face. Its \emph{holonomy around an edge} $e$ of $D$ is the unit
\[
   \mhol_e=\frac{\vth_{f_{\rm L}(e)}}{\vth_{f_{\rm R}(e)}}\ \in\kk^\times .
\]
\end{definition}

This is the multiplicative form of the potentials of Lemma~\ref{lem.faces}: if $k\in\BZ^{\CF}$ is a potential of a flow $z$ and $t$ is a unit of $\kk$, then $\vth_f=t^{k_f}$ is a multiplicative potential, with $\mhol_e=t^{z_e}$ by \eqref{eq.potential}. In general $\mhol_e$ is the product of the units $\vth_f^{\,m(e)_f}$, where $m(e)$ is the vector \eqref{eq.m}, so that \eqref{eq.vertexrelation} gives, at every crossing,
\begin{equation}\label{eq.mvertex}
   \mhol_{o_j}\,\mhol_{u_j}=\mhol_{o'_j}\,\mhol_{u'_j} .
\end{equation}

\emph{The local system.} For a loop $\gamma$ in $X_G$ we put
\begin{equation}\label{eq.holonomy}
   \rho_\vth(\gamma)=\prod_{f\in\CF}\vth_f^{\,-\lk(\partial f,\gamma)}\ \in\kk^\times .
\end{equation}
Since the linking number with a fixed cycle of $G$ is additive on the products of loops and invariant under their deformations in $X_G$, this is a homomorphism from the fundamental group of $X_G$ to $\kk^\times$. We denote by $\CL_\vth$ the local system of free $\kk$-modules of rank one on $X_G$ with this holonomy. By Lemma~\ref{lem.words}~(iii), the holonomy of the loop of a word is the product of the units of its letters,
\begin{equation}\label{eq.holword}
   \rho_\vth\bigl(\pth_{f_1}^{d_1}\cdots\pth_{f_n}^{d_n}\bigr)=\prod_{i=1}^n\vth_{f_i}^{\,d_i} :
\end{equation}
the transport of $\CL_\vth$ along $\pth_f$, from $\star_+$ to $\star_-$, is the multiplication by $\vth_f$, when the fibers at $\star_+$ and at $\star_-$ are identified with $\kk$ in a suitable way. For instance $\mhol_e$ is the holonomy of the meridian $\pth_{f_{\rm L}(e)}\pth_{f_{\rm R}(e)}^{-1}$ of the edge $e$, which explains its name. Multiplying all the $\vth_f$ by the same unit does not change $\CL_\vth$, and by Remark~\ref{rem.dehn} every local system of rank one on $X_G$ is of this form.

\emph{Chains.} We shall use the chains of $X_G$ with coefficients in $\CL_\vth$ which are carried by the paths $\pth_f$. They form the following complex, which we define directly.

\begin{definition}[lifts]\label{def.chains}
Let $\vth$ be a multiplicative potential. We denote by $C_\vth$ the free $\kk$-module with a basis $(\tpth_f)_{f\in\CF}$ indexed by the faces, by $C^0$ the free $\kk$-module with the basis $\tstar_+,\tstar_-$, and by $\partial\colon C_\vth\to C^0$ the $\kk$-linear map with
\[
   \partial\,\tpth_f=\vth_f\,\tstar_--\tstar_+ .
\]
The \emph{lift} of the letter $\pth_f$ which starts at $a\,\tstar_+$, where $a\in\kk$, is the chain $a\,\tpth_f$, and it ends at $a\vth_f\,\tstar_-$. The lift of the letter $\pth_f^{-1}$ which starts at $a\,\tstar_-$ is the chain $-a\vth_f^{-1}\,\tpth_f$, and it ends at $a\vth_f^{-1}\,\tstar_+$. The lift of a word which starts at a given multiple of $\tstar_+$ or of $\tstar_-$ is the sum of the lifts of its letters, each of them starting where the previous one ends.
\end{definition}

So the boundary of the lift of a word is its end minus its starting point, and by \eqref{eq.holword} the lift of a loop $\gamma$ which starts at $a\,\tstar_\pm$ ends at $\rho_\vth(\gamma)\,a\,\tstar_\pm$. Here is the meaning of this definition. Let $\tstar_+$ and $\tstar_-$ be generators of the fibers of $\CL_\vth$ at the base points, chosen so that the transport along the paths $\pth_f$ is as in \eqref{eq.holword}; then $\tpth_f$ is the path $\pth_f$ together with the parallel section of $\CL_\vth$ along it which starts at $\tstar_+$, a singular $1$-chain with coefficients in $\CL_\vth$ whose boundary is $\vth_f\tstar_--\tstar_+$, and the lift of a word is the chain of its path with the parallel section which starts at the given point of the fiber. When $\kk$ is the group ring of a free abelian group which is generated by the units $\vth_f$, this is the language of the covering of $X_G$ defined by $\rho_\vth$, whose group of deck transformations is this group: $\tstar_\pm$ are points over $\star_\pm$, $\tpth_f$ is the lift of the path $\pth_f$ which starts at $\tstar_+$, and it ends at the image of $\tstar_-$ under the deck transformation $\vth_f$.

\emph{Chains as linear forms on the flows.} The linking number pairs the loops of $X_G$ with the cycles of $G$, that is, with the flows. The lifted chains are paired with the flows in two ways, according to the side of $\Sigma$ from which the coefficients are read.

\begin{definition}\label{def.iota}
For a chain $v=\sum_fv_f\,\tpth_f\in C_\vth$ and a vector $k\in\kk^{\CF}$ we put
\[
   \iota(v)(k)=-\sum_{f\in\CF}v_f\,k_f,\qquad \iota'(v)(k)=-\sum_{f\in\CF}v_f\,\vth_f\,k_f .
\]
We denote by $C^-_\vth$ the submodule of the chains $v$ with $\partial v\in\kk\,\tstar_-$, and by $C^+_\vth$ the submodule of the chains with $\partial v\in\kk\,\tstar_+$.
\end{definition}

Since $\partial v=\bigl(\sum_fv_f\vth_f\bigr)\tstar_--\bigl(\sum_fv_f\bigr)\tstar_+$, the chains of $C^-_\vth$ are those with $\sum_fv_f=0$, and the chains of $C^+_\vth$ those with $\sum_fv_f\vth_f=0$. The lift of a word from $\star_-$ to $\star_-$ belongs to $C^-_\vth$, and the lift of a word from $\star_+$ to $\star_+$ belongs to $C^+_\vth$.

\begin{lemma}[chains are linear forms on the flows]\label{lem.iota}
\begin{enumerate}\renewcommand{\theenumi}{\roman{enumi}}
\item If $v\in C^-_\vth$, then $\iota(v)(k)$ only depends on the flow $z=\sum_fk_f\,\partial f$, and we write $\iota(v)(z)$ for it. This defines an isomorphism $\iota$ from $C^-_\vth$ onto the module of the linear forms $\FF_D\to\kk$. In the same way $\iota'$ is an isomorphism from $C^+_\vth$ onto this module.
\item If $\vth_f=1$ for all $f$ and $\gamma$ is a loop based at $\star_-$, resp.\ at $\star_+$, which is the path of a word, and $\tilde\gamma$ is its lift which starts at $\tstar_-$, resp.\ at $\tstar_+$, then $\iota(\tilde\gamma)(z)=\lk(z,\gamma)$, resp.\ $\iota'(\tilde\gamma)(z)=\lk(z,\gamma)$, for every flow $z$.
\item Let $e$ be an edge of $D$, with $f=f_{\rm L}(e)$ and $g=f_{\rm R}(e)$. If $v$ is the lift of the loop $\pth_f^{-1}\pth_g$ which starts at $a\,\tstar_-$, then $\iota(v)(z)=a\,\vth_f^{-1}\,z_e$. If $v$ is the lift of the loop $\pth_f\,\pth_g^{-1}$ which starts at $a\,\tstar_+$, then $\iota'(v)(z)=-a\,\vth_f\,z_e$.
\end{enumerate}
\end{lemma}
\begin{proof}
(i) Two vectors $k\in\BZ^{\CF}$ with the same flow $\sum_fk_f\,\partial f$ differ by a constant vector (Lemma~\ref{lem.faces}~(ii)), and $\iota(v)$ vanishes on the constant vectors if and only if $\sum_fv_f=0$. Moreover the linear forms on $\FF_D\cong\BZ^{\CF}/\BZ\1$ with values in $\kk$ are the forms $k\mapsto-\sum_fv_fk_f$ on $\BZ^{\CF}$ with $v\in\kk^{\CF}$ and $\sum_fv_f=0$, each of them for exactly one $v$. For $\iota'$ the argument is the same, the condition being $\sum_fv_f\vth_f=0$.

(ii) For $\vth=1$ the lift of the letter $\pth_f^{d}$ is $d\,\tpth_f$, wherever it starts, so that $\iota(\tilde\gamma)(k)=\iota'(\tilde\gamma)(k)=-\sum_id_i\,k_{f_i}$ for the word $\pth_{f_1}^{d_1}\cdots\pth_{f_n}^{d_n}$; this is $\lk(z,\gamma)$ by Lemma~\ref{lem.words}~(iii).

(iii) The lift of $\pth_f^{-1}\pth_g$ which starts at $a\,\tstar_-$ is $-a\vth_f^{-1}\,\tpth_f+a\vth_f^{-1}\,\tpth_g$, and $\iota$ of it is $k\mapsto a\vth_f^{-1}(k_f-k_g)$; and $k_f-k_g=z_e$ by \eqref{eq.potential}. The lift of $\pth_f\,\pth_g^{-1}$ which starts at $a\,\tstar_+$ is $a\,\tpth_f-a\vth_f\vth_g^{-1}\,\tpth_g$, and $\iota'$ of it is $k\mapsto-a\vth_f(k_f-k_g)$.
\end{proof}

So $\iota$ and $\iota'$ are twisted versions of the linking number, and Lemma~\ref{lem.iota}~(i) is a form of Alexander duality between the complement of $G$ and $G$. Geometrically, let $z$ be a flow with potential $k$; it is the boundary of the $2$-chain $\sum_fk_f\,\Delta_f$. The chain $\tpth_f$ crosses this $2$-chain once, downward, through $\Delta_f$; its coefficient is $1$ before this passage, in $X_+$, and $\vth_f$ after it, in $X_-$. So $\iota(v)(z)$ is the intersection number of the $2$-chain with $v$ in which every passage of $v$ through a disk is counted with the coefficient which $v$ has on the upper side of the disk, and $\iota'(v)(z)$ is this number with the coefficients on the lower side. By (iii), the linear form $z\mapsto z_e$, the flow through an edge, is a lifted meridian of this edge, up to a unit.

\subsection{The meridian of a co-core and the crossing forms}\label{ssec.meridian}
We now apply the dictionary of Lemma~\ref{lem.iota} to the crossing forms. Each co-core has a meridian, a loop of $X_G$ which goes once around it, and the result of this subsection is that the two crossing forms of the crossing $j$ are two lifts of this one loop, which start on the two sides of $\Sigma$ (Proposition~\ref{prop.readings}). We first write the meridian as a word in the four faces around the crossing.

\begin{figure}[ht]
\centering
\begin{tikzpicture}[>=Stealth,scale=0.94,every node/.style={font=\small}]
\def\meridian{%
  \draw[blue,semithick,postaction={decorate},decoration={markings,
     mark=at position 0.125 with {\arrow{<}},mark=at position 0.375 with {\arrow{<}},
     mark=at position 0.625 with {\arrow{<}},mark=at position 0.875 with {\arrow{<}}}] (0,0) circle (0.85);
  \node[below] at (0,-2) {$u_j$}; \node[above] at (0,2) {$u'_j$};}
\def\points{%
  \fill[blue] (0,-0.85) circle (1.7pt); \fill[blue] (0,0.85) circle (1.7pt);
  \draw[blue,fill=white] (-0.85,0) circle (1.7pt); \draw[blue,fill=white] (0.85,0) circle (1.7pt);
  \node[blue,font=\scriptsize] at (0.22,-1.05) {$0$}; \node[blue,font=\scriptsize] at (-1.05,0.22) {$1$};
  \node[blue,font=\scriptsize] at (0.22,1.05) {$2$}; \node[blue,font=\scriptsize] at (1.05,0.22) {$3$};}
\begin{scope}
  \meridian
  \draw[thick] (0,-2) -- (0,-0.17); \draw[thick,->] (0,0.17) -- (0,2);
  \draw[thick,->] (-2.2,0) -- (2.2,0);
  \points
  \node[left] at (-2.2,0) {$\hat e^+_j=o_j$}; \node[right] at (2.2,0) {$o'_j=\hat e^-_j$};
  \node at (-1.45,-1.3) {$f^1_j=f^\downarrow_j$}; \node at (-1.45,1.3) {$f^2_j=f^\leftarrow_j$};
  \node at (1.45,1.3) {$f^3_j=f^\uparrow_j$}; \node at (1.45,-1.3) {$f^4_j=f^\rightarrow_j$};
  \node at (0,-2.85) {$\ve_j=1$};
\end{scope}
\begin{scope}[xshift=8.6cm]
  \meridian
  \draw[thick] (0,-2) -- (0,-0.17); \draw[thick,->] (0,0.17) -- (0,2);
  \draw[thick,->] (2.2,0) -- (-2.2,0);
  \points
  \node[right] at (2.2,0) {$o_j=\hat e^-_j$}; \node[left] at (-2.2,0) {$\hat e^+_j=o'_j$};
  \node at (-1.45,-1.3) {$f^1_j=f^\leftarrow_j$}; \node at (-1.45,1.3) {$f^2_j=f^\uparrow_j$};
  \node at (1.45,1.3) {$f^3_j=f^\rightarrow_j$}; \node at (1.45,-1.3) {$f^4_j=f^\downarrow_j$};
  \node at (0,-2.85) {$\ve_j=-1$};
\end{scope}
\end{tikzpicture}
\caption{The crossing $j$ seen from above, drawn with the under-strand pointing upward, for the two signs: its four edges, the four faces $f^1_j,\dots,f^4_j$ of \eqref{eq.fourfaces}, and the meridian $\mer_j$ of the co-core $\alpha_j$, which is perpendicular to the figure at the crossing. At the full dots $0$ and $2$ the meridian passes above the under-strand, at the hollow dots $1$ and $3$ it passes below the over-strand.}
\label{fig.meridian}
\end{figure}
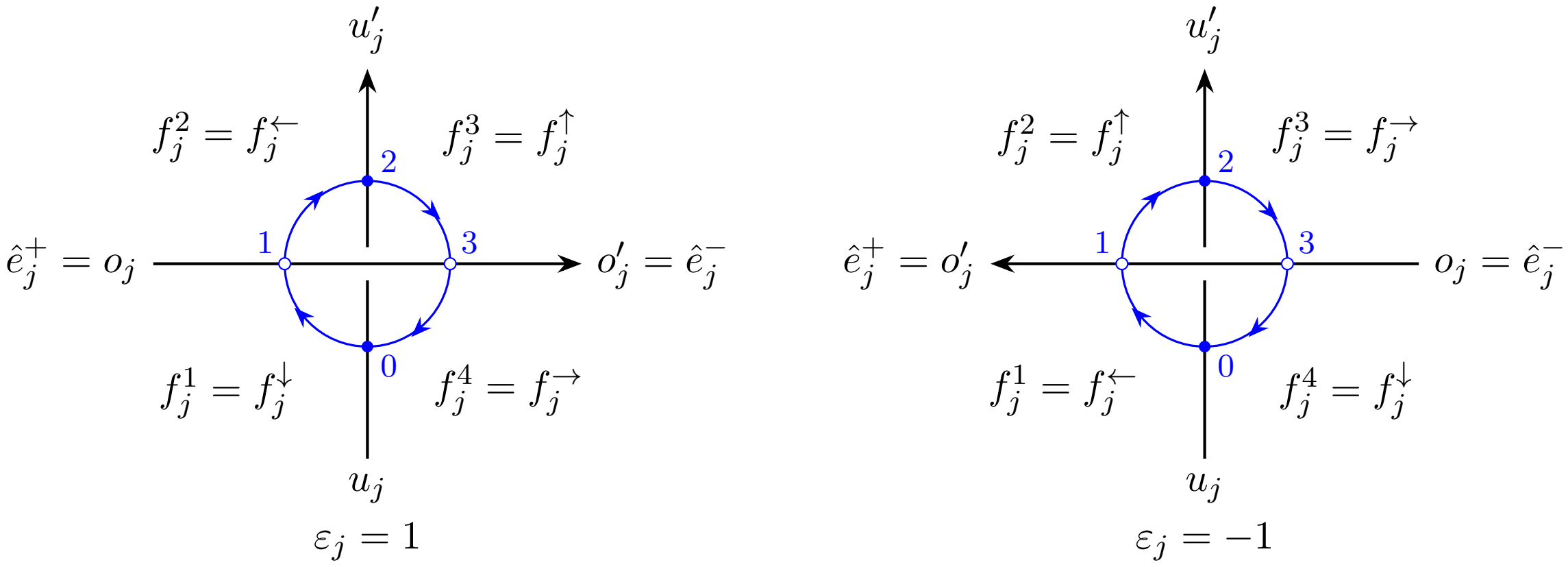

Let $j$ be a crossing. We denote by
\begin{equation}\label{eq.fourfaces}
   f^1_j=f_{\rm L}(u_j),\qquad f^2_j=f_{\rm L}(u'_j),\qquad f^3_j=f_{\rm R}(u'_j),\qquad f^4_j=f_{\rm R}(u_j)
\end{equation}
the four faces around the crossing, in the clockwise order when the diagram is seen from above, starting on the left of the incoming edge of the under-strand (Figure~\ref{fig.meridian}); they need not be distinct. In the notation of Section~\ref{ssec.flowlattice},
\[
   (f^1_j,f^2_j,f^3_j,f^4_j)=\begin{cases}(f^\downarrow_j,f^\leftarrow_j,f^\uparrow_j,f^\rightarrow_j)&\text{if }\ve_j=1,\\[2pt] (f^\leftarrow_j,f^\uparrow_j,f^\rightarrow_j,f^\downarrow_j)&\text{if }\ve_j=-1,\end{cases}
\]
by \eqref{eq.mcrossing}. The under-strand separates $f^1_j,f^2_j$ from $f^3_j,f^4_j$, and the over-strand separates $f^1_j,f^4_j$ from $f^2_j,f^3_j$. As in Proposition~\ref{prop.otherside}, $\hat e^+_j$ is the edge of the over-strand which lies on the left of the under-strand, between $f^1_j$ and $f^2_j$, and $\hat e^-_j$ the one which lies on its right, between $f^3_j$ and $f^4_j$: thus $\hat e^+_j=o_j$ and $\hat e^-_j=o'_j$ if $\ve_j=1$, and $\hat e^+_j=o'_j$ and $\hat e^-_j=o_j$ if $\ve_j=-1$.

\begin{definition}[the meridian of a co-core]\label{def.meridian}
The \emph{meridian} of the co-core $\alpha_j$ is the loop based at $\star_+$
\begin{equation}\label{eq.merword}
   \mer_j=\pth_{f^1_j}\,\pth_{f^2_j}^{-1}\,\pth_{f^3_j}\,\pth_{f^4_j}^{-1} .
\end{equation}
\end{definition}

\begin{lemma}\label{lem.meridian}
Let $\mer^\circ_j$ be a horizontal circle centered at the middle point of the co-core $\alpha_j$, of a radius $r<r_0$ so small that $\chi(r)>\frac12$, traveled clockwise as seen from above. It is a loop in $X_G$ which passes above the under-strand at two points, over $u_j$ and over $u'_j$, and below the over-strand at two points, under $\hat e^+_j$ and under $\hat e^-_j$. Starting from its point over $u_j$, joined to $\star_+$ by a path in $X_+$, it is homotopic to $\mer_j$. For every flow $z$,
\[
   \lk(z,\mer_j)=z_{u'_j}-z_{u_j},
\]
which is minus the coefficient of the co-core $\alpha_j$ in the cycle \eqref{eq.flowcycle}.
\end{lemma}
\begin{proof}
The circle lies at the height $\frac12h_0$, over the circle of radius $r$ around the crossing in $P$. Over this circle the link has four points: two points of the under-strand, at the height $0$, and two points of the over-strand, at the height $h_0\chi(r)>\frac12h_0$. So $\mer^\circ_j$ does not meet $G$; its two points over the under-strand lie above $\Sigma$, in $X_+$, and its two points under the over-strand lie below $\Sigma$, in $X_-$. These four points cut $\mer^\circ_j$ into four open arcs, each of which lies over an open quadrant of the disk $\CB_j$, hence over an open disk contained in one face. Turning clockwise from the incoming half-strand of the under-strand, one meets the quadrants of $f^1_j$, $f^2_j$, $f^3_j$, $f^4_j$ in this order, separated by the half-strands $\hat e^+_j$, $u'_j$, $\hat e^-_j$ (Figure~\ref{fig.meridian}). So the four arcs go from $X_+$ to $X_-$ over $f^1_j$, from $X_-$ to $X_+$ over $f^2_j$, from $X_+$ to $X_-$ over $f^3_j$, and from $X_-$ to $X_+$ over $f^4_j$. Since $X_+$ and $X_-$ are open, each arc, shortened slightly at both ends, is a path as in Lemma~\ref{lem.words}~(ii), up to its direction; so the four arcs are homotopic to the four letters of $\mer_j$, when their ends are joined to the base points inside $X_+$ and $X_-$. By Lemma~\ref{lem.words}~(iii) and \eqref{eq.potential}, $\lk(z,\mer_j)=-(k_{f^1_j}-k_{f^4_j})+(k_{f^2_j}-k_{f^3_j})=-z_{u_j}+z_{u'_j}$.
\end{proof}

So $\mer_j$ goes once around the co-core, between the two strands of the crossing (Figure~\ref{fig.XG}, right). The complement of $L$ is obtained from $X_G$ by putting the co-cores back, and in the complement of $L$ a loop can be moved across $\alpha_j$, which changes its class in $\pi_1(X_G)$ by a conjugate of $\mer_j^{\pm1}$. So a local system of rank one on $X_G$ extends to the complement of $L$ if and only if its holonomy around every co-core is trivial.

\emph{The lifted meridian.} Let $\vth$ be a multiplicative potential. For a crossing $j$ we write $\vth_i=\vth_{f^i_j}$ for $i=1,\dots,4$, and
\begin{equation}\label{eq.omega}
   \omega_j=\rho_\vth(\mer_j)=\frac{\vth_1\vth_3}{\vth_2\vth_4}=\frac{\mhol_{u_j}}{\mhol_{u'_j}}=\frac{\mhol_{o'_j}}{\mhol_{o_j}}
\end{equation}
for the holonomy around the co-core; the last two equalities follow from $\mhol_{u_j}=\vth_1/\vth_4$, $\mhol_{u'_j}=\vth_2/\vth_3$ and \eqref{eq.mvertex}. We lift the meridian traveled again and again, that is, the infinite word $\cdots\mer_j\,\mer_j\,\mer_j\cdots$, in which we number the places between the letters by the integers: the place $0$ is the beginning of a copy of $\mer_j$, so that the even places are at $\star_+$ and the odd places at $\star_-$. Geometrically, the places $0$ and $2$ modulo $4$ are the points of the circle $\mer^\circ_j$ of Lemma~\ref{lem.meridian} over $u_j$ and over $u'_j$, and the places $1$ and $3$ modulo $4$ are its points under $\hat e^+_j$ and under $\hat e^-_j$ (Figures~\ref{fig.XG} and~\ref{fig.meridian}). The lift which starts at $\tstar_+$ at the place $0$ passes at the place $i$ through a point $\star_j(i)$, which is a multiple of $\tstar_+$ if $i$ is even and of $\tstar_-$ if $i$ is odd; we denote its coefficient by $\mathsf c_j(i)\in\kk^\times$. Thus $\star_j(i)=\mathsf c_j(i)\,\tstar_\pm$, with
\begin{equation}\label{eq.liftedpoints}
   \mathsf c_j(0)=1,\qquad \mathsf c_j(1)=\vth_1,\qquad \mathsf c_j(2)=\frac{\vth_1}{\vth_2},\qquad \mathsf c_j(3)=\frac{\vth_1\vth_3}{\vth_2},\qquad \mathsf c_j(i+4)=\omega_j\,\mathsf c_j(i) .
\end{equation}
For integers $i\le l$ we denote by $\mer_j(i,l)\in C_\vth$ the part of this lift between the places $i$ and $l$. By Definition~\ref{def.chains},
\begin{equation}\label{eq.pieces}
\begin{gathered}
   \mer_j(0,1)=\tpth_{f^1_j},\qquad \mer_j(1,2)=-\frac{\vth_1}{\vth_2}\,\tpth_{f^2_j},\qquad \mer_j(2,3)=\frac{\vth_1}{\vth_2}\,\tpth_{f^3_j},\qquad \mer_j(3,4)=-\omega_j\,\tpth_{f^4_j},\\
   \mer_j(i+4,l+4)=\omega_j\,\mer_j(i,l),\qquad \mer_j(i,l)=\mer_j(i,i')+\mer_j(i',l)\quad(i\le i'\le l),
\end{gathered}
\end{equation}
and $\partial\,\mer_j(i,l)=\star_j(l)-\star_j(i)$. We call $\mer_j(i,i+2)$ a \emph{half-turn} and $\mer_j(i,i+4)$ a \emph{turn} of the lifted meridian. A turn which starts at an odd place belongs to $C^-_\vth$, and a turn which starts at an even place to $C^+_\vth$; so do the half-turns. Two turns which start a half-turn apart are different chains: by \eqref{eq.pieces},
\begin{equation}\label{eq.turndiff}
   \mer_j(i+2,i+6)-\mer_j(i,i+4)=\mer_j(i+4,i+6)-\mer_j(i,i+2)=(\omega_j-1)\,\mer_j(i,i+2) .
\end{equation}

\begin{proposition}[the meridian of a co-core reads the crossing forms]\label{prop.readings}
Let $\vth$ be a multiplicative potential and $j$ a crossing; we write $\ve=\ve_j$, $\omega=\omega_j$, $\mer=\mer_j$, $\hat e^\pm=\hat e^\pm_j$, and $u,u'$ for $u_j,u'_j$.
\begin{enumerate}\renewcommand{\theenumi}{\roman{enumi}}
\item The half-turns are, up to units, the flows through the four edges of the crossing: for every flow $z$,
\begin{align*}
   \iota\bigl(\mer(1,3)\bigr)(z)&=\mhol_{\hat e^+}^{-\ve}\,z_{u'}, & \iota\bigl(\mer(3,5)\bigr)(z)&=-\omega\,z_{u},\\
   \iota'\bigl(\mer(0,2)\bigr)(z)&=\ve\,\vth_1\,z_{\hat e^+}, & \iota'\bigl(\mer(2,4)\bigr)(z)&=-\ve\,\frac{\vth_1}{\mhol_{u'}}\,z_{\hat e^-} .
\end{align*}
\item The turns are, up to units, the crossing forms of the two flavors, each of them with two values of its variable. With the coefficients $\mathsf c=\mathsf c_j$ of \eqref{eq.liftedpoints},
\begin{align*}
   \iota\bigl(\mer(1,5)\bigr)&=\mathsf c(4)\,\pi_j^{(\xs)}\quad\text{with }\xs=\mhol_{\hat e^-}^{-\ve}, & \iota\bigl(\mer(3,7)\bigr)&=\mathsf c(4)\,\pi_j^{(\xs)}\quad\text{with }\xs=\mhol_{\hat e^+}^{-\ve},\\
   \iota'\bigl(\mer(0,4)\bigr)&=\mathsf c(2+\ve)\,\tpi_j^{(\ys)}\quad\text{with }\ys=\mhol_{u'}^{\ve}, & \iota'\bigl(\mer(2,6)\bigr)&=\mathsf c(4-\ve)\,\tpi_j^{(\ys)}\quad\text{with }\ys=\mhol_{u}^{\ve}.
\end{align*}
Here $\mathsf c(4)=\omega$, and $\mathsf c(2+\ve)=\mathsf c(4-\ve)=\vth_1/\mhol_{u'}$ if $\ve=1$, whereas $\mathsf c(2+\ve)=\vth_1$ and $\mathsf c(4-\ve)=\omega\,\vth_1$ if $\ve=-1$. In each of the four formulas the unit is the coefficient of the lifted point at the place, inside the turn, where the meridian passes above the incoming under-edge $u_j$, resp.\ below the outgoing over-edge $o'_j$.
\end{enumerate}
\end{proposition}
\begin{proof}
(i) Let $k$ be a potential of $z$, and write $k_i=k_{f^i_j}$. By \eqref{eq.pieces} and Definition~\ref{def.iota},
\[
   \iota\bigl(\mer(1,3)\bigr)(z)=\frac{\vth_1}{\vth_2}\,(k_2-k_3),\qquad \iota\bigl(\mer(3,5)\bigr)(z)=\omega\,(k_4-k_1),
\]
since $\mer(3,5)=\mer(3,4)+\omega\,\mer(0,1)$. By \eqref{eq.fourfaces} and \eqref{eq.potential}, $k_2-k_3=z_{u'}$ and $k_1-k_4=z_u$. The faces $f^1_j$ and $f^2_j$ are the two faces adjacent to $\hat e^+$: if $\ve=1$, then $\hat e^+=o_j$ has $f^2_j$ on its left and $f^1_j$ on its right, so that $\vth_1/\vth_2=\mhol_{\hat e^+}^{-1}$ and $k_2-k_1=z_{\hat e^+}$; if $\ve=-1$, then $\hat e^+=o'_j$ has $f^1_j$ on its left and $f^2_j$ on its right, so that $\vth_1/\vth_2=\mhol_{\hat e^+}$ and $k_2-k_1=-z_{\hat e^+}$ (Figure~\ref{fig.meridian}). In both cases $\vth_1/\vth_2=\mhol_{\hat e^+}^{-\ve}$ and $k_2-k_1=\ve\,z_{\hat e^+}$; this gives the first line. In the same way $f^3_j$ and $f^4_j$ are the two faces adjacent to $\hat e^-$, and
\begin{equation}\label{eq.eminus}
   \frac{\vth_3}{\vth_4}=\mhol_{\hat e^-}^{\ve},\qquad k_4-k_3=-\ve\,z_{\hat e^-} .
\end{equation}
For the second line, Definition~\ref{def.iota} and \eqref{eq.pieces} give
\begin{align*}
   \iota'\bigl(\mer(0,2)\bigr)(z)&=-\vth_1k_1+\frac{\vth_1}{\vth_2}\,\vth_2k_2=\vth_1\,(k_2-k_1),\\
   \iota'\bigl(\mer(2,4)\bigr)(z)&=-\frac{\vth_1}{\vth_2}\,\vth_3k_3+\omega\,\vth_4k_4=\frac{\vth_1\vth_3}{\vth_2}\,(k_4-k_3),
\end{align*}
and $\vth_3/\vth_2=\mhol_{u'}^{-1}$.

(ii) By \eqref{eq.omega}, by the proof of (i) and by \eqref{eq.eminus}, $\omega=(\vth_1/\vth_2)(\vth_3/\vth_4)=\mhol_{\hat e^+}^{-\ve}\,\mhol_{\hat e^-}^{\ve}$. Adding the two formulas of the first line of (i),
\[
   \iota\bigl(\mer(1,5)\bigr)(z)=\mhol_{\hat e^+}^{-\ve}\,z_{u'}-\omega\,z_u=\omega\bigl(\mhol_{\hat e^-}^{-\ve}\,z_{u'}-z_u\bigr),
\]
and since $\mer(3,7)=\mer(3,5)+\omega\,\mer(1,3)$, $\iota\bigl(\mer(3,7)\bigr)(z)=\omega\bigl(\mhol_{\hat e^+}^{-\ve}\,z_{u'}-z_u\bigr)$; and $\omega=\mathsf c(4)$. For the second flavor, the second line of (i) gives
\[
   \iota'\bigl(\mer(0,4)\bigr)(z)=\ve\,\vth_1\Bigl(z_{\hat e^+}-\frac{z_{\hat e^-}}{\mhol_{u'}}\Bigr),\qquad \iota'\bigl(\mer(2,6)\bigr)(z)=\ve\,\vth_1\Bigl(\omega\,z_{\hat e^+}-\frac{z_{\hat e^-}}{\mhol_{u'}}\Bigr),
\]
because $\mer(2,6)=\mer(2,4)+\omega\,\mer(0,2)$. If $\ve=1$, then $(\hat e^+,\hat e^-)=(o_j,o'_j)$, and these are $\frac{\vth_1}{\mhol_{u'}}\bigl(\mhol_{u'}\,z_{o_j}-z_{o'_j}\bigr)$ and $\frac{\vth_1}{\mhol_{u'}}\bigl(\omega\,\mhol_{u'}\,z_{o_j}-z_{o'_j}\bigr)$, where $\omega\,\mhol_{u'}=\mhol_u$ by \eqref{eq.omega}, and $\vth_1/\mhol_{u'}=\vth_1\vth_3/\vth_2=\mathsf c(3)$. If $\ve=-1$, then $(\hat e^+,\hat e^-)=(o'_j,o_j)$, and they are $\vth_1\bigl(\mhol_{u'}^{-1}z_{o_j}-z_{o'_j}\bigr)$ and $\omega\,\vth_1\bigl(\mhol_u^{-1}z_{o_j}-z_{o'_j}\bigr)$, where $\vth_1=\mathsf c(1)$ and $\omega\,\vth_1=\mathsf c(5)$. The places $0$ and $4$ lie over $u_j$, and the place under $o'_j$ is $3$ modulo $4$ if $\ve=1$ and $1$ modulo $4$ if $\ve=-1$; this gives the last statement.
\end{proof}

In words: the first crossing form $\pi_j$ is the meridian of the co-core, lifted from a point below the over-strand, and its variable is the holonomy around an edge of the over-strand, to the power $-\ve_j$; the second crossing form $\tpi_j$ is the same meridian, lifted from a point above the under-strand, and its variable is the holonomy around an edge of the under-strand, to the power $\ve_j$. The two flavors, which Theorem~\ref{thm.deformedseifert} places on the two sides of the plane of the diagram, are here the two sides of $\Sigma$ from which one closed curve is lifted. The holonomies around the two edges of the over-strand differ by the factor $\omega_j^{\pm1}$, by \eqref{eq.omega}, and so do those of the under-strand; this is why each form occurs with two values of its variable, according to the starting point of the turn. They coincide if and only if $\omega_j=1$.

\begin{remark}[the Alexander point]\label{rem.alexpoint}
Let $t$ be a unit of $\kk$ and $\vth_f=t^{w(f)}$, where $w$ is the winding function of Lemma~\ref{lem.winding}. Then $\mhol_e=t$ for every edge $e$, by \eqref{eq.wh}, and $\omega_j=1$ for every crossing: the local system $\CL_\vth$ extends to the complement of $L$, where it is the local system which sends the meridians of all the components of $L$ to the same unit, the one which defines the Alexander polynomial. All the turns $\mer_j(i,i+4)$ with $i$ odd are then equal to one closed chain, which reads $\pi_j$ at $x_j=t^{-\ve_j}$, and those with $i$ even are equal to one closed chain, which reads $\tpi_j$ at $y_j=t^{\ve_j}$, up to units. By Theorem~\ref{thm.flows} the determinants of these forms in a basis of $\FF_{D,\hate}$ are $\pm P_{D,\hate}(t^{-\ve})$ and $\pm\tP_{D,\hate}(t^{\ve})$, that is, up to sign, the Alexander determinants $\det A_{D,\hate}(t^{-1})$ and $\det B_{D,\hate}(t)$ of Section~\ref{sec.alexander}.
\end{remark}

\subsection{Interleaved chains and the coupled crossing forms}\label{ssec.interleaved}
The coupled crossing form $\Psi_j$ of Section~\ref{ssec.mixed} is a bilinear form on the lattice of flows. In the dictionary of Lemma~\ref{lem.iota} it therefore corresponds to a chain in the tensor product of two copies of the complex of Definition~\ref{def.chains}, and in this subsection we determine this chain. It is not the product of two turns of the lifted meridian, which gives a form of rank one. It is a sum of two such products, in which the turn of the first factor depends on the half-turn of the second one: an \emph{interleaved chain}.

For $v\in C^-_\vth$ and $v'\in C^+_\vth$ we denote by $(\iota\otimes\iota')(v\otimes v')$ the bilinear form $(z,z')\mapsto\iota(v)(z)\,\iota'(v')(z')$ on $\FF_D$, with values in $\kk$, and we extend this notation linearly to the tensor product $C^-_\vth\otimes_{\kk}C^+_\vth$. By Proposition~\ref{prop.readings}, the product of two turns of the lifted meridian of $\alpha_j$, the first one starting at an odd place and the second one at an even place, gives a form $\pi_j^{(\xs)}\otimes\tpi_j^{(\ys)}$ of rank one, up to a unit. The coupled crossing forms come from a different combination.

\begin{definition}[interleaved chains]\label{def.interleaved}
Let $\vth$ be a multiplicative potential and $j$ a crossing. The \emph{squares} of the crossing $j$ are the chains $\mer_j(i,i+4)\otimes\mer_j(i',i'+4)$ with $i$ odd and $i'$ even, and its two \emph{interleaved chains} are
\begin{align*}
   \CZ^\out_j&=\mer_j(3,7)\otimes\mer_j(0,2)+\mer_j(1,5)\otimes\mer_j(2,4),\\
   \CZ^\inn_j&=\mer_j(1,5)\otimes\mer_j(0,2)+\mer_j(-1,3)\otimes\mer_j(2,4),
\end{align*}
elements of $C^-_\vth\otimes_{\kk}C^+_\vth$.
\end{definition}

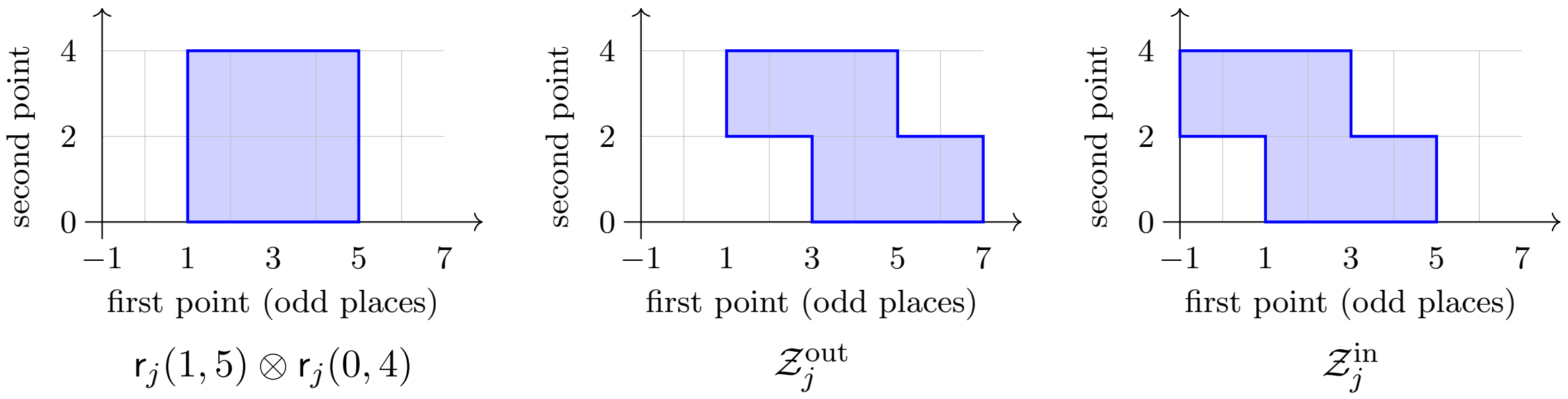
\begin{figure}[ht]
\centering
\begin{tikzpicture}[scale=0.42,every node/.style={font=\scriptsize}]
\def\grid{%
  \draw[black!25,very thin] (-1,0) grid[xstep=2,ystep=2] (7,4);
  \draw[->] (-1.4,0) -- (7.9,0);
  \draw[->] (-1,-0.4) -- (-1,5.0);
  \foreach \x in {-1,1,3,5,7} \node[below] at (\x,-0.25) {$\x$};
  \foreach \y in {0,2,4} \node[left] at (-1.25,\y) {$\y$};
  \node[below] at (3,-1.2) {first point (odd places)};
  \node[rotate=90] at (-2.9,2) {second point};}
\begin{scope}
  \fill[blue!18] (1,0) rectangle (5,4);
  \grid
  \draw[thick,blue] (1,0) rectangle (5,4);
  \node at (3,-3.4) {\small $\mer_j(1,5)\otimes\mer_j(0,4)$};
\end{scope}
\begin{scope}[xshift=12.6cm]
  \fill[blue!18] (3,0) rectangle (7,2); \fill[blue!18] (1,2) rectangle (5,4);
  \grid
  \draw[thick,blue] (3,0) -- (7,0) -- (7,2) -- (5,2) -- (5,4) -- (1,4) -- (1,2) -- (3,2) -- cycle;
  \node at (3,-3.4) {\small $\CZ^\out_j$};
\end{scope}
\begin{scope}[xshift=25.2cm]
  \fill[blue!18] (1,0) rectangle (5,2); \fill[blue!18] (-1,2) rectangle (3,4);
  \grid
  \draw[thick,blue] (1,0) -- (5,0) -- (5,2) -- (3,2) -- (3,4) -- (-1,4) -- (-1,2) -- (1,2) -- cycle;
  \node at (3,-3.4) {\small $\CZ^\inn_j$};
\end{scope}
\end{tikzpicture}
\caption{The plane of the pairs of places of the lifted meridian of a co-core: the first coordinate is a place below the over-strand, the second one a place above the under-strand. A rectangle $[i,l]\times[i',l']$ stands for the chain $\mer_j(i,l)\otimes\mer_j(i',l')$. Shaded: a square and the two interleaved chains of Definition~\ref{def.interleaved}.}
\label{fig.staircase}
\end{figure}

Think of two points which move on the lifted meridian, the first one between places below the over-strand and the second one between places above the under-strand (Figure~\ref{fig.staircase}). In a square each of them makes one turn, independently of the other one. In an interleaved chain the second point also makes the turn $\mer_j(0,4)=\mer_j(0,2)+\mer_j(2,4)$, but the turn of the first point depends on where the second point is. In $\CZ^\out_j$: while the second point makes the half-turn from the place $0$ to the place $2$, under the edge $\hat e^+_j$ of the over-strand, the first point makes the turn which starts at the place $3$, below the other edge $\hat e^-_j$; while the second point makes the half-turn from $2$ to $4$, under $\hat e^-_j$, the first point makes the turn which starts at the place $1$, below $\hat e^+_j$. In $\CZ^\inn_j$ the two turns of the first point start at the places $1$ and $-1$: it is the same rule, shifted by a half-turn. In the plane of the pairs of places, a square is a fundamental domain for the translations of the two points by one turn, and the interleaved chains are two other fundamental domains, in the shape of a staircase.

The forms which these chains define involve the holonomies around the four edges of the crossing, through the units
\begin{equation}\label{eq.xyvalues}
   \xs^\out_j=\mhol_{o'_j}^{-\ve_j},\quad \ys^\out_j=\mhol_{u'_j}^{\ve_j},\qquad \xs^\inn_j=\mhol_{o_j}^{-\ve_j},\quad \ys^\inn_j=\mhol_{u_j}^{\ve_j},\qquad \sv_j=\Bigl(\frac{\mhol_{u'_j}}{\mhol_{o_j}}\Bigr)^{\ve_j}=\Bigl(\frac{\mhol_{u_j}}{\mhol_{o'_j}}\Bigr)^{\ve_j} :
\end{equation}
the first two are the holonomies around the outgoing edges of the two strands, with the exponents $-\ve_j$ and $\ve_j$, and the next two those around the incoming edges; the two expressions of $\sv_j$ agree by \eqref{eq.mvertex}. By \eqref{eq.omega},
\begin{equation}\label{eq.inout}
   \xs^\inn_j=\omega_j^{\ve_j}\,\xs^\out_j,\qquad \ys^\inn_j=\omega_j^{\ve_j}\,\ys^\out_j,\qquad \sv_j=\xs^\inn_j\,\ys^\out_j=\xs^\out_j\,\ys^\inn_j .
\end{equation}

\begin{proposition}[the interleaved chains read the coupled crossing forms]\label{prop.interleaved}
Let $\vth$ be a multiplicative potential and $j$ a crossing; we write $\omega=\omega_j$, $\mer=\mer_j$ and $\vth_1=\vth_{f^1_j}$.
\begin{enumerate}\renewcommand{\theenumi}{\roman{enumi}}
\item The interleaved chains differ from a square by a multiple of $\omega-1$:
\begin{align*}
   \CZ^\out_j&=\mer(1,5)\otimes\mer(0,4)+(\omega-1)\,\mer(1,3)\otimes\mer(0,2),\\
   \CZ^\inn_j&=\mer(1,5)\otimes\mer(0,4)-(\omega-1)\,\mer(-1,1)\otimes\mer(2,4).
\end{align*}
In particular, if $\omega_j=1$, then all the squares and the two interleaved chains of the crossing $j$ are one and the same chain.
\item The forms of the interleaved chains are, up to units, the forms of Definition~\ref{def.coupledforms}:
\[
   (\iota\otimes\iota')\bigl(\CZ^\out_j\bigr)=\us^\out_j\,\Psi_j\bigl[\sv_j;\xs^\out_j,\ys^\out_j\bigr],\qquad (\iota\otimes\iota')\bigl(\CZ^\inn_j\bigr)=\us^\inn_j\,\Psi_j\bigl[\sv_j;\xs^\inn_j,\ys^\inn_j\bigr],
\]
with the units $\us^\out_j=\mathsf c_j(4)\,\mathsf c_j(2+\ve_j)$ and $\us^\inn_j=\mathsf c_j(4-\ve_j)$ of Proposition~\ref{prop.readings}~(ii): $\us^\out_j=\omega\,\vth_1/\mhol_{u'_j}$ and $\us^\inn_j=\vth_1/\mhol_{u'_j}$ if $\ve_j=1$, and $\us^\out_j=\us^\inn_j=\omega\,\vth_1$ if $\ve_j=-1$. We also write $\us(\CZ^\out_j)=\us^\out_j$ and $\us(\CZ^\inn_j)=\us^\inn_j$.
\end{enumerate}
\end{proposition}
\begin{proof}
(i) By \eqref{eq.turndiff}, $\mer(3,7)=\mer(1,5)+(\omega-1)\,\mer(1,3)$ and $\mer(-1,3)=\mer(1,5)-(\omega-1)\,\mer(-1,1)$; we substitute this in Definition~\ref{def.interleaved} and use $\mer(0,2)+\mer(2,4)=\mer(0,4)$. If $\omega=1$, the chain $\mer(i,i+4)$ only depends on the parity of $i$, by \eqref{eq.turndiff}.

(ii) Write $\ve=\ve_j$ and $(o,u,o',u')$ for the four edges of the crossing. By Proposition~\ref{prop.readings} and since $\mer(-1,3)=\omega^{-1}\mer(3,7)$, the two forms are
\begin{align*}
   (\iota\otimes\iota')\bigl(\CZ^\out_j\bigr)(z,z')&=\ve\,\omega\,\vth_1\Bigl(\bigl(\mhol_{\hat e^+}^{-\ve}z_{u'}-z_u\bigr)\,z'_{\hat e^+}-\frac1{\mhol_{u'}}\bigl(\mhol_{\hat e^-}^{-\ve}z_{u'}-z_u\bigr)\,z'_{\hat e^-}\Bigr),\\
   (\iota\otimes\iota')\bigl(\CZ^\inn_j\bigr)(z,z')&=\ve\,\vth_1\Bigl(\omega\,\bigl(\mhol_{\hat e^-}^{-\ve}z_{u'}-z_u\bigr)\,z'_{\hat e^+}-\frac1{\mhol_{u'}}\bigl(\mhol_{\hat e^+}^{-\ve}z_{u'}-z_u\bigr)\,z'_{\hat e^-}\Bigr).
\end{align*}
Suppose that $\ve=1$, so that $\hat e^+=o$ and $\hat e^-=o'$. The first form is
\[
   \frac{\omega\,\vth_1}{\mhol_{u'}}\Bigl(\frac{\mhol_{u'}}{\mhol_o}\,z_{u'}z'_o-\frac1{\mhol_{o'}}\,z_{u'}z'_{o'}-\mhol_{u'}\,z_uz'_o+z_uz'_{o'}\Bigr)=\us^\out_j\,\Psi_j\bigl[\sv_j;\xs^\out_j,\ys^\out_j\bigr](z,z')
\]
by \eqref{eq.Psigeneral} and \eqref{eq.xyvalues}, and the second one is
\[
   \frac{\vth_1}{\mhol_{u'}}\Bigl(\frac{\omega\,\mhol_{u'}}{\mhol_{o'}}\,z_{u'}z'_o-\frac1{\mhol_o}\,z_{u'}z'_{o'}-\omega\,\mhol_{u'}\,z_uz'_o+z_uz'_{o'}\Bigr)=\us^\inn_j\,\Psi_j\bigl[\sv_j;\xs^\inn_j,\ys^\inn_j\bigr](z,z'),
\]
because $\omega\,\mhol_{u'}=\mhol_u$ and $\omega\,\mhol_{u'}/\mhol_{o'}=\mhol_{u'}/\mhol_o$ by \eqref{eq.omega}. Suppose that $\ve=-1$, so that $\hat e^+=o'$ and $\hat e^-=o$. The first form is
\[
   \omega\,\vth_1\Bigl(\frac{\mhol_o}{\mhol_{u'}}\,z_{u'}z'_o-\mhol_{o'}\,z_{u'}z'_{o'}-\frac1{\mhol_{u'}}\,z_uz'_o+z_uz'_{o'}\Bigr)=\us^\out_j\,\Psi_j\bigl[\sv_j;\xs^\out_j,\ys^\out_j\bigr](z,z'),
\]
and the second one is
\[
   \omega\,\vth_1\Bigl(\frac{\mhol_{o'}}{\omega\,\mhol_{u'}}\,z_{u'}z'_o-\mhol_o\,z_{u'}z'_{o'}-\frac1{\omega\,\mhol_{u'}}\,z_uz'_o+z_uz'_{o'}\Bigr)=\us^\inn_j\,\Psi_j\bigl[\sv_j;\xs^\inn_j,\ys^\inn_j\bigr](z,z'),
\]
because $\mhol_{o'}/\omega=\mhol_o$ and $\omega\,\mhol_{u'}=\mhol_u$.
\end{proof}

Recall that the form $\Psi_j[\sv;\xs,\ys]$ has rank one when $\sv=\xs\ys$. By \eqref{eq.inout}, $\sv_j-\xs^\out_j\ys^\out_j=(\omega_j^{\ve_j}-1)\,\xs^\out_j\ys^\out_j$: the forms of the interleaved chains have rank two exactly when the holonomy $\omega_j$ around the co-core is not trivial.

\emph{The pairing.} Let $\hate$ be an edge of $D$, and let $\CF_{\hate}$ be the set of the $c$ faces which are not adjacent to $\hate$. For an element $\CZ=\sum_{f,g\in\CF}\CZ_{fg}\,\tpth_f\otimes\tpth_g$ of $C_\vth\otimes_{\kk}C_\vth$ we denote by $\CZ|_{\hate}=(\CZ_{fg})_{f,g\in\CF_{\hate}}$ the $c\times c$ matrix of its coefficients on the faces of $\CF_{\hate}$, and for $c$ such elements we put
\[
   \lb\CZ_1,\dots,\CZ_c\rb_{\hate}=\bigl\lb\,\CZ_1|_{\hate},\dots,\CZ_c|_{\hate}\,\bigr\rb\ \in\kk,
\]
the mixed discriminant of Definition~\ref{def.mixed}.

\begin{lemma}\label{lem.chainpairing}
If $\CZ_1,\dots,\CZ_c\in C^-_\vth\otimes_{\kk}C^+_\vth$, then
\[
   \bigl\lb\,(\iota\otimes\iota')(\CZ_1),\dots,(\iota\otimes\iota')(\CZ_c)\,\bigr\rb_{\hate}=\Bigl(\prod_{f\in\CF_{\hate}}\vth_f\Bigr)\,\lb\CZ_1,\dots,\CZ_c\rb_{\hate} .
\]
\end{lemma}
\begin{proof}
By Lemma~\ref{lem.faces}~(iii), the boundaries $\partial f$, $f\in\CF_{\hate}$, form a basis of $\FF_{D,\hate}$, and the unit vector $e_f$ is a potential of $\partial f$. So $\iota(v)(\partial f)=-v_f$ for $v\in C^-_\vth$ and $\iota'(v')(\partial g)=-\vth_g\,v'_g$ for $v'\in C^+_\vth$, by Definition~\ref{def.iota}. Hence the matrix of the form $(\iota\otimes\iota')(\CZ)$ in this basis is $\bigl(\CZ_{fg}\,\vth_g\bigr)_{f,g\in\CF_{\hate}}$, the product of $\CZ|_{\hate}$ and of the diagonal matrix of the units $\vth_g$, and the claim follows from \eqref{eq.mixedproperties}.
\end{proof}

\begin{corollary}[all the local systems of rank one]\label{cor.anypotential}
Let $\vth$ be a multiplicative potential, let $\hate$ be an edge of $D$, and for every crossing $j$ let $\CZ_j$ be one of the two interleaved chains $\CZ^\out_j$, $\CZ^\inn_j$, with its unit $\us(\CZ_j)$ of Proposition~\ref{prop.interleaved}, and let $(\xs_j,\ys_j)$ be the pair $(\xs^\out_j,\ys^\out_j)$ or $(\xs^\inn_j,\ys^\inn_j)$ of \eqref{eq.xyvalues}, accordingly. Then
\[
   \Bigl(\prod_{f\in\CF_{\hate}}\vth_f\Bigr)\,\lb\CZ_1,\dots,\CZ_c\rb_{\hate}=\Bigl(\prod_{j=1}^c\us(\CZ_j)\Bigr)\,(-1)^{\mu-1}\,\BE'\bigl[P_{D,\hate}(x)\,\tP_{D,\hate}(y)\bigr],
\]
where $\BE'$ is the linear map of Theorem~\ref{thm.mixed}~\textup{(i)} with $x_j\mapsto\xs_j$, $y_j\mapsto\ys_j$ and $x_jy_j\mapsto\sv_j$.
\end{corollary}
\begin{proof}
Write $\Psi'_j=\Psi_j[\sv_j;\xs_j,\ys_j]$. By Lemma~\ref{lem.chainpairing}, by Proposition~\ref{prop.interleaved}~(ii), by the multilinearity of the mixed discriminant and by Theorem~\ref{thm.mixed}~(i), in this order,
\begin{align*}
   \Bigl(\prod_{f\in\CF_{\hate}}\vth_f\Bigr)\,\lb\CZ_1,\dots,\CZ_c\rb_{\hate}
   &=\bigl\lb\,(\iota\otimes\iota')(\CZ_1),\dots,(\iota\otimes\iota')(\CZ_c)\,\bigr\rb_{\hate}\\
   &=\bigl\lb\,\us(\CZ_1)\,\Psi'_1,\dots,\us(\CZ_c)\,\Psi'_c\,\bigr\rb_{\hate}\\
   &=\Bigl(\prod_{j=1}^c\us(\CZ_j)\Bigr)\,\lb\Psi'_1,\dots,\Psi'_c\rb_{\hate}\\
   &=\Bigl(\prod_{j=1}^c\us(\CZ_j)\Bigr)\,(-1)^{\mu-1}\,\BE'\bigl[P_{D,\hate}(x)\,\tP_{D,\hate}(y)\bigr].\qedhere
\end{align*}
\end{proof}

\begin{remark}[evaluations and the coupling]\label{rem.dictionary}
Corollary~\ref{cor.anypotential} is a dictionary between the local systems of rank one on $X_G$ and the linear maps $\BE'$ which can be applied to the product $P_{D,\hate}\tP_{D,\hate}$ of the two edge polynomials. Two families of local systems stand out.

(i) \emph{Evaluations.} The map $\BE'$ is an evaluation of the crossing variables, that is, $\sv_j=\xs_j\ys_j$ for all $j$, if and only if $\omega_j=1$ for all $j$, by \eqref{eq.inout}. By \eqref{eq.omega} this means that the holonomy $\mhol_e$ is the same for the two edges of each strand at every crossing, hence that it is constant along every component of $L$: these are the local systems which extend to the complement of the link. Then the chains $\CZ_j$ are squares, their forms have rank one, and the mixed discriminant is the product of two determinants (Lemma~\ref{lem.cauchybinet}), namely $(-1)^{\mu-1}P_{D,\hate}(\xs)\,\tP_{D,\hate}(\ys)$ at $\xs_j=\mhol_{o_j}^{-\ve_j}$ and $\ys_j=\mhol_{u_j}^{\ve_j}$.

(ii) \emph{The rule $x_jy_j\mapsto1$.} We have $\sv_j=1$ if and only if $\mhol_{u'_j}=\mhol_{o_j}$, and then $\mhol_{u_j}=\mhol_{o'_j}$: the holonomy is the same for the two edges $o_j,u'_j$ of the Seifert circle $A_j$, and for the two edges $u_j,o'_j$ of the Seifert circle $B_j$ (Figure~\ref{fig.crossing}). By \eqref{eq.mcrossing}, for both signs $\sv_j=\vth_{f^\downarrow_j}/\vth_{f^\uparrow_j}$: this is the holonomy of the loop $\pth_{f^\downarrow_j}\pth_{f^\uparrow_j}^{-1}$, which passes under the two strands of the crossing and comes back over them. This loop goes once around the band $b_j$ of the canonical surface; it is a meridian of the core $\kappa_j$ of this band (Section~\ref{ssec.statesurface}). So the rule $x_jy_j\mapsto1$ of the coupling $\BE_D$ says that the holonomy around the band of the crossing $j$ is trivial, whereas the rule $x_jy_j\mapsto x_jy_j$ of an evaluation says that the holonomy around its co-core is trivial. If $\sv_j=1$ for all $j$, the holonomy $\mhol_e$ is constant along every Seifert circle of $D$, and $\vth$ is constant on every region of the complement of the Seifert circles (Section~\ref{ssec.regions}): the local system is determined by one unit for each Seifert circle, as if it lived on the complement of the Seifert circles instead of the complement of the link.

(iii) Both conditions hold for all $j$ if and only if $\mhol_e$ is the same unit $t$ for all the edges of $D$: this is the Alexander point of Remark~\ref{rem.alexpoint}.
\end{remark}

\subsection{The two-colored local system and the Links--Gould polynomial}\label{ssec.twocolored}
We now choose the local system for which the forms of Corollary~\ref{cor.anypotential} are the coupled crossing forms $\Psi_j=\Psi_j[1;t_0^{\ve_j},t_1^{\ve_j}]$ of Definition~\ref{def.coupledforms}. Recall that every crossing $j$ joins two different Seifert circles of $D$: the circle $A_j$, which contains the edges $o_j$ and $u'_j$, and the circle $B_j$, which contains $u_j$ and $o'_j$ (Figure~\ref{fig.crossing}). A \emph{two-coloring} of the Seifert circles is a map $\mathrm{col}$ from the set of the Seifert circles of $D$ to $\{0,1\}$ such that $\mathrm{col}(A_j)\ne\mathrm{col}(B_j)$ for every crossing $j$; the color of an edge is the color of its Seifert circle.

\begin{lemma}[the two-colored potential]\label{lem.twocolored}
\begin{enumerate}\renewcommand{\theenumi}{\roman{enumi}}
\item There are exactly two two-colorings of the Seifert circles of $D$: the parity of the level $h(C)$ of Definition~\ref{def.level}, and the opposite coloring.
\item Let $\mathrm{col}$ be a two-coloring. There is a multiplicative potential $\vth$ with values in $\BZ[t_0^{\pm1},t_1^{\pm1}]$, unique up to a common factor, such that
\[
   \mhol_e=t_1\quad\text{if }\mathrm{col}(e)=1,\qquad\qquad \mhol_e=t_0^{-1}\quad\text{if }\mathrm{col}(e)=0,
\]
and the common factor can be chosen so that its values are monomials in $t_0$ and $t_1$. With this choice we call it the \emph{two-colored potential} of $\mathrm{col}$.
\item For the two-colored potential, $\sv_j=1$ for every crossing $j$, and
\begin{align*}
   &\text{if }\mathrm{col}(B_j)=0:&\xs^\out_j&=t_0^{\ve_j},&\ys^\out_j&=t_1^{\ve_j},&\omega_j&=(t_0t_1)^{-1};\\
   &\text{if }\mathrm{col}(A_j)=0:&\xs^\inn_j&=t_0^{\ve_j},&\ys^\inn_j&=t_1^{\ve_j},&\omega_j&=t_0t_1 .
\end{align*}
We say that the crossing $j$ is \emph{of type out} in the first case, in which its outgoing over-edge $o'_j$ has the color $0$, and \emph{of type in} in the second case, in which its incoming over-edge $o_j$ has the color $0$.
\item All the crossings of $D$ have the same type if and only if the diagram $D$ is alternating.
\end{enumerate}
\end{lemma}

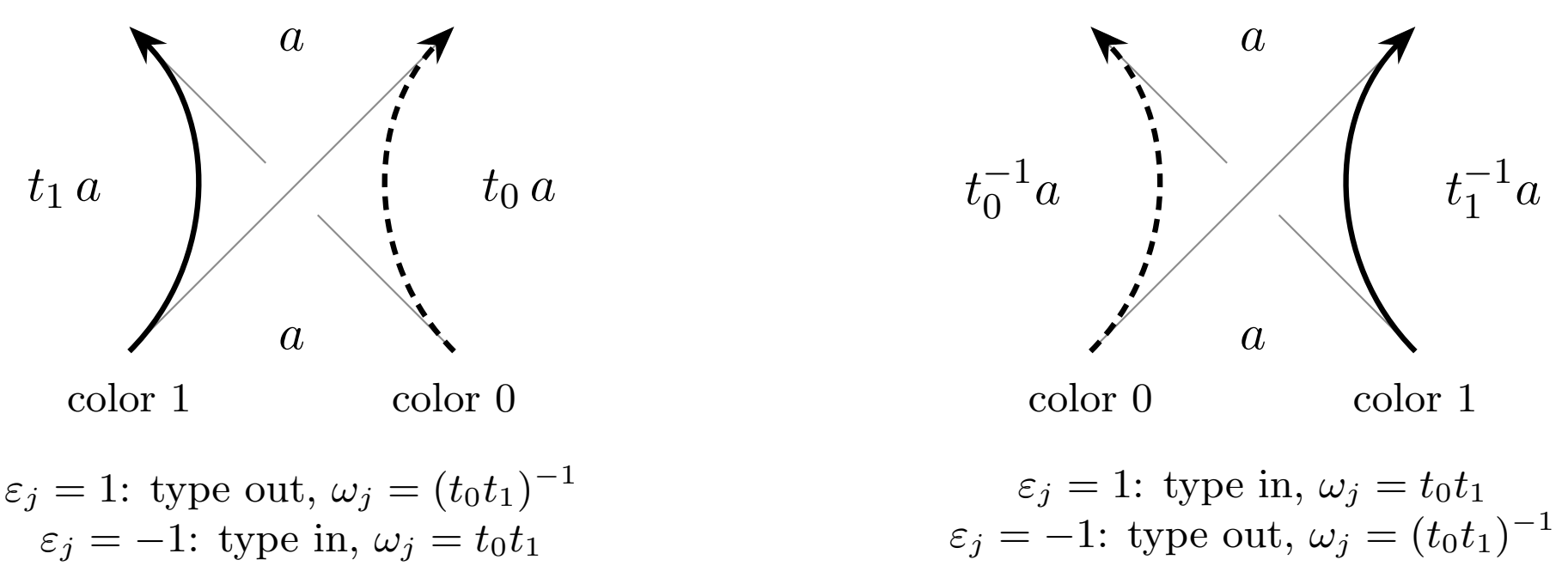
\begin{figure}[ht]
\centering
\begin{tikzpicture}[>=Stealth,scale=1.0,every node/.style={font=\small}]
\def\crossingpic{%
  \draw[black!45] (1.25,-1.25) -- (0.2,-0.2); \draw[black!45,->] (-0.2,0.2) -- (-1.25,1.25);
  \draw[black!45,->] (-1.25,-1.25) -- (1.25,1.25);}
\begin{scope}
  \crossingpic
  \draw[very thick,->] (-1.25,-1.25) to[out=45,in=-45] (-1.25,1.25);
  \draw[very thick,densely dashed,->] (1.25,-1.25) to[out=135,in=-135] (1.25,1.25);
  \node at (0,1.15) {$a$}; \node at (0,-1.15) {$a$}; \node at (-1.75,0) {$t_1\,a$}; \node at (1.75,0) {$t_0\,a$};
  \node[font=\scriptsize] at (-1.25,-1.6) {color $1$}; \node[font=\scriptsize] at (1.25,-1.6) {color $0$};
  \node[font=\scriptsize,align=center] at (0,-2.5) {$\ve_j=1$: type out, $\omega_j=(t_0t_1)^{-1}$\\ $\ve_j=-1$: type in, $\omega_j=t_0t_1$};
\end{scope}
\begin{scope}[xshift=7.4cm]
  \crossingpic
  \draw[very thick,densely dashed,->] (-1.25,-1.25) to[out=45,in=-45] (-1.25,1.25);
  \draw[very thick,->] (1.25,-1.25) to[out=135,in=-135] (1.25,1.25);
  \node at (0,1.15) {$a$}; \node at (0,-1.15) {$a$}; \node at (-1.85,0) {$t_0^{-1}a$}; \node at (1.85,0) {$t_1^{-1}a$};
  \node[font=\scriptsize] at (-1.25,-1.6) {color $0$}; \node[font=\scriptsize] at (1.25,-1.6) {color $1$};
  \node[font=\scriptsize,align=center] at (0,-2.5) {$\ve_j=1$: type in, $\omega_j=t_0t_1$\\ $\ve_j=-1$: type out, $\omega_j=(t_0t_1)^{-1}$};
\end{scope}
\end{tikzpicture}
\caption{The values of the two-colored potential on the four faces at a crossing $j$, drawn with both strands pointing upward (in gray, for $\ve_j=1$; for $\ve_j=-1$ the other strand is the over-strand). The two arcs of the oriented smoothing belong to the two Seifert circles of the crossing, of color $1$ (solid) and $0$ (dashed); $a$ is the common value on the faces $f^\uparrow_j$ and $f^\downarrow_j$. If $\ve_j=1$ the circle on the left is $A_j$, and if $\ve_j=-1$ it is $B_j$ (Figure~\ref{fig.crossing}).}
\label{fig.twocolored}
\end{figure}

\begin{proof}
(i) By Lemma~\ref{lem.level}, $h(B_j)=h(A_j)+\ve_j$, so the parity of $h$ is a two-coloring. The Seifert graph $\Gamma_D$ of Section~\ref{ssec.smoothings}, whose vertices are the Seifert circles and whose edges are the crossings, is connected because $D$ is, and a two-coloring is determined by the color of one vertex of this graph.

(ii) Let $\1_{(1)}\in\BZ^E$ be the indicator vector of the set of the edges of color $1$, and let $\1_{(0)}=\1-\1_{(1)}$ be that of the set of the edges of color $0$. At every crossing one incoming and one outgoing edge lie on $A_j$, and the other two on $B_j$; so $\1_{(1)}$ and $\1_{(0)}$ are flows: they are the sums of the Seifert circles of each color. Let $k^{(1)}$ and $k^{(0)}$ be potentials of these two flows (Lemma~\ref{lem.faces}~(ii)). By \eqref{eq.potential}, the multiplicative potential $\vth_f=t_1^{\,k^{(1)}_f}\,t_0^{-k^{(0)}_f}$ has the holonomy $\mhol_e=t_1$ around the edges of color $1$ and $\mhol_e=t_0^{-1}$ around the edges of color $0$. If two potentials have the same holonomies, their quotient takes the same value on the two faces adjacent to every edge, hence on all the faces.

(iii) We have $\mhol_{o_j}=\mhol_{u'_j}$ and $\mhol_{u_j}=\mhol_{o'_j}$, hence $\sv_j=1$ by \eqref{eq.xyvalues}. If $\mathrm{col}(B_j)=0$, then $\mhol_{o'_j}=\mhol_{u_j}=t_0^{-1}$ and $\mhol_{o_j}=\mhol_{u'_j}=t_1$, so that $\xs^\out_j=t_0^{\ve_j}$, $\ys^\out_j=t_1^{\ve_j}$ and $\omega_j=\mhol_{o'_j}/\mhol_{o_j}=(t_0t_1)^{-1}$. If $\mathrm{col}(A_j)=0$, then $\mhol_{o_j}=\mhol_{u'_j}=t_0^{-1}$ and $\mhol_{u_j}=\mhol_{o'_j}=t_1$, so that $\xs^\inn_j=t_0^{\ve_j}$, $\ys^\inn_j=t_1^{\ve_j}$ and $\omega_j=t_0t_1$.

(iv) By (i) we may assume that the color is the parity of $h$; then the type of $j$ is given by the parity of $h(u_j)=h(B_j)$. Let $e$ be an edge from a crossing $j$ to a crossing $k$. If $e=o'_j=u_k$, then $h(u_k)=h(o'_j)=h(u_j)$ by \eqref{eq.levels}. If $e=u'_j=o_k$, then $h(u_k)=h(o_k)+\ve_k=h(u'_j)+\ve_k=h(u_j)-\ve_j+\ve_k$. In both cases $j$ and $k$ have the same type. If $e=o'_j=o_k$, then $h(u_k)=h(o_k)+\ve_k=h(u_j)+\ve_k$, and if $e=u'_j=u_k$, then $h(u_k)=h(u'_j)=h(u_j)-\ve_j$; in these two cases $j$ and $k$ have different types. Since $D$ is connected, all the crossings have the same type if and only if every edge is an edge of the over-strand at exactly one of its ends, that is, if and only if $D$ is alternating.
\end{proof}

Figure~\ref{fig.twocolored} shows the two-colored potential near a crossing. It only changes when one crosses a Seifert circle: from the right to the left of a circle of color $1$ it is multiplied by $t_1$, and from the right to the left of a circle of color $0$ by $t_0^{-1}$.

\begin{theorem}[the Links--Gould polynomial in the complement of the co-cores]\label{thm.cocores}
Let $D$ be a diagram of a link $L$ and $\hate$ an edge of $D$. Let $\mathrm{col}$ be a two-coloring of the Seifert circles of $D$ and $\vth$ its two-colored potential, and for every crossing $j$ let
\[
   \CZ_j=\CZ^\out_j\quad\text{if $j$ is of type out,}\qquad\qquad \CZ_j=\CZ^\inn_j\quad\text{if $j$ is of type in.}
\]
\begin{enumerate}\renewcommand{\theenumi}{\roman{enumi}}
\item For every crossing $j$, the form of the interleaved lift $\CZ_j$ of the meridian of the co-core $\alpha_j$ is the coupled crossing form, up to a unit: $(\iota\otimes\iota')(\CZ_j)=\us(\CZ_j)\,\Psi_j$, where $\us(\CZ_j)$ is the unit $\us^\out_j$ or $\us^\inn_j$ of Proposition~\ref{prop.interleaved}, a monomial in $t_0$ and $t_1$.
\item The Links--Gould polynomial is the mixed discriminant of the interleaved chains of the $c$ co-cores, up to a unit:
\[
   \lb\CZ_1,\dots,\CZ_c\rb_{\hate}=\Bigl(\prod_{j=1}^c\us(\CZ_j)\Bigr)\Bigl(\prod_{f\in\CF_{\hate}}\vth_f\Bigr)^{-1}\LG(L;t_0,t_1),
\]
and the factor in front of $\LG$ is a monomial $t_0^{a}t_1^{b}$.
\item The holonomy of the local system $\CL_\vth$ around the co-core $\alpha_j$ is $(t_0t_1)^{-1}$ if $j$ is of type out and $t_0t_1$ if $j$ is of type in. After the specialization $t_1=t_0^{-1}$ it extends to the complement of $L$, where it is the local system of the Alexander polynomial, with the holonomy $t_0^{-1}$ around every edge; the chains $\CZ_j$ become squares, and $\lb\CZ_1,\dots,\CZ_c\rb_{\hate}$ becomes a product of two determinants.
\end{enumerate}
\end{theorem}
\begin{proof}
(i) follows from Proposition~\ref{prop.interleaved}~(ii) and Lemma~\ref{lem.twocolored}~(iii), since $\Psi_j=\Psi_j[1;t_0^{\ve_j},t_1^{\ve_j}]$. (ii) By Lemma~\ref{lem.chainpairing}, by (i), by the multilinearity of the mixed discriminant and by Theorem~\ref{thm.mixed}~(ii), in this order,
\begin{align*}
   \Bigl(\prod_{f\in\CF_{\hate}}\vth_f\Bigr)\,\lb\CZ_1,\dots,\CZ_c\rb_{\hate}
   &=\bigl\lb\,\us(\CZ_1)\,\Psi_1,\dots,\us(\CZ_c)\,\Psi_c\,\bigr\rb_{\hate}\\
   &=\Bigl(\prod_{j=1}^c\us(\CZ_j)\Bigr)\,\lb\Psi_1,\dots,\Psi_c\rb_{\hate}
   =\Bigl(\prod_{j=1}^c\us(\CZ_j)\Bigr)\,\LG(L;t_0,t_1).
\end{align*}
The values $\vth_f$ of the two-colored potential are monomials by Lemma~\ref{lem.twocolored}~(ii), and so are the units $\us(\CZ_j)$, which are products of values of the potential and of their inverses (Proposition~\ref{prop.interleaved}~(ii)). (iii) The holonomies are those of Lemma~\ref{lem.twocolored}~(iii). At $t_1=t_0^{-1}$ they are trivial, $\mhol_e=t_0^{-1}$ for every edge, and the rest is Proposition~\ref{prop.interleaved}~(i), Remark~\ref{rem.dictionary}~(i) and Remark~\ref{rem.alexpoint}.
\end{proof}

So the two flavors of the determinant formula are the two sides of $\Sigma$ from which the meridians of the co-cores are lifted, the two variables are, up to an inversion, the holonomies around the Seifert circles of the two colors, namely $t_1$ and $t_0^{-1}$, and the coupling is the interleaving of the two lifts, which matters exactly because the local system does not extend across the co-cores: its holonomy around each of them is $(t_0t_1)^{\pm1}$. At the free fermion point $t_0t_1=1$ of Section~\ref{sec.formula} this holonomy is trivial, the interleaving disappears, and one is left with the square of the Alexander polynomial.

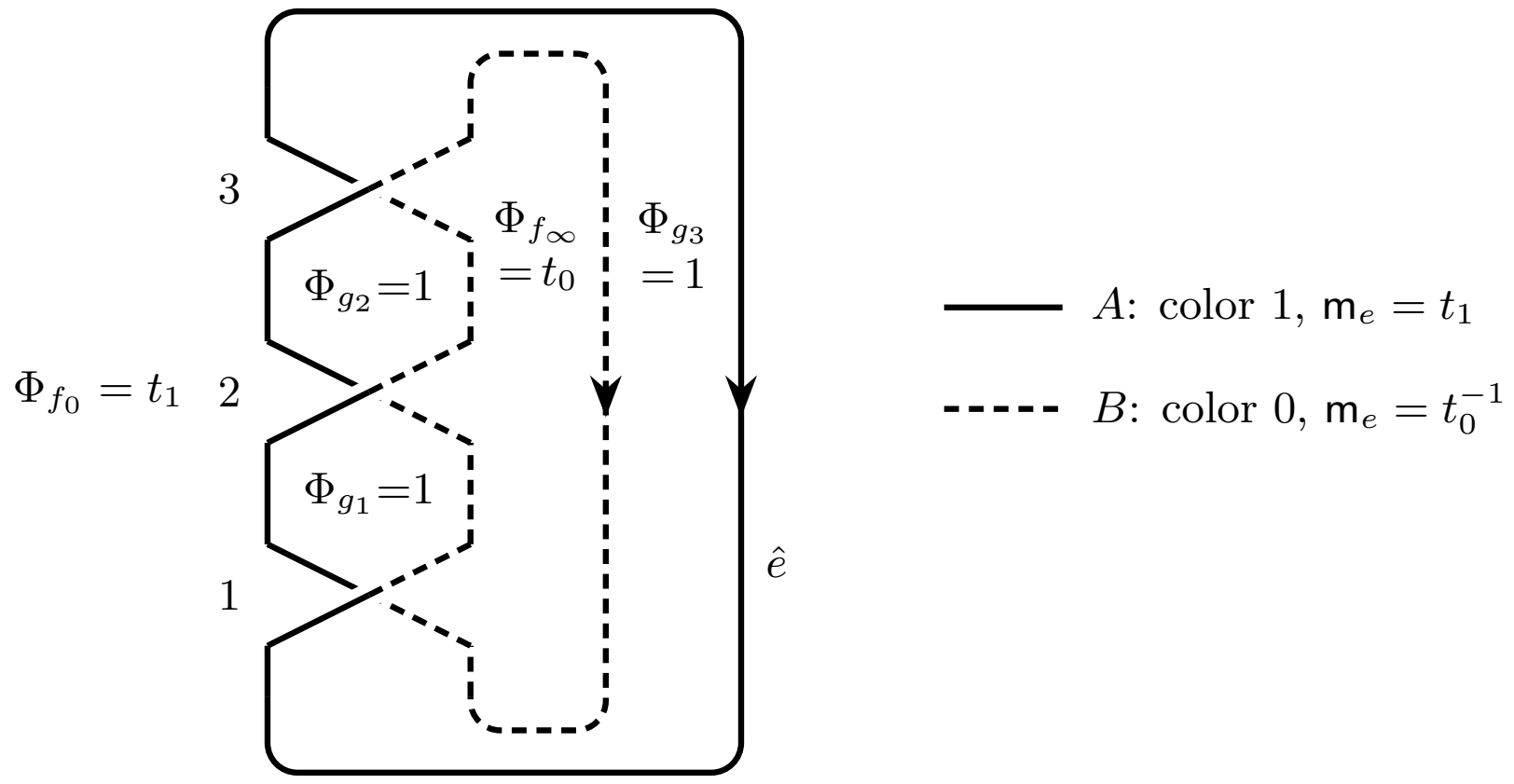
\begin{figure}[ht]
\centering
\begin{tikzpicture}[>=Stealth,scale=0.6,every node/.style={font=\small},cone/.style={very thick},czero/.style={very thick,densely dashed}]
\draw[cone] (0,0) -- (0,0.6); \draw[cone] (0,1.8) -- (0,3.0); \draw[cone] (0,4.2) -- (0,5.4); \draw[cone] (0,6.6) -- (0,7.2);
\draw[cone,rounded corners=6pt,->] (0,7.2) -- (0,8.1) -- (5.6,8.1) -- (5.6,3.3);
\draw[cone,rounded corners=6pt] (5.6,3.5) -- (5.6,-0.9) -- (0,-0.9) -- (0,0);
\draw[czero] (2.4,0) -- (2.4,0.6); \draw[czero] (2.4,1.8) -- (2.4,3.0); \draw[czero] (2.4,4.2) -- (2.4,5.4); \draw[czero] (2.4,6.6) -- (2.4,7.2);
\draw[czero,rounded corners=6pt,->] (2.4,7.2) -- (2.4,7.6) -- (4.0,7.6) -- (4.0,3.3);
\draw[czero,rounded corners=6pt] (4.0,3.5) -- (4.0,-0.4) -- (2.4,-0.4) -- (2.4,0);
\foreach \y in {1.2,3.6,6.0} {
  \draw[czero] (2.4,{\y-0.6}) -- (1.2,\y); \draw[cone] (1.2,\y) -- (0,{\y+0.6});
  \draw[white,line width=5pt] (0.45,{\y-0.375}) -- (1.95,{\y+0.375});
  \draw[cone] (0,{\y-0.6}) -- (1.2,\y); \draw[czero] (1.2,\y) -- (2.4,{\y+0.6}); }
\node[font=\scriptsize] at (-0.45,1.2) {$1$}; \node[font=\scriptsize] at (-0.45,3.6) {$2$}; \node[font=\scriptsize] at (-0.45,6.0) {$3$};
\node[font=\scriptsize] at (1.2,2.4) {$\vth_{g_1}{=}1$}; \node[font=\scriptsize] at (1.2,4.8) {$\vth_{g_2}{=}1$}; \node[font=\scriptsize] at (4.8,5.6) {$\vth_{g_3}$}; \node[font=\scriptsize] at (4.8,5.0) {${=}\,1$};
\node[font=\scriptsize] at (3.2,5.6) {$\vth_{f_\infty}$}; \node[font=\scriptsize] at (3.2,5.0) {${=}\,t_0$}; \node[font=\scriptsize] at (-2.0,3.6) {$\vth_{f_0}=t_1$};
\node[font=\scriptsize,right] at (5.65,1.6) {$\hate$};
\draw[cone] (8.0,4.6) -- (9.4,4.6); \node[font=\scriptsize,right] at (9.5,4.6) {$A$: color $1$, $\mhol_e=t_1$};
\draw[czero] (8.0,3.4) -- (9.4,3.4); \node[font=\scriptsize,right] at (9.5,3.4) {$B$: color $0$, $\mhol_e=t_0^{-1}$};
\end{tikzpicture}
\caption{The diagram of the trefoil of Example~\ref{ex.trefoilcocores}, with its two Seifert circles $A$ (solid, color $1$) and $B$ (dashed, color $0$), which exchange the strands of the link at every crossing, and the values of the two-colored potential on the five faces. The holonomy around the co-core of every crossing is $\vth_{g_{j-1}}\vth_{g_j}/(\vth_{f_0}\vth_{f_\infty})=(t_0t_1)^{-1}$.}
\label{fig.trefoilcocores}
\end{figure}

\begin{example}[the trefoil]\label{ex.trefoilcocores}
Let $D$ be the diagram of the trefoil of Examples~\ref{ex.trefoilflows} and~\ref{ex.trefoilmixed}, the closure of $\sigma_1^3$, with the faces $f_0,g_1,g_2,g_3,f_\infty$ and the cut edge $\hate=u'_3=o_1$ (Figure~\ref{fig.trefoilcocores}). It has two Seifert circles: $A=A_1=A_2=A_3$, which contains the edges $o_j$ and $u'_j$ and bounds the face $f_0$, and $B=B_1=B_2=B_3$, which contains the edges $u_j$ and $o'_j$ and bounds the face $f_\infty$. We take $\mathrm{col}(A)=1$ and $\mathrm{col}(B)=0$. The three bigons $g_1,g_2,g_3$ lie between the two circles, $f_0$ lies on the left of $A$ and $f_\infty$ on the right of $B$, so that the two-colored potential is
\[
   \vth_{g_1}=\vth_{g_2}=\vth_{g_3}=1,\qquad \vth_{f_0}=t_1,\qquad \vth_{f_\infty}=t_0 ,
\]
with $\mhol_e=t_1$ on the edges of $A$ and $\mhol_e=t_0^{-1}$ on the edges of $B$. The three crossings are positive and of type out. At the crossing $j$ the four faces \eqref{eq.fourfaces} are $(f^1_j,f^2_j,f^3_j,f^4_j)=(g_{j-1},f_0,g_j,f_\infty)$, with $g_0=g_3$, and $\omega_j=\omega=(t_0t_1)^{-1}$. By \eqref{eq.pieces} the lifted meridian of the co-core $\alpha_j$ consists of the chains
\[
   \mer_j(0,1)=\tpth_{g_{j-1}},\qquad \mer_j(1,2)=-t_1^{-1}\,\tpth_{f_0},\qquad \mer_j(2,3)=t_1^{-1}\,\tpth_{g_j},\qquad \mer_j(3,4)=-\omega\,\tpth_{f_\infty}
\]
and of their multiples by the powers of $\omega$. The faces which are not adjacent to $\hate$ are $g_1,g_2,f_\infty$. Writing the coefficients of a chain on these three faces as a row, we find for $j=2$
\begin{align*}
   \mer_2(0,2)&:\ (1,0,0), & \mer_2(2,4)&:\ (0,t_1^{-1},-\omega),\\
   \mer_2(1,5)&:\ (\omega,t_1^{-1},-\omega), & \mer_2(3,7)&:\ (\omega,\omega t_1^{-1},-\omega),
\end{align*}
and similarly for $j=1,3$. So the matrices of the interleaved chains $\CZ_j=\CZ^\out_j$ are
\begin{gather*}
   \CZ_1|_{\hate}=\begin{pmatrix}t_1^{-1}\\0\\-\omega\end{pmatrix}\begin{pmatrix}t_1^{-1}&0&-\omega\end{pmatrix},\qquad
   \CZ_3|_{\hate}=\begin{pmatrix}0\\ \omega\\-\omega\end{pmatrix}\begin{pmatrix}0&1&-\omega\end{pmatrix},\\
   \CZ_2|_{\hate}=\begin{pmatrix}\omega\\ \omega t_1^{-1}\\-\omega\end{pmatrix}\begin{pmatrix}1&0&0\end{pmatrix}+\begin{pmatrix}\omega\\ t_1^{-1}\\-\omega\end{pmatrix}\begin{pmatrix}0&t_1^{-1}&-\omega\end{pmatrix} .
\end{gather*}
Multiplying the columns by $\vth_{g_1}=\vth_{g_2}=1$ and $\vth_{f_\infty}=t_0$ and using $\omega t_0t_1=1$, one finds the matrices $\omega t_1^{-1}\,\Psi_j$, where $\Psi_1,\Psi_2,\Psi_3$ are the matrices of Example~\ref{ex.trefoilmixed}, in accordance with Theorem~\ref{thm.cocores}~(i) and with $\us^\out_j=\omega\,\vth_1/\mhol_{u'_j}=\omega\,t_1^{-1}$. Hence
\[
   \lb\CZ_1,\CZ_2,\CZ_3\rb_{\hate}=\frac{(\omega\,t_1^{-1})^3}{t_0}\,\LG(3_1;t_0,t_1)=t_0^{-4}\,t_1^{-6}\,\LG(3_1;t_0,t_1).
\]
\end{example}

\begin{remark}[other chains]\label{rem.variants}
(i) \emph{Squares.} If the interleaved chain $\CZ_j$ of Theorem~\ref{thm.cocores} is replaced, for every $j$, by the square $\mer_j(i,i+4)\otimes\mer_j(i',i'+4)$ which reads $\pi_j$ at $x_j=t_0^{\ve_j}$ and $\tpi_j$ at $y_j=t_1^{\ve_j}$ (Proposition~\ref{prop.readings}), the mixed discriminant becomes, up to a unit, the product $P_{D,\hate}(t_0^{\ve})\,\tP_{D,\hate}(t_1^{\ve})=\det A_{D,\hate}(t_0)\,\det B_{D,\hate}(t_1)$ of two Alexander determinants, by Lemma~\ref{lem.cauchybinet}. So the Links--Gould polynomial and the product $\Delta_L(t_0)\Delta_L(t_1)$ are computed from the same lifted meridians, for the same local system; the only difference is the interleaving.

(ii) \emph{The other interleaved chain.} If $\CZ^\out_j$ and $\CZ^\inn_j$ are exchanged in Theorem~\ref{thm.cocores}, the forms become $\Psi_j[1;t_1^{-\ve_j},t_0^{-\ve_j}]$ up to units, by Lemma~\ref{lem.twocolored}~(iii) and \eqref{eq.inout}, and the mixed discriminant becomes $\LG(L;t_1^{-1},t_0^{-1})$ up to a unit. By \eqref{eq.mirror} this is the Links--Gould polynomial of the mirror image of $L$. The same happens if one keeps the chains and exchanges the two colors.

(iii) \emph{Alternating diagrams.} By Lemma~\ref{lem.twocolored}~(iv), if $D$ is alternating all the crossings have the same type, and for one of the two colorings all the chains are $\CZ^\out_j$ and all the holonomies around the co-cores are $(t_0t_1)^{-1}$. For a diagram which is not alternating both types occur, for either coloring.
\end{remark}

\subsection{The role of the plane of the diagram}\label{ssec.vertical}
Theorem~\ref{thm.cocores} computes $\LG(L)$ from the following objects.
\begin{enumerate}\renewcommand{\theenumi}{\alph{enumi}}
\item The canonical surface $F_D$ with its decomposition into disks and bands, and a coloring of its disks with two colors. They determine the graph $G$, which is the link with the co-cores of the bands, the space $X_G$, and the local system: its holonomy is trivial around the bands, and around the boundary of a disk it is $t_1$ or $t_0^{-1}$, according to the color of the disk (Remark~\ref{rem.dictionary}~(ii)).
\item For every band, the meridian of its co-core, with its two points above the under-strand and its two points below the over-strand, and paths which join the first two points to $\star_+$ and the last two to $\star_-$. In Lemma~\ref{lem.meridian} these paths lie in $X_+$ and in $X_-$, the two sides of the surface $\Sigma$ which is made of the faces of the diagram. Since $X_+$ and $X_-$ are contractible, this determines them up to homotopy.
\item The edge $\hate$, through the faces which are adjacent to it (Lemma~\ref{lem.chainpairing}).
\end{enumerate}
The objects (a) only depend on the surface and on its decomposition, up to isotopy. The paths of (b) and the edge of (c) depend on the diagram: they are given by the vertical direction. One would like to replace them by data which only depend on (a); this would turn Theorem~\ref{thm.cocores} into a definition of the Links--Gould polynomial from a Seifert surface with a decomposition into disks and bands. We do not know how to do this, and the following results show that the paths cannot simply be forgotten.

\emph{Dependence on the paths.} Recall the coefficients $\mathsf c_j(i)$ of the lifted points, defined in \eqref{eq.liftedpoints}.

\begin{proposition}[changing a path]\label{prop.paths}
Let $\vth$ be a multiplicative potential and $j$ a crossing; write $\omega=\omega_j$, $\mer=\mer_j$ and $\mathsf c=\mathsf c_j$.
\begin{enumerate}\renewcommand{\theenumi}{\roman{enumi}}
\item Let $\gamma$ be a loop based at $\star_+$, the path of a word, with holonomy $\mathsf g=\rho_\vth(\gamma)$, and let $\tilde\gamma\in C^+_\vth$ be its lift which starts at $\tstar_+$. Suppose that the path which joins the meridian of $\alpha_j$ to $\star_+$ is changed by $\gamma$, so that the loops based at $\star_+$ of which the chains $\mer(i,l)$ with $i,l$ even are lifts are conjugated by $\gamma$. Then $\mer(i,l)$ is replaced by $\mathsf g\,\mer(i,l)+\bigl(\mathsf c(i)-\mathsf c(l)\bigr)\,\tilde\gamma$ for $i,l$ even, and the interleaved chains are replaced by
\begin{gather*}
   \mathsf g\,\CZ^\out_j+(1-\omega)\,\bigl(\mer(1,5)+(\mathsf c(2)-1)\,\mer(1,3)\bigr)\otimes\tilde\gamma\qquad\text{and}\\
   \mathsf g\,\CZ^\inn_j+(1-\omega)\,\bigl(\mer(-1,3)+(\mathsf c(2)-1)\,\mer(-1,1)\bigr)\otimes\tilde\gamma .
\end{gather*}
\item Let $\gamma$ be a loop based at $\star_-$ with holonomy $\mathsf g$, and let $\tilde\gamma\in C^-_\vth$ be its lift which starts at $\tstar_-$. If the path which joins the meridian of $\alpha_j$ to $\star_-$ is changed by $\gamma$, then $\mer(i,l)$ is replaced by $\mathsf g\,\mer(i,l)+\bigl(\mathsf c(i)-\mathsf c(l)\bigr)\,\tilde\gamma$ for $i,l$ odd, and the interleaved chains are replaced by
\begin{gather*}
   \mathsf g\,\CZ^\out_j+(1-\omega)\,\tilde\gamma\otimes\bigl(\mathsf c(3)\,\mer(0,2)+\mathsf c(1)\,\mer(2,4)\bigr)\qquad\text{and}\\
   \mathsf g\,\CZ^\inn_j+(1-\omega)\,\tilde\gamma\otimes\bigl(\mathsf c(1)\,\mer(0,2)+\mathsf c(-1)\,\mer(2,4)\bigr).
\end{gather*}
\end{enumerate}
\end{proposition}
\begin{proof}
Let $\gamma'$ be a loop based at $\star=\star_+$ or $\star_-$, the path of a word, and suppose that its lift $\tilde\gamma'$ which starts at $a\,\tstar$ ends at $b\,\tstar$. The lift of $\gamma\,\gamma'\,\gamma^{-1}$ which starts at $a\,\tstar$ consists of the lift $a\,\tilde\gamma$ of $\gamma$, which ends at $\mathsf ga\,\tstar$, of the lift $\mathsf g\,\tilde\gamma'$ of $\gamma'$, which ends at $\mathsf gb\,\tstar$, and of the lift $-b\,\tilde\gamma$ of $\gamma^{-1}$: it is $\mathsf g\,\tilde\gamma'+(a-b)\,\tilde\gamma$. For $\mer(i,l)$ we have $a=\mathsf c(i)$ and $b=\mathsf c(l)$. In (i), Definition~\ref{def.interleaved} with the new chains $\mer(0,2)$ and $\mer(2,4)$ gives $\mathsf g\,\CZ^\out_j+\bigl((1-\mathsf c(2))\,\mer(3,7)+(\mathsf c(2)-\omega)\,\mer(1,5)\bigr)\otimes\tilde\gamma$, since $\mathsf c(0)=1$ and $\mathsf c(4)=\omega$, and by \eqref{eq.turndiff}
\[
   (1-\mathsf c(2))\,\mer(3,7)+(\mathsf c(2)-\omega)\,\mer(1,5)=(1-\omega)\,\mer(1,5)+(1-\mathsf c(2))(\omega-1)\,\mer(1,3) ;
\]
the computation for $\CZ^\inn_j$ is the same, with $\mer(1,5)$, $\mer(-1,3)$, $\mer(-1,1)$ in place of $\mer(3,7)$, $\mer(1,5)$, $\mer(1,3)$. In (ii) the new turns are $\mathsf g\,\mer(i,i+4)+(1-\omega)\,\mathsf c(i)\,\tilde\gamma$ for $i$ odd, and we substitute them in Definition~\ref{def.interleaved}.
\end{proof}

So a change of one path of the crossing $j$ replaces $\CZ_j$ by $\mathsf g\,\CZ_j+(1-\omega_j)\,\CZ'_j$, where $\CZ'_j$ is the tensor product of two chains, and by the multilinearity of the mixed discriminant it replaces $\lb\CZ_1,\dots,\CZ_c\rb_{\hate}$ by
\[
   \mathsf g\,\lb\CZ_1,\dots,\CZ_c\rb_{\hate}+(1-\omega_j)\,\lb\CZ_1,\dots,\CZ'_j,\dots,\CZ_c\rb_{\hate} .
\]
For the two-colored potential the first term is $\LG(L)$ up to a unit, and $1-\omega_j$ is $1-t_0t_1$ up to a unit. So the dependence on the paths is entirely in a term which is divisible by $1-t_0t_1$: modulo $1-t_0t_1$, where one computes the square of the Alexander polynomial, there is no dependence. There is no reason for the second term to vanish, and the next example shows that it does not.

\emph{Surfaces in other positions.} To test the role of the paths we need surfaces with disks and bands which are placed differently with respect to a diagram: the bands are still crossings of the diagram, with vertical co-cores, but the diagram has other crossings.

\begin{definition}\label{def.recipe}
Let $\CJ$ be a set of crossings of $D$. Smoothing the crossings of $\CJ$ along the orientations and keeping the other crossings gives a link diagram, whose components we call the \emph{$\CJ$-circles}; every edge of $D$ lies on one of them. Suppose that the $\CJ$-circles can be colored with the colors $0$ and $1$ in such a way that the two arcs of the smoothing of every crossing of $\CJ$ lie on two different $\CJ$-circles, of different colors, and let $\mathrm{col}$ be such a coloring. Let $\vth$ be the multiplicative potential, with values in $\BZ[t_0^{\pm1},t_1^{\pm1}]$, such that $\mhol_e=t_1$ if the $\CJ$-circle of $e$ has the color $1$ and $\mhol_e=t_0^{-1}$ if it has the color $0$. For an edge $\hate$ of $D$ we put
\[
   \CV(D,\CJ,\mathrm{col};\hate)=(-1)^{\mu-1}\,\BE_{D,\CJ}\bigl[P_{D,\hate}(x)\,\tP_{D,\hate}(y)\bigr]\ \in\BZ[t_0^{\pm1},t_1^{\pm1}],
\]
where $\BE_{D,\CJ}$ is the linear map which is $\BE_D$ on the variables of the crossings of $\CJ$, that is, $x_j\mapsto t_0^{\ve_j}$, $y_j\mapsto t_1^{\ve_j}$, $x_jy_j\mapsto1$ for $j\in \CJ$, and which is an evaluation of the variables of the other crossings: $x_k\mapsto \mhol_{o_k}^{-\ve_k}$ and $y_k\mapsto \mhol_{u_k}^{\ve_k}$ for $k\notin \CJ$.
\end{definition}

The potential $\vth$ exists and is unique up to a common factor, as in Lemma~\ref{lem.twocolored}~(ii): a $\CJ$-circle follows a strand through the crossings which are not in $\CJ$ and passes from $o_j$ to $u'_j$, or from $u_j$ to $o'_j$, at the crossings of $\CJ$, so that the sum of the edges of the $\CJ$-circles of a given color is a flow. By Corollary~\ref{cor.anypotential}, $\CV(D,\CJ,\mathrm{col};\hate)$ is, up to a unit, the mixed discriminant $\lb\CZ_1,\dots,\CZ_c\rb_{\hate}$ of the interleaved chains for this potential, chosen as in Theorem~\ref{thm.cocores} at the crossings of $\CJ$; at a crossing $k\notin \CJ$ the holonomy is the same along each strand, so that $\omega_k=1$ and the two interleaved chains coincide with the squares (Proposition~\ref{prop.interleaved}~(i)), with $\sv_k=\mhol_{o_k}^{-\ve_k}\mhol_{u_k}^{\ve_k}$.

Here is the geometric meaning of the definition. Suppose that the $\CJ$-circles bound disjoint disks in $\BR^3$, and let $F$ be the surface made of these disks and of half-twisted bands at the crossings of $\CJ$. Then $F$ is a Seifert surface of $L$, with a decomposition into disks and bands and with a coloring of its disks; the crossings of $\CJ$ are its bands, and the other crossings are places where two pieces of its boundary pass over each other. The local system $\CL_\vth$ extends across the co-cores of the crossings which are not in $\CJ$, and it is the two-colored local system of $F$ in the sense of (a) above. So $\CV(D,\CJ,\mathrm{col};\hate)$ is the construction of Theorem~\ref{thm.cocores} for the surface $F$, in which the paths of (b) are the vertical paths of the diagram $D$ and the edge of (c) is the edge $\hate$ of $D$. If $\CJ$ is the set of all the crossings, then $F=F_D$ and $\CV(D,\CJ,\mathrm{col};\hate)=\LG(L)$.

The first test is positive: an isotopy of the surface which creates two crossings between pieces of its boundary, away from the bands, does not change the value.

\begin{proposition}[a move of type II]\label{prop.typeII}
Let $D$ be a closed braid diagram, and let $D'$ be obtained from $D$ by inserting $\sigma_i\sigma_i^{-1}$ or $\sigma_i^{-1}\sigma_i$ at some place of the braid word. Let $\CJ$ be the set of the crossings of $D$, viewed as crossings of $D'$, and let $\hate$ be an edge of $D'$ outside the bigon. Then the $\CJ$-circles are the Seifert circles of $D$, and for each of their two colorings
\[
   \CV(D',\CJ,\mathrm{col};\hate)=\LG(L;t_0,t_1).
\]
\end{proposition}
\begin{proof}
Let $a$ and $b$ be the two new crossings. They have the same over-strand and the same under-strand, and opposite signs. Since $a,b\notin \CJ$, the $\CJ$-circles follow these two strands through $a$ and $b$, so they are the Seifert circles of $D$, and $\mhol_{o_a}=\mhol_{o_b}$ and $\mhol_{u_a}=\mhol_{u_b}$. Hence $\BE_{D',\CJ}$ sends $x_ax_b$ and $y_ay_b$ to $1$. By Proposition~\ref{prop.hollowhandle}, $P_{D',\hate}=P_{D,\hate}+(x_ax_b-1)R$ and $\tP_{D',\hate}=\tP_{D,\hate}+(y_ay_b-1)\tilde R$, where $P_{D,\hate}$, $\tP_{D,\hate}$, $R$ and $\tilde R$ do not contain the variables of $a$ and $b$. Since $\BE_{D',\CJ}$ is an evaluation of the variables of $a$ and $b$ and is $\BE_D$ on the other variables, $\BE_{D',\CJ}[P_{D',\hate}\tP_{D',\hate}]=\BE_D[P_{D,\hate}\tP_{D,\hate}]$, and the claim follows from Theorem~\ref{thm.HK}.
\end{proof}

\begin{figure}[ht]
\centering
\begin{tikzpicture}[scale=0.92,every node/.style={font=\small}]
\def\bandpicture{%
  \draw[thick] (-2.2,0) -- (2.2,0); \node[above] at (-1.5,0) {under-strand};
  \draw[very thick] (0,0) -- (0,1.1); \node[right] at (0.02,0.55) {$\alpha_j$};
  \fill (0,1.1) circle (2.4pt); \node[right] at (0.15,1.15) {over-strand};
  \fill (0,2.5) circle (1.5pt) node[above] {$\star_+$}; \fill (0,-3) circle (1.5pt) node[below] {$\star_-$};
  \draw[blue,semithick] (0,2.5) -- (0,1.1);}
\begin{scope}
  \bandpicture
  \draw[blue,semithick] (0,-3) -- (0,0);
  \draw[fill=white] (-1,-1.7) circle (2.6pt); \node[left,font=\scriptsize] at (-1.1,-1.7) {strand};
  \draw[->,densely dotted] (-0.8,-1.7) -- (0.85,-1.7);
\end{scope}
\begin{scope}[xshift=7.6cm]
  \bandpicture
  \draw[red,semithick] (0,-3) -- (0,0);
  \draw[blue,semithick,dashed] (0,-3) .. controls (2.3,-2.6) and (2.3,-0.9) .. (0,0);
  \draw[fill=white] (1,-1.7) circle (2.6pt);
  \node[red,left] at (0,-0.9) {vertical}; \node[blue,right] at (1.75,-1.7) {carried along};
\end{scope}
\end{tikzpicture}
\caption{A strand of the link (hollow dot) is pushed across the vertical line of the crossing $j$, below the under-strand, seen from the side: the under-strand of the crossing lies in the plane of the figure, the over-strand and the strand which is pushed cross it. The vertical segments are the paths which join the meridian of the co-core $\alpha_j$ to the two base points, up to homotopy. Left: before the isotopy. Right: after it; the image of the lower path under the isotopy (dashed) and the vertical path of the new diagram differ by a meridian of this strand.}
\label{fig.push}
\end{figure}
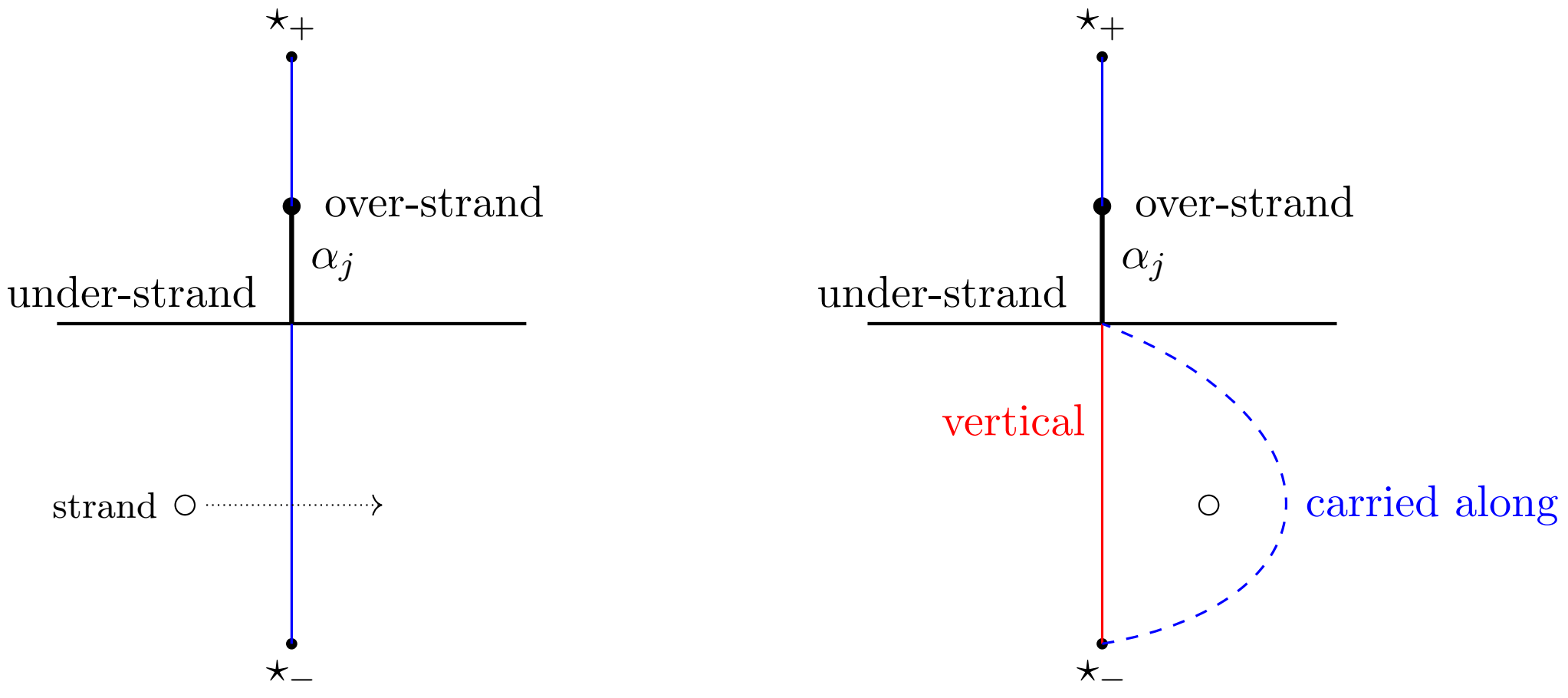

The second test is negative. Let $w$ be a braid word on $n$ strands and $D$ its closure, so that $F_D$ consists of $n$ disks, one for each position of the braid, and of one band for each letter. We replace one letter $\sigma_i^{\ve}$ of $w$, where $\ve=\pm1$ is the sign of its crossing, by the braid word
\begin{equation}\label{eq.conjugate}
   \sigma_k^{\ve'}\,\sigma_i^{\ve'}\ \sigma_k^{\ve}\ \sigma_i^{-\ve'}\,\sigma_k^{-\ve'},\qquad k=i\pm1,\quad \ve'=\pm1,
\end{equation}
which represents the same braid, by the relation $\sigma_k\sigma_i\sigma_k=\sigma_i\sigma_k\sigma_i$. Let $D'$ be the closure of the new word, and let $\CJ$ be the set of the letters of $w$, the letter which was replaced being now the middle letter $\sigma_k^{\ve}$ of \eqref{eq.conjugate}. In $D'$ the strand next to the two strands of the crossing has been pushed across the crossing, over it or under it according to $\ve'$. The $\CJ$-circles of $D'$ are the Seifert circles of $D$, one of which has been pushed across a band: the surface $F$ of $(D',\CJ)$ is the surface $F_D$, with the same disks and the same bands, moved by an isotopy. But the strand has crossed the vertical line through the crossing, so that the vertical path of this band in $D'$ is not the image under the isotopy of its vertical path in $D$: the two differ by a meridian of the strand (Figure~\ref{fig.push}).

\begin{example}[a disk bounded by the unknot]\label{ex.unknotdisc}
Let $w=\sigma_1\sigma_2\sigma_3$ on four strands. Its closure is a diagram of the unknot, and its canonical surface consists of four disks and three bands: it is a disk. Replacing the letter $\sigma_2$ by \eqref{eq.conjugate} with $k=3$ and $\ve'=1$ gives the braid word
\[
   w'=\sigma_1\cdot\sigma_3\sigma_2\,\sigma_3\,\sigma_2^{-1}\sigma_3^{-1}\cdot\sigma_3 ,
\]
and $\CJ$ consists of the letters $1$, $4$ and $7$ of $w'$: the fourth strand has been pushed under the second band. Let $D'$ be the closure of $w'$ and let $\hate$ be the edge which closes the first strand of the braid; it does not touch the four new crossings. For the coloring in which the disk of the first position has the color $1$,
\[
   \CV(D',\CJ,\mathrm{col};\hate)=1+t_0-t_0^2t_1=1+t_0\,(1-t_0t_1),
\]
whereas $\LG=1$ for the unknot. If the lower path of the band of the letter $4$ is changed by the meridian of the edge of the fourth strand which runs beside this band, as in Proposition~\ref{prop.paths}~(ii), the value becomes a unit. For the other coloring the value at this edge is the unit $t_0t_1$.
\end{example}

So the construction, with the vertical paths of a diagram in which the bands are crossings, is not invariant under the isotopies of a surface with disks and bands. This is not a peculiarity of the example. For $60$ random examples of the same kind (braid words with at most five letters on three or four strands, random choices of the letter, of $k$ and of $\ve'$) we computed the values for all the edges $\hate$ which do not touch the four new crossings. With the vertical paths of $D'$ the value is $\LG(L)$ up to a unit for all these edges in $33$ examples, and not in the $27$ other ones. With the path of the band changed by the meridian of the strand which was pushed, by the formulas of Proposition~\ref{prop.paths}, the value is $\LG(L)$ up to a unit for all these edges in all the $60$ examples, as it should be: the construction with the paths carried along is the image of the construction for $D$ under the isotopy.

\begin{remark}[coarser decompositions]\label{rem.coarser}
The canonical surface has other decompositions into disks and bands. Let $\CJ$ be a set of crossings of $D$ whose complement contains no cycle of the Seifert graph, so that the bands of the crossings outside $\CJ$ and the disks which they join form disjoint disks. Then $F_D$ is also the union of these bigger disks and of the bands of $\CJ$; for instance a chain of $p$ bands and $p-1$ small disks between two disks becomes one band with $p$ half-twists, whose co-core is placed at one of the $p$ crossings. When the $\CJ$-circles can be two-colored, Definition~\ref{def.recipe} applies. For the diagram with five crossings of the knot $5_2$ with the planar diagram code $X_{1,5,2,4}\,X_{3,9,4,8}\,X_{5,1,6,10}\,X_{7,3,8,2}\,X_{9,7,10,6}$, the crossings $1,3,5$ form a band with three half-twists between two disks, so that $\CJ=\{2,4,5\}$, $\{2,3,4\}$ or $\{1,2,4\}$. For each of these three sets, for each of the two colorings and for each of the ten edges $\hate$, the value $\CV(D,\CJ,\mathrm{col};\hate)$ is not $\LG(5_2)$ up to a unit; as $\hate$ varies it takes $7$, $10$ and $7$ different values up to units; and changing the two paths of the twisted band by the meridian of an arbitrary edge of the diagram, in either direction, as in Proposition~\ref{prop.paths}, never gives $\LG(5_2)$ either. More generally, for $97$ diagrams with at most eight crossings and all their sets $\CJ$ as above, a decomposition which is coarser than that of Seifert's algorithm gave $\LG(L)$ for all the edges $\hate$ in each of the $47$ cases in which all the crossings outside $\CJ$ are \emph{nugatory}, that is, have two opposite corners in the same face, and in only three of the $237$ other cases; and when the signs of the crossings outside $\CJ$ which are not nugatory do not add up to zero, it never did ($196$ cases).
\end{remark}

In conclusion, Theorem~\ref{thm.cocores} expresses $\LG(L)$ by objects which make sense for every Seifert surface with a decomposition into disks and bands and a coloring of its disks: the complement of the link and of the co-cores of the bands, the local system which is trivial around the bands and has the holonomies $t_1$ and $t_0^{-1}$ around the boundaries of the disks, and the interleaved lifts of the meridians of the co-cores. But it also uses the paths which join these meridians to two base points, and an edge; in the theorem these are given by the plane of a diagram for which the surface is the canonical surface. The results of this subsection show that these data cannot be ignored: the value depends on the paths, through a term which is divisible by $1-t_0t_1$ (Proposition~\ref{prop.paths}); the prescription ``the vertical paths of a diagram in which the bands are crossings'' is compatible with the moves of type II, which do not cross the paths (Proposition~\ref{prop.typeII}), but not with the isotopies which do (Example~\ref{ex.unknotdisc}); and for the decompositions which are not those of Seifert's algorithm it does not give $\LG$ (Remark~\ref{rem.coarser}). So the theorem is a formula for $\LG$ in terms of a diagram, written in the complement of the link and of the co-cores. It would become a definition of $\LG$ from a canonical surface if the paths and the edge $\hate$ could be characterized in terms of the surface, of its disks and of its bands. Two limitations are already visible. Modulo $1-t_0t_1$ everything is canonical and gives $\Delta_L^2$, so that such a characterization has to be sensitive to more than the homology of the infinite cyclic covering. And not every Seifert surface is isotopic to a canonical surface: the complement of a canonical surface is a handlebody, which is not the case for all Seifert surfaces of minimal genus \cite{Moriah87,Livingston88,KobayashiKobayashi96}.

\part{Faces: Kashaev's matrix, the Gordon--Litherland forms and the Kauffman model}\label{part.faces}

In this part the crossing forms are written in the basis of the boundaries of the faces of the diagram, the basis of Alexander's original matrix (Section~\ref{ssec.facematrices}). Section~\ref{sec.faces} compares the matrices so obtained with Kashaev's matrix. Section~\ref{sec.minusone} studies the point where all the crossing variables are equal to $-1$, at which the forms of Gordon and Litherland of all the state surfaces of Section~\ref{sec.statesurfaces} appear. Section~\ref{sec.kauffman} expands the two determinants over the Kauffman states of the diagram. Section~\ref{sec.minusone} uses Section~\ref{sec.faces} only through the winding numbers and the definitions of Kashaev's matrix and of its deformation; Section~\ref{sec.kauffman} uses the winding numbers of Section~\ref{ssec.windings} and, for alternating diagrams, the coloring and the signature formula of Section~\ref{ssec.altminusone}.

\section{Kashaev's matrix}\label{sec.faces}

In the basis of the boundaries of the faces, the crossing forms are Alexander's matrix (Theorem~\ref{thm.faceform}). More recently, Kashaev \cite{Kashaev21} attached to a diagram a \emph{symmetric} matrix $\tau_D$ indexed by the faces, whose entries lie in $\BZ[2\zeta]$, $2\zeta =t^{1/2}+t^{-1/2}$, and conjectured that its signature computes the Levine--Tristram signature of $L$. The conjecture was proved by Liu \cite{Liu25}, by constructing a Seifert surface for $L\sqcup rL$ whose first homology is generated by the faces and the vertices of $D$ and whose Seifert form has $\tau_D$ as a Schur complement; Cimasoni and Ferretti \cite{CimasoniFerretti24} proved that $\tau_D$ is a Gram matrix $K(D)^{\T}SK(D)$ of Kauffman's crossing-by-face matrix $K(D)$ with respect to the diagonal matrix $S$ of the crossing signs, and both papers show that the determinant of $\tau_D$ with two adjacent faces removed is $\Delta_L(t)^2$ up to a unit \cite[Theorem~1]{CimasoniFerretti24}, \cite[Corollary~3.3]{Liu25}.

Since $\LG(L;t,t^{-1})=\Delta_L(t)^2$, this suggests a relation between Kashaev's matrix and the determinant formula, and the purpose of this section is to make it precise. Kashaev's matrix is, up to the winding numbers of the faces, the Gram matrix of the face matrices of the two flavors at the free fermion point (Theorem~\ref{thm.kashaev}); in the language of Section~\ref{sec.flows}, it is the matrix of the bilinear form $\sum_j\ve_j\,\pi_j\otimes\tpi_j$ of Remark~\ref{rem.flows}~(iii) on the boundaries of the faces, at $x=t^{\ve}$ and $y=t^{-\ve}$. When the variables are not specialized, the same matrix is a two-variable deformation $\tau_D(x,y)$ of Kashaev's matrix whose coupling by $\BE_D$ is $\LG$ (Corollary~\ref{cor.deformedkashaev}). The deformation breaks the symmetry of $\tau_D$: the rows of $\tau_D(x,y)$ belong to the first flavor and its columns to the second one, and $\BE_D$ couples them at every crossing. The value of Kashaev's matrix at $t=-1$, where it becomes the direct sum of the two Goeritz matrices of the diagram, is the subject of Section~\ref{sec.minusone}.

\subsection{Winding numbers}\label{ssec.windings}
The comparison of Kashaev's matrix with the face matrices involves the winding numbers of the diagram around its faces. In this subsection we relate them to the levels of the Seifert circles (Lemma~\ref{lem.winding}), and we show that, together with the levels, they turn the face matrix of the first flavor at $x=t^{\ve}$ into the face matrix of the second flavor at $y=t^{-\ve}$ (Lemma~\ref{lem.facelevel}).

We use the notation of Sections~\ref{ssec.flowlattice} and~\ref{ssec.facematrices}: the set $\CF$ of the faces of $D$, the faces $f_{\rm L}(e)$ and $f_{\rm R}(e)$ on the two sides of an edge $e$, the four faces $f^\uparrow_j,f^\leftarrow_j,f^\downarrow_j,f^\rightarrow_j$ at a crossing, and the face matrices $Q_D(x)$ and $\tQ_D(y)$, whose rows $a_j(x_j)$ and $b_j(y_j)$ are given by the table \eqref{eq.facetable}.

\begin{lemma}[winding numbers]\label{lem.winding}
There is a function $w\colon\CF\to\BZ$, unique up to an additive constant, such that $w(f_{\rm L}(e))=w(f_{\rm R}(e))+1$ for every edge $e$. If $D$ is drawn in the plane as in Section~\ref{ssec.level} and $w$ is normalized by $w=0$ on the unbounded face, then $w(f)$ is the winding number of $D$ around the points of $f$, and
\begin{equation}\label{eq.wh}
   w(f_{\rm L}(e))=-h(e),\qquad w(f_{\rm R}(e))=-h(e)-1\qquad\text{for every edge $e$},
\end{equation}
where $h$ is the level function of Definition~\ref{def.level}. At every crossing $j$,
\begin{equation}\label{eq.wcrossing}
   w(f^\leftarrow_j)=w(f^\uparrow_j)+1=w(f^\downarrow_j)+1=w(f^\rightarrow_j)+2 .
\end{equation}
\end{lemma}
\begin{proof}
The winding number of the oriented curves of $D$ around a point increases by one when the point crosses a strand from its right side to its left side, and the winding number around the unbounded face is $0$; so the winding number has the stated property, and uniqueness up to a constant follows from the connectedness of the dual graph. Let $e$ be an edge, lying on the Seifert circle $C$. The winding number around a point $p$ is $-\sum s(C')$, the sum running over the Seifert circles $C'$ enclosing $p$, since a clockwise circle ($s(C')=1$) winds $-1$ times around the points it encloses, and smoothing the crossings does not change winding numbers. If $s(C)=-1$, the left side of $C$ is its interior, the circles enclosing a point just to the left of $e$ are $C$ and the circles enclosing $C$, and $w(f_{\rm L}(e))=1-\sum_{C'}s(C')=-h(C)$ by Definition~\ref{def.level}, the sum running over the circles $C'\ne C$ enclosing $C$; if $s(C)=1$, the left side of $C$ is its exterior, the circles enclosing a point just to the left of $e$ are those enclosing $C$, and $w(f_{\rm L}(e))=-\sum_{C'}s(C')=-h(C)$ again. This proves \eqref{eq.wh}, and \eqref{eq.wcrossing} follows from $w(f_{\rm L})=w(f_{\rm R})+1$ applied to the two upward strands.
\end{proof}

At the Alexander point the two face matrices are related by the levels, as the two edge matrices are (Lemma~\ref{lem.level}~(i)).

\begin{lemma}\label{lem.facelevel}
$\tQ_D(t^{-\ve})=\diag\bigl(t^{-h(u_j)-(1+\ve_j)/2}\bigr)_{j}\ Q_D(t^{\ve})\ \diag\bigl(t^{-w(f)}\bigr)_{f\in\CF}$.
\end{lemma}
\begin{proof}
Write the entries in the order $(f^\uparrow_j,f^\leftarrow_j,f^\downarrow_j,f^\rightarrow_j)$, put $k=w(f^\leftarrow_j)$, and use \eqref{eq.facetable} and \eqref{eq.wcrossing}. For $\ve_j=1$, the row $a_j(t)^{\T}\diag(t^{-w})$ is $(-t\cdot t^{1-k},\,t\cdot t^{-k},\,-t^{1-k},\,t^{2-k})=t^{1-k}(-t,1,-1,t)$ and $b_j(t^{-1})^{\T}=(-1,t^{-1},-t^{-1},1)=t^{-1}(-t,1,-1,t)$; so $b_j(t^{-1})^{\T}=t^{k-2}\,a_j(t)^{\T}\diag(t^{-w})$, and $k-2=w(f^\rightarrow_j)=w(f_{\rm R}(u_j))=-h(u_j)-1$ by \eqref{eq.wh}. For $\ve_j=-1$, the row $a_j(t^{-1})^{\T}\diag(t^{-w})$ is $(t^{-1}\cdot t^{1-k},\,-t^{-k},\,t^{1-k},\,-t^{-1}\cdot t^{2-k})=t^{-k}(1,-1,t,-t)$ and $b_j(t)^{\T}=(1,-1,t,-t)$; so $b_j(t)^{\T}=t^{k}\,a_j(t^{-1})^{\T}\diag(t^{-w})$, and $k=w(f^\leftarrow_j)=w(f_{\rm L}(u_j))=-h(u_j)$. In both cases the exponent is $-h(u_j)-(1+\ve_j)/2$.
\end{proof}

\subsection{Kashaev's matrix at the free fermion point}\label{ssec.kashaev}
We recall the definition of Kashaev's matrix, and we prove the main result of this section: after a conjugation by the winding numbers, Kashaev's matrix is the product of the face matrices of the two flavors at the free fermion point, with the signs of the crossings in between (Theorem~\ref{thm.kashaev}). The formula for its determinant and the description of its kernel follow (Corollary~\ref{cor.kashaev}).

Let $\zeta$ be an indeterminate; it is Kashaev's variable, written $x$ in \cite{Kashaev21, Liu25, CimasoniFerretti24}.

\begin{definition}[{Kashaev \cite[Section~2]{Kashaev21}}]\label{def.kashaev}
The \emph{Kashaev matrix} of $D$ (the letter $\tau$ is Kashaev's; it should not be confused with the trees $\tau$ of Section~\ref{ssec.alternating}) is the symmetric $\CF\times\CF$ matrix with entries in $\BZ[\zeta]$
\[
   \tau_D=\sum_{j=1}^c\ve_j\,\tau_j,
\]
where $\tau_j$ is zero outside the four faces of the crossing $j$ and is given on them, in the basis $(f^\uparrow_j,f^\leftarrow_j,f^\downarrow_j,f^\rightarrow_j)$, by the block
\[
   \begin{pmatrix} 2\zeta^2-1 & \zeta  & 1 & \zeta\\ \zeta  & 1 & \zeta  & 1\\ 1 & \zeta  & 2\zeta^2-1 & \zeta\\ \zeta  & 1 & \zeta  & 1\end{pmatrix}
\]
(when faces coincide, the corresponding rows and columns are added). The diagonal entries are $2\zeta^2-1$ at $f^\uparrow_j$ and $f^\downarrow_j$ and $1$ at $f^\leftarrow_j$ and $f^\rightarrow_j$; the entry is $1$ between the two opposite faces $f^\uparrow_j,f^\downarrow_j$ and between $f^\leftarrow_j,f^\rightarrow_j$, and $\zeta$ between two faces adjacent along an edge of the crossing.
\end{definition}

The blocks $\tau_j$ have entries $\zeta$ and $2\zeta^2-1$, which do not lie in $\BZ[2\zeta]$, but their sum $\tau_D$ does \cite[Remark~2.2]{Liu25}: two faces adjacent along an edge receive $\pm \zeta$ from each of the two crossings at the ends of the edge, in total $0$ or $\pm2\zeta$, and a face has as many corners of type $f^\uparrow$ as of type $f^\downarrow$ (see the proof of Theorem~\ref{thm.kashaev} below), so that its diagonal entry is an even multiple of $2\zeta^2-1$, i.e.\ a multiple of $(2\zeta)^2-2$, plus an integer.

Kashaev proved that $\tau_D$, together with the writhe, is invariant under Reidemeister moves up to a modified $S$-equivalence, and conjectured that $\sign\tau_D-\mathrm{wr}(D)=2\sigma_L(\omega)$ for $\omega\in S^1\sm\{1\}$ and $2\zeta =\omega^{1/2}+\bar\omega^{1/2}$, where $\sigma_L$ is the Levine--Tristram signature \cite{Kashaev21}; this is Liu's theorem \cite[Theorem~3.1]{Liu25}. We set $2\zeta =t^{1/2}+t^{-1/2}$, so that $2\zeta^2-1=\frac12(t+t^{-1})$, and we write $\tau_D(t)$ for the matrix with entries in $\BZ[t^{\pm1/2}]$ so obtained; for $\omega\in S^1$ and $t=\omega$ this is the matrix of Kashaev's conjecture. Let
\begin{equation}\label{eq.W}
   \WW_t=\diag\bigl((-t^{1/2})^{w(f)}\bigr)_{f\in\CF}\qquad\text{and}\qquad \DD_\ve=\diag(\ve_1,\dots,\ve_c),
\end{equation}
where $w$ is the winding function of Lemma~\ref{lem.winding}; a change of the additive constant in $w$ multiplies $\WW_t$ by a scalar and does not affect the conjugation by $\WW_t$.

\begin{theorem}[Kashaev's matrix is the Gram matrix of the two flavors]\label{thm.kashaev}
For every diagram $D$,
\[
   \tau_D(t)=\WW_t^{-1}\ Q_D(t^{\ve})^{\T}\,\DD_\ve\,\tQ_D(t^{-\ve})\ \WW_t
   =\WW_t^{-1}\Bigl(\sum_{j=1}^c\ve_j\,a_j(t^{\ve_j})\,b_j(t^{-\ve_j})^{\T}\Bigr)\WW_t .
\]
\end{theorem}
The proof shows where the symmetry comes from. Conjugated by $\WW_t$, the product $\ve_ja_j(t^{\ve_j})b_j(t^{-\ve_j})^{\T}$ of the two relations of the crossing $j$ is a \emph{corner matrix} supported on the four faces of $j$, which is Kashaev's block $\ve_j\tau_j(t)$ up to a correction $\vartheta\Delta_j$ with $\vartheta=\frac12(t^{1/2}-t^{-1/2})$; the corrections cancel globally, $\sum_j\Delta_j=0$, because along the boundary of a face the corners at which both edges leave the crossing alternate with those at which both edges enter it, and because every edge leaves one crossing and enters another.

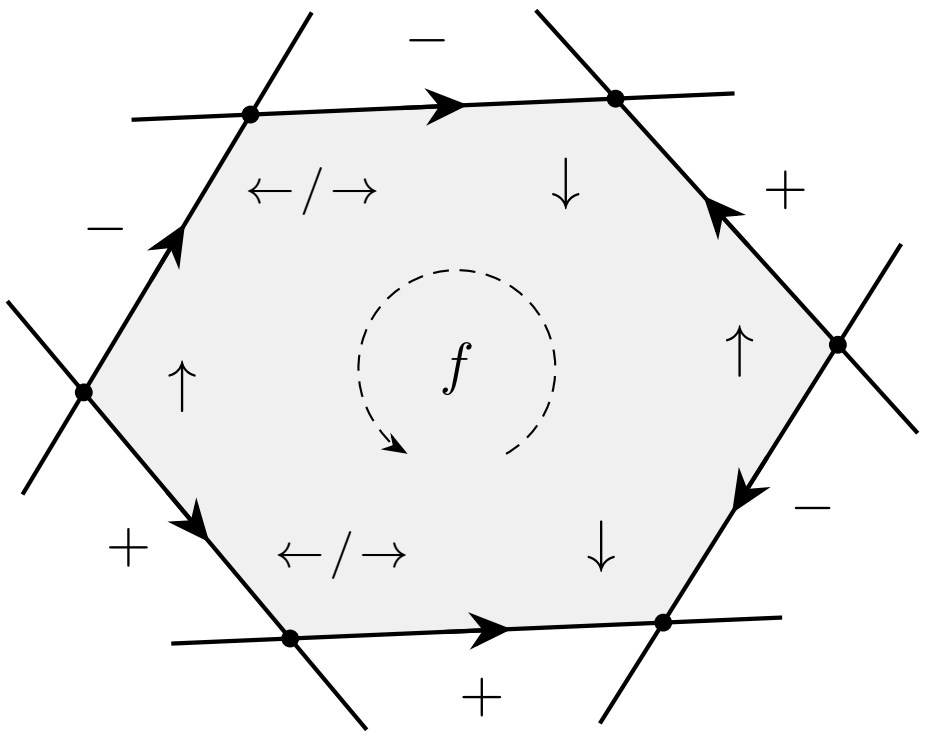
\begin{figure}[ht]
\centering
\begin{tikzpicture}[>=Stealth,scale=1.05,every node/.style={font=\small}]
\fill[gray!12] (-2.300,-0.200) -- (-1.000,-1.750) -- (1.350,-1.650) -- (2.450,0.100) -- (1.050,1.650) -- (-1.250,1.550) -- cycle;
\draw[thick] (-2.782,0.375) -- (-0.518,-2.325);
\draw[thick] (-1.749,-1.782) -- (2.099,-1.618);
\draw[thick] (0.951,-2.285) -- (2.849,0.735);
\draw[thick] (2.953,-0.457) -- (0.547,2.207);
\draw[thick] (1.799,1.683) -- (-1.999,1.517);
\draw[thick] (-0.864,2.193) -- (-2.686,-0.843);
\draw[thick,-{Stealth[length=8pt,width=7pt]}] (-1.779,-0.822) -- (-1.509,-1.144);
\node at (-2.019,-1.176) {$+$};
\draw[thick,-{Stealth[length=8pt,width=7pt]}] (-0.025,-1.709) -- (0.395,-1.691);
\node at (0.207,-2.119) {$+$};
\draw[thick,-{Stealth[length=8pt,width=7pt]}] (2.006,-0.606) -- (1.783,-0.961);
\node at (2.291,-0.928) {$-$};
\draw[thick,-{Stealth[length=8pt,width=7pt]}] (1.884,0.727) -- (1.603,1.038);
\node at (2.119,1.076) {$+$};
\draw[thick,-{Stealth[length=8pt,width=7pt]}] (-0.300,1.591) -- (0.120,1.610);
\node at (-0.138,2.018) {$-$};
\draw[thick,-{Stealth[length=8pt,width=7pt]}] (-1.878,0.504) -- (-1.662,0.864);
\node at (-2.165,0.830) {$-$};
\node at (-1.681,-0.161) {$\uparrow$};
\fill (-2.300,-0.200) circle (1.6pt);
\node at (-0.674,-1.223) {{\scriptsize$\leftarrow\!/\!\rightarrow$}};
\fill (-1.000,-1.750) circle (1.6pt);
\node at (0.959,-1.169) {$\downarrow$};
\fill (1.350,-1.650) circle (1.6pt);
\node at (1.831,0.061) {$\uparrow$};
\fill (2.450,0.100) circle (1.6pt);
\node at (0.736,1.116) {$\downarrow$};
\fill (1.050,1.650) circle (1.6pt);
\node at (-0.859,1.069) {{\scriptsize$\leftarrow\!/\!\rightarrow$}};
\fill (-1.250,1.550) circle (1.6pt);
\draw[densely dashed,->] (0.360,-0.587) arc (-60:240:0.620);
\node at (0.050,-0.050) {$f$};
\end{tikzpicture}
\caption{The boundary of a face $f$, followed counterclockwise; the over-strands and the under-strands are not distinguished. An edge is marked $+$ if it is traveled in its own direction and $-$ otherwise. At a corner of type $\uparrow$, where both edges leave the crossing, the sign changes from $-$ to $+$; at a corner of type $\downarrow$, where both edges enter it, it changes from $+$ to $-$; and at the corners of the types $\leftarrow$ and $\rightarrow$ it does not change. So the corners of the types $\uparrow$ and $\downarrow$ alternate along the boundary, and a face has as many corners of each of these two types (proofs of Theorem~\ref{thm.kashaev} and of Proposition~\ref{prop.omega}).}
\label{fig.facewalk}
\end{figure}

\begin{proof}
Let $\Xi=\sum_j\ve_ja_j(t^{\ve_j})b_j(t^{-\ve_j})^{\T}$. For faces $f,f'$, the entry $(f,f')$ of $\WW_t^{-1}\Xi\WW_t$ is $(-t^{1/2})^{w(f')-w(f)}\Xi_{ff'}$, and $\Xi_{ff'}$ is the sum over the crossings $j$ of $\ve_j(a_j)_f(b_j)_{f'}$; by \eqref{eq.facetable}, $(a_j)_f$ is the sum of the entries of $a_j$ at the corners of $j$ which belong to the face $f$, and likewise for $(b_j)_{f'}$. So $\WW_t^{-1}\Xi\WW_t=\sum_jc_j$, where $c_j$ is zero outside the four faces of $j$ and is given on them, in the basis $(f^\uparrow_j,f^\leftarrow_j,f^\downarrow_j,f^\rightarrow_j)$, by the \emph{corner matrix} with entries $\ve_j(-t^{1/2})^{w(p')-w(p)}(a_j)_p(b_j)_{p'}$, $p,p'\in\{\uparrow,\leftarrow,\downarrow,\rightarrow\}$, the windings $w(p)$ being given by \eqref{eq.wcrossing}. From \eqref{eq.facetable} one computes
\begin{equation}\label{eq.cornermatrix}
\begin{aligned}
   c_j&=\begin{pmatrix} t & t^{1/2} & 1 & t^{1/2}\\ t^{1/2} & 1 & t^{-1/2} & 1\\ 1 & t^{-1/2} & t^{-1} & t^{-1/2}\\ t^{1/2} & 1 & t^{-1/2} & 1\end{pmatrix}&&\text{if $\ve_j=1$},\\[2pt]
   c_j&=-\begin{pmatrix} t^{-1} & t^{-1/2} & 1 & t^{-1/2}\\ t^{-1/2} & 1 & t^{1/2} & 1\\ 1 & t^{1/2} & t & t^{1/2}\\ t^{-1/2} & 1 & t^{1/2} & 1\end{pmatrix}&&\text{if $\ve_j=-1$}.
\end{aligned}
\end{equation}
For instance, for $\ve_j=1$ the entry $(\uparrow,\leftarrow)$ is $(-t)(t^{-1})(-t^{1/2})^{w(\leftarrow)-w(\uparrow)}=(-1)(-t^{1/2})=t^{1/2}$, and the entry $(\leftarrow,\rightarrow)$ is $t\cdot1\cdot(-t^{1/2})^{-2}=1$; for $\ve_j=-1$ the entry $(\downarrow,\downarrow)$ is $-(1)(t)=-t$. In both cases $c_j=\ve_j\,u_ju_j^{\T}$, where $u_j$ is the column vector with the coordinates $(t^{\ve_j/2},\,1,\,t^{-\ve_j/2},\,1)$ in this basis; this makes the sixteen entries easy to check, and it shows that the corner matrix is symmetric. The corner matrix $c_j$ is Kashaev's block $\ve_j\tau_j$ in which the entries $2\zeta^2-1$ at $\uparrow$ and $\downarrow$ have been replaced by $t^{\ve_j}$ and $t^{-\ve_j}$, and the entries $\zeta$ adjacent to $\uparrow$, resp.\ to $\downarrow$, by $t^{\ve_j/2}$, resp.\ $t^{-\ve_j/2}$. With $\vartheta=\frac12(t^{1/2}-t^{-1/2})$, so that $t^{\pm1/2}=\zeta \pm\vartheta$ and $t^{\pm1}=(2\zeta^2-1)\pm2\zeta \vartheta$, this reads
\begin{equation}\label{eq.cornerdiff}
\begin{aligned}
   c_j&=\ve_j\tau_j(t)+\vartheta\,\Delta_j,\\
   \Delta_j&=(f^\uparrow_j-f^\downarrow_j)(f^\leftarrow_j+f^\rightarrow_j)^{\T}+(f^\leftarrow_j+f^\rightarrow_j)(f^\uparrow_j-f^\downarrow_j)^{\T}+2\zeta\bigl(f^\uparrow_jf^{\uparrow\T}_j-f^\downarrow_jf^{\downarrow\T}_j\bigr),
\end{aligned}
\end{equation}
for both signs of the crossing (for $\ve_j=-1$, replacing $t$ by $t^{-1}$ changes the sign of $\vartheta$ and of $c_j$ and fixes $\tau_j$ and $\Delta_j$). It remains to show that $\sum_j\Delta_j=0$.

The last term of $\Delta_j$: the coefficient of $e_fe_f^{\T}$ in $\sum_j\bigl(f^\uparrow_jf^{\uparrow\T}_j-f^\downarrow_jf^{\downarrow\T}_j\bigr)$ is the number of corners of $f$ at which $f$ is the face $f^\uparrow$ (both edges of the corner leave the crossing) minus the number of corners at which it is $f^\downarrow$ (both edges enter the crossing). Walk once around the boundary of $f$ and record, for each boundary edge, whether it is traversed along ($+$) or against ($-$) its orientation. At a corner of type $\uparrow$ the walk arrives against an edge and leaves along the other ($-+$), at a corner of type $\downarrow$ it changes from $+$ to $-$, and at the corners of type $\leftarrow,\rightarrow$ the sign does not change. Along a closed walk the changes $-+$ and $+-$ alternate (Figure~\ref{fig.facewalk}), so the two numbers are equal and the last terms of the $\Delta_j$ sum to $0$.

The first two terms of $\Delta_j$: $(f^\uparrow_j)(f^\leftarrow_j+f^\rightarrow_j)^{\T}$ is the sum of $f^\uparrow_jf^{\leftarrow\T}_j$ and $f^\uparrow_jf^{\rightarrow\T}_j$, the pairs of faces adjacent along the two outgoing edges of $j$, and $f^\downarrow_j(f^\leftarrow_j+f^\rightarrow_j)^{\T}$ is the sum over the two incoming edges. For an edge $e$ put $s_e=e_{f_{\rm L}(e)}e_{f_{\rm R}(e)}^{\T}+e_{f_{\rm R}(e)}e_{f_{\rm L}(e)}^{\T}$. The sum of the first two terms of $\Delta_j$ is $\sum_{e\text{ leaves }j}s_e-\sum_{e\text{ enters }j}s_e$: for the outgoing left edge $e$ of $j$, which runs between $f^\leftarrow_j$ and $f^\uparrow_j$, $f^\uparrow_jf^{\leftarrow\T}_j+f^\leftarrow_jf^{\uparrow\T}_j=s_e$; in the same way the outgoing right edge, which runs between $f^\uparrow_j$ and $f^\rightarrow_j$, gives the two terms which contain $f^\uparrow_j$ and $f^\rightarrow_j$, and the two incoming edges, which run between $f^\downarrow_j$ and $f^\leftarrow_j$ and between $f^\downarrow_j$ and $f^\rightarrow_j$, give the four terms which contain $f^\downarrow_j$, with the sign $-$. Summing over $j$, every edge leaves one crossing and enters one, so $\sum_j\Delta_j=\sum_e(s_e-s_e)=0$.
\end{proof}

\begin{corollary}\label{cor.kashaev}
Let $\hate$ be an edge of $D$ and let $\tilde\tau_D(t)$ be Kashaev's matrix with the rows and columns of the two faces adjacent to $\hate$ removed.
\begin{enumerate}\renewcommand{\theenumi}{\roman{enumi}}
\item Let $\ve_D=\prod_j\ve_j$ be the product of the signs of the crossings and $\Delta_L(t)=\nabla_L(t^{1/2}-t^{-1/2})$ the Conway-normalized Alexander polynomial. Then
\[
   \det\tilde\tau_D(t)=\ve_D\,(-1)^{\mu-1}\,P_{D,\hate}(t^{\ve})\,\tP_{D,\hate}(t^{-\ve})=\ve_D\,\Delta_L(t)^2
\]
\textup(up to a unit $\pm t^m$ this is \cite[Theorem~1]{CimasoniFerretti24}, \cite[Corollary~3.3]{Liu25}\textup).
\item The kernel of $\tau_D(t)$ over $\BQ(t^{1/2})$ contains the two vectors $v_\pm=\bigl((-t^{\pm1/2})^{w(f)}\bigr)_{f\in\CF}$, which are the images under $\WW_t$ of the vectors $\1$ and $(t^{-w(f)})_{f\in\CF}$ of $\ker Q_D(t^{\ve})$; they span the kernel if $\Delta_L\ne0$ \textup(\cite[Corollary~3.8]{Liu25}\textup).
\end{enumerate}
\end{corollary}
\begin{proof}
(i) Conjugation by the diagonal matrix $\WW_t$ commutes with the removal of rows and columns of the same index, and removing the two columns $f_{\rm L}(\hate),f_{\rm R}(\hate)$ from $Q_D$ and $\tQ_D$ in the product of Theorem~\ref{thm.kashaev} gives $\tilde\tau_D(t)=\WW_t'^{-1}\,Q_{D,\hate}(t^{\ve})^{\T}\DD_\ve\tQ_{D,\hate}(t^{-\ve})\,\WW'_t$ with $\WW'_t$ the restricted diagonal matrix. Taking determinants, $\det\tilde\tau_D(t)=\ve_D\det Q_{D,\hate}(t^{\ve})\det\tQ_{D,\hate}(t^{-\ve})$, and Theorem~\ref{thm.faceform} gives the first equality. Lemma~\ref{lem.alexpoly} gives the second one: since $\Delta_L(t^{-1})=(-1)^{\mu-1}\Delta_L(t)$,
\begin{align*}
   P_{D,\hate}(t^{\ve})&=(-1)^{\mu-1}\,t^{m_{D,\hate}/2}\,\Delta_L(t),\\
   \tP_{D,\hate}(t^{-\ve})&=(-1)^{\mu-1}\,t^{-m_{D,\hate}/2}\,\Delta_L(t^{-1})=t^{-m_{D,\hate}/2}\,\Delta_L(t).
\end{align*}
(ii) By \eqref{eq.facetable} and \eqref{eq.wcrossing}, $a_j(t^{\ve_j})$ vanishes on $\1$ (each $m(e)$ does) and on $(t^{-w(f)})_f$: for $\ve_j=1$ and $k=w(f^\leftarrow_j)$ this is $-t\cdot t^{1-k}+t\cdot t^{-k}-t^{1-k}+t^{2-k}=0$, and for $\ve_j=-1$ it is $t^{-1}t^{1-k}-t^{-k}+t^{1-k}-t^{-1}t^{2-k}=0$. Since $\tau_D(t)$ is symmetric, its kernel equals its left kernel, which by Theorem~\ref{thm.kashaev} contains $\WW_t\ker Q_D(t^{\ve})$: if $Q_D(t^{\ve})k=0$, then
\[
   (\WW_tk)^{\T}\tau_D(t)=k^{\T}Q_D(t^{\ve})^{\T}\DD_\ve\tQ_D(t^{-\ve})\WW_t=0 .
\] If $\Delta_L\ne0$, then $Q_{D,\hate}(t^{\ve})$ is invertible by (i), so $\dim\ker Q_D(t^{\ve})=2$ and $\dim\ker\tau_D(t)=2$ as well.
\end{proof}

\begin{remark}\label{rem.kashaev}
(i) Theorem~\ref{thm.kashaev} is a two-flavor form of the factorization $\tau_D=K(D)^{\T}SK(D)$ of \cite[Proposition~3.1]{CimasoniFerretti24}, in which $K(D)$ is the crossing-by-face matrix of Kauffman's state model \cite{KauffmanFKT} and $S=\DD_\ve$. Indeed, by Lemma~\ref{lem.facelevel} the second factor $\tQ_D(t^{-\ve})$ is the first one rescaled by the levels and the windings, and since $\WW_t^{-1}=\diag(t^{-w(f)})\WW_t=\diag\bigl((-t^{-1/2})^{w(f)}\bigr)$, the theorem can be rewritten as the Gram matrix of a single matrix,
\[
   \tau_D(t)=K(t)^{\T}\,\DD_\ve\,K(t),\qquad K(t)=\diag\bigl(t^{-h(u_j)/2-(1+\ve_j)/4}\bigr)_j\ Q_D(t^{\ve})\ \diag\bigl((-t^{-1/2})^{w(f)}\bigr)_f ,
\]
Alexander's matrix with its rows and columns rescaled by monomials in $t^{1/2}$; this is the role played by Kauffman's matrix in \cite{CimasoniFerretti24}. The proof of the theorem explains the mechanism: the product $\ve_ja_jb_j^{\T}$ of the two relations of the crossing $j$ (under-strand, over-strand) is Kashaev's block up to a correction $\vartheta\Delta_j$, and the corrections cancel globally, by the alternation of the corners of a face and because every edge leaves one crossing and enters another.

(ii) The two null vectors of (ii) are, in the face picture, the right null vector $\1$ and the left null vector $\ell_t$ of the Alexander matrix of Lemma~\ref{lem.level}: by \eqref{eq.wh} the potential $k=(t^{-w(f)})_f$ satisfies $k_{f_{\rm L}(e)}-k_{f_{\rm R}(e)}=t^{h(e)}(1-t)$, i.e.\ $mk=(1-t)\ell_t^{\T}$. The grading of the faces by the winding number, which Liu introduces in order to identify the Schur complement of her Seifert form with $\tau_D$ \cite[proof of Theorem~3.1]{Liu25}, is thus the grading by the levels of Section~\ref{ssec.regions}: a face $f$ lies in a single region $\Omega$ of the complement of the Seifert circles, and $w(f)=-H_\Omega-1$ by \eqref{eq.wh} and Lemma~\ref{lem.regions}.
\end{remark}

\subsection{The deformed Kashaev matrix}\label{ssec.deformedkashaev}
Theorem~\ref{thm.kashaev} invites the following definition. In it the crossing variables are no longer specialized, and the determinant formula becomes a statement about one matrix which is indexed by the faces (Corollary~\ref{cor.deformedkashaev}).

\begin{definition}\label{def.deformedkashaev}
The \emph{deformed Kashaev matrix} of $D$ is the $\CF\times\CF$ matrix with entries in $\BZ[x,y]$
\[
   \tau_D(x,y)=Q_D(x)^{\T}\,\DD_\ve\,\tQ_D(y)=\sum_{j=1}^c\ve_j\,a_j(x_j)\,b_j(y_j)^{\T},
\]
and $\tau_{D,\hate}(x,y)=Q_{D,\hate}(x)^{\T}\DD_\ve\tQ_{D,\hate}(y)$ is the matrix with the rows and columns of the two faces adjacent to $\hate$ removed.
\end{definition}
In the notation of Remark~\ref{rem.flows}~(iii), the entry of $\tau_D(x,y)$ at two faces $f,f'$ is the value $\sum_j\ve_j\,\pi_j(\partial f)\,\tpi_j(\partial f')$ of the coupled form on the boundaries of the two faces, and $\tau_{D,\hate}(x,y)$ is the matrix of this form in the basis of the faces of the lattice $\FF_{D,\hate}$.

The matrix $\tau_D(x,y)$ is the sum over the crossings of matrices of rank one, each supported on the four faces of its crossing and bilinear in the two variables of the crossing: by \eqref{eq.facetable}, in the basis $(f^\uparrow_j,f^\leftarrow_j,f^\downarrow_j,f^\rightarrow_j)$,
\begin{equation}\label{eq.localblocks}
\begin{aligned}
   \ve_ja_jb_j^{\T}&=\begin{pmatrix} x_j & -x_jy_j & x_jy_j & -x_j\\ -x_j & x_jy_j & -x_jy_j & x_j\\ 1 & -y_j & y_j & -1\\ -1 & y_j & -y_j & 1\end{pmatrix}&&(\ve_j=1),\\[2pt]
   \ve_ja_jb_j^{\T}&=\begin{pmatrix} -x_j & x_j & -x_jy_j & x_jy_j\\ 1 & -1 & y_j & -y_j\\ -1 & 1 & -y_j & y_j\\ x_j & -x_j & x_jy_j & -x_jy_j\end{pmatrix}&&(\ve_j=-1),
\end{aligned}
\end{equation}
the second block being the first one with $\leftarrow$ and $\rightarrow$ exchanged and the sign changed. The variable $x_j$ sits in the rows $f^\uparrow_j,f^\leftarrow_j$ for a positive crossing, i.e.\ in the rows of the faces adjacent to the outgoing under-edge $u'_j$, and $y_j$ in the columns of the faces adjacent to the incoming over-edge $o_j$.

\begin{corollary}[Kashaev form of the determinant formula]\label{cor.deformedkashaev}
For every diagram $D$ of $L$ and every edge $\hate$ of $D$,
\[
   \LG(L;t_0,t_1)=\ve_D\ \BE_D\bigl[\det\tau_{D,\hate}(x,y)\bigr],\qquad\text{where}\quad \ve_D=\prod_{j=1}^c\ve_j .
\]
At the free fermion point $x=t^{\ve}$, $y=t^{-\ve}$, the deformed Kashaev matrix is Kashaev's matrix up to the windings: $\tau_D(t^{\ve},t^{-\ve})=\WW_t\,\tau_D(t)\,\WW_t^{-1}$.
\end{corollary}
\begin{proof}
By Theorem~\ref{thm.faceform},
\[
   \det\tau_{D,\hate}(x,y)=\ve_D\det Q_{D,\hate}(x)\det\tQ_{D,\hate}(y)=\ve_D\,(-1)^{\mu-1}P_{D,\hate}(x)\,\tP_{D,\hate}(y),
\]
and Theorem~\ref{thm.HK} applies. The second statement is Theorem~\ref{thm.kashaev}.
\end{proof}

\begin{example}[the trefoil]\label{ex.trefoilkashaev}
Let $D$ be the standard diagram of the trefoil of Example~\ref{ex.trefoilflows}. Its five faces (Figure~\ref{fig.trefoilfaces}) are the unbounded face $f_0$ (outside $C_1$), the face $f_\infty$ inside $C_2$, and the three bigons $g_1,g_2,g_3$ between the strands, $g_j$ lying between the crossings $j$ and $j+1$ (indices modulo $3$), so that $g_j$ is bounded by the edges $u'_j=o_{j+1}$ and $o'_j=u_{j+1}$; the windings are $w(f_0)=0$, $w(g_j)=-1$, $w(f_\infty)=-2$. At the crossing $j$ (all crossings are positive) $f^\uparrow_j=g_j$, $f^\downarrow_j=g_{j-1}$, $f^\leftarrow_j=f_0$ and $f^\rightarrow_j=f_\infty$, so that $a_j=x_j(f_0-g_j)+(f_\infty-g_{j-1})$ and $b_j=y_j(f_0-g_{j-1})+(f_\infty-g_j)$, and in the basis $(f_0,g_1,g_2,g_3,f_\infty)$
\[
   \tau_D(x,y)=\begin{pmatrix}
   x_1y_1+x_2y_2+x_3y_3 & -x_1-x_2y_2 & -x_2-x_3y_3 & -x_3-x_1y_1 & x_1+x_2+x_3\\
   -x_1y_1-y_2 & x_1+y_2 & 1 & x_1y_1 & -x_1-1\\
   -x_2y_2-y_3 & x_2y_2 & x_2+y_3 & 1 & -x_2-1\\
   -x_3y_3-y_1 & 1 & x_3y_3 & x_3+y_1 & -x_3-1\\
   y_1+y_2+y_3 & -y_2-1 & -y_3-1 & -y_1-1 & 3
   \end{pmatrix}.
\]
At the free fermion point $x_j=t$, $y_j=t^{-1}$, the entries become $3$, $-t-1$, $3t$, $-1-t^{-1}$, $t+t^{-1}$, $1$ and $3t^{-1}$, and the conjugation by $\WW_t=\diag(1,-t^{-1/2},-t^{-1/2},-t^{-1/2},t^{-1})$ turns $-t-1$ and $-1-t^{-1}$ into $t^{1/2}+t^{-1/2}=2\zeta$ and $3t,3t^{-1}$ into $3$: Kashaev's matrix of the trefoil has diagonal entries $3,\,2\zeta^2-1+2\zeta^2-1,\dots,3$, entries $2\zeta$ between $f_0$ or $f_\infty$ and a bigon (adjacent along one edge, i.e.\ at two corners), $1$ between two bigons and $3$ between $f_0$ and $f_\infty$ (opposite at three crossings). For the cut edge $\hate=u'_3=o_1$, the faces $f_0$ and $g_3$ are removed, and $\det\tau_{D,\hate}(x,y)$ is a polynomial in six variables of degree one in each, whose image under $\BE_D$ is $\LG(3_1)$.
\end{example}

\begin{remark}\label{rem.deformedkashaev}
(i) \emph{What the deformation does.} Kashaev's matrix is symmetric, with an entry of $\tau_j$ for each pair of corners of a crossing. The deformed matrix is a sum of rank-one blocks $\ve_ja_jb_j^{\T}$ in which the rows belong to the first flavor and the columns to the second one: the symmetry is broken, and it is restored at the free fermion point only after the conjugation by the windings (Lemma~\ref{lem.facelevel}). The coupling $\BE_D$ acts on the variables of one crossing at a time, i.e.\ on one block at a time, exactly as it acts on the deformed Seifert matrices of Corollary~\ref{cor.seifertformula}. By Theorem~\ref{thm.mixed} this action can be carried out block by block, before the determinant is taken: it replaces the block $\ve_ja_jb_j^{\T}$ of rank one by the matrix of the coupled crossing form $\Psi_j$ in the basis of the faces, which has rank at most two and entries in $\BZ[t_0^{\pm1},t_1^{\pm1}]$, and it replaces the determinant by the mixed discriminant of the $c$ blocks; the sign $\ve_D$ disappears. In Liu's proof of Kashaev's conjecture \cite{Liu25}, $\tau_D$ is the Schur complement, over the classes of the vertices, of the Seifert form $Q=(1-\omega)A+(1-\bar\omega)A^{\T}$ of a Seifert surface for $L\sqcup rL$ whose first homology is generated by the faces and the vertices of $D$ (and the vertex block of $Q$ is diagonal). The two variables of a crossing thus sit exactly on the generator which Liu eliminates. The block matrix $\bigl(\begin{smallmatrix}0&Q_D(x)^{\T}\\ \tQ_D(y)&-\DD_\ve\end{smallmatrix}\bigr)$, indexed by the faces and the crossings, has $\tau_D(x,y)$ as its Schur complement over the crossings and determinant $(-1)^{c+\mu-1}P_{D,\hate}(x)\tP_{D,\hate}(y)$ after the removal of the two faces of $\hate$; whether it is, in a suitable basis, a deformed Seifert matrix of Liu's surface in the sense of Section~\ref{sec.seifert}, with the crossing variables distributed over the vertex classes, is a natural question which we leave open.

(ii) \emph{Signatures.} The Hermitian structure of $\tau_D$ has no counterpart for $\LG$. The transpose of $\tau_D(x,y)$ exchanges the two flavors, which by Lemma~\ref{lem.inverse} amounts to passing to the inverse diagram, i.e.\ to the mirror image of the reversed link; correspondingly, the mirror image $L^*$ of $L$ changes $q$ into $q^{-1}$ and $\LG(L^*;t_0,t_1)=\LG(L;t_0^{-1},t_1^{-1})$, so that on the torus $|t_0|=|t_1|=1$ the complex conjugate of $\LG(L)$ is $\LG(L^*)$: the polynomial $\LG(L)$ is real on the torus only up to the mirror image (it is not real for the trefoil, for $5_2$ or for $6_1$, and it is real for the amphichiral knots $4_1$ and $6_3$). Moreover $\BE_D$ is not an evaluation. It can be carried inside the determinant, but what one obtains at a point $(t_0,t_1)$ is not one matrix of which a signature could be taken: it is the pencil $\sum_j\lambda_j\Psi_j$ of the coupled crossing forms, of which $\LG(L)$ is the mixed discriminant (Theorem~\ref{thm.mixed}). The matrix of the form $\sum_j\Psi_j$ is not Hermitian on the torus in general (it is not for the trefoil, by Example~\ref{ex.trefoilmixed}), and its determinant is not an invariant of the link (Proposition~\ref{prop.kink} and Example~\ref{ex.pencilunknot}). The real quantity which $\LG$ does provide on the torus is $\LG(L)\overline{\LG(L)}=\LG(L)\LG(L^*)$, the Links--Gould polynomial of the connected sum of $L$ with its mirror image. On alternating diagrams the signature $\sigma_L(-1)$ does appear in the state sum for $\LG$, by a different route: see Proposition~\ref{prop.signature} and Corollary~\ref{cor.extremalstates}.

(iii) \emph{The Kinoshita--Terasaka knot again.} By Corollary~\ref{cor.kashaev}~(i), $\det\tilde\tau_D=\ve_D=\pm1$ for every diagram of the Kinoshita--Terasaka knot, since its Alexander polynomial is $1$; yet the coupling $\BE_D[\det\tau_{D,\hate}(x,y)]$ of the deformed matrix is not $1$, by Remark~\ref{rem.KT}. The crossing variables of the deformed Kashaev matrix carry information which the Seifert form of Liu's surface, like every Seifert form, does not.
\end{remark}

\section{The point \texorpdfstring{$x=-1$}{x = -1}: Goeritz matrices and the forms of Gordon and Litherland}\label{sec.minusone}

At the Alexander point $x_j=t^{\ve_j}$ the crossing variables depend on the orientation of the link, through the signs of the crossings. For exactly one value of $t$ they do not: at $t=-1$ all the crossing variables are equal to $-1$. At this point the crossing forms are
\begin{equation}\label{eq.piminusone}
   \pi_j\big|_{x_j=-1}=-U_j,\qquad \tpi_j\big|_{y_j=-1}=-O_j ,
\end{equation}
where $U_j(z)=z_{u'_j}+z_{u_j}$ and $O_j(z)=z_{o_j}+z_{o'_j}$ are the linear forms of \eqref{eq.UO}, and the edge polynomials are, up to sign, the determinant $|\Delta_L(-1)|$ of the link (Lemma~\ref{lem.alexpoly}). The oriented smoothing is no longer distinguished: the forms $U_j$ and $O_j$ are related in the same way to the transfers of all the $2^c$ smoothings of the diagram, and each smoothing $\vs$ comes with a spanning surface of the link, orientable or not, its state surface $F_\vs$ (Section~\ref{sec.statesurfaces}).

In this section we describe the lattice of flows at $x=-1$. The objects which appear are classical: the matrices of Goeritz \cite{Goeritz33}, the forms of Gordon and Litherland of the spanning surfaces of a link \cite{GordonLitherland}, and the value at $t=-1$ of Kashaev's matrix, which was studied by Cimasoni and Ferretti \cite{CimasoniFerretti24}. The lattice of flows puts them in one frame, in which the statements are uniform and the proofs are short. The results are the following.
\begin{itemize}
\item Each crossing has two transfers $w^0_j$ and $w^1_j$, one for each of its two smoothings, and all the forms of this section are written in terms of them. For every smoothing $\vs$, the forms of Gordon and Litherland of the two state surfaces $F_\vs$ and $F^\uparrow_\vs$, below and above the plane of the diagram, are $\Goe_\vs+\Stk_\vs$ and $\Goe_\vs-\Stk_\vs$, where $\Goe_\vs$ is a Goeritz form and $\Stk_\vs$ records the nesting of the circles of the smoothing (Proposition~\ref{prop.formsminusone}); the two are exchanged by an involution of the homology of the surface (Proposition~\ref{prop.J}).
\item At $x=y=-1$ the coupled form of the two flavors is symmetric. In the basis of the faces it is the direct sum of the two Goeritz matrices of the diagram, and Kashaev's matrix at $t=-1$ is congruent to it by a diagonal matrix of signs, for every diagram (Theorem~\ref{thm.goeritz}).
\item The Gram matrices of the forms of Gordon and Litherland of all the state surfaces present one and the same group, the quotient of $\BZ^c$ by the image of the lattice of flows under $(U_1,\dots,U_c)$, whose order is the determinant of the link (Theorem~\ref{thm.determinant}).
\item The signature of the link is the signature of the form of $F_\vs$ minus $\sum_{\vs_j=1}\ve_j$, for every smoothing $\vs$ (Theorem~\ref{thm.signature}). This is the formula of Gordon and Litherland for the state surfaces. The proof given here is elementary: changing the smoothing of one crossing $j$ adds to the form the square $\ve_j\,U_j\otimes U_j$ of a linear form, and changes its rank by exactly one. For the checkerboard smoothings one obtains the formula of Gordon and Litherland for the Goeritz matrices, and from it the theorem of Cimasoni and Ferretti on the signature of Kashaev's matrix at $t=-1$ (Corollaries~\ref{cor.GLgoeritz} and~\ref{cor.kashaevsignature}); for alternating diagrams one obtains Traczyk's formula for the signature (Corollary~\ref{cor.altdefinite}).
\end{itemize}

Throughout this section $D$ is a diagram of a link $L$ with $c$ crossings, drawn in the plane $P$ as in Section~\ref{ssec.statesurface}, and $\hate$ is an edge of $D$. We use the lattice of flows $\FF_D$ and the potentials of the flows (Lemma~\ref{lem.faces}), the smoothings $\vs\in\{0,1\}^c$, their circles $z^C$ and their transfers $w^\vs_j$ (Section~\ref{ssec.smoothings}), the covered and the free quadrants of a crossing, the state surfaces $F_\vs$ and $F^\uparrow_\vs$ and the band signs $\theta_j=\ve_j(-1)^{\vs_j}$ (Section~\ref{sec.statesurfaces}). For a symmetric bilinear form $Q$ on a free abelian group of finite rank, or for a symmetric matrix, $\sign Q$ denotes its signature, the number of positive squares minus the number of negative squares in a diagonalization of $Q$ over $\BR$, and $\rk Q$ its rank; they do not change when $Q$ is regarded as a form on the quotient by a subgroup of its radical.

\subsection{The two transfers of a crossing}\label{ssec.twotransfers}
A crossing has two smoothings, hence two transfers. In this subsection we express the two forms $U_j$ and $O_j$, to which the crossing forms reduce at $x=-1$, as combinations of these two transfers, and we describe the transfers in terms of the potentials of the flows (Lemma~\ref{lem.cornertransfers}). All the computations of this section rest on these formulas.

Let $j$ be a crossing of $D$. By Definition~\ref{def.transfer}, the transfers at $j$ of the two constant smoothings, the oriented smoothing $\vs=0=(0,\dots,0)$ and the smoothing $\vs=1=(1,\dots,1)$, are the linear forms on $\FF_D$
\begin{equation}\label{eq.w01}
   w^0_j(z)=z_{o_j}-z_{u'_j}=z_{o'_j}-z_{u_j},\qquad\qquad w^1_j(z)=-z_{o'_j}-z_{u'_j}=-z_{o_j}-z_{u_j},
\end{equation}
where the second expressions follow from the flow condition \eqref{eq.flow}. Thus $w^0_j$ is the form $w_j$ of the previous sections, and $-w^1_j(z)$ is the total flow of $z$ through the crossing $j$. The transfer of an arbitrary smoothing $\vs$ at the crossing $j$ is one of these two forms: $w^\vs_j=w^{\vs_j}_j$. We write
\[
   \bar\vs=1-\vs=(1-\vs_1,\dots,1-\vs_c)
\]
for the \emph{complementary smoothing}, which smooths every crossing in the other way. Its arcs lie in the free quadrants of $\vs$, its transfer $w^{\bar\vs}_j$ at $j$ is the other one of the two forms $w^0_j,w^1_j$, and its band signs are $\theta_j(\bar\vs)=-\theta_j(\vs)$. The two forms \eqref{eq.UO} are combinations of the two transfers: for every smoothing $\vs$,
\begin{equation}\label{eq.UOw}
   U_j=-w^0_j-w^1_j=-\bigl(w^\vs_j+w^{\bar\vs}_j\bigr),\qquad\qquad O_j=w^0_j-w^1_j=(-1)^{\vs_j}\bigl(w^\vs_j-w^{\bar\vs}_j\bigr).
\end{equation}
Indeed, by \eqref{eq.w01}, $-w^0_j(z)-w^1_j(z)=(z_{u'_j}-z_{o_j})+(z_{o_j}+z_{u_j})=U_j(z)$ and $w^0_j(z)-w^1_j(z)=(z_{o_j}-z_{u'_j})+(z_{o'_j}+z_{u'_j})=O_j(z)$; the expressions in terms of $\vs$ follow because $(w^\vs_j,w^{\bar\vs}_j)$ is $(w^0_j,w^1_j)$ if $\vs_j=0$ and $(w^1_j,w^0_j)$ if $\vs_j=1$. (Solved for $w^\vs_j$, the two identities are \eqref{eq.transferUO}.)

In terms of a potential, the two transfers of a crossing are the differences of the potential across the crossing in its two directions, and for a given smoothing they are distributed over the free and the covered quadrants as follows.

\begin{lemma}[transfers and potentials]\label{lem.cornertransfers}
Let $z\in\FF_D$ be a flow, let $k$ be a potential of $z$ (Lemma~\ref{lem.faces}), and let $j$ be a crossing.
\begin{enumerate}\renewcommand{\theenumi}{\roman{enumi}}
\item $w^0_j(z)=\ve_j\bigl(k(f^\uparrow_j)-k(f^\downarrow_j)\bigr)$ and $w^1_j(z)=k(f^\rightarrow_j)-k(f^\leftarrow_j)$.
\item Let $\vs$ be a smoothing and $\theta_j=\ve_j(-1)^{\vs_j}$. Each of the two half-strands $u'_j$ and $u_j$ of the under-strand of $j$ separates a covered quadrant from a free quadrant. Let $c'_j$ and $\phi'_j$ be the faces of the covered and of the free quadrant next to $u'_j$, and let $c_j$ and $\phi_j$ be the faces of the covered and of the free quadrant next to $u_j$ \textup(Figure~\ref{fig.cornersmoothing}\textup). Then
\[
   w^\vs_j(z)=\theta_j\bigl(k(\phi'_j)-k(\phi_j)\bigr)\qquad\text{and}\qquad w^{\bar\vs}_j(z)=-\theta_j\bigl(k(c'_j)-k(c_j)\bigr).
\]
\end{enumerate}
\end{lemma}
\begin{proof}
(i) We use \eqref{eq.potential} and \eqref{eq.mcrossing}, and we write $k^\uparrow$ for $k(f^\uparrow_j)$, and so on. If $\ve_j=1$, then $z_{o_j}=k^\leftarrow-k^\downarrow$, $z_{u'_j}=k^\leftarrow-k^\uparrow$ and $z_{o'_j}=k^\uparrow-k^\rightarrow$; hence $w^0_j(z)=z_{o_j}-z_{u'_j}=k^\uparrow-k^\downarrow$ and $w^1_j(z)=-z_{o'_j}-z_{u'_j}=k^\rightarrow-k^\leftarrow$. If $\ve_j=-1$, then $z_{o_j}=k^\downarrow-k^\rightarrow$, $z_{u'_j}=k^\uparrow-k^\rightarrow$ and $z_{o'_j}=k^\leftarrow-k^\uparrow$; hence $w^0_j(z)=k^\downarrow-k^\uparrow$ and $w^1_j(z)=k^\rightarrow-k^\leftarrow$.

(ii) By \eqref{eq.mcrossing}, the half-strand $u'_j$ separates the quadrants of $f^\leftarrow_j$ and $f^\uparrow_j$ if $\ve_j=1$, and those of $f^\uparrow_j$ and $f^\rightarrow_j$ if $\ve_j=-1$; the half-strand $u_j$ separates the quadrants of $f^\downarrow_j$ and $f^\rightarrow_j$ if $\ve_j=1$, and those of $f^\leftarrow_j$ and $f^\downarrow_j$ if $\ve_j=-1$. The covered quadrants are those of $f^\leftarrow_j$ and $f^\rightarrow_j$ if $\vs_j=0$, and those of $f^\uparrow_j$ and $f^\downarrow_j$ if $\vs_j=1$ (Section~\ref{ssec.statesurface}). So the four faces are given by the following table, in which the arrows stand for the faces $f^\leftarrow_j$, $f^\uparrow_j$, $f^\rightarrow_j$, $f^\downarrow_j$:
\begin{equation}\label{eq.cornertable}
\begin{array}{cc|c|cccc}
   \ve_j&\vs_j&\theta_j&c'_j&\phi'_j&c_j&\phi_j\\ \hline
   1&0&1&\leftarrow&\uparrow&\rightarrow&\downarrow\\
   1&1&-1&\uparrow&\leftarrow&\downarrow&\rightarrow\\
   -1&0&-1&\rightarrow&\uparrow&\leftarrow&\downarrow\\
   -1&1&1&\uparrow&\rightarrow&\downarrow&\leftarrow
\end{array}
\end{equation}
In the first line, $w^\vs_j=w^0_j$ and $w^{\bar\vs}_j=w^1_j$, and by (i) $w^0_j(z)=k^\uparrow-k^\downarrow=\theta_j\bigl(k(\phi'_j)-k(\phi_j)\bigr)$ and $w^1_j(z)=k^\rightarrow-k^\leftarrow=-\theta_j\bigl(k(c'_j)-k(c_j)\bigr)$. In the second line, $w^\vs_j=w^1_j$ and $w^{\bar\vs}_j=w^0_j$, and $w^1_j(z)=k^\rightarrow-k^\leftarrow=-\bigl(k(\phi'_j)-k(\phi_j)\bigr)$, $w^0_j(z)=k^\uparrow-k^\downarrow=k(c'_j)-k(c_j)$. In the third line, $w^0_j(z)=k^\downarrow-k^\uparrow=-\bigl(k(\phi'_j)-k(\phi_j)\bigr)$ and $w^1_j(z)=k^\rightarrow-k^\leftarrow=k(c'_j)-k(c_j)$. In the fourth line, $w^1_j(z)=k^\rightarrow-k^\leftarrow=k(\phi'_j)-k(\phi_j)$ and $w^0_j(z)=k^\downarrow-k^\uparrow=-\bigl(k(c'_j)-k(c_j)\bigr)$. In all four cases these are the two formulas of the lemma.
\end{proof}

\begin{figure}[ht]
\centering
\begin{tikzpicture}[>=Stealth,scale=0.82,every node/.style={font=\small}]
\def\covLR{\fill[gray!25] (0,0) -- (135:1.25) arc (135:225:1.25) -- cycle; \fill[gray!25] (0,0) -- (-45:1.25) arc (-45:45:1.25) -- cycle;
   \draw[densely dashed] (135:1.25) arc (135:225:1.25); \draw[densely dashed] (-45:1.25) arc (-45:45:1.25);}
\def\covUD{\fill[gray!25] (0,0) -- (45:1.25) arc (45:135:1.25) -- cycle; \fill[gray!25] (0,0) -- (225:1.25) arc (225:315:1.25) -- cycle;
   \draw[densely dashed] (45:1.25) arc (45:135:1.25); \draw[densely dashed] (225:1.25) arc (225:315:1.25);}
\def\pos{\draw (1.45,-1.45) -- (0.2,-0.2); \draw[->] (-0.2,0.2) -- (-1.45,1.45); \draw[->] (-1.45,-1.45) -- (1.45,1.45);
   \node[font=\scriptsize] at (-1.6,-1.7) {$o_j$}; \node[font=\scriptsize] at (1.6,1.75) {$o'_j$};
   \node[font=\scriptsize] at (1.6,-1.7) {$u_j$}; \node[font=\scriptsize] at (-1.6,1.75) {$u'_j$};}
\def\negc{\draw (-1.45,-1.45) -- (-0.2,-0.2); \draw[->] (0.2,0.2) -- (1.45,1.45); \draw[->] (1.45,-1.45) -- (-1.45,1.45);
   \node[font=\scriptsize] at (1.6,-1.7) {$o_j$}; \node[font=\scriptsize] at (-1.6,1.75) {$o'_j$};
   \node[font=\scriptsize] at (-1.6,-1.7) {$u_j$}; \node[font=\scriptsize] at (1.6,1.75) {$u'_j$};}
\begin{scope}
\covLR \pos
\node at (-0.8,0) {$c'_j$}; \node at (0,0.8) {$\phi'_j$}; \node at (0.8,0) {$c_j$}; \node at (0,-0.8) {$\phi_j$};
\node at (0,-2.5) {$\ve_j=1,\ \vs_j=0$}; \node at (0,-3.1) {$\theta_j=1$};
\end{scope}
\begin{scope}[xshift=4.75cm]
\covUD \pos
\node at (0,0.8) {$c'_j$}; \node at (-0.8,0) {$\phi'_j$}; \node at (0,-0.8) {$c_j$}; \node at (0.8,0) {$\phi_j$};
\node at (0,-2.5) {$\ve_j=1,\ \vs_j=1$}; \node at (0,-3.1) {$\theta_j=-1$};
\end{scope}
\begin{scope}[xshift=9.5cm]
\covLR \negc
\node at (0.8,0) {$c'_j$}; \node at (0,0.8) {$\phi'_j$}; \node at (-0.8,0) {$c_j$}; \node at (0,-0.8) {$\phi_j$};
\node at (0,-2.5) {$\ve_j=-1,\ \vs_j=0$}; \node at (0,-3.1) {$\theta_j=-1$};
\end{scope}
\begin{scope}[xshift=14.25cm]
\covUD \negc
\node at (0,0.8) {$c'_j$}; \node at (0.8,0) {$\phi'_j$}; \node at (0,-0.8) {$c_j$}; \node at (-0.8,0) {$\phi_j$};
\node at (0,-2.5) {$\ve_j=-1,\ \vs_j=1$}; \node at (0,-3.1) {$\theta_j=1$};
\end{scope}
\end{tikzpicture}
\caption{The four corners of a crossing $j$ with respect to a smoothing $\vs$, in the four cases of the table \eqref{eq.cornertable}. The covered quadrants are shaded and the arcs of the smoothing are dashed; $c'_j$ and $\phi'_j$ are the faces of the covered and of the free quadrant next to the outgoing half-strand $u'_j$ of the under-strand, $c_j$ and $\phi_j$ those next to its incoming half-strand $u_j$. The covered quadrant next to $u'_j$ lies on the left of $u'_j$ exactly when $\theta_j=1$, as in \eqref{eq.thetaside}. The transfer of $\vs$ is, up to the sign $\theta_j$, the difference of the potential between the free corners $\phi'_j$ and $\phi_j$, and the transfer of the complementary smoothing is the difference between the covered corners $c'_j$ and $c_j$ (Lemma~\ref{lem.cornertransfers}).}
\label{fig.cornersmoothing}
\end{figure}

So the transfer of a smoothing at a crossing is, up to the sign $\theta_j$, the difference of the potential between the two \emph{free} corners of the crossing, and the transfer of the complementary smoothing is the difference between the two \emph{covered} corners. This agrees with Proposition~\ref{prop.statehomology}: a flow which is the boundary of a face $f$ runs through the band of the crossing $j$ of $F_\vs$ only if it passes from one strand of the crossing to the other one at a corner which the band does not cover.

\subsection{The Goeritz form, the stacking form and Kashaev's form}\label{ssec.formsminusone}
All the forms of this section are built from the products of the two transfers of the crossings. They are symmetric because of the following identity.

\begin{lemma}\label{lem.basicidentity}
For all flows $z,z'\in\FF_D$,
\[
   \sum_{j=1}^c\ve_j\Bigl(w^0_j(z)\,w^1_j(z')-w^1_j(z)\,w^0_j(z')\Bigr)=0 .
\]
\end{lemma}
\begin{proof}
Let $k$ and $k'$ be potentials of $z$ and $z'$. For two faces $p,q$ put $a(p,q)=k(p)k'(q)-k(q)k'(p)$, so that $a(q,p)=-a(p,q)$, and for an edge $e$ put $A(e)=a\bigl(f_{\rm L}(e),f_{\rm R}(e)\bigr)$. Fix a crossing $j$ and write $k^\uparrow=k(f^\uparrow_j)$, and so on. By Lemma~\ref{lem.cornertransfers}~(i),
\begin{align*}
   \ve_j\Bigl(w^0_j(z)\,w^1_j(z')-w^1_j(z)\,w^0_j(z')\Bigr)
   &=(k^\uparrow-k^\downarrow)(k'^\rightarrow-k'^\leftarrow)-(k^\rightarrow-k^\leftarrow)(k'^\uparrow-k'^\downarrow)\\
   &=a(\leftarrow,\uparrow)+a(\uparrow,\rightarrow)-a(\leftarrow,\downarrow)-a(\downarrow,\rightarrow),
\end{align*}
where the second line is obtained by expanding the two products and collecting the eight terms, and where the arrows stand for the four faces of the crossing. By \eqref{eq.mcrossing}, for both signs of the crossing, the two edges which leave $j$ have the faces $(f_{\rm L},f_{\rm R})=(f^\leftarrow_j,f^\uparrow_j)$ and $(f^\uparrow_j,f^\rightarrow_j)$ on their two sides, and the two edges which enter $j$ have the faces $(f^\leftarrow_j,f^\downarrow_j)$ and $(f^\downarrow_j,f^\rightarrow_j)$. So the term of the crossing $j$ is
\[
   \sum_{e\text{ leaves }j}A(e)-\sum_{e\text{ enters }j}A(e).
\]
Every edge of $D$ leaves exactly one crossing and enters exactly one crossing; so in the sum over $j$ every $A(e)$ occurs once with each sign, and the sum vanishes.
\end{proof}

The identity is a form of the fact that two closed curves on the sphere have intersection number zero: if $z$ and $z'$ pass straight through the crossing $j$, with multiplicity one, along its two different strands, the term of $j$ is $\pm2$, and it vanishes if they pass along the same strand.

\begin{definition}[the forms at $x=-1$]\label{def.formsminusone}
Let $\vs$ be a smoothing of $D$, with complementary smoothing $\bar\vs$ and band signs $\theta_j=\ve_j(-1)^{\vs_j}$. We define five bilinear forms on $\FF_D$ with integer values:
\begin{align*}
   \Goe_\vs&=-\sum_{j=1}^c\theta_j\,w^\vs_j\otimes w^\vs_j &&\text{(the \emph{Goeritz form} of the smoothing $\vs$),}\\
   \Stk_\vs&=-\sum_{j=1}^c\theta_j\,w^{\bar\vs}_j\otimes w^\vs_j &&\text{(the \emph{stacking form} of $\vs$),}\\
   \GLf_\vs&=\Goe_\vs+\Stk_\vs,\qquad \GLf^\uparrow_\vs=\Goe_\vs-\Stk_\vs,&&\\
   \Kash&=\sum_{j=1}^c\ve_j\,U_j\otimes O_j &&\text{(\emph{Kashaev's form}).}
\end{align*}
Here $(\lambda\otimes\lambda')(z,z')=\lambda(z)\lambda'(z')$ for two linear forms $\lambda,\lambda'$.
\end{definition}

The names are justified by the following proposition and by Theorem~\ref{thm.goeritz}.

\begin{proposition}[the forms at $x=-1$]\label{prop.formsminusone}
Let $\vs$ be a smoothing of $D$.
\begin{enumerate}\renewcommand{\theenumi}{\roman{enumi}}
\item \emph{The forms of Gordon and Litherland.} For all flows $z,z'$,
\[
   \GLf_\vs(z,z')=\sum_{j=1}^c\theta_j\,U_j(z)\,w^\vs_j(z')\qquad\text{and}\qquad \GLf^\uparrow_\vs(z,z')=-\sum_{j=1}^c\ve_j\,O_j(z)\,w^\vs_j(z') .
\]
Consequently $\GLf_\vs(z,z')=\CG_{F_\vs}([z],[z'])$ and $\GLf^\uparrow_\vs(z,z')=\CG_{F^\uparrow_\vs}([z],[z'])$: the forms $\GLf_\vs$ and $\GLf^\uparrow_\vs$ are the forms of Gordon and Litherland of the state surfaces of $\vs$ below and above the plane, evaluated on the classes of the flows.
\item \emph{Symmetry.} The five forms $\Goe_\vs$, $\Stk_\vs$, $\GLf_\vs$, $\GLf^\uparrow_\vs$ and $\Kash$ are symmetric.
\item \emph{Kashaev's form.} For every smoothing $\vs$,
\begin{equation}\label{eq.Kashw}
   \Kash=\sum_{j=1}^c\ve_j\bigl(w^1_j\otimes w^1_j-w^0_j\otimes w^0_j\bigr)=\Goe_\vs+\Goe_{\bar\vs}=\GLf_\vs+\GLf_{\bar\vs}=\GLf^\uparrow_\vs+\GLf^\uparrow_{\bar\vs},
\end{equation}
and $\Stk_{\bar\vs}=-\Stk_\vs$. Moreover $\Kash$ is the value at $x=y=-1$ of the coupled form $\tau=\sum_j\ve_j\,\pi_j\otimes\tpi_j$ of Remark~\ref{rem.flows}~\textup{(iii)}; on the boundaries of the faces it is the deformed Kashaev matrix $\tau_D(-1,-1)$ of Definition~\ref{def.deformedkashaev}.
\item \emph{Radicals.} The circles $z^C$ of $\vs$ lie in the radicals of $\Goe_\vs$, $\Stk_\vs$, $\GLf_\vs$ and $\GLf^\uparrow_\vs$; so these four forms are defined on $H_1(F_\vs;\BZ)=\FF_D/\langle z^C\rangle$. The alternating sum of the Seifert circles
\[
   z^0=\sum_C(-1)^{h(C)}\,C,
\]
where $h$ is the level of Definition~\ref{def.level}, satisfies $w^0_j(z^0)=w^1_j(z^0)=0$ for every crossing $j$; it lies in the radical of all the forms of Definition~\ref{def.formsminusone}, and it is a combination of the circles of every smoothing.
\item \emph{Changing one smoothing.} Let $j$ be a crossing with $\vs_j=0$ and let $\vs'$ be the smoothing which differs from $\vs$ at the crossing $j$ only, where $\vs'_j=1$. Then
\[
   \GLf_{\vs'}=\GLf_\vs+\ve_j\,U_j\otimes U_j\qquad\text{and}\qquad \GLf^\uparrow_{\vs'}=\GLf^\uparrow_\vs+\ve_j\,O_j\otimes O_j .
\]
Hence $\GLf_\vs=\GLf_0+\sum_{\vs_j=1}\ve_j\,U_j\otimes U_j$ for every smoothing $\vs$.
\item \emph{The oriented smoothing.} On $H_1(F_D;\BZ)$ the form $\GLf_0$ is the symmetrized Seifert form $V_F+V_F^{\T}$ of the canonical surface, and $\GLf^\uparrow_0$ is the symmetrized Seifert form $V_{F^\uparrow}+V_{F^\uparrow}^{\T}$ of the canonical surface above the plane.
\end{enumerate}
\end{proposition}
\begin{proof}
(i) By \eqref{eq.UOw}, $\theta_j\,U_j\otimes w^\vs_j=-\theta_j\bigl(w^\vs_j+w^{\bar\vs}_j\bigr)\otimes w^\vs_j$, and the sum over $j$ is $\Goe_\vs+\Stk_\vs$. Likewise $-\ve_j\,O_j\otimes w^\vs_j=-\ve_j(-1)^{\vs_j}\bigl(w^\vs_j-w^{\bar\vs}_j\bigr)\otimes w^\vs_j=-\theta_j\,w^\vs_j\otimes w^\vs_j+\theta_j\,w^{\bar\vs}_j\otimes w^\vs_j$, and the sum over $j$ is $\Goe_\vs-\Stk_\vs$. The two right-hand sides are the forms of Gordon and Litherland by Corollary~\ref{cor.GLform} and by Proposition~\ref{prop.otherside}~(i).

(ii) The form $\Goe_\vs$ is a combination of squares of linear forms. For the stacking form,
\[
   \Stk_\vs(z,z')-\Stk_\vs(z',z)=-\sum_j\theta_j\Bigl(w^{\bar\vs}_j(z)\,w^\vs_j(z')-w^\vs_j(z)\,w^{\bar\vs}_j(z')\Bigr).
\]
If $\vs_j=0$, then $\theta_j=\ve_j$, $w^\vs_j=w^0_j$ and $w^{\bar\vs}_j=w^1_j$, and the term of $j$ is $\ve_j\bigl(w^0_j(z)w^1_j(z')-w^1_j(z)w^0_j(z')\bigr)$. If $\vs_j=1$, then $\theta_j=-\ve_j$, $w^\vs_j=w^1_j$ and $w^{\bar\vs}_j=w^0_j$, and the term of $j$ is the same expression. So the difference is the left-hand side of Lemma~\ref{lem.basicidentity}, which is zero. Hence $\Stk_\vs$ is symmetric, and so are $\GLf_\vs$ and $\GLf^\uparrow_\vs$. For Kashaev's form, by \eqref{eq.UOw},
\[
   \ve_j\,U_j\otimes O_j=-\ve_j\bigl(w^0_j+w^1_j\bigr)\otimes\bigl(w^0_j-w^1_j\bigr)=\ve_j\bigl(w^1_j\otimes w^1_j-w^0_j\otimes w^0_j\bigr)+\ve_j\bigl(w^0_j\otimes w^1_j-w^1_j\otimes w^0_j\bigr),
\]
and the sum over $j$ of the last terms vanishes by Lemma~\ref{lem.basicidentity}. This proves the first equality of \eqref{eq.Kashw}, which shows that $\Kash$ is symmetric.

(iii) Since $\theta_j(\bar\vs)=-\theta_j$, the term of $j$ in $\Goe_\vs+\Goe_{\bar\vs}$ is $-\theta_j\bigl(w^\vs_j\otimes w^\vs_j-w^{\bar\vs}_j\otimes w^{\bar\vs}_j\bigr)$, which is $\ve_j\bigl(w^1_j\otimes w^1_j-w^0_j\otimes w^0_j\bigr)$ both for $\vs_j=0$ and for $\vs_j=1$. Next, $\Stk_{\bar\vs}=-\sum_j\theta_j(\bar\vs)\,w^\vs_j\otimes w^{\bar\vs}_j=\sum_j\theta_j\,w^\vs_j\otimes w^{\bar\vs}_j$, which is $-\Stk_\vs$ with its two arguments exchanged, hence $-\Stk_\vs$ by (ii). So $\GLf_\vs+\GLf_{\bar\vs}=\Goe_\vs+\Goe_{\bar\vs}+\Stk_\vs+\Stk_{\bar\vs}=\Kash$, and in the same way $\GLf^\uparrow_\vs+\GLf^\uparrow_{\bar\vs}=\Kash$. By \eqref{eq.piminusone}, $\ve_j\,\pi_j\otimes\tpi_j=\ve_j\,U_j\otimes O_j$ at $x_j=y_j=-1$, and by Definition~\ref{def.deformedkashaev} the entry of $\tau_D(x,y)$ at two faces $f,f'$ is $\sum_j\ve_j\,\pi_j(\partial f)\,\tpi_j(\partial f')$.

(iv) By Lemma~\ref{lem.circles}~(i), $w^\vs(z^C)=0$ for every circle $C$ of $\vs$. So $\Goe_\vs(z,z^C)=\Stk_\vs(z,z^C)=0$ for every flow $z$, and by the symmetry $z^C$ lies in the radicals of the two forms, hence in those of $\GLf_\vs$ and $\GLf^\uparrow_\vs$. At a crossing $j$ the edges $o_j$ and $u'_j$ lie on the Seifert circle $A_j$, the edges $u_j$ and $o'_j$ lie on $B_j$, and $h(B_j)=h(A_j)+\ve_j$ (Lemma~\ref{lem.level}); so $z^0_{o_j}=z^0_{u'_j}=-z^0_{u_j}=-z^0_{o'_j}$, and \eqref{eq.w01} gives $w^0_j(z^0)=0$ and $w^1_j(z^0)=0$. All the forms of Definition~\ref{def.formsminusone} are combinations of products of these linear forms, by \eqref{eq.UOw}; so $z^0$ lies in their radicals. Finally $w^\vs(z^0)=0$ for every smoothing $\vs$, and $z^0$ is a combination of the circles of $\vs$ by Lemma~\ref{lem.circles}~(i).

(v) By (i), the forms $\GLf_\vs$ and $\GLf_{\vs'}$ differ only in the term of the crossing $j$, which is $\ve_j\,U_j\otimes w^0_j$ for $\vs$ and $-\ve_j\,U_j\otimes w^1_j$ for $\vs'$. The difference is $-\ve_j\,U_j\otimes\bigl(w^0_j+w^1_j\bigr)=\ve_j\,U_j\otimes U_j$ by \eqref{eq.UOw}. For the surfaces above the plane the term of $j$ is $-\ve_j\,O_j\otimes w^0_j$ for $\vs$ and $-\ve_j\,O_j\otimes w^1_j$ for $\vs'$, and the difference is $\ve_j\,O_j\otimes\bigl(w^0_j-w^1_j\bigr)=\ve_j\,O_j\otimes O_j$. The last statement follows by changing the smoothings of the crossings with $\vs_j=1$ one at a time, starting from the oriented smoothing.

(vi) For $\vs=0$ we have $\theta_j=\ve_j$ and $w^\vs_j=w_j$, and $U_j(z)=z_{e^+_j}+z_{e^-_j}$, $O_j(z)=z_{\hat e^+_j}+z_{\hat e^-_j}$ in the notation of \eqref{eq.eplus} and of Proposition~\ref{prop.otherside}. So, by (i), by Corollary~\ref{cor.seifertform} and by Proposition~\ref{prop.otherside}~(ii),
\begin{align*}
   \GLf_0(z,z')&=\sum_j\ve_j\bigl(z_{e^+_j}+z_{e^-_j}\bigr)w_j(z')=V_F\bigl([z],[z']\bigr)+V_F\bigl([z'],[z]\bigr),\\
   \GLf^\uparrow_0(z,z')&=-\sum_j\ve_j\bigl(z_{\hat e^-_j}+z_{\hat e^+_j}\bigr)w_j(z')=V_{F^\uparrow}\bigl([z],[z']\bigr)+V_{F^\uparrow}\bigl([z'],[z]\bigr).
\end{align*}
\end{proof}

With Lemma~\ref{lem.cornertransfers}, the forms can be read from the potentials. Let $z$ and $z'$ be flows with potentials $k$ and $k'$. In the notation of Lemma~\ref{lem.cornertransfers}~(ii), and since $\theta_j^2=1$,
\begin{equation}\label{eq.formspot}
\begin{aligned}
   \Goe_\vs(z,z')&=-\sum_{j=1}^c\theta_j\,\bigl(k(\phi'_j)-k(\phi_j)\bigr)\bigl(k'(\phi'_j)-k'(\phi_j)\bigr),\\
   \Stk_\vs(z,z')&=\ \ \sum_{j=1}^c\theta_j\,\bigl(k(c'_j)-k(c_j)\bigr)\bigl(k'(\phi'_j)-k'(\phi_j)\bigr),\\
   \Kash(z,z')&=\ \ \sum_{j=1}^c\ve_j\Bigl(\bigl(k^\leftarrow-k^\rightarrow\bigr)\bigl(k'^\leftarrow-k'^\rightarrow\bigr)-\bigl(k^\uparrow-k^\downarrow\bigr)\bigl(k'^\uparrow-k'^\downarrow\bigr)\Bigr),
\end{aligned}
\end{equation}
where, in the term of the crossing $j$ of the last formula, $k^\leftarrow=k(f^\leftarrow_j)$ and so on; this formula follows from \eqref{eq.Kashw} and Lemma~\ref{lem.cornertransfers}~(i). The Goeritz form only involves the free corners of the crossings, and the stacking form pairs the covered corners of each crossing with its free corners.

\subsection{Goeritz matrices and Kashaev's matrix at \texorpdfstring{$t=-1$}{t = -1}}\label{ssec.goeritz}
We now compute the forms of Section~\ref{ssec.formsminusone} in the basis of the faces, for the two checkerboard smoothings of the diagram. Their Goeritz forms are the two classical Goeritz matrices, their stacking forms vanish, and Kashaev's form is the direct sum of the two Goeritz matrices; this determines Kashaev's matrix at $t=-1$ for every diagram (Theorem~\ref{thm.goeritz}). We first recall the checkerboard surfaces and the Goeritz matrices.

Color the faces of $D$ black and white in such a way that the two faces adjacent to an edge have different colors; this is possible, for instance by the parity of the winding numbers of Lemma~\ref{lem.winding}, and the coloring is unique up to the exchange of the two colors. At every crossing $j$ the faces $f^\leftarrow_j$ and $f^\rightarrow_j$ have one color and the faces $f^\uparrow_j$ and $f^\downarrow_j$ have the other one, because two consecutive corners of a crossing are separated by an edge. As in Remark~\ref{rem.checkerboardbasis}, let $\vs_{\rm b}$ be the smoothing whose covered quadrants are the black corners, and let $\vs_{\rm w}=\bar\vs_{\rm b}$ be the complementary smoothing, whose covered quadrants are the white corners:
\[
   (\vs_{\rm b})_j=0\ \text{ if $f^\leftarrow_j$ and $f^\rightarrow_j$ are black},\qquad(\vs_{\rm b})_j=1\ \text{ if $f^\uparrow_j$ and $f^\downarrow_j$ are black.}
\]
Every black face $f$ contains exactly one circle of $\vs_{\rm b}$, which follows the boundary of $f$ and cuts off its corners, and $\partial f$ is, up to sign, the flow of this circle. The state surface $F_{\vs_{\rm b}}$ consists of a disk for each black face and of a band over the two black corners of each crossing: it is the \emph{black checkerboard surface} of $D$, and $F_{\vs_{\rm w}}$ is the white one.

For a color $v\in\{{\rm b},{\rm w}\}$ and a crossing $j$ put
\begin{equation}\label{eq.eta}
   \eta_v(j)=\begin{cases}\ \ \ve_j&\text{if the two faces of color $v$ at $j$ are $f^\uparrow_j$ and $f^\downarrow_j$,}\\ -\ve_j&\text{if they are $f^\leftarrow_j$ and $f^\rightarrow_j$.}\end{cases}
\end{equation}
These are the signs $\eta_v(c)$ of \cite[Section~2.2 and Figure~3]{CimasoniFerretti24}; they depend on the crossing and on the coloring, but not on the orientation of the link (Figure~\ref{fig.goeritz}). (In the notation of that paper, $t_v(c)=1$ in the first case and $t_v(c)=-1$ in the second one, and $\eta_v(c)\,t_v(c)$ is the sign of the crossing \cite[Lemma~3.2]{CimasoniFerretti24}.) Since the two faces of one color at $j$ are $f^\uparrow_j$ and $f^\downarrow_j$ exactly when the two faces of the other color are $f^\leftarrow_j$ and $f^\rightarrow_j$, we have $\eta_{\rm b}(j)=-\eta_{\rm w}(j)$. The \emph{Goeritz matrix} of $D$ for the color $v$ is the symmetric integer matrix $G_v(D)$, indexed by the faces of color $v$, with the entries
\[
   G_v(D)_{f,g}=\sum_{j}\eta_v(j)\quad(f\ne g),\qquad\qquad G_v(D)_{f,f}=-\sum_{g\ne f}G_v(D)_{f,g},
\]
where the first sum runs over the crossings $j$ at which the two faces of color $v$ are $f$ and $g$ \cite{Goeritz33, GordonLitherland}; we follow the conventions of \cite[Section~2.2]{CimasoniFerretti24}, in which the matrix is not reduced, so that the sum of its columns is zero. Equivalently, if $p_j$ and $q_j$ denote the two faces of color $v$ at the crossing $j$, then for all vectors $k,k'$ indexed by the faces of color $v$
\begin{equation}\label{eq.goeritzform}
   k^{\T}G_v(D)\,k'=-\sum_{j=1}^c\eta_v(j)\,\bigl(k(p_j)-k(q_j)\bigr)\bigl(k'(p_j)-k'(q_j)\bigr).
\end{equation}
Indeed, for the unit vectors of two faces $f\ne g$ the right-hand side is the sum of $\eta_v(j)$ over the crossings with $\{p_j,q_j\}=\{f,g\}$, and for $f=g$ it is minus the sum of $\eta_v(j)$ over the crossings at which exactly one of the two faces $p_j,q_j$ is $f$.

\begin{figure}[ht]
\centering
\begin{tikzpicture}[>=Stealth,scale=0.72,every node/.style={font=\small}]
\def\shUD{\fill[gray!30] (0,0) -- (45:1.9) arc (45:135:1.9) -- cycle; \fill[gray!30] (0,0) -- (225:1.9) arc (225:315:1.9) -- cycle;}
\def\shLR{\fill[gray!30] (0,0) -- (135:1.9) arc (135:225:1.9) -- cycle; \fill[gray!30] (0,0) -- (-45:1.9) arc (-45:45:1.9) -- cycle;}
\def\pos{\draw[thick] (1.45,-1.45) -- (0.22,-0.22); \draw[thick,->] (-0.22,0.22) -- (-1.45,1.45); \draw[thick,->] (-1.45,-1.45) -- (1.45,1.45);}
\def\negc{\draw[thick] (-1.45,-1.45) -- (-0.22,-0.22); \draw[thick,->] (0.22,0.22) -- (1.45,1.45); \draw[thick,->] (1.45,-1.45) -- (-1.45,1.45);}
\begin{scope}
\shUD \pos
\node at (0,-2.45) {$\ve_j=1$}; \node at (0,-3.15) {$\eta_v(j)=1$}; \node at (0,-3.85) {$t_v=1$};
\end{scope}
\begin{scope}[xshift=4.6cm]
\shLR \pos
\node at (0,-2.45) {$\ve_j=1$}; \node at (0,-3.15) {$\eta_v(j)=-1$}; \node at (0,-3.85) {$t_v=-1$};
\end{scope}
\begin{scope}[xshift=9.2cm]
\shUD \negc
\node at (0,-2.45) {$\ve_j=-1$}; \node at (0,-3.15) {$\eta_v(j)=-1$}; \node at (0,-3.85) {$t_v=1$};
\end{scope}
\begin{scope}[xshift=13.8cm]
\shLR \negc
\node at (0,-2.45) {$\ve_j=-1$}; \node at (0,-3.15) {$\eta_v(j)=1$}; \node at (0,-3.85) {$t_v=-1$};
\end{scope}
\end{tikzpicture}
\caption{The signs $\eta_v(j)$ of \eqref{eq.eta}; the two faces of color $v$ at the crossing are shaded. In the first and in the third picture they are $f^\uparrow_j$ and $f^\downarrow_j$, and $\eta_v(j)=\ve_j$; in the second and in the fourth one they are $f^\leftarrow_j$ and $f^\rightarrow_j$, and $\eta_v(j)=-\ve_j$. The sign does not depend on the orientations of the strands: the fourth picture is the first one, turned by a quarter of a turn, with the orientation of one of the two strands reversed, and the third one is obtained in the same way from the second one. The numbers $t_v=t_v(c)$ are those of \cite{CimasoniFerretti24}.}
\label{fig.goeritz}
\end{figure}

\begin{theorem}[Goeritz matrices and Kashaev's matrix at $t=-1$]\label{thm.goeritz}
Let $D$ be a diagram with a checkerboard coloring of its faces.
\begin{enumerate}\renewcommand{\theenumi}{\roman{enumi}}
\item \emph{The band signs.} For every crossing $j$, $\theta_j(\vs_{\rm b})=\eta_{\rm w}(j)$ and $\theta_j(\vs_{\rm w})=\eta_{\rm b}(j)$.
\item \emph{The checkerboard surfaces.} The stacking form of $\vs_{\rm b}$ vanishes, so that $\GLf_{\vs_{\rm b}}=\GLf^\uparrow_{\vs_{\rm b}}=\Goe_{\vs_{\rm b}}$, and
\[
   \Goe_{\vs_{\rm b}}(\partial f,\partial g)=\begin{cases}G_{\rm w}(D)_{f,g}&\text{if $f$ and $g$ are white,}\\ 0&\text{if $f$ or $g$ is black.}\end{cases}
\]
Thus the form of Gordon and Litherland of the black checkerboard surface is given by the Goeritz matrix of the white faces, the boundaries of the black faces being zero in the homology of the surface; and the same holds with the two colors exchanged.
\item \emph{Kashaev's form.} $\Kash(\partial f,\partial g)=G_v(D)_{f,g}$ if the faces $f$ and $g$ have the same color $v$, and $\Kash(\partial f,\partial g)=0$ if they have different colors. In other words, with the faces ordered by their colors,
\[
   \tau_D(-1,-1)=\bigl(\Kash(\partial f,\partial g)\bigr)_{f,g\in\CF}=G_{\rm w}(D)\oplus G_{\rm b}(D).
\]
\item \emph{Kashaev's matrix.} Let $\tau_D|_{\zeta=0}$ be Kashaev's matrix \textup(Definition~\ref{def.kashaev}\textup) at $\zeta=0$, that is, at $t=-1$, and let $\mathsf P=\diag\bigl((-1)^{\lfloor w(f)/2\rfloor}\bigr)_{f\in\CF}$, where $w$ is the winding function of Lemma~\ref{lem.winding}. Then
\[
   \mathsf P\ \tau_D|_{\zeta=0}\ \mathsf P=G_{\rm w}(D)\oplus G_{\rm b}(D)
\]
for every diagram $D$.
\end{enumerate}
\end{theorem}
\begin{proof}
(i) If $(\vs_{\rm b})_j=0$, the black faces at $j$ are $f^\leftarrow_j,f^\rightarrow_j$ and the white ones are $f^\uparrow_j,f^\downarrow_j$; so $\eta_{\rm w}(j)=\ve_j=\theta_j(\vs_{\rm b})$. If $(\vs_{\rm b})_j=1$, the white faces at $j$ are $f^\leftarrow_j,f^\rightarrow_j$ and $\eta_{\rm w}(j)=-\ve_j=\theta_j(\vs_{\rm b})$. The second statement follows, since $\theta_j(\vs_{\rm w})=-\theta_j(\vs_{\rm b})$ and $\eta_{\rm b}(j)=-\eta_{\rm w}(j)$.

(ii) For the smoothing $\vs_{\rm b}$ the covered quadrants of a crossing $j$ are its black corners and the free quadrants are its white corners: in the notation of Lemma~\ref{lem.cornertransfers}~(ii), $c'_j$ and $c_j$ are the two black faces at $j$, and $\phi'_j$ and $\phi_j$ the two white ones. Let $f,g$ be two faces; the flows $\partial f$ and $\partial g$ have the potentials $e_f$ and $e_g$. By \eqref{eq.formspot} and (i),
\[
   \Goe_{\vs_{\rm b}}(\partial f,\partial g)=-\sum_j\eta_{\rm w}(j)\,\bigl(e_f(\phi'_j)-e_f(\phi_j)\bigr)\bigl(e_g(\phi'_j)-e_g(\phi_j)\bigr).
\]
This is zero if $f$ or $g$ is black, and it is $G_{\rm w}(D)_{f,g}$ if both are white, by \eqref{eq.goeritzform}. By \eqref{eq.formspot} again, $\Stk_{\vs_{\rm b}}(\partial f,\partial g)$ is a sum of products of a factor $e_f(c'_j)-e_f(c_j)$, which vanishes if $f$ is white, and of a factor $e_g(\phi'_j)-e_g(\phi_j)$, which vanishes if $g$ is black. So $\Stk_{\vs_{\rm b}}(\partial f,\partial g)=0$ unless $f$ is black and $g$ is white; and in that case it is zero as well, because $\Stk_{\vs_{\rm b}}(\partial f,\partial g)=\Stk_{\vs_{\rm b}}(\partial g,\partial f)$ by Proposition~\ref{prop.formsminusone}~(ii). Since the boundaries of the faces generate $\FF_D$, the stacking form vanishes, and $\GLf_{\vs_{\rm b}}=\GLf^\uparrow_{\vs_{\rm b}}=\Goe_{\vs_{\rm b}}$ by Definition~\ref{def.formsminusone}. The interpretation is Proposition~\ref{prop.formsminusone}~(i).

(iii) By \eqref{eq.Kashw}, $\Kash=\Goe_{\vs_{\rm b}}+\Goe_{\vs_{\rm w}}$. By (ii), the first form is $G_{\rm w}(D)$ on the white faces and vanishes if one of the two faces is black; by (ii) with the two colors exchanged, the second one is $G_{\rm b}(D)$ on the black faces and vanishes if one of the two faces is white. The matrix of $\Kash$ on the boundaries of the faces is $\tau_D(-1,-1)$ by Proposition~\ref{prop.formsminusone}~(iii).

(iv) At $\zeta=0$ the block of the crossing $j$ in Definition~\ref{def.kashaev} is, in the basis $(f^\uparrow_j,f^\leftarrow_j,f^\downarrow_j,f^\rightarrow_j)$,
\[
   \tau_j\big|_{\zeta=0}=\begin{pmatrix}-1&0&1&0\\0&1&0&1\\1&0&-1&0\\0&1&0&1\end{pmatrix},
\]
so that, for all $k,k'\in\BZ^{\CF}$, with $k^\uparrow=k(f^\uparrow_j)$ and so on,
\[
   k^{\T}\,\tau_D|_{\zeta=0}\,k'=\sum_{j=1}^c\ve_j\Bigl(\bigl(k^\leftarrow+k^\rightarrow\bigr)\bigl(k'^\leftarrow+k'^\rightarrow\bigr)-\bigl(k^\uparrow-k^\downarrow\bigr)\bigl(k'^\uparrow-k'^\downarrow\bigr)\Bigr).
\]
By \eqref{eq.wcrossing}, $w(f^\leftarrow_j)=w(f^\rightarrow_j)+2$ and $w(f^\uparrow_j)=w(f^\downarrow_j)$. Hence the diagonal entries of $\mathsf P$ at $f^\leftarrow_j$ and at $f^\rightarrow_j$ are opposite, those at $f^\uparrow_j$ and at $f^\downarrow_j$ are equal, and
\[
   (\mathsf Pk)^\leftarrow+(\mathsf Pk)^\rightarrow=\pm\bigl(k^\leftarrow-k^\rightarrow\bigr),\qquad(\mathsf Pk)^\uparrow-(\mathsf Pk)^\downarrow=\pm\bigl(k^\uparrow-k^\downarrow\bigr),
\]
with signs which depend on the crossing but not on $k$. Therefore $(\mathsf Pk)^{\T}\tau_D|_{\zeta=0}(\mathsf Pk')$ is the right-hand side of the last formula of \eqref{eq.formspot}, that is, $\Kash(z,z')$ for the flows $z=\sum_fk_f\,\partial f$ and $z'=\sum_fk'_f\,\partial f$. So $\mathsf P\,\tau_D|_{\zeta=0}\,\mathsf P$ is the matrix of $\Kash$ on the boundaries of the faces, and (iii) applies.
\end{proof}

\begin{remark}[comparison with Cimasoni and Ferretti]\label{rem.CF}
Parts (iii) and (iv) make precise the relation between Kashaev's matrix and the Goeritz matrices which is the basis of \cite[Section~3.2]{CimasoniFerretti24}. There it is shown that Kashaev's matrix at $\zeta=0$ splits as a direct sum $\tau_{\rm w}(D)\oplus\tau_{\rm b}(D)$ over the two colors, and that $\tau_v(D)$ is very close to $G_v(D)$: the off-diagonal entries of the two matrices differ by signs $t_{fg}=\pm1$, where $t_{fg}$ is the common value of $t_v(c)$ at the crossings $c$ at which the faces $f$ and $g$ of color $v$ meet. The two matrices are then shown to be congruent for a special diagram, with a suitable coloring, which is sufficient for the purpose of that paper because every link has a special diagram. By (iv), $t_{fg}=(-1)^{(w(f)-w(g))/2}$ is the product of the diagonal entries of $\mathsf P$ at $f$ and at $g$, and the two matrices are congruent, by a diagonal matrix of signs, for every diagram. In the lattice of flows the signs are absent from the start: the conjugation by the winding numbers which relates Kashaev's matrix to the coupled form of the two flavors (Theorem~\ref{thm.kashaev}, whose matrix $\WW_t$ is $\diag\bigl((-\sqrt{-1})^{w(f)}\bigr)$ at $t=-1$) removes them. The splitting itself is the identity $\Kash=\Goe_{\vs_{\rm b}}+\Goe_{\vs_{\rm w}}$, a special case of $\Kash=\Goe_\vs+\Goe_{\bar\vs}$: Kashaev's form is the sum of the Goeritz forms of any two complementary smoothings.
\end{remark}

\subsection{The regions of a smoothing and the two sides of the plane}\label{ssec.sides}
The form of Gordon and Litherland of the state surface below the plane is $\Goe_\vs+\Stk_\vs$, and that of the state surface above the plane is $\Goe_\vs-\Stk_\vs$ (Proposition~\ref{prop.formsminusone}). In this subsection we describe the two summands in terms of the regions which the circles of the smoothing cut out in the plane. The Goeritz form is a sum of forms which belong to the regions; the stacking form couples the regions which are adjacent; and an involution of the homology of the surface, given by the parity of the depth of the regions, exchanges the forms of the two surfaces (Proposition~\ref{prop.J}).

\begin{figure}[ht]
\centering
\begin{tikzpicture}[>=Stealth,every node/.style={font=\small}]
\begin{scope}[scale=0.95]
\fill[gray!22] (0,0) circle (2.7);
\fill[white] (0,0) circle (1.8);
\fill[gray!22] (0,0) circle (0.9);
\draw[thick] (0,0) circle (2.7); \draw[thick] (0,0) circle (1.8); \draw[thick] (0,0) circle (0.9);
\draw[line width=3pt,gray!75] (212:1.8) -- (212:2.7); \draw[line width=3pt,gray!75] (158:1.8) -- (158:2.7);
\draw[line width=3pt,gray!75] (186:0.9) -- (186:1.8); \draw[line width=3pt,gray!75] (132:0.9) -- (132:1.8);
\node[font=\scriptsize] at (217:2.95) {$1$}; \node[font=\scriptsize] at (153:2.95) {$3$};
\node[font=\scriptsize] at (196:1.35) {$2$}; \node[font=\scriptsize] at (122:1.35) {$4$};
\node[font=\scriptsize] at (60:2.92) {$C_1$}; \node[font=\scriptsize] at (60:2.02) {$C_2$}; \node[font=\scriptsize] at (60:1.12) {$C_3$};
\node[font=\scriptsize] at (0,0.2) {$\Omega_3$}; \node[font=\scriptsize] at (0,-0.22) {$d=3$};
\node[font=\scriptsize] at (0,-1.35) {$\Omega_2,\ d=2$};
\node[font=\scriptsize] at (0,-2.25) {$\Omega_1,\ d=1$};
\node[font=\scriptsize] at (2.75,-2.6) {$\Omega_0,\ d=0$};
\node at (0,-3.55) {(a)};
\end{scope}
\begin{scope}[xshift=8.3cm,scale=1.05]
\fill[gray!14] (135:2.6) -- (135:1.1) arc (135:45:1.1) -- (45:2.6) arc (45:135:2.6);
\fill[gray!14] (225:2.6) -- (225:1.1) arc (225:315:1.1) -- (315:2.6) arc (315:225:2.6);
\fill[gray!14] (0,0) circle (1.1);
\fill[gray!45] (0,0) -- (135:1.1) arc (135:225:1.1) -- cycle;
\fill[gray!45] (0,0) -- (-45:1.1) arc (-45:45:1.1) -- cycle;
\draw[very thick] (135:2.6) -- (135:1.1) arc (135:225:1.1) -- (225:2.6);
\draw[very thick] (45:2.6) -- (45:1.1) arc (45:-45:1.1) -- (-45:2.6);
\draw (135:1.1) -- (-45:1.1); \draw[white,line width=3pt] (225:0.3) -- (45:0.3); \draw (225:1.1) -- (45:1.1);
\draw[gray,thin] (45:1.1) arc (45:135:1.1); \draw[gray,thin] (225:1.1) arc (225:315:1.1);
\node[font=\scriptsize] at (-2.2,0) {$f^\circ\subset\Omega'$}; \node[font=\scriptsize] at (2.25,0) {$f'^\circ\subset\Omega''$};
\node[font=\scriptsize] at (0,1.85) {$g^\circ\subset\Omega(j)$}; \node[font=\scriptsize] at (0,-1.85) {$g'^\circ\subset\Omega(j)$};
\node[font=\scriptsize] at (-0.62,0) {$f$}; \node[font=\scriptsize] at (0.62,0) {$f'$};
\node[font=\scriptsize] at (0,0.62) {$g$}; \node[font=\scriptsize] at (0,-0.62) {$g'$};
\node[font=\scriptsize] at (132:2.85) {$C$}; \node[font=\scriptsize] at (48:2.85) {$C'$}; \node[font=\scriptsize] at (97:1.32) {$\CB_j$};
\node at (0,-3.2) {(b)};
\end{scope}
\end{tikzpicture}
\caption{(a) The regions of the oriented smoothing of the closure of $(\sigma_1\sigma_2^{-1})^2$, a diagram of the figure-eight knot (Example~\ref{ex.figeightminusone}). The Seifert circles $C_1,C_2,C_3$ are drawn as concentric circles and the bands of the four crossings in gray. The regions $\Omega_0,\dots,\Omega_3$ have the depths $0,\dots,3$; the crossings $1$ and $3$ lie in $\Omega_1$, the crossings $2$ and $4$ in $\Omega_2$, and the involution $J$ is $-1$ on the homology of the regions of odd depth (shaded). (b) A crossing $j$ and the circles $C$ and $C'$ which contain the two arcs of its smoothing (thick). The region $\Omega(j)$ (light gray) contains the disk $\CB_j$. The faces $g,g'$ of the free quadrants have $g^\circ,g'^\circ\subset\Omega(j)$; the faces $f,f'$ of the covered quadrants (dark) have $f^\circ$ and $f'^\circ$ on the other side of $C$ and of $C'$, in the regions $\Omega'$ and $\Omega''$ adjacent to $\Omega(j)$ (Lemma~\ref{lem.smoothingregions}~(i)). The circles $C$ and $C'$ may coincide, and so may the faces $g,g'$, or $f,f'$.}
\label{fig.smoothingregions}
\end{figure}
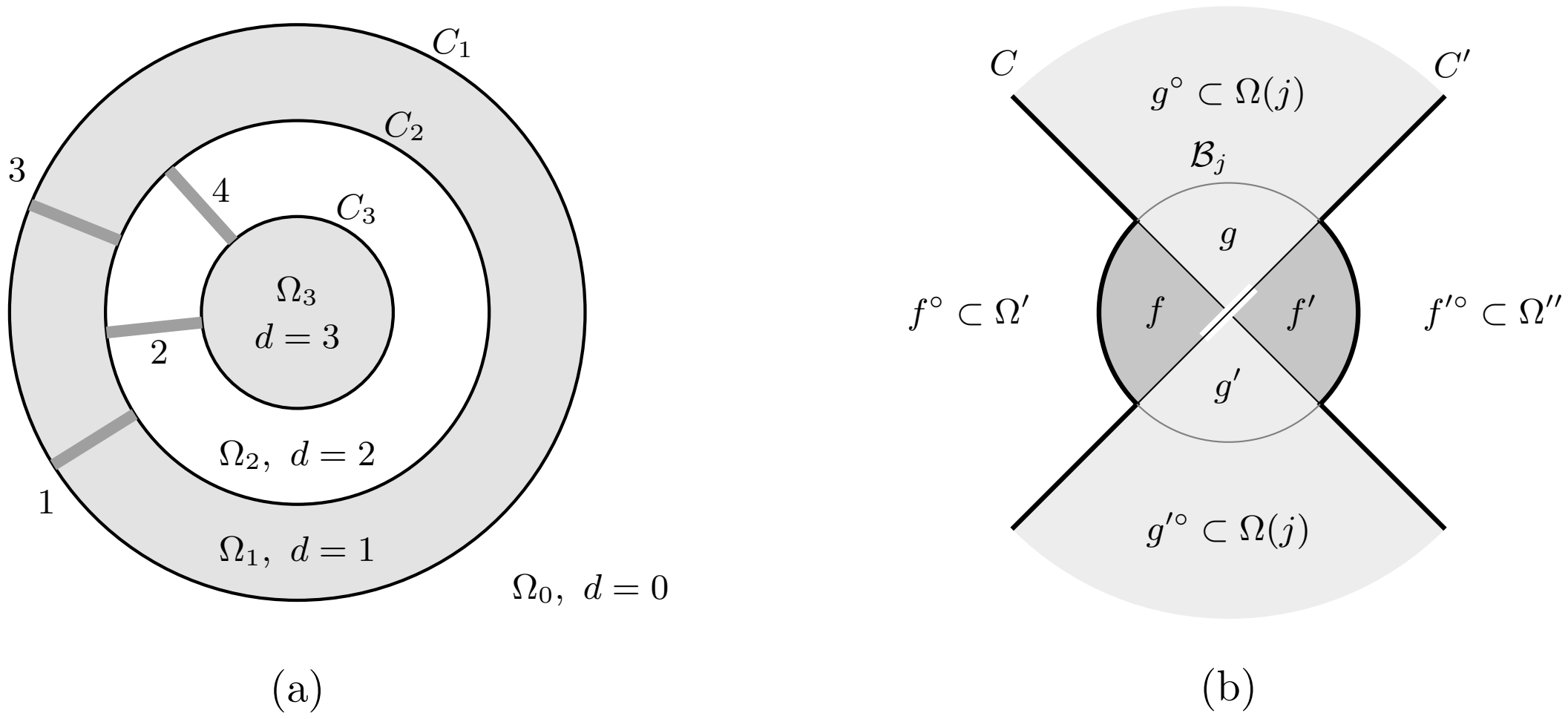

Let $\vs$ be a smoothing of $D$. Its circles are pairwise disjoint simple closed curves in the plane $P$, made of the edges of $D$ outside the disks $\CB_j$ and of the arcs of the smoothing on the circles $\partial\CB_j$ (Section~\ref{ssec.statesurface}). The connected components of the complement in $P$ of the union of the circles are the \emph{regions} of $\vs$; for the oriented smoothing these are the regions of Section~\ref{ssec.regions}. Every circle $C$ lies in the boundary of exactly two regions, the one directly inside $C$ and the one directly outside it; we call two such regions \emph{adjacent}. The \emph{depth} $d(\Omega)$ of a region $\Omega$ is the number of circles of $\vs$ which enclose it; the depths of two adjacent regions differ by one (Figure~\ref{fig.smoothingregions}~(a)).

The crossings and the faces of $D$ are distributed among the regions as follows.\label{not.regions} The interior of the disk $\CB_j$ does not meet the circles, so it is contained in one region $\Omega(j)$; we say that the crossing $j$ \emph{lies in} $\Omega(j)$. For a face $f$, let $f^\circ$ be the set of the points of $f$ which do not belong to a covered quadrant of a crossing: it is the face $f$ with its covered corners cut off along the arcs of the smoothing. The set $f^\circ$ is connected, because cutting a corner off a disk leaves a disk, and it does not meet the circles; so it is contained in one region $\Omega_f$. For a region $\Omega$ we denote by $\CF_\Omega$ the set of the faces $f$ with $\Omega_f=\Omega$, and by $\1_\Omega\in\BZ^{\CF}$ the indicator vector of this set.

\begin{lemma}[regions, faces and circles]\label{lem.smoothingregions}
Let $\vs$ be a smoothing of $D$.
\begin{enumerate}\renewcommand{\theenumi}{\roman{enumi}}
\item If a face $f$ has a free corner at a crossing $j$, then $\Omega_f=\Omega(j)$. If $f$ has a covered corner at $j$, then $\Omega_f$ and $\Omega(j)$ are adjacent: they are the two regions on the two sides of the circle which contains the arc of the smoothing in this corner \textup(Figure~\ref{fig.smoothingregions}~\textup{(b)}\textup).
\item Let $C$ be a circle of $\vs$, traveled counterclockwise. Then $z^C=\sum\partial f$, the sum over the faces $f$ such that $f^\circ$ lies inside $C$.
\item A vector $k\in\BZ^{\CF}$ is the potential of a combination of circles of $\vs$ if and only if it is a combination of the vectors $\1_\Omega$. Consequently
\[
   H_1(F_\vs;\BZ)=\FF_D\big/\langle z^C\rangle\ \cong\ \BZ^{\CF}\big/\langle\1_\Omega\rangle=\bigoplus_{\Omega}H_\Omega,\qquad H_\Omega=\BZ^{\CF_\Omega}\big/\BZ\,\1_\Omega,
\]
where the sum runs over the regions of $\vs$ and $H_\Omega$ is the subgroup of $H_1(F_\vs;\BZ)$ generated by the classes of the boundaries of the faces of $\CF_\Omega$.
\end{enumerate}
\end{lemma}
\begin{proof}
(i) Let $Q$ be a quadrant of the crossing $j$ which is a corner of $f$. The points of the interior of $Q$ belong to $f$ and to the interior of $\CB_j$, hence to $\Omega(j)$. If $Q$ is free, these points belong to $f^\circ$, so $\Omega_f=\Omega(j)$. If $Q$ is covered, let $a$ be the arc of the smoothing on the boundary of $Q$ and $C$ the circle which contains it. The points of $f$ which lie just outside $\CB_j$, near $a$, belong to $f^\circ$, and they are separated from the interior of $Q$ by the arc $a$. The two sides of an arc of $C$ belong to the two regions adjacent to $C$; so these are $\Omega_f$ and $\Omega(j)$.

(ii) Let $k\in\BZ^{\CF}$ be the indicator vector of the set of the faces $f$ such that $f^\circ$ lies inside $C$, and let $z=\sum_fk_f\,\partial f$. Let $e$ be an edge of $D$, and let $C_e$ be the circle of $\vs$ which contains the part of $e$ outside the disks. Near a point of this part, the points on the left of $e$ belong to $f_{\rm L}(e)^\circ$ and the points on its right to $f_{\rm R}(e)^\circ$. If $C_e\ne C$, these points lie on the same side of $C$, so that $k_{f_{\rm L}(e)}=k_{f_{\rm R}(e)}$ and $z_e=0$ by \eqref{eq.potential}. If $C_e=C$, exactly one of the two sides of $e$ lies inside $C$: it is the left side if the counterclockwise direction of $C$ is the direction of $e$, and then $z_e=1$; otherwise $z_e=-1$. In all cases $z_e=z^C_e$.

(iii) By (ii) and by Lemma~\ref{lem.faces}~(ii), the potentials of the combinations of the circles are the combinations of the constant vector $\1$ and of the vectors $\1_{{\rm in}\,C}=\sum_{\Omega\text{ inside }C}\1_\Omega$. These are combinations of the vectors $\1_\Omega$, and $\1=\sum_\Omega\1_\Omega$. Conversely, if $\Omega$ is a bounded region, with outer boundary circle $C_0$ and inner boundary circles $C_1,\dots,C_r$, then $\1_\Omega=\1_{{\rm in}\,C_0}-\sum_{i=1}^r\1_{{\rm in}\,C_i}$; and for the unbounded region, $\1_\Omega=\1-\sum\1_{{\rm in}\,C}$, the sum over the circles which are enclosed by no other circle. This proves the first statement. Every face belongs to exactly one set $\CF_\Omega$, so $\BZ^{\CF}=\bigoplus_\Omega\BZ^{\CF_\Omega}$, with $\1_\Omega\in\BZ^{\CF_\Omega}$; and the map which sends a vector $k$ to the flow $\sum_fk_f\,\partial f$ induces an isomorphism of $\BZ^{\CF}/\BZ\1$ onto $\FF_D$ by Lemma~\ref{lem.faces}~(ii), which maps the classes of the vectors $\1_\Omega$ onto the subgroup generated by the circles. The decomposition follows.
\end{proof}

We let $J$ be the involution of $H_1(F_\vs;\BZ)$ which is the multiplication by $(-1)^{d(\Omega)}$ on $H_\Omega$, for every region $\Omega$. (If another face is chosen as the unbounded one, the depths change, and $J$ is replaced by $J$ or by $-J$.)

\begin{proposition}[the two sides of the plane]\label{prop.J}
Let $\vs$ be a smoothing of $D$, and regard the forms of Definition~\ref{def.formsminusone} as forms on $H_1(F_\vs;\BZ)$.
\begin{enumerate}\renewcommand{\theenumi}{\roman{enumi}}
\item \emph{Blocks.} The Goeritz form is block diagonal for the decomposition of Lemma~\ref{lem.smoothingregions}: $\Goe_\vs(a,b)=0$ if $a\in H_\Omega$ and $b\in H_{\Omega'}$ with $\Omega\ne\Omega'$. The stacking form only couples adjacent regions: $\Stk_\vs(a,b)=0$ if $a\in H_\Omega$ and $b\in H_{\Omega'}$, unless $\Omega$ and $\Omega'$ are adjacent.
\item \emph{The involution.} $\Goe_\vs(Ja,Jb)=\Goe_\vs(a,b)$ and $\Stk_\vs(Ja,Jb)=-\Stk_\vs(a,b)$ for all $a,b$. Hence
\[
   \CG_{F^\uparrow_\vs}(a,b)=\CG_{F_\vs}(Ja,Jb)\qquad\text{for all }a,b\in H_1(F_\vs;\BZ):
\]
the forms of Gordon and Litherland of the state surfaces above and below the plane are isometric, by the involution $J$.
\item \emph{Circles with crossings on one side.} If no circle of $\vs$ has crossings on both of its sides, that is, if all the crossings lie in one region, then $\Stk_\vs=0$ and $\CG_{F^\uparrow_\vs}=\CG_{F_\vs}=\Goe_\vs$. This is the case for the two checkerboard smoothings.
\end{enumerate}
\end{proposition}
\begin{proof}
(i) It suffices to consider the classes of the boundaries of two faces $f\in\CF_\Omega$ and $g\in\CF_{\Omega'}$, whose potentials are $e_f$ and $e_g$. By \eqref{eq.formspot}, a crossing $j$ contributes to $\Goe_\vs(\partial f,\partial g)$ only if both $f$ and $g$ are faces of free quadrants of $j$, and then $\Omega=\Omega'=\Omega(j)$ by Lemma~\ref{lem.smoothingregions}~(i). It contributes to $\Stk_\vs(\partial f,\partial g)$ only if $f$ is the face of a covered quadrant and $g$ the face of a free quadrant of $j$, and then $\Omega'=\Omega(j)$ and $\Omega$ is adjacent to it.

(ii) On $H_\Omega\times H_{\Omega'}$, replacing $(a,b)$ by $(Ja,Jb)$ multiplies a bilinear form by $(-1)^{d(\Omega)+d(\Omega')}$, which is $1$ if $\Omega=\Omega'$ and $-1$ if $\Omega$ and $\Omega'$ are adjacent. So (i) gives the two identities, and then $\GLf^\uparrow_\vs(a,b)=\Goe_\vs(a,b)-\Stk_\vs(a,b)=\GLf_\vs(Ja,Jb)$, which is the statement by Proposition~\ref{prop.formsminusone}~(i).

(iii) A circle $C$ has crossings on both sides if there are two crossings, each with an arc on $C$, which lie in the two different regions adjacent to $C$. If there is no such circle, any two crossings which have an arc on a common circle lie in the same region; since the state graph is connected, all the crossings then lie in one region $\Omega^*$. By the proof of (i), for two faces $f,g$ the number $\Stk_\vs(\partial f,\partial g)$ is zero unless $g$ has a free corner at some crossing, so that $\Omega_g=\Omega^*$, and $f$ has a covered corner there, so that $\Omega_f\ne\Omega^*$. Since $\Stk_\vs$ is symmetric, the same holds with $f$ and $g$ exchanged, and the two conditions exclude each other; so $\Stk_\vs$ vanishes on the boundaries of the faces, which generate $\FF_D$. For the smoothing $\vs_{\rm b}$ of Section~\ref{ssec.goeritz}, the circle of a black face $f$ is the boundary of $f^\circ$, and all the crossings lie on the other side of it.
\end{proof}

So the passage from one side of the plane to the other changes the sign of the stacking form and nothing else. For the oriented smoothing this determines how the Seifert forms of the two canonical surfaces of Section~\ref{ssec.otherside} are related.

\begin{corollary}[the two canonical Seifert surfaces]\label{cor.twoseifert}
Let $V_F$ and $V_{F^\uparrow}$ be the Seifert forms of the canonical surfaces $F_D$ and $F^\uparrow_D$ on $H_1(F_D;\BZ)=H_1(F^\uparrow_D;\BZ)=\bigoplus_\Omega H_\Omega$, the sum over the regions of the complement of the Seifert circles. Then
\[
   V_{F^\uparrow}(a,b)=\begin{cases}\ \ V_F(a,b)&\text{if $a$ and $b$ belong to the same summand $H_\Omega$,}\\ -V_F(b,a)&\text{if $a\in H_\Omega$ and $b\in H_{\Omega'}$ with $\Omega\ne\Omega'$.}\end{cases}
\]
\end{corollary}
\begin{proof}
By Theorem~\ref{thm.deformedseifert}~(i) and~(iii), the two surfaces have the same antisymmetrized Seifert form:
\begin{align*}
   V_F(a,b)-V_F(b,a)&=\sum_j\bigl(V_j-V'_j\bigr)(z)\,w_j(z')\\
   &=\sum_j\bigl(V^\uparrow_j-V'^{\uparrow}_j\bigr)(z)\,w_j(z')=V_{F^\uparrow}(a,b)-V_{F^\uparrow}(b,a)
\end{align*}
for the classes $a,b$ of two flows $z,z'$. By Proposition~\ref{prop.formsminusone}~(vi) and Proposition~\ref{prop.J}~(ii), $V_{F^\uparrow}(a,b)+V_{F^\uparrow}(b,a)=V_F(Ja,Jb)+V_F(Jb,Ja)$, which is $V_F(a,b)+V_F(b,a)$ if $a$ and $b$ belong to the same summand $H_\Omega$, and $-V_F(a,b)-V_F(b,a)$ if they belong to the summands of two adjacent regions. Adding the two identities gives $V_{F^\uparrow}(a,b)=V_F(a,b)$ in the first case and $V_{F^\uparrow}(a,b)=-V_F(b,a)$ in the second one. Finally let $a$ and $b$ be the classes of the boundaries of two faces $f\in\CF_\Omega$ and $g\in\CF_{\Omega'}$, where $\Omega\ne\Omega'$ are not adjacent. By Theorem~\ref{thm.deformedseifert}~(i), a crossing $j$ contributes to $V_F(a,b)$ or to $V_{F^\uparrow}(a,b)$ only if $w_j(\partial g)\ne0$, which requires that $g$ has a free corner at $j$ (Lemma~\ref{lem.cornertransfers}), and if one of the four coordinates of $\partial f$ at the edges of the crossing is not zero, which requires that $f$ has a corner at $j$; by Lemma~\ref{lem.smoothingregions}~(i), $\Omega'=\Omega(j)$ and $\Omega$ is $\Omega(j)$ or adjacent to it, which is excluded. So $V_F(a,b)=V_{F^\uparrow}(a,b)=0$, and likewise $V_F(b,a)=0$.
\end{proof}

Thus the two canonical surfaces have the same Seifert form on the homology of each region, and opposite and transposed Seifert forms between different regions. This is the algebraic form of the description of Section~\ref{ssec.otherside}: each of the two surfaces is obtained from the surfaces of the regions by Murasugi sums along the disks of the Seifert circles which have crossings on both sides, and they differ by the side of the disk on which each summand is placed. More generally, for every smoothing $\vs$ the restriction of $\Goe_\vs$ to $H_\Omega$ is a Goeritz form in the sense of \eqref{eq.goeritzform}, with the faces of $\CF_\Omega$ in the role of the faces of one color and the signs $\theta_j$ of the crossings which lie in $\Omega$, and the stacking form describes how these pieces are assembled.

\subsection{The determinant}\label{ssec.determinant}
We now show that the Gram matrices of the forms $\GLf_\vs$, for all the smoothings $\vs$, are obtained from one and the same $c\times c$ matrix: the matrix of the forms $U_1,\dots,U_c$ on the lattice of flows, which is the matrix of the crossing forms at $x=-1$. We need the description of the transfers by the state graph.

Recall from Section~\ref{ssec.statehomology} the state graph $\Gamma_\vs$ of a smoothing $\vs$: its vertices are the circles of $\vs$, and it has one edge for each crossing $j$, which joins the circle $C'_j$ containing the arc of the smoothing at the end of the half-strand $u'_j$ to the circle $C''_j$ containing the arc at the end of $u_j$; these two circles may coincide. Note that the edge $u'_j$ of $D$ lies on $C'_j$ and the edge $u_j$ on $C''_j$. We orient the edge $j$ of $\Gamma_\vs$ from $C'_j$ to $C''_j$, and we let
\[
   Z_1(\Gamma_\vs)=\Bigl\{n\in\BZ^c:\ \sum_{j:\,C''_j=C}n_j=\sum_{j:\,C'_j=C}n_j\ \text{ for every circle $C$ of $\vs$}\Bigr\}
\]
be its cycle space; it is a subgroup of $\BZ^c$ of rank $c-s_\vs+1$, because $\Gamma_\vs$ is connected. For the oriented smoothing, $\Gamma_0$ is the Seifert graph $\Gamma_D$ and $Z_1(\Gamma_0)$ is the cycle space $H$ of Section~\ref{ssec.smoothings}.

\begin{lemma}[transfers and the state graph]\label{lem.stategraph}
Let $\vs$ be a smoothing of $D$, let $\nu_\vs=c-s_\vs+1$, and let $\CT\subseteq\{1,\dots,c\}$ be the set of the edges of a spanning tree of $\Gamma_\vs$.
\begin{enumerate}\renewcommand{\theenumi}{\roman{enumi}}
\item The map $w^\vs\colon\FF_D\to\BZ^c$ induces an isomorphism of $H_1(F_\vs;\BZ)=\FF_D/\langle z^C\rangle$ onto $Z_1(\Gamma_\vs)$.
\item If $z_1,\dots,z_{\nu_\vs}$ are flows whose classes form a basis of $\FF_D/\langle z^C\rangle$, then the vectors $w^\vs(z_1),\dots,w^\vs(z_{\nu_\vs})$ and the unit vectors $e_j$, $j\in\CT$, form a basis of $\BZ^c$.
\item Let $\hat C$ be a circle of $\vs$. The square matrix $\bigl(U_j(z^C)\bigr)_{C\ne\hat C,\ j\in\CT}$, of size $s_\vs-1$, has determinant $\pm1$; and so has the matrix $\bigl(O_j(z^C)\bigr)_{C\ne\hat C,\ j\in\CT}$.
\end{enumerate}
\end{lemma}
\begin{proof}
(i) By Section~\ref{ssec.statehomology}, the surface $F_\vs$ retracts by deformation onto a copy of the graph $\Gamma_\vs$, whose edges are the cores of the bands; so the inclusion induces an isomorphism of $H_1(\Gamma_\vs;\BZ)$ onto $H_1(F_\vs;\BZ)$. The first homology group of a graph is the group of its $1$-cycles, and the $1$-cycle of a closed walk $b$ in $\Gamma_\vs$ is $\sum_jw^\vs_j(b)\,e_j\in Z_1(\Gamma_\vs)$: by definition, $w^\vs_j(b)$ is the number of times $b$ runs through the edge $j$ in the direction of its orientation, from the half-band of $u'_j$ to the half-band of $u_j$, minus the number of times it runs through it in the other direction. By Proposition~\ref{prop.statehomology}, $w^\vs(b)=w^\vs(z)$ for every flow $z$ in the class of $b$. Hence the isomorphism $H_1(F_\vs;\BZ)\to Z_1(\Gamma_\vs)$ is induced by $w^\vs$.

(ii) For every edge $j\notin\CT$ there is exactly one cycle $\gamma_j\in Z_1(\Gamma_\vs)$ whose coordinate at $j$ is $1$ and whose coordinates at the other edges which are not in $\CT$ are $0$: the edge $j$ followed by the path in the tree from its end back to its origin (and $\gamma_j=e_j$ if $j$ is a loop). A cycle is determined by its coordinates at the edges which are not in $\CT$, because a cycle supported in a tree is zero; so the $\nu_\vs$ cycles $\gamma_j$ form a basis of $Z_1(\Gamma_\vs)$, and together with the vectors $e_i$, $i\in\CT$, they form a basis of $\BZ^c$, the matrix of these vectors being triangular with ones on the diagonal when the edges which are not in $\CT$ are put first. By (i) the vectors $w^\vs(z_a)$ form another basis of $Z_1(\Gamma_\vs)$.

(iii) Let $j\in\CT$; it is not a loop, so $C'_j\ne C''_j$. Since $u'_j$ lies on $C'_j$ and $u_j$ on $C''_j$, the number $U_j(z^C)=z^C_{u'_j}+z^C_{u_j}$ is $z^{C}_{u'_j}=\pm1$ if $C=C'_j$, it is $z^C_{u_j}=\pm1$ if $C=C''_j$, and it is $0$ for the other circles. So the matrix is, up to the signs of its entries, the incidence matrix of the tree $\CT$ with the row of the vertex $\hat C$ removed. We show by induction on $s_\vs$ that such a matrix has determinant $\pm1$. If $s_\vs=1$ it is empty. Otherwise the tree has at least two leaves, hence a leaf $C_1\ne\hat C$; the row of $C_1$ has exactly one nonzero entry, $\pm1$, in the column of the edge $j_1$ of the tree at $C_1$. Expanding the determinant along this row gives $\pm1$ times the determinant of the matrix of the same kind for the tree obtained by removing the vertex $C_1$ and the edge $j_1$. The same argument applies to $O_j(z^C)=z^C_{o_j}+z^C_{o'_j}$, because the edges $o_j$ and $o'_j$ lie on the two different arcs of the smoothing of the crossing $j$, hence one of them on $C'_j$ and the other one on $C''_j$.
\end{proof}

Let $n(L)$ denote the corank\label{not.nL} of $V+V^{\T}$, where $V$ is a Seifert matrix of $L$; it is an invariant of $L$, the first Betti number of the double covering of $S^3$ branched along $L$, and $n(L)=0$ if and only if $\Delta_L(-1)\ne0$ \cite[Chapter~9]{Lickorish}.

\begin{theorem}[determinant and cokernel]\label{thm.determinant}
Let $\vs$ be a smoothing of $D$ and $\nu_\vs=c-s_\vs+1$. Let $z_1,\dots,z_{\nu_\vs}\in\FF_{D,\hate}$ be flows whose classes form a basis of $H_1(F_\vs;\BZ)=\FF_D/\langle z^C\rangle$, and let
\[
   G_\vs=\bigl(\GLf_\vs(z_a,z_b)\bigr)_{a,b=1,\dots,\nu_\vs}
\]
be the corresponding Gram matrix of the form of Gordon and Litherland of the state surface $F_\vs$. Let $U=(U_1,\dots,U_c)\colon\FF_D\to\BZ^c$.
\begin{enumerate}\renewcommand{\theenumi}{\roman{enumi}}
\item $\det G_\vs=\pm P_{D,\hate}(-1,\dots,-1)$. Consequently $|\det G_\vs|=|\Delta_L(-1)|$ is the determinant of the link, for every smoothing.
\item The abelian group $\BZ^{\nu_\vs}/G_\vs\BZ^{\nu_\vs}$ presented by $G_\vs$ is isomorphic to $\BZ^c/U(\FF_D)$. In particular it does not depend on the smoothing, up to isomorphism; it is the group presented by the symmetrized Seifert matrix $V_F+V_F^{\T}$ of the canonical surface, which is the first homology group of the double covering of $S^3$ branched along $L$ \textup{\cite[Theorem~9.1]{Lickorish}}.
\item $\rk\GLf_\vs=\nu_\vs-n(L)=c+1-s_\vs-n(L)$.
\item The same statements hold for the Gram matrices of the forms $\GLf^\uparrow_\vs$ of the state surfaces above the plane.
\end{enumerate}
\end{theorem}
\begin{proof}
The idea is to multiply the matrix of the forms $U_1,\dots,U_c$ in a basis adapted to the smoothing, whose determinant is the edge polynomial at $x=-1$ up to sign, by a matrix of determinant $\pm1$ which is built from the transfers. The product is block triangular, and its first diagonal block is the Gram matrix $G_\vs$. Let $\hat C$ be the circle of $\vs$ which contains $\hate$. By Proposition~\ref{prop.adapted}, the flows $z_1,\dots,z_{\nu_\vs}$ and the flows $z^C$ of the circles $C\ne\hat C$ form a basis of $\FF_{D,\hate}$. Let $\CU$ be the $c\times c$ integer matrix whose rows are the vectors $U(z)\in\BZ^c$ for $z$ in this basis. By \eqref{eq.piminusone}, $-\CU$ is the matrix $\bigl(\pi_j(z_i)\bigr)$ of Theorem~\ref{thm.flows} at $x_1=\dots=x_c=-1$; so
\[
   \det\CU=\pm P_{D,\hate}(-1,\dots,-1).
\]
Let $\CT$ be a spanning tree of the state graph $\Gamma_\vs$, let $\DD_\theta=\diag(\theta_1,\dots,\theta_c)$, and let $\CW$ be the $c\times c$ matrix whose columns are the vectors $\DD_\theta\,w^\vs(z_b)$, $b=1,\dots,\nu_\vs$, followed by the unit vectors $e_j$, $j\in\CT$. By Lemma~\ref{lem.stategraph}~(ii), $\det\CW=\pm1$. We compute the product $\CU\CW$. The entry in the row of a basis flow $z$ and in the column $b\le\nu_\vs$ is $\sum_jU_j(z)\,\theta_j\,w^\vs_j(z_b)=\GLf_\vs(z,z_b)$, by Proposition~\ref{prop.formsminusone}~(i); it is the entry $(a,b)$ of $G_\vs$ if $z=z_a$, and it is zero if $z=z^C$, because the circles lie in the radical of $\GLf_\vs$ (Proposition~\ref{prop.formsminusone}~(iv)). The entry in the row of $z$ and in the column of $j\in\CT$ is $U_j(z)$. So
\begin{equation}\label{eq.UWblock}
   \CU\,\CW=\begin{pmatrix}G_\vs&*\\0&\CU_{\CT}\end{pmatrix},\qquad \CU_{\CT}=\bigl(U_j(z^C)\bigr)_{C\ne\hat C,\ j\in\CT},
\end{equation}
and $\det\CU_{\CT}=\pm1$ by Lemma~\ref{lem.stategraph}~(iii).

(i) Taking determinants in \eqref{eq.UWblock}, $\det G_\vs=\pm\det\CU=\pm P_{D,\hate}(-1,\dots,-1)$. By Lemma~\ref{lem.alexpoly} at $t=-1$, where $x_j=t^{\ve_j}=-1$ for every crossing, the absolute value of this integer is $|\Delta_L(-1)|$.

(ii) The row space of $\CU$ is $U(\FF_{D,\hate})$. Now $\FF_D=\FF_{D,\hate}\oplus\BZ z^0$, because the coordinate of $z^0$ at $\hate$ is $\pm1$, and $U(z^0)=0$ by Proposition~\ref{prop.formsminusone}~(iv) and \eqref{eq.UOw}; so the row space of $\CU$ is $U(\FF_D)$. The quotient of $\BZ^c$ by the row space of a square integer matrix is isomorphic to the quotient by its column space, and its isomorphism class does not change when the matrix is multiplied on either side by an integer matrix of determinant $\pm1$ (both facts follow from the Smith normal form). Multiplying the matrix \eqref{eq.UWblock} on the right by $\bigl(\begin{smallmatrix}1&0\\0&\CU_{\CT}^{-1}\end{smallmatrix}\bigr)$ and then subtracting from its first $\nu_\vs$ rows combinations of the last ones gives $\bigl(\begin{smallmatrix}G_\vs&0\\0&1\end{smallmatrix}\bigr)$. Hence $\BZ^c/U(\FF_D)\cong\BZ^{\nu_\vs}/G_\vs\BZ^{\nu_\vs}$. For the oriented smoothing, $G_0$ is a matrix of $V_F+V_F^{\T}$ by Proposition~\ref{prop.formsminusone}~(vi).

(iii) By (ii), $\nu_\vs-\rk G_\vs$ is the rank of the abelian group $\BZ^c/U(\FF_D)$, which does not depend on $\vs$; for $\vs=0$ it is the corank $n(L)$ of $V_F+V_F^{\T}$.

(iv) By Proposition~\ref{prop.J}~(ii), the Gram matrix of $\GLf^\uparrow_\vs$ in the basis $[z_1],\dots,[z_{\nu_\vs}]$ is $\mathsf J^{\T}G_\vs\mathsf J$, where $\mathsf J$ is the matrix of the involution $J$ in this basis, and $\det\mathsf J=\pm1$.
\end{proof}

\begin{remark}\label{rem.Mminusone}
(i) \emph{The deformed Seifert matrices at $x=-1$.} For the oriented smoothing, \eqref{eq.UWblock} is the block decomposition of the proof of Corollary~\ref{cor.seifert}~(iii) at $x=-1$, written in a basis of the deformed cycle space, and the deformed Seifert matrices of Definition~\ref{def.deformedseifertmatrices} are given by Theorem~\ref{thm.alexpoint} at $t=-1$: for the basis of fundamental cycles of a spanning tree $\CT$ of the Seifert graph,
\[
   M_{D,\hate}(-1,\dots,-1)=-\bigl(V_F+V_F^{\T}\bigr)\DD'_\ve,\qquad M^b_{D,\hate}(-1,\dots,-1)=-\DD'_\ve\bigl(V_{F^\uparrow}+V_{F^\uparrow}^{\T}\bigr),
\]
where $\DD'_\ve=\diag(\ve_i)_{i\notin\CT}$, and the deformed cycle space at this point is $\DD_\ve H$, by Lemma~\ref{lem.Yt}. Theorem~\ref{thm.determinant} extends this to the other smoothings: at $x=-1$ the matrix of the crossing forms in a basis adapted to \emph{any} smoothing $\vs$ reduces, after the elimination of the circles, to the form of Gordon and Litherland of the state surface $F_\vs$. For general $x$ there is a matrix of the crossing forms in a basis adapted to $\vs$ as well (Proposition~\ref{prop.adapted}), but we do not know a surface interpretation of it when $\vs$ is not the oriented smoothing.

(ii) \emph{The second flavor.} The same computation with the forms $O_j$ in place of the $U_j$, and with the columns $-\DD_\ve\,w^\vs(z_b)$, expresses $\pm\tP_{D,\hate}(-1,\dots,-1)$ as the determinant of a Gram matrix of $\GLf^\uparrow_\vs$, by Proposition~\ref{prop.formsminusone}~(i) and Lemma~\ref{lem.stategraph}~(iii). So at $x=-1$, as at every other point (Theorem~\ref{thm.deformedseifert}), the first flavor belongs to the surfaces below the plane and the second one to the surfaces above it.
\end{remark}

\subsection{The signature}\label{ssec.signature}
We now compute the signatures of the forms of Gordon and Litherland of all the state surfaces. The method is to pass from one smoothing to another one by changing the smoothing of one crossing at a time: each step adds to the form the square of a linear form, with a sign, and changes its rank by exactly one, and this determines the change of the signature.

Let $\sigma_L=\sign(V+V^{\T})$ be the signature of the link $L$, where $V$ is a Seifert matrix of $L$ in the conventions of Section~\ref{ssec.localization}; it does not depend on the Seifert surface, and with these conventions the signature of the right-handed trefoil is $-2$ (Example~\ref{ex.trefoilseifert}).

\begin{theorem}[the signature of the state surfaces]\label{thm.signature}
For every smoothing $\vs$ of a diagram $D$ of a link $L$,
\[
   \sign\GLf_\vs=\sign\GLf^\uparrow_\vs=\sigma_L+\sum_{j:\,\vs_j=1}\ve_j .
\]
In other words, the signature of the form of Gordon and Litherland of the state surface $F_\vs$, or of $F^\uparrow_\vs$, minus the sum of the signs of the crossings which are not smoothed along the orientation, is the signature of the link.
\end{theorem}

The proof uses two lemmas. The first one is a general fact about quadratic forms, and the second one is the only place where the planarity of the diagram enters.

\begin{lemma}\label{lem.rankone}
Let $Q$ be a symmetric bilinear form on a real vector space of finite dimension, $\lambda$ a linear form, $\ve=\pm1$, and $Q'=Q+\ve\,\lambda\otimes\lambda$. If $\rk Q'=\rk Q\pm1$, then $\sign Q'=\sign Q+\ve$.
\end{lemma}
\begin{proof}
Replacing $Q,Q'$ by $-Q,-Q'$ if necessary, we may assume that $\ve=1$. Let $p$ and $q$ be the maximal dimensions of subspaces on which $Q$ is positive definite and negative definite, so that $\rk Q=p+q$ and $\sign Q=p-q$, and define $p',q'$ in the same way for $Q'$. Since $Q'(v,v)=Q(v,v)+\lambda(v)^2\ge Q(v,v)$, the form $Q'$ is positive definite wherever $Q$ is, and $Q$ is negative definite wherever $Q'$ is: $p'\ge p$ and $q\ge q'$. On the hyperplane $\ker\lambda$ the two forms coincide; so if $Q'$ is positive definite on a subspace $W$, then $Q$ is positive definite on $W\cap\ker\lambda$, whose dimension is at least $\dim W-1$, and $p\ge p'-1$; in the same way $q'\ge q-1$. Thus $p'-p\in\{0,1\}$ and $q'-q\in\{-1,0\}$. If the rank changes by exactly one, then $(p'-p)+(q'-q)=\pm1$, so that $(p'-p,q'-q)$ is $(1,0)$ or $(0,-1)$, and in both cases $(p'-q')-(p-q)=1$.
\end{proof}

\begin{figure}[ht]
\centering
\begin{tikzpicture}[>=Stealth,scale=0.74,every node/.style={font=\small}]
\def\discbase{\draw[gray!60] (0,0) circle (0.8);
  \foreach \a/\l in {45/1,135/2,225/3,315/4} {\fill (\a:0.8) circle (1.4pt); \node[font=\scriptsize] at (\a:0.42) {$p_\l$};}}
\def\arcsvs{\draw[very thick] (45:0.8) arc (45:135:0.8); \draw[very thick] (225:0.8) arc (225:315:0.8);}
\def\arcsvsp{\draw[very thick] (135:0.8) arc (135:225:0.8); \draw[very thick] (315:0.8) arc (315:405:0.8);}
\def\betas{\draw[thick] (135:0.8) .. controls (-2.6,2.7) and (2.6,2.7) .. (45:0.8); \draw[thick] (315:0.8) .. controls (2.6,-2.7) and (-2.6,-2.7) .. (225:0.8);
   \node[font=\scriptsize] at (0,2.25) {$\beta$}; \node[font=\scriptsize] at (0,-2.25) {$\beta'$};}
\begin{scope}
\betas \discbase \arcsvs \node at (0,-3.0) {$\vs$: two circles};
\end{scope}
\begin{scope}[xshift=4.6cm]
\betas \discbase \arcsvsp \node at (0,-3.0) {$\vs'$: one circle};
\end{scope}
\draw[<->] (1.75,0) -- (2.85,0);
\def\sides{\draw[thick] (135:0.8) .. controls (-2.7,2.6) and (-2.7,-2.6) .. (225:0.8); \draw[thick] (315:0.8) .. controls (2.7,-2.6) and (2.7,2.6) .. (45:0.8);}
\begin{scope}[xshift=9.7cm]
\sides \discbase \arcsvs \node at (0,-3.0) {$\vs$: one circle};
\end{scope}
\begin{scope}[xshift=14.5cm]
\sides \discbase \arcsvsp \node at (0,-3.0) {$\vs'$: two circles};
\end{scope}
\draw[<->] (11.6,0) -- (12.6,0);
\end{tikzpicture}
\caption{Changing the smoothing of one crossing $j$ (Lemma~\ref{lem.flipcircles}). The four end points $p_1,\dots,p_4$ of the half-strands lie on the circle $\partial\CB_j$; the arcs of $\vs$ at $j$ are $p_1p_2$ and $p_3p_4$, those of $\vs'$ are $p_2p_3$ and $p_4p_1$ (thick), and the rest of the circles lies outside $\CB_j$. Left: the two arcs of $\vs$ lie on two circles, which are merged into one circle of $\vs'$. Right: the two arcs of $\vs$ lie on one circle, which is split into two circles of $\vs'$; the two paths outside the disk cannot join $p_1$ to $p_3$ and $p_2$ to $p_4$, because they are disjoint.}
\label{fig.flip}
\end{figure}
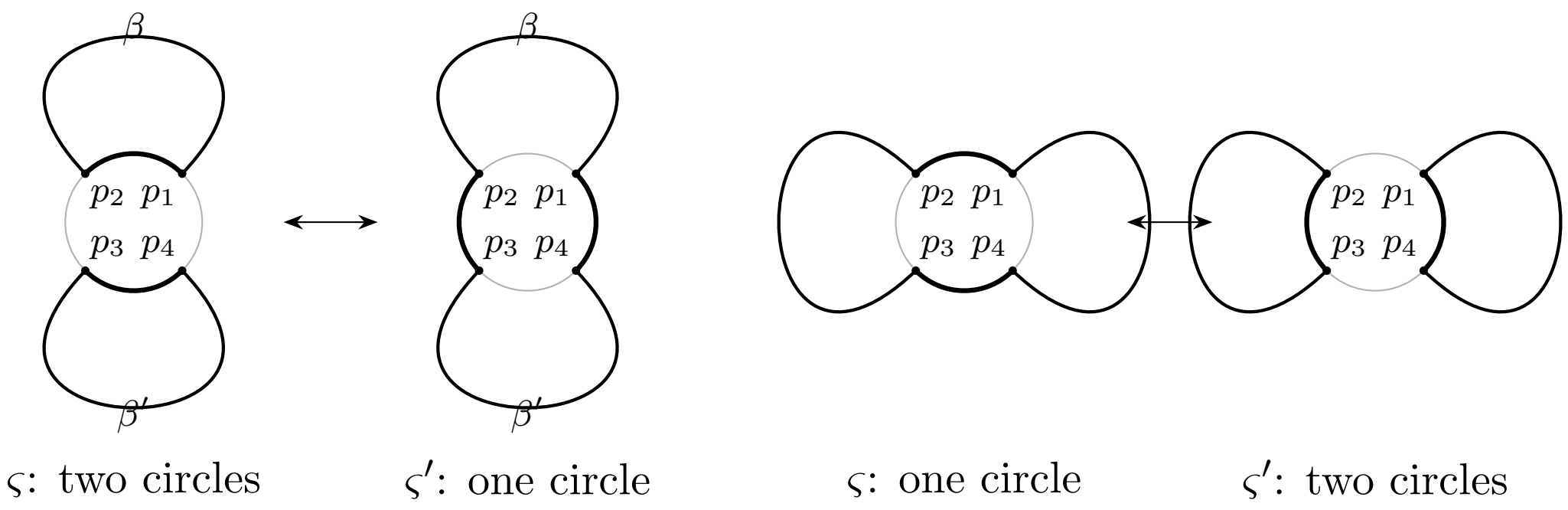

\begin{lemma}\label{lem.flipcircles}
Let $\vs$ and $\vs'$ be two smoothings of $D$ which differ at exactly one crossing. Then $s_{\vs'}=s_\vs\pm1$.
\end{lemma}
\begin{proof}
Let $j$ be the crossing at which they differ, and let $p_1,p_2,p_3,p_4$ be the four end points of the half-strands of $j$ on the circle $\partial\CB_j$, in their cyclic order. We may assume that the arcs of $\vs$ at $j$ are the arcs $p_1p_2$ and $p_3p_4$ of $\partial\CB_j$, and that those of $\vs'$ are $p_2p_3$ and $p_4p_1$ (Figure~\ref{fig.flip}). Outside the open disk $\CB_j^\circ$ the circles of $\vs$ and of $\vs'$ coincide. If the two arcs of $\vs$ lie on two different circles, these consist of the arc $p_1p_2$ and a path $\beta$ from $p_2$ to $p_1$ outside $\CB_j^\circ$, and of the arc $p_3p_4$ and a path $\beta'$ from $p_4$ to $p_3$; in $\vs'$ they are replaced by the single circle formed by $\beta$, the arc $p_1p_4$, $\beta'$ and the arc $p_3p_2$, and $s_{\vs'}=s_\vs-1$. If the two arcs lie on the same circle $C$, the rest of $C$ consists of two disjoint paths outside $\CB_j^\circ$ which join the four points in pairs. These paths lie in the complement of $\CB_j^\circ$ in the sphere, which is a closed disk with the four points on its boundary in the cyclic order $p_1,p_2,p_3,p_4$; so they cannot join $p_2$ to $p_4$ and $p_1$ to $p_3$, by the Jordan curve theorem, and since $C$ is connected they do not join $p_1$ to $p_2$. Hence they join $p_2$ to $p_3$ and $p_4$ to $p_1$, and in $\vs'$ the circle $C$ is replaced by two circles, each formed by one of these paths and one of the arcs $p_2p_3$, $p_4p_1$: $s_{\vs'}=s_\vs+1$.
\end{proof}

\begin{proof}[Proof of Theorem~\ref{thm.signature}]
We regard all the forms $\GLf_\vs$ as forms on the same real vector space $\FF_D\otimes\BR$; their signatures and their ranks are those of the forms induced on the homology of the state surfaces, since the circles lie in the radicals (Proposition~\ref{prop.formsminusone}~(iv)). Let $\vs$ be a smoothing with $\vs_j=0$ and let $\vs'$ be the smoothing which differs from $\vs$ at the crossing $j$ only. By Proposition~\ref{prop.formsminusone}~(v), $\GLf_{\vs'}=\GLf_\vs+\ve_j\,U_j\otimes U_j$. By Theorem~\ref{thm.determinant}~(iii) and Lemma~\ref{lem.flipcircles},
\[
   \rk\GLf_{\vs'}-\rk\GLf_\vs=s_\vs-s_{\vs'}=\pm1 .
\]
So Lemma~\ref{lem.rankone} gives $\sign\GLf_{\vs'}=\sign\GLf_\vs+\ve_j$. Starting from the oriented smoothing and changing the smoothings of the crossings one at a time, we obtain
\[
   \sign\GLf_\vs=\sign\GLf_0+\sum_{j:\,\vs_j=1}\ve_j\qquad\text{for every smoothing $\vs$.}
\] On $H_1(F_D;\BZ)$ the form $\GLf_0$ is $V_F+V_F^{\T}$ (Proposition~\ref{prop.formsminusone}~(vi)), whose signature is $\sigma_L$ because $F_D$ is a Seifert surface of $L$. Finally $\sign\GLf^\uparrow_\vs=\sign\GLf_\vs$, because the two forms are isometric (Proposition~\ref{prop.J}~(ii)).
\end{proof}

\begin{corollary}[Gordon and Litherland \cite{GordonLitherland}]\label{cor.GLgoeritz}
Let $D$ be a diagram of $L$ with a checkerboard coloring, let $v$ be one of the two colors, and let
\[
   \mu_v(D)=-\sum_j\eta_v(j)=\sum_j\ve_j ,
\]
both sums running over the crossings $j$ at which the two faces of color $v$ are $f^\leftarrow_j$ and $f^\rightarrow_j$. Then $\sign G_v(D)-\mu_v(D)=\sigma_L$.
\end{corollary}
\begin{proof}
The two sums agree by \eqref{eq.eta}. Take $v={\rm w}$. By Theorem~\ref{thm.goeritz}~(ii), the matrix of the form $\GLf_{\vs_{\rm b}}$ on the boundaries of all the faces is $G_{\rm w}(D)\oplus0$, so that $\sign G_{\rm w}(D)=\sign\GLf_{\vs_{\rm b}}$; and $(\vs_{\rm b})_j=1$ exactly when the black faces at $j$ are $f^\uparrow_j,f^\downarrow_j$, that is, when the white ones are $f^\leftarrow_j,f^\rightarrow_j$. So this is Theorem~\ref{thm.signature} for the smoothing $\vs_{\rm b}$; for $v={\rm b}$ exchange the two colors.
\end{proof}

This is the formula \cite[(3)]{CimasoniFerretti24}, in which $\mu_v(D)$ is written as the sum of $-\eta_v(c)$ over the crossings with $t_v(c)=-1$. Since every crossing has its faces $f^\leftarrow_j,f^\rightarrow_j$ of exactly one of the two colors,
\begin{equation}\label{eq.musum}
   \mu_{\rm w}(D)+\mu_{\rm b}(D)=\sum_{j=1}^c\ve_j=\mathrm{wr}(D)
\end{equation}
is the writhe of the diagram \cite[Lemma~3.3]{CimasoniFerretti24}.

\begin{corollary}[Cimasoni and Ferretti \cite{CimasoniFerretti24}]\label{cor.kashaevsignature}
For every diagram $D$ of a link $L$,
\[
   \sign\bigl(\tau_D|_{\zeta=0}\bigr)=\sign\Kash=2\,\sigma_L+\mathrm{wr}(D),
\]
where $\tau_D|_{\zeta=0}$ is Kashaev's matrix at $t=-1$. Moreover $\sign\GLf_\vs+\sign\GLf_{\bar\vs}=\sign\Kash$ for every smoothing $\vs$.
\end{corollary}
\begin{proof}
By Theorem~\ref{thm.goeritz}~(iii) and~(iv), $\sign(\tau_D|_{\zeta=0})=\sign\Kash=\sign G_{\rm w}(D)+\sign G_{\rm b}(D)$, which is $2\sigma_L+\mu_{\rm w}(D)+\mu_{\rm b}(D)$ by Corollary~\ref{cor.GLgoeritz}; then use \eqref{eq.musum}. By Theorem~\ref{thm.signature}, $\sign\GLf_\vs+\sign\GLf_{\bar\vs}=2\sigma_L+\sum_j\ve_j$, since every crossing has $\vs_j=1$ or $\bar\vs_j=1$.
\end{proof}

The first statement is \cite[Theorem~1~(ii)]{CimasoniFerretti24}: it is the case $\omega=-1$ of Kashaev's conjecture $\sign\tau_D-\mathrm{wr}(D)=2\sigma_L(\omega)$, which was proved in general by Liu \cite{Liu25}. The proof above follows the same path as that of Cimasoni and Ferretti, through the two Goeritz matrices; the differences are that the congruence of Theorem~\ref{thm.goeritz}~(iv) holds for every diagram, so that special diagrams are not needed, and that the formula of Corollary~\ref{cor.GLgoeritz} is proved here directly from the definition of the signature by the Seifert form of the canonical surface, with the lattice of flows in the place of the four-dimensional arguments of \cite{GordonLitherland}. The last statement shows that Kashaev's form plays the same role for every pair of complementary smoothings as for the two checkerboard smoothings, although it is a direct sum only for the latter.

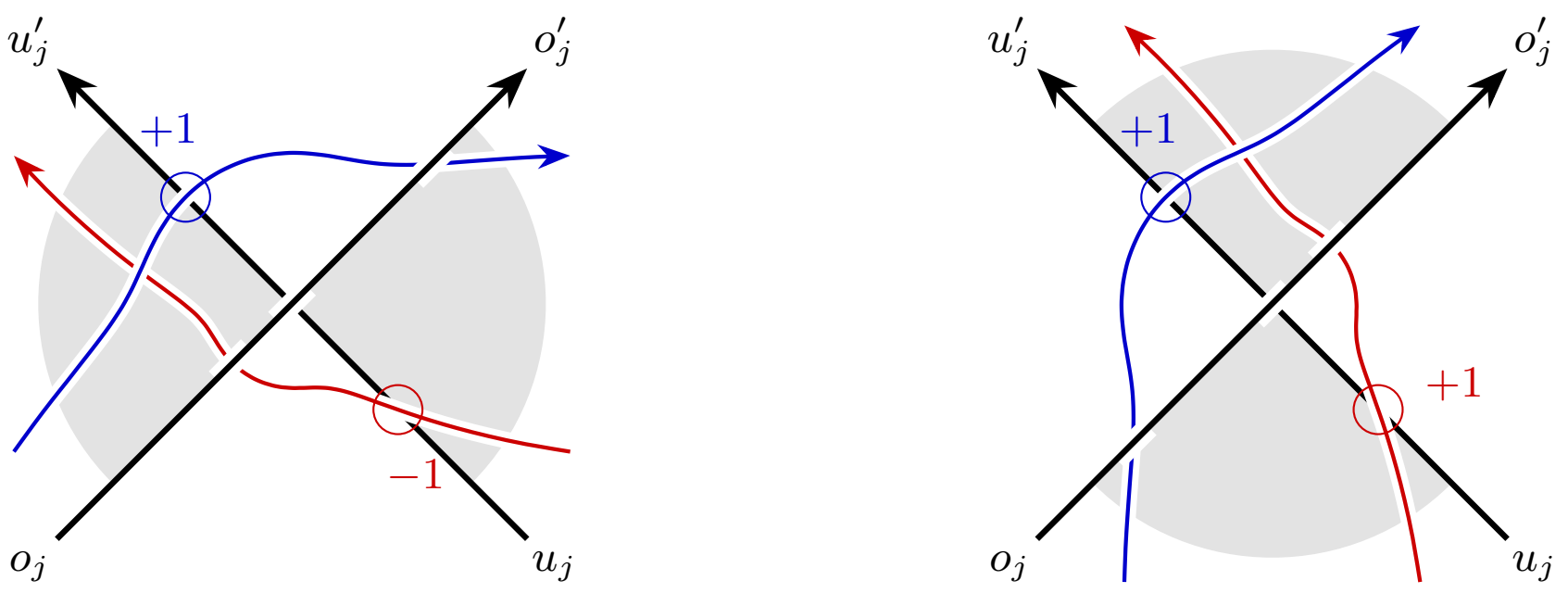
\begin{figure}[ht]
\centering
\begin{tikzpicture}[>=Stealth,scale=1.25,every node/.style={font=\small}]
\begin{scope}[xshift=0.00cm]
\fill[gray!22] (0,0) -- (135:1.45) arc (135:225:1.45) -- cycle; \fill[gray!22] (0,0) -- (-45:1.45) arc (-45:45:1.45) -- cycle;
\draw[very thick,->] (315:1.9) -- (135:1.9);
\draw[white,line width=3.4pt] plot[smooth] coordinates {(1.589,-0.845) (1.335,-0.802) (1.121,-0.756) (0.944,-0.711) (0.796,-0.668) (0.674,-0.629) (0.574,-0.594) (0.491,-0.565) (0.422,-0.540) (0.364,-0.520) (0.316,-0.505) (0.274,-0.494) (0.237,-0.487) (0.205,-0.482) (0.175,-0.480) (0.146,-0.479) (0.119,-0.479) (0.093,-0.480) (0.068,-0.481) (0.042,-0.481) (0.017,-0.482) (-0.008,-0.481) (-0.034,-0.479) (-0.059,-0.477) (-0.083,-0.473) (-0.108,-0.468) (-0.132,-0.461) (-0.156,-0.454) (-0.180,-0.445) (-0.203,-0.435) (-0.225,-0.424) (-0.247,-0.411) (-0.268,-0.398) (-0.289,-0.383) (-0.309,-0.368) (-0.327,-0.351) (-0.345,-0.333) (-0.362,-0.315) (-0.379,-0.296) (-0.394,-0.276) (-0.409,-0.255) (-0.423,-0.234) (-0.436,-0.213) (-0.450,-0.191) (-0.464,-0.169) (-0.479,-0.146) (-0.495,-0.123) (-0.514,-0.100) (-0.536,-0.075) (-0.563,-0.049) (-0.596,-0.021) (-0.635,0.011) (-0.684,0.048) (-0.742,0.091) (-0.813,0.143) (-0.898,0.207) (-0.999,0.286) (-1.117,0.385) (-1.254,0.507) (-1.411,0.658) (-1.589,0.845)};
\draw[red!80!black,thick,->] plot[smooth] coordinates {(1.589,-0.845) (1.335,-0.802) (1.121,-0.756) (0.944,-0.711) (0.796,-0.668) (0.674,-0.629) (0.574,-0.594) (0.491,-0.565) (0.422,-0.540) (0.364,-0.520) (0.316,-0.505) (0.274,-0.494) (0.237,-0.487) (0.205,-0.482) (0.175,-0.480) (0.146,-0.479) (0.119,-0.479) (0.093,-0.480) (0.068,-0.481) (0.042,-0.481) (0.017,-0.482) (-0.008,-0.481) (-0.034,-0.479) (-0.059,-0.477) (-0.083,-0.473) (-0.108,-0.468) (-0.132,-0.461) (-0.156,-0.454) (-0.180,-0.445) (-0.203,-0.435) (-0.225,-0.424) (-0.247,-0.411) (-0.268,-0.398) (-0.289,-0.383) (-0.309,-0.368) (-0.327,-0.351) (-0.345,-0.333) (-0.362,-0.315) (-0.379,-0.296) (-0.394,-0.276) (-0.409,-0.255) (-0.423,-0.234) (-0.436,-0.213) (-0.450,-0.191) (-0.464,-0.169) (-0.479,-0.146) (-0.495,-0.123) (-0.514,-0.100) (-0.536,-0.075) (-0.563,-0.049) (-0.596,-0.021) (-0.635,0.011) (-0.684,0.048) (-0.742,0.091) (-0.813,0.143) (-0.898,0.207) (-0.999,0.286) (-1.117,0.385) (-1.254,0.507) (-1.411,0.658) (-1.589,0.845)};
\draw[white,line width=3.4pt] plot[smooth] coordinates {(-1.589,-0.845) (-1.475,-0.688) (-1.374,-0.555) (-1.285,-0.443) (-1.210,-0.347) (-1.145,-0.264) (-1.090,-0.192) (-1.043,-0.128) (-1.004,-0.070) (-0.970,-0.017) (-0.942,0.033) (-0.917,0.080) (-0.895,0.126) (-0.875,0.170) (-0.855,0.213) (-0.836,0.256) (-0.817,0.297) (-0.797,0.338) (-0.776,0.379) (-0.754,0.418) (-0.730,0.456) (-0.705,0.494) (-0.678,0.530) (-0.649,0.564) (-0.619,0.597) (-0.587,0.629) (-0.553,0.659) (-0.518,0.687) (-0.481,0.713) (-0.443,0.737) (-0.404,0.759) (-0.363,0.779) (-0.322,0.797) (-0.280,0.813) (-0.237,0.827) (-0.193,0.838) (-0.149,0.847) (-0.105,0.854) (-0.060,0.858) (-0.015,0.861) (0.030,0.861) (0.075,0.859) (0.120,0.855) (0.165,0.850) (0.210,0.844) (0.256,0.836) (0.302,0.828) (0.348,0.820) (0.396,0.812) (0.446,0.805) (0.499,0.799) (0.557,0.795) (0.619,0.793) (0.689,0.793) (0.769,0.796) (0.859,0.801) (0.964,0.809) (1.086,0.818) (1.228,0.828) (1.395,0.838) (1.589,0.845)};
\draw[blue!80!black,thick,->] plot[smooth] coordinates {(-1.589,-0.845) (-1.475,-0.688) (-1.374,-0.555) (-1.285,-0.443) (-1.210,-0.347) (-1.145,-0.264) (-1.090,-0.192) (-1.043,-0.128) (-1.004,-0.070) (-0.970,-0.017) (-0.942,0.033) (-0.917,0.080) (-0.895,0.126) (-0.875,0.170) (-0.855,0.213) (-0.836,0.256) (-0.817,0.297) (-0.797,0.338) (-0.776,0.379) (-0.754,0.418) (-0.730,0.456) (-0.705,0.494) (-0.678,0.530) (-0.649,0.564) (-0.619,0.597) (-0.587,0.629) (-0.553,0.659) (-0.518,0.687) (-0.481,0.713) (-0.443,0.737) (-0.404,0.759) (-0.363,0.779) (-0.322,0.797) (-0.280,0.813) (-0.237,0.827) (-0.193,0.838) (-0.149,0.847) (-0.105,0.854) (-0.060,0.858) (-0.015,0.861) (0.030,0.861) (0.075,0.859) (0.120,0.855) (0.165,0.850) (0.210,0.844) (0.256,0.836) (0.302,0.828) (0.348,0.820) (0.396,0.812) (0.446,0.805) (0.499,0.799) (0.557,0.795) (0.619,0.793) (0.689,0.793) (0.769,0.796) (0.859,0.801) (0.964,0.809) (1.086,0.818) (1.228,0.828) (1.395,0.838) (1.589,0.845)};
\draw[white,line width=4.5pt] (45:-0.130) -- (45:0.130);
\draw[white,line width=4.5pt] (225:0.350) -- (225:0.610);
\draw[white,line width=4.5pt] (45:0.998) -- (45:1.258);
\draw[very thick,->] (225:1.9) -- (45:1.9);
\node[font=\scriptsize] at (225:2.12) {$o_j$}; \node[font=\scriptsize] at (45:2.12) {$o'_j$};
\node[font=\scriptsize] at (315:2.12) {$u_j$}; \node[font=\scriptsize] at (135:2.12) {$u'_j$};
\draw[red!80!black] (0.605,-0.605) circle (0.14);
\draw[blue!80!black] (-0.608,0.608) circle (0.14);
\node[font=\scriptsize,red!80!black] at (0.705,-0.985) {$-1$};
\node[font=\scriptsize,blue!80!black] at (-0.708,0.988) {$+1$};
\node at (0,-2.75) {$\vs_j=0$, $\theta_j=1$};
\end{scope}
\begin{scope}[xshift=5.60cm]
\fill[gray!22] (0,0) -- (45:1.45) arc (45:135:1.45) -- cycle; \fill[gray!22] (0,0) -- (225:1.45) arc (225:315:1.45) -- cycle;
\draw[very thick,->] (315:1.9) -- (135:1.9);
\draw[white,line width=3.4pt] plot[smooth] coordinates {(0.845,-1.589) (0.802,-1.335) (0.756,-1.121) (0.711,-0.944) (0.668,-0.796) (0.629,-0.674) (0.594,-0.574) (0.565,-0.491) (0.540,-0.422) (0.520,-0.364) (0.505,-0.316) (0.494,-0.274) (0.487,-0.237) (0.482,-0.205) (0.480,-0.175) (0.479,-0.146) (0.479,-0.119) (0.480,-0.093) (0.481,-0.068) (0.481,-0.042) (0.482,-0.017) (0.481,0.008) (0.479,0.034) (0.477,0.059) (0.473,0.083) (0.468,0.108) (0.461,0.132) (0.454,0.156) (0.445,0.180) (0.435,0.203) (0.424,0.225) (0.411,0.247) (0.398,0.268) (0.383,0.289) (0.368,0.309) (0.351,0.327) (0.333,0.345) (0.315,0.362) (0.296,0.379) (0.276,0.394) (0.255,0.409) (0.234,0.423) (0.213,0.436) (0.191,0.450) (0.169,0.464) (0.146,0.479) (0.123,0.495) (0.100,0.514) (0.075,0.536) (0.049,0.563) (0.021,0.596) (-0.011,0.635) (-0.048,0.684) (-0.091,0.742) (-0.143,0.813) (-0.207,0.898) (-0.286,0.999) (-0.385,1.117) (-0.507,1.254) (-0.658,1.411) (-0.845,1.589)};
\draw[red!80!black,thick,->] plot[smooth] coordinates {(0.845,-1.589) (0.802,-1.335) (0.756,-1.121) (0.711,-0.944) (0.668,-0.796) (0.629,-0.674) (0.594,-0.574) (0.565,-0.491) (0.540,-0.422) (0.520,-0.364) (0.505,-0.316) (0.494,-0.274) (0.487,-0.237) (0.482,-0.205) (0.480,-0.175) (0.479,-0.146) (0.479,-0.119) (0.480,-0.093) (0.481,-0.068) (0.481,-0.042) (0.482,-0.017) (0.481,0.008) (0.479,0.034) (0.477,0.059) (0.473,0.083) (0.468,0.108) (0.461,0.132) (0.454,0.156) (0.445,0.180) (0.435,0.203) (0.424,0.225) (0.411,0.247) (0.398,0.268) (0.383,0.289) (0.368,0.309) (0.351,0.327) (0.333,0.345) (0.315,0.362) (0.296,0.379) (0.276,0.394) (0.255,0.409) (0.234,0.423) (0.213,0.436) (0.191,0.450) (0.169,0.464) (0.146,0.479) (0.123,0.495) (0.100,0.514) (0.075,0.536) (0.049,0.563) (0.021,0.596) (-0.011,0.635) (-0.048,0.684) (-0.091,0.742) (-0.143,0.813) (-0.207,0.898) (-0.286,0.999) (-0.385,1.117) (-0.507,1.254) (-0.658,1.411) (-0.845,1.589)};
\draw[white,line width=3.4pt] plot[smooth] coordinates {(-0.845,-1.589) (-0.838,-1.395) (-0.828,-1.228) (-0.818,-1.086) (-0.809,-0.964) (-0.801,-0.859) (-0.796,-0.769) (-0.793,-0.689) (-0.793,-0.619) (-0.795,-0.557) (-0.799,-0.499) (-0.805,-0.446) (-0.812,-0.396) (-0.820,-0.348) (-0.828,-0.302) (-0.836,-0.256) (-0.844,-0.210) (-0.850,-0.165) (-0.855,-0.120) (-0.859,-0.075) (-0.861,-0.030) (-0.861,0.015) (-0.858,0.060) (-0.854,0.105) (-0.847,0.149) (-0.838,0.193) (-0.827,0.237) (-0.813,0.280) (-0.797,0.322) (-0.779,0.363) (-0.759,0.404) (-0.737,0.443) (-0.713,0.481) (-0.687,0.518) (-0.659,0.553) (-0.629,0.587) (-0.597,0.619) (-0.564,0.649) (-0.530,0.678) (-0.494,0.705) (-0.456,0.730) (-0.418,0.754) (-0.379,0.776) (-0.338,0.797) (-0.297,0.817) (-0.256,0.836) (-0.213,0.855) (-0.170,0.875) (-0.126,0.895) (-0.080,0.917) (-0.033,0.942) (0.017,0.970) (0.070,1.004) (0.128,1.043) (0.192,1.090) (0.264,1.145) (0.347,1.210) (0.443,1.285) (0.555,1.374) (0.688,1.475) (0.845,1.589)};
\draw[blue!80!black,thick,->] plot[smooth] coordinates {(-0.845,-1.589) (-0.838,-1.395) (-0.828,-1.228) (-0.818,-1.086) (-0.809,-0.964) (-0.801,-0.859) (-0.796,-0.769) (-0.793,-0.689) (-0.793,-0.619) (-0.795,-0.557) (-0.799,-0.499) (-0.805,-0.446) (-0.812,-0.396) (-0.820,-0.348) (-0.828,-0.302) (-0.836,-0.256) (-0.844,-0.210) (-0.850,-0.165) (-0.855,-0.120) (-0.859,-0.075) (-0.861,-0.030) (-0.861,0.015) (-0.858,0.060) (-0.854,0.105) (-0.847,0.149) (-0.838,0.193) (-0.827,0.237) (-0.813,0.280) (-0.797,0.322) (-0.779,0.363) (-0.759,0.404) (-0.737,0.443) (-0.713,0.481) (-0.687,0.518) (-0.659,0.553) (-0.629,0.587) (-0.597,0.619) (-0.564,0.649) (-0.530,0.678) (-0.494,0.705) (-0.456,0.730) (-0.418,0.754) (-0.379,0.776) (-0.338,0.797) (-0.297,0.817) (-0.256,0.836) (-0.213,0.855) (-0.170,0.875) (-0.126,0.895) (-0.080,0.917) (-0.033,0.942) (0.017,0.970) (0.070,1.004) (0.128,1.043) (0.192,1.090) (0.264,1.145) (0.347,1.210) (0.443,1.285) (0.555,1.374) (0.688,1.475) (0.845,1.589)};
\draw[white,line width=4.5pt] (45:-0.130) -- (45:0.130);
\draw[white,line width=4.5pt] (45:0.350) -- (45:0.610);
\draw[white,line width=4.5pt] (225:0.998) -- (225:1.258);
\draw[very thick,->] (225:1.9) -- (45:1.9);
\node[font=\scriptsize] at (225:2.12) {$o_j$}; \node[font=\scriptsize] at (45:2.12) {$o'_j$};
\node[font=\scriptsize] at (315:2.12) {$u_j$}; \node[font=\scriptsize] at (135:2.12) {$u'_j$};
\draw[red!80!black] (0.605,-0.605) circle (0.14);
\draw[blue!80!black] (-0.608,0.608) circle (0.14);
\node[font=\scriptsize,red!80!black] at (1.040,-0.460) {$+1$};
\node[font=\scriptsize,blue!80!black] at (-0.708,0.988) {$+1$};
\node at (0,-2.75) {$\vs_j=1$, $\theta_j=-1$};
\end{scope}
\end{tikzpicture}
\caption{The link $L$ (black) and its push-off $L'$ into the state surface near a positive crossing $j$, for the two smoothings of $j$ (Remark~\ref{rem.euler}); the covered quadrants are shaded. The push-off of the under-strand (red) and that of the over-strand (blue) lie over the covered quadrants; they have been moved slightly away from the co-core, around which they turn, so that the crossings of the picture are distinct. Each of them passes once over the under-strand (circles) and once under the over-strand. The push-off of the over-strand crosses over the under-strand with the sign $\ve_j=1$, and the push-off of the under-strand with the sign $-\theta_j$: the two signs cancel for the oriented smoothing, and they add up to $2\ve_j$ for the other one.}
\label{fig.euler}
\end{figure}

\begin{remark}[the Euler number]\label{rem.euler}
Gordon and Litherland prove that, for every spanning surface $F$ of $L$, the signature of the link is the signature of the form $\CG_F$ plus a correction term which depends only on the way the surface approaches its boundary: in our conventions it is $-\frac12\lk(L,L')$, where $L'$ is the link $L$ pushed into the surface $F$ and oriented parallel to $L$; this is one half of the Euler number of the surface in \cite{GordonLitherland}. For the state surfaces Theorem~\ref{thm.signature} is this statement, with
\[
   \tfrac12\,\lk(L,L')=\sum_{j:\,\vs_j=1}\ve_j .
\]
This value can be seen directly in the model of Section~\ref{ssec.statesurface} (Figure~\ref{fig.euler}). The curve $L'$ lies on the surface and runs parallel to $L$; along the under-strand of a crossing $j$ it lies above the plane, over the covered quadrant next to $u_j$ before the crossing and over the covered quadrant next to $u'_j$ after it, and it passes from one to the other along the co-core, slightly above its lower end; along the over-strand it lies slightly below the over-strand, and it passes below its highest point. So in the projection $L'$ passes over $L$ twice near each crossing: the push-off of the over-strand passes over the under-strand, which is a crossing of sign $\ve_j$; and the push-off of the under-strand passes over the under-strand itself, from its right side to its left side if $\theta_j=1$ and from its left side to its right side if $\theta_j=-1$, by \eqref{eq.thetaside}, which is a crossing of sign $-\theta_j$. With the convention of Section~\ref{ssec.spanning}, $\lk(L,L')=\sum_j(\ve_j-\theta_j)=2\sum_{\vs_j=1}\ve_j$.
\end{remark}

\subsection{Alternating diagrams}\label{ssec.altminusone}
A diagram is \emph{alternating} if every edge is the over-strand at one of its ends and the under-strand at the other one, that is, if every edge is $o'_j=u_k$ or $u'_j=o_k$ for some crossings $j,k$ (Section~\ref{ssec.alternating}). For such a diagram the checkerboard coloring can be read from the signs of the crossings, and the forms of the two checkerboard surfaces are definite.

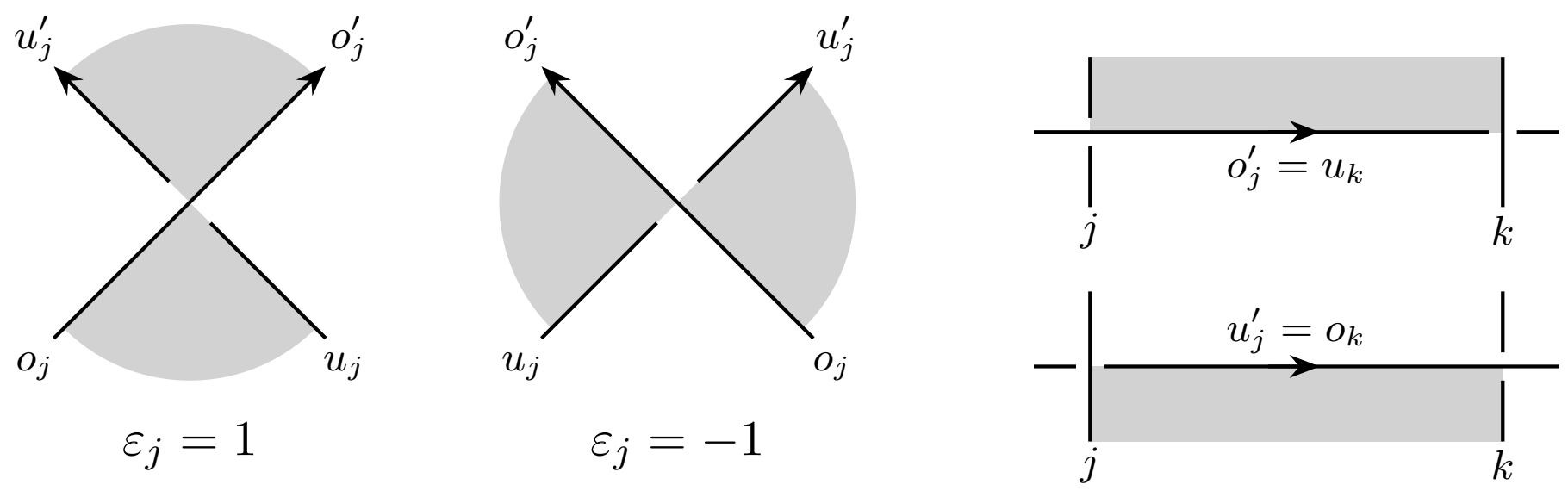
\begin{figure}[ht]
\centering
\begin{tikzpicture}[>=Stealth,scale=0.72,every node/.style={font=\small}]
\def\shUD{\fill[gray!35] (0,0) -- (45:1.9) arc (45:135:1.9) -- cycle; \fill[gray!35] (0,0) -- (225:1.9) arc (225:315:1.9) -- cycle;}
\def\shLR{\fill[gray!35] (0,0) -- (135:1.9) arc (135:225:1.9) -- cycle; \fill[gray!35] (0,0) -- (-45:1.9) arc (-45:45:1.9) -- cycle;}
\begin{scope}
\shUD
\draw[thick] (1.45,-1.45) -- (0.22,-0.22); \draw[thick,->] (-0.22,0.22) -- (-1.45,1.45); \draw[thick,->] (-1.45,-1.45) -- (1.45,1.45);
\node[font=\scriptsize] at (-1.65,-1.75) {$o_j$}; \node[font=\scriptsize] at (1.7,1.75) {$o'_j$};
\node[font=\scriptsize] at (1.65,-1.75) {$u_j$}; \node[font=\scriptsize] at (-1.65,1.75) {$u'_j$};
\node at (0,-2.6) {$\ve_j=1$};
\end{scope}
\begin{scope}[xshift=5.2cm]
\shLR
\draw[thick] (-1.45,-1.45) -- (-0.22,-0.22); \draw[thick,->] (0.22,0.22) -- (1.45,1.45); \draw[thick,->] (1.45,-1.45) -- (-1.45,1.45);
\node[font=\scriptsize] at (1.65,-1.75) {$o_j$}; \node[font=\scriptsize] at (-1.65,1.75) {$o'_j$};
\node[font=\scriptsize] at (-1.65,-1.75) {$u_j$}; \node[font=\scriptsize] at (1.7,1.75) {$u'_j$};
\node at (0,-2.6) {$\ve_j=-1$};
\end{scope}
\begin{scope}[xshift=9.6cm,yshift=0.75cm]
\fill[gray!35] (0,0) rectangle (4.4,0.8);
\draw[thick] (0,-0.8) -- (0,-0.15); \draw[thick] (0,0.15) -- (0,0.8);
\draw[thick] (-0.6,0) -- (4.25,0); \draw[thick] (4.55,0) -- (5.0,0);
\draw[thick] (4.4,-0.8) -- (4.4,0.8);
\draw[thick,->] (1.9,0) -- (2.45,0);
\node[font=\scriptsize] at (2.2,-0.38) {$o'_j=u_k$}; \node[font=\scriptsize] at (0,-1.05) {$j$}; \node[font=\scriptsize] at (4.4,-1.05) {$k$};
\end{scope}
\begin{scope}[xshift=9.6cm,yshift=-1.75cm]
\fill[gray!35] (0,-0.8) rectangle (4.4,0);
\draw[thick] (0,-0.8) -- (0,0.8);
\draw[thick] (-0.6,0) -- (-0.15,0); \draw[thick] (0.15,0) -- (5.0,0);
\draw[thick] (4.4,-0.8) -- (4.4,-0.15); \draw[thick] (4.4,0.15) -- (4.4,0.8);
\draw[thick,->] (1.9,0) -- (2.45,0);
\node[font=\scriptsize] at (2.2,0.38) {$u'_j=o_k$}; \node[font=\scriptsize] at (0,-1.05) {$j$}; \node[font=\scriptsize] at (4.4,-1.05) {$k$};
\end{scope}
\end{tikzpicture}
\caption{The coloring of an alternating diagram (Lemma~\ref{lem.altcoloring}). Left: the black corners of a positive and of a negative crossing. Right: an edge $o'_j=u_k$, which leaves the crossing $j$ as its over-strand and arrives at the crossing $k$ as its under-strand, has the black face on its left; an edge $u'_j=o_k$ has the black face on its right.}
\label{fig.altcolor}
\end{figure}

\begin{lemma}[the coloring of an alternating diagram]\label{lem.altcoloring}
Let $D$ be an alternating diagram. Call a corner of a crossing $j$ \emph{black} if it is one of the corners of $f^\uparrow_j,f^\downarrow_j$ and $\ve_j=1$, or one of the corners of $f^\leftarrow_j,f^\rightarrow_j$ and $\ve_j=-1$, and \emph{white} otherwise. Then all the corners of a face have the same color, and the two faces adjacent to an edge have different colors: this defines a checkerboard coloring of $D$, in which a face is black if and only if it lies on the left of an edge $o'_j=u_k$ or on the right of an edge $u'_j=o_k$ \textup(Figure~\ref{fig.altcolor}\textup). For this coloring,
\[
   \theta_j(\vs_{\rm b})=-1\qquad\text{and}\qquad\theta_j(\vs_{\rm w})=1\qquad\text{for every crossing $j$.}
\]
\end{lemma}
\begin{proof}
By \eqref{eq.mcrossing}, for both signs of the crossing $j$, the corners on the left of $o'_j$ and of $u_j$ are black and the corners on their right are white, while the corners on the left of $u'_j$ and of $o_j$ are white and the corners on their right are black. Since $D$ is alternating, every edge is $o'_j=u_k$ or $u'_j=o_k$ for some $j,k$; so the two corners on the left of an edge, at its two ends, have the same color, and the two corners on its right have the other color. Two consecutive corners of a face are adjacent to the same side of the edge between them; hence all the corners of a face have the same color, the two faces adjacent to an edge have different colors, and the black faces are those which lie on the left of an edge $o'_j=u_k$ or on the right of an edge $u'_j=o_k$. The covered quadrants of $\vs_{\rm b}$ are the black corners: $(\vs_{\rm b})_j=1$ if $\ve_j=1$ and $(\vs_{\rm b})_j=0$ if $\ve_j=-1$, so that $\theta_j(\vs_{\rm b})=\ve_j(-1)^{(\vs_{\rm b})_j}=-1$ in both cases, and $\theta_j(\vs_{\rm w})=-\theta_j(\vs_{\rm b})$.
\end{proof}

\begin{corollary}[alternating diagrams]\label{cor.altdefinite}
Let $D$ be an alternating diagram of a link $L$, colored as in Lemma~\ref{lem.altcoloring}, with $F_{\rm b}$ black faces, $F_{\rm w}$ white faces, $c_+$ positive and $c_-$ negative crossings.
\begin{enumerate}\renewcommand{\theenumi}{\roman{enumi}}
\item For all flows $z,z'$,
\[
   \GLf_{\vs_{\rm b}}(z,z')=\sum_{j=1}^cw^{\vs_{\rm b}}_j(z)\,w^{\vs_{\rm b}}_j(z')\qquad\text{and}\qquad\GLf_{\vs_{\rm w}}(z,z')=-\sum_{j=1}^cw^{\vs_{\rm w}}_j(z)\,w^{\vs_{\rm w}}_j(z').
\]
So the form of Gordon and Litherland of the black checkerboard surface is positive definite, of rank $c+1-F_{\rm b}$, and that of the white one is negative definite, of rank $c+1-F_{\rm w}$: by means of the transfers they are the restriction of the standard scalar product of $\BZ^c$ to the cycle space $Z_1(\Gamma_{\vs_{\rm b}})$, and the restriction of its opposite to $Z_1(\Gamma_{\vs_{\rm w}})$.
\item $\Delta_L(-1)\ne0$, and
\[
   \sigma_L=c_-+1-F_{\rm b}=F_{\rm w}-1-c_+ .
\]
\end{enumerate}
\end{corollary}
\begin{proof}
(i) By Theorem~\ref{thm.goeritz}~(ii), $\GLf_{\vs_{\rm b}}=\Goe_{\vs_{\rm b}}=-\sum_j\theta_j(\vs_{\rm b})\,w^{\vs_{\rm b}}_j\otimes w^{\vs_{\rm b}}_j$, and all the signs are $-1$ by Lemma~\ref{lem.altcoloring}. The map $w^{\vs_{\rm b}}$ is injective on $H_1(F_{\vs_{\rm b}};\BZ)$, with image $Z_1(\Gamma_{\vs_{\rm b}})$ (Lemma~\ref{lem.stategraph}~(i)); so the form is positive definite on this group, whose rank is $c+1-F_{\rm b}$ because the circles of $\vs_{\rm b}$ correspond to the black faces. For the white surface, $\GLf_{\vs_{\rm w}}=\Goe_{\vs_{\rm w}}=-\sum_j\theta_j(\vs_{\rm w})\,w^{\vs_{\rm w}}_j\otimes w^{\vs_{\rm w}}_j$ by the same theorem with the two colors exchanged, and now all the signs are $1$; the map $w^{\vs_{\rm w}}$ is injective on $H_1(F_{\vs_{\rm w}};\BZ)$, with image $Z_1(\Gamma_{\vs_{\rm w}})$, and the circles of $\vs_{\rm w}$ correspond to the white faces. So this form is negative definite, of rank $c+1-F_{\rm w}$.

(ii) A definite form has full rank, so $n(L)=0$ by Theorem~\ref{thm.determinant}~(iii), that is, $\Delta_L(-1)\ne0$. By Theorem~\ref{thm.signature} for the smoothing $\vs_{\rm b}$, whose crossings with $(\vs_{\rm b})_j=1$ are the positive ones, $\sigma_L=\sign\GLf_{\vs_{\rm b}}-c_+=(c+1-F_{\rm b})-c_+=c_-+1-F_{\rm b}$, the signature of a positive definite form being its rank. The second expression for $\sigma_L$ follows from $F_{\rm b}+F_{\rm w}=c+2$, which is Euler's formula for the diagram.
\end{proof}

Part (ii) is Traczyk's formula for the signature of an alternating link \cite[Theorem~2]{Traczyk04}, which reads $\sigma_L=-\frac12(c_+-c_-)+\frac12(F_{\rm w}-F_{\rm b})$ in our conventions for the colors and for the sign of the signature; the diagram need not be reduced. We shall use it in Section~\ref{ssec.statestrees}, where the number $F_{\rm b}-1-c_-=-\sigma_L$ counts the vertical corners of every Kauffman state of an alternating diagram.

\subsection{The lattice of the under-strands}\label{ssec.underlattice}
By Theorem~\ref{thm.determinant}, the forms of all the state surfaces are governed by one sublattice of $\BZ^c$, the image
\[
   \Lambda_D=U(\FF_D)\subseteq\BZ^c
\]
of the lattice of flows under the map $U=(U_1,\dots,U_c)$, $U_j(z)=z_{u'_j}+z_{u_j}$: the lattice of the vectors which record, crossing by crossing, the flow along the under-strands. We describe it in the basis of the faces. Fix a checkerboard coloring, let $\Gamma_{\rm b}=\Gamma_{\vs_{\rm b}}$ be the state graph of the black checkerboard smoothing, and write $\theta_j=\theta_j(\vs_{\rm b})$ and $\DD_\theta=\diag(\theta_1,\dots,\theta_c)$. The vertices of $\Gamma_{\rm b}$ are the black faces, and the edge of the crossing $j$ goes from the black face $c'_j$ at the corner next to $u'_j$ to the black face $c_j$ at the corner next to $u_j$ (Lemma~\ref{lem.cornertransfers}); this is the \emph{Tait graph} of the black faces. Besides its cycle space $Z_1(\Gamma_{\rm b})\subseteq\BZ^c$ we consider its \emph{lattice of cuts} $B^1(\Gamma_{\rm b})\subseteq\BZ^c$, which is generated by the vectors
\[
   \delta(f)=\sum_{j=1}^c\bigl(\ind{c_j=f}-\ind{c'_j=f}\bigr)\,e_j,\qquad f\text{ a black face;}
\]
it has rank $F_{\rm b}-1$, where $F_{\rm b}$ is the number of black faces, and it is the orthogonal complement of $Z_1(\Gamma_{\rm b})$ for the standard scalar product of $\BZ^c$.

\begin{proposition}[the lattice of the under-strands]\label{prop.underlattice}
Let $D$ be a diagram of $L$ with a checkerboard coloring.
\begin{enumerate}\renewcommand{\theenumi}{\roman{enumi}}
\item $U(\partial g)=-w^{\vs_{\rm b}}(\partial g)$ for every white face $g$, and these vectors generate $Z_1(\Gamma_{\rm b})$; $U(\partial f)=-\DD_\theta\,\delta(f)$ for every black face $f$. Hence
\[
   \Lambda_D=Z_1(\Gamma_{\rm b})+\DD_\theta\,B^1(\Gamma_{\rm b}).
\]
\item The kernel of $U$ contains the flow $z^0$ of Proposition~\ref{prop.formsminusone}~\textup{(iv)}, which is, up to sign, the sum of the boundaries of the faces of one color. The group $\BZ^c/\Lambda_D$ is isomorphic to the first homology group of the double covering of $S^3$ branched along $L$. If $\Delta_L(-1)\ne0$, then $\ker U=\BZ z^0$, the sum in \textup{(i)} is direct, and $\Lambda_D$ has the index $|\Delta_L(-1)|$ in $\BZ^c$.
\item If $D$ is alternating and colored as in Lemma~\ref{lem.altcoloring}, then $\DD_\theta=-1$ and
\[
   \Lambda_D=Z_1(\Gamma_{\rm b})\oplus B^1(\Gamma_{\rm b})
\]
is the sum of the lattice of the integral flows and of the lattice of the integral cuts of the Tait graph.
\end{enumerate}
\end{proposition}
\begin{proof}
(i) By \eqref{eq.UOw}, $U_j=-w^{\vs_{\rm b}}_j-w^{\vs_{\rm w}}_j$. Let $g$ be a white face. Its boundary is, up to sign, a circle of $\vs_{\rm w}$, so $w^{\vs_{\rm w}}(\partial g)=0$ by Lemma~\ref{lem.circles}~(i) and $U(\partial g)=-w^{\vs_{\rm b}}(\partial g)$. The classes of the boundaries of the white faces generate $H_1(F_{\vs_{\rm b}};\BZ)$, because the boundaries of all the faces generate $\FF_D$ and those of the black faces are circles of $\vs_{\rm b}$; so their images under $w^{\vs_{\rm b}}$ generate $Z_1(\Gamma_{\rm b})$ by Lemma~\ref{lem.stategraph}~(i). Let $f$ be a black face. Then $w^{\vs_{\rm b}}(\partial f)=0$, and by Lemma~\ref{lem.cornertransfers}~(ii) for the smoothing $\vs_{\rm b}$ and the potential $e_f$ of $\partial f$,
\[
   U_j(\partial f)=-w^{\vs_{\rm w}}_j(\partial f)=\theta_j\bigl(e_f(c'_j)-e_f(c_j)\bigr)=-\theta_j\,\delta(f)_j .
\]
Since the boundaries of the faces generate $\FF_D$, the lattice $\Lambda_D$ is generated by these two families of vectors.

(ii) We have $U(z^0)=0$ by Proposition~\ref{prop.formsminusone}~(iv) and \eqref{eq.UOw}. By \eqref{eq.wh}, $z^0_e=(-1)^{h(e)}=k_{f_{\rm L}(e)}-k_{f_{\rm R}(e)}$ for the vector $k$ with $k_f=1$ if $w(f)$ is even and $k_f=0$ if $w(f)$ is odd; so $z^0$ is the sum of the boundaries of the faces with an even winding number, which are the faces of one color. The second statement is Theorem~\ref{thm.determinant}~(ii). If $\Delta_L(-1)\ne0$, the matrix $\CU$ of the proof of that theorem is invertible over $\BQ$; so $U$ is injective on $\FF_{D,\hate}$, its kernel on $\FF_D=\FF_{D,\hate}\oplus\BZ z^0$ is $\BZ z^0$, and the index of $\Lambda_D$ is $|\det\CU|=|\Delta_L(-1)|$. Then $\Lambda_D$ has rank $c$, which is the sum of the ranks $c+1-F_{\rm b}$ and $F_{\rm b}-1$ of the two lattices of (i); so their sum is direct.

(iii) By Lemma~\ref{lem.altcoloring}, $\theta_j=-1$ for every crossing, so that $\DD_\theta=-1$ and $\DD_\theta\,B^1(\Gamma_{\rm b})=B^1(\Gamma_{\rm b})$; by (i), $\Lambda_D=Z_1(\Gamma_{\rm b})+B^1(\Gamma_{\rm b})$. The sum is direct by (ii), because $\Delta_L(-1)\ne0$ for an alternating link (Corollary~\ref{cor.altdefinite}~(ii)); this is also clear from the fact that the two lattices are orthogonal to each other for the standard scalar product.
\end{proof}

For an alternating diagram $\Lambda_D$ is, up to the signs of the coordinates at the negative crossings, the lattice in which the state sum of Section~\ref{ssec.tait} is written and which is studied in Section~\ref{ssec.lattice} (Remark~\ref{rem.branched}~(ii)); for a general diagram the signs $\theta_j$, which are the signs $\eta_{\rm w}(j)$ of the Goeritz matrix (Theorem~\ref{thm.goeritz}~(i)), twist the lattice of cuts, and the two summands are orthogonal for the form $\DD_\theta$ instead of the standard one. The second flavor has the lattice $O(\FF_D)$ of the over-strands, with the same properties.

\subsection{Examples}\label{ssec.exminusone}
We compute the forms of this section for the standard diagrams of the trefoil and of the figure-eight knot, and we check the theorems on the determinant and on the signature on them.

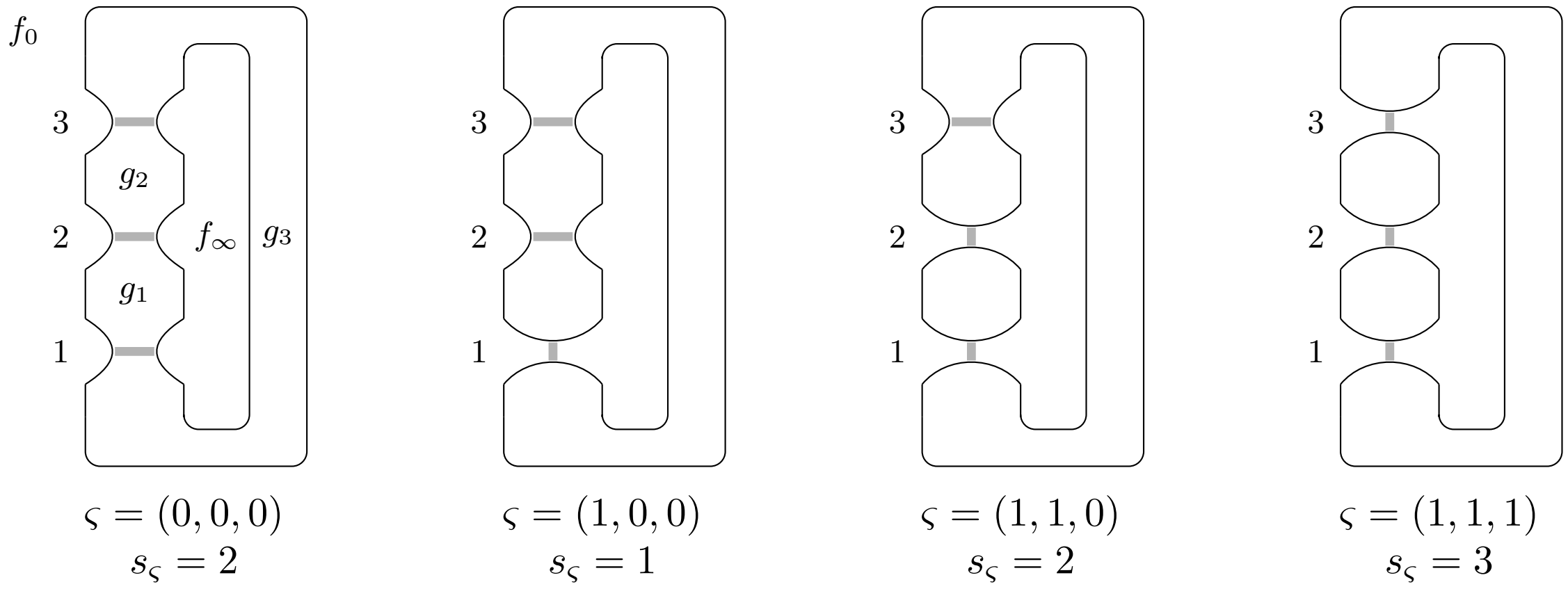
\begin{figure}[ht]
\centering
\begin{tikzpicture}[scale=0.78,every node/.style={font=\small}]
\def\smcross#1#2{%
  \ifnum#2=0
    \draw (0,{#1-0.4}) .. controls (0.44,{#1-0.12}) and (0.44,{#1+0.12}) .. (0,{#1+0.4});
    \draw (1.2,{#1-0.4}) .. controls (0.76,{#1-0.12}) and (0.76,{#1+0.12}) .. (1.2,{#1+0.4});
    \draw[line width=2.4pt,gray!60] (0.36,#1) -- (0.84,#1);
  \else
    \draw (0,{#1-0.4}) .. controls (0.3,{#1-0.04}) and (0.9,{#1-0.04}) .. (1.2,{#1-0.4});
    \draw (0,{#1+0.4}) .. controls (0.3,{#1+0.04}) and (0.9,{#1+0.04}) .. (1.2,{#1+0.4});
    \draw[line width=2.4pt,gray!60] (0.6,{#1-0.11}) -- (0.6,{#1+0.11});
  \fi}
\def\trefoilframe{%
  \draw (0,0.2) -- (0,0.6); \draw (0,1.4) -- (0,2.0); \draw (0,2.8) -- (0,3.4); \draw (0,4.2) -- (0,4.6);
  \draw (1.2,0.2) -- (1.2,0.6); \draw (1.2,1.4) -- (1.2,2.0); \draw (1.2,2.8) -- (1.2,3.4); \draw (1.2,4.2) -- (1.2,4.6);
  \draw[rounded corners=4pt] (0,4.6) -- (0,5.2) -- (2.7,5.2) -- (2.7,-0.4) -- (0,-0.4) -- (0,0.2);
  \draw[rounded corners=4pt] (1.2,4.6) -- (1.2,4.75) -- (2.0,4.75) -- (2.0,0.05) -- (1.2,0.05) -- (1.2,0.2);
  \node[font=\scriptsize] at (-0.3,1.0) {$1$}; \node[font=\scriptsize] at (-0.3,2.4) {$2$}; \node[font=\scriptsize] at (-0.3,3.8) {$3$};}
\begin{scope}
\trefoilframe \smcross{1.0}{0} \smcross{2.4}{0} \smcross{3.8}{0}
\node[font=\scriptsize] at (0.6,1.7) {$g_1$}; \node[font=\scriptsize] at (0.6,3.1) {$g_2$}; \node[font=\scriptsize] at (2.35,2.4) {$g_3$};
\node[font=\scriptsize] at (1.6,2.4) {$f_\infty$}; \node[font=\scriptsize] at (-0.75,4.9) {$f_0$};
\node at (1.2,-1.0) {$\vs=(0,0,0)$}; \node at (1.2,-1.6) {$s_\vs=2$};
\end{scope}
\begin{scope}[xshift=5.1cm]
\trefoilframe \smcross{1.0}{1} \smcross{2.4}{0} \smcross{3.8}{0}
\node at (1.2,-1.0) {$\vs=(1,0,0)$}; \node at (1.2,-1.6) {$s_\vs=1$};
\end{scope}
\begin{scope}[xshift=10.2cm]
\trefoilframe \smcross{1.0}{1} \smcross{2.4}{1} \smcross{3.8}{0}
\node at (1.2,-1.0) {$\vs=(1,1,0)$}; \node at (1.2,-1.6) {$s_\vs=2$};
\end{scope}
\begin{scope}[xshift=15.3cm]
\trefoilframe \smcross{1.0}{1} \smcross{2.4}{1} \smcross{3.8}{1}
\node at (1.2,-1.0) {$\vs=(1,1,1)$}; \node at (1.2,-1.6) {$s_\vs=3$};
\end{scope}
\end{tikzpicture}
\caption{Four smoothings of the standard diagram of the trefoil, with the crossings $1,2,3$ from bottom to top as in Figure~\ref{fig.trefoilfaces}. The circles of the smoothing are drawn in black and the bands of the state surface are indicated in gray; the faces are labeled in the first picture. The state surfaces are the canonical Seifert surface, two nonorientable surfaces and a M\"obius band; their forms of Gordon and Litherland are given in Example~\ref{ex.trefoilminusone}.}
\label{fig.trefoilsmoothings}
\end{figure}

\begin{example}[the trefoil]\label{ex.trefoilminusone}
Let $D$ be the standard diagram of the right-handed trefoil (Example~\ref{ex.trefoilflows} and Figure~\ref{fig.trefoilfaces}). Its faces are $f_0$, $f_\infty$ and the bigons $g_1,g_2,g_3$; at the crossing $j$, $f^\leftarrow_j=f_0$, $f^\rightarrow_j=f_\infty$, $f^\uparrow_j=g_j$ and $f^\downarrow_j=g_{j-1}$, with $g_0=g_3$, and all the crossings are positive. For a flow $z$ with potential $k$, Lemma~\ref{lem.cornertransfers}~(i) and \eqref{eq.UOw} give
\begin{gather*}
   w^0_j(z)=k(g_j)-k(g_{j-1}),\qquad w^1_j(z)=k(f_\infty)-k(f_0),\\
   U_j(z)=k(f_0)-k(f_\infty)-k(g_j)+k(g_{j-1}) .
\end{gather*}
\emph{Kashaev's form.} By \eqref{eq.Kashw}, in the basis $(\partial f_0,\partial f_\infty,\partial g_1,\partial g_2,\partial g_3)$,
\[
   \Kash=\begin{pmatrix}3&-3\\-3&3\end{pmatrix}\oplus\begin{pmatrix}-2&1&1\\1&-2&1\\1&1&-2\end{pmatrix}.
\]
The two blocks are the Goeritz matrices of the two colors: for the color of $f_0$ and $f_\infty$, which are the faces $f^\leftarrow_j,f^\rightarrow_j$ of every crossing, $\eta=-1$ and $\mu=3$; for the color of the bigons, $\eta=1$ and $\mu=0$. Their signatures are $1$ and $-2$, and Corollary~\ref{cor.GLgoeritz} gives $\sigma_L=1-3=-2-0=-2$; the signature of $\Kash$ is $-1=2\sigma_L+\mathrm{wr}(D)$, as in Corollary~\ref{cor.kashaevsignature}.

\emph{The state surfaces.} Figure~\ref{fig.trefoilsmoothings} shows the smoothings $(0,0,0)$, $(1,0,0)$, $(1,1,0)$ and $(1,1,1)$; by the symmetry of the diagram the other four smoothings behave like the second and the third one. The form $\GLf_0$ vanishes on $\partial f_0$ and $\partial f_\infty$, which are the Seifert circles, and it is the second block of $\Kash$ on the bigons (Theorem~\ref{thm.goeritz}~(ii)); the other forms are obtained from it by adding the squares $U_j\otimes U_j$ one at a time (Proposition~\ref{prop.formsminusone}~(v)), where $U_j$ takes the values $1,-1,-1,1$ on $\partial f_0,\partial f_\infty,\partial g_j,\partial g_{j-1}$ and vanishes on the third bigon. One finds the following Gram matrices, in the indicated bases of the homology of the surfaces.
\[
\begin{array}{c|c|c|c|c|c}
   \vs&s_\vs&\text{basis of }H_1(F_\vs;\BZ)&G_\vs&\det G_\vs&\sign G_\vs\\ \hline
   (0,0,0)&2&[\partial g_1],[\partial g_2]&\begin{pmatrix}-2&1\\1&-2\end{pmatrix}&3&-2\rule{0pt}{20pt}\\[10pt]
   (1,0,0)&1&[\partial f_0],[\partial g_1],[\partial g_2]&\begin{pmatrix}1&-1&0\\-1&-1&1\\0&1&-2\end{pmatrix}&3&-1\\[14pt]
   (1,1,0)&2&[\partial f_0],[\partial g_2]&\begin{pmatrix}2&-1\\-1&-1\end{pmatrix}&-3&0\\[10pt]
   (1,1,1)&3&[\partial f_0]&(3)&3&1
\end{array}
\]
The first surface is the canonical Seifert surface, of genus one, and the last one is a M\"obius band with three half twists, whose core has the class $[\partial f_0]$; the two others are nonorientable surfaces with Euler characteristic $-2$ and $-1$. In all cases $|\det G_\vs|=3=|\Delta_L(-1)|$, the rank is $4-s_\vs$, and the signature minus the number of the crossings with $\vs_j=1$ is $\sigma_L=-2$, as in Theorems~\ref{thm.determinant} and~\ref{thm.signature}; each change of a smoothing raises the signature by one, and it changes the rank by $+1$, $-1$ and $-1$. For $\vs=(1,0,0)$ the single circle has crossings on both sides: its two regions are that of the faces $f_0$ and $f_\infty$, in which the crossing $1$ lies, and that of the three bigons, in which the crossings $2$ and $3$ lie. So $J$ changes the sign of $[\partial f_0]$ and fixes $[\partial g_1]$ and $[\partial g_2]$, up to a global sign, and the form of the state surface above the plane is obtained from $G_\vs$ by changing the sign of the entries $-1$ of its first row and of its first column (Proposition~\ref{prop.J}).
\end{example}

\begin{example}[the figure-eight knot]\label{ex.figeightminusone}
Let $D$ be the closure of $(\sigma_1\sigma_2^{-1})^2$ (Example~\ref{ex.figeight} and Figures~\ref{fig.diagrams} and~\ref{fig.circles}). For the oriented smoothing the three Seifert circles are nested, $C_1$ enclosing $C_2$ and $C_2$ enclosing $C_3$. There are four regions (Figure~\ref{fig.smoothingregions}~(a)): the unbounded one $\Omega_0$, of depth $0$; the region $\Omega_1$ between $C_1$ and $C_2$, of depth $1$, in which the positive crossings $1$ and $3$ lie; the region $\Omega_2$ between $C_2$ and $C_3$, of depth $2$, in which the negative crossings $2$ and $4$ lie; and the interior $\Omega_3$ of $C_3$, of depth $3$. The circle $C_2$ has crossings on both sides. The groups $H_{\Omega_1}$ and $H_{\Omega_2}$ have rank one, and they are generated by the classes of the cycles $\gamma_1=e_3-e_1$ and $\gamma_2=e_4-e_2$ of Example~\ref{ex.figeight}. In this basis
\[
   \Goe_0=\begin{pmatrix}-2&0\\0&2\end{pmatrix},\qquad \Stk_0=\begin{pmatrix}0&-1\\-1&0\end{pmatrix},\qquad \GLf_0=\begin{pmatrix}-2&-1\\-1&2\end{pmatrix},\qquad \GLf^\uparrow_0=\begin{pmatrix}-2&1\\1&2\end{pmatrix}:
\]
the Goeritz form is the sum of the forms $(-2)$ and $(2)$ of a positive and of a negative Hopf band, one for each region, the canonical surface is the plumbing of these two bands, and the stacking form is the coupling between them, which changes sign when the surface is pushed to the other side of the plane. The involution is $J=\diag(-1,1)$. Accordingly, the Seifert matrices of the two canonical surfaces, $V_F=\bigl(\begin{smallmatrix}-1&-1\\0&1\end{smallmatrix}\bigr)$ and $V_{F^\uparrow}=\bigl(\begin{smallmatrix}-1&0\\1&1\end{smallmatrix}\bigr)$, have the same diagonal entries, and $V_{F^\uparrow}(\gamma_1,\gamma_2)=-V_F(\gamma_2,\gamma_1)=0$, $V_{F^\uparrow}(\gamma_2,\gamma_1)=-V_F(\gamma_1,\gamma_2)=1$, as in Corollary~\ref{cor.twoseifert}. The determinant of $\GLf_0$ is $-5$ and its signature is $0=\sigma_L$.

The diagram is alternating. In the coloring of Lemma~\ref{lem.altcoloring} the black checkerboard smoothing is $\vs_{\rm b}=(1,0,1,0)$ and the white one is $\vs_{\rm w}=(0,1,0,1)$; there are $F_{\rm b}=3$ black and $F_{\rm w}=3$ white faces, and in suitable bases the forms of the two checkerboard surfaces are
\[
   \GLf_{\vs_{\rm b}}=\begin{pmatrix}3&1\\1&2\end{pmatrix}\qquad\text{and}\qquad\GLf_{\vs_{\rm w}}=\begin{pmatrix}-3&2\\2&-3\end{pmatrix},
\]
which are positive and negative definite, of determinant $5$. Their signatures $2$ and $-2$, corrected by $\ve_1+\ve_3=2$ and by $\ve_2+\ve_4=-2$, give $\sigma_L=0=c_-+1-F_{\rm b}$ (Corollary~\ref{cor.altdefinite}).
\end{example}

\section{A two-flavor Kauffman state sum and positivity on alternating diagrams}\label{sec.kauffman}

Kauffman's state model \cite{KauffmanFKT} is the expansion of the determinant of Alexander's matrix over the bijections between the crossings and the faces of the diagram. By Theorem~\ref{thm.faceform} both edge polynomials are such determinants, with the crossing variables inserted: they are the determinants of the crossing forms in the basis of the faces of the lattice of flows. Expanding them turns the determinant formula into a sum over \emph{pairs} of Kauffman states, one for each flavor, which interact through $\BE_D$ at the crossings (Proposition~\ref{prop.statesum}). Grouped according to the two markers which they put at every crossing, these pairs are the Kauffman states with multiplicity two of a state model of Ballinger \cite{BallingerTalk}, whose invariant for $N=2$ is thereby identified with $\LG$ (Corollary~\ref{cor.ballinger}).

The two specializations of $\LG$ recalled in the introduction, which hold for every link, appear as two regimes of this model: at $t_0t_1=1$ the two states of a pair are independent (Remark~\ref{rem.statesum}), and at $t_0t_1=-1$ the pairs of two equal states give the whole sum, which is Kauffman's state sum for $\Delta_L(t_0^2)$, the other pairs canceling. The second regime is treated in Section~\ref{ssec.signs}, together with the exact signs of the states, on which its statement rests.

For an alternating diagram the signs of the states obey Kasteleyn's rule for planar bipartite graphs, and the state sum has no cancellation (Theorem~\ref{thm.kauffmanalt}). This is the theorem of \cite{HK}, that $\LG$ of an alternating link is alternating, proved in Kauffman's model instead of Crowell's. The states of the two models correspond (Corollary~\ref{cor.statestrees}), and the two extremal trees of Section~\ref{ssec.alternating} become one extremal state. Finally, on an alternating diagram the crossings at which a state sits between the two outgoing or between the two incoming edges count the signature $\sigma_L$ of the link (Proposition~\ref{prop.signature}, which rests on the signature formula of Section~\ref{sec.minusone}); consequently the top coefficient of a fibered alternating link is $q^{-\sigma_L}$, and the Hopf invariant of its fiber is $\frac12(2g(L)+\mu-1+\sigma_L)$ (Corollary~\ref{cor.extremalstates}).

Section~\ref{ssec.tait} writes the state sum of an alternating diagram on its Tait graph: a state is a spanning tree, its two flavors are two orientations of the graph, and the weight of a pair of states is read off the inner products of these orientations with each other and with the orientation defined by the over-strands (Theorem~\ref{thm.tait}). The latter is the sum of an integral flow, carried by the positive crossings, and of an integral cut, carried by the negative ones (Proposition~\ref{prop.omega}), and in these terms $-\sigma_L$ is a correction term of the double branched cover (Remark~\ref{rem.branched}).

Section~\ref{ssec.gluing} gives a formula for the polynomial of a link which is cut into two tangles by a Conway sphere, in terms of three closures of each tangle (Theorem~\ref{thm.gluing}); on the Tait graph it deals with the pieces which are attached to the rest of the graph along two vertices (Corollary~\ref{cor.twosum}), and its proof uses the representation theory instead of the determinant formula.

The last subsection discusses to what extent the polynomial of an alternating link is computed from the lattice of its two checkerboard surfaces, which is the lattice of the over-strands of Section~\ref{ssec.underlattice}: the states are the $\mathrm{Spin}^c$ structures of the double branched cover, and the state sum is a sum over the pairs of them (Corollary~\ref{cor.latticeknots} and Proposition~\ref{prop.spinc}).

\subsection{States}\label{ssec.states}
We introduce the Kauffman states of a diagram, in a notation which treats the two signs of a crossing at once, and we expand the two edge polynomials over them: every state contributes one monomial to each of the two polynomials, with a sign (Proposition~\ref{prop.statesum}). The determinant formula then becomes a sum over the ordered pairs of states.

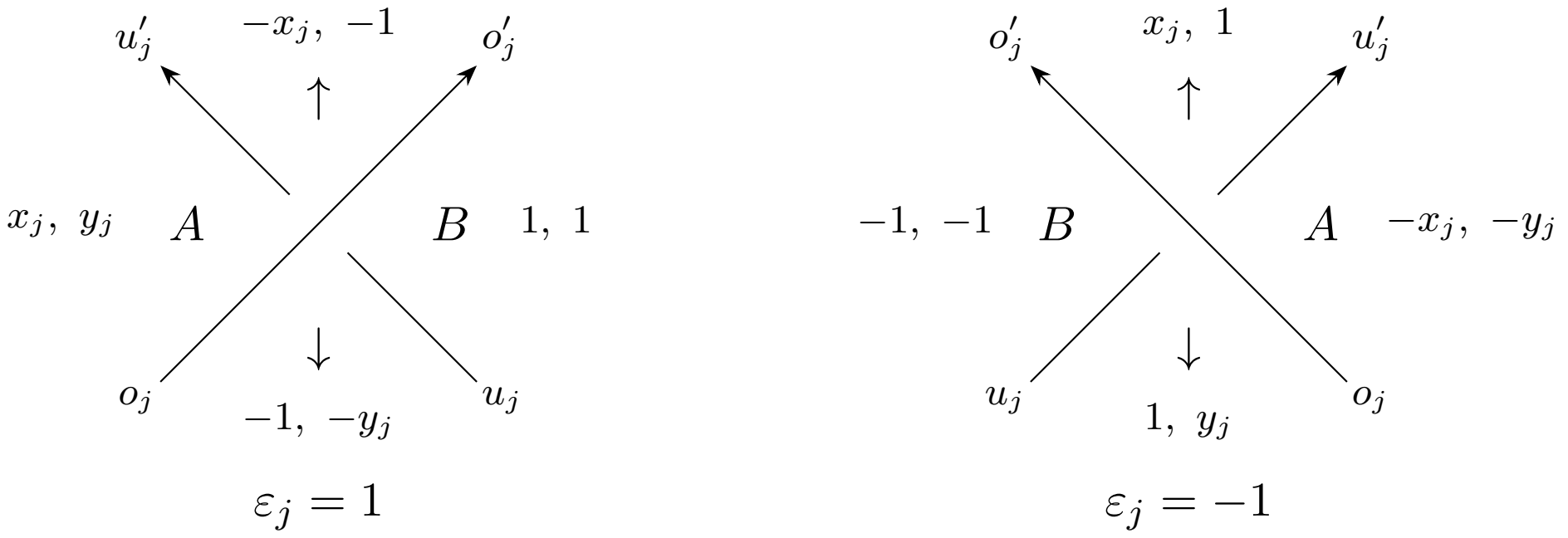
\begin{figure}[ht]
\centering
\begin{tikzpicture}[>=Stealth,scale=1.05,every node/.style={font=\small}]
\begin{scope}
\draw (1.2,-1.2) -- (0.22,-0.22); \draw[->] (-0.22,0.22) -- (-1.2,1.2);
\draw[->] (-1.2,-1.2) -- (1.2,1.2);
\node[below left,font=\scriptsize] at (-1.1,-1.1) {$o_j$}; \node[above right,font=\scriptsize] at (1.1,1.1) {$o'_j$};
\node[below right,font=\scriptsize] at (1.1,-1.1) {$u_j$}; \node[above left,font=\scriptsize] at (-1.1,1.1) {$u'_j$};
\node at (0,0.95) {$\uparrow$}; \node at (0,-0.95) {$\downarrow$};
\node at (-1.0,0) {$A$}; \node at (1.0,0) {$B$};
\node[font=\scriptsize] at (0,1.5) {$-x_j,\ -1$}; \node[font=\scriptsize] at (0,-1.5) {$-1,\ -y_j$};
\node[font=\scriptsize] at (-1.95,0) {$x_j,\ y_j$}; \node[font=\scriptsize] at (1.8,0) {$1,\ 1$};
\node at (0,-2.15) {$\ve_j=1$};
\end{scope}
\begin{scope}[xshift=6.6cm]
\draw (-1.2,-1.2) -- (-0.22,-0.22); \draw[->] (0.22,0.22) -- (1.2,1.2);
\draw[->] (1.2,-1.2) -- (-1.2,1.2);
\node[below right,font=\scriptsize] at (1.1,-1.1) {$o_j$}; \node[above left,font=\scriptsize] at (-1.1,1.1) {$o'_j$};
\node[below left,font=\scriptsize] at (-1.1,-1.1) {$u_j$}; \node[above right,font=\scriptsize] at (1.1,1.1) {$u'_j$};
\node at (0,0.95) {$\uparrow$}; \node at (0,-0.95) {$\downarrow$};
\node at (-1.0,0) {$B$}; \node at (1.0,0) {$A$};
\node[font=\scriptsize] at (0,1.5) {$x_j,\ 1$}; \node[font=\scriptsize] at (0,-1.5) {$1,\ y_j$};
\node[font=\scriptsize] at (-2.0,0) {$-1,\ -1$}; \node[font=\scriptsize] at (2.15,0) {$-x_j,\ -y_j$};
\node at (0,-2.15) {$\ve_j=-1$};
\end{scope}
\end{tikzpicture}
\caption{The corners of a crossing: the vertical corners $\uparrow$, $\downarrow$, the $A$-corner between $o_j$ and $u'_j$ and the $B$-corner between $u_j$ and $o'_j$. Next to each corner, the coefficients $a^\diamond_j,\ b^\diamond_j$ of \eqref{eq.cornercoeff}, which are the entries of \eqref{eq.facetable}. The crossing belongs to $X(\mathfrak s)$ if the corner of the state $\mathfrak s$ is $\uparrow$ or $A$, and to $Y(\mathfrak s)$ if it is $\downarrow$ or $A$.}
\label{fig.corners}
\end{figure}

Let $D$ be a diagram and $\hate$ an edge of $D$; we use the notation of Section~\ref{ssec.facematrices}. A crossing $j$ has four corners, which we denote by $\uparrow,\leftarrow,\downarrow,\rightarrow$ as in Figure~\ref{fig.faces}, the corner $\diamond$ lying in the face $f^\diamond_j$. We call $\uparrow$ and $\downarrow$, the corners between the two outgoing and between the two incoming edges, the \emph{vertical} corners. One of the two horizontal corners lies between $o_j$ and $u'_j$ and the other between $u_j$ and $o'_j$; we call them the \emph{$A$-corner} and the \emph{$B$-corner} of $j$, because $o_j,u'_j$ lie on the Seifert circle $A_j$ and $u_j,o'_j$ on $B_j$ (Figure~\ref{fig.crossing}). Thus the $A$-corner is $\leftarrow$ if $\ve_j=1$ and $\rightarrow$ if $\ve_j=-1$ (Figure~\ref{fig.corners}), and the two corners adjacent to each of the four edges of the crossing are
\begin{equation}\label{eq.corneredges}
   u'_j:\ \{\uparrow,A\},\qquad o_j:\ \{\downarrow,A\},\qquad u_j:\ \{\downarrow,B\},\qquad o'_j:\ \{\uparrow,B\}.
\end{equation}
In this notation the table \eqref{eq.facetable} is the same for the two signs of the crossing: the rows of the face matrices are $a_j(x_j)=\sum_\diamond a^\diamond_j\,e_{f^\diamond_j}$ and $b_j(y_j)=\sum_\diamond b^\diamond_j\,e_{f^\diamond_j}$, the sums running over the four corners of $j$, with
\begin{equation}\label{eq.cornercoeff}
   a^\diamond_j=\ve_j\,(-1)^{\ind{\diamond\text{ vertical}}}\,x_j^{\ind{\diamond\in\{\uparrow,A\}}},\qquad
   b^\diamond_j=\ve_j\,(-1)^{\ind{\diamond\text{ vertical}}}\,y_j^{\ind{\diamond\in\{\downarrow,A\}}} .
\end{equation}
So the coefficients of the two flavors have the same signs, the variable $x_j$ sits at the two corners adjacent to the outgoing under-edge $u'_j$, and $y_j$ at the two corners adjacent to the incoming over-edge $o_j$.

\begin{definition}\label{def.state}
A \emph{Kauffman state} of $(D,\hate)$ is a map $\mathfrak s$ which assigns to every crossing $j$ one of its four corners $\mathfrak s(j)$, in such a way that $j\mapsto f^{\mathfrak s(j)}_j$ is a bijection from the set of crossings onto $\CF\sm\{f_{\rm L}(\hate),f_{\rm R}(\hate)\}$: every face other than the two faces adjacent to $\hate$ contains exactly one of the chosen corners, and these two faces contain none. We write $\sgn(\mathfrak s)=\pm1$ for the sign of this bijection with respect to the orders of the crossings and of the faces which are used for the rows and the columns of $Q_{D,\hate}$ and $\tQ_{D,\hate}$, and we put
\begin{gather*}
   X(\mathfrak s)=\{j:\mathfrak s(j)\in\{\uparrow,A\}\},\qquad Y(\mathfrak s)=\{j:\mathfrak s(j)\in\{\downarrow,A\}\},\\
   V(\mathfrak s)=\{j:\mathfrak s(j)\text{ is vertical}\},
\end{gather*}
so that $V(\mathfrak s)=X(\mathfrak s)\triangle Y(\mathfrak s)$, with $X(\mathfrak s)\sm Y(\mathfrak s)=\{j:\mathfrak s(j)={\uparrow}\}$ and $Y(\mathfrak s)\sm X(\mathfrak s)=\{j:\mathfrak s(j)={\downarrow}\}$.
\end{definition}

These are the states of \cite{KauffmanFKT} for the diagram $D$ with its two stars in the faces adjacent to $\hate$, the corner $\mathfrak s(j)$ being the state marker at the crossing $j$. By \eqref{eq.corneredges}, $X(\mathfrak s)$ is the set of crossings at which the marker is adjacent to the outgoing under-edge, and $Y(\mathfrak s)$ the set of those at which it is adjacent to the incoming over-edge. A change of the orders multiplies all the signs $\sgn(\mathfrak s)$ by the same sign. Figure~\ref{fig.trefoilstates} shows the three states of the standard diagram of the trefoil, which is worked out in Example~\ref{ex.states}. As in Section~\ref{ssec.epsdegree} we write $x^S=\prod_{j\in S}x_j$, $y^S=\prod_{j\in S}y_j$ and $\ve(S)=\sum_{j\in S}\ve_j$ for a set $S$ of crossings, and we recall that $|S|\equiv\ve(S)$ modulo $2$.

\begin{figure}[ht]
\centering
\begin{tikzpicture}[>=Stealth,scale=0.6]
\def\trefoilbase{%
   \draw (0,0) -- (0,0.5); \draw[->] (0,1.5) -- (0,2.5); \draw[->] (0,3.5) -- (0,4.5); \draw (0,5.5) -- (0,6);
   \draw[rounded corners=5pt] (0,6) -- (0,7.1) -- (4.0,7.1) -- (4.0,-1.1) -- (0,-1.1) -- (0,0);
   \draw (1.7,0) -- (1.7,0.5); \draw[->] (1.7,1.5) -- (1.7,2.5); \draw[->] (1.7,3.5) -- (1.7,4.5); \draw (1.7,5.5) -- (1.7,6);
   \draw[rounded corners=5pt] (1.7,6) -- (1.7,6.5) -- (2.9,6.5) -- (2.9,-0.5) -- (1.7,-0.5) -- (1.7,0);
   \foreach \y in {1,3,5} { \draw (1.7,{\y-0.5}) -- (0,{\y+0.5}); \draw[preaction={draw=white,line width=3pt}] (0,{\y-0.5}) -- (1.7,{\y+0.5}); }
   \foreach \y/\j in {1/1,3/2,5/3} \node[font=\scriptsize] at (-0.3,\y) {$\j$};
   \node[font=\scriptsize] at (-0.9,6.6) {$\bigstar$}; \node[font=\scriptsize] at (3.45,3) {$\bigstar$};}
\def\mk#1#2{\fill[red] (0.85,#1) ++(#2) circle (4.2pt);}
\begin{scope}
\trefoilbase \mk{1}{0,0.48} \mk{3}{0,0.48} \mk{5}{0.5,0}
\node[font=\small] at (1.6,-1.9) {$\mathfrak s_1=(\uparrow,\uparrow,B)$};
\node[font=\small] at (1.6,-2.8) {$X=\{1,2\},\ Y=\emptyset$};
\node[font=\small] at (1.6,-3.7) {$b=0$};
\end{scope}
\begin{scope}[xshift=8cm]
\trefoilbase \mk{1}{0,0.48} \mk{3}{0.5,0} \mk{5}{0,-0.48}
\node[font=\small] at (1.6,-1.9) {$\mathfrak s_2=(\uparrow,B,\downarrow)$};
\node[font=\small] at (1.6,-2.8) {$X=\{1\},\ Y=\{3\}$};
\node[font=\small] at (1.6,-3.7) {$b=1$};
\end{scope}
\begin{scope}[xshift=16cm]
\trefoilbase \mk{1}{0.5,0} \mk{3}{0,-0.48} \mk{5}{0,-0.48}
\node[font=\small] at (1.6,-1.9) {$\mathfrak s_3=(B,\downarrow,\downarrow)$};
\node[font=\small] at (1.6,-2.8) {$X=\emptyset,\ Y=\{2,3\}$};
\node[font=\small] at (1.6,-3.7) {$b=2$};
\end{scope}
\end{tikzpicture}
\caption{The three Kauffman states of the standard diagram of the trefoil for the cut edge $\hate=u'_3=o_1$ (Example~\ref{ex.states}; the faces are named in Figure~\ref{fig.trefoilfaces}). The two faces adjacent to $\hate$ are marked by stars, and the marker of a state at a crossing is the dot in one of its four corners: every other face contains exactly one marker. The crossing belongs to $X(\mathfrak s)$ if the marker is adjacent to the outgoing under-edge, which here is the upper left one, and to $Y(\mathfrak s)$ if it is adjacent to the incoming over-edge, the lower left one; $b(\mathfrak s)$ counts the markers between the two incoming edges. Thus $P_{D,\hate}=x_1x_2-x_1+1$ and $\tP_{D,\hate}=1-y_3+y_2y_3$.}
\label{fig.trefoilstates}
\end{figure}

\begin{proposition}[the two-flavor state sum]\label{prop.statesum}
Let $D$ be a diagram of a link $L$ with $\mu$ components, $\hate$ an edge of $D$ and $\ve_D=\prod_j\ve_j$. In the following sums $\mathfrak s$ and $\mathfrak s'$ run over the Kauffman states of $(D,\hate)$, and $b(\mathfrak s)$ is the number of crossings $j$ with $\mathfrak s(j)={\downarrow}$, i.e.\ at which the marker lies between the two incoming edges.
\begin{enumerate}\renewcommand{\theenumi}{\roman{enumi}}
\item For every state $\mathfrak s$,
\[
   \ve(X(\mathfrak s))+\ve(Y(\mathfrak s))=m_{D,\hate},\qquad\text{where}\quad m_{D,\hate}=\sum_{f\ne f_{\rm L}(\hate),f_{\rm R}(\hate)}w(f)\ +\ \sum_{j=1}^c\Bigl(h(u_j)+\frac{1+\ve_j}2\Bigr),
\]
$w$ being the winding function of Lemma~\ref{lem.winding}, normalized by \eqref{eq.wh}, and $h$ the level function of Definition~\ref{def.level}. Moreover $|V(\mathfrak s)|\equiv m_{D,\hate}\equiv\mu-1$ modulo $2$ for every state.
\item The edge polynomials are
\[
   P_{D,\hate}(x)=(-1)^{\mu-1}\sum_{\mathfrak s}(-1)^{b(\mathfrak s)}\,x^{X(\mathfrak s)},\qquad \tP_{D,\hate}(y)=\sum_{\mathfrak s}(-1)^{b(\mathfrak s)}\,y^{Y(\mathfrak s)} .
\]
The signs of the states are $\sgn(\mathfrak s)=\eta\,(-1)^{b(\mathfrak s)}$, where the sign $\eta=\pm1$ depends only on the orders of the crossings and of the faces, and the determinants of the face matrices are
\begin{align*}
   \det Q_{D,\hate}(x)&=\ve_D\,(-1)^{m_{D,\hate}}\sum_{\mathfrak s}\sgn(\mathfrak s)\,x^{X(\mathfrak s)}=\ve_D\,\eta\,P_{D,\hate}(x),\\
   \det\tQ_{D,\hate}(y)&=\ve_D\,(-1)^{m_{D,\hate}}\sum_{\mathfrak s}\sgn(\mathfrak s)\,y^{Y(\mathfrak s)}=(-1)^{\mu-1}\,\ve_D\,\eta\,\tP_{D,\hate}(y).
\end{align*}
\item The Links--Gould polynomial is
\[
   \LG(L;t_0,t_1)=\sum_{(\mathfrak s,\mathfrak s')}(-1)^{b(\mathfrak s)+b(\mathfrak s')}\prod_{j\in X(\mathfrak s)\sm Y(\mathfrak s')}t_0^{\ve_j}\ \prod_{j\in Y(\mathfrak s')\sm X(\mathfrak s)}t_1^{\ve_j},
\]
the sum running over the ordered pairs of states.
\end{enumerate}
\end{proposition}
\begin{proof}
(i) Let $w_j$ be the winding number of the face of the $B$-corner of $j$. For every corner $\diamond$ of $j$,
\[
   \ve_j\bigl(\ind{\diamond\in\{\uparrow,A\}}+\ind{\diamond\in\{\downarrow,A\}}\bigr)=w(f^\diamond_j)-w_j:
\]
the left-hand side is $2\ve_j$, $\ve_j$, $\ve_j$, $0$ at the corners $A$, $\uparrow$, $\downarrow$, $B$, and by \eqref{eq.wcrossing} the winding numbers of the faces at $\leftarrow$, $\uparrow$, $\downarrow$, $\rightarrow$ are $k$, $k-1$, $k-1$, $k-2$, the $A$-corner being $\leftarrow$ if $\ve_j=1$ and $\rightarrow$ if $\ve_j=-1$. Summing over the crossings, for a state $\mathfrak s$,
\[
   \ve(X(\mathfrak s))+\ve(Y(\mathfrak s))=\sum_jw\bigl(f^{\mathfrak s(j)}_j\bigr)-\sum_jw_j=\sum_{f\ne f_{\rm L}(\hate),f_{\rm R}(\hate)}w(f)-\sum_jw_j,
\]
since $j\mapsto f^{\mathfrak s(j)}_j$ is a bijection onto these faces. By \eqref{eq.mcrossing} the face of the $B$-corner is $f_{\rm R}(u_j)$ if $\ve_j=1$ and $f_{\rm L}(u_j)$ if $\ve_j=-1$, so that $w_j=-h(u_j)-\frac{1+\ve_j}2$ by \eqref{eq.wh}. Finally $|V(\mathfrak s)|=|X(\mathfrak s)\triangle Y(\mathfrak s)|\equiv|X(\mathfrak s)|+|Y(\mathfrak s)|\equiv\ve(X(\mathfrak s))+\ve(Y(\mathfrak s))$ modulo $2$. The congruence $m_{D,\hate}\equiv\mu-1$ is part of the determination of the signs, in Section~\ref{ssec.signs} (Proposition~\ref{prop.signs}~(i)).

(ii) \emph{The face matrices.} By \eqref{eq.cornercoeff}, the row $j$ of $Q_{D,\hate}(x)$ is the sum, over the corners $\diamond$ of $j$ which do not lie in the two removed faces, of $a^\diamond_j$ times the unit row vector of the face $f^\diamond_j$. The determinant is a linear function of each row, so it is the sum, over all the choices of one such corner $\mathfrak s(j)$ at every crossing, of $\prod_ja^{\mathfrak s(j)}_j$ times the determinant of the matrix whose row $j$ is the unit vector of $f^{\mathfrak s(j)}_j$. This last determinant is $\sgn(\mathfrak s)$ if the faces $f^{\mathfrak s(j)}_j$ are pairwise distinct, i.e.\ if $\mathfrak s$ is a state, and $0$ otherwise; and for a state $\prod_ja^{\mathfrak s(j)}_j=\ve_D\,(-1)^{|V(\mathfrak s)|}\,x^{X(\mathfrak s)}=\ve_D\,(-1)^{m_{D,\hate}}\,x^{X(\mathfrak s)}$ by \eqref{eq.cornercoeff} and (i). The same computation with the coefficients $b^\diamond_j$, for which $\prod_jb^{\mathfrak s(j)}_j=\ve_D\,(-1)^{|V(\mathfrak s)|}\,y^{Y(\mathfrak s)}$, gives the expansion of $\det\tQ_{D,\hate}(y)$.

\emph{The signs of the states.} By the Clock Theorem of \cite{KauffmanFKT}, any two states are related by a sequence of moves each of which exchanges the faces of the markers of two crossings joined by an edge; such a move changes $\sgn(\mathfrak s)$ and changes by one the number $b(\mathfrak s)$ of crossings at which the marker lies between the two incoming edges (the \emph{black holes} of \cite{KauffmanFKT}), since among the four corners adjacent to an edge exactly one lies between two incoming edges. Hence
\begin{equation}\label{eq.blackholes}
   \sgn(\mathfrak s)=\eta\,(-1)^{b(\mathfrak s)},\qquad b(\mathfrak s)=|\{j:\mathfrak s(j)={\downarrow}\}|,
\end{equation}
with a sign $\eta=\pm1$ which does not depend on the state.

\emph{The edge polynomials.} By Theorem~\ref{thm.faceform} there is a sign $\kappa=\pm1$ such that
\[
   P_{D,\hate}=\kappa\det Q_{D,\hate}\qquad\text{and}\qquad\tP_{D,\hate}=(-1)^{\mu-1}\kappa\det\tQ_{D,\hate} .
\]
With the expansions of the two determinants and with \eqref{eq.blackholes} this gives
\begin{equation}\label{eq.PPstates}
   P_{D,\hate}(x)=\beta_a\sum_{\mathfrak s}(-1)^{b(\mathfrak s)}\,x^{X(\mathfrak s)},\qquad \tP_{D,\hate}(y)=\beta_b\sum_{\mathfrak s}(-1)^{b(\mathfrak s)}\,y^{Y(\mathfrak s)},\qquad \beta_a\beta_b=(-1)^{\mu-1},
\end{equation}
where $\beta_a=\kappa\,\ve_D\,(-1)^{m_{D,\hate}}\,\eta$ and $\beta_b=(-1)^{\mu-1}\beta_a$. So the two state sums have the stated form up to one common sign, the relation between the signs of the two flavors being that of Theorem~\ref{thm.flows}. This sign is determined in Section~\ref{ssec.signs}: $\beta_a=(-1)^{\mu-1}$ and $\beta_b=1$ (Proposition~\ref{prop.signs}). For knots it is fixed by the values of the edge polynomials at the point where all the variables are $1$ and by the descending diagrams, and for links by induction on the number of components, a diagram with one more component being the smoothing of a crossing. The expressions of the two determinants by the edge polynomials follow, since $(-1)^{m_{D,\hate}}=(-1)^{\mu-1}$ by (i).

(iii) By (ii), Theorem~\ref{thm.HK} and the definition \eqref{eq.ED} of $\BE_D$,
\[
   \LG(L;t_0,t_1)=(-1)^{\mu-1}\,\BE_D\bigl[P_{D,\hate}(x)\,\tP_{D,\hate}(y)\bigr]=\sum_{(\mathfrak s,\mathfrak s')}(-1)^{b(\mathfrak s)+b(\mathfrak s')}\ \BE_D\bigl[x^{X(\mathfrak s)}y^{Y(\mathfrak s')}\bigr],
\]
which is the formula. Only the product $\beta_a\beta_b=(-1)^{\mu-1}$ of the two signs of \eqref{eq.PPstates} enters here, and not their individual values.
\end{proof}

\begin{remark}\label{rem.statesum}
(i) \emph{Kauffman's model.} At the Alexander point, (ii) is Kauffman's state sum for the Alexander polynomial: $\sum_{\mathfrak s}(-1)^{b(\mathfrak s)}\,t^{\ve(X(\mathfrak s))}=(-1)^{\mu-1}P_{D,\hate}(t^{\ve})$ is $\Delta_L(t)$ up to a unit by Lemma~\ref{lem.alexpoly}, a state being weighted by $t^{\ve_j}$ at the crossings where its marker is adjacent to the outgoing under-edge; the unit is $t^{m_{D,\hate}/2}$ for the Conway-normalized Alexander polynomial (Proposition~\ref{prop.signs}~(iv)). By (i), the second flavor is the same state sum read backward, $\ve(Y(\mathfrak s))=m_{D,\hate}-\ve(X(\mathfrak s))$, with the same signs. The signs are local as well: by \eqref{eq.blackholes} the sign of a state is $(-1)^{b(\mathfrak s)}$, where $b(\mathfrak s)$ counts the markers which lie between two incoming edges. Thus (iii) is a state sum for $\LG$ with local weights and local signs.

(ii) \emph{The interaction, and the point $t_0t_1=1$.} Since $\BE_D[x^Sy^T]=t_0^{\ve(S)}\,t_1^{\ve(T)}\,(t_0t_1)^{-\ve(S\cap T)}$, the weight of a pair $(\mathfrak s,\mathfrak s')$ in (iii) is the product of the weights $t_0^{\ve(X(\mathfrak s))}$ and $t_1^{\ve(Y(\mathfrak s'))}$ of two independent states, one for each flavor, corrected by the factor $(t_0t_1)^{-\ve_j}$ at every crossing $j\in X(\mathfrak s)\cap Y(\mathfrak s')$: the two states interact at the crossings where the marker of the first one is adjacent to the outgoing under-edge and the marker of the second one to the incoming over-edge. At the free fermion point $t_0=t$, $t_1=t^{-1}$ there is no interaction, and by (i) the sum is $t^{-m_{D,\hate}}\bigl(\sum_{\mathfrak s}\sgn(\mathfrak s)\,t^{\ve(X(\mathfrak s))}\bigr)^2$, the square of Kauffman's state sum, in agreement with Lemma~\ref{lem.alexpoly}.

(iii) \emph{The diagonal.} The pairs $(\mathfrak s,\mathfrak s)$ contribute
\[
   \sum_{\mathfrak s}\ \prod_{j:\,\mathfrak s(j)={\uparrow}}t_0^{\ve_j}\prod_{j:\,\mathfrak s(j)={\downarrow}}t_1^{\ve_j}
\]
to the sum of (iii): only the vertical corners carry a weight, $t_0^{\ve_j}$ between the two outgoing edges and $t_1^{\ve_j}$ between the two incoming ones, and all the coefficients are positive. At $t_0t_1=-1$ this diagonal is Kauffman's state sum for $\Delta_L(t_0^2)$, and it is the whole sum: the specialization $\LG(L;t_0,-t_0^{-1})=\Delta_L(t_0^2)$, which holds for every link \cite{DWIL}, says that the pairs with $\mathfrak s\ne\mathfrak s'$ cancel there (Corollary~\ref{cor.diagonal}). In this sense the two specializations of $\LG$ recalled in the introduction are the two extreme regimes of the state sum: at $t_0t_1=1$ the two states of a pair are independent, and at $t_0t_1=-1$ they coincide.

(iv) \emph{The Kashaev form.} The pairs of states are the terms of the expansion of the determinant of the deformed Kashaev matrix $\tau_{D,\hate}(x,y)=Q_{D,\hate}(x)^{\T}\DD_\ve\tQ_{D,\hate}(y)$ of Corollary~\ref{cor.deformedkashaev}: since $\det\tau_{D,\hate}(x,y)=\ve_D\det Q_{D,\hate}(x)\det\tQ_{D,\hate}(y)$, a pair $(\mathfrak s,\mathfrak s')$ selects, in the block \eqref{eq.localblocks} of each crossing, the entry in the row of the corner of $\mathfrak s$ and in the column of the corner of $\mathfrak s'$, and $\BE_D$ acts on the selected entries, crossing by crossing. By Theorem~\ref{thm.mixed} this can be said without $\BE_D$: the state sum (iii) is the expansion \eqref{eq.mixedexplicit} of the mixed discriminant $\lb\Psi_1,\dots,\Psi_c\rb_{\hate}$ of the coupled crossing forms in the basis of the boundaries of the faces, the entries being expanded over the corners. The factor of the crossing $j$ in the term of a pair of states is, up to sign, the coefficient of $\Psi_j$ which the two markers select. Let $e$ be the edge of the under-strand which is adjacent to the marker of $\mathfrak s$, and $e'$ the edge of the over-strand which is adjacent to the marker of $\mathfrak s'$; by \eqref{eq.Psi} the factor is $1$ if $(e,e')=(u'_j,o_j)$, $t_0^{\ve_j}$ if $(e,e')=(u'_j,o'_j)$, $t_1^{\ve_j}$ if $(e,e')=(u_j,o_j)$, and $1$ if $(e,e')=(u_j,o'_j)$. The interaction of (ii) is the coefficient $1$ of $z_{u'_j}z'_{o_j}$ in $\Psi_j$, which stands in the place of the product $t_0^{\ve_j}t_1^{\ve_j}$ of the two weights.

(v) \emph{Trees and states.} For closed braids, \cite[Section~6.4]{HK} expands the two edge polynomials over spanning trees by the matrix-tree theorem and obtains a state sum for $\LG$ over pairs of trees, in which a tree has the weight $x_j$, $1-x_j$ or $1$ at each crossing. These trees are in general more numerous than the Kauffman states (eleven against five for the closure of $\sigma_1^5$) and their weights are not monomials; in the expansion over Kauffman states each state is one monomial of the edge polynomial, with a sign.
\end{remark}

\subsection{The signs of the states, and the point \texorpdfstring{$t_0t_1=-1$}{t0 t1 = -1}}\label{ssec.signs}
In this subsection we determine the signs of Proposition~\ref{prop.statesum}: the global signs of the state sums of the two edge polynomials, with which the state sum for $\LG$ has no global sign for any link, and the exact form of Kauffman's state sum for the Conway polynomial (Proposition~\ref{prop.signs}). For knots these signs are fixed by the values at the point where all the variables are $1$ and by one diagram for which everything is known, the descending one. For links the edge polynomials vanish at that point, and we argue by induction on the number of components: a diagram with one more component is obtained by smoothing a crossing, and the edge polynomials of the smoothed diagram are determined by those of the diagram and of the diagram with the crossing changed (Lemma~\ref{lem.smoothing}). With the signs at hand we then describe the form which the first of the two specializations of $\LG$ recalled in the introduction,
\[
   \LG(L;t_0,-t_0^{-1})=\Delta_L(t_0^2),
\]
a theorem for every oriented link $L$ \cite{DWIL}, takes in the two-flavor state sum: at $t_0t_1=-1$ the pairs $(\mathfrak s,\mathfrak s)$ of two equal states contribute Kauffman's state sum for $\Delta_L(t_0^2)$, and the theorem is the statement that all the other pairs cancel (Corollary~\ref{cor.diagonal}).

We keep the notation of Section~\ref{ssec.states}: the Kauffman states $\mathfrak s$ of $(D,\hate)$, the sets $X(\mathfrak s)$, $Y(\mathfrak s)$ and $V(\mathfrak s)$ of Definition~\ref{def.state}, the number $m_{D,\hate}=\ve(X(\mathfrak s))+\ve(Y(\mathfrak s))$ of Proposition~\ref{prop.statesum}~(i), and the number $b(\mathfrak s)$ of black holes of \eqref{eq.blackholes}, i.e.\ of crossings at which the marker of $\mathfrak s$ is the corner $\downarrow$ between the two incoming edges. The Alexander polynomial is normalized as follows. Let $\nabla_L(z)\in\BZ[z]$ be the Conway polynomial of the oriented link $L$, which is characterized by $\nabla_O=1$ for the unknot $O$ and by the skein relation
\begin{equation}\label{eq.conwayskein}
   \nabla_{L_+}(z)-\nabla_{L_-}(z)=z\,\nabla_{L_0}(z)
\end{equation}
for three links which have diagrams differing at one crossing only, positive for $L_+$, negative for $L_-$ and smoothed according to the orientation for $L_0$ \cite{KauffmanFKT}, \cite[Chapter~8]{Lickorish}; it vanishes on split links, and $\nabla_L(-z)=(-1)^{\mu-1}\nabla_L(z)$. Then $\Delta_L(t)=\nabla_L(t^{1/2}-t^{-1/2})$, so that
\[
   \Delta_L(t_0^2)=\nabla_L(t_0-t_0^{-1})
\]
is a Laurent polynomial in $t_0$, which is $t_0-t_0^{-1}$ for the positive Hopf link; for a knot $K$, $\Delta_K$ is the Alexander polynomial with $\Delta_K(t)=\Delta_K(t^{-1})$ and $\Delta_K(1)=1$. It is with this normalization that the specialization displayed above holds for the links with an even number of components, where the sign matters.

\begin{proposition}[the signs of the states]\label{prop.signs}
Let $D$ be a diagram of a link $L$ with $\mu$ components and let $\hate$ be an edge of $D$. In the following sums $\mathfrak s$ and $\mathfrak s'$ run over the Kauffman states of $(D,\hate)$.
\begin{enumerate}\renewcommand{\theenumi}{\roman{enumi}}
\item $|V(\mathfrak s)|\equiv m_{D,\hate}\equiv\mu-1$ modulo $2$ for every state $\mathfrak s$.
\item The signs $\beta_a$ and $\beta_b$ of the proof of Proposition~\ref{prop.statesum}~\textup{(ii)} are $(-1)^{\mu-1}$ and $1$: the state sums of the two edge polynomials are
\[
   P_{D,\hate}(x)=(-1)^{\mu-1}\sum_{\mathfrak s}(-1)^{b(\mathfrak s)}\,x^{X(\mathfrak s)},\qquad \tP_{D,\hate}(y)=\sum_{\mathfrak s}(-1)^{b(\mathfrak s)}\,y^{Y(\mathfrak s)} .
\]
\item Consequently the state sum for $\LG$ of Proposition~\ref{prop.statesum}~\textup{(iii)} holds with no global sign.
\item At the Alexander point,
\[
   \sum_{\mathfrak s}(-1)^{b(\mathfrak s)}\,t^{\ve(X(\mathfrak s))-\ve(Y(\mathfrak s))}=\nabla_L(t-t^{-1}),
\]
that is, $P_{D,\hate}\bigl((t^2)^{\ve}\bigr)=\tP_{D,\hate}\bigl((t^2)^{\ve}\bigr)=(-1)^{\mu-1}\,t^{m_{D,\hate}}\,\nabla_L(t-t^{-1})$.
\end{enumerate}
\end{proposition}

The left-hand side of (iv) is Kauffman's state sum: a state has the weight $t^{\ve_j}$ at the crossings where its marker lies between the two outgoing edges, $-t^{-\ve_j}$ at those where it lies between the two incoming ones, and $1$ at the others. So (iv) is, in our conventions, Kauffman's theorem that this state sum is the Conway polynomial \cite{KauffmanFKT}, and the proof given below is the skein-theoretic one. Parts (i)--(iii) concern the crossing variables and the two flavors. For a knot the proposition says that both state sums hold without a global sign, that $m_{D,\hate}$ is even, and that $P_{D,\hate}(t^{\ve})=\tP_{D,\hate}(t^{\ve})=t^{m_{D,\hate}/2}\,\Delta_K(t)$.

\begin{figure}[ht]
\centering
\begin{tikzpicture}[>=Stealth,scale=1.0,every node/.style={font=\small}]
\begin{scope}
\draw (1.2,-1.2) -- (0.22,-0.22); \draw[->] (-0.22,0.22) -- (-1.2,1.2);
\draw[->] (-1.2,-1.2) -- (1.2,1.2);
\node[below left,font=\scriptsize] at (-1.1,-1.1) {$o$}; \node[above right,font=\scriptsize] at (1.1,1.1) {$o'$};
\node[below right,font=\scriptsize] at (1.1,-1.1) {$u$}; \node[above left,font=\scriptsize] at (-1.1,1.1) {$u'$};
\node at (0,0.75) {$\uparrow$}; \node at (0,-0.75) {$\downarrow$};
\node at (-0.8,0) {$A$}; \node at (0.8,0) {$B$};
\node[font=\scriptsize] at (0,1.3) {$t$}; \node[font=\scriptsize] at (0,-1.3) {$-t^{-1}$};
\node[font=\scriptsize] at (-1.4,0) {$1$}; \node[font=\scriptsize] at (1.4,0) {$1$};
\node at (0,-2.2) {$D$};
\end{scope}
\begin{scope}[xshift=4.7cm]
\draw (-1.2,-1.2) -- (-0.22,-0.22); \draw[->] (0.22,0.22) -- (1.2,1.2);
\draw[->] (1.2,-1.2) -- (-1.2,1.2);
\node[below left,font=\scriptsize] at (-1.1,-1.1) {$o$}; \node[above right,font=\scriptsize] at (1.1,1.1) {$o'$};
\node[below right,font=\scriptsize] at (1.1,-1.1) {$u$}; \node[above left,font=\scriptsize] at (-1.1,1.1) {$u'$};
\node at (0,0.75) {$\uparrow$}; \node at (0,-0.75) {$\downarrow$};
\node at (-0.8,0) {$B$}; \node at (0.8,0) {$A$};
\node[font=\scriptsize] at (0,1.3) {$t^{-1}$}; \node[font=\scriptsize] at (0,-1.3) {$-t$};
\node[font=\scriptsize] at (-1.4,0) {$1$}; \node[font=\scriptsize] at (1.4,0) {$1$};
\node at (0,-2.2) {$D'$};
\end{scope}
\begin{scope}[xshift=9.4cm]
\draw[->] (-1.2,-1.2) to[out=45,in=-45] (-1.2,1.2);
\draw[->] (1.2,-1.2) to[out=135,in=-135] (1.2,1.2);
\node[below left,font=\scriptsize] at (-1.1,-1.1) {$o$}; \node[above right,font=\scriptsize] at (1.1,1.1) {$o'$};
\node[below right,font=\scriptsize] at (1.1,-1.1) {$u$}; \node[above left,font=\scriptsize] at (-1.1,1.1) {$u'$};
\node at (0,0) {$f_0$};
\node[font=\scriptsize] at (0,1.3) {$f^\uparrow_j$}; \node[font=\scriptsize] at (0,-1.3) {$f^\downarrow_j$};
\draw[densely dotted] (-0.62,0) -- (-0.28,0); \draw[densely dotted] (0.28,0) -- (0.62,0);
\node at (0,-2.2) {$D_0$};
\end{scope}
\end{tikzpicture}
\caption{The three diagrams of Lemma~\ref{lem.smoothing} near the crossing $j$, drawn for $\ve_j=1$ in $D$, with the edges $o,u,o',u'$ of $j$ in $D$ and the names of the corners. The smoothing joins the faces $f^\uparrow_j$ and $f^\downarrow_j$ of the two vertical corners into one face $f_0$ of $D_0$, so that the states of $D_0$ are the states of $D$, or of $D'$, whose marker at $j$ is vertical. Next to each corner of $j$, the weight of a marker in Kauffman's state sum of Proposition~\ref{prop.signs}~(iv): the weights in $D$ and $D'$ differ by $\ve_j(t-t^{-1})$ at the vertical corners and are equal at the horizontal ones.}
\label{fig.smoothing}
\end{figure}
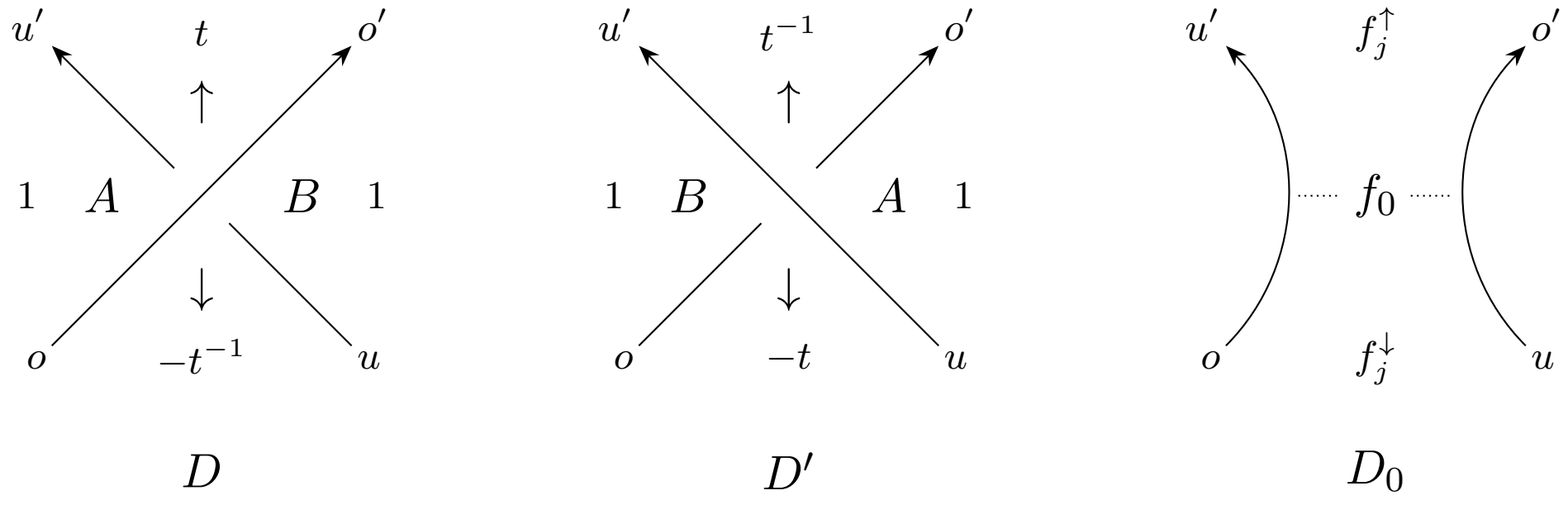

\begin{lemma}[smoothing a crossing]\label{lem.smoothing}
Let $D$ be a diagram with $c\ge2$ crossings and let $j$ be a crossing of $D$ whose four edges are pairwise distinct. Let $D'$ be the diagram obtained from $D$ by changing the crossing $j$, and $D_0$ the one obtained by smoothing $j$ according to the orientation (Figures~\ref{fig.crossing} and~\ref{fig.smoothing}); we assume that $D_0$ is connected. The crossings $k\ne j$ are crossings of the three diagrams, with the same signs and the same corners, and we use the same variables for them.
\begin{enumerate}\renewcommand{\theenumi}{\alph{enumi}}
\item Let $\hate$ be an edge of $D$ and $\hate_0$ the edge of $D_0$ which contains it. The restriction to the crossings $k\ne j$ is a bijection from the states $\mathfrak s$ of $(D,\hate)$ for which $\mathfrak s(j)$ is a vertical corner onto the states of $(D_0,\hate_0)$. The states of $(D',\hate)$ are those of $(D,\hate)$, with the same vertical corners; the $A$-corner of $j$ in $D$ is its $B$-corner in $D'$.
\item If $\hate$ is not one of the four edges of $j$, then
\[
   P_{D_0,\hate}=P_{D',\hate}\big|_{x_j=0}-\frac{\partial P_{D,\hate}}{\partial x_j},\qquad \tP_{D_0,\hate}=\tP_{D',\hate}\big|_{y_j=0}-\frac{\partial\tP_{D,\hate}}{\partial y_j}.
\]
\end{enumerate}
\end{lemma}
\begin{proof}
(a) The smoothing joins the faces of $D$ at the two vertical corners of $j$ and does not change the other faces. By Euler's formula $D$ has $c+2$ faces and the connected diagram $D_0$ has $c+1$; so the faces $f^\uparrow_j$ and $f^\downarrow_j$ are distinct, and the faces of $D_0$ are those of $D$ with these two replaced by their union $f_0$. An edge of $j$ is adjacent to one vertical and to one horizontal corner of $j$ by \eqref{eq.corneredges}; hence, whether $\hate$ is an edge of $j$ or not, the two faces of $D_0$ adjacent to $\hate_0$ are the images of the two faces of $D$ adjacent to $\hate$, and since they are distinct, at most one of $f^\uparrow_j$, $f^\downarrow_j$ is adjacent to $\hate$. Let $\mathfrak s_0$ be a state of $(D_0,\hate_0)$. If $f_0$ is not adjacent to $\hate_0$, exactly one of the corners of $\mathfrak s_0$ lies in $f_0$, hence in one of the faces $f^\uparrow_j$, $f^\downarrow_j$ of $D$, and we let $\mathfrak s(j)$ be the vertical corner of $j$ in the other one. If $f_0$ is adjacent to $\hate_0$, no corner of $\mathfrak s_0$ lies in $f_0$, and we let $\mathfrak s(j)$ be the vertical corner of $j$ in the one of the faces $f^\uparrow_j$, $f^\downarrow_j$ which is not adjacent to $\hate$. With $\mathfrak s(k)=\mathfrak s_0(k)$ for $k\ne j$ this defines a state of $(D,\hate)$, the only one with a vertical corner at $j$ which restricts to $\mathfrak s_0$; and the restriction of such a state of $(D,\hate)$ is a state of $(D_0,\hate_0)$. The last sentence of (a) holds because the faces, the corners and the orientations are those of the curve underlying $D$, while the $A$-corner is the horizontal corner between the incoming over-edge and the outgoing under-edge.

(b) Write $o,u,o',u'$ for the edges $o_j,u_j,o'_j,u'_j$ of $D$; in $D'$ the over-strand of $j$ goes from $u$ to $u'$ and the under-strand from $o$ to $o'$, and in $D_0$ the edges $o$ and $u'$ are parts of one edge, and so are $u$ and $o'$. The rows of the edges other than $u'$ and $o'$ are the same in $1-T^a_D(x)$ and in $1-T^a_{D'}(x)$, because these edges leave crossings other than $j$. For two row vectors $r,r'$ let $[r\,;r']$ be the determinant of the matrix obtained from $(1-T^a_D(x))_{\hate}$ by replacing the row of $u'$ by $r$ and the row of $o'$ by $r'$. By \eqref{eq.Ta}, with $w=e_o-e_u$,
\[
   P_{D,\hate}=\bigl[(e_{u'}-e_o)^{\T}+x_jw^{\T}\,;\,(e_{o'}-e_o)^{\T}\bigr],\qquad P_{D',\hate}=\bigl[(e_{u'}-e_u)^{\T}\,;\,(e_{o'}-e_u)^{\T}-x_jw^{\T}\bigr],
\]
and we claim that
\[
   P_{D_0,\hate}=\bigl[(e_{u'}-e_o)^{\T}\,;\,(e_{o'}-e_u)^{\T}\bigr].
\]
Indeed, in the matrix on the right, adding the column of $u'$ to the column of $o$ and the column of $o'$ to the column of $u$ turns the rows of $u'$ and $o'$ into $e_{u'}^{\T}$ and $e_{o'}^{\T}$; removing these two rows and the columns of $u'$ and $o'$ does not change the determinant, and leaves the matrix $(1-T^a_{D_0}(x))_{\hate}$, in which $o$ and $u$ stand for the two edges of $D_0$ which contain them. Since $[r\,;r']$ is bilinear and vanishes when $r'$ is a multiple of $r$,
\[
   P_{D',\hate}\big|_{x_j=0}=\bigl[(e_{u'}-e_o)^{\T}+w^{\T}\,;\,(e_{o'}-e_u)^{\T}\bigr]=P_{D_0,\hate}+\bigl[w^{\T}\,;\,(e_{o'}-e_u)^{\T}\bigr]
\]
and
\[
   \frac{\partial P_{D,\hate}}{\partial x_j}=\bigl[w^{\T}\,;\,(e_{o'}-e_o)^{\T}\bigr]=\bigl[w^{\T}\,;\,(e_{o'}-e_u)^{\T}-w^{\T}\bigr]=\bigl[w^{\T}\,;\,(e_{o'}-e_u)^{\T}\bigr],
\]
which proves the first identity. The second one is proved in the same way, with the rows and the columns exchanged. By \eqref{eq.Tb} the columns of $o$ and $u$ of $1-T^b(y)$ are $e_o-e_{u'}-y_jv$ and $e_u-e_{u'}$ for $D$, with $v=e_{o'}-e_{u'}$, then $e_o-e_{o'}$ and $e_u-e_{o'}+y_jv$ for $D'$, and $e_o-e_{u'}$ and $e_u-e_{o'}$ for $D_0$, the other columns being the same for $D$ and $D'$: writing $\langle c\,;c'\rangle$ for the determinant with these two columns replaced by $c$ and $c'$,
\begin{align*}
   \tP_{D',\hate}\big|_{y_j=0}&=\langle e_o-e_{u'}-v\,;\,e_u-e_{o'}\rangle=\tP_{D_0,\hate}-\langle v\,;\,e_u-e_{o'}\rangle,\\
   \frac{\partial\tP_{D,\hate}}{\partial y_j}&=-\langle v\,;\,e_u-e_{u'}\rangle=-\langle v\,;\,e_u-e_{o'}+v\rangle=-\langle v\,;\,e_u-e_{o'}\rangle .\qedhere
\end{align*}
\end{proof}

\begin{proof}[Proof of Proposition~\ref{prop.signs}]
\emph{Knots.} Let $D$ be a diagram of a knot $K$, let $e_0=\hate,e_1,\dots,e_{2c-1}$ be the edges of $D$ in the order in which the knot traverses them, and put $e_{2c}=e_0$. At the point where all the variables are $1$, the edge matrices $T^a_D$ and $T^b_D$ are both equal to $\sum_j(E_{o'_jo_j}+E_{u'_ju_j})$ by \eqref{eq.Ta} and \eqref{eq.Tb}, whose row $e_k$ has the single entry $1$ in the column $e_{k-1}$; so the matrices $(1-T^a_D)_{\hate}$ and $(1-T^b_D)_{\hate}$ are lower unitriangular in the order $e_1,\dots,e_{2c-1}$, and both edge polynomials take the value $1$ at this point. By \eqref{eq.PPstates} and \eqref{eq.blackholes} there are signs $\beta_a,\beta_b=\pm1$ such that $P_{D,\hate}(x)=\beta_a\sum_{\mathfrak s}(-1)^{b(\mathfrak s)}x^{X(\mathfrak s)}$ and $\tP_{D,\hate}(y)=\beta_b\sum_{\mathfrak s}(-1)^{b(\mathfrak s)}y^{Y(\mathfrak s)}$. Evaluating at the same point gives $\beta_aB=\beta_bB=1$ with $B=\sum_{\mathfrak s}(-1)^{b(\mathfrak s)}$; hence $B=\pm1$ and $\beta_a=\beta_b=B$, and we have to prove that $B=1$.

The states of $(D,\hate)$ and their black holes depend only on the oriented curve underlying $D$, and not on which strand is the over-strand at each crossing; hence so does $B$, and we may change the crossings of $D$. We make $D$ \emph{descending from $\hate$}: traveling along the knot from $\hate$, every crossing is first met along its over-strand. If the two passages of the knot through the crossing $j$ are from $e_{p-1}$ to $e_p$ and from $e_{p'-1}$ to $e_{p'}$, with $1\le p<p'\le2c$, then $o_j=e_{p-1}$, $o'_j=e_p$, $u_j=e_{p'-1}$ and $u'_j=e_{p'}$; so the entries of $T^b_D(y)$, which by \eqref{eq.Tb} are in the positions $(o'_j,o_j)$, $(u'_j,o_j)$ and $(u'_j,u_j)$, lie strictly below the diagonal in the order $e_1,\dots,e_{2c-1}$. Thus $\tP_{D,\hate}(y)=1$ identically for the descending diagram. On the other hand the constant term of $\tP_{D,\hate}(y)$ is $B$ times the sum of $(-1)^{b(\mathfrak s)}$ over the states with $Y(\mathfrak s)=\emptyset$, and such a state has no corner $\downarrow$, so that $b(\mathfrak s)=0$. So $B$ times a nonnegative integer is $1$, and $B=1$. This is (ii) for knots. The sign $B$ cannot be seen in $\LG$, where it occurs squared, which is why it has to be determined on a diagram for which the edge polynomials are known.

By (ii) and Proposition~\ref{prop.statesum}~(i),
\[
   \tP_{D,\hate}(t^{\ve})=\sum_{\mathfrak s}(-1)^{b(\mathfrak s)}\,t^{\ve(Y(\mathfrak s))}=t^{m_{D,\hate}}\sum_{\mathfrak s}(-1)^{b(\mathfrak s)}\,t^{-\ve(X(\mathfrak s))}=t^{m_{D,\hate}}\,P_{D,\hate}(t^{-\ve}),
\]
and $\tP_{D,\hate}(t^{\ve})=P_{D,\hate}(t^{\ve})$ by Lemma~\ref{lem.level}~(i). By Lemma~\ref{lem.alexpoly}, $P_{D,\hate}(t^{\ve})=\pm t^a\Delta_K(t)$ for some integer $a$; the symmetry just obtained gives $m_{D,\hate}=2a$, and the value $1$ at $t=1$ gives the sign $+$. So $m_{D,\hate}$ is even, which is (i) for knots since $|V(\mathfrak s)|\equiv m_{D,\hate}$ modulo $2$ (Proposition~\ref{prop.statesum}~(i)), and $P_{D,\hate}(t^{\ve})=t^{m_{D,\hate}/2}\Delta_K(t)$, which is (iv) for knots.

\emph{Links: (i) and (ii).} We argue by induction on $\mu$. Let $D_0$ be a diagram of a link with $\mu+1\ge2$ components and $\hate$ an edge of $D_0$, and assume (i) and (ii) for the links with $\mu$ components. Since $D_0$ is connected, it has a crossing $k$ whose two strands belong to different components. No edge both leaves and enters $k$: such an edge would be a whole component, which would meet the rest of the diagram at $k$ only, while two components cross an even number of times. So $\hate$ is not one of the two edges which leave $k$, or it is not one of the two edges which enter $k$. In the first case let $D_+$ and $D_-$ be the diagrams obtained from $D_0$ by making the two edges which leave $k$ cross once more, next to $k$, at a new crossing $j$ which is positive in $D_+$ and negative in $D_-$ (Figure~\ref{fig.fusion}); in the second case we do the same with the two edges which enter $k$. Thus $D_+$ and $D_-$ differ at the crossing $j$, and $D_0$ is the smoothing of both at $j$. They are diagrams of links with $\mu$ components, because the new crossing exchanges two strands which belong to different components of $D_0$. Each of the two edges of $D_0$ is cut by $j$ into two edges of $D_\pm$; these are the four edges of $j$, they are pairwise distinct, and $\hate$, which is an edge of $D_\pm$, is not one of them. So Lemma~\ref{lem.smoothing} applies to $D=D_+$ and $D'=D_-$.

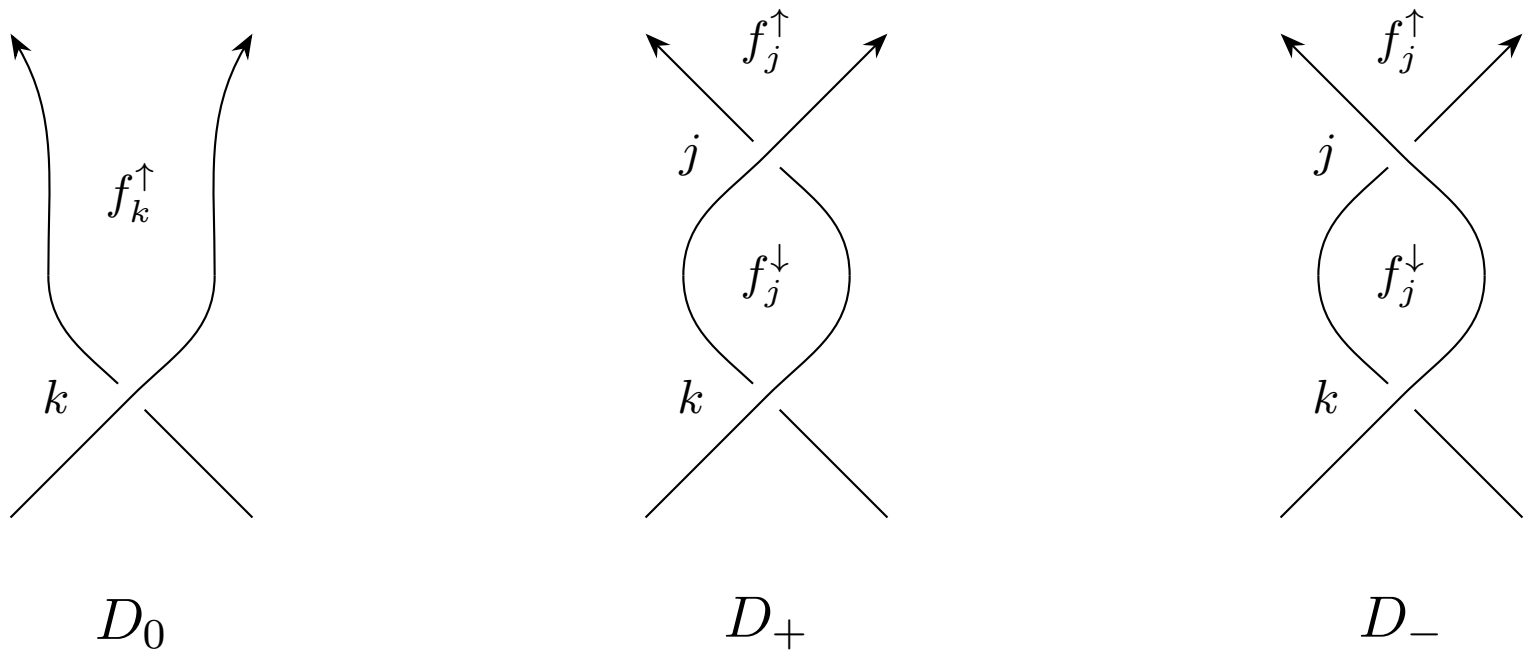
\begin{figure}[ht]
\centering
\begin{tikzpicture}[>=Stealth,scale=1.0,every node/.style={font=\small}]
\begin{scope}
\draw (0.8,-0.8) -- (0,0) to[out=135,in=-90] (-0.55,0.8);
\begin{scope}\clip (0,0) circle (0.3); \draw[white,line width=7pt] (-0.8,-0.8) -- (0,0) to[out=45,in=-90] (0.55,0.8);\end{scope}
\draw (-0.8,-0.8) -- (0,0) to[out=45,in=-90] (0.55,0.8);
\draw[->] (-0.55,0.8) to[out=90,in=-60] (-0.8,2.4);
\draw[->] (0.55,0.8) to[out=90,in=-120] (0.8,2.4);
\node[font=\scriptsize] at (-0.5,0) {$k$};
\node[font=\scriptsize] at (0,1.35) {$f^\uparrow_k$};
\node at (0,-1.5) {$D_0$};
\end{scope}
\begin{scope}[xshift=4.2cm]
\draw (0.8,-0.8) -- (0,0) to[out=135,in=-90] (-0.55,0.8);
\begin{scope}\clip (0,0) circle (0.3); \draw[white,line width=7pt] (-0.8,-0.8) -- (0,0) to[out=45,in=-90] (0.55,0.8);\end{scope}
\draw (-0.8,-0.8) -- (0,0) to[out=45,in=-90] (0.55,0.8);
\draw[->] (0.55,0.8) to[out=90,in=-45] (0,1.6) -- (-0.8,2.4);
\begin{scope}\clip (0,1.6) circle (0.3); \draw[white,line width=7pt] (-0.55,0.8) to[out=90,in=-135] (0,1.6) -- (0.8,2.4);\end{scope}
\draw[->] (-0.55,0.8) to[out=90,in=-135] (0,1.6) -- (0.8,2.4);
\node[font=\scriptsize] at (-0.5,0) {$k$}; \node[font=\scriptsize] at (-0.5,1.6) {$j$};
\node[font=\scriptsize] at (0,0.8) {$f^\downarrow_j$}; \node[font=\scriptsize] at (0,2.35) {$f^\uparrow_j$};
\node at (0,-1.5) {$D_+$};
\end{scope}
\begin{scope}[xshift=8.4cm]
\draw (0.8,-0.8) -- (0,0) to[out=135,in=-90] (-0.55,0.8);
\begin{scope}\clip (0,0) circle (0.3); \draw[white,line width=7pt] (-0.8,-0.8) -- (0,0) to[out=45,in=-90] (0.55,0.8);\end{scope}
\draw (-0.8,-0.8) -- (0,0) to[out=45,in=-90] (0.55,0.8);
\draw[->] (-0.55,0.8) to[out=90,in=-135] (0,1.6) -- (0.8,2.4);
\begin{scope}\clip (0,1.6) circle (0.3); \draw[white,line width=7pt] (0.55,0.8) to[out=90,in=-45] (0,1.6) -- (-0.8,2.4);\end{scope}
\draw[->] (0.55,0.8) to[out=90,in=-45] (0,1.6) -- (-0.8,2.4);
\node[font=\scriptsize] at (-0.5,0) {$k$}; \node[font=\scriptsize] at (-0.5,1.6) {$j$};
\node[font=\scriptsize] at (0,0.8) {$f^\downarrow_j$}; \node[font=\scriptsize] at (0,2.35) {$f^\uparrow_j$};
\node at (0,-1.5) {$D_-$};
\end{scope}
\end{tikzpicture}
\caption{The induction step in the proof of Proposition~\ref{prop.signs}. The two strands of the crossing $k$ of $D_0$, here drawn positive, belong to different components; the two edges which leave $k$ are made to cross once more, at a new crossing $j$, positive in $D_+$ and negative in $D_-$. The face $f^\uparrow_k$ of $D_0$ is cut into the faces $f^\downarrow_j$ and $f^\uparrow_j$ of $D_\pm$, the diagram $D_0$ is the smoothing of $D_\pm$ at $j$, and $D_\pm$ has one component less than $D_0$.}
\label{fig.fusion}
\end{figure}

For a state $\mathfrak s$ of $(D_+,\hate)$ let $b_j(\mathfrak s)$ be the number of crossings $k\ne j$ with $\mathfrak s(k)={\downarrow}$, and for a corner $\diamond$ of $j$ put
\[
   \Sigma^\diamond=\sum_{\mathfrak s(j)=\diamond}(-1)^{b_j(\mathfrak s)}\,x^{X(\mathfrak s)\sm\{j\}},\qquad \widetilde\Sigma^\diamond=\sum_{\mathfrak s(j)=\diamond}(-1)^{b_j(\mathfrak s)}\,y^{Y(\mathfrak s)\sm\{j\}},
\]
the sums running over the states $\mathfrak s$ of $(D_+,\hate)$ with $\mathfrak s(j)=\diamond$; here $A$ and $B$ denote the horizontal corners of $j$ in $D_+$. By Definition~\ref{def.state}, by Lemma~\ref{lem.smoothing}~(a) and by the induction hypothesis,
\begin{align*}
   (-1)^{\mu-1}P_{D_+,\hate}&=x_j\bigl(\Sigma^\uparrow+\Sigma^A\bigr)+\Sigma^B-\Sigma^\downarrow, & (-1)^{\mu-1}P_{D_-,\hate}&=x_j\bigl(\Sigma^\uparrow+\Sigma^B\bigr)+\Sigma^A-\Sigma^\downarrow,\\
   \tP_{D_+,\hate}&=y_j\bigl(\widetilde\Sigma^A-\widetilde\Sigma^\downarrow\bigr)+\widetilde\Sigma^B+\widetilde\Sigma^\uparrow, & \tP_{D_-,\hate}&=y_j\bigl(\widetilde\Sigma^B-\widetilde\Sigma^\downarrow\bigr)+\widetilde\Sigma^A+\widetilde\Sigma^\uparrow .
\end{align*}
Hence, by Lemma~\ref{lem.smoothing}~(b),
\begin{align*}
   P_{D_0,\hate}&=(-1)^{\mu-1}\Bigl(\bigl(\Sigma^A-\Sigma^\downarrow\bigr)-\bigl(\Sigma^\uparrow+\Sigma^A\bigr)\Bigr)=(-1)^{\mu}\bigl(\Sigma^\uparrow+\Sigma^\downarrow\bigr),\\
   \tP_{D_0,\hate}&=\bigl(\widetilde\Sigma^A+\widetilde\Sigma^\uparrow\bigr)-\bigl(\widetilde\Sigma^A-\widetilde\Sigma^\downarrow\bigr)=\widetilde\Sigma^\uparrow+\widetilde\Sigma^\downarrow .
\end{align*}
By Lemma~\ref{lem.smoothing}~(a) the states of $(D_+,\hate)$ with a vertical corner at $j$ are the states $\mathfrak s_0$ of $(D_0,\hate)$, with $b_j(\mathfrak s)=b(\mathfrak s_0)$, $X(\mathfrak s)\sm\{j\}=X(\mathfrak s_0)$ and $Y(\mathfrak s)\sm\{j\}=Y(\mathfrak s_0)$; so $\Sigma^\uparrow+\Sigma^\downarrow=\sum_{\mathfrak s_0}(-1)^{b(\mathfrak s_0)}x^{X(\mathfrak s_0)}$ and $\widetilde\Sigma^\uparrow+\widetilde\Sigma^\downarrow=\sum_{\mathfrak s_0}(-1)^{b(\mathfrak s_0)}y^{Y(\mathfrak s_0)}$, which is (ii) for $(D_0,\hate)$. Moreover $|V(\mathfrak s_0)|=|V(\mathfrak s)|-1\equiv\mu$ modulo $2$, and $(D_0,\hate)$ has at least one state \cite{KauffmanFKT}; this is (i).

\emph{(iii)} follows from (ii), Theorem~\ref{thm.HK} and \eqref{eq.ED}, the two signs $(-1)^{\mu-1}$ canceling.

\emph{Links: (iv).} Put $\Phi_{D,\hate}(t)=\sum_{\mathfrak s}(-1)^{b(\mathfrak s)}\,t^{\ve(X(\mathfrak s))-\ve(Y(\mathfrak s))}$. Since $\ve(X(\mathfrak s))-\ve(Y(\mathfrak s))=2\ve(X(\mathfrak s))-m_{D,\hate}=m_{D,\hate}-2\ve(Y(\mathfrak s))$, (ii) gives
\begin{equation}\label{eq.omegaP}
   P_{D,\hate}\bigl((t^2)^{\ve}\bigr)=(-1)^{\mu-1}\,t^{m_{D,\hate}}\,\Phi_{D,\hate}(t),\qquad \tP_{D,\hate}\bigl((t^2)^{\ve}\bigr)=t^{m_{D,\hate}}\,\Phi_{D,\hate}(t^{-1}),
\end{equation}
and the two are equal by Lemma~\ref{lem.level}~(i). So it suffices to prove that $\Phi_{D,\hate}(t)=\nabla_L(t-t^{-1})$, which for knots was proved above. A state contributes to $\Phi_{D,\hate}$ the product of the weights $t^{\ve_j}$, $-t^{-\ve_j}$ and $1$ of its markers. Let $j$ be a crossing of $D$ whose two strands belong to different components, and let $D'$ and $D_0$ be as in Lemma~\ref{lem.smoothing}; the four edges of $j$ are pairwise distinct, as above, and $D_0$ is connected: otherwise the two vertical corners of $j$ would lie in the same face of $D$, by the count of the faces in the proof of Lemma~\ref{lem.smoothing}, and a closed curve through these two corners which meets $D$ at $j$ only would meet each of the two components exactly once. The states of $D$ and $D'$ have the same weights at the crossings $k\ne j$, and their weights at $j$ differ by $\ve_j(t-t^{-1})$ if the marker of $j$ is vertical, and are both $1$ otherwise. By Lemma~\ref{lem.smoothing}~(a),
\[
   \Phi_{D,\hate}-\Phi_{D',\hate}=\ve_j\,(t-t^{-1})\,\Phi_{D_0,\hate_0} .
\]
Now let $\mu\ge2$ and assume (iv) for the links with $\mu-1$ components. Then $\Phi_{D_0,\hate_0}(t)=\nabla_{L_0}(t-t^{-1})$, where $L_0$ is the link of $D_0$, and comparing with \eqref{eq.conwayskein} shows that $\Phi_{D,\hate}(t)-\nabla_L(t-t^{-1})$ does not change when a crossing between two different components is changed. Changing the crossings at which one given component passes under another one, we arrive at a diagram of a split link $L''$. For this diagram both terms vanish: $\nabla_{L''}=0$, and $\Phi_{D,\hate}=0$ by \eqref{eq.omegaP} and Lemma~\ref{lem.alexpoly}, since the Alexander polynomial of a split link is zero. Hence $\Phi_{D,\hate}(t)=\nabla_L(t-t^{-1})$.
\end{proof}

\begin{corollary}[the state sum at $t_0t_1=-1$]\label{cor.diagonal}
Let $D$ be a diagram of a link $L$ and $\hate$ an edge of $D$. The pairs $(\mathfrak s,\mathfrak s)$ of two equal states contribute
\[
   \sum_{\mathfrak s}\ \prod_{j:\,\mathfrak s(j)={\uparrow}}t_0^{\ve_j}\prod_{j:\,\mathfrak s(j)={\downarrow}}t_1^{\ve_j}
\]
to the state sum of Proposition~\ref{prop.statesum}~(iii), and at $t_1=-t_0^{-1}$ this contribution is $\nabla_L(t_0-t_0^{-1})=\Delta_L(t_0^2)$. Consequently the identity $\LG(L;t_0,-t_0^{-1})=\Delta_L(t_0^2)$ is equivalent to the cancellation of the pairs of two distinct states at $t_0t_1=-1$:
\[
   \sum_{\mathfrak s\ne\mathfrak s'}(-1)^{b(\mathfrak s)+b(\mathfrak s')}\prod_{j\in X(\mathfrak s)\sm Y(\mathfrak s')}t_0^{\ve_j}\ \prod_{j\in Y(\mathfrak s')\sm X(\mathfrak s)}\bigl(-t_0^{-1}\bigr)^{\ve_j}=0 .
\]
By \cite{DWIL} this holds for every link.
\end{corollary}
\begin{proof}
For $\mathfrak s'=\mathfrak s$ the sign is $1$, and $X(\mathfrak s)\sm Y(\mathfrak s)$ and $Y(\mathfrak s)\sm X(\mathfrak s)$ are the sets of crossings with $\mathfrak s(j)={\uparrow}$ and with $\mathfrak s(j)={\downarrow}$ (Definition~\ref{def.state}). Only the vertical corners carry a weight, $t_0^{\ve_j}$ between the two outgoing edges and $t_1^{\ve_j}$ between the two incoming ones, and all the coefficients are positive. Since $(-t_0^{-1})^{\ve_j}=-t_0^{-\ve_j}$, the pair $(\mathfrak s,\mathfrak s)$ contributes $(-1)^{b(\mathfrak s)}\,t_0^{\ve(X(\mathfrak s))-\ve(Y(\mathfrak s))}$ at $t_1=-t_0^{-1}$, and the sum of these terms is $\nabla_L(t_0-t_0^{-1})$ by Proposition~\ref{prop.signs}~(iv). The rest follows, the whole sum being $\LG(L;t_0,-t_0^{-1})$.
\end{proof}

So at $t_0t_1=-1$ the diagonal of the state sum is the whole sum, and it is Kauffman's state sum for the Conway-normalized Alexander polynomial, in the variable $t_0^2$. In this sense the two specializations of $\LG$ are the two extreme regimes of the state sum: at $t_0t_1=1$ the two states of a pair are independent (Remark~\ref{rem.statesum}~(ii)), and at $t_0t_1=-1$ they coincide. In Corollary~\ref{cor.diagonal} the cancellation is a consequence of the theorem of \cite{DWIL}; a proof of the cancellation in the state model itself would give another proof of that theorem. Corollary~\ref{cor.ballinger} gives such a proof: the pairs of two distinct states cancel within each superposition (Remark~\ref{rem.ballinger}~(ii)).

\begin{example}[the trefoil and the Hopf link]\label{ex.signs}
For the standard diagram of the trefoil, Example~\ref{ex.states} gives $\LG(3_1;t_0,t_1)=1-(t_0+t_1)(1-t_0)(1-t_1)$, in which the three pairs $(\mathfrak s,\mathfrak s)$ contribute $t_0^2+t_0t_1+t_1^2$, which is $t_0^2-1+t_0^{-2}=\Delta_{3_1}(t_0^2)$ at $t_1=-t_0^{-1}$, and the six other pairs contribute $(1+t_0t_1)(1-t_0-t_1)$.

Let $D$ be the closure of the $2$-braid $\sigma_1^2$, a diagram of the positive Hopf link with two positive crossings, and $\hate=o_1=u'_2$. There are two states, with the corners $(B,\downarrow)$ and $(\uparrow,B)$ at the crossings $1,2$; so $X(\mathfrak s)=\emptyset,\{1\}$, $Y(\mathfrak s)=\{2\},\emptyset$ and $b(\mathfrak s)=1,0$, and $m_{D,\hate}=1$. Proposition~\ref{prop.signs}~(ii) reads $P_{D,\hate}(x)=-(-1+x_1)$ and $\tP_{D,\hate}(y)=-y_2+1$, and (iii) reads
\[
   \LG(L;t_0,t_1)=t_0+t_1-1-t_0t_1=-(1-t_0)(1-t_1):
\]
the two pairs $(\mathfrak s,\mathfrak s)$ contribute $t_0+t_1$, which is $t_0-t_0^{-1}=\nabla_L(t_0-t_0^{-1})$ at $t_1=-t_0^{-1}$, and the two other pairs contribute $-(1+t_0t_1)$.
\end{example}

\subsection{Kauffman states with multiplicity: Ballinger's model}\label{ssec.ballinger}
In the state sum of Proposition~\ref{prop.statesum}~(iii) the two states of a pair play different roles, and the weight of a pair is not symmetric in them. If the pairs are grouped according to the two markers which they put at every crossing, whichever of the two states puts them, the sum takes a symmetric form: its terms are indexed by Kauffman states with multiplicity two, and the interaction of the two flavors becomes a factor $1+t_0t_1$ for every cycle along which the two states of a pair differ. In this form the state sum is the case $N=2$ of a state model due to William Ballinger \cite{BallingerTalk}, which we first describe.

For every $N\ge1$, Ballinger's model is a sum over the \emph{Kauffman states with multiplicity $N$} of $(D,\hate)$: the assignments of a multiplicity in $\{0,\dots,N\}$ to every corner of every crossing such that the multiplicities add up to $N$ at every crossing and in every face other than $f_{\rm L}(\hate)$ and $f_{\rm R}(\hate)$, and vanish in these two faces. The corners of positive multiplicity, labeled by their multiplicities, are the edges of a plane graph, and a state is weighted by the evaluation of this graph as a web of $\mathfrak{sl}_N$ and by a product of local weights at the crossings, which are monomials in two variables $s$ and $q$. For $N=1$ the states are the Kauffman states, and the sum is Kauffman's state sum for the Alexander polynomial. According to \cite{BallingerTalk} the sum is an invariant $Q_N(L)$ of the link, which is conjectured there to be the Links--Gould invariant $\LG^{N,1}$ of $U_q\mathfrak{gl}(N|1)$ after a change of variables. We refer to \cite{BallingerTalk} for the general definition and spell it out for $N=2$, where it is elementary.

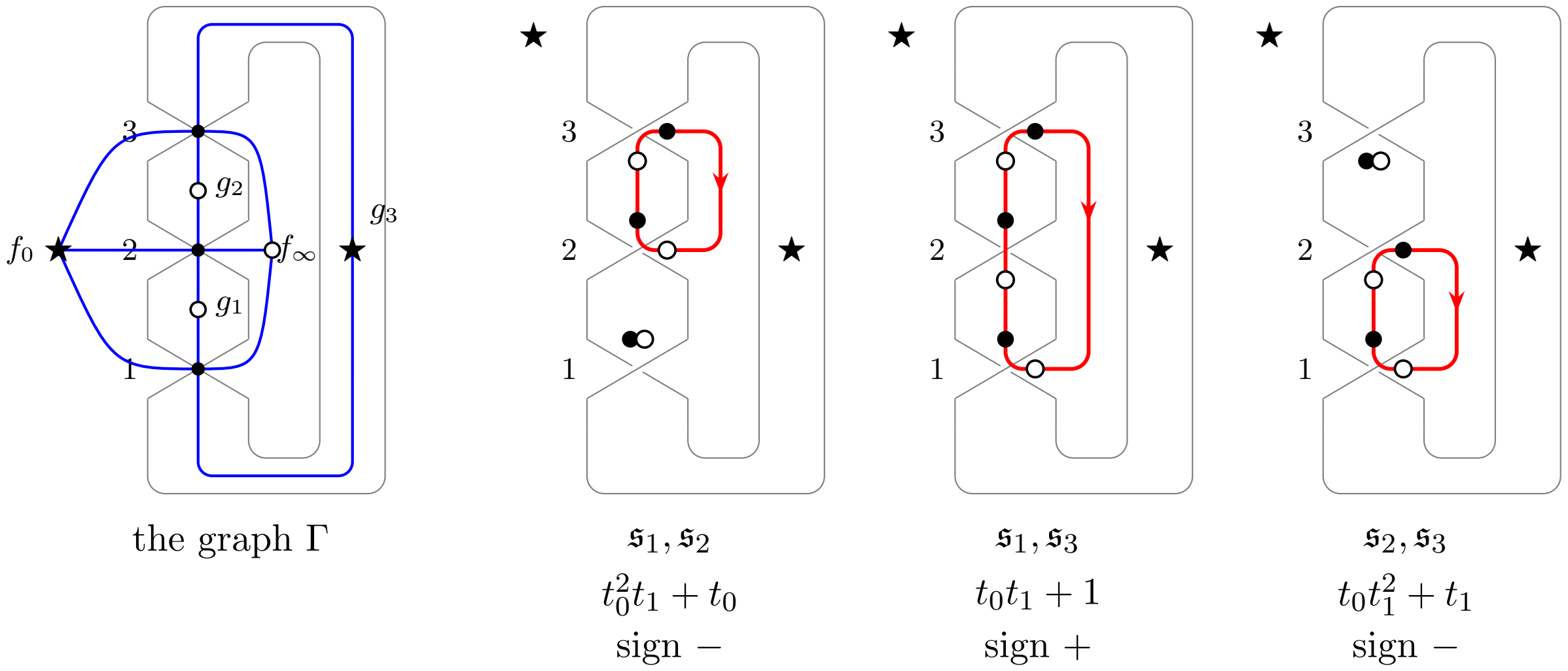
\begin{figure}[ht]
\centering
\begin{tikzpicture}[>=Stealth,scale=0.6]
\def\trefoilbase{%
   \draw[gray] (0,0) -- (0,0.5) (0,1.5) -- (0,2.5) (0,3.5) -- (0,4.5) (0,5.5) -- (0,6);
   \draw[gray,rounded corners=5pt] (0,6) -- (0,7.1) -- (4.0,7.1) -- (4.0,-1.1) -- (0,-1.1) -- (0,0);
   \draw[gray] (1.7,0) -- (1.7,0.5) (1.7,1.5) -- (1.7,2.5) (1.7,3.5) -- (1.7,4.5) (1.7,5.5) -- (1.7,6);
   \draw[gray,rounded corners=5pt] (1.7,6) -- (1.7,6.5) -- (2.9,6.5) -- (2.9,-0.5) -- (1.7,-0.5) -- (1.7,0);
   \foreach \y in {1,3,5} { \draw[gray] (1.7,{\y-0.5}) -- (0,{\y+0.5}); \draw[gray,preaction={draw=white,line width=3pt}] (0,{\y-0.5}) -- (1.7,{\y+0.5}); }
   \foreach \y/\j in {1/1,3/2,5/3} \node[font=\scriptsize] at (-0.3,\y) {$\j$};}
\def\stars{\node[font=\scriptsize] at (-0.9,6.6) {$\bigstar$}; \node[font=\scriptsize] at (3.45,3) {$\bigstar$};}
\def\mka#1#2{\fill[black] (0.85,#1) ++(#2) circle (4pt);}
\def\mkb#1#2{\draw[black,fill=white,line width=0.7pt] (0.85,#1) ++(#2) circle (4pt);}
\begin{scope}
\trefoilbase
\coordinate (f0) at (-1.5,3); \coordinate (g1) at (0.85,2); \coordinate (g2) at (0.85,4); \coordinate (g3) at (3.45,3); \coordinate (fi) at (2.1,3);
\begin{scope}[blue,line width=0.8pt]
\foreach \y in {1,3,5} { \draw (0.85,\y) .. controls (-0.6,\y) .. (f0); \draw (0.85,\y) .. controls (1.9,\y) .. (fi); }
\draw (0.85,1) -- (g1) -- (0.85,3) -- (g2) -- (0.85,5);
\draw[rounded corners=4pt] (0.85,5) -- (0.85,6.8) -- (3.45,6.8) -- (g3);
\draw[rounded corners=4pt] (0.85,1) -- (0.85,-0.8) -- (3.45,-0.8) -- (g3);
\end{scope}
\foreach \y in {1,3,5} \fill (0.85,\y) circle (3.2pt);
\foreach \p in {g1,g2,fi} \draw[fill=white,line width=0.7pt] (\p) circle (3.6pt);
\node[font=\scriptsize] at (f0) {$\bigstar$}; \node[font=\scriptsize] at (g3) {$\bigstar$};
\node[font=\scriptsize,left] at (-1.65,3) {$f_0$}; \node[font=\scriptsize,right] at (3.5,3.6) {$g_3$};
\node[font=\scriptsize,right] at (0.9,2.05) {$g_1$}; \node[font=\scriptsize,right] at (0.9,4.05) {$g_2$}; \node[font=\scriptsize] at (2.52,3.02) {$f_\infty$};
\node[font=\small] at (1.4,-1.9) {the graph $\Gamma$};
\end{scope}
\begin{scope}[xshift=7.4cm]
\trefoilbase \stars
\draw[red,line width=1.1pt,rounded corners=5pt,postaction={decorate,decoration={markings,mark=at position 0.62 with {\arrow{Stealth[length=6pt]}}}}] (0.85,3) -- (0.85,5) -- (2.25,5) -- (2.25,3) -- cycle;
\mka{1}{-0.12,0.5} \mkb{1}{0.12,0.5}
\mka{3}{0,0.5} \mkb{3}{0.5,0}
\mka{5}{0.5,0} \mkb{5}{0,-0.5}
\node[font=\small] at (1.4,-1.9) {$\mathfrak s_1,\mathfrak s_2$};
\node[font=\small] at (1.4,-2.8) {$t_0^2t_1+t_0$};
\node[font=\small] at (1.4,-3.7) {sign $-$};
\end{scope}
\begin{scope}[xshift=13.6cm]
\trefoilbase \stars
\draw[red,line width=1.1pt,rounded corners=5pt,postaction={decorate,decoration={markings,mark=at position 0.62 with {\arrow{Stealth[length=6pt]}}}}] (0.85,1) -- (0.85,5) -- (2.25,5) -- (2.25,1) -- cycle;
\mka{1}{0,0.5} \mkb{1}{0.5,0}
\mka{3}{0,0.5} \mkb{3}{0,-0.5}
\mka{5}{0.5,0} \mkb{5}{0,-0.5}
\node[font=\small] at (1.4,-1.9) {$\mathfrak s_1,\mathfrak s_3$};
\node[font=\small] at (1.4,-2.8) {$t_0t_1+1$};
\node[font=\small] at (1.4,-3.7) {sign $+$};
\end{scope}
\begin{scope}[xshift=19.8cm]
\trefoilbase \stars
\draw[red,line width=1.1pt,rounded corners=5pt,postaction={decorate,decoration={markings,mark=at position 0.62 with {\arrow{Stealth[length=6pt]}}}}] (0.85,1) -- (0.85,3) -- (2.25,3) -- (2.25,1) -- cycle;
\mka{1}{0,0.5} \mkb{1}{0.5,0}
\mka{3}{0.5,0} \mkb{3}{0,-0.5}
\mka{5}{-0.12,-0.5} \mkb{5}{0.12,-0.5}
\node[font=\small] at (1.4,-1.9) {$\mathfrak s_2,\mathfrak s_3$};
\node[font=\small] at (1.4,-2.8) {$t_0t_1^2+t_1$};
\node[font=\small] at (1.4,-3.7) {sign $-$};
\end{scope}
\end{tikzpicture}
\caption{Left: the graph $\Gamma$ of the trefoil diagram of Figure~\ref{fig.trefoilstates}. Its vertices are the three crossings and one point in each of the five faces, and its twelve edges are the corners; its faces are quadrilaterals, one disk $R_e$ around each edge $e$ of the diagram. A Kauffman state is a perfect matching of $\Gamma$ with the two starred vertices removed. Right: the three states with multiplicity two which are not the superposition of a state with itself. Each of them is the superposition of two ordered pairs $(\mathfrak s,\mathfrak s')$ and $(\mathfrak s',\mathfrak s)$; for the first pair the markers of $\mathfrak s$ are black and those of $\mathfrak s'$ white. The corners chosen once form one cycle (red), oriented as in Corollary~\ref{cor.ballinger}~\textup{(ii)} for the first pair; the cut edge is on its left, so that the monomial of the second pair is that of the first one divided by $t_0t_1$. Below each picture, the sum of the two monomials, which is a monomial times $1+t_0t_1$, and the sign $(-1)^{b(\mathfrak s)+b(\mathfrak s')}$.}
\label{fig.gammagraph}
\end{figure}

A \emph{state with multiplicity two} of $(D,\hate)$ is a map $x$ which assigns to every crossing $j$ an unordered pair $x(j)$ of corners of $j$, possibly twice the same corner, in such a way that every face other than $f_{\rm L}(\hate)$ and $f_{\rm R}(\hate)$ contains exactly two of the $2c$ chosen corners, and these two faces contain none. The \emph{superposition} of two Kauffman states $\mathfrak s$ and $\mathfrak s'$ is the state with multiplicity two $j\mapsto\{\mathfrak s(j),\mathfrak s'(j)\}$. Let $\Gamma$ be the graph on $S^2$ whose vertices are the crossings of $D$\label{not.Gamma} and one point $v_f$ in every face $f$, and whose edges are the corners (Figure~\ref{fig.gammagraph}): the edge of a corner of the crossing $j$ which lies in the face $f$ is an arc from $j$ to $v_f$ inside $f$, leaving $j$ through this corner, and distinct arcs meet at their ends only. The corners which $x$ chooses exactly once form a subgraph $\Gamma_x$ of $\Gamma$ in which every vertex has degree two or zero, that is, a disjoint union of cycles; let $k(x)$ be the number of these cycles. In Ballinger's model the corners chosen twice carry the trivial representation of $\mathfrak{sl}_2$, and each cycle of $\Gamma_x$ is a closed curve which carries the fundamental representation and evaluates to $[2]=q+q^{-1}$. The local weight of the crossing $j$ is
\begin{equation}\label{eq.ellweight}
   \ell(x,j)=\Bigl((-1)^{d}\,s^{\frac12(a-d)}\,q^{\binom a2+\binom d2+\frac h2}\Bigr)^{\ve_j},
\end{equation}
where $a$ and $d$ are the numbers of times the corners $\uparrow$ and $\downarrow$ occur in $x(j)$, and $h$ is the number of horizontal corners which occur exactly once in $x(j)$. So, for $N=2$,
\begin{equation}\label{eq.Q2}
   Q_2(D,\hate)=\sum_x\ [2]^{k(x)}\,\prod_{j=1}^c\ell(x,j),
\end{equation}
the sum running over the states with multiplicity two of $(D,\hate)$.

\begin{corollary}[Ballinger's model \cite{BallingerTalk}]\label{cor.ballinger}
Let $D$ be a diagram of a link $L$ and $\hate$ an edge of $D$.
\begin{enumerate}\renewcommand{\theenumi}{\roman{enumi}}
\item Every state $x$ with multiplicity two of $(D,\hate)$ is the superposition of exactly $2^{k(x)}$ ordered pairs of Kauffman states. Each cycle of $\Gamma_x$ passes through at least two crossings, and the pairs with superposition $x$ are obtained by distributing the corners of each cycle alternately between the two states, in one of two ways.
\item The pairs $(\mathfrak s,\mathfrak s')$ with superposition $x$ have the same sign $(-1)^{b(\mathfrak s)+b(\mathfrak s')}$, and their monomials in the state sum of Proposition~\ref{prop.statesum}~\textup{(iii)} have the same degree \eqref{eq.KTdegree}. Exchanging the two states along one cycle $\gamma$ of $\Gamma_x$ multiplies the monomial of a pair by $(t_0t_1)^{\pm1}$: if $\gamma$ is oriented so that it enters its crossings through the corners of $\mathfrak s'$ and leaves them through the corners of $\mathfrak s$, the factor is $t_0t_1$ if $\hate$ lies on the right of $\gamma$, and $(t_0t_1)^{-1}$ if it lies on its left.
\item Let $t_0=sq$ and $t_1=s^{-1}q$. The terms of the $2^{k(x)}$ pairs with superposition $x$ in the state sum of $\LG$ add up to $[2]^{k(x)}\prod_j\ell(x,j)$. Consequently
\[
   \LG(L;\,sq,\,s^{-1}q)=\sum_x\ [2]^{k(x)}\,\prod_{j=1}^c\ell(x,j)=Q_2(D,\hate):
\]
for $N=2$, Ballinger's invariant is the Links--Gould invariant.
\end{enumerate}
\end{corollary}

\begin{proof}
For two corners $c,c'$ of a crossing $j$ put
\begin{equation}\label{eq.wpair}
   w_j(c,c')=\begin{cases} t_0^{\ve_j}&\text{if }c\in\{\uparrow,A\}\text{ and }c'\in\{\uparrow,B\},\\ t_1^{\ve_j}&\text{if }c\in\{\downarrow,B\}\text{ and }c'\in\{\downarrow,A\},\\ 1&\text{otherwise.}\end{cases}
\end{equation}
{\sloppy By Definition~\ref{def.state} the monomial of a pair $(\mathfrak s,\mathfrak s')$ in Proposition~\ref{prop.statesum}~(iii) is $\prod_jw_j\bigl(\mathfrak s(j),\mathfrak s'(j)\bigr)$, with the coefficient $(-1)^{b(\mathfrak s)+b(\mathfrak s')}$.\par}

\emph{The pairs with a given superposition.} The corners chosen by $x$, counted with their multiplicities, form a graph on the crossings and on the faces other than $f_{\rm L}(\hate)$ and $f_{\rm R}(\hate)$ in which every vertex has degree two. Its components are the corners chosen twice and the cycles of $\Gamma_x$, which have even length because $\Gamma$ is bipartite. A cycle has exactly two perfect matchings. Choosing one of them on every cycle of $\Gamma_x$ and adding the corners chosen twice gives a Kauffman state $\mathfrak s$; the other matchings and the corners chosen twice form a Kauffman state $\mathfrak s'$; the superposition of $\mathfrak s$ and $\mathfrak s'$ is $x$, and every pair with superposition $x$ is obtained in this way exactly once. So there are $2^{k(x)}$ such pairs. Their sign is $(-1)^{b(x)}$, where $b(x)$ is the number of corners $\downarrow$ chosen by $x$, counted with multiplicity. The degree \eqref{eq.KTdegree} of the monomial of a pair is $\ve(X(\mathfrak s))-\ve(Y(\mathfrak s'))$, which is $\ve(X(\mathfrak s))+\ve(X(\mathfrak s'))-m_{D,\hate}$ by Proposition~\ref{prop.statesum}~(i), and $\ve(X(\mathfrak s))+\ve(X(\mathfrak s'))$ is the sum of the signs $\ve_j$ over the corners $\uparrow$ and $A$ chosen by $x$, counted with multiplicity. So the sign and the degree depend on $x$ only, and the monomials of two pairs with superposition $x$ differ by a power of $t_0t_1$.

\emph{One cycle.} Let $(\mathfrak s,\mathfrak s')$ be a pair with superposition $x$, and let $\gamma$ be a cycle of $\Gamma_x$, of length $2\ell$ with $\ell\ge1$. It passes through crossings $j_1,\dots,j_\ell$ and through the points of faces $F_1,\dots,F_\ell$, which we number so that the corner $c_k=\mathfrak s(j_k)$ lies in $F_k$ and the corner $c'_k=\mathfrak s'(j_k)$ in $F_{k-1}$, the indices being taken modulo $\ell$; and we orient $\gamma$ from $F_{k-1}$ through $j_k$ to $F_k$, so that it enters $j_k$ through $c'_k$ and leaves it through $c_k$. Then $\gamma$ is a simple closed curve on $S^2$ which meets $D$ at the crossings $j_k$ only. It does not meet $\hate$; let $\Theta$ be the component of $S^2\sm\gamma$ which does not contain $\hate$. Let $m$ be the number of crossings in $\Theta$, $n$ the number of faces contained in $\Theta$, and $h_k$ the number of ends of edges of $D$ at $j_k$ which lie in $\Theta$. A crossing in $\Theta$ has its corners in faces which are contained in $\Theta$ or are among the $F_k$, and $\mathfrak s$ maps the crossings of $\gamma$ onto the faces $F_k$; so $\mathfrak s$ maps the crossings in $\Theta$ injectively into the faces contained in $\Theta$. A face contained in $\Theta$ is not adjacent to $\hate$, so it contains the corner $\mathfrak s(j)$ of some crossing $j$, which lies in $\Theta$ and is not one of the $j_k$. Hence $n=m$. The closed disk $\overline\Theta$ is decomposed into the $m+\ell$ crossings, the $2m+\frac12\sum_kh_k$ edges of $D$ in $\Theta$ and the $\ell$ arcs of $\gamma$ between its crossings, and the $n$ faces and the $\ell$ parts of the faces $F_k$ which lie in $\Theta$; since $D$ is connected, the union of $\gamma$ and of the part of $D$ in $\overline\Theta$ is connected, and the Euler characteristic of the disk gives
\begin{equation}\label{eq.eulercycle}
   1=(m+\ell)-\bigl(2m+\tfrac12\textstyle\sum_kh_k+\ell\bigr)+(n+\ell),\qquad\text{that is,}\qquad \textstyle\sum_kh_k=2\ell-2 .
\end{equation}
For two different corners $c',c$ of a crossing let $\tau(c',c)$ be $1$ if $c$ is the counterclockwise neighbor of $c'$, $-1$ if it is its clockwise neighbor, and $0$ if the two corners are opposite. A curve which enters a crossing through $c'$ and leaves it through $c$ has $2+\tau(c',c)$ of the four ends of edges at this crossing on its left, and $2-\tau(c',c)$ on its right (Figure~\ref{fig.turningnumber}). In particular $1\le h_k\le3$, so that \eqref{eq.eulercycle} excludes $\ell=1$: a cycle of $\Gamma_x$ passes through at least two crossings, which completes the proof of (i). Moreover $h_k=2+\tau(c'_k,c_k)$ if $\Theta$ is on the left of $\gamma$, and $h_k=2-\tau(c'_k,c_k)$ if it is on its right, so that by \eqref{eq.eulercycle}
\begin{equation}\label{eq.turning}
   \sum_{k=1}^\ell\tau(c'_k,c_k)=-2\ \text{ if $\hate$ lies on the right of $\gamma$},\qquad \sum_{k=1}^\ell\tau(c'_k,c_k)=2\ \text{ if it lies on its left.}
\end{equation}
Now exchange the two states along $\gamma$. This divides the monomial of the pair by
\[
   R_\gamma=\prod_{k=1}^\ell\frac{w_{j_k}(c_k,c'_k)}{w_{j_k}(c'_k,c_k)},
\]
which is a power of $t_0t_1$ by the first part of the proof. An inspection of \eqref{eq.wpair} shows that for two different corners $c',c$ of a crossing $j$, of either sign, the quotient $w_j(c,c')/w_j(c',c)$ is $t_0$ or $t_1$ if $\tau(c',c)=1$, according as the end of edge between the two corners is outgoing or incoming; it is the inverse of $t_0$ or of $t_1$ if $\tau(c',c)=-1$; and it is $1$ or $(t_0/t_1)^{\pm1}$ if $\tau(c',c)=0$. So its total degree in $t_0$ and $t_1$ is $\tau(c',c)$. Hence $R_\gamma$ has the total degree $\sum_k\tau(c'_k,c_k)$, and by \eqref{eq.turning} it is $(t_0t_1)^{-1}$ if $\hate$ lies on the right of $\gamma$, and $t_0t_1$ if it lies on its left. This is (ii).

\begin{figure}[ht]
\centering
\begin{tikzpicture}[>=Stealth,scale=0.95]
\def\crossing{%
   \draw[->] (0.7,-0.7) -- (-0.7,0.7); \draw[->,preaction={draw=white,line width=4pt}] (-0.7,-0.7) -- (0.7,0.7);}
\begin{scope}
\crossing
\draw[red,line width=1.1pt,rounded corners=3pt,postaction={decorate,decoration={markings,mark=at position 0.3 with {\arrow{Stealth[length=6pt]}},mark=at position 0.85 with {\arrow{Stealth[length=6pt]}}}}] (0,-1.1) -- (0,0) -- (1.1,0);
\node[font=\scriptsize,red,left] at (0,-0.95) {$c'$}; \node[font=\scriptsize,red,above] at (0.95,0) {$c$};
\node[font=\small] at (0,-1.7) {$\tau(c',c)=1$};
\node[font=\small] at (0,-2.3) {three ends on the left};
\end{scope}
\begin{scope}[xshift=5cm]
\crossing
\draw[red,line width=1.1pt,postaction={decorate,decoration={markings,mark=at position 0.25 with {\arrow{Stealth[length=6pt]}},mark=at position 0.85 with {\arrow{Stealth[length=6pt]}}}}] (0,-1.1) -- (0,1.1);
\node[font=\scriptsize,red,left] at (0,-0.95) {$c'$}; \node[font=\scriptsize,red,left] at (0,0.95) {$c$};
\node[font=\small] at (0,-1.7) {$\tau(c',c)=0$};
\node[font=\small] at (0,-2.3) {two ends on the left};
\end{scope}
\begin{scope}[xshift=10cm]
\crossing
\draw[red,line width=1.1pt,rounded corners=3pt,postaction={decorate,decoration={markings,mark=at position 0.3 with {\arrow{Stealth[length=6pt]}},mark=at position 0.85 with {\arrow{Stealth[length=6pt]}}}}] (0,-1.1) -- (0,0) -- (-1.1,0);
\node[font=\scriptsize,red,right] at (0,-0.95) {$c'$}; \node[font=\scriptsize,red,above] at (-0.95,0) {$c$};
\node[font=\small] at (0,-1.7) {$\tau(c',c)=-1$};
\node[font=\small] at (0,-2.3) {one end on the left};
\end{scope}
\end{tikzpicture}
\caption{A curve $\gamma$ (red) which enters a crossing through the corner $c'$ and leaves it through the corner $c$, drawn for a positive crossing and for $c'={\downarrow}$: it has $2+\tau(c',c)$ of the four ends of edges of the crossing on its left. In the first picture the end between the two corners is the incoming under-edge, and $w_j(c,c')/w_j(c',c)=t_1$; in the third one it is the incoming over-edge, and the quotient is $t_1^{-1}$; in the second one the two corners are opposite, and the quotient is $1$.}
\label{fig.turningnumber}
\end{figure}
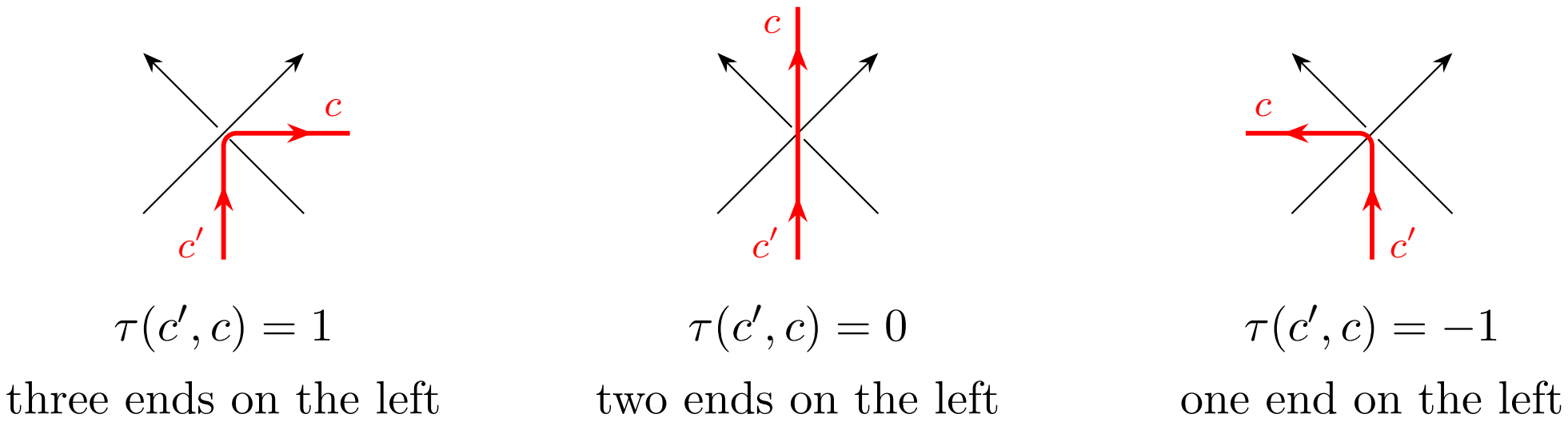

\emph{The weights.} Let $t_0=sq$ and $t_1=s^{-1}q$. If $x(j)=\{c,c'\}$, then
\begin{equation}\label{eq.geometricmean}
   \ell(x,j)^2=w_j(c,c')\,w_j(c',c),
\end{equation}
and the sign of $\ell(x,j)$ is $(-1)^d$, where $d$ is the number of corners $\downarrow$ among $c$ and $c'$: up to this sign, the weight \eqref{eq.ellweight} is the geometric mean of the weights of the two orders of the two markers. This is read off the following table, in which $\ve_j=1$; for $\ve_j=-1$ all the entries are replaced by their inverses.
\[
\begin{array}{c|c|c|c}
   x(j)=\{c,c'\} & \ell(x,j) & \text{in the variables }t_0,t_1 & w_j(c,c'),\ w_j(c',c)\\ \hline
   \{\uparrow,\uparrow\} & sq & t_0 & t_0\\
   \{\downarrow,\downarrow\} & s^{-1}q & t_1 & t_1\\
   \{A,A\},\ \{B,B\} & 1 & 1 & 1\\
   \{\uparrow,A\} & s^{1/2}q^{1/2} & t_0^{1/2} & 1,\ t_0\\
   \{\uparrow,B\} & s^{1/2}q^{1/2} & t_0^{1/2} & t_0,\ 1\\
   \{\uparrow,\downarrow\} & -1 & -1 & 1,\ 1\\
   \{A,B\} & q & (t_0t_1)^{1/2} & t_0,\ t_1\\
   \{\downarrow,A\} & -s^{-1/2}q^{1/2} & -t_1^{1/2} & t_1,\ 1\\
   \{\downarrow,B\} & -s^{-1/2}q^{1/2} & -t_1^{1/2} & 1,\ t_1
\end{array}
\]
The signs of the weights multiply to $(-1)^{b(x)}$, the sign of the pairs. A crossing at which $x$ chooses twice the corner $c$ has the weight $\ell(x,j)=w_j(c,c)$, its factor in the monomial of every pair with superposition $x$. For a cycle $\gamma$ of $\Gamma_x$, let $W_\gamma$ and $t_0t_1W_\gamma=q^2W_\gamma$ be the products of the factors of the crossings of $\gamma$ in the monomials of the pairs, for the two ways of distributing the corners of $\gamma$ between the two states, as in (ii). By \eqref{eq.geometricmean} the product of the numbers $|\ell(x,j)|$ over the crossings of $\gamma$ is $(W_\gamma\cdot q^2W_\gamma)^{1/2}=q\,W_\gamma$, and
\[
   [2]\cdot q\,W_\gamma=W_\gamma+q^2\,W_\gamma .
\]
Taking the product over the cycles of $\Gamma_x$ and over the crossings at which $x$ chooses twice the same corner, we find that $[2]^{k(x)}\prod_j\ell(x,j)$ is the sum of the terms of the $2^{k(x)}$ pairs with superposition $x$. The sum over $x$ is the sum over all the ordered pairs of Kauffman states, which is $\LG(L;t_0,t_1)$ by Proposition~\ref{prop.statesum}~(iii), and it is $Q_2(D,\hate)$ by \eqref{eq.Q2}.
\end{proof}

The model of Corollary~\ref{cor.ballinger} is Ballinger's. He has presented the invariants $Q_N$ in several talks, among them \cite{BallingerTalk}, and we learned of them from these presentations and from conversations with him. It was not known that $Q_N$ is the Links--Gould invariant for any $N\ge2$: this is conjectured in \cite{BallingerTalk}, and Corollary~\ref{cor.ballinger} proves it for $N=2$, by showing that the state sum which the determinant formula of \cite{HK} gives for $\LG$ is Ballinger's sum. With his model, Ballinger plans to prove three properties of the invariants $Q_N$: the genus bound for every Seifert surface of the link, the multiplicativity of the coefficient of extreme degree under arbitrary Murasugi sums, and the monicity on fibered links. For the Links--Gould invariant these are the theorems of \cite{Kohli-Tahar, LNVdV25}, of \cite{ALvdV26} and of \cite{LNVdV25} which were recalled in the introduction and in Section~\ref{sec.degree}. From the determinant formula, Section~\ref{sec.degree} obtains the first two for the canonical surface of a diagram and the third for alternating links and for the closures of homogeneous braids, and Section~\ref{ssec.beyond} shows why the deformed Seifert matrices do not give more. So the approach of this monograph does ultimately reach these three properties in full generality, modulo Ballinger's work: the determinant formula gives the state sum of Proposition~\ref{prop.statesum}, this state sum is his model for $N=2$, and the three properties then follow from what Ballinger plans to prove for his model.

\begin{remark}\label{rem.ballinger}
(i) \emph{The variables.} In the variables $t_0=p^{-1}q$, $t_1=pq$ of Lemma~\ref{lem.pq}, the variables of Ballinger's model are $s=p^{-1}$ and $q$: the exponent of $s$ is the degree \eqref{eq.KTdegree}, to which the genus bound applies, and the exponent of $q$ is the total degree in $t_0$ and $t_1$. In the variables $t_0,t_1$ the term of a state $x$ with multiplicity two is
\[
   [2]^{k(x)}\prod_{j=1}^c\ell(x,j)=(-1)^{b(x)}\,W(x)\,(1+t_0t_1)^{k(x)},
\]
where $W(x)$ is the monomial of lowest total degree among the monomials of the pairs with superposition $x$.

(ii) \emph{The points $t_0t_1=\pm1$.} At $t_0t_1=-1$ the terms of the states with $k(x)\ge1$ vanish, and the states with $k(x)=0$ are the superpositions of a Kauffman state with itself, whose terms are the diagonal of Remark~\ref{rem.statesum}~(iii). So the pairs of two distinct states cancel at $t_0t_1=-1$, and they do so within each superposition. This is a proof, in the state model, of the cancellation of Corollary~\ref{cor.diagonal}, hence another proof of the specialization $\LG(L;t_0,-t_0^{-1})=\Delta_L(t_0^2)$ of \cite{DWIL}. At $t_0t_1=1$ the factor $(1+t_0t_1)^{k(x)}$ is the number of pairs with superposition $x$, in agreement with Remark~\ref{rem.statesum}~(ii).

(iii) \emph{Pairs of states and colorings.} The $2^{k(x)}$ pairs with superposition $x$ are the colorings of the cycles of $\Gamma_x$ by the two states, and $[2]^{k(x)}=(q+q^{-1})^{k(x)}$ is a sum over these colorings: by Corollary~\ref{cor.ballinger}~(ii), the monomial of a pair is $\prod_j|\ell(x,j)|$ times $q^{-1}$ for every cycle which has $\hate$ on its right and $q$ for every cycle which has $\hate$ on its left, the cycles being oriented as in that statement. So the state sum of Proposition~\ref{prop.statesum}~(iii) is the sum \eqref{eq.Q2} in which the evaluation of the web is expanded over its colorings. The weight of a pair is a product over the crossings which is not symmetric in the two states; the weight $\prod_j\ell(x,j)$ is symmetric, and the difference is the contribution $q^{\pm1}$ of each cycle, which Corollary~\ref{cor.ballinger}~(ii) distributes over its crossings.
\end{remark}

\subsection{Alternating diagrams: Kasteleyn's rule}\label{ssec.kasteleyn}
In the rest of this section, except in Section~\ref{ssec.gluing}, the diagrams are alternating: every edge is an edge of the under-strand at exactly one of its two ends, and an edge of the over-strand at the other one.

\emph{A guide to the rest of the section.} The five subsections which follow use more graph theory and more lattice theory than the rest of this monograph. We define the vocabulary as it comes, in Definitions~\ref{def.graphs}, \ref{def.lattices}, \ref{def.flowcut}, \ref{def.directed} and~\ref{def.spinc}, and each subsection begins with a description of its content. In short:
\begin{itemize}
\item Section~\ref{ssec.kasteleyn}: on an alternating diagram all the terms of the state sum have the same sign, once the signs of the variables are changed, so that nothing cancels (Theorem~\ref{thm.kauffmanalt}).
\item Section~\ref{ssec.statestrees}: the states are the trees of \cite{HK}, and the signature of the link can be read off every single state (Corollary~\ref{cor.statestrees} and Proposition~\ref{prop.signature}).
\item Section~\ref{ssec.tait}: the states are the spanning trees of a plane graph whose edges are the crossings, the Tait graph; the two flavors of a state become two ways of orienting the edges of this graph, and the monomial of a pair of states is computed from these orientations, edge by edge (Theorem~\ref{thm.tait}).
\item Section~\ref{ssec.gluing}: a formula for the polynomial of a link which is cut into two pieces by a sphere meeting it in four points (Theorem~\ref{thm.gluing}). It holds for all links, and it can be read independently of the rest.
\item Section~\ref{ssec.lattice}: how much of the state sum only depends on the lattice of the integral flows and of the integral cuts of the Tait graph, which is an invariant of the link up to mutation; this leads to a proof, by the Links--Gould polynomial, of two conjectures of Tait (Theorem~\ref{thm.extreme} and Corollary~\ref{cor.crossings}).
\end{itemize}

\emph{This subsection.} For a general diagram the signs $(-1)^{b(\mathfrak s)}$ of Proposition~\ref{prop.statesum} depend on the state, and the terms of the state sum cancel each other in part. We show that for an alternating diagram the sign of a state is determined, up to a global sign, by the number of crossings in $X(\mathfrak s)$, and also by the number of crossings in $Y(\mathfrak s)$ (Lemma~\ref{lem.kasteleyn}). Consequently, when $t_0$ and $t_1$ are replaced by $-t_0$ and $-t_1$, all the pairs of states contribute to $\LG$ with the same sign (Theorem~\ref{thm.kauffmanalt}). This is a case of a classical sign rule for determinants, due to Kasteleyn, which we first recall together with some vocabulary.

\begin{definition}[graphs, trees and matchings]\label{def.graphs}
A \emph{graph} has finitely many vertices and edges. Every edge joins two vertices, its \emph{ends}; the two ends may coincide, in which case the edge is a \emph{loop}, and several edges may join the same two vertices. A \emph{path} is a sequence of edges in which each edge has an end in common with the next one, and a \emph{cycle} is a closed path which passes only once through each of its vertices. The graph is \emph{connected} if any two vertices are joined by a path.
\begin{itemize}
\item A \emph{plane graph} is a graph drawn on the sphere $S^2$, in such a way that its edges meet at their ends only; its \emph{faces} are the components of the complement. The \emph{dual} of a connected plane graph $\CG$ is the plane graph $\CG^*$ which has one vertex in each face of $\CG$ and one edge $e^*$ for each edge $e$ of $\CG$: the edge $e^*$ crosses $e$ once, meets no other edge of $\CG$, and joins the faces on the two sides of $e$.
\item A \emph{spanning tree} of a connected graph is a set $T$ of edges which contains no cycle and which connects all the vertices; it has one edge less than the graph has vertices. If one vertex $r$, the \emph{root}, is chosen, every other vertex is joined to $r$ by exactly one path in $T$, and each edge of $T$ has one end which is \emph{nearer to} the root and one end which is \emph{farther from} it. If $T$ is a spanning tree of a connected plane graph $\CG$, the edges $e^*$ with $e\notin T$ form a spanning tree of $\CG^*$, the \emph{dual tree} of $T$.
\item A graph is \emph{bipartite} if its vertices are colored with two colors in such a way that every edge joins two vertices of different colors. A \emph{perfect matching} of a graph is a set of edges such that every vertex is an end of exactly one of them.
\end{itemize}
\end{definition}

Here is the example which matters in this subsection. Let $\Gamma$ be the graph of Section~\ref{ssec.ballinger}: its vertices are the crossings of $D$ and one point $v_f$ in each face $f$, and its edges are the corners of the crossings, the edge of a corner joining its crossing to the point of its face (Figure~\ref{fig.gammagraph}, left). It is a plane bipartite graph, the two colors being ``crossing'' and ``face''. By Definition~\ref{def.state}, a Kauffman state of $(D,\hate)$ is exactly a perfect matching of the graph obtained from $\Gamma$ by removing the points of the two faces adjacent to $\hate$: every crossing is matched with one of its four corners, and every remaining face with one of its corners.

\emph{Kasteleyn's rule.} Consider a bipartite graph with $c$ vertices of each color, which we call the rows and the columns, and suppose that every edge $e$ carries a weight $w_e$ and a sign $\varkappa_e=\pm1$. Let $K$ be the $c\times c$ matrix whose entry in the row $r$ and in the column $s$ is the sum of the numbers $\varkappa_e\,w_e$ over the edges $e$ which join $r$ to $s$. A perfect matching $M$ matches every row with a column, and expanding the determinant gives
\[
   \det K=\sum_{M}\ \sgn(M)\,\prod_{e\in M}\varkappa_e\,w_e ,
\]
where $M$ runs over the perfect matchings and $\sgn(M)$ is the sign of the bijection from the rows to the columns which $M$ defines. In general the terms of this sum have different signs and cancel in part. Kasteleyn's theorem \cite{Kasteleyn63} says that for a \emph{plane} bipartite graph the signs $\varkappa_e$ can be chosen in such a way that $\sgn(M)\prod_{e\in M}\varkappa_e$ is the same for all the perfect matchings, so that $\pm\det K$ is the sum of the weights $\prod_{e\in M}w_e$ of the matchings, without any cancellation; and that for this it is enough that, around every face which is bounded by $2\ell$ edges, the product of the signs be $(-1)^{\ell+1}$. When all the faces are quadrilaterals, as is the case for $\Gamma$, the condition is that the product of the four signs around every face be $-1$.

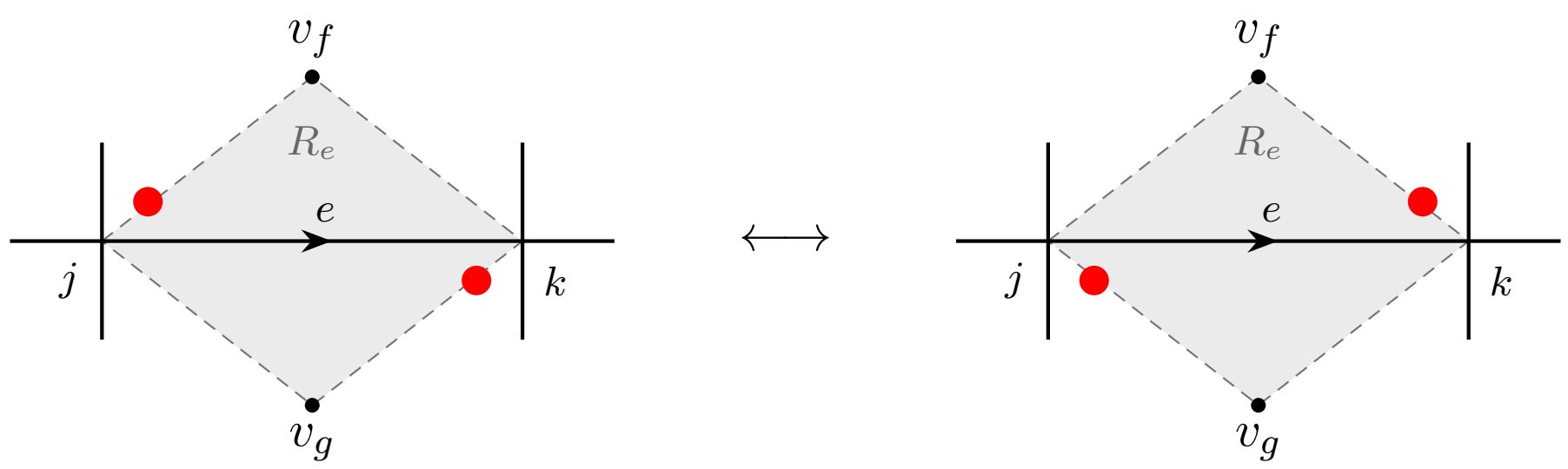
\begin{figure}[ht]
\centering
\begin{tikzpicture}[>=Stealth,scale=1.0,every node/.style={font=\small}]
\def\clockbase{%
  \fill[black!8] (-1.6,0) -- (0,1.25) -- (1.6,0) -- (0,-1.25) -- cycle;
  \draw[black!55,densely dashed] (-1.6,0) -- (0,1.25) -- (1.6,0) -- (0,-1.25) -- cycle;
  \draw[thick,->] (-1.6,0) -- (0.15,0); \draw[thick] (0,0) -- (1.6,0);
  \draw[thick] (-1.6,-0.75) -- (-1.6,0.75); \draw[thick] (1.6,-0.75) -- (1.6,0.75);
  \draw[thick] (-2.3,0) -- (-1.6,0); \draw[thick] (1.6,0) -- (2.3,0);
  \fill (0,1.25) circle (1.6pt); \fill (0,-1.25) circle (1.6pt);
  \node[above] at (0,1.25) {$v_f$}; \node[below] at (0,-1.25) {$v_g$};
  \node[font=\scriptsize] at (0.1,0.22) {$e$};
  \node[font=\scriptsize] at (-1.85,-0.3) {$j$}; \node[font=\scriptsize] at (1.85,-0.3) {$k$};
  \node[font=\scriptsize,black!60] at (0.0,0.75) {$R_e$};}
\begin{scope}
  \clockbase
  \fill[red] (-1.25,0.3) circle (3.2pt); \fill[red] (1.25,-0.3) circle (3.2pt);
\end{scope}
\node at (3.6,0) {$\longleftrightarrow$};
\begin{scope}[xshift=7.2cm]
  \clockbase
  \fill[red] (-1.25,-0.3) circle (3.2pt); \fill[red] (1.25,0.3) circle (3.2pt);
\end{scope}
\end{tikzpicture}
\caption{A clock move. The edge $e$ of $D$ joins the crossings $j$ and $k$ and separates the faces $f$ and $g$. The four corners adjacent to $e$, two at each of its ends, are four edges of the graph $\Gamma$ (dashed): they bound the quadrilateral $R_e$, which is a face of $\Gamma$. The two states have the same markers (dots) at the other crossings; at $j$ and $k$ their markers are two opposite sides of the quadrilateral in the first state, and the other two sides in the second one.}
\label{fig.clockmove}
\end{figure}

By the proof of Proposition~\ref{prop.statesum}~(ii), the face matrices $Q_{D,\hate}(x)$ and $\tQ_{D,\hate}(y)$ are matrices $K$ of this kind for the graph $\Gamma$ without two vertices: the rows are the crossings, the columns are the faces, and the entry of a corner is its coefficient \eqref{eq.cornercoeff}. Lemma~\ref{lem.kasteleyn} below says that on an alternating diagram, after the change of variables $x_j\mapsto-x_j$, resp.\ $y_j\mapsto-y_j$, all the states have the same sign; in its proof the sign of this change of variables at a corner is called $\varkappa$, resp.\ $\tilde\varkappa$, and the condition \eqref{eq.kastcond} which makes the proof work is Kasteleyn's condition for the quadrilaterals of $\Gamma$. We do not use Kasteleyn's theorem itself, but the following theorem of Kauffman, which was already used in the proof of Proposition~\ref{prop.statesum}. A \emph{clock move} changes a Kauffman state at two crossings $j$ and $k$ which are joined by an edge $e$ of $D$, when the markers of $j$ and $k$ are corners adjacent to $e$ which lie in the two faces adjacent to $e$, one in each: it replaces each of the two markers by the corner of the same crossing which is adjacent to $e$ in the other face (Figure~\ref{fig.clockmove}). The \emph{Clock Theorem} \cite{KauffmanFKT} says that any two Kauffman states of $(D,\hate)$ are related by a sequence of clock moves.

\begin{lemma}\label{lem.kasteleyn}
Let $D$ be an alternating diagram and $\hate$ an edge of $D$. Then $\sgn(\mathfrak s)\,(-1)^{|X(\mathfrak s)|}$ takes the same value on all the states $\mathfrak s$ of $(D,\hate)$, and so does $\sgn(\mathfrak s)\,(-1)^{|Y(\mathfrak s)|}$. Equivalently, the number $a(\mathfrak s)=|X(\mathfrak s)\cap Y(\mathfrak s)|$ of crossings at which the marker of $\mathfrak s$ is the $A$-corner has the same parity for all the states of $(D,\hate)$.
\end{lemma}
\begin{proof}
For a corner $\diamond$ of a crossing put $\varkappa(\diamond)=-1$ if $\diamond\in\{\uparrow,A\}$, i.e.\ if $\diamond$ is adjacent to the outgoing under-edge, and $\varkappa(\diamond)=1$ otherwise; then
\begin{equation}\label{eq.etaX}
   \sgn(\mathfrak s)\,(-1)^{|X(\mathfrak s)|}=\sgn(\mathfrak s)\prod_{j=1}^c\varkappa(\mathfrak s(j))
\end{equation}
for every state $\mathfrak s$. By \eqref{eq.corneredges}, the product of $\varkappa$ over the two corners adjacent to an edge at one of its ends is $+1$ if the edge is $u_j$ or $u'_j$ there, and $-1$ if it is $o_j$ or $o'_j$. Since $D$ is alternating, every edge is an over-edge at exactly one of its ends, so that
\begin{equation}\label{eq.kastcond}
   \prod_{\diamond\text{ adjacent to }e}\varkappa(\diamond)=-1\qquad\text{for every edge $e$ of $D$,}
\end{equation}
the product running over the four corners adjacent to $e$, two at each end.

Let $\Gamma$ be the graph of Section~\ref{ssec.ballinger}, drawn on $S^2$. Its vertices are the crossings of $D$ and one point $v_f$ in every face $f$. Its edges are in bijection with the corners: the edge of a corner of the crossing $j$ which lies in the face $f$ is an arc from $j$ to $v_f$ inside $f$, leaving $j$ through this corner, and distinct arcs meet at their ends only. Such arcs exist because the faces of the connected diagram $D$ are open disks. The graph $\Gamma$ is connected, and its arcs cut every face of $D$ into sectors, one for each side of an edge of $D$ on the boundary of the face. Hence every component of $S^2\sm\Gamma$ is an open disk $R_e$ which contains exactly one edge $e$ of $D$, and $R_e$ is bounded by the closed walk formed by the edges of $\Gamma$ of the four corners adjacent to $e$. We identify the corners with the edges of $\Gamma$. A state of $(D,\hate)$ is then a perfect matching of the graph obtained from $\Gamma$ by removing the two vertices $v_{f_{\rm L}(\hate)}$ and $v_{f_{\rm R}(\hate)}$.

By the Clock Theorem of \cite{KauffmanFKT}, which was used in the proof of Proposition~\ref{prop.statesum}~(ii), any two states are related by a sequence of moves each of which exchanges the faces of the markers of two crossings joined by an edge $e$. Such a move replaces two of the four corners adjacent to $e$, one at each end of $e$, by the other two, so that the two states differ along the boundary of the disk $R_e$ (Figures~\ref{fig.gammagraph} and~\ref{fig.clockmove}). It changes $\sgn(\mathfrak s)$, because the bijection from the crossings to the faces which the state defines is composed with a transposition. It also changes the sign of $\prod_j\varkappa(\mathfrak s(j))$: by \eqref{eq.kastcond}, the product of $\varkappa$ over the two corners which are removed is the opposite of its product over the two corners which replace them. Hence \eqref{eq.etaX} takes the same value on all the states, which is the first claim.

For the second claim replace $\varkappa$ by the function $\tilde\varkappa$ which is $-1$ at the corners $\downarrow$ and $A$, i.e.\ at the corners adjacent to the incoming over-edge, and $1$ at the other two. Then $\sgn(\mathfrak s)(-1)^{|Y(\mathfrak s)|}=\sgn(\mathfrak s)\prod_j\tilde\varkappa(\mathfrak s(j))$, the product of $\tilde\varkappa$ over the two corners adjacent to an edge at one of its ends is $+1$ for $o_j$ and $o'_j$ and $-1$ for $u_j$ and $u'_j$, so that \eqref{eq.kastcond} holds for $\tilde\varkappa$, and the rest of the proof is unchanged.

Finally $|Y(\mathfrak s)|=b(\mathfrak s)+a(\mathfrak s)$ and $|X(\mathfrak s)|=|V(\mathfrak s)|-b(\mathfrak s)+a(\mathfrak s)$ by Definition~\ref{def.state}, where $(-1)^{|V(\mathfrak s)|}=(-1)^{m_{D,\hate}}$ does not depend on the state (Proposition~\ref{prop.statesum}~(i)). Since $\sgn(\mathfrak s)=\eta\,(-1)^{b(\mathfrak s)}$ by \eqref{eq.blackholes}, we get
\[
   \sgn(\mathfrak s)\,(-1)^{|Y(\mathfrak s)|}=\eta\,(-1)^{a(\mathfrak s)}\qquad\text{and}\qquad \sgn(\mathfrak s)\,(-1)^{|X(\mathfrak s)|}=\eta\,(-1)^{m_{D,\hate}}(-1)^{a(\mathfrak s)},
\]
so that each of the first two statements is equivalent to the last one.
\end{proof}

The condition \eqref{eq.kastcond} is Kasteleyn's condition for the planar bipartite graph $\Gamma$, all of whose faces are quadrilaterals \cite{Kasteleyn63}; for a reduced diagram of a knot, the graph $\Gamma$ without the two removed vertices is the graph of the dimer model of Cohen, Dasbach and Russell for the Alexander polynomial, whose dimer coverings, that is, perfect matchings, are Kauffman's states \cite{CohenDasbachRussell}. By \eqref{eq.cornercoeff} the entries of the face matrices at $x_j=-1$ and at $y_j=-1$ have the signs $\ve_j(-1)^{\ind{\diamond\text{ vertical}}}\varkappa(\diamond)$ and $\ve_j(-1)^{\ind{\diamond\text{ vertical}}}\tilde\varkappa(\diamond)$, and the product of the first two factors over the four corners adjacent to an edge is $1$. So an alternating diagram is one for which Alexander's matrix at $t=-1$ is a Kasteleyn matrix of $\Gamma$, and this holds for the face matrices of both flavors, the crossing variables playing the role of the weights of the edges of $\Gamma$. In the proof above the lemma follows from \eqref{eq.kastcond} and the Clock Theorem of \cite{KauffmanFKT}, as in the proof of Proposition~\ref{prop.statesum}~(ii): two states related by a clock move differ along the boundary of one disk $R_e$.

\begin{theorem}[positivity of the state sum]\label{thm.kauffmanalt}
Let $D$ be an alternating diagram of a link $L$ with $\mu$ components, and let $\hate$ be any edge of $D$.
\begin{enumerate}\renewcommand{\theenumi}{\roman{enumi}}
\item Let $\sigma_{D,\hate}=(-1)^{a(\mathfrak s)}$, where $a(\mathfrak s)$ is the number of crossings at which the marker of a state $\mathfrak s$ of $(D,\hate)$ is the $A$-corner; this sign does not depend on the state (Lemma~\ref{lem.kasteleyn}). Then
\[
   P_{D,\hate}(x)=\sigma_{D,\hate}\sum_{\mathfrak s}(-1)^{|X(\mathfrak s)|}\,x^{X(\mathfrak s)}\qquad\text{and}\qquad \tP_{D,\hate}(y)=\sigma_{D,\hate}\sum_{\mathfrak s}(-1)^{|Y(\mathfrak s)|}\,y^{Y(\mathfrak s)},
\]
the sums running over the Kauffman states of $(D,\hate)$. The number of these states is $|\Delta_L(-1)|$, the determinant of $L$.
\item The Links--Gould polynomial is the sum over the ordered pairs of states
\[
   (-1)^{\mu-1}\,\LG(L;-t_0,-t_1)=\sum_{(\mathfrak s,\mathfrak s')}\ \prod_{j\in X(\mathfrak s)\sm Y(\mathfrak s')}t_0^{\ve_j}\ \prod_{j\in Y(\mathfrak s')\sm X(\mathfrak s)}t_1^{\ve_j},
\]
a sum of $\Delta_L(-1)^2$ monomials. In particular the coefficient of $t_0^at_1^b$ in $\LG(L)$ is zero or has the sign $(-1)^{\mu-1+a+b}$.
\end{enumerate}
\end{theorem}
\begin{proof}
(i) By Definition~\ref{def.state}, $|Y(\mathfrak s)|=b(\mathfrak s)+a(\mathfrak s)$ and $|X(\mathfrak s)|=|V(\mathfrak s)|-b(\mathfrak s)+a(\mathfrak s)$, and $|V(\mathfrak s)|\equiv\mu-1$ modulo $2$; so the two identities are those of Proposition~\ref{prop.statesum}~(ii), in which $(-1)^{b(\mathfrak s)}=\sigma_{D,\hate}(-1)^{|Y(\mathfrak s)|}=(-1)^{\mu-1}\sigma_{D,\hate}(-1)^{|X(\mathfrak s)|}$ for every state, by Lemma~\ref{lem.kasteleyn}. Let $n$ be the number of states. The point $x_j=y_j=-1$ for all $j$ is the Alexander point at $t=-1$, whatever the signs of the crossings, and there the two identities give $P_{D,\hate}=\tP_{D,\hate}=\sigma_{D,\hate}\,n$, in agreement with Lemma~\ref{lem.level}~(i), by which the two edge polynomials have the same value at the Alexander point: $P_{D,\hate}(t^{\ve})=\det A_{D,\hate}(t)=\det B_{D,\hate}(t)=\tP_{D,\hate}(t^{\ve})$. Hence $n=|P_{D,\hate}((-1)^{\ve})|=|\Delta_L(-1)|$ by Lemma~\ref{lem.alexpoly}.

(ii) By (i), $P_{D,\hate}(x)\,\tP_{D,\hate}(y)=\sum_{(\mathfrak s,\mathfrak s')}(-1)^{|X(\mathfrak s)|+|Y(\mathfrak s')|}\,x^{X(\mathfrak s)}y^{Y(\mathfrak s')}$, so by Theorem~\ref{thm.HK} and \eqref{eq.ED}
\[
   \LG(L;t_0,t_1)=(-1)^{\mu-1}\sum_{(\mathfrak s,\mathfrak s')}(-1)^{|X(\mathfrak s)|+|Y(\mathfrak s')|}\prod_{j\in X(\mathfrak s)\sm Y(\mathfrak s')}t_0^{\ve_j}\ \prod_{j\in Y(\mathfrak s')\sm X(\mathfrak s)}t_1^{\ve_j}.
\]
Now $|X(\mathfrak s)|+|Y(\mathfrak s')|\equiv|X(\mathfrak s)\sm Y(\mathfrak s')|+|Y(\mathfrak s')\sm X(\mathfrak s)|$ modulo $2$ and $(-t)^{\ve_j}=-t^{\ve_j}$, so replacing $t_0,t_1$ by $-t_0,-t_1$ removes all the signs. The last statement follows, since $(-t_0)^a(-t_1)^b=(-1)^{a+b}t_0^at_1^b$.
\end{proof}

Part (ii) is the theorem of \cite{HK} that the Links--Gould polynomial of an alternating link is alternating, and since a connected diagram has at least one Kauffman state \cite{KauffmanFKT}, the sum is not empty. At the Alexander point, (i) is the classical fact that the states of an alternating diagram all contribute to Kauffman's state sum at $t=-1$ with the same sign; what the argument above adds is that the sign rule holds with the crossing variables, for the two flavors and with the same global sign, which is what the coupling $\BE_D$ requires. The proof of \cite{HK} follows Crowell \cite{Crowell59} instead, and the two expansions have the same terms.

\subsection{States, trees and the signature}\label{ssec.statestrees}
\emph{This subsection.} We compare the Kauffman states of an alternating diagram with the trees of \cite{HK}, which were used in Section~\ref{ssec.alternating} (Corollary~\ref{cor.statestrees}), and we show that every state knows the signature of the link: the crossings at which its marker lies between the two outgoing or between the two incoming edges, counted with their signs, add up to $-\sigma_L$ (Proposition~\ref{prop.signature}). For a fibered alternating link this gives the top coefficient of $\LG$ and the Hopf invariant of the fiber in terms of the genus, the number of components and the signature (Corollary~\ref{cor.extremalstates}). The following objects were introduced earlier.
\begin{itemize}
\item \emph{The sets of a state} (Definition~\ref{def.state} and Figure~\ref{fig.corners}): $X(\mathfrak s)$ is the set of the crossings at which the marker of $\mathfrak s$ is adjacent to the outgoing under-edge, $Y(\mathfrak s)$ the set of those at which it is adjacent to the incoming over-edge, and $V(\mathfrak s)$ the set of those at which the marker is vertical, that is, lies between the two outgoing or between the two incoming edges. For a set $S$ of crossings, $\ve(S)$ is the sum of the signs of its elements.
\item \emph{The trees of \cite{HK}} (Section~\ref{ssec.alternating}): for an edge $\hate=o'_{j_1}=u_{k_0}$, an under-tree $\tau$ chooses at every crossing one of the two edges of the under-strand, an over-tree $\tau'$ one of the two edges of the over-strand, and $S(\tau)$ and $S'(\tau')$ are the sets of the crossings at which the chosen edge is $u_j$, resp.\ $o'_j$; the edge polynomials are sums over these trees, by \eqref{eq.HKtrees}.
\item \emph{The two colors} (Lemma~\ref{lem.altcoloring} and Figure~\ref{fig.altcolor}): the faces of an alternating diagram are colored black and white, two faces which are adjacent to the same edge having different colors, and every corner has the color of its face. The black corners of a positive crossing are its two vertical corners, and those of a negative crossing are its $A$-corner and its $B$-corner. We write $F_{\rm b}$ for the number of black faces.
\item \emph{The signature, the top coefficient and the Hopf invariant.} $\sigma_L$ is the signature of $L$, with the convention of Section~\ref{ssec.signature}. In the variables $t_0=q/p$, $t_1=pq$ of Lemma~\ref{lem.pq}, $\mathrm{top}(L,q)$ is the coefficient of $p^{2g(L)+\mu-1}$ in $\LG(L)$, a Laurent polynomial in $q$. If $L$ is fibered, with fiber $F$, then $\mathrm{top}(L,q)$ is a monomial $q^{\nu-2\lambda(\xi_F)}$, where $\nu=2g(L)+\mu-1$ and the integer $\lambda(\xi_F)$ is the Hopf invariant of the plane field defined by the fibration \cite{LNVdV25}; we only use it through this formula.
\end{itemize}

\begin{corollary}[Kauffman states and spanning trees]\label{cor.statestrees}
Let $D$ be an alternating diagram and let $\hate=o'_{j_1}=u_{k_0}$ be an edge which leaves a crossing as its over-strand, as in Section~\ref{ssec.alternating}. Then $\mathfrak s\mapsto X(\mathfrak s)$ is a bijection from the Kauffman states of $(D,\hate)$ onto the sets $S(\tau)$ of the under-trees, and $\mathfrak s\mapsto Y(\mathfrak s)$ is a bijection from the states onto the sets $S'(\tau')$ of the over-trees. Under these two bijections the sum of Theorem~\ref{thm.kauffmanalt}~(ii) over the pairs of states is, term by term, the sum over the pairs of an under-tree and an over-tree given by \eqref{eq.HKtrees} and Theorem~\ref{thm.HK}, and the sign $\sigma$ of \eqref{eq.HKtrees} is the sign $\sigma_{D,\hate}$ of Theorem~\ref{thm.kauffmanalt}~(i).
\end{corollary}
\begin{proof}
By Theorem~\ref{thm.kauffmanalt}~(i) and \eqref{eq.HKtrees}, $\sum_{\mathfrak s}(-1)^{|X(\mathfrak s)|}x^{X(\mathfrak s)}=\pm\sum_\tau(-1)^{|S(\tau)|}x^{S(\tau)}$. Comparing the coefficients of $x^S$, the number of states with $X(\mathfrak s)=S$ is $\pm$ the number of under-trees with $S(\tau)=S$, which is $0$ or $1$ since $\tau\mapsto S(\tau)$ is injective; so the sign is $+$ and the two numbers are equal; hence $\sigma=\sigma_{D,\hate}$, and $\mathfrak s\mapsto X(\mathfrak s)$ is a bijection onto the sets $S(\tau)$. For the second flavor, Theorem~\ref{thm.kauffmanalt}~(i) and \eqref{eq.HKtrees} give $\sum_{\mathfrak s}(-1)^{|Y(\mathfrak s)|}y^{Y(\mathfrak s)}=\sum_{\tau'}(-1)^{|S'(\tau')|}y^{S'(\tau')}$, with the same sign $\sigma=\sigma_{D,\hate}$ on both sides, and $\tau'\mapsto S'(\tau')$ is injective; so the same comparison of the coefficients shows that $\mathfrak s\mapsto Y(\mathfrak s)$ is a bijection onto the sets $S'(\tau')$. Finally, in both sums the monomial of a pair only depends on the two sets of crossings, $X(\mathfrak s)$ and $Y(\mathfrak s')$ in one sum, $S(\tau)$ and $S'(\tau')$ in the other one; so the two sums agree term by term.
\end{proof}

A Kauffman state is a spanning tree of the checkerboard graph of the diagram together with the complementary tree of the dual graph \cite{KauffmanFKT} (Section~\ref{ssec.tait} gives the sets $X(\mathfrak s)$ and $Y(\mathfrak s)$ in these terms), and the trees of \cite{HK} are those of Crowell's proof \cite{Crowell59}; Corollary~\ref{cor.statestrees} matches the two kinds of states through the sets of crossings, and not only through their weights $t^{\ve(X(\mathfrak s))}=t^{\ve(S(\tau))}$, but it does not describe the bijection. For special alternating links a weight-preserving bijection between the states of Crowell's model and those of Kauffman's model was constructed by Hafner, M\'esz\'aros and Vidinas by means of Eulerian tours \cite[Theorem~6.7]{HMV24}; we have not compared it with the present one.

In Kauffman's model an under-tree and an over-tree which correspond to the same state differ exactly at the crossings where the marker of the state is vertical, and on an alternating diagram the signed number of these crossings does not depend on the state: it is the signature. Recall that $\sigma_L$ denotes the signature of $L$, with the conventions of Section~\ref{ssec.signature}, for which the right-handed trefoil has $\sigma_L=-2$; it should not be confused with the sign $\sigma_{D,\hate}$ of Theorem~\ref{thm.kauffmanalt}.

\begin{proposition}[the vertical corners count the signature]\label{prop.signature}
Let $D$ be an alternating diagram of $L$ with $c_-$ negative crossings, colored as in Lemma~\ref{lem.altcoloring}, with $F_{\rm b}$ black faces, and let $\hate$ be an edge of $D$. Then, for every Kauffman state $\mathfrak s$ of $(D,\hate)$,
\[
   \ve(V(\mathfrak s))=F_{\rm b}-1-c_-=-\sigma_L .
\]
Consequently $\ve\bigl(X(\mathfrak s)\cap Y(\mathfrak s)\bigr)=\frac12\bigl(m_{D,\hate}+\sigma_L\bigr)$ for every state, and the sign of Theorem~\ref{thm.kauffmanalt}~\textup{(i)} is
\[
   \sigma_{D,\hate}=(-1)^{(m_{D,\hate}+\sigma_L)/2}.
\]
\end{proposition}
\begin{proof}
By Lemma~\ref{lem.altcoloring}, the corners of a black face are black: they are vertical corners of positive crossings or horizontal corners of negative crossings; and the corners of a white face are white. A state $\mathfrak s$ has one corner in every face other than $f_{\rm L}(\hate)$ and $f_{\rm R}(\hate)$, one of which is black; so it has $F_{\rm b}-1$ black corners. On the other hand its black corners are its vertical corners at the positive crossings and its horizontal corners at the negative ones, so their number is $|V_+|+(c_--|V_-|)=\ve(V(\mathfrak s))+c_-$, where $V_\pm$ is the set of the positive, resp.\ negative, crossings in $V(\mathfrak s)$. This proves the first equality. The second one is Corollary~\ref{cor.altdefinite}~(ii), that is, Traczyk's formula for the signature of an alternating link \cite{Traczyk04}.

For the last statement, $\ve(X(\mathfrak s))+\ve(Y(\mathfrak s))=\ve(V(\mathfrak s))+2\,\ve\bigl(X(\mathfrak s)\cap Y(\mathfrak s)\bigr)$ and the left-hand side is $m_{D,\hate}$ (Proposition~\ref{prop.statesum}~(i)); and $a(\mathfrak s)=|X(\mathfrak s)\cap Y(\mathfrak s)|\equiv\ve\bigl(X(\mathfrak s)\cap Y(\mathfrak s)\bigr)$ modulo $2$.
\end{proof}

\begin{corollary}[the extremal states]\label{cor.extremalstates}
Let $D$ be an alternating diagram of $L$, let $\hate$ be any edge of $D$, $\nu=c-s+1$, and let $N$ be the absolute value of the leading coefficient of $\Delta_L$.
\begin{enumerate}\renewcommand{\theenumi}{\roman{enumi}}
\item The states $\mathfrak s$ of $(D,\hate)$ for which $\ve(X(\mathfrak s))$ is smallest are the states for which $\ve(Y(\mathfrak s))$ is largest. There are exactly $N$ of them, and
\[
   \mathrm{top}(L,q)=\sum_{(\mathfrak s,\mathfrak s')}q^{\,\ve(X(\mathfrak s)\triangle Y(\mathfrak s'))},
\]
the sum running over the $N^2$ ordered pairs of such states; each of the $N$ pairs $(\mathfrak s,\mathfrak s)$ contributes $q^{-\sigma_L}$.
\item $L$ is fibered if and only if there is a unique such state $\mathfrak s_0$, and then, with $F$ the fiber,
\[
   \mathrm{top}(L,q)=q^{-\sigma_L}\qquad\text{and}\qquad \lambda(\xi_F)=\tfrac12\bigl(\nu+\sigma_L\bigr)=g(L)+\tfrac12\bigl(\mu-1+\sigma_L\bigr).
\]
If $\hate=o'_{j_1}=u_{k_0}$, the two extremal trees of Theorem~\ref{thm.alternating} are the two flavors of this state, $S(\tau_0)=X(\mathfrak s_0)$ and $S'(\tau'_0)=Y(\mathfrak s_0)$, so that $\ve(S(\tau_0)\triangle S'(\tau'_0))=-\sigma_L$.
\end{enumerate}
\end{corollary}
\begin{proof}
As in the proof of Theorem~\ref{thm.alternating}, $D$ is tight, so that Theorem~\ref{thm.tight} applies to $(D,\hate)$, and $\nu=2g(L)+\mu-1$. By Theorem~\ref{thm.kauffmanalt}~(i), the coefficients of the edge polynomials in the monomial basis \eqref{eq.monomials} are
\[
   P_S=\sigma_{D,\hate}\,(-1)^{|S|}\,\bigl|\{\mathfrak s:X(\mathfrak s)=S\}\bigr|,\qquad \tP_T=\sigma_{D,\hate}\,(-1)^{|T|}\,\bigl|\{\mathfrak s:Y(\mathfrak s)=T\}\bigr| .
\]
In particular the $\ve$-degrees of the monomials of $P_{D,\hate}$ are the numbers $\ve(X(\mathfrak s))$, whose minimum is therefore $\delta-\nu$, and those of $\tP_{D,\hate}$ are the numbers $\ve(Y(\mathfrak s))$, whose maximum is $\delta$ (Theorem~\ref{thm.tight}~(i)); by Proposition~\ref{prop.statesum}~(i) the minimum and the maximum are attained at the same states. The coefficient of $t^{\delta-\nu}$ in $P_{D,\hate}(t^{\ve})$ is $\sigma_{D,\hate}(-1)^{\delta-\nu}$ times the number of these states; it is an extreme coefficient of $\Delta_L$ up to sign (Theorem~\ref{thm.tight}~(i)), so this number is $N$. In the formula of Theorem~\ref{thm.tight}~(iii) each pair $(\mathfrak s,\mathfrak s')$ of these states contributes $(-1)^{\mu-1}(-1)^{|X(\mathfrak s)|+|Y(\mathfrak s')|}\,q^{\ve(X(\mathfrak s)\triangle Y(\mathfrak s'))}$, and $(-1)^{|X(\mathfrak s)|+|Y(\mathfrak s')|}=(-1)^{(\delta-\nu)+\delta}=(-1)^{\nu}=(-1)^{\mu-1}$; for $\mathfrak s'=\mathfrak s$ the exponent is $\ve(V(\mathfrak s))=-\sigma_L$ by Proposition~\ref{prop.signature}. This proves (i). By Murasugi's theorem \cite{Murasugi63}, $L$ is fibered if and only if $N=1$, and then (i) gives $\mathrm{top}(L,q)=q^{-\sigma_L}$; the formula for the Hopf invariant follows from $\mathrm{top}(L,q)=q^{\nu-2\lambda(\xi_F)}$ \cite[Corollary~1]{LNVdV25}, as in the proof of Theorem~\ref{thm.alternating}. Finally, if $\hate=o'_{j_1}=u_{k_0}$, then by Corollary~\ref{cor.statestrees} the set $S(\tau_0)$, the only set $S(\tau)$ of smallest $\ve$-degree, is $X(\mathfrak s_0)$, and $S'(\tau'_0)=Y(\mathfrak s_0)$ likewise.
\end{proof}

So for a fibered alternating link the top coefficient of $\LG$ and the Hopf invariant of the fiber are determined by the genus, the number of components and the signature. For instance, for the positive and the negative Hopf link, whose fibers are the two Hopf bands, $\nu=1$ and $\sigma_L=\mp1$, so that $\mathrm{top}(L,q)=q^{\pm1}$ and $\lambda(\xi_F)=0$ and $1$ respectively; and for the four knots of Example~\ref{ex.alttrees}, $\sigma_L=-2,0,-2,0$ and $\lambda(\xi_F)=0,1,1,2$.

\begin{example}[the trefoil and the figure-eight knot]\label{ex.states}
Let $D$ be the standard diagram of the trefoil, with the faces $f_0,g_1,g_2,g_3,f_\infty$ of Example~\ref{ex.trefoilkashaev} and the cut edge $\hate=u'_3=o_1$, so that the faces $f_0$ and $g_3$ are removed. At the crossing $j$ the corners $\uparrow$, $\downarrow$, $A$, $B$ lie in $g_j$, $g_{j-1}$, $f_0$ and $f_\infty$. There are three states (Figure~\ref{fig.trefoilstates}), listed with the corners of the crossings $1,2,3$ and with the sign for the order $g_1,g_2,f_\infty$ of the faces:
\[
\begin{array}{c|c|c|c|c|c|c}
   \mathfrak s & \text{corners} & \text{faces} & X(\mathfrak s) & Y(\mathfrak s) & V(\mathfrak s) & \sgn(\mathfrak s)\\ \hline
   \mathfrak s_1 & (\uparrow,\uparrow,B) & (g_1,g_2,f_\infty) & \{1,2\} & \emptyset & \{1,2\} & 1\\
   \mathfrak s_2 & (\uparrow,B,\downarrow) & (g_1,f_\infty,g_2) & \{1\} & \{3\} & \{1,3\} & -1\\
   \mathfrak s_3 & (B,\downarrow,\downarrow) & (f_\infty,g_1,g_2) & \emptyset & \{2,3\} & \{2,3\} & 1
\end{array}
\]
Thus $\det Q_{D,\hate}(x)=1-x_1+x_1x_2=P_{D,\hate}(x)$ (Example~\ref{ex.trefoil}), $\det\tQ_{D,\hate}(y)=1-y_3+y_2y_3$, the signs are $(-1)^{b(\mathfrak s)}$ as in \eqref{eq.blackholes}, with $\eta=1$, no state has an $A$-corner, so that $\sigma_{D,\hate}=1$, and $\ve(X(\mathfrak s))+\ve(Y(\mathfrak s))=2=m_{D,\hate}$ (the three faces have the windings $-1,-1,-2$, and $h(u_j)=1$). The black faces are the three bigons, and $\ve(V(\mathfrak s))=2=F_{\rm b}-1-c_-=-\sigma_L$ for the three states. The nine pairs of states give
\[
   \LG(3_1;-t_0,-t_1)=1+t_0+t_1+t_0^2+2t_0t_1+t_1^2+t_0^2t_1+t_0t_1^2=1+(t_0+t_1)(1+t_0)(1+t_1),
\]
in which the three pairs $(\mathfrak s,\mathfrak s)$ contribute $t_0^2+t_0t_1+t_1^2$, and the six other pairs $(1+t_0t_1)(1+t_0+t_1)$, which vanishes at $t_0t_1=-1$ (Corollary~\ref{cor.diagonal}). In the model of Corollary~\ref{cor.ballinger} these six pairs are three states with multiplicity two, the superpositions of $\mathfrak s_1$ and $\mathfrak s_2$, of $\mathfrak s_1$ and $\mathfrak s_3$ and of $\mathfrak s_2$ and $\mathfrak s_3$; each of them has one cycle, and they contribute $t_0(1+t_0t_1)$, $1+t_0t_1$ and $t_1(1+t_0t_1)$ (Figure~\ref{fig.gammagraph}). The extremal state is $\mathfrak s_3$: $\mathrm{top}(3_1,q)=q^2$ and $\lambda(\xi_F)=0$.

For the closure of $\sigma_1\sigma_2^{-1}\sigma_1\sigma_2^{-1}$, the diagram of the figure-eight knot of Example~\ref{ex.alttrees}, with $\hate=u_1$, there are five states: their corners at the crossings $1,2,3,4$ are $(A,B,B,A)$, $(A,A,\downarrow,\downarrow)$, $(\uparrow,\uparrow,A,A)$, $(\uparrow,A,A,\downarrow)$ and $(A,\uparrow,\downarrow,A)$, with
\begin{align*}
   X(\mathfrak s)&=\{1,4\},\ \{1,2\},\ \{1,2,3,4\},\ \{1,2,3\},\ \{1,2,4\},\\
   Y(\mathfrak s)&=\{1,4\},\ \{1,2,3,4\},\ \{3,4\},\ \{2,3,4\},\ \{1,3,4\}
\end{align*}
respectively; these are the monomials of $P_{D,u_1}$ and $\tP_{D,u_1}$, and $m_{D,\hate}=0$. The sets $V(\mathfrak s)$ are $\emptyset$, $\{3,4\}$, $\{1,2\}$, $\{1,4\}$ and $\{2,3\}$, all of $\ve$-degree $0=-\sigma_L$. The extremal state is the last one: its sets $X=\{1,2,4\}$ and $Y=\{1,3,4\}$ are the sets $S(\tau_0)$ and $S'(\tau'_0)$ of the two extremal trees of Example~\ref{ex.alttrees}, so that $\mathrm{top}(4_1,q)=1$ and $\lambda(\xi_F)=1$. The same happens for the other two knots of that example.
\end{example}

\begin{remark}[the diagonal of an alternating diagram]\label{rem.diagonalalt}
In the variables $t_0=q/p$, $t_1=pq$, a pair $(\mathfrak s,\mathfrak s')$ contributes $p^{\,\ve(Y(\mathfrak s'))-\ve(X(\mathfrak s))}q^{\,\ve(X(\mathfrak s)\triangle Y(\mathfrak s'))}$ to the sum of Theorem~\ref{thm.kauffmanalt}~(ii), by Lemma~\ref{lem.pq}. By Propositions~\ref{prop.statesum}~(i) and~\ref{prop.signature}, the diagonal pairs $(\mathfrak s,\mathfrak s)$ of an alternating diagram contribute
\[
   q^{-\sigma_L}D(p),\qquad D(p)=\sum_{\mathfrak s}p^{\,m_{D,\hate}-2\ve(X(\mathfrak s))} ,
\]
and by Theorem~\ref{thm.kauffmanalt}~(i), $D(p)=\sigma_{D,\hate}\,p^{m_{D,\hate}}P_{D,\hate}\bigl((-p^{-2})^{\ve}\bigr)$: the diagonal is the Alexander polynomial of $L$ in the variable $-p^{-2}$, with positive coefficients, placed in the $q$-degree $-\sigma_L$. Moreover $D(p)=D(p^{-1})$: by Theorem~\ref{thm.kauffmanalt}~(i) at $x_j=y_j=-t^{\ve_j}$ and Lemma~\ref{lem.level}~(i), $\sum_{\mathfrak s}t^{\ve(X(\mathfrak s))}=\sum_{\mathfrak s}t^{\ve(Y(\mathfrak s))}=t^{m_{D,\hate}}\sum_{\mathfrak s}t^{-\ve(X(\mathfrak s))}$.

The state sum of Proposition~\ref{prop.statesum}~(iii) has the same diagonal. At $q^2=-1$, where $t_0t_1=-1$ and $t_0-t_0^{-1}=q\,(p+p^{-1})$, this diagonal is the whole of that sum, by Corollary~\ref{cor.diagonal}:
\[
   q^{-\sigma_L}\sum_{\mathfrak s}p^{\,m_{D,\hate}-2\ve(X(\mathfrak s))}=\nabla_L\bigl(q\,(p+p^{-1})\bigr)\qquad(q^2=-1).
\]
At $p=1$ and $q=\sqrt{-1}$ the left-hand side is $q^{-\sigma_L}$ times the number of states, which is $|\Delta_L(-1)|$ (Theorem~\ref{thm.kauffmanalt}~(i)): this is the classical relation between the determinant and the signature,
\[
   \nabla_L\bigl(2\sqrt{-1}\bigr)=\bigl(\sqrt{-1}\bigr)^{-\sigma_L}\,|\Delta_L(-1)|,
\]
which for a knot reads $\Delta_K(-1)=(-1)^{\sigma_K/2}\,|\Delta_K(-1)|$.
\end{remark}

\begin{remark}[diagrams which are not alternating]\label{rem.nonalt}
The hypothesis of Lemma~\ref{lem.kasteleyn} enters through \eqref{eq.kastcond} only: for an edge which is an edge of the under-strand at both of its ends, or at neither, the product in \eqref{eq.kastcond} is $+1$, and the signs $\sgn(\mathfrak s)(-1)^{|X(\mathfrak s)|}$ then depend on the state. For the closure of $\sigma_1^3\sigma_2\sigma_1^{-1}\sigma_2$, the diagram of $5_2$ of Remark~\ref{rem.loose}, there are eleven states for every cut edge, nine with one sign and two with the other, in agreement with $|\Delta_{5_2}(-1)|=7$, and $\ve(V(\mathfrak s))$ takes the values $2$ and $4$. For the closure of $(\sigma_1\sigma_2)^3$ there are sixteen states for every cut edge, and some of them have the same set $X(\mathfrak s)$ and opposite signs, so that they cancel in $P_{D,\hate}$. The state sum of Proposition~\ref{prop.statesum}~(iii) holds for these diagrams, with cancellations.
\end{remark}

\subsection{The state sum on the Tait graph}\label{ssec.tait}

\emph{This subsection.} So far a Kauffman state was a choice of corners. On an alternating diagram it is also a spanning tree of a plane graph, the \emph{Tait graph} $\CG$ of the diagram, whose vertices are the black faces and whose edges are the crossings \cite{KauffmanFKT}. We rewrite the state sum of Theorem~\ref{thm.kauffmanalt} on this graph. Each of the two flavors $X(\mathfrak s)$ and $Y(\mathfrak s)$ of a state becomes an \emph{orientation} of the graph, that is, a choice of a direction for every edge, which we record as a vector with coordinates $\pm1$. The diagram provides a reference orientation, defined by the over-strands. The monomial of a pair of states is then a product over the edges, in which the factor of an edge only depends on the directions which the two orientations give to this edge, and its two degrees are scalar products of these vectors (Theorem~\ref{thm.tait}). The reference orientation is the sum of an integral flow, carried by the positive crossings, and of an integral cut, carried by the negative ones (Proposition~\ref{prop.omega}); flows and cuts are defined in Definition~\ref{def.flowcut}. This places the state sum in the lattice from which the Heegaard Floer homology of the double branched cover of an alternating link is computed, and in these terms the number $-\sigma_L$ of Proposition~\ref{prop.signature} is a correction term (Remark~\ref{rem.branched}).

\begin{figure}[ht]
\centering
\begin{tikzpicture}[>=Stealth,scale=1.05,every node/.style={font=\small}]
\begin{scope}
\fill[black!13] (0,0) -- (1.2,1.2) -- (-1.2,1.2) -- cycle; \fill[black!13] (0,0) -- (1.2,-1.2) -- (-1.2,-1.2) -- cycle;
\draw (1.2,-1.2) -- (0.22,-0.22); \draw[->] (-0.22,0.22) -- (-1.2,1.2);
\draw[->] (-1.2,-1.2) -- (1.2,1.2);
\draw[->,red,line width=1pt] (0,-0.95) -- (0,0.95);
\draw[->,blue,dashed,line width=0.8pt] (0.95,0) -- (-0.95,0);
\node[font=\scriptsize] at (0.52,1.02) {head}; \node[font=\scriptsize] at (0.45,-1.02) {tail};
\node[font=\scriptsize] at (-1.55,0) {head}; \node[font=\scriptsize] at (1.5,0) {tail};
\node[font=\scriptsize] at (0,1.45) {$\uparrow$}; \node[font=\scriptsize] at (0,-1.45) {$\downarrow$};
\node[font=\scriptsize] at (-2.15,0) {$A$}; \node[font=\scriptsize] at (2.05,0) {$B$};
\node at (0,-2.05) {$\ve_j=1$};
\end{scope}
\begin{scope}[xshift=7.0cm]
\fill[black!13] (0,0) -- (1.2,1.2) -- (1.2,-1.2) -- cycle; \fill[black!13] (0,0) -- (-1.2,1.2) -- (-1.2,-1.2) -- cycle;
\draw (-1.2,-1.2) -- (-0.22,-0.22); \draw[->] (0.22,0.22) -- (1.2,1.2);
\draw[->] (1.2,-1.2) -- (-1.2,1.2);
\draw[->,red,line width=1pt] (0.95,0) -- (-0.95,0);
\draw[->,blue,dashed,line width=0.8pt] (0,0.95) -- (0,-0.95);
\node[font=\scriptsize] at (-1.55,0.22) {head}; \node[font=\scriptsize] at (1.5,0.22) {tail};
\node[font=\scriptsize] at (0.52,-1.02) {head}; \node[font=\scriptsize] at (0.45,1.02) {tail};
\node[font=\scriptsize] at (0,1.45) {$\uparrow$}; \node[font=\scriptsize] at (0,-1.45) {$\downarrow$};
\node[font=\scriptsize] at (-2.15,0) {$B$}; \node[font=\scriptsize] at (2.05,0) {$A$};
\node at (0,-2.05) {$\ve_j=-1$};
\end{scope}
\end{tikzpicture}
\caption{The two edges of a crossing $j$ of an alternating diagram, with the orientations of \eqref{eq.headstails}; the crossing is drawn as in Figure~\ref{fig.corners}, with both strands pointing upward, and its black corners are shaded. The edge of $j$ in the Tait graph $\CG$ (solid, red) joins the two black corners through the crossing, and its edge in the dual graph $\CG^*$ (dashed, blue) joins the two white corners. Each of them goes from its tail to its head. The over-strand crosses the edge of $\CG$ from the left to the right of this edge.}
\label{fig.headstails}
\end{figure}
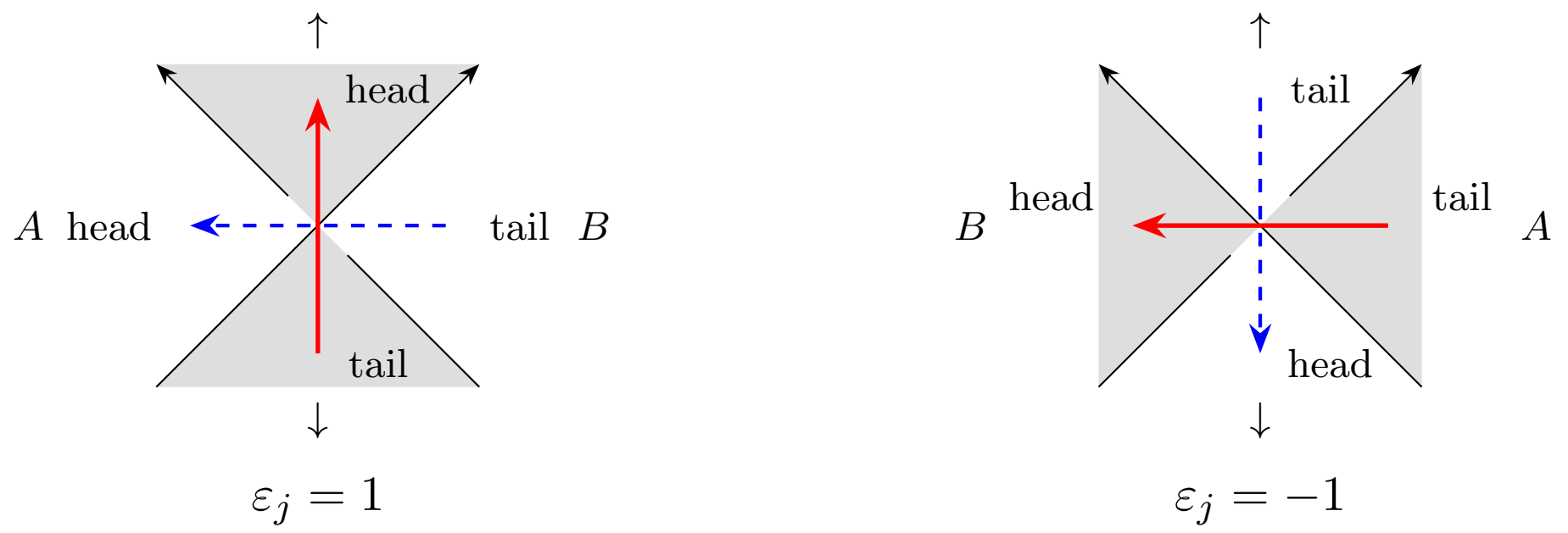

Let $D$ be an alternating diagram, with its black and white faces and corners as in Lemma~\ref{lem.altcoloring}: the black corners of a crossing $j$ are $\uparrow$ and $\downarrow$ if $\ve_j=1$, and its $A$-corner and its $B$-corner if $\ve_j=-1$; a face is black if its corners are. The \emph{Tait graph} $\CG$ of $D$ has the black faces as vertices\label{not.tait} and one edge for each crossing $j$, which joins the faces of the two black corners of $j$ through the crossing. Its dual $\CG^*$ has the white faces as vertices, the edge of $j$ joining the faces of its two white corners. Both are connected graphs drawn on $S^2$, possibly with loops and multiple edges, and every face of $\CG$ is an open disk which contains exactly one white face of $D$. We orient the two edges of every crossing, and we call \emph{heads} of $j$ the two corners at which they end and \emph{tails} of $j$ the two corners at which they start:
\begin{equation}\label{eq.headstails}
\begin{array}{c|cc|cc}
   & \multicolumn{2}{c|}{\text{edge of $j$ in $\CG$}} & \multicolumn{2}{c}{\text{edge of $j$ in $\CG^*$}}\\
   & \text{tail} & \text{head} & \text{tail} & \text{head}\\ \hline
   \ve_j=1 & \downarrow & \uparrow & B & A\\
   \ve_j=-1 & A & B & \uparrow & \downarrow
\end{array}
\end{equation}
These orientations are drawn in Figure~\ref{fig.headstails}. For both signs, the over-strand of $j$ crosses the edge of $j$ in $\CG$ from its left to its right, and the edge of $j$ in $\CG^*$ is this edge turned counterclockwise by a quarter of a turn; at a positive crossing the edge of $\CG$ points in the direction of the two strands. We call this orientation of $\CG$ its \emph{over-orientation}.

\begin{definition}[the lattice of the crossings and its scalar product]\label{def.lattices}
\begin{enumerate}\renewcommand{\theenumi}{\roman{enumi}}
\item We denote by $\BZ^c$ the group of the vectors $a=(a_1,\dots,a_c)$ with integer coordinates, one coordinate for each crossing of $D$, that is, for each edge of $\CG$. The \emph{scalar product} of two such vectors is
\[
   \langle a,b\rangle=\sum_{j=1}^ca_j\,b_j\ \in\BZ .
\]
This is the usual scalar product, for which the $c$ crossings form an orthonormal basis, the one which was already used in Sections~\ref{ssec.altminusone} and~\ref{ssec.underlattice}; it is the only pairing between vectors of $\BZ^c$ which is used in Sections~\ref{ssec.tait} to~\ref{ssec.lattice}. Thus $\langle a,a\rangle=\sum_ja_j^2$, and $a$ and $b$ are \emph{orthogonal} if $\langle a,b\rangle=0$.
\item A \emph{lattice} is a subgroup $N$ of $\BZ^c$, together with the restriction of the scalar product. It is a free abelian group, and its \emph{rank} is the number of elements of a basis. Its \emph{orthogonal complement} is the lattice of the vectors $a\in\BZ^c$ such that $\langle a,x\rangle=0$ for all $x\in N$. If $N$ has rank $c$, its \emph{index} is the number of elements of the finite group $\BZ^c/N$. Two lattices are \emph{isometric} if there is an isomorphism of groups from one onto the other which preserves the scalar products.
\item An \emph{orientation} of $\CG$ is a choice of a direction for each of its edges. We record it by the vector $\chi\in\{\pm1\}^c$ whose coordinate $\chi_j$ is $1$ if the edge of $j$ has its over-orientation, and $-1$ if it has the opposite direction. For two orientations $\chi$ and $\chi'$, the number $\langle\chi,\chi'\rangle$ is the number of the edges to which they give the same direction minus the number of the edges to which they give opposite directions.
\item We write $\1=(1,\dots,1)$ and $\ve=(\ve_1,\dots,\ve_c)$, and $\1_+$ and $\1_-$ for the vectors which are $1$ at the positive, resp.\ at the negative, crossings and $0$ at the other ones, so that $\1=\1_++\1_-$ and $\ve=\1_+-\1_-$. The \emph{writhe} of $D$ is $\mathrm{wr}(D)=\sum_j\ve_j=\langle\ve,\1\rangle$. For a vector $a$, the number $\langle a,\1\rangle$ is the sum of its coordinates, and $\langle a,\ve\rangle$ is the sum of its coordinates at the positive crossings minus the sum of its coordinates at the negative ones.
\end{enumerate}
\end{definition}

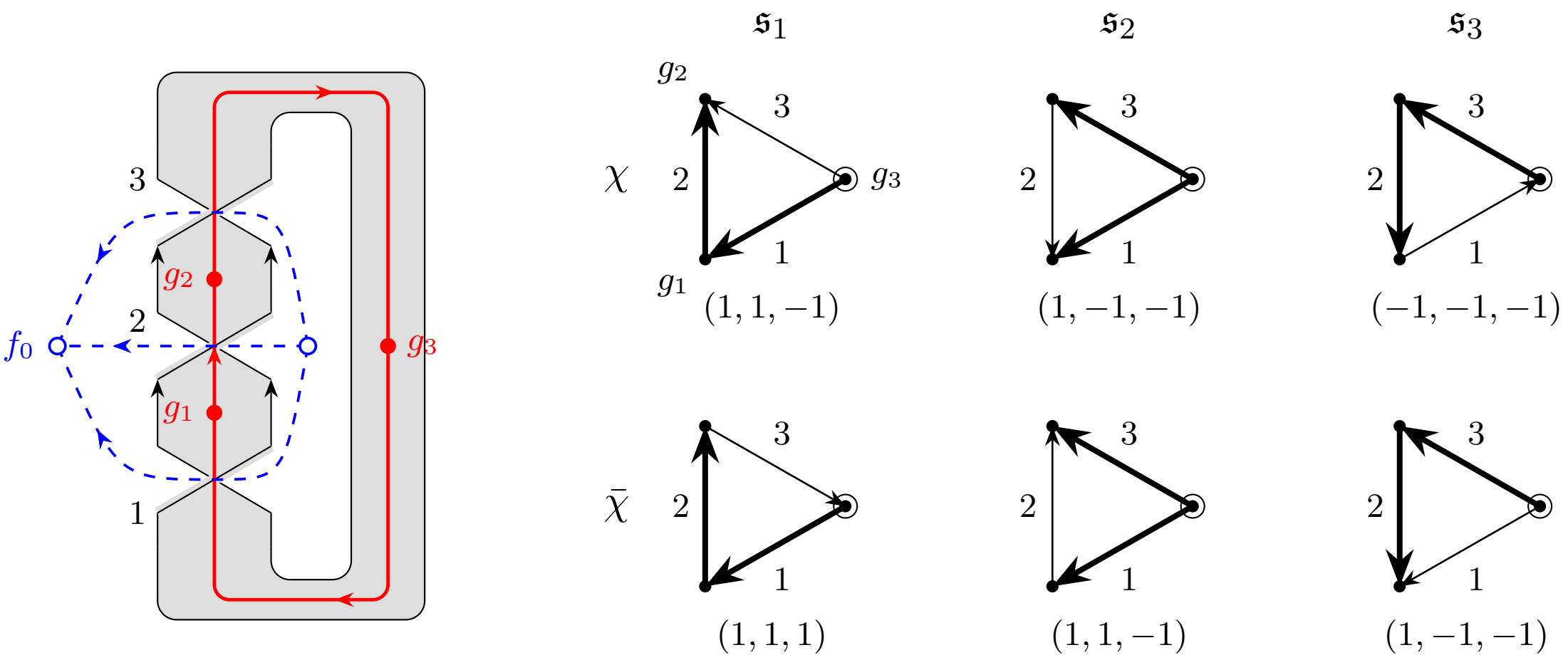
\begin{figure}[ht]
\centering
\begin{tikzpicture}[>=Stealth,scale=0.62]
\begin{scope}
\fill[black!13] (0.85,1) -- (1.7,1.5) -- (1.7,2.5) -- (0.85,3) -- (0,2.5) -- (0,1.5) -- cycle;
\fill[black!13] (0.85,3) -- (1.7,3.5) -- (1.7,4.5) -- (0.85,5) -- (0,4.5) -- (0,3.5) -- cycle;
\fill[black!13] (0.85,5) -- (0,5.5) [rounded corners=5pt] -- (0,7.1) -- (4.0,7.1) -- (4.0,-1.1) -- (0,-1.1) [sharp corners] -- (0,0.5) -- (0.85,1) -- (1.7,0.5) [rounded corners=5pt] -- (1.7,-0.5) -- (2.9,-0.5) -- (2.9,6.5) -- (1.7,6.5) [sharp corners] -- (1.7,5.5) -- cycle;
\draw (0,0) -- (0,0.5); \draw[->] (0,1.5) -- (0,2.5); \draw[->] (0,3.5) -- (0,4.5); \draw (0,5.5) -- (0,6);
\draw[rounded corners=5pt] (0,6) -- (0,7.1) -- (4.0,7.1) -- (4.0,-1.1) -- (0,-1.1) -- (0,0);
\draw (1.7,0) -- (1.7,0.5); \draw[->] (1.7,1.5) -- (1.7,2.5); \draw[->] (1.7,3.5) -- (1.7,4.5); \draw (1.7,5.5) -- (1.7,6);
\draw[rounded corners=5pt] (1.7,6) -- (1.7,6.5) -- (2.9,6.5) -- (2.9,-0.5) -- (1.7,-0.5) -- (1.7,0);
\foreach \y in {1,3,5} { \draw (1.7,{\y-0.5}) -- (0,{\y+0.5}); \draw[preaction={draw=black!13,line width=3pt}] (0,{\y-0.5}) -- (1.7,{\y+0.5}); }
\coordinate (g1) at (0.85,2); \coordinate (g2) at (0.85,4); \coordinate (g3) at (3.45,3); \coordinate (f0) at (-1.5,3); \coordinate (fi) at (2.25,3);
\begin{scope}[red,line width=0.9pt,decoration={markings,mark=at position 0.5 with {\arrow{Stealth[length=5pt]}}}]
\draw[postaction={decorate}] (g1) -- (g2);
\draw[rounded corners=4pt,postaction={decorate}] (g3) -- (3.45,-0.8) -- (0.85,-0.8) -- (g1);
\draw[rounded corners=4pt,postaction={decorate}] (g2) -- (0.85,6.8) -- (3.45,6.8) -- (g3);
\end{scope}
\begin{scope}[blue,dashed,line width=0.7pt,decoration={markings,mark=at position 0.78 with {\arrow{Stealth[length=5pt]}}}]
\foreach \y in {1,3,5} \draw[postaction={decorate}] (fi) .. controls (1.9,\y) .. (0.85,\y) .. controls (-0.6,\y) .. (f0);
\end{scope}
\foreach \p in {g1,g2,g3} \fill[red] (\p) circle (3.4pt);
\foreach \p in {f0,fi} \draw[blue,fill=white,line width=0.7pt] (\p) circle (3.4pt);
\node[font=\scriptsize,red,left] at (0.8,2.0) {$g_1$}; \node[font=\scriptsize,red,left] at (0.8,4.0) {$g_2$}; \node[font=\scriptsize,red,right] at (3.5,3) {$g_3$};
\node[font=\scriptsize,blue,left] at (-1.6,3) {$f_0$};
\node[font=\scriptsize] at (-0.3,0.5) {$1$}; \node[font=\scriptsize] at (-0.3,3.38) {$2$}; \node[font=\scriptsize] at (-0.3,5.5) {$3$};
\end{scope}
\def\tri#1#2#3#4#5#6{
   \coordinate (a) at (0,0); \coordinate (b) at (0,2.4); \coordinate (c) at (2.1,1.2);
   \ifnum#1=1 \draw[->,line width=#4] (c) -- (a); \else \draw[->,line width=#4] (a) -- (c); \fi
   \ifnum#2=1 \draw[->,line width=#5] (a) -- (b); \else \draw[->,line width=#5] (b) -- (a); \fi
   \ifnum#3=1 \draw[->,line width=#6] (b) -- (c); \else \draw[->,line width=#6] (c) -- (b); \fi
   \fill (a) circle (2.6pt); \fill (b) circle (2.6pt); \fill (c) circle (2.6pt); \draw (c) circle (5pt);
   \node[font=\scriptsize,below] at (1.15,0.5) {$1$}; \node[font=\scriptsize,left] at (0,1.2) {$2$}; \node[font=\scriptsize,above] at (1.15,1.9) {$3$};}
\def\thick{1.5pt}\def\thin{0.5pt}
\begin{scope}[xshift=8.2cm,yshift=4.3cm]
\tri{1}{1}{-1}{\thick}{\thick}{\thin} \node[font=\small,anchor=east] at (-0.9,1.2) {$\chi$};
\node[font=\scriptsize,below left] at (0,0) {$g_1$}; \node[font=\scriptsize,above left] at (0,2.4) {$g_2$}; \node[font=\scriptsize,right] at (2.25,1.2) {$g_3$};
\node[font=\small] at (1.0,3.5) {$\mathfrak s_1$};
\node[font=\scriptsize] at (1.0,-0.75) {$(1,1,-1)$};
\end{scope}
\begin{scope}[xshift=13.4cm,yshift=4.3cm]
\tri{1}{-1}{-1}{\thick}{\thin}{\thick}
\node[font=\small] at (1.0,3.5) {$\mathfrak s_2$};
\node[font=\scriptsize] at (1.0,-0.75) {$(1,-1,-1)$};
\end{scope}
\begin{scope}[xshift=18.6cm,yshift=4.3cm]
\tri{-1}{-1}{-1}{\thin}{\thick}{\thick}
\node[font=\small] at (1.0,3.5) {$\mathfrak s_3$};
\node[font=\scriptsize] at (1.0,-0.75) {$(-1,-1,-1)$};
\end{scope}
\begin{scope}[xshift=8.2cm,yshift=-0.6cm]
\tri{1}{1}{1}{\thick}{\thick}{\thin} \node[font=\small,anchor=east] at (-0.9,1.2) {$\bar\chi$};
\node[font=\scriptsize] at (1.0,-0.75) {$(1,1,1)$};
\end{scope}
\begin{scope}[xshift=13.4cm,yshift=-0.6cm]
\tri{1}{1}{-1}{\thick}{\thin}{\thick}
\node[font=\scriptsize] at (1.0,-0.75) {$(1,1,-1)$};
\end{scope}
\begin{scope}[xshift=18.6cm,yshift=-0.6cm]
\tri{1}{-1}{-1}{\thin}{\thick}{\thick}
\node[font=\scriptsize] at (1.0,-0.75) {$(1,-1,-1)$};
\end{scope}
\end{tikzpicture}
\caption{Left: the standard diagram of the trefoil with its black faces $g_1,g_2,g_3$ (shaded), its Tait graph $\CG$ (red) and the dual graph $\CG^*$ (blue, dashed; its two vertices are $f_0$, on the left, and $f_\infty$), with the orientations of \eqref{eq.headstails}: the edge of the crossing $j$ in $\CG$ goes from $g_{j-1}$ to $g_j$, the over-strand crossing it from left to right, and its edge in $\CG^*$ goes from $f_\infty$ to $f_0$. Right: the three Kauffman states of Figure~\ref{fig.trefoilstates} as spanning trees of $\CG$ (thick), with the orientations $\chi(\mathfrak s)$ and $\bar\chi(\mathfrak s)$ of Definition~\ref{def.chi}; the triangle is drawn with $g_1$ at the bottom, $g_2$ at the top and the root $r_{\rm b}=g_3$ (circled) on the right. In $\chi(\mathfrak s)$ the edges of the tree point away from the root, and $\bar\chi(\mathfrak s)$ is obtained by reversing the edge outside the tree. The crossing $j$ belongs to $X(\mathfrak s)$ if $\chi_j(\mathfrak s)=1$, the crossings being positive, and to $Y(\mathfrak s)$ if $\bar\chi_j(\mathfrak s)=-1$ (Theorem~\ref{thm.tait}).}
\label{fig.taittrefoil}
\end{figure}

\begin{figure}[ht]
\centering
\begin{tikzpicture}[>=Stealth,scale=0.9,every node/.style={font=\small},
   mid arrow/.style={postaction={decorate,decoration={markings,mark=at position 0.56 with {\arrow{Stealth}}}}}]
\def\graphpts{\coordinate (a) at (-1.3,0); \coordinate (b) at (1.3,0); \coordinate (c) at (0,-1.9);}
\begin{scope}
\graphpts
\draw[mid arrow] (a) to[bend right=50] (b); \draw[mid arrow] (b) to[bend right=50] (a);
\draw[mid arrow] (c) -- (b); \draw[mid arrow] (c) -- (a);
\fill (a) circle (2pt) (b) circle (2pt) (c) circle (2pt);
\node[font=\scriptsize,left] at (a) {$v_1$}; \node[font=\scriptsize,right] at (b) {$v_2$}; \node[font=\scriptsize,below] at (c) {$v_3$};
\node[font=\scriptsize] at (0,-0.78) {$1$}; \node[font=\scriptsize] at (0,0.78) {$3$};
\node[font=\scriptsize] at (1.0,-1.05) {$2$}; \node[font=\scriptsize] at (-1.0,-1.05) {$4$};
\node at (0,-3.0) {the over-orientation $\1$};
\end{scope}
\begin{scope}[xshift=5.6cm]
\graphpts
\draw[mid arrow,line width=1.6pt] (a) to[bend right=50] (b); \draw[mid arrow,line width=1.6pt] (b) to[bend right=50] (a);
\draw[black!45] (c) -- (b); \draw[black!45] (c) -- (a);
\fill (a) circle (2pt) (b) circle (2pt) (c) circle (2pt);
\node[font=\scriptsize] at (0,-0.82) {$1$}; \node[font=\scriptsize] at (0,0.82) {$1$};
\node[font=\scriptsize,black!60] at (1.0,-1.05) {$0$}; \node[font=\scriptsize,black!60] at (-1.0,-1.05) {$0$};
\node at (0,-3.0) {the flow $\1_+=(1,0,1,0)$};
\end{scope}
\begin{scope}[xshift=11.2cm]
\graphpts
\draw[black!45] (a) to[bend right=50] (b); \draw[black!45] (b) to[bend right=50] (a);
\draw[mid arrow,line width=1.6pt] (c) -- (b); \draw[mid arrow,line width=1.6pt] (c) -- (a);
\fill (a) circle (2pt) (b) circle (2pt) (c) circle (2pt);
\node[font=\scriptsize,left] at (a) {$\psi=1$}; \node[font=\scriptsize,right] at (b) {$\psi=1$}; \node[font=\scriptsize,below] at (c) {$\psi=0$};
\node[font=\scriptsize,black!60] at (0,-0.78) {$0$}; \node[font=\scriptsize,black!60] at (0,0.78) {$0$};
\node[font=\scriptsize] at (1.0,-1.05) {$1$}; \node[font=\scriptsize] at (-1.0,-1.05) {$1$};
\node at (0,-3.0) {the cut $\1_-=(0,1,0,1)$};
\end{scope}
\end{tikzpicture}
\caption{Flows and cuts on the Tait graph of the diagram of the figure-eight knot of Example~\ref{ex.tait}, which has three vertices and four edges; the crossings $1$ and $3$ are positive, and $2$ and $4$ are negative. Left: the graph with its over-orientation and the numbering of its edges. Middle: the vector $\1_+$, which is $1$ on the edges of the positive crossings, is a flow: at $v_1$ and at $v_2$ one unit comes in and one goes out. Right: the vector $\1_-$ is the cut of the function $\psi$ whose values are written at the vertices: it is $1$ on the two edges along which $\psi$ increases by one. Their sum is the over-orientation $\1$ (Proposition~\ref{prop.omega}).}
\label{fig.flowcut}
\end{figure}
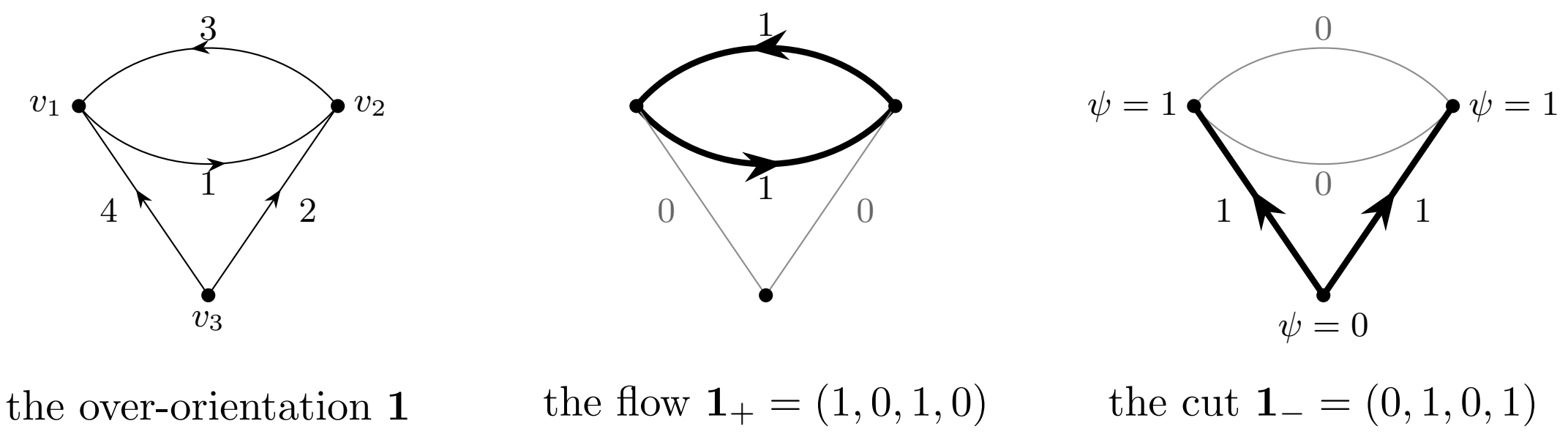

Thus the over-orientation is the vector $\1$, and $-\ve=(-\ve_1,\dots,-\ve_c)$ is the \emph{under-orientation}, for which the under-strand of every crossing crosses its edge from left to right. Figure~\ref{fig.taittrefoil} shows the Tait graph of the trefoil with its over-orientation; the Tait graph of the figure-eight knot is drawn in Figures~\ref{fig.flowcut} and~\ref{fig.fig8smoothings}.

\begin{definition}\label{def.chi}
Let $\mathfrak s$ be a Kauffman state of $(D,\hate)$. We denote by $T(\mathfrak s)$ the set of the crossings $j$ whose corner $\mathfrak s(j)$ is black, and by $\chi(\mathfrak s),\bar\chi(\mathfrak s)\in\{\pm1\}^c$ the vectors with the coordinates
\[
   \chi_j(\mathfrak s)=\begin{cases}\ \ 1&\text{if $\mathfrak s(j)$ is a head of $j$,}\\ -1&\text{if $\mathfrak s(j)$ is a tail of $j$,}\end{cases}\qquad\qquad
   \bar\chi_j(\mathfrak s)=\begin{cases}\ \ \chi_j(\mathfrak s)&\text{if $j\in T(\mathfrak s)$,}\\ -\chi_j(\mathfrak s)&\text{if $j\notin T(\mathfrak s)$.}\end{cases}
\]
\end{definition}

\begin{theorem}[the state sum on the Tait graph]\label{thm.tait}
Let $D$ be an alternating diagram of a link $L$ with $\mu$ components, let $\hate$ be an edge of $D$, and let $r_{\rm b}$ and $r_{\rm w}$ be the black one and the white one of the two faces adjacent to $\hate$.
\begin{enumerate}\renewcommand{\theenumi}{\roman{enumi}}
\item The map $\mathfrak s\mapsto T(\mathfrak s)$ is a bijection from the Kauffman states of $(D,\hate)$ onto the spanning trees of $\CG$, and the crossings which are not in $T(\mathfrak s)$ are the edges of a spanning tree of $\CG^*$. If $j\in T(\mathfrak s)$, the corner $\mathfrak s(j)$ lies in the end of the edge of $j$ which is the farther from $r_{\rm b}$ in the tree $T(\mathfrak s)$; if $j\notin T(\mathfrak s)$, it lies in the end of the edge of $j$ which is the farther from $r_{\rm w}$ in the dual tree. Thus $\chi_j(\mathfrak s)=1$ if and only if the oriented edge of $j$, in the tree to which it belongs, points away from the root $r_{\rm b}$, resp.\ $r_{\rm w}$.
\item For every state $\mathfrak s$,
\begin{gather*}
   X(\mathfrak s)=\{j:\chi_j(\mathfrak s)=\ve_j\},\qquad\qquad Y(\mathfrak s)=\{j:\bar\chi_j(\mathfrak s)=-1\},\\
   V(\mathfrak s)=\{j\in T(\mathfrak s):\ve_j=1\}\cup\{j\notin T(\mathfrak s):\ve_j=-1\},
\end{gather*}
and for all states $\mathfrak s$ and $\mathfrak s'$,
\begin{equation}\label{eq.chidegrees}
\begin{aligned}
   \ve\bigl(X(\mathfrak s)\bigr)&=\tfrac12\bigl(\mathrm{wr}(D)+\langle\chi(\mathfrak s),\1\rangle\bigr),\\
   \ve\bigl(Y(\mathfrak s')\bigr)&=\tfrac12\bigl(\mathrm{wr}(D)-\langle\bar\chi(\mathfrak s'),\ve\rangle\bigr),\\
   \ve\bigl(X(\mathfrak s)\triangle Y(\mathfrak s')\bigr)&=\tfrac12\bigl(\mathrm{wr}(D)+\langle\chi(\mathfrak s),\bar\chi(\mathfrak s')\rangle\bigr).
\end{aligned}
\end{equation}
\item The Links--Gould polynomial is the sum over the ordered pairs of states, that is, of spanning trees of $\CG$,
\begin{equation}\label{eq.taitsum}
   (-1)^{\mu-1}\,\LG(L;-t_0,-t_1)=\sum_{(\mathfrak s,\mathfrak s')}\ \prod_{j=1}^c m_j\bigl(\chi_j(\mathfrak s),\bar\chi_j(\mathfrak s')\bigr),
\end{equation}
with the local weights
\[
\begin{array}{c|cccc}
   \bigl(\chi_j(\mathfrak s),\bar\chi_j(\mathfrak s')\bigr) & (1,1) & (-1,-1) & (1,-1) & (-1,1)\\ \hline
   m_j\ \textup{ if }\ve_j=1 & t_0 & t_1 & 1 & 1\\
   m_j\ \textup{ if }\ve_j=-1 & 1 & 1 & t_1^{-1} & t_0^{-1}
\end{array}
\]
In the variables $t_0=p^{-1}q$, $t_1=pq$ of Lemma~\ref{lem.pq}, the pair $(\mathfrak s,\mathfrak s')$ contributes the monomial $p^{\,\ve(Y(\mathfrak s'))-\ve(X(\mathfrak s))}\,q^{\,\ve(X(\mathfrak s)\triangle Y(\mathfrak s'))}$, whose degrees in $p$ and in $q$ are
\[
   -\tfrac12\bigl(\langle\chi(\mathfrak s),\1\rangle+\langle\bar\chi(\mathfrak s'),\ve\rangle\bigr)\qquad\text{and}\qquad\tfrac12\bigl(\mathrm{wr}(D)+\langle\chi(\mathfrak s),\bar\chi(\mathfrak s')\rangle\bigr) .
\]
\end{enumerate}
\end{theorem}
\begin{proof}
(i) This is the classical description of the states by trees \cite{KauffmanFKT}; we recall the argument because we need the position of the markers. There are three steps: $T(\mathfrak s)$ is a spanning tree; the markers lie at the far ends of the edges; and every spanning tree comes from exactly one state.

\emph{The tree of a state.} Let $\mathfrak s$ be a state and $T=T(\mathfrak s)$. Every black face other than $r_{\rm b}$ contains the corner of exactly one crossing, which belongs to $T$ and whose edge has an end in this face, and $r_{\rm b}$ contains none; so $|T|=F_{\rm b}-1$, and in the same way the $c-|T|=F_{\rm w}-1$ crossings outside $T$ are matched with the white faces other than $r_{\rm w}$. Suppose that $\ell\ge1$ edges of $T$ form a cycle of $\CG$. This cycle is a simple closed curve $\gamma$ on $S^2$ which passes through $\ell$ black faces and $\ell$ crossings and meets $D$ at these crossings only. Let $\Theta$ be the component of $S^2\sm\gamma$ which does not contain $r_{\rm w}$; it does not contain $r_{\rm b}$ either, because $\hate$ lies on the boundary of both faces and does not meet $\gamma$. Let $n$, $m$ and $k$ be the numbers of black faces, of crossings and of white faces contained in $\Theta$. The closure of $\Theta$ is a disk, decomposed by $\CG$ into $n+\ell$ vertices, $m+\ell$ edges and $k$ faces, so that $(n+\ell)-(m+\ell)+k=1$. On the other hand the $n+k$ faces of $D$ contained in $\Theta$ contain the corners of $n+k$ distinct crossings. None of these crossings is on $\gamma$: the corner of $\mathfrak s$ at a crossing of $\gamma$ is black, and it lies in one of the black faces through which $\gamma$ passes. So they are contained in $\Theta$ and $m\ge n+k=m+1$, a contradiction. Hence $T$ has no cycle, and since it has $F_{\rm b}-1$ edges it is a spanning tree of $\CG$.

\emph{The position of the markers.} Direct every edge of $T$ toward the face which contains the corner of $\mathfrak s$: every vertex other than $r_{\rm b}$ is the end point of exactly one directed edge, and $r_{\rm b}$ of none, so along the path of $T$ from $r_{\rm b}$ to any vertex all the edges are directed forward. That is, $\mathfrak s(j)$ lies in the end of the edge of $j$ which is the farther from $r_{\rm b}$. The same argument with the two colors exchanged applies to the crossings outside $T$, to $\CG^*$ and to $r_{\rm w}$.

\emph{From a tree to a state.} Conversely, let $T$ be a spanning tree of $\CG$. By planar duality the crossings outside $T$ are the edges of a spanning tree of the dual graph, and choosing at every crossing the corner in the farther end of its edge, in $T$ or in the dual tree, defines a state $\mathfrak s$ with $T(\mathfrak s)=T$; by what precedes it is the only one. The last sentence of (i) restates the position of the corners, the head of an edge being the corner at which it ends.

(ii) By \eqref{eq.headstails} the corners $\uparrow$ and $A$ are the two heads of $j$ if $\ve_j=1$ and its two tails if $\ve_j=-1$; so, by Definition~\ref{def.state}, $j\in X(\mathfrak s)$ if and only if $\chi_j(\mathfrak s)=\ve_j$. For both signs, the corners $\downarrow$ and $A$ are the tail of the edge of $j$ in $\CG$ and the head of its edge in $\CG^*$; so $j\in Y(\mathfrak s)$ if and only if $\chi_j(\mathfrak s)=-1$ and $j\in T(\mathfrak s)$, or $\chi_j(\mathfrak s)=1$ and $j\notin T(\mathfrak s)$, that is, if and only if $\bar\chi_j(\mathfrak s)=-1$. The vertical corners of a positive crossing are black and those of a negative crossing are white, which gives $V(\mathfrak s)$. Consequently, with $\chi=\chi(\mathfrak s)$ and $\bar\chi'=\bar\chi(\mathfrak s')$, the indicator functions of $X(\mathfrak s)$ and of $Y(\mathfrak s')$ are $j\mapsto\frac12(1+\ve_j\chi_j)$ and $j\mapsto\frac12(1-\bar\chi'_j)$, hence the first two identities of \eqref{eq.chidegrees}, and
\[
   \ve\bigl(X(\mathfrak s)\cap Y(\mathfrak s')\bigr)=\tfrac14\sum_j\ve_j(1+\ve_j\chi_j)(1-\bar\chi'_j)=\tfrac14\bigl(\mathrm{wr}(D)+\langle\chi,\1\rangle-\langle\bar\chi',\ve\rangle-\langle\chi,\bar\chi'\rangle\bigr);
\]
the third one follows, since $\ve(X\triangle Y')=\ve(X)+\ve(Y')-2\ve(X\cap Y')$.

(iii) By (ii), $j\in X(\mathfrak s)\sm Y(\mathfrak s')$ if and only if $\chi_j(\mathfrak s)=\ve_j$ and $\bar\chi_j(\mathfrak s')=1$, and $j\in Y(\mathfrak s')\sm X(\mathfrak s)$ if and only if $\chi_j(\mathfrak s)=-\ve_j$ and $\bar\chi_j(\mathfrak s')=-1$. So the monomial of the pair $(\mathfrak s,\mathfrak s')$ in Theorem~\ref{thm.kauffmanalt}~(ii) is the product of the weights $m_j$ of the table. By Lemma~\ref{lem.pq} this monomial is $p^{\,\ve(Y(\mathfrak s'))-\ve(X(\mathfrak s))}q^{\,\ve(X(\mathfrak s)\triangle Y(\mathfrak s'))}$ in the variables $p,q$, and \eqref{eq.chidegrees} gives its two degrees.
\end{proof}

So the two flavors of a state are two orientations of its Tait graph. In the first one, $\chi(\mathfrak s)$, the edges of the tree $T(\mathfrak s)$ point away from the root $r_{\rm b}$ and the other edges are oriented so that their dual edges point away from $r_{\rm w}$ in the dual tree; up to the conventions for the dual graph, this is the orientation which Bernardi associates with a spanning tree of a rooted graph embedded in a surface \cite[Section~3.2]{Bernardi08}, in the case of the sphere. The second one, $\bar\chi(\mathfrak s)$, is obtained from the first one by reversing the edges outside the tree. The weight of a pair of states is local: a positive crossing contributes $t_0$ or $t_1$ when the first orientation of $\mathfrak s$ and the second orientation of $\mathfrak s'$ agree on its edge, and a negative crossing contributes $t_0^{-1}$ or $t_1^{-1}$ when they differ. In particular the $q$-degree of the pair is the number of positive crossings at which the two orientations agree minus the number of negative crossings at which they differ, and by \eqref{eq.chidegrees} the two orientations are compared with the over-orientation and with the under-orientation only through the $p$-degree, in which the two states do not interact.

The vectors of \eqref{eq.chidegrees} live in the lattice $\BZ^c$ of Definition~\ref{def.lattices}, which contains two classical sublattices.

\begin{definition}[integral flows and integral cuts]\label{def.flowcut}
Let $\CG$ be oriented by its over-orientation, and for a crossing $j$ let $f^+_j$ and $f^-_j$ be the black faces at which the edge of $j$ ends and starts, that is, the faces of its head and of its tail.
\begin{enumerate}\renewcommand{\theenumi}{\roman{enumi}}
\item A vector $x\in\BZ^c$ is an \emph{integral flow} of $\CG$ if at every vertex $v$ what flows in is what flows out: $\sum_{j:f^+_j=v}x_j=\sum_{j:f^-_j=v}x_j$. This is the notion of Definition~\ref{def.flows}, for the graph $\CG$ instead of the diagram $D$. For instance, a cycle of $\CG$, traveled in one of its two directions, defines a flow: $x_j=1$ on the edges which are traveled along their orientation, $x_j=-1$ on those which are traveled against it, and $x_j=0$ on the other edges. The integral flows form a lattice $\mathrm{Flow}(\CG)$.
\item A vector $k\in\BZ^c$ is an \emph{integral cut} of $\CG$ if there is an integer-valued function $\psi$ on the set of the vertices, that is, of the black faces, such that $k_j=\psi(f^+_j)-\psi(f^-_j)$ for every $j$: the coordinate of an edge is the increase of $\psi$ along this edge. For instance, if $R$ is a set of vertices and $\psi$ is equal to $1$ on $R$ and to $0$ elsewhere, then $k_j=1$ on the edges which enter $R$, $k_j=-1$ on the edges which leave $R$, and $k_j=0$ on the other edges. The integral cuts form a lattice $\mathrm{Cut}(\CG)$.
\end{enumerate}
\end{definition}

A flow $x$ is orthogonal to a cut $k$: collecting the terms of each vertex,
\[
   \langle x,k\rangle=\sum_jx_j\bigl(\psi(f^+_j)-\psi(f^-_j)\bigr)=\sum_v\psi(v)\Bigl(\sum_{j:f^+_j=v}x_j-\sum_{j:f^-_j=v}x_j\Bigr)=0 ,
\]
and taking for $\psi$ the function which is $1$ at one vertex and $0$ at the other ones shows that $\mathrm{Flow}(\CG)$ is exactly the orthogonal complement of $\mathrm{Cut}(\CG)$ in $\BZ^c$. Since $\CG$ is connected, a cut determines its function $\psi$ up to a constant; so the rank of $\mathrm{Cut}(\CG)$ is $F_{\rm b}-1$, and that of $\mathrm{Flow}(\CG)$ is $c-F_{\rm b}+1$. The sum of the two lattices is direct and has rank $c$, but it is not the whole of $\BZ^c$ in general: its index is the number of spanning trees of $\CG$ \cite[Lemma~1 and Propositions~1 and~2]{BdlHN97}, which is $|\Delta_L(-1)|$ by Theorems~\ref{thm.tait}~(i) and~\ref{thm.kauffmanalt}~(i). For the trefoil, for instance, $\CG$ is a triangle whose three edges form a directed cycle (Figure~\ref{fig.taittrefoil}): $\mathrm{Flow}(\CG)$ consists of the multiples of $(1,1,1)$, $\mathrm{Cut}(\CG)$ of the vectors $a$ with $a_1+a_2+a_3=0$, and a vector $a\in\BZ^3$ belongs to $\mathrm{Flow}(\CG)\oplus\mathrm{Cut}(\CG)$ if and only if $a_1+a_2+a_3$ is divisible by $3$; the index is $3$, the number of the spanning trees of the triangle and of the states of Figure~\ref{fig.trefoilstates}. So a vector of $\{\pm1\}^c$ does not belong to $\mathrm{Cut}(\CG)\oplus\mathrm{Flow}(\CG)$ in general. The over-orientation does, and its two components are the positive and the negative crossings (Figure~\ref{fig.flowcut}).

\begin{proposition}[the over-orientation is a flow plus a cut]\label{prop.omega}
Let $D$ be an alternating diagram, and let $\1_+,\1_-\in\BZ^c$ be the indicator vectors of the set of positive crossings and of the set of negative crossings, so that $\1=\1_++\1_-$ and $\ve=\1_+-\1_-$. Then
\[
   \1_+\in\mathrm{Flow}(\CG)\qquad\text{and}\qquad\1_-\in\mathrm{Cut}(\CG);
\]
more precisely, if $v_0$ is a black face, $\1_-$ is the vector of the function $\psi=\frac12\bigl(w-w(v_0)\bigr)$ on the black faces, where $w$ is the winding function of Lemma~\ref{lem.winding}. Consequently the over-orientation $\1$ and the under-orientation $-\ve$ belong to $\mathrm{Cut}(\CG)\oplus\mathrm{Flow}(\CG)$, and
\begin{equation}\label{eq.latticesignature}
   \langle\1_+,\1_+\rangle-\rk\,\mathrm{Flow}(\CG)=\rk\,\mathrm{Cut}(\CG)-\langle\1_-,\1_-\rangle=F_{\rm b}-1-c_-=-\sigma_L .
\end{equation}
\end{proposition}
\begin{proof}
Let $v$ be a black face. It is an open disk, and its boundary is a closed walk along which the edges of $D$ and the corners of $v$ alternate; we follow it with $v$ on the left. By the proof of Lemma~\ref{lem.altcoloring}, the black corner adjacent to an end of an edge lies on the left of the edge if the edge is $o'_j$ or $u_j$ at this end, and on its right if it is $u'_j$ or $o_j$. So, by \eqref{eq.corneredges}, at a corner $\uparrow$ of a positive crossing $j$ the walk arrives along $u'_j$, against the orientation of this edge, and leaves along $o'_j$, following its orientation; at a corner $\downarrow$ of a positive crossing it arrives along $u_j$, following the orientation, and leaves along $o_j$, against it; and at a corner of a negative crossing, which is adjacent to $o_j$ and $u'_j$ or to $u_j$ and $o'_j$, it arrives and leaves either both times following the orientation or both times against it. Hence the corners of $v$ at which the walk passes from an edge traversed against its orientation to an edge traversed along it are the corners $\uparrow$ of the positive crossings, and those at which it passes from an edge traversed along its orientation to an edge traversed against it are the corners $\downarrow$ of the positive crossings. The two kinds alternate along the closed walk (compare Figure~\ref{fig.facewalk}), so $v$ contains as many corners $\uparrow$ as corners $\downarrow$ of positive crossings. By \eqref{eq.headstails} this is the equality $\sum_{j:f^+_j=v}(\1_+)_j=\sum_{j:f^-_j=v}(\1_+)_j$, and $\1_+\in\mathrm{Flow}(\CG)$.

By \eqref{eq.wcrossing} the winding numbers of the faces at the corners $\leftarrow$, $\uparrow$, $\downarrow$, $\rightarrow$ of a crossing are $k$, $k-1$, $k-1$, $k-2$. If $\ve_j=1$, the edge of $j$ in $\CG$ goes from $\downarrow$ to $\uparrow$, between two faces with the same winding number; if $\ve_j=-1$, it goes from the $A$-corner $\rightarrow$ to the $B$-corner $\leftarrow$, and the winding number increases by $2$. So $w(f^+_j)-w(f^-_j)=2\,(\1_-)_j$ for every $j$. Since $\CG$ is connected, $w-w(v_0)$ is even on the black faces, $\psi$ is integer-valued, and its vector is $\1_-$.

Finally $\langle\1_\pm,\1_\pm\rangle=c_\pm$ with $c_++c_-=c$, the two ranks are $c-F_{\rm b}+1$ and $F_{\rm b}-1$, and $F_{\rm b}-1-c_-=-\sigma_L$ by Proposition~\ref{prop.signature}.
\end{proof}

The next remark relates these lattices to objects of three-dimensional topology\label{not.glossary} which are not used in the proofs of this monograph; here is what the words mean. The two \emph{checkerboard surfaces} of $D$ are the state surfaces of Section~\ref{sec.statesurfaces} which are made of the black faces, resp.\ of the white faces, joined by a half-twisted band at every crossing; their forms of Gordon and Litherland were studied in Section~\ref{sec.minusone}. The \emph{double branched cover} $\Sigma(L)$ is the closed oriented $3$-manifold which covers the sphere $S^3$ twice outside $L$ and once along $L$; its first homology group is the quotient of $\BZ^c$ by the lattice of the under-strands, a finite group with $|\Delta_L(-1)|$ elements when this number is not zero (Proposition~\ref{prop.underlattice}). The \emph{Heegaard Floer homology} of Ozsv\'ath and Szab\'o attaches to a closed oriented $3$-manifold $Y$ whose first homology group is finite, and to each of its $\mathrm{Spin}^c$ structures $\mathfrak t$, a rational number $d(Y,\mathfrak t)$, its \emph{correction term}; the $\mathrm{Spin}^c$ structures are explained in Definition~\ref{def.spinc} below, and when $H_1(Y;\BZ)$ has odd order, as for the double branched cover of a knot, one of them is distinguished: it is the one which comes from the unique spin structure of $Y$, and we denote it by $\mathfrak t_0$. For the double branched cover of an alternating link the correction terms are computed from the lattice of flows of the Tait graph \cite{OSbranched}, by means of its covectors: a \emph{covector} of a lattice $N$ is a linear form $N\to\BZ$, for instance the form $x\mapsto\langle a,x\rangle$ defined by a vector $a\in\BZ^c$, and it is \emph{characteristic} if its value at every $x\in N$ has the same parity as $\langle x,x\rangle$ (Definition~\ref{def.spinc}). We use these notions through the results we will quote.

\begin{remark}[the lattice of the two checkerboard surfaces]\label{rem.branched}
(i) \emph{The signs of the crossings.} Modulo $2$, Proposition~\ref{prop.omega} says that the positive crossings form an even subgraph of $\CG$, every black face having an even number of corners at positive crossings, and that the negative crossings form a cut of $\CG$, the one between the black faces at which $\psi$ is even and those at which it is odd. For a knot $K$ this determines the signs. Indeed $d=|\Delta_K(-1)|$ is odd and $d\,\BZ^c\subseteq\mathrm{Cut}(\CG)\oplus\mathrm{Flow}(\CG)$, so every vector of $(\BZ/2)^c$ is in a unique way the sum of the reduction of an integral cut and of the reduction of an integral flow. Thus the set of positive crossings of an alternating diagram of a knot is the only even subgraph of its Tait graph whose complement is a cut: the signs of the crossings are determined by the cycle space of $\CG$.

(ii) \emph{The Goeritz lattices.} The graph $\CG$ is the state graph $\Gamma_{\vs_{\rm b}}$ of the black checkerboard smoothing of Section~\ref{ssec.underlattice}, and $\CG^*$ is that of the white one; the over-orientation of the edge of $j$ is its orientation in $\Gamma_{\vs_{\rm b}}$ if $\ve_j=-1$ and the opposite one if $\ve_j=1$, and the orientation of Section~\ref{ssec.underlattice} is the under-orientation. With the restriction of the standard form of $\BZ^c$, the lattice $\mathrm{Flow}(\CG)$ is the form of Gordon and Litherland of the black checkerboard surface, which retracts onto $\CG$, and $\mathrm{Cut}(\CG)$, which is the lattice of flows of $\CG^*$, is the opposite of that of the white one (Corollary~\ref{cor.altdefinite}~(i); \cite{GordonLitherland, Greene13}). In the lattice of flows of the diagram the two appear together: by \eqref{eq.UOw}, $O_j(\partial g)=\ve_j\,U_j(\partial g)$ for a white face $g$ and $O_j(\partial f)=-\ve_j\,U_j(\partial f)$ for a black face $f$, so that, by Proposition~\ref{prop.underlattice}, the map $O=(O_1,\dots,O_c)$ sends the boundary of a black face $f$ to the cut of $\CG$ at the vertex $f$ and the boundaries of the white faces onto generators of $\mathrm{Flow}(\CG)$. Thus
\[
   \mathrm{Flow}(\CG)\oplus\mathrm{Cut}(\CG)=O(\FF_D)
\]
is the lattice of the over-strands of the diagram, and its image under the change of the signs of the coordinates at the negative crossings is the lattice $\Lambda_D=U(\FF_D)$ of the under-strands. The vectors of $\{\pm1\}^c$, that is, the orientations of $\CG$, induce the short characteristic covectors of these two lattices \cite[Section~3.3]{Greene13}, and it is from these covectors that Ozsv\'ath and Szab\'o compute the correction terms $d(\Sigma(L),\mathfrak t)$ of the Heegaard Floer homology of the double branched cover $\Sigma(L)$ \cite{OSbranched}, \cite[Theorem~4.7]{Greene13}. So the state sum \eqref{eq.taitsum} is expressed in the lattice of the two checkerboard surfaces: it involves the two vectors of each state, the vectors $\1_\pm$, and their inner products. What depends on the diagram beyond this lattice is the choice of the vectors $\chi(\mathfrak s)$ and $\bar\chi(\mathfrak s)$, which uses the vertices and the faces of $\CG$ and the root. How much of this is determined by the lattice is discussed in Remark~\ref{rem.pieces}~(ii) and in Section~\ref{ssec.lattice}.

(iii) \emph{The signature as a correction term.} For an alternating knot $K$, Manolescu and Owens proved that $2\,d(\Sigma(K),\mathfrak t_0)=-\sigma_K/2$, where $\mathfrak t_0$ is the spin structure of $\Sigma(K)$ \cite[Theorem~1.2]{ManolescuOwens}; their signature is ours, the right-handed trefoil having signature $-2$. With Proposition~\ref{prop.signature} and \eqref{eq.latticesignature} this gives, for every state $\mathfrak s$ of an alternating diagram of $K$,
\[
   \ve(V(\mathfrak s))=-\sigma_K=4\,d(\Sigma(K),\mathfrak t_0)=\langle\1_+,\1_+\rangle-\rk\,\mathrm{Flow}(\CG) .
\]
In particular the top coefficient of a fibered alternating knot is $\mathrm{top}(K,q)=q^{4d(\Sigma(K),\mathfrak t_0)}$, by Corollary~\ref{cor.extremalstates}. The right-hand side has the form of the lattice formula for the correction terms: in \cite{OSbranched}, \cite[Theorem~4.7]{Greene13} they are, up to the conventions for the signs, the numbers $\frac14(\langle\chi,\chi\rangle-\rk)$ for the shortest characteristic covectors $\chi$ of the lattice of flows in each class modulo twice the lattice. The covector of $\mathrm{Flow}(\CG)$ induced by the over-orientation is the pairing with $\1_+$, because $\1_-$ is orthogonal to the flows; it is characteristic, since $\langle\1_+,x\rangle=\langle\1,x\rangle\equiv\langle x,x\rangle$ modulo $2$ for $x\in\mathrm{Flow}(\CG)$, and it comes from an element of the lattice itself.
\end{remark}

\begin{example}[the trefoil and the figure-eight knot on their Tait graphs]\label{ex.tait}
For the trefoil of Example~\ref{ex.states}, the black faces are the bigons $g_1,g_2,g_3$, and $\CG$ is the triangle whose edge $j$ goes from $g_{j-1}$ to $g_j$, with $g_0=g_3$ (Figure~\ref{fig.taittrefoil}). The over-orientation is a directed cycle: $\1=\1_+$ is a flow, and $\langle\1_+,\1_+\rangle-\rk\,\mathrm{Flow}(\CG)=3-1=2=-\sigma_L$. The roots are $r_{\rm b}=g_3$ and $r_{\rm w}=f_0$, and the three edges of $\CG^*$ go from $f_\infty$ to $f_0$. The states $\mathfrak s_1,\mathfrak s_2,\mathfrak s_3$ are the spanning trees $\{1,2\}$, $\{1,3\}$, $\{2,3\}$, with
\begin{align*}
   \chi(\mathfrak s)&=(1,1,-1),\ (1,-1,-1),\ (-1,-1,-1),\\
   \bar\chi(\mathfrak s)&=(1,1,1),\ (1,1,-1),\ (1,-1,-1)
\end{align*}
respectively, which gives back the sets $X(\mathfrak s)$ and $Y(\mathfrak s)$ of the table of Example~\ref{ex.states}. The numbers $\frac12\bigl(3+\langle\chi(\mathfrak s_a),\bar\chi(\mathfrak s_b)\rangle\bigr)$, $b=1,2,3$, are $2,3,2$ for $a=1$, then $1,2,3$ for $a=2$ and $0,1,2$ for $a=3$: these are the $q$-degrees of the nine monomials of $\LG(3_1;-t_0,-t_1)$.

For the diagram of the figure-eight knot of the same example, $\CG$ has three vertices. The positive crossings $1$ and $3$ are two edges which join two of them in opposite directions, and the negative crossings $2$ and $4$ go from the third vertex to these two. So $\1_+$ is a directed cycle of length two, $\1_-$ is the cut around the third vertex, and both sides of \eqref{eq.latticesignature} are $2-2=0=-\sigma_L$. The spanning trees are the five pairs of crossings other than $\{1,3\}$, and the five states listed in Example~\ref{ex.states} are the trees $\{2,4\}$, $\{2,3\}$, $\{1,4\}$, $\{1,2\}$ and $\{3,4\}$, in this order.
\end{example}

\subsection{Gluing along a Conway sphere}\label{ssec.gluing}

\emph{This subsection.} A sphere which meets a link in four points cuts it into two pieces, called tangles. We show that the Links--Gould polynomial of the link is determined by six polynomials: those of the three links which are obtained by closing the first tangle in three standard ways, and the same for the second tangle (Theorem~\ref{thm.gluing}). This is the analogue for $\LG$ of a classical formula of Conway for the Alexander--Conway polynomial (Remark~\ref{rem.conwaygluing}). One consequence is that $\LG$ does not change under mutation, the operation which turns one of the two tangles by a half-turn before it is glued back (Corollary~\ref{cor.mutation}). The subsection can be read independently of the rest of the section; it is placed here for the following reason.

The state sum of Theorem~\ref{thm.tait} is written on the plane graph $\CG$, with its vertices, its faces and a root, whereas the lattices of Remark~\ref{rem.branched} determine a reduced alternating diagram only up to mutation \cite{Greene13}: they do not see how two pieces of $\CG$ which have two vertices in common are glued. In this subsection we show that the Links--Gould polynomial of a link which is cut into two tangles by a Conway sphere is an explicit bilinear function of the polynomials of three closures of each tangle (Theorem~\ref{thm.gluing}). This holds for all links; for algebraic links, which are built from rational tangles by such gluings, Ishii computed $\LG$ recursively by skein relations \cite{Ishii04}. The formula contains the invariance of $\LG$ under mutation \cite{DeWitKauffmanLinks}, and for an alternating diagram the closures can be taken alternating: their Tait graphs are the two pieces of $\CG$, each of them unchanged, with its two marked vertices identified, or with one or two new edges which join them (Corollary~\ref{cor.twosum} and Remark~\ref{rem.pieces}). Unlike the rest of this section, the proof does not use the determinant formula, but the representation theory which defines $\LG$.

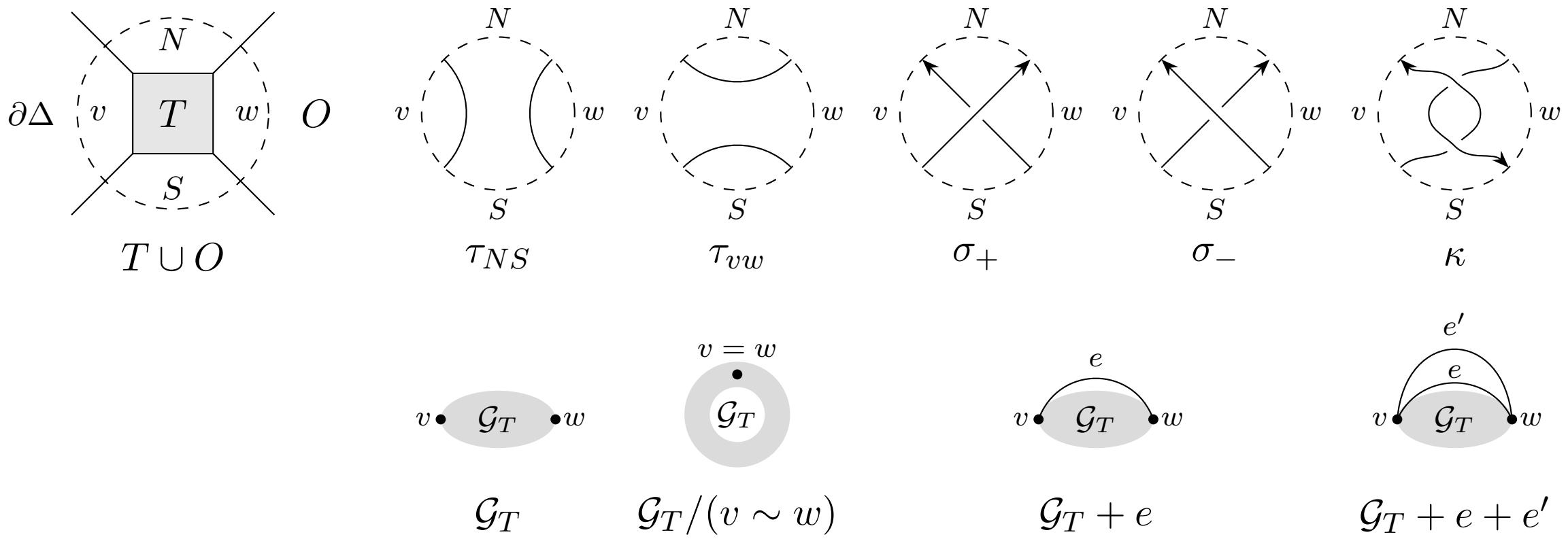
\begin{figure}[ht]
\centering
\begin{tikzpicture}[>=Stealth,scale=0.93,every node/.style={font=\small}]
\begin{scope}
\draw[dashed] (0,0) circle (1);
\foreach \a in {45,135,225,315} \draw (\a:0.55) -- (\a:1.5);
\draw[fill=black!10] (-0.42,-0.42) rectangle (0.42,0.42); \node at (0,0) {$T$};
\node[font=\scriptsize] at (180:0.78) {$v$}; \node[font=\scriptsize] at (90:0.78) {$N$};
\node[font=\scriptsize] at (0:0.78) {$w$}; \node[font=\scriptsize] at (270:0.78) {$S$};
\node at (1.5,0) {$O$}; \node[font=\scriptsize] at (-1.48,0) {$\partial\Delta$};
\node at (0,-1.5) {$T\cup O$};
\end{scope}
\newcommand{\glstd}{\draw[dashed] (0,0) circle (0.8);
 \node[font=\tiny] at (180:1.0) {$v$}; \node[font=\tiny] at (90:1.0) {$N$}; \node[font=\tiny] at (0:1.0) {$w$}; \node[font=\tiny] at (270:1.0) {$S$};}
\begin{scope}[xshift=3.4cm]
\glstd
\draw (225:0.8) to[out=45,in=-45] (135:0.8);
\draw (315:0.8) to[out=135,in=-135] (45:0.8);
\node at (0,-1.5) {$\tau_{NS}$};
\end{scope}
\begin{scope}[xshift=5.9cm]
\glstd
\draw (225:0.8) to[out=45,in=135] (315:0.8);
\draw (135:0.8) to[out=-45,in=-135] (45:0.8);
\node at (0,-1.5) {$\tau_{vw}$};
\end{scope}
\begin{scope}[xshift=8.4cm]
\glstd
\draw[->] (315:0.8) -- (135:0.8);
\begin{scope}\clip (0,0) circle (0.2); \draw[white,line width=5pt] (225:0.8) -- (45:0.8);\end{scope}
\draw[->] (225:0.8) -- (45:0.8);
\node at (0,-1.5) {$\sigma_+$};
\end{scope}
\begin{scope}[xshift=10.9cm]
\glstd
\draw[->] (225:0.8) -- (45:0.8);
\begin{scope}\clip (0,0) circle (0.2); \draw[white,line width=5pt] (315:0.8) -- (135:0.8);\end{scope}
\draw[->] (315:0.8) -- (135:0.8);
\node at (0,-1.5) {$\sigma_-$};
\end{scope}
\begin{scope}[xshift=13.4cm]
\glstd
\draw[->] (45:0.8) to[out=-135,in=45] (0,0.34) to[out=-135,in=90] (-0.27,0) to[out=-90,in=135] (0,-0.34) to[out=-45,in=135] (315:0.8);
\begin{scope}\clip (0,0.34) circle (0.15); \draw[white,line width=5pt] (0.27,0) to[out=90,in=-45] (0,0.34) to[out=135,in=-45] (135:0.8);\end{scope}
\draw[->] (225:0.8) to[out=45,in=-135] (0,-0.34) to[out=45,in=-90] (0.27,0) to[out=90,in=-45] (0,0.34) to[out=135,in=-45] (135:0.8);
\begin{scope}\clip (0,-0.34) circle (0.15); \draw[white,line width=5pt] (-0.27,0) to[out=-90,in=135] (0,-0.34) to[out=-45,in=135] (315:0.8);
 \draw (-0.27,0) to[out=-90,in=135] (0,-0.34) to[out=-45,in=135] (315:0.8);\end{scope}
\node at (0,-1.5) {$\kappa$};
\end{scope}
\newcommand{\glpiece}{\fill[black!14] (0,0) ellipse (0.6 and 0.3); \node[font=\scriptsize] at (0,0) {$\CG_T$};
 \fill (-0.6,0) circle (1.5pt); \fill (0.6,0) circle (1.5pt);
 \node[font=\tiny] at (-0.78,0) {$v$}; \node[font=\tiny] at (0.8,0) {$w$};}
\begin{scope}[yshift=-3.2cm]
\begin{scope}[xshift=3.4cm]
\glpiece
\node at (0,-1.0) {$\CG_T$};
\end{scope}
\begin{scope}[xshift=5.9cm]
\draw[black!14,line width=7pt] (0,0.05) circle (0.42);
\fill (0,0.47) circle (1.5pt); \node[font=\tiny] at (0,0.72) {$v=w$}; \node[font=\scriptsize] at (0,0.05) {$\CG_T$};
\node at (0,-1.0) {$\CG_T/(v\sim w)$};
\end{scope}
\begin{scope}[xshift=9.65cm]
\glpiece
\draw (-0.6,0) to[out=70,in=110,looseness=1.3] (0.6,0); \node[font=\tiny] at (0,0.62) {$e$};
\node at (0,-1.0) {$\CG_T+e$};
\end{scope}
\begin{scope}[xshift=13.4cm]
\glpiece
\draw (-0.6,0) to[out=65,in=115,looseness=1.2] (0.6,0);
\draw (-0.6,0) to[out=85,in=95,looseness=2.1] (0.6,0);
\node[font=\tiny] at (0,0.5) {$e$}; \node[font=\tiny] at (0,1.0) {$e'$};
\node at (0,-1.0) {$\CG_T+e+e'$};
\end{scope}
\end{scope}
\end{tikzpicture}
\caption{The disk $\Delta$ with the inside tangle $T$ and the outside tangle $O$, and the standard tangles, drawn inside $\Delta$; the arcs $v$, $N$, $w$, $S$ of $\partial\Delta$ are on the left, at the top, on the right and at the bottom. The tangles $\sigma_\pm$ are drawn for a disk of parallel type whose entry points are the ends of $S$, and the clasp $\kappa$ for a disk of antiparallel type; its two crossings are positive. Second row: for an alternating diagram in which $v$ and $w$ are black faces, the Tait graphs of the diagrams obtained by closing $T$ with the standard tangles above them, placed outside $\Delta$ (Corollary~\ref{cor.twosum}); for $\sigma_\pm$, the graph of the one of the two diagrams which is alternating.}
\label{fig.gluing}
\end{figure}

\emph{Tangles.}\label{not.tangles} A \emph{Conway sphere} for a link $L\subset S^3$ is a sphere which meets $L$ transversally in four points. It cuts $S^3$ into two balls and $L$ into two \emph{tangles}, the parts of $L$ in the two balls; each tangle consists of two arcs with their ends on the sphere, and possibly of closed components. In a diagram this is seen as follows. Let $D\subset S^2$ be a diagram of an oriented link $L$, and let $\Delta\subset S^2$ be a disk whose boundary meets $D$ transversally in four points which are not crossings. The parts $T$ and $O$ of $D$ inside and outside $\Delta$ are diagrams of the two tangles into which $L$ is cut by a Conway sphere; we write $D=T\cup O$, and we call $T$ the \emph{inside} tangle and $O$ the \emph{outside} tangle. The four points cut $\partial\Delta$ into four arcs, which we denote by $v$, $N$, $w$, $S$ in their cyclic order (Figure~\ref{fig.gluing}). The link enters $\Delta$ at two of the four points and leaves it at the other two. We say that $\Delta$ is of \emph{parallel type} if the two entry points are the ends of one of the four arcs, and of \emph{antiparallel type} if the entry points and the exit points alternate on $\partial\Delta$.

\emph{Standard tangles.} These are the following diagrams in $\Delta$, with the same four ends (Figure~\ref{fig.gluing}). The trivial tangles $\tau_{NS}$ and $\tau_{vw}$ consist of two disjoint arcs without crossing: the complement of $\tau_{NS}$ in $\Delta$ has a component which contains the arcs $N$ and $S$, and that of $\tau_{vw}$ a component which contains $v$ and $w$. The tangles $\sigma_+$ and $\sigma_-$ consist of two arcs which cross once. The \emph{clasp} $\kappa$ is $\tau_{NS}$ with a full twist: its two arcs cross twice, and the arcs $v$ and $w$ lie in two different components of its complement, each of which has a corner at both crossings. A standard tangle has at most one orientation which enters and leaves $\Delta$ at the same points as $D$, and we use the following ones. If $\Delta$ is of parallel type: the one of the two trivial tangles which has such an orientation, which we denote by $\tau_0$, and the tangles $\sigma_+$ and $\sigma_-$, whose crossing is positive, resp.\ negative. If $\Delta$ is of antiparallel type: $\tau_{NS}$, $\tau_{vw}$ and $\kappa$, whose two arcs then run in opposite directions, the twist being the one for which the two crossings are positive. The standard tangles in the disk $S^2\sm\mathrm{int}\,\Delta$ are defined in the same way. For a standard tangle $s$ we denote by $T\cup s$ the diagram formed by $T$ and by $s$ placed outside $\Delta$, and by $s\cup O$ the diagram formed by $s$ placed inside $\Delta$ and by $O$.

\begin{theorem}[gluing along a Conway sphere]\label{thm.gluing}
Let $L$ be the link of the diagram $D=T\cup O$, and let $H_\pm=-(1-t_0^{\pm1})(1-t_1^{\pm1})$ be the Links--Gould polynomials of the positive and of the negative Hopf link.
\begin{enumerate}\renewcommand{\theenumi}{\roman{enumi}}
\item If $\Delta$ is of parallel type, let $L_0,L_+,L_-$ be the Links--Gould polynomials of the links of the diagrams $T\cup\tau_0$, $T\cup\sigma_+$, $T\cup\sigma_-$, and $M_0,M_+,M_-$ those of $\tau_0\cup O$, $\sigma_+\cup O$, $\sigma_-\cup O$. Then
\begin{multline*}
   (H_++H_-)\,\LG(L)=(L_+-L_-)(M_+-M_-)+H_-\,(L_0M_++L_+M_0)\\ +H_+\,(L_0M_-+L_-M_0)-H_+H_-\,L_0M_0 ,
\end{multline*}
where $H_++H_-=-(1-t_0)(1-t_1)(1+t_0t_1)\,t_0^{-1}t_1^{-1}$.
\item If $\Delta$ is of antiparallel type, let $L_{NS},L_{vw},L_\kappa$ be the Links--Gould polynomials of the links of $T\cup\tau_{NS}$, $T\cup\tau_{vw}$, $T\cup\kappa$, and $M_{NS},M_{vw},M_\kappa$ those of $\tau_{NS}\cup O$, $\tau_{vw}\cup O$, $\kappa\cup O$. Then
\[
   \LG(L)=L_{NS}M_{vw}+L_{vw}M_{NS}+\frac{\delta(T)\,\delta(O)}{H_+\,(1+t_0t_1)},
\]
where $\delta(T)=L_\kappa-L_{NS}-H_+L_{vw}$ and $\delta(O)=M_\kappa-M_{NS}-H_+M_{vw}$.
\end{enumerate}
\end{theorem}

The proof uses the definition of $\LG$ by the representation theory of the quantum supergroup $U_q\mathfrak{gl}(2|1)$, which is not used elsewhere in this monograph. Only three properties of this definition are needed, which are recalled in the first paragraph of the proof, together with the fact that the endomorphisms of the tensor square of the representation which commute with the quantum supergroup form a space of dimension three.

\begin{proof}
We use the definition of $\LG$ by the construction of Reshetikhin and Turaev \cite{LinksGould, DeWitKauffmanLinks}. Let $U=U_q\mathfrak{gl}(2|1)$ and let $V$ be the four-dimensional representation of $U$ which defines $\LG$. The construction associates a morphism of $U$-modules to every diagram of an oriented tangle whose ends lie on two horizontal lines: it goes from the tensor product of the modules of the lower ends to that of the upper ends, the module of an end being $V$ if the strand is oriented upward there and $V^*$ otherwise. It is compatible with the composition and with the juxtaposition of diagrams, and a diagram with one lower and one upper end whose closure is a diagram of a link $L'$ has the morphism $\LG(L')\,\mathrm{id}_V$. Both sides of (i), and both sides of (ii) multiplied by $H_+(1+t_0t_1)$, are Laurent polynomials in $t_0$ and $t_1$, so it is enough to prove (i) and (ii) for generic values of the parameters.

Suppose that $\Delta$ is of parallel type. After an isotopy of $S^2$ the disk $\Delta$ is a rectangle with the two entry points on its lower side and the two exit points on its upper side. Then every diagram $T'$ in $\Delta$ which enters and leaves $\Delta$ at the same points as $D$ has a morphism $F(T')$ in $E=\mathrm{End}_U(V\otimes V)$. Cut $O$ at one point. This turns $T'\cup O$ into a diagram with one lower and one upper end which contains the rectangle $\Delta$, and its morphism is obtained from $F(T')$ by tensoring and composing with the morphisms of the pieces of $O$. So there is a linear form $\ell_O$ on $E$, which does not depend on $T'$, such that
\begin{equation}\label{eq.linearform}
   \LG(T'\cup O)=\ell_O\bigl(F(T')\bigr)\qquad\text{for all such }T' .
\end{equation}
For generic values of the parameters, $V\otimes V$ is the direct sum of three simple modules which are pairwise not isomorphic \cite{DeWitKauffmanLinks, DWIL}, so that $\dim E=3$. Let $(s_1,s_2,s_3)=(\tau_0,\sigma_+,\sigma_-)$, and let $G$ be the matrix with the entries $G_{ab}=\LG(s_a\cup s_b)$, where $s_a$ is placed inside and $s_b$ outside $\Delta$. The diagram $\tau_0\cup\tau_0$ consists of two disjoint circles, $\tau_0\cup\sigma_\pm$ and $\sigma_\pm\cup\tau_0$ are diagrams of the unknot with one crossing, $\sigma_+\cup\sigma_+$ and $\sigma_-\cup\sigma_-$ are the diagrams with two crossings of the two Hopf links (Example~\ref{ex.signs} and Remark~\ref{rem.deformedkashaev}~(ii)), and $\sigma_\pm\cup\sigma_\mp$ is a diagram of the unlink with two components. Since $\LG$ vanishes on split links,
\[
   G=\begin{pmatrix}0&1&1\\1&H_+&0\\1&0&H_-\end{pmatrix},\qquad \det G=-(H_++H_-)\neq0 .
\]
If $\sum_ac_aF(s_a)=0$, then \eqref{eq.linearform} for $O=s_b$ gives $\sum_ac_aG_{ab}=0$ for $b=1,2,3$, hence $c=0$: the three morphisms $F(s_a)$ are a basis of $E$. Write $F(T)=\sum_ac_aF(s_a)$. By \eqref{eq.linearform} for the outside tangles $s_b$ and $O$,
\[
   \LG(T\cup s_b)=\sum_ac_a\,G_{ab}\quad(b=1,2,3)\qquad\text{and}\qquad\LG(L)=\sum_ac_a\,\LG(s_a\cup O) .
\]
Solving the first three equations for $c$ gives $\LG(L)=\lambda^{\T}G^{-1}\mu$ with $\lambda=(L_0,L_+,L_-)^{\T}$ and $\mu=(M_0,M_+,M_-)^{\T}$. This is (i), because
\[
   -(H_++H_-)\,G^{-1}=\begin{pmatrix}H_+H_-&-H_-&-H_+\\-H_-&-1&1\\-H_+&1&-1\end{pmatrix} .
\]

Suppose that $\Delta$ is of antiparallel type. We now take for $\Delta$ a rectangle with one entry point and one exit point on each of its horizontal sides, and $E=\mathrm{End}_U(V\otimes V^*)$. This space also has dimension three: by the duality, $\mathrm{End}_U(V\otimes V^*)\cong\mathrm{Hom}_U(V\otimes V^*\otimes V,V)$ and $\mathrm{End}_U(V\otimes V)\cong\mathrm{Hom}_U(V\otimes V\otimes V^*,V)$, and the braiding of the last two factors identifies the two spaces on the right. The same argument applies to $(s_1,s_2,s_3)=(\tau_{NS},\tau_{vw},\kappa)$, for which
\[
   G=\begin{pmatrix}0&1&H_+\\1&0&1\\H_+&1&(3+t_0t_1)H_+\end{pmatrix},\qquad \det G=-H_+\,(1+t_0t_1)\neq0 .
\]
Indeed $\tau_{NS}\cup\tau_{NS}$ and $\tau_{vw}\cup\tau_{vw}$ consist of two disjoint circles, $\tau_{NS}\cup\tau_{vw}$ is a circle, $\tau_{NS}\cup\kappa$ is the diagram of the positive Hopf link with two crossings, $\tau_{vw}\cup\kappa$ is a diagram of the unknot, and $\kappa\cup\kappa$ is the diagram with four positive crossings of the torus link $T(2,4)$, oriented as the boundary of an annulus, for which the determinant formula \eqref{eq.intro.formula} gives the polynomial $(3+t_0t_1)H_+$. Then $\LG(L)=\lambda^{\T}G^{-1}\mu$, with $\lambda=(L_{NS},L_{vw},L_\kappa)^{\T}$ and $\mu=(M_{NS},M_{vw},M_\kappa)^{\T}$, is (ii), because
\[
   -H_+\,(1+t_0t_1)\,G^{-1}=\begin{pmatrix}-1&-(2+t_0t_1)H_+&1\\-(2+t_0t_1)H_+&-H_+^2&H_+\\1&H_+&-1\end{pmatrix} .\qedhere
\]
\end{proof}

\begin{remark}[the Conway polynomial]\label{rem.conwaygluing}
On the curve $t_0t_1=-1$ we have $\LG(L)=\nabla_L(z)$ with $z=t_0-t_0^{-1}$, and $H_\pm=\pm z$ (Section~\ref{ssec.signs}). For a disk of antiparallel type the polynomial $\delta(T)$ vanishes on this curve. Indeed, changing one of the two crossings of $\kappa$ turns $T\cup\kappa$ into a diagram of the link of $T\cup\tau_{NS}$, and smoothing it gives a diagram of the link of $T\cup\tau_{vw}$, so that $\delta(T)$ is there the difference of the two sides of the skein relation \eqref{eq.conwayskein}. Hence $\delta(T)$ and $\delta(O)$ are divisible by $1+t_0t_1$, the last term of (ii) vanishes on the curve, and (ii) becomes Conway's formula for the union of two tangles \cite{Conway70}:
\[
   \nabla_{T\cup O}=\nabla_{T\cup\tau_{NS}}\,\nabla_{\tau_{vw}\cup O}+\nabla_{T\cup\tau_{vw}}\,\nabla_{\tau_{NS}\cup O} .
\]
The last term of (ii) is what $\LG$ adds to this formula, and it is the product of an invariant of $T$ and of an invariant of $O$: since $\delta$ vanishes for the two trivial tangles and $\delta(\kappa)=(1+t_0t_1)H_+$, the quotient $\delta(T)/((1+t_0t_1)H_+)$ is the coefficient of $F(\kappa)$ when $F(T)$ is written in the basis of the proof. For a disk of parallel type both sides of (i) vanish on the curve.
\end{remark}

\begin{figure}[ht]
\centering
\begin{tikzpicture}[>=Stealth,scale=0.95,every node/.style={font=\small}]
\def\tanglebase{%
  \draw[dashed] (0,0) circle (1);
  \foreach \a in {45,135,225,315} \draw (\a:0.55) -- (\a:1.45);
  \draw[fill=black!10] (-0.42,-0.42) rectangle (0.42,0.42);}
\def\letter{\textsf{\bfseries R}}
\begin{scope}
  \tanglebase \node at (0,0) {\letter};
  \node at (0,-2.0) {the tangle $T$};
\end{scope}
\begin{scope}[xshift=4.1cm]
  \tanglebase \node[rotate=180] at (0,0) {\letter};
  \draw[->,thick] (80:1.75) arc (80:10:1.75);
  \draw[->,thick] (260:1.75) arc (260:190:1.75);
  \node at (0,-2.0) {axis perpendicular};
  \node at (0,-2.45) {to the figure};
\end{scope}
\begin{scope}[xshift=8.2cm]
  \tanglebase \node[yscale=-1] at (0,0) {\letter};
  \draw[densely dotted,thick] (-1.9,0) -- (-0.5,0); \draw[densely dotted,thick] (0.5,0) -- (1.9,0);
  \node at (0,-2.0) {horizontal axis};
\end{scope}
\begin{scope}[xshift=12.3cm]
  \tanglebase \node[xscale=-1] at (0,0) {\letter};
  \draw[densely dotted,thick] (0,-1.6) -- (0,-0.5); \draw[densely dotted,thick] (0,0.5) -- (0,1.6);
  \node at (0,-2.0) {vertical axis};
\end{scope}
\end{tikzpicture}
\caption{A tangle $T$ in the disk $\Delta$, shown by a box with the letter R, and its images under the three rotations of angle $\pi$ of the ball which preserve the set of the four points of the Conway sphere: the rotation about the axis perpendicular to the plane of the diagram, and the rotations about two axes in this plane (dotted). For the last two, the diagram of the rotated tangle is the mirror image of the diagram of $T$ in the axis, in which all the crossings have been changed.}
\label{fig.mutation}
\end{figure}
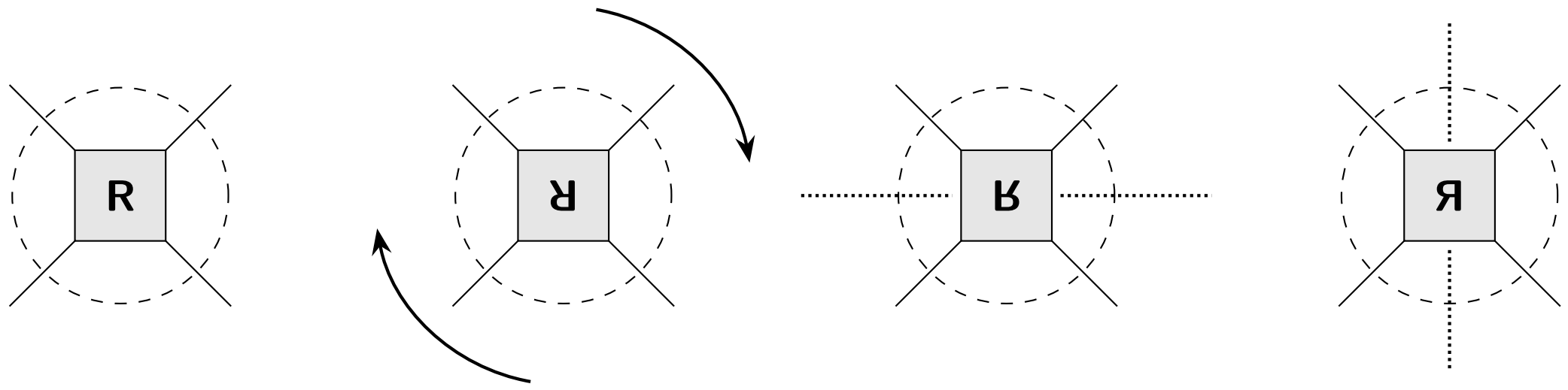

\emph{Mutation.} Recall that a \emph{mutation} of a link $L$ consists in cutting $L$ along a Conway sphere, turning one of the two balls, with its tangle, by a rotation of angle $\pi$ which maps the set of the four points to itself, and gluing it back; the result is a \emph{mutant} of $L$. There are three such rotations, about three perpendicular axes (Figure~\ref{fig.mutation}). A link and its mutants are not isotopic in general, but they are not distinguished by many invariants, among them the Alexander polynomial and the Jones polynomial.

\begin{corollary}[mutation \cite{DeWitKauffmanLinks}]\label{cor.mutation}
Let $\rho$ be one of the three rotations by $\pi$ of the ball bounded by the Conway sphere which preserve the set of the four points, and let $\rho T$ be the image of $T$, with all its components reversed if $\rho$ exchanges the entry points with the exit points. Then $\LG(\rho T\cup O)=\LG(T\cup O)$.
\end{corollary}

\begin{proof}
The rotation $\rho$ extends to a rotation of $S^3$ which preserves the sphere, and the standard tangles outside $\Delta$ are invariant under it up to isotopy, as unoriented tangles. So for each of them the link of $\rho T\cup s$ is the image of the link of $T\cup s$, with all its components reversed if those of $\rho T$ have been; and $\LG$ does not change when all the components of a link are reversed \cite{DeWitKauffmanLinks}. Thus $T$ and $\rho T$ have the same three polynomials in Theorem~\ref{thm.gluing}, and $O$ has not changed.
\end{proof}

\begin{corollary}[the Tait graph]\label{cor.twosum}
Let $D=T\cup O$ be an alternating diagram such that the arcs $v$, $N$, $w$, $S$ lie in four different faces of $D$, those of $v$ and $w$ being black; we denote these two faces by $v$ and $w$ again. Let $\CG_T$ and $\CG_O$ be the subgraphs of the Tait graph $\CG$ formed by the vertices $v$ and $w$ and by the crossings inside, resp.\ outside, $\Delta$ with their ends. Then $\CG=\CG_T\cup\CG_O$ and $\CG_T\cap\CG_O=\{v,w\}$, and the closures of $T$ have the following Tait graphs (Figure~\ref{fig.gluing}).
\begin{enumerate}\renewcommand{\theenumi}{\roman{enumi}}
\item The diagrams $T\cup\tau_{NS}$ and $T\cup\tau_{vw}$ are alternating, with the Tait graphs $\CG_T$ and $\CG_T/(v\sim w)$, the graph $\CG_T$ with $v$ and $w$ identified.
\item If $\Delta$ is of antiparallel type, $T\cup\kappa$ is alternating, with the Tait graph $\CG_T+e+e'$, where $e$ and $e'$ are two new edges which join $v$ and $w$.
\item If $\Delta$ is of parallel type, let $\ve=1$ if the two entry points are the ends of $v$ or of $w$, and $\ve=-1$ if they are the ends of $N$ or of $S$. Then $\tau_0$ is $\tau_{vw}$ if $\ve=1$ and $\tau_{NS}$ if $\ve=-1$; the diagram $T\cup\sigma_\ve$ is alternating, with the Tait graph $\CG_T+e$; and $T\cup\sigma_{-\ve}$ is this diagram with the crossing $e$ changed.
\end{enumerate}
The same holds for the diagrams $s\cup O$ and the graph $\CG_O$.
\end{corollary}

\begin{proof}
A black face of $D$ other than $v$ and $w$ does not meet $\partial\Delta$, so it lies inside or outside $\Delta$ together with the crossings at its corners; this gives the decomposition of $\CG$. Since the four arcs lie in four different faces, $T$ is connected and contains the four points (two components of $T$ would be separated by an arc in $\Delta$ which joins two of the four arcs of $\partial\Delta$), so the diagrams $T\cup s$ are connected. Placing a standard tangle outside $\Delta$ changes neither the crossings of $T$ nor the colors of their corners, and a connected diagram with a checkerboard coloring is alternating, with the coloring of Section~\ref{ssec.tait}, if and only if the black corners of every positive crossing are the vertical ones and those of every negative crossing the horizontal ones. In $T\cup\tau_{NS}$ the arcs $v$ and $w$ lie in two different black faces, and in $T\cup\tau_{vw}$ in the same one; this is (i). Every crossing of $\kappa$ and of $\sigma_\pm$ has a corner in the face of $v$ and a corner in the face of $w$, so it is an edge which joins $v$ and $w$. In $\kappa$ the two arcs run in opposite directions, with $v$ and $w$ on their two sides, so these are the vertical corners, and the crossings have to be positive; this is (ii). In $\sigma_\pm$ they are the vertical corners if the link runs from $v$ to $w$ or from $w$ to $v$ in $\Delta$, that is, if the entry points are the ends of $v$ or of $w$, and the horizontal ones otherwise; in the first case the two arcs of $\tau_0$ join the ends of $v$ to the ends of $w$, so that $\tau_0=\tau_{vw}$, and in the second case $\tau_0=\tau_{NS}$ in the same way. This is (iii).
\end{proof}

\begin{remark}[the pieces of the Tait graph]\label{rem.pieces}
(i) \emph{Alternating closures.} For a disk of parallel type the diagram $T\cup\sigma_{-\ve}$ of Corollary~\ref{cor.twosum} is not alternating. It can be replaced by the alternating diagram $T\cup\sigma_\ve^2$, where $\sigma_\pm^2$ is the tangle with two crossings of the same sign in a row; its Tait graph is $\CG_T$ with two new edges between $v$ and $w$, in series if $\ve=1$ and in parallel if $\ve=-1$. Indeed, the morphism $R=F(\sigma_+)$ satisfies $(R-t_0)(R-t_1)(R+1)=0$ \cite{DeWitKauffmanLinks, GHKKST}, and $F(\sigma_-)=R^{-1}$, so that by linearity, as in \eqref{eq.linearform} with the roles of the two disks exchanged,
\[
   t_0t_1\,\LG(T\cup\sigma_-)=-\LG(T\cup\sigma_+^2)+(t_0+t_1-1)\,\LG(T\cup\sigma_+)+(t_0+t_1-t_0t_1)\,\LG(T\cup\tau_0),
\]
and the same relation holds with $\sigma_+$ and $\sigma_-$ exchanged and $t_0,t_1$ replaced by $t_0^{-1},t_1^{-1}$.

(ii) \emph{The plane graph and the lattice.} By Remark~\ref{rem.branched}~(ii), the state sum of Theorem~\ref{thm.tait} uses the plane graph $\CG$ and not only its lattice of flows. Recall that two graphs are \emph{$2$-isomorphic} if there is a bijection between their sets of edges which maps the cycles of one onto the cycles of the other, and that a graph with more than $k$ vertices is \emph{$k$-connected} if it remains connected whenever fewer than $k$ of its vertices are removed. For a reduced diagram this lattice determines $\CG$ up to $2$-isomorphism \cite[Theorem~1.2]{Greene13}, which for a $2$-connected graph means up to the twists of the pieces attached to the rest of the graph along two vertices \cite{Whitney33}: such a piece is detached, turned over, and attached again with its two vertices exchanged; each of these twists is a mutation of the diagram \cite{Greene13}, by one of the rotations of Corollary~\ref{cor.mutation}. Corollary~\ref{cor.twosum} shows how $\LG$ ignores them: the pieces $\CG_T$ and $\CG_O$ enter through the polynomials of three graphs each, in which the way they are glued does not appear. The theorem can be applied at several disjoint Conway spheres one after the other; a piece which is attached to the rest of the graph along $k$ pairs of vertices then enters through the polynomials of its $3^k$ closures. If $\CG$ is $3$-connected there is no such piece, and $\CG$ has only one embedding in the sphere, up to the homeomorphisms of the sphere \cite{Whitney32}; its vertices and its faces are then determined by its cycle space, and what remains to be chosen in Theorem~\ref{thm.tait} is the root. The terms of the state sum are discussed in Section~\ref{ssec.lattice}.
\end{remark}

\subsection{The lattice of an alternating link}\label{ssec.lattice}

\emph{This subsection.} Up to isometry, the lattice $\mathrm{Flow}(\CG)\oplus\mathrm{Cut}(\CG)\subset\BZ^c$ of the last two subsections only depends on the alternating link, and even only on its double branched cover (Remark~\ref{rem.branched}~(ii)), whereas the state sum of Theorem~\ref{thm.tait} uses more: the plane graph $\CG$, and the edge $\hate$, through the roots $r_{\rm b}$ and $r_{\rm w}$ of the trees. In this subsection we ask how much of the state sum can be expressed in the lattice alone. The motivation is that a formula for $\LG$ in this lattice would be a formula in terms of a three-dimensional object, which does not depend on any diagram. The answer has four parts.
\begin{itemize}
\item For a knot the lattice determines the polynomial, because it determines the diagram up to mutation (Corollary~\ref{cor.latticeknots}).
\item The orientations of the Tait graph fall into classes, which are the $\mathrm{Spin}^c$ structures of the double branched cover (Definition~\ref{def.spinc}), and the states are in bijection with these classes, in two ways, one for each flavor. So the state sum \eqref{eq.taitsum} is a sum over the ordered pairs of classes; but the term of a pair is computed from two orientations which represent the two classes, and these representatives depend on the root (Proposition~\ref{prop.spinc}). The dependence is explicit: every class determines a way of smoothing all the crossings of the diagram, and its representatives only change when the root passes from one circle of this smoothing to another one (Corollary~\ref{cor.classsmoothing}).
\item For every class there is one pair of classes whose term does not depend on the root, and is written in the lattice (Proposition~\ref{prop.nat}). These $|\Delta_L(-1)|$ terms contain the terms of the largest and of the smallest total degree in $t_0$ and $t_1$, which are $c_+$ and $-c_-$ for a reduced alternating diagram with $c_+$ positive and $c_-$ negative crossings; this gives a proof, by the Links--Gould polynomial, of Tait's conjectures on the crossing number and on the writhe (Theorem~\ref{thm.extreme} and Corollary~\ref{cor.crossings}). Here a diagram is \emph{reduced} if it has no \emph{nugatory} crossing, that is, no crossing which has two opposite corners in the same face.
\item The other terms do depend on the root, and two computer experiments indicate that this dependence cannot be removed by a simple device (Remark~\ref{rem.rootdependence}).
\end{itemize}

\begin{corollary}[alternating knots]\label{cor.latticeknots}
Let $D$ and $D'$ be reduced alternating diagrams of two knots $K$ and $K'$, with the Tait graphs $\CG$ and $\CG'$. If the lattices $\mathrm{Flow}(\CG)$ and $\mathrm{Flow}(\CG')$ are isometric, then $\LG(K)=\LG(K')$. In particular $\LG(K)=\LG(K')$ if the double branched covers of $K$ and $K'$ are homeomorphic by a homeomorphism which preserves the orientations.
\end{corollary}

\begin{proof}
The Tait graph of a reduced diagram has no \emph{bridge}, that is, no edge whose removal disconnects the graph, so $\CG$ and $\CG'$ are $2$-isomorphic by \cite[Theorem~1.2]{Greene13}, and $D'$ is obtained from $D$ by a sequence of mutations of unoriented diagrams by \cite[Corollary~4.5]{Greene13}. For a knot the orientation of a mutant is determined up to the reversal of the knot, which does not change $\LG$ \cite{DeWitKauffmanLinks}; so each of these mutations is one of the rotations of Corollary~\ref{cor.mutation}, and $\LG(K)=\LG(K')$. If the double branched covers are homeomorphic as oriented manifolds, their Heegaard Floer homology groups are isomorphic, so $D$ and $D'$ are mutants by \cite[Theorem~1.1]{Greene13}, and the same argument applies.
\end{proof}

The same argument applies to every invariant of knots which does not change under mutation and under reversal. What is particular to $\LG$ is that the state sum \eqref{eq.taitsum} and Theorem~\ref{thm.gluing} compute it from the lattice, in a sense which we now make precise. We need the following vocabulary.

\begin{definition}[characteristic covectors, torsors, $\mathrm{Spin}^c$ structures]\label{def.spinc}
\begin{enumerate}\renewcommand{\theenumi}{\roman{enumi}}
\item Let $N\subseteq\BZ^c$ be a lattice (Definition~\ref{def.lattices}). A \emph{covector} of $N$ is a linear form $\xi\colon N\to\BZ$; for instance, every vector $a\in\BZ^c$ defines the covector $x\mapsto\langle a,x\rangle$ of $N$, which we denote by $\langle a,\cdot\,\rangle$. A covector $\xi$ is \emph{characteristic} if $\xi(x)\equiv\langle x,x\rangle$ modulo $2$ for all $x\in N$. If $\xi$ is characteristic and $y\in N$, then $\xi+2\langle y,\cdot\,\rangle$ is characteristic, and we consider the characteristic covectors \emph{modulo $2N$}, that is, up to the addition of the covectors $2\langle y,\cdot\,\rangle$ with $y\in N$.
\item Let $A$ be an abelian group. A \emph{torsor} over $A$ is a nonempty set $\CS$ on which $A$ acts in such a way that for any two elements $\mathfrak t,\mathfrak t'$ of $\CS$ there is exactly one $a\in A$ which sends $\mathfrak t$ to $\mathfrak t'$. So a torsor has as many elements as the group, and it becomes a copy of the group as soon as one of its elements is chosen as an origin, but no element is preferred; an affine space, with the action of its vector space of translations, is an example. The map $\mathfrak t\mapsto a\cdot\mathfrak t$, for a fixed $a\in A$, is a \emph{translation} of the torsor.
\item Let $Y$ be a closed oriented $3$-manifold. A \emph{$\mathrm{Spin}^c$ structure} on $Y$ is an additional structure on the tangent bundle of $Y$. By a theorem of Turaev \cite{Turaev97} it can be described concretely as a class of nowhere vanishing vector fields on $Y$, two such fields being in the same class if they can be deformed into each other, through nowhere vanishing vector fields, outside a ball. All that matters here is the following. The $\mathrm{Spin}^c$ structures of $Y$ form a torsor $\mathrm{Spin}^c(Y)$ over the group $H_1(Y;\BZ)$; in particular they are finitely many if this group is finite. And the invariants of Heegaard Floer theory, for instance the correction terms $d(Y,\mathfrak t)$, are indexed by the elements $\mathfrak t$ of this set. For the double branched cover $Y=\Sigma(L)$ of an alternating link, the group $H_1(Y;\BZ)$ is isomorphic to $\BZ^c/\Lambda$, where $\Lambda=\mathrm{Flow}(\CG)\oplus\mathrm{Cut}(\CG)$ (Remark~\ref{rem.branched}~(ii) and Proposition~\ref{prop.underlattice}), and $\mathrm{Spin}^c(Y)$ is identified with the set of the classes of Proposition~\ref{prop.spinc} below \cite{OSbranched}, \cite[Theorem~4.7]{Greene13}. A reader who is not familiar with these structures can take this as a definition: in this section, a $\mathrm{Spin}^c$ structure of $\Sigma(L)$ is a class of vectors of $\BZ^c$ with odd coordinates, two vectors $\chi$ and $\chi'$ being in the same class if $\frac12(\chi-\chi')$ belongs to $\Lambda$. Nothing else about $\mathrm{Spin}^c$ structures is used in the proofs.
\end{enumerate}
\end{definition}

\begin{definition}[directed paths, cycles and cuts]\label{def.directed}
Let $\chi\in\{\pm1\}^c$ be an orientation of $\CG$ (Definition~\ref{def.lattices}). It also orients the dual graph $\CG^*$: the edge of $j$ in $\CG^*$ has the orientation of \eqref{eq.headstails} if $\chi_j=1$ and the opposite one if $\chi_j=-1$, so that it is always the oriented edge of $j$ in $\CG$ turned counterclockwise by a quarter of a turn (Figure~\ref{fig.heads}).
\begin{enumerate}\renewcommand{\theenumi}{\roman{enumi}}
\item A \emph{directed path} of $\CG$ is a path all of whose edges are traveled along their orientation, and a \emph{directed cycle} is a cycle with this property; a loop is a directed cycle.
\item For a set $R$ of vertices of $\CG$ we denote by $\partial R$ the set of the edges with exactly one end in $R$. It is a \emph{directed cut} if it is not empty and if all its edges point toward $R$.
\item The orientation $\chi$ is \emph{acyclic} if it has no directed cycle. A \emph{source} is a vertex at which no edge ends.
\item \emph{Reversing} a directed cycle, or a directed cut, means changing the direction of all its edges, that is, the signs of the corresponding coordinates of $\chi$.
\end{enumerate}
The same words are used for the graph $\CG^*$ with the orientation which $\chi$ gives to it.
\end{definition}

\begin{figure}[ht]
\centering
\begin{tikzpicture}[>=Stealth,scale=0.9,every node/.style={font=\scriptsize},
   mid arrow/.style={postaction={decorate,decoration={markings,mark=at position 0.6 with {\arrow{Stealth}}}}}]
\def\tri#1#2#3{%
   \coordinate (a) at (0,0); \coordinate (b) at (0,1.5); \coordinate (c) at (1.3,0.75);
   \ifnum#1=1 \draw[mid arrow] (c) -- (a); \else \draw[mid arrow] (a) -- (c); \fi
   \ifnum#2=1 \draw[mid arrow] (a) -- (b); \else \draw[mid arrow] (b) -- (a); \fi
   \ifnum#3=1 \draw[mid arrow] (b) -- (c); \else \draw[mid arrow] (c) -- (b); \fi
   \fill (a) circle (1.8pt) (b) circle (1.8pt) (c) circle (1.8pt);}
\begin{scope}
\tri{1}{1}{1} \node at (0.6,-0.5) {$+{+}+$}; \node at (0.6,-1.05) {$\bar\chi(\mathfrak s_1)$};
\node[left] at (0,0) {$g_1$}; \node[left] at (0,1.5) {$g_2$}; \node[right] at (1.3,0.75) {$g_3$};
\node at (0.85,0.2) {$1$}; \node at (-0.2,0.75) {$2$}; \node at (0.85,1.3) {$3$};
\end{scope}
\begin{scope}[xshift=2.5cm] \tri{-1}{-1}{-1} \node at (0.6,-0.5) {$-{-}-$}; \node at (0.6,-1.05) {$\chi(\mathfrak s_3)$}; \end{scope}
\draw[black!50] (4.3,-1.9) -- (4.3,2.0);
\begin{scope}[xshift=4.9cm] \tri{1}{1}{-1} \node at (0.6,-0.5) {$+{+}-$}; \node at (0.6,-1.05) {$\chi(\mathfrak s_1)$}; \node at (0.6,-1.55) {${}=\bar\chi(\mathfrak s_2)$}; \end{scope}
\begin{scope}[xshift=7.0cm] \tri{1}{-1}{1} \node at (0.6,-0.5) {$+{-}+$}; \end{scope}
\begin{scope}[xshift=9.1cm] \tri{-1}{1}{1} \node at (0.6,-0.5) {$-{+}+$}; \end{scope}
\draw[black!50] (10.9,-1.9) -- (10.9,2.0);
\begin{scope}[xshift=11.5cm] \tri{1}{-1}{-1} \node at (0.6,-0.5) {$+{-}-$}; \node at (0.6,-1.05) {$\chi(\mathfrak s_2)$}; \node at (0.6,-1.55) {${}=\bar\chi(\mathfrak s_3)$}; \end{scope}
\begin{scope}[xshift=13.6cm] \tri{-1}{1}{-1} \node at (0.6,-0.5) {$-{+}-$}; \end{scope}
\begin{scope}[xshift=15.7cm] \tri{-1}{-1}{1} \node at (0.6,-0.5) {$-{-}+$}; \end{scope}
\end{tikzpicture}
\caption{The eight orientations of the Tait graph of the trefoil, a triangle with the vertices $g_1,g_2,g_3$ and the edges $1,2,3$ (Figure~\ref{fig.taittrefoil}), with the signs of their coordinates: the sign of an edge is $+$ if it has its over-orientation, in which the three edges form the directed cycle $g_3\to g_1\to g_2\to g_3$. Here $\mathrm{Flow}(\CG)\oplus\mathrm{Cut}(\CG)$ consists of the vectors whose sum of coordinates is divisible by $3$, so that the class of an orientation is given by the number of its signs $+$ modulo $3$. There are three classes, separated by the vertical lines: the two directed cycles, which are exchanged by reversing the cycle; and two classes of three acyclic orientations each, in which one passes from an orientation to another one by reversing the directed cut formed by the two edges at its source. These are the three $\mathrm{Spin}^c$ structures of the double branched cover of the trefoil, whose first homology group has three elements. Below four of the orientations we indicate the states of Example~\ref{ex.tait}, for the root $r_{\rm b}=g_3$, of which they are the first orientation $\chi(\mathfrak s)$ or the second orientation $\bar\chi(\mathfrak s)$: each class contains the first orientation of exactly one state and the second orientation of exactly one state, and one orientation can play the two roles, for two different states.}
\label{fig.triangleclasses}
\end{figure}
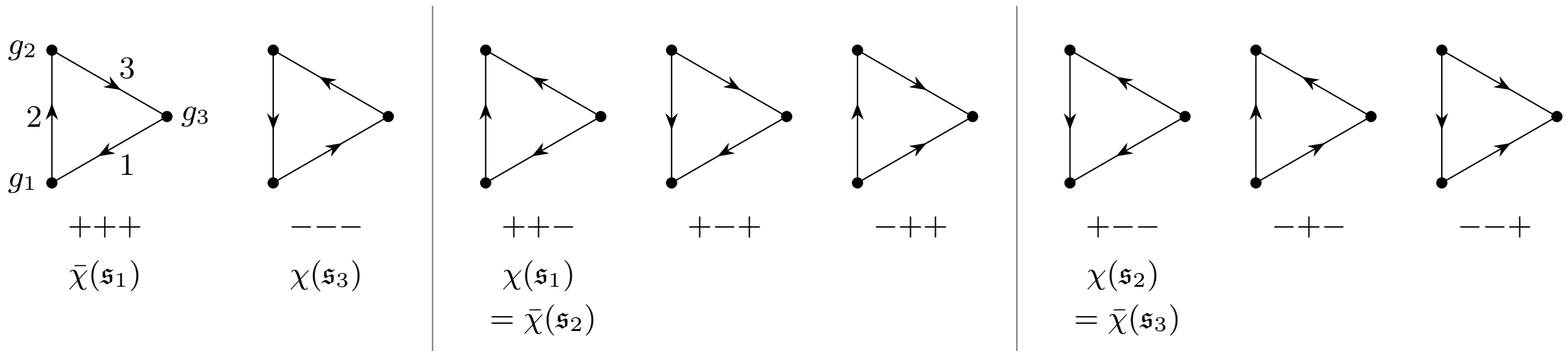

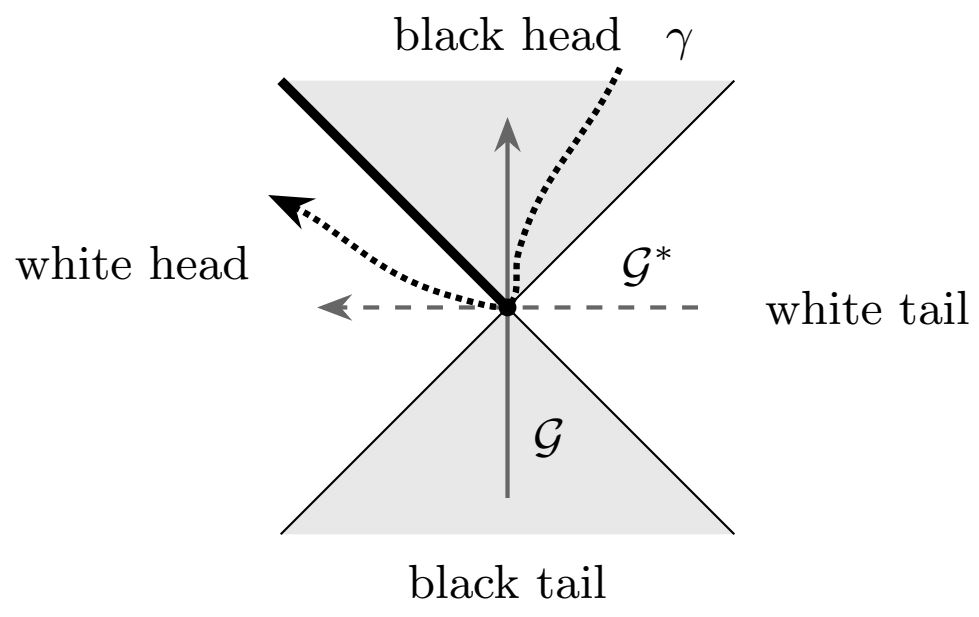
\begin{figure}[ht]
\centering
\begin{tikzpicture}[>=Stealth,scale=1.2]
\fill[black!9] (0,0) -- (1.25,1.25) -- (-1.25,1.25) -- cycle;
\fill[black!9] (0,0) -- (1.25,-1.25) -- (-1.25,-1.25) -- cycle;
\draw (0,0) -- (1.25,1.25); \draw (0,0) -- (1.25,-1.25); \draw (0,0) -- (-1.25,-1.25);
\draw[line width=1.8pt] (0,0) -- (-1.25,1.25);
\draw[->,thick,black!60] (0,-1.05) -- (0,1.05);
\draw[->,thick,black!60,dashed] (1.05,0) -- (-1.05,0);
\draw[->,very thick,densely dotted] (0.62,1.32) to[out=-115,in=75] (0.06,0.3) to[out=-105,in=45] (0,0) to[out=180,in=-15] (-0.3,0.06) to[out=165,in=-25] (-1.32,0.62);
\fill (0,0) circle (1.4pt);
\node[font=\scriptsize,right] at (0.02,-0.72) {$\CG$};
\node[font=\scriptsize,above] at (0.78,0) {$\CG^*$};
\node[font=\scriptsize] at (0.95,1.45) {$\gamma$};
\node[font=\scriptsize,above] at (0,1.3) {black head};
\node[font=\scriptsize,left] at (-1.3,0.25) {white head};
\node[font=\scriptsize,below] at (0,-1.3) {black tail};
\node[font=\scriptsize,right] at (1.3,0) {white tail};
\end{tikzpicture}
\caption{A crossing $j$ with the edges of $j$ in $\CG$ and in $\CG^*$, for an orientation $\chi$; the black faces are shaded. The two heads of $j$ are adjacent corners, separated by one of the four ends of edges of $D$ at $j$ (thick). A curve $\gamma$ which passes from the black head to the white head has this end on its right, and the three other ones on its left.}
\label{fig.heads}
\end{figure}

Here is, in plain words, what Proposition~\ref{prop.spinc} below says. The $2^c$ orientations of the Tait graph are grouped into classes: two orientations $\chi$ and $\chi'$ are in the same class if $\frac12(\chi-\chi')$, a vector with coordinates $0$ and $\pm1$, is the sum of an integral flow and of an integral cut; by Lemma~\ref{lem.minty} below this means that $\chi'$ is obtained from $\chi$ by reversing directed cycles and directed cuts, one after the other. There are as many classes as spanning trees of $\CG$, that is, as states, and these classes are the $\mathrm{Spin}^c$ structures of $\Sigma(L)$. Every state $\mathfrak s$ has two orientations $\chi(\mathfrak s)$ and $\bar\chi(\mathfrak s)$, one for each flavor (Definition~\ref{def.chi}), and each class contains the orientation $\chi(\mathfrak s)$ of exactly one state $\mathfrak s$ and the orientation $\bar\chi(\mathfrak s')$ of exactly one state $\mathfrak s'$. So a state is a way of choosing an orientation in a class, and the state sum, which is a sum over the pairs of states, becomes a sum over the pairs of classes. Figure~\ref{fig.triangleclasses} shows the classes for the trefoil.

\begin{proposition}[the states and the $\mathrm{Spin}^c$ structures]\label{prop.spinc}
Let $D$ be an alternating diagram of a link $L$, let $\hate$ be an edge of $D$, and let $\Lambda=\mathrm{Flow}(\CG)\oplus\mathrm{Cut}(\CG)$. We say that two vectors $\chi,\chi'$ of $\BZ^c$ with odd coordinates are \emph{equivalent} if $\frac12(\chi-\chi')\in\Lambda$, and we denote by $[\chi]$ the class of $\chi$.
\begin{enumerate}\renewcommand{\theenumi}{\roman{enumi}}
\item The classes form a torsor over the group $\BZ^c/\Lambda$, with $|\Delta_L(-1)|$ elements. Restricting the linear forms $\langle\chi,\cdot\,\rangle$ to $\mathrm{Flow}(\CG)$ identifies it with the set of the characteristic covectors of $\mathrm{Flow}(\CG)$ modulo $2\,\mathrm{Flow}(\CG)$, and the same holds for $\mathrm{Cut}(\CG)$. This is the set of the $\mathrm{Spin}^c$ structures of the double branched cover $\Sigma(L)$ \cite{OSbranched}, \cite[Theorem~4.7]{Greene13}.
\item The maps $\mathfrak s\mapsto[\chi(\mathfrak s)]$ and $\mathfrak s\mapsto[\bar\chi(\mathfrak s)]$ are bijections from the set of the states of $(D,\hate)$ onto the set of the classes.
\item The vector $\chi(\mathfrak s)$ is the only vector of $\{\pm1\}^c$ in its class for which every black face is the end of a directed path of $\CG$ which starts at $r_{\rm b}$, and every white face the end of a directed path of $\CG^*$ which starts at $r_{\rm w}$. The vector $\bar\chi(\mathfrak s)$ is the only vector of $\{\pm1\}^c$ in its class for which every black face is the end of a directed path of $\CG$ which starts at $r_{\rm b}$, and every white face the origin of a directed path of $\CG^*$ which ends at $r_{\rm w}$.
\item If $\hate$ is replaced by another edge, and the states are identified with the spanning trees of $\CG$ by Theorem~\ref{thm.tait}~(i), the two bijections of (ii) are composed with translations of the torsor.
\end{enumerate}
Consequently, if $u_{\mathfrak t}=\chi(\mathfrak s)$ and $\bar u_{\mathfrak t}=\bar\chi(\mathfrak s')$ denote the vectors of the two states with $[\chi(\mathfrak s)]=\mathfrak t$ and $[\bar\chi(\mathfrak s')]=\mathfrak t$, the state sum \eqref{eq.taitsum} is the sum over the ordered pairs of classes
\begin{equation}\label{eq.spincsum}
   (-1)^{\mu-1}\,\LG(L;-t_0,-t_1)=\sum_{(\mathfrak t,\mathfrak t')}\ t_0^{\,a(\mathfrak t,\mathfrak t')}\,t_1^{\,b(\mathfrak t,\mathfrak t')},\qquad
   \begin{aligned} a(\mathfrak t,\mathfrak t')&=\tfrac14\langle u_{\mathfrak t}+\ve,\,\bar u_{\mathfrak t'}+\1\rangle,\\ b(\mathfrak t,\mathfrak t')&=\tfrac14\langle u_{\mathfrak t}-\ve,\,\bar u_{\mathfrak t'}-\1\rangle .\end{aligned}
\end{equation}
\end{proposition}

\begin{proof}
(i) The group $\BZ^c$ acts simply transitively on the vectors with odd coordinates by $a\cdot\chi=\chi+2a$, and the classes are the orbits of $\Lambda$. So they form a torsor over $\BZ^c/\Lambda$, and the index of $\Lambda$ in $\BZ^c$ is the number of spanning trees of $\CG$, which is $|\Delta_L(-1)|$ (Theorems~\ref{thm.tait}~(i) and~\ref{thm.kauffmanalt}~(i); compare Proposition~\ref{prop.underlattice}). If $\chi$ has odd coordinates and $x\in\mathrm{Flow}(\CG)$, then $\langle\chi,x\rangle\equiv\langle\1,x\rangle\equiv\langle x,x\rangle$ modulo $2$, so the restriction of $\langle\chi,\cdot\,\rangle$ is characteristic. The restrictions of $\langle\chi,\cdot\,\rangle$ and of $\langle\chi',\cdot\,\rangle$ differ by $2\langle y,\cdot\,\rangle$ with $y\in\mathrm{Flow}(\CG)$ if and only if $\frac12(\chi-\chi')-y$ is orthogonal to $\mathrm{Flow}(\CG)$, that is, belongs to $\mathrm{Cut}(\CG)$, which is a \emph{primitive} sublattice of $\BZ^c$: a vector of $\BZ^c$ which has a nonzero multiple in $\mathrm{Cut}(\CG)$ belongs to $\mathrm{Cut}(\CG)$. So the map from the classes to the characteristic covectors modulo $2\,\mathrm{Flow}(\CG)$ is injective, and both sets have as many elements as $\CG$ has spanning trees \cite{BdlHN97}. The same argument applies to $\mathrm{Cut}(\CG)$.

(ii) and (iii). Let us say that a vector of $\{\pm1\}^c$ is connected from $r_{\rm b}$ if every black face is the end of a directed path of $\CG$ which starts at $r_{\rm b}$, connected from $r_{\rm w}$ if every white face is the end of a directed path of $\CG^*$ which starts at $r_{\rm w}$, and connected to $r_{\rm w}$ if every white face is the origin of a directed path of $\CG^*$ which ends at $r_{\rm w}$. By Theorem~\ref{thm.tait}~(i), $\chi(\mathfrak s)$ is connected from $r_{\rm b}$ and from $r_{\rm w}$; and $\bar\chi(\mathfrak s)$, which differs from it on the edges outside $T(\mathfrak s)$, is connected from $r_{\rm b}$ and to $r_{\rm w}$.

We first show that if two equivalent vectors $\chi,\chi'\in\{\pm1\}^c$ are connected from $r_{\rm b}$, then $\frac12(\chi-\chi')\in\mathrm{Flow}(\CG)$. Write $\frac12(\chi-\chi')=x+k$ with $x\in\mathrm{Flow}(\CG)$ and $k=\bigl(\psi(f^+_j)-\psi(f^-_j)\bigr)_j\in\mathrm{Cut}(\CG)$, and suppose that the function $\psi$ is not constant. Let $S$ be the set of the black faces at which it is maximal, and let $\partial S\neq\emptyset$ be the set of the edges of $\CG$ with exactly one end in $S$. The vector $\kappa\in\mathrm{Cut}(\CG)$ of the indicator function of $S$ is supported on $\partial S$, and $\kappa_j=1$ if $f^+_j\in S$, $\kappa_j=-1$ if $f^-_j\in S$. Every edge of $\partial S$ contributes at least $1$ to $\langle k,\kappa\rangle$, so $\langle\frac12(\chi-\chi'),\kappa\rangle=\langle k,\kappa\rangle\ge|\partial S|$. On the other hand $\langle\chi,\kappa\rangle=2n_\chi-|\partial S|$, where $n_\chi$ is the number of edges of $\partial S$ which $\chi$ orients toward $S$. Hence $n_\chi-n_{\chi'}\ge|\partial S|$: all the edges of $\partial S$ point toward $S$ for $\chi$, and away from $S$ for $\chi'$. If $r_{\rm b}\in S$, no directed path of $\chi$ leaves $S$, and if $r_{\rm b}\notin S$, no directed path of $\chi'$ enters $S$; in both cases one of the two vectors is not connected from $r_{\rm b}$. So $\psi$ is constant and $k=0$. The same argument applies to two vectors which are connected to a root, and to the graph $\CG^*$, whose lattices of cuts and of flows are $\mathrm{Flow}(\CG)$ and $\mathrm{Cut}(\CG)$. Therefore, if $\chi$ and $\chi'$ are equivalent, connected from $r_{\rm b}$, and both connected from $r_{\rm w}$ or both connected to $r_{\rm w}$, then $\frac12(\chi-\chi')\in\mathrm{Flow}(\CG)\cap\mathrm{Cut}(\CG)=0$. This is the uniqueness in (iii), and it shows that $\chi(\mathfrak s)=\chi(\mathfrak s')$ if $[\chi(\mathfrak s)]=[\chi(\mathfrak s')]$, and the same for $\bar\chi$.

We now show that a state is determined by its vector $\chi(\mathfrak s)$. For an orientation $\chi$, we call heads of a crossing $j$ the corners at which the edges of $j$ in $\CG$ and in $\CG^*$ end. They are adjacent, and separated by one of the four ends of edges of $D$ at $j$ (Figure~\ref{fig.heads}); and by Definition~\ref{def.chi} the corner $\mathfrak s(j)$ is a head of $j$ for the orientation $\chi(\mathfrak s)$. Suppose that $\mathfrak s\neq\mathfrak s'$ are two states with $\chi(\mathfrak s)=\chi(\mathfrak s')$. Both are bijections from the crossings onto the faces other than $r_{\rm b}$ and $r_{\rm w}$, so there are $\ell\ge2$ distinct crossings $j_1,\dots,j_\ell$ and distinct faces $F_1,\dots,F_\ell$ with $\mathfrak s(j_k)=F_k$ and $\mathfrak s'(j_k)=F_{k+1}$, the indices being taken modulo $\ell$. The faces $F_k$ and $F_{k+1}$ contain the two heads of $j_k$, so one is black and the other one white, and the colors alternate. Joining in every $F_k$ the corners at $j_{k-1}$ and at $j_k$ by an arc, we obtain a simple closed curve $\gamma$ which meets $D$ at the crossings $j_k$ only. It does not meet $\hate$, so $r_{\rm b}$ and $r_{\rm w}$ lie in the same component of $S^2\sm\gamma$; let $\Theta$ be the other one. Let $m$ be the number of crossings in $\Theta$ and $n$ the number of faces contained in $\Theta$. A crossing in $\Theta$ has its corners in faces which are contained in $\Theta$ or are among the $F_k$, and $\mathfrak s$ maps the crossings of $\gamma$ onto the $F_k$; so $\mathfrak s$ maps the crossings in $\Theta$ bijectively onto the faces contained in $\Theta$, and $n=m$. The closed disk $\overline\Theta$ is decomposed into the $m+\ell$ crossings, the edges of $D$ in $\Theta$ and the $\ell$ arcs of $\gamma$, and the $n$ faces and the $\ell$ parts of the faces $F_k$ which lie in $\Theta$. If $h_k\in\{1,3\}$ is the number of ends of edges of $D$ at $j_k$ which lie in $\Theta$, there are $2m+\frac12\sum_kh_k$ edges of $D$ in $\Theta$, and the Euler characteristic of the disk gives
\[
   1=(m+\ell)-\bigl(2m+\tfrac12\textstyle\sum_kh_k+\ell\bigr)+(n+\ell),\qquad\text{that is,}\qquad \textstyle\sum_kh_k=2\ell-2 .
\]
But follow $\gamma$ through $F_1,j_1,F_2,j_2,\dots$ The end of edge which separates the two heads of $j_k$ is on the right of $\gamma$ if $F_k$ is black and on its left if $F_k$ is white (Figure~\ref{fig.heads}). As the colors alternate, it is on the side of $\Theta$ for exactly half of the crossings $j_k$, so that $\sum_kh_k=\frac\ell2+\frac{3\ell}2=2\ell$, a contradiction. For $\bar\chi$ the argument is the same: for the orientation $\bar\chi(\mathfrak s)$, the corner $\mathfrak s(j)$ is the head of the edge of $j$ in $\CG$ or the tail of the edge of $j$ in $\CG^*$, and these two corners are also adjacent.

So $\mathfrak s\mapsto[\chi(\mathfrak s)]$ and $\mathfrak s\mapsto[\bar\chi(\mathfrak s)]$ are injective. By Theorem~\ref{thm.tait}~(i) there are as many states as spanning trees of $\CG$, that is, as classes; this proves (ii).

(iv) Let $\hate'$ be another edge, with the roots $r'_{\rm b}$ and $r'_{\rm w}$, and let $\mathfrak s$ and $\mathfrak s'$ be the states of $(D,\hate)$ and of $(D,\hate')$ with the same tree $T$. By Theorem~\ref{thm.tait}~(i) the vectors $\chi(\mathfrak s)$ and $\chi(\mathfrak s')$ differ exactly on the path of $T$ from $r_{\rm b}$ to $r'_{\rm b}$ and on the path of the dual tree from $r_{\rm w}$ to $r'_{\rm w}$, and $\frac12(\chi(\mathfrak s)-\chi(\mathfrak s'))=\pi+\pi^*$, where $\pi,\pi^*\in\{0,\pm1\}^c$ are these two oriented paths. For another tree, $\pi$ changes by a closed path of $\CG$, which is an element of $\mathrm{Flow}(\CG)$, and $\pi^*$ by a closed path of $\CG^*$, which is an element of $\mathrm{Cut}(\CG)$. So the class of $\pi+\pi^*$ in $\BZ^c/\Lambda$ does not depend on $T$. For $\bar\chi$ the same holds with $\pi-\pi^*$.

Finally, the exponents $a$ and $b$ in \eqref{eq.spincsum} are those of the table of Theorem~\ref{thm.tait}~(iii): for $\ve_j=1$ the crossing $j$ contributes $\frac14(\chi_j+1)(\bar\chi'_j+1)$ to the exponent of $t_0$ and $\frac14(\chi_j-1)(\bar\chi'_j-1)$ to that of $t_1$, and for $\ve_j=-1$ it contributes $\frac14(\chi_j-1)(\bar\chi'_j+1)$ and $\frac14(\chi_j+1)(\bar\chi'_j-1)$.
\end{proof}

For the vectors $\chi(\mathfrak s)$, statements (ii) and (iv) can also be deduced from the results of Bernardi and of Baker and Wang on the spanning trees of a plane graph: the orientations $\chi(\mathfrak s)$ are those of \cite{Bernardi08}, and by \cite{BakerWang} the resulting action of the group $\BZ^c/\Lambda$ on the spanning trees does not depend on the root.

It remains to see how the representatives $u_{\mathfrak t}$ and $\bar u_{\mathfrak t}$ depend on the root, and which terms of \eqref{eq.spincsum} do not depend on it.

Let $\chi\in\{\pm1\}^c$. All the directed paths, cycles and cuts below are taken with respect to the orientation $\chi$ of $\CG$ and of $\CG^*$; a directed cycle may be a loop. As in the proof of Proposition~\ref{prop.spinc}, we call \emph{black head} and \emph{black tail} of a crossing $j$ the corners at which its edge in $\CG$ ends and starts, and \emph{white head} and \emph{white tail} those at which its edge in $\CG^*$ ends and starts; around the crossing they come in the order black head, white head, black tail, white tail, counterclockwise (Figure~\ref{fig.heads}). We say that $\chi$ is \emph{connected from} a black face $v$ if every black face is the end of a directed path of $\CG$ which starts at $v$, \emph{connected from} a white face $f$ if every white face is the end of a directed path of $\CG^*$ which starts at $f$, and \emph{connected to} $f$ if every white face is the origin of a directed path of $\CG^*$ which ends at $f$; for $v=r_{\rm b}$ and $f=r_{\rm w}$ these are the notions of the proof of Proposition~\ref{prop.spinc}. Let $C_{\rm b}(\chi)$ be the set of the black faces from which $\chi$ is connected, $C_{\rm w}(\chi)$ the set of the white faces from which it is connected, and $C'_{\rm w}(\chi)$ the set of the white faces to which it is connected. For an edge $e$ of $D$ we write $r_{\rm b}(e)$ and $r_{\rm w}(e)$ for the black and the white face adjacent to $e$, and we put
\begin{align*}
   E(\chi)&=\bigl\{e:\ r_{\rm b}(e)\in C_{\rm b}(\chi)\ \text{ and }\ r_{\rm w}(e)\in C_{\rm w}(\chi)\bigr\},\\
   \bar E(\chi)&=\bigl\{e:\ r_{\rm b}(e)\in C_{\rm b}(\chi)\ \text{ and }\ r_{\rm w}(e)\in C'_{\rm w}(\chi)\bigr\}.
\end{align*}
By Proposition~\ref{prop.spinc}~(ii) and~(iii) and the first lines of their proof, applied to the root $e$,
\begin{equation}\label{eq.rootsofchi}
\begin{aligned}
   e\in E(\chi)\ &\iff\ \chi=\chi(\mathfrak s)\ \text{ for a state $\mathfrak s$ of $(D,e)$},\\
   e\in\bar E(\chi)\ &\iff\ \chi=\bar\chi(\mathfrak s)\ \text{ for a state $\mathfrak s$ of $(D,e)$}:
\end{aligned}
\end{equation}
the vectors of the states of $(D,e)$ have these connectivity properties, every class contains one of them, and a class contains only one vector of $\{\pm1\}^c$ with these properties. Thus $E(\chi)$ is the set of the roots for which $\chi$ is the representative $u_{\mathfrak t}$ of its class $\mathfrak t=[\chi]$, and $\bar E(\chi)$ the set of those for which it is $\bar u_{\mathfrak t}$. Each of them may be empty; and for every class $\mathfrak t$, every edge of $D$ belongs to $E(\chi)$ for exactly one vector $\chi\in\{\pm1\}^c$ of $\mathfrak t$, and to $\bar E(\chi)$ for exactly one.

\begin{lemma}[cyclic crossings]\label{lem.minty}
Let $\chi\in\{\pm1\}^c$.
\begin{enumerate}\renewcommand{\theenumi}{\roman{enumi}}
\item For every crossing $j$, exactly one of the following holds: the edge of $j$ in $\CG$ lies on a directed cycle of $\CG$; the edge of $j$ in $\CG^*$ lies on a directed cycle of $\CG^*$. In the first case we say that $j$ is \emph{cyclic} for $\chi$, and we denote by $A(\chi)$ the set of the cyclic crossings. A crossing which is not cyclic lies on a \emph{directed cut}: a set $\partial R$, formed by the edges of $\CG$ with exactly one end in a set $R$ of black faces, all of which point toward $R$. No crossing lies both on a directed cycle and on a directed cut of $\CG$.
\item If $\chi'\in\{\pm1\}^c$ is equivalent to $\chi$, then $\chi'$ is obtained from $\chi$ by reversing successively directed cycles and directed cuts, and $A(\chi')=A(\chi)$.
\end{enumerate}
\end{lemma}
\begin{proof}
(i) Let $\partial R$ be a directed cut. Follow the boundary of a face of $\CG$ with this face on the left; it is a closed walk of $\CG$, which enters and leaves $R$ alternately. When it enters $R$ it follows an edge of $\partial R$ in its direction, and the edge of this crossing in $\CG^*$, which is this edge turned counterclockwise, ends at the white face $f$ contained in the face of $\CG$; when it leaves $R$, the edge in $\CG^*$ starts at $f$. So in the subgraph of $\CG^*$ formed by the edges of the crossings of $\partial R$ every white face has as many incoming as outgoing edges, and this subgraph is a union of directed cycles: the edge in $\CG^*$ of every crossing of a directed cut of $\CG$ lies on a directed cycle of $\CG^*$. Now let the edge of $j$ in $\CG$ go from the black face $v$ to the black face $v'$, and let $R$ be the set of the black faces which are the ends of the directed paths starting at $v'$. If $v\notin R$, no edge leaves $R$, so $\partial R$ is a directed cut which contains the edge of $j$, and the edge of $j$ in $\CG^*$ lies on a directed cycle. If $v\in R$, the edge of $j$ lies on a directed cycle of $\CG$, hence on a simple one, $Z$, which is a simple closed curve on $S^2$. The edges in $\CG^*$ of the crossings of $Z$ all cross $Z$ from its right to its left; since a closed walk of $\CG^*$ crosses $Z$ as many times in each direction, none of these edges lies on a directed cycle of $\CG^*$. This proves the dichotomy, and the last sentence of (i) follows from it and from the first step.

(ii) Let $d=\frac12(\chi-\chi')=x+k$ with $x\in\mathrm{Flow}(\CG)$ and $k\in\mathrm{Cut}(\CG)$. The coordinates of $d$ are $0$ and $\pm1$, and $d_j=\chi_j=-\chi'_j$ if $d_j\neq0$. We argue by induction on the number of the nonzero coordinates of $d$. If $k\neq0$, the first step of the proof of the uniqueness in Proposition~\ref{prop.spinc}~(iii), which up to that point uses only that $\chi$ and $\chi'$ are equivalent, gives a set $S$ of black faces with $\partial S\neq\emptyset$ such that all the edges of $\partial S$ point toward $S$ for $\chi$ and away from $S$ for $\chi'$. So $\partial S$ is a directed cut of $\chi$, and if $\kappa\in\mathrm{Cut}(\CG)$ is the vector of the indicator function of $S$, then $\kappa_j=\chi_j=d_j$ for $j\in\partial S$, and $\chi-2\kappa$ is the vector $\chi$ with this directed cut reversed. If $k=0$ and $d\neq0$, the edges $j$ with $d_j\neq0$, oriented by $\chi$, form a subgraph in which every black face has as many incoming as outgoing edges, because $d$ is a flow; this subgraph contains a directed cycle of $\chi$, and if $\zeta\in\mathrm{Flow}(\CG)$ is the vector of this cycle, with $\zeta_j=\chi_j$ on the cycle and $\zeta_j=0$ elsewhere, then $\chi-2\zeta$ is $\chi$ with this directed cycle reversed. In both cases the new vector is equivalent to $\chi'$ and differs from it at fewer crossings. Finally, by (i), a directed cut does not meet a directed cycle; so the reversal of a directed cycle, which is undone by the reversal of a directed cycle, does not change the set of the directed cuts, nor the set of the crossings which lie on them, which is the complement of $A(\chi)$. In the same way the reversal of a directed cut does not change the directed cycles.
\end{proof}

Part (i) is Minty's lemma \cite{Minty66} for the plane graph $\CG$ and its dual, and by (ii) the classes of Proposition~\ref{prop.spinc}, restricted to $\{\pm1\}^c$, are the classes of the cycle--cocycle reversing system of Gioan \cite{Gioan07}. For a class $\mathfrak t$ we write $A(\mathfrak t)$ for the set $A(\chi)$ of any vector $\chi\in\{\pm1\}^c$ of $\mathfrak t$, and we call it the \emph{cyclic part} of $\mathfrak t$. It is determined by the lattice with its coordinates: $j\notin A(\mathfrak t)$ if and only if there is a vector $\kappa\in\mathrm{Cut}(\CG)$ with coordinates $0,\pm1$ such that $\kappa_j\neq0$ and $\chi-2\kappa\in\{\pm1\}^c$, for one, or for every, vector $\chi\in\{\pm1\}^c$ of $\mathfrak t$. Indeed, the vector of the indicator function of $R$ has these properties if $j$ lies on the directed cut $\partial R$; conversely, $\chi-2\kappa\in\{\pm1\}^c$ means that $\kappa_i=\chi_i$ whenever $\kappa_i\neq0$, so that $\kappa$ is the vector of a function $\psi$ on the black faces which increases by one along the edges $i$ with $\kappa_i\neq0$, oriented by $\chi$, and is constant along the other edges; if $n$ is its value at the origin of the edge of $j$, the edges between $\{\psi\le n\}$ and $\{\psi>n\}$ form a directed cut which contains this edge.

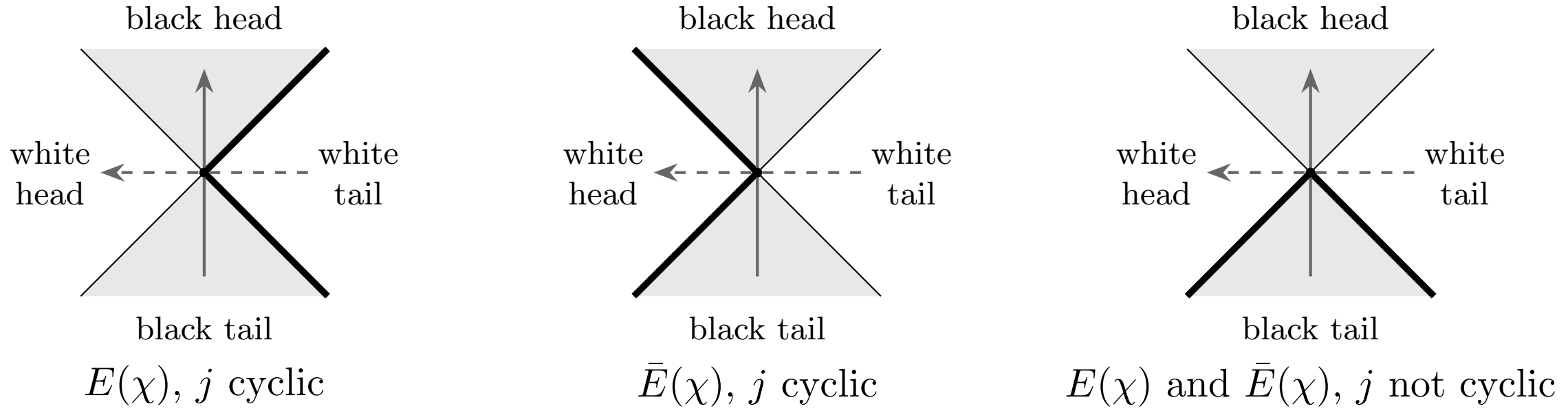
\begin{figure}[ht]
\centering
\begin{tikzpicture}[>=Stealth,scale=0.95]
\foreach \k/\lab/\ea/\eb in {0/{$E(\chi)$, $j$ cyclic}/{1.25,1.25}/{1.25,-1.25}, 1/{$\bar E(\chi)$, $j$ cyclic}/{-1.25,1.25}/{-1.25,-1.25}, 2/{$E(\chi)$ and $\bar E(\chi)$, $j$ not cyclic}/{-1.25,-1.25}/{1.25,-1.25}}{
\begin{scope}[xshift=\k*5.6cm]
\fill[black!9] (0,0) -- (1.25,1.25) -- (-1.25,1.25) -- cycle;
\fill[black!9] (0,0) -- (1.25,-1.25) -- (-1.25,-1.25) -- cycle;
\draw (0,0) -- (1.25,1.25); \draw (0,0) -- (1.25,-1.25); \draw (0,0) -- (-1.25,-1.25); \draw (0,0) -- (-1.25,1.25);
\draw[line width=1.8pt] (\ea) -- (0,0) -- (\eb);
\draw[->,thick,black!60] (0,-1.05) -- (0,1.05);
\draw[->,thick,black!60,dashed] (1.05,0) -- (-1.05,0);
\fill (0,0) circle (1.4pt);
\node[font=\scriptsize,above] at (0,1.3) {black head};
\node[font=\scriptsize,below] at (0,-1.3) {black tail};
\node[font=\scriptsize,left,align=center] at (-1.0,0) {white\\head};
\node[font=\scriptsize,right,align=center] at (1.0,0) {white\\tail};
\node[font=\small] at (0,-2.15) {\lab};
\end{scope}}
\end{tikzpicture}
\caption{A crossing $j$ with its edges in $\CG$ (solid arrow) and in $\CG^*$ (dashed arrow) for an orientation $\chi$, as in Figure~\ref{fig.heads}. In thick, the two ends of edges of $D$ at $j$ which belong to the curve of $\chi$, when this curve passes through $j$: it turns around the white tail (first flavor) or around the white head (second flavor) if $j$ is cyclic, and around the black tail if it is not.}
\label{fig.turn}
\end{figure}

\begin{proposition}[the curve of an orientation]\label{prop.curve}
Let $\chi\in\{\pm1\}^c$ with $E(\chi)\neq\emptyset$.
\begin{enumerate}\renewcommand{\theenumi}{\roman{enumi}}
\item At every crossing $j$, either none of the four ends of edges of $D$ at $j$ is an end of an edge of $E(\chi)$, or exactly two of them are: the two ends adjacent to the white tail of $j$ if $j$ is cyclic, and the two ends adjacent to the black tail of $j$ if it is not (Figure~\ref{fig.turn}).
\item The edges of $E(\chi)$ form one simple closed curve on $S^2$. The faces of $C_{\rm b}(\chi)$ lie on one side of it, and the faces of $C_{\rm w}(\chi)$ on the other side.
\end{enumerate}
The same holds for $\bar E(\chi)$ when it is not empty, with the white head in place of the white tail in \textup{(i)} and $C'_{\rm w}(\chi)$ in place of $C_{\rm w}(\chi)$ in \textup{(ii)}.
\end{proposition}

\begin{proof}
(i) Let $e\in E(\chi)$ have an end at $j$. This end separates a black corner and a white corner of $j$, which lie in $v=r_{\rm b}(e)\in C_{\rm b}(\chi)$ and in $f=r_{\rm w}(e)\in C_{\rm w}(\chi)$; let $v'$ and $f'$ be the faces of the other black corner and of the other white corner of $j$. Suppose that $j$ is cyclic. Then $v'\in C_{\rm b}(\chi)$, because $v$ is the end of a directed path which starts at $v'$. By Lemma~\ref{lem.minty}~(i) the edge of $j$ in $\CG^*$ does not lie on a directed cycle of $\CG^*$. It cannot go from $f'$ to $f$, since $f'$ is the end of a directed path which starts at $f$; so it goes from $f$ to $f'$, the corner in $f$ is the white tail of $j$, and $f'\notin C_{\rm w}(\chi)$, since otherwise $f$ would be the end of a directed path which starts at $f'$. The four ends at $j$ separate the corners in $(v,f)$, $(f,v')$, $(v',f')$ and $(f',v)$, so exactly the first two are ends of edges of $E(\chi)$: the two ends adjacent to the white tail. If $j$ is not cyclic, exchange the two colors: the edge of $j$ in $\CG^*$ lies on a directed cycle, so $f'\in C_{\rm w}(\chi)$; the edge of $j$ in $\CG$ goes from $v$ to $v'$, and $v'\notin C_{\rm b}(\chi)$; and the ends of edges of $E(\chi)$ are the two ends adjacent to the corner in $v$, which is the black tail. (If $v=v'$ the crossing is cyclic, and if $f=f'$ it is not; the argument covers these cases.)

(ii) By (i), the edges of $E(\chi)$ form pairwise disjoint simple closed curves $\gamma_1,\dots,\gamma_m$. Along each of them the white faces adjacent to its edges lie on one side and the black faces on the other: at a crossing where the curve turns around a white corner, its two edges there are adjacent to the face of this corner, which is on the inner side of the turn, and their black faces are on the outer side; and similarly for a black corner. Two faces of $C_{\rm b}(\chi)$ are joined by a path of $\CG$ whose edges are cyclic crossings between faces of $C_{\rm b}(\chi)$; if one of the curves passes through such a crossing, it turns around a white corner, and it does not separate the two black corners. So all the faces of $C_{\rm b}(\chi)$ lie in one component $U$ of the complement of $\gamma_1\cup\dots\cup\gamma_m$ in $S^2$, and in the same way all the faces of $C_{\rm w}(\chi)$ lie in one component $V$. Every curve $\gamma_i$ has a face of $C_{\rm b}(\chi)$ on one side and a face of $C_{\rm w}(\chi)$ on the other, so it lies on the boundary of both $U$ and $V$, which are on different sides of it. If $m\ge2$, the curve $\gamma_2$ is contained in one of the two disks bounded by $\gamma_1$, and the one of $U$, $V$ which lies in the other disk does not have $\gamma_2$ on its boundary. So $m=1$.

For $\bar E(\chi)$ the argument is the same: a face of $C'_{\rm w}(\chi)$ is the end of directed paths starting at all the white faces, so that at a cyclic crossing the edge of $j$ in $\CG^*$ goes from $f'$ to $f$, and the corner in $f$ is the white head.
\end{proof}

\begin{corollary}[the smoothing of a class]\label{cor.classsmoothing}
Let $\mathfrak t$ be a class and $A=A(\mathfrak t)$ its cyclic part. Let $s_A(D)$ be the union of disjoint circles obtained from $D$ by replacing every crossing $j$ by two arcs, which join the two ends of edges adjacent to each white corner of $j$ if $j\in A$, and the two ends adjacent to each black corner of $j$ if $j\notin A$.
\begin{enumerate}\renewcommand{\theenumi}{\roman{enumi}}
\item The curves $E(\chi)$ of the vectors $\chi\in\{\pm1\}^c$ of $\mathfrak t$ with $E(\chi)\neq\emptyset$ are the circles of $s_A(D)$: every edge of $D$ lies on exactly one of them, and the two which pass through a crossing are different. The same holds for the curves $\bar E(\chi)$. Consequently the representatives $u_{\mathfrak t}$ and $\bar u_{\mathfrak t}$ for a root $\hate$ depend on $\hate$ only through the circle of $s_A(D)$ on which $\hate$ lies.
\item The number of these circles is $k+b_1$, where $k$ is the number of connected components and $b_1$ the first Betti number of the graph $\CG_A$ formed by all the black faces and the edges of the crossings of $A$.
\item If $\mathfrak t=[\1]$ is the class of the over-orientation, then $A$ is the set of the positive crossings, $s_A(D)$ is the smoothing of $D$ along its orientation, and its circles are the Seifert circles of $D$.
\end{enumerate}
\end{corollary}
\begin{proof}
(i) By \eqref{eq.rootsofchi} and Proposition~\ref{prop.spinc}~(ii), every edge of $D$ lies on $E(\chi)$ for exactly one vector $\chi\in\{\pm1\}^c$ of $\mathfrak t$. By Lemma~\ref{lem.minty}~(ii) these vectors have the same cyclic part $A$, so by Proposition~\ref{prop.curve} each of the curves $E(\chi)$ is a circle of $s_A(D)$, all the circles are obtained, and a curve which passes through a crossing contains only two of the four ends at this crossing. The same applies to $\bar E(\chi)$, since the two arcs of $s_A(D)$ at a crossing of $A$ turn around its two white corners, its white head and its white tail.

(ii) The circles of $s_A(D)$ are the boundary components of a regular neighborhood in $S^2$ of the plane graph $\CG_A$: the black faces, joined by a band at every crossing of $A$. A connected plane graph with first Betti number $b$ has $b+1$ faces, so this boundary has $k+b_1$ components.

(iii) By Proposition~\ref{prop.omega}, $\1_+$ is a flow with coordinates $0$ and $1$; so for the orientation $\1$ the edges of the positive crossings form a subgraph with as many incoming as outgoing edges at every black face, and each of them lies on a directed cycle. By the same proposition $\1_-$ is the vector of an integer-valued function $\psi$ on the black faces; so the edge of a negative crossing goes from a face with $\psi=n$ to a face with $\psi=n+1$, the edges of the positive crossings join faces with the same value of $\psi$, and every negative crossing lies on the directed cut between $\{\psi\le n\}$ and $\{\psi>n\}$. Hence $A$ is the set of the positive crossings, by Lemma~\ref{lem.minty}~(i). The white corners of a positive crossing and the black corners of a negative crossing are its $A$-corner and its $B$-corner; so at every crossing $s_A(D)$ joins $o_j$ to $u'_j$ and $u_j$ to $o'_j$, by \eqref{eq.corneredges}, and its circles are the Seifert circles.
\end{proof}

So every $\mathrm{Spin}^c$ structure of $\Sigma(L)$ comes with a smoothing of the diagram, that is, with a state of the Kauffman bracket \cite{Kauffman87}, and the vectors of the class which occur as representatives are in bijection with the circles of this smoothing, in two ways, one for each flavor; the other vectors of the class are never the orientation of a state. For the class of the over-orientation, whose covector on $\mathrm{Flow}(\CG)$ comes from an element of the lattice (Remark~\ref{rem.branched}~(iii)), the smoothing is the one which gives the Seifert circles. The curves of Proposition~\ref{prop.curve} also make explicit the face potential of Proposition~\ref{prop.potential}, which becomes a sum of indicator functions of disks with coefficients of constant sign. In particular they give a form of the state sum without a root.

\begin{corollary}[the face potential of an alternating diagram]\label{cor.altpotential}
For $\chi\in\{\pm1\}^c$ with $E(\chi)\neq\emptyset$ let $\CR(\chi)$ be the set of the faces of $D$ which lie on the same side of the curve $E(\chi)$ as the faces of $C_{\rm b}(\chi)$, and let $S(\chi)=\{j:\chi_j=\ve_j\}$. If $\bar E(\chi)\neq\emptyset$ let $\bar\CR(\chi)$ be the set of the faces which lie on the same side of $\bar E(\chi)$ as the faces of $C_{\rm b}(\chi)$, and let $Y(\chi)=\{j:\chi_j=-1\}$. Let
\[
   \mathcal K(x)=\sum_{\chi:\,E(\chi)\neq\emptyset}x^{S(\chi)}\,\1_{\CR(\chi)}\in\BZ[x]^{\CF},\qquad\qquad \widetilde{\mathcal K}(y)=\sum_{\chi:\,\bar E(\chi)\neq\emptyset}y^{Y(\chi)}\,\1_{\bar\CR(\chi)}\in\BZ[y]^{\CF},
\]
where $\1_{\CR}$ is the indicator vector of a set $\CR$ of faces.
\begin{enumerate}\renewcommand{\theenumi}{\roman{enumi}}
\item For an edge $e$ of $D$ let $s(e)=1$ if the black face adjacent to $e$ is on its left and $s(e)=-1$ if it is on its right, and let $\sigma_{D,e}$ be the sign of Theorem~\ref{thm.kauffmanalt}~\textup{(i)}. Then the sign $\theta=s(e)\,\sigma_{D,e}$ does not depend on the edge $e$, and the vectors $k(x)=\theta\,\mathcal K(-x)$ and $\tilde k(y)=\theta\,\widetilde{\mathcal K}(-y)$ satisfy, for every edge $e$ of $D$,
\[
   P_{D,e}(x)=k_{f_{\rm L}(e)}(x)-k_{f_{\rm R}(e)}(x)\qquad\text{and}\qquad \tP_{D,e}(y)=\tilde k_{f_{\rm L}(e)}(y)-\tilde k_{f_{\rm R}(e)}(y).
\]
Thus $k(x)$ is the face potential of Proposition~\ref{prop.potential}, up to a multiple of $\1$.
\item The Links--Gould polynomial is the sum over the ordered pairs of orientations of the Tait graph
\[
   2c\,(-1)^{\mu-1}\,\LG(L;-t_0,-t_1)=\sum_{(\chi,\chi')}\bigl|E(\chi)\cap\bar E(\chi')\bigr|\ \prod_{j=1}^cm_j(\chi_j,\chi'_j),
\]
where $m_j$ is the local weight of Theorem~\ref{thm.tait}~(iii) and $|E(\chi)\cap\bar E(\chi')|$ is the number of edges of $D$ common to the two curves.
\end{enumerate}
\end{corollary}
\begin{proof}
(i) The two faces adjacent to an edge $e$ are separated by the curve $E(\chi)$ if and only if $e$ lies on it, and then $r_{\rm b}(e)\in C_{\rm b}(\chi)\subseteq \CR(\chi)$ and $r_{\rm w}(e)\notin \CR(\chi)$, by Proposition~\ref{prop.curve}~(ii). So $\mathcal K_{r_{\rm b}(e)}-\mathcal K_{r_{\rm w}(e)}$ is the sum of $x^{S(\chi)}$ over the vectors $\chi$ with $e\in E(\chi)$. By \eqref{eq.rootsofchi} these are the vectors $\chi(\mathfrak s)$ of the states of $(D,e)$, and $S(\chi(\mathfrak s))=X(\mathfrak s)$ by Theorem~\ref{thm.tait}~(ii); hence $\mathcal K_{r_{\rm b}(e)}(x)-\mathcal K_{r_{\rm w}(e)}(x)=\sum_{\mathfrak s}x^{X(\mathfrak s)}$, and in the same way $\widetilde{\mathcal K}_{r_{\rm b}(e)}(y)-\widetilde{\mathcal K}_{r_{\rm w}(e)}(y)=\sum_{\mathfrak s}y^{Y(\mathfrak s)}$, the sums running over the states of $(D,e)$. By Theorem~\ref{thm.kauffmanalt}~(i), $P_{D,e}(x)=\sigma_{D,e}\sum_{\mathfrak s}(-1)^{|X(\mathfrak s)|}x^{X(\mathfrak s)}$ and $\tP_{D,e}(y)=\sigma_{D,e}\sum_{\mathfrak s}(-1)^{|Y(\mathfrak s)|}y^{Y(\mathfrak s)}$. At the point $x_j=-1$ for all $j$, which is the Alexander point at $t=-1$, this gives $P_{D,e}=\sigma_{D,e}\,|\Delta_L(-1)|\neq0$. On the other hand, by the argument of the proof of Lemma~\ref{lem.understrand} applied to the matrix $A_D(-1)$, the vector $(P_{D,e})_e$ at this point is proportional to the row vector which spans the left kernel of $A_D(-1)$, and this vector is $\bigl((-1)^{h(e)}\bigr)_e$ by Lemma~\ref{lem.level}~(ii). So $\sigma_{D,e}(-1)^{h(e)}$ does not depend on $e$. Now $h(e)=-w(f_{\rm L}(e))$ by \eqref{eq.wh}, and the winding function $w$ has one parity on the black faces (Proposition~\ref{prop.omega}) and the other one on the white faces (Lemma~\ref{lem.winding}). Hence $(-1)^{h(e)}s(e)$ does not depend on $e$ either, so that $\theta=s(e)\,\sigma_{D,e}$ does not depend on $e$, and
\[
   k_{f_{\rm L}(e)}(x)-k_{f_{\rm R}(e)}(x)=\theta\,s(e)\bigl(\mathcal K_{r_{\rm b}(e)}(-x)-\mathcal K_{r_{\rm w}(e)}(-x)\bigr)=\theta\,s(e)\,\sigma_{D,e}\,P_{D,e}(x)=P_{D,e}(x).
\]
The same computation applies to the second flavor, with the same signs.

(ii) By \eqref{eq.rootsofchi}, for every edge $e$ the right-hand side of \eqref{eq.taitsum} for the root $e$ is the sum of the products $\prod_jm_j(\chi_j,\chi'_j)$ over the pairs $(\chi,\chi')$ such that $e\in E(\chi)\cap\bar E(\chi')$. Summing over the $2c$ edges gives the formula.
\end{proof}

We now turn to the terms of \eqref{eq.spincsum}. For $\chi\in\{\pm1\}^c$ let $\chi^\natural\in\{\pm1\}^c$ be the vector which is equal to $-\chi$ on the cyclic part $A(\chi)$ and to $\chi$ outside it: the orientation $\chi$ with all its directed cycles reversed.

\begin{proposition}[the terms which do not depend on the root]\label{prop.nat}
Let $\chi\in\{\pm1\}^c$, with class $\mathfrak t=[\chi]$ and cyclic part $A=A(\mathfrak t)$, and let $J(\mathfrak t)$ be the set of the positive crossings which are not in $A$ and of the negative crossings which are in $A$.
\begin{enumerate}\renewcommand{\theenumi}{\roman{enumi}}
\item $A(\chi^\natural)=A$, and the class of $\chi^\natural$ depends only on $\mathfrak t$; we denote it by $\mathfrak t^\natural$. The map $\mathfrak t\mapsto\mathfrak t^\natural$ is an involution of the set of the classes.
\item $\bar E(\chi^\natural)=E(\chi)$. Hence, for every root $\hate$, the second representative of $\mathfrak t^\natural$ is $\bar u_{\mathfrak t^\natural}=(u_{\mathfrak t})^\natural$.
\item The monomial of the pair $(\chi,\chi^\natural)$ is
\[
   \prod_{j=1}^cm_j(\chi_j,\chi^\natural_j)=t_0^{\,a}\,t_1^{\,b}\qquad\text{with}\qquad a+b=\ve(J(\mathfrak t))=c_+-|A|,\qquad b-a=-\sum_{j\in J(\mathfrak t)}\chi_j ,
\]
and it is the same for all the vectors $\chi\in\{\pm1\}^c$ of $\mathfrak t$; we denote it by $M(\mathfrak t)$.
\end{enumerate}
Consequently, for every root, the term of the pair of classes $(\mathfrak t,\mathfrak t^\natural)$ in \eqref{eq.spincsum} is $M(\mathfrak t)$, which in the variables $p,q$ of Lemma~\ref{lem.pq} is $p^{\,b-a}q^{\,c_+-|A(\mathfrak t)|}$. These are $|\Delta_L(-1)|$ of the $\Delta_L(-1)^2$ terms of \eqref{eq.spincsum}, and they are determined by the lattice $\BZ^c\supset\mathrm{Flow}(\CG)\oplus\mathrm{Cut}(\CG)$ with its coordinates and by the signs of the crossings.
\end{proposition}
\begin{proof}
(i) A directed cycle of $\chi$ consists of cyclic crossings, and it is a directed cycle of $\chi^\natural$, in the opposite direction; a directed cut of $\chi$ contains no crossing of $A$, by Lemma~\ref{lem.minty}~(i), and it is a directed cut of $\chi^\natural$. So for $\chi^\natural$ the crossings of $A$ lie on directed cycles and the other ones on directed cuts, and $A(\chi^\natural)=A$ by Lemma~\ref{lem.minty}~(i); in particular $(\chi^\natural)^\natural=\chi$. If $\chi_1$ is obtained from $\chi$ by reversing a directed cycle, with vector $\zeta$, then $A(\chi_1)=A$ and $\chi_1^\natural=\chi^\natural+2\zeta$ is obtained from $\chi^\natural$ by reversing its directed cycle $-\zeta$; if $\chi_1=\chi-2\kappa$ is obtained by reversing a directed cut, with vector $\kappa$, then $\chi_1^\natural=\chi^\natural-2\kappa$. So $\chi_1^\natural$ is equivalent to $\chi^\natural$ in both cases, and by Lemma~\ref{lem.minty}~(ii) the class of $\chi^\natural$ depends only on that of $\chi$.

(ii) Let $v\in C_{\rm b}(\chi)$ and let $v'$ be a black face. A directed path of $\chi$ from $v$ to $v'$ becomes a directed path of $\chi^\natural$ if each of its cyclic edges, from $v_1$ to $v_2$ say, is replaced by the rest of a directed cycle of $\chi$ through this edge, which goes from $v_2$ to $v_1$ for $\chi$, hence from $v_1$ to $v_2$ for $\chi^\natural$. So $C_{\rm b}(\chi)\subseteq C_{\rm b}(\chi^\natural)$, and the two sets are equal, because $(\chi^\natural)^\natural=\chi$. By Lemma~\ref{lem.minty}~(i), the vector $-\chi^\natural$ is equal to $-\chi$ on the crossings whose edge in $\CG^*$ lies on a directed cycle of $\CG^*$, and to $\chi$ on the other ones; so the same argument, applied to $\CG^*$, gives $C_{\rm w}(\chi)=C_{\rm w}(-\chi^\natural)$, and $C_{\rm w}(-\chi^\natural)=C'_{\rm w}(\chi^\natural)$, since reversing all the edges exchanges the origins and the ends of the directed paths. Hence $\bar E(\chi^\natural)=E(\chi)$. If $u=u_{\mathfrak t}$ for the root $\hate$, then $\hate\in E(u)=\bar E(u^\natural)$, and $u^\natural$ is the second representative of its class $\mathfrak t^\natural$.

(iii) By the table of Theorem~\ref{thm.tait}~(iii), a crossing outside $A$, at which the two vectors are equal, contributes $t_0$ or $t_1$ if it is positive, according to whether $\chi_j=1$ or $-1$, and $1$ if it is negative; a crossing of $A$, at which they are opposite, contributes $1$ if it is positive, and $t_1^{-1}$ or $t_0^{-1}$ if it is negative, according to whether $\chi_j=1$ or $-1$. This gives $a+b$ and $b-a$, and $\ve(J(\mathfrak t))=c_+-|A|$. Let $\1_J\in\BZ^c$ be the indicator vector of $J(\mathfrak t)$. If $\zeta$ is the vector of a directed cycle of $\chi$, it is supported in $A$ and $\langle\zeta,\1_J\rangle=\langle\zeta,\1_-\rangle=0$, because $\1_-$ is a cut (Proposition~\ref{prop.omega}); if $\kappa$ is the vector of a directed cut, it is supported outside $A$ and $\langle\kappa,\1_J\rangle=\langle\kappa,\1_+\rangle=0$, because $\1_+$ is a flow. So $b-a=-\langle\chi,\1_J\rangle$ does not change under the reversals of Lemma~\ref{lem.minty}~(ii).

The last statements follow from (ii) and (iii): the pairs $(\mathfrak t,\mathfrak t^\natural)$ are as many as the classes; and $\mathfrak t^\natural$ and $M(\mathfrak t)$ are defined by a vector $\chi\in\{\pm1\}^c$ of $\mathfrak t$, by the signs of the crossings and by the cyclic part $A(\mathfrak t)$, which was described in the lattice after Lemma~\ref{lem.minty}.
\end{proof}

In the next theorem an orientation is acyclic if it has no directed cycle, and a source is a vertex at which no edge ends (Definition~\ref{def.directed}); every acyclic orientation has at least one source. The \emph{Tutte polynomial} $T_{\CG}(x,y)$ of a graph $\CG$ is a polynomial in two variables with nonnegative integer coefficients which contains many numerical invariants of the graph; for instance $T_{\CG}(1,1)$ is the number of its spanning trees. We only use two of its values: by a theorem of Greene and Zaslavsky \cite{GreeneZaslavsky83}, $T_{\CG}(1,0)$ is the number of the acyclic orientations of a connected graph $\CG$ which have a given vertex as their only source, and for a connected plane graph $T_{\CG}(0,1)=T_{\CG^*}(1,0)$.

\begin{theorem}[the two extreme $q$-degrees]\label{thm.extreme}
Let $D$ be a diagram of a link $L$, with $c_+$ positive and $c_-$ negative crossings.
\begin{enumerate}\renewcommand{\theenumi}{\roman{enumi}}
\item Every monomial $t_0^at_1^b$ of $\LG(L)$ satisfies $-c_-\le a+b\le c_+$: in the variables $p,q$ of Lemma~\ref{lem.pq}, the degrees in $q$ of the monomials of $\LG(L)$ lie in $[-c_-,c_+]$.
\item Let $D$ be alternating, and let $v_0$ be a black face and $f_0$ a white face of $D$. The part of degree $c_+$ in $q$ of $(-1)^{\mu-1}\LG(L;-t_0,-t_1)$ is
\[
   \sum_u\ t_0^{\,n(u)}\,t_1^{\,c_+-n(u)}=q^{\,c_+}\sum_u\ p^{\,c_+-2n(u)},
\]
the sum running over the acyclic orientations $u$ of $\CG$ whose only source is $v_0$, where $n(u)$ is the number of positive crossings $j$ with $u_j=1$. It is the sum of the monomials $M(\mathfrak t)$ of the classes with $A(\mathfrak t)=\emptyset$. The part of degree $-c_-$ in $q$ is
\[
   \sum_u\ t_0^{\,-(c_--n'(u))}\,t_1^{\,-n'(u)}=q^{\,-c_-}\sum_u\ p^{\,c_--2n'(u)},
\]
the sum running over the orientations $u$ of $\CG$ for which the orientation of $\CG^*$ is acyclic with $f_0$ as its only source, where $n'(u)$ is the number of negative crossings $j$ with $u_j=1$. It is the sum of the monomials $M(\mathfrak t)$ of the classes whose cyclic part is the set of all the crossings. The numbers of terms of these two sums are the values $T_{\CG}(1,0)$ and $T_{\CG}(0,1)$ of the Tutte polynomial of $\CG$.
\item If $D$ is a reduced alternating diagram, these two parts are not zero: $c_+$ is the largest and $-c_-$ the smallest degree in $q$ of a monomial of $\LG(L)$.
\end{enumerate}
\end{theorem}
\begin{proof}
(i) By Lemma~\ref{lem.pq}, $\LG(L)=(-1)^{\mu-1}\sum_{S,T}P_S\tP_T\,\BE_D[x^Sy^T]$, and $\BE_D[x^Sy^T]$ is a monomial of degree $\ve(S\triangle T)$ in $q$; and $-c_-\le\ve(S\triangle T)\le c_+$ for all sets of crossings $S$ and $T$.

(ii) Let $\hate$ be an edge of $D$ with $r_{\rm b}(\hate)=v_0$. By Theorem~\ref{thm.tait}~(iii) the monomial of a pair of states $(\mathfrak s,\mathfrak s')$ of $(D,\hate)$ has the degree $\frac12\bigl(\mathrm{wr}(D)+\langle\chi(\mathfrak s),\bar\chi(\mathfrak s')\rangle\bigr)\le\frac12(\mathrm{wr}(D)+c)=c_+$ in $q$, with equality if and only if $\chi(\mathfrak s)=\bar\chi(\mathfrak s')$. So, by \eqref{eq.rootsofchi}, the pairs of degree $c_+$ are given by the vectors $u$ with $\hate\in E(u)\cap\bar E(u)$, that is, with $v_0\in C_{\rm b}(u)$ and $r_{\rm w}(\hate)\in C_{\rm w}(u)\cap C'_{\rm w}(u)$; and the monomial of such a pair is $t_0^{n(u)}t_1^{c_+-n(u)}$, by the table of Theorem~\ref{thm.tait}~(iii). If $C_{\rm w}(u)\cap C'_{\rm w}(u)\neq\emptyset$, any two white faces are joined by a directed path of $\CG^*$, every edge of $\CG^*$ lies on a directed cycle, and by Lemma~\ref{lem.minty}~(i) the orientation $u$ of $\CG$ is acyclic; since every black face is the end of a directed path starting at $v_0$, this face is its only source. Conversely, if $u$ is acyclic with $v_0$ as its only source, then every black face is the end of a directed path starting at $v_0$; every edge of $\CG^*$ lies on a directed cycle, by Lemma~\ref{lem.minty}~(i); and since $\CG^*$ is connected, any two white faces are joined by a directed path, so that $C_{\rm w}(u)$ and $C'_{\rm w}(u)$ contain all the white faces. These orientations are the vectors of $\{\pm1\}^c$ with empty cyclic part which are connected from $v_0$; they are the representatives $u_{\mathfrak t}$, for the root $\hate$, of the classes with $A(\mathfrak t)=\emptyset$, for which $u^\natural=u$ and $\mathfrak t^\natural=\mathfrak t$, so that the monomial of $u$ is $M(\mathfrak t)$ by Proposition~\ref{prop.nat}. The number of the acyclic orientations of a connected graph with a given vertex as their only source is the value $T(1,0)$ of its Tutte polynomial, by a theorem of Greene and Zaslavsky \cite{GreeneZaslavsky83}.

The other end is the same statement with the two colors exchanged. The degree in $q$ is at least $\frac12(\mathrm{wr}(D)-c)=-c_-$, with equality if and only if $\chi(\mathfrak s)=-\bar\chi(\mathfrak s')$. A vector $u$ satisfies $\hate\in E(u)\cap\bar E(-u)$ if and only if $r_{\rm b}(\hate)$ is the origin and the end of directed paths of $\CG$ to and from all the black faces, and $r_{\rm w}(\hate)\in C_{\rm w}(u)$; as above, this means that the orientation of $\CG^*$ is acyclic with $r_{\rm w}(\hate)$ as its only source. The monomial of the pair $(u,-u)$ is $t_0^{-(c_--n'(u))}t_1^{-n'(u)}$. These vectors are the representatives of the classes whose cyclic part is the set of all the crossings, for which $u^\natural=-u$; their number is $T_{\CG^*}(1,0)=T_{\CG}(0,1)$; and we take for $\hate$ an edge with $r_{\rm w}(\hate)=f_0$.

(iii) The numbers $T_{\CG}(1,0)$ and $T_{\CG}(0,1)$ are not zero if and only if $\CG$, resp.\ $\CG^*$, has no loop: a connected graph without loops has an acyclic orientation whose only source is a given vertex. A loop of $\CG$ or of $\CG^*$ is a crossing with two opposite corners in the same face, that is, a nugatory crossing. Since the coefficients of $(-1)^{\mu-1}\LG(L;-t_0,-t_1)$ are nonnegative, the two parts are not zero.
\end{proof}

\begin{corollary}[the crossings of an alternating link]\label{cor.crossings}
Let $D$ be a reduced alternating diagram of a link $L$, with $c_+(D)$ positive and $c_-(D)$ negative crossings. Then every diagram $D'$ of $L$ satisfies
\[
   c_+(D')\ \ge\ c_+(D)\qquad\text{and}\qquad c_-(D')\ \ge\ c_-(D).
\]
In particular $D$ has the smallest number of crossings among the diagrams of $L$; all the reduced alternating diagrams of $L$ have the same number of crossings and the same writhe; and this writhe is zero if $L$ is isotopic to its mirror image.
\end{corollary}
\begin{proof}
By Theorem~\ref{thm.extreme}~(iii) and~(i), $c_+(D)$ is the largest degree in $q$ of the monomials of $\LG(L)$, which is at most $c_+(D')$, and likewise $-c_-(D)\ge-c_-(D')$. Here $D'$ has to be a diagram in the sense of Section~\ref{ssec.diagrams}, connected and with at least one crossing. Every diagram of $L$ is connected, because $L$ is not split: $\LG(L)\neq0$ by Theorem~\ref{thm.kauffmanalt}, the connected graph $\CG$ having at least one spanning tree, whereas $\LG$ vanishes on split links \cite[Lemma~3.2]{GHKKST}. And $L$ is not the unknot, whose polynomial is $1$: the largest degree in $q$ of $\LG(L)$ is $c_+(D)$, the smallest one is $-c_-(D)$, and $D$ has a crossing. Finally, the mirror image of $D$ is a reduced alternating diagram of the mirror image of $L$, with $c_-(D)$ positive and $c_+(D)$ negative crossings; so $c_+(D)=c_-(D)$ if $L$ is isotopic to its mirror image.
\end{proof}

The last statements are the first two conjectures of Tait on alternating links: a reduced alternating diagram has the smallest number of crossings, and its writhe is an invariant of the link, which is zero if the link is amphichiral, that is, isotopic to its mirror image. They were proved by Kauffman, Murasugi and Thistlethwaite with the Jones polynomial, for the number of crossings \cite{Kauffman87, Murasugi87, Thistlethwaite87}, and, for the writhe, by Thistlethwaite with the Kauffman polynomial \cite{Thistlethwaite88} and by Murasugi \cite{Murasugi87b}. In the proof above the Links--Gould polynomial plays the role of the Jones polynomial: the bound (i) is a property of the determinant formula, for every diagram, and its sharpness for the reduced alternating diagrams comes from the positivity of the state sum. The degree in $q$ sees the positive and the negative crossings separately, which is why the writhe comes with the number of crossings. With Theorems~\ref{thm.tight} and~\ref{thm.alternating}, the theorem places the support of $\LG$ of a link with a reduced alternating diagram $D$ in a rectangle of the plane of the exponents of $p$ and $q$, all four sides of which are attained: the degree in $p$ runs over $[-\nu,\nu]$ with $\nu=2g(L)+\mu-1$, and the degree in $q$ over $[-c_-(D),c_+(D)]$. For a diagram which is not alternating the bound (i) need not be attained: for the closure of the braid $(\sigma_1\sigma_2)^4$, a diagram with eight positive crossings of the torus knot $T(3,4)$, which is the knot $8_{19}$ up to mirror image and has crossing number eight, the degrees in $q$ of $\LG$ run from $0$ to $7$.

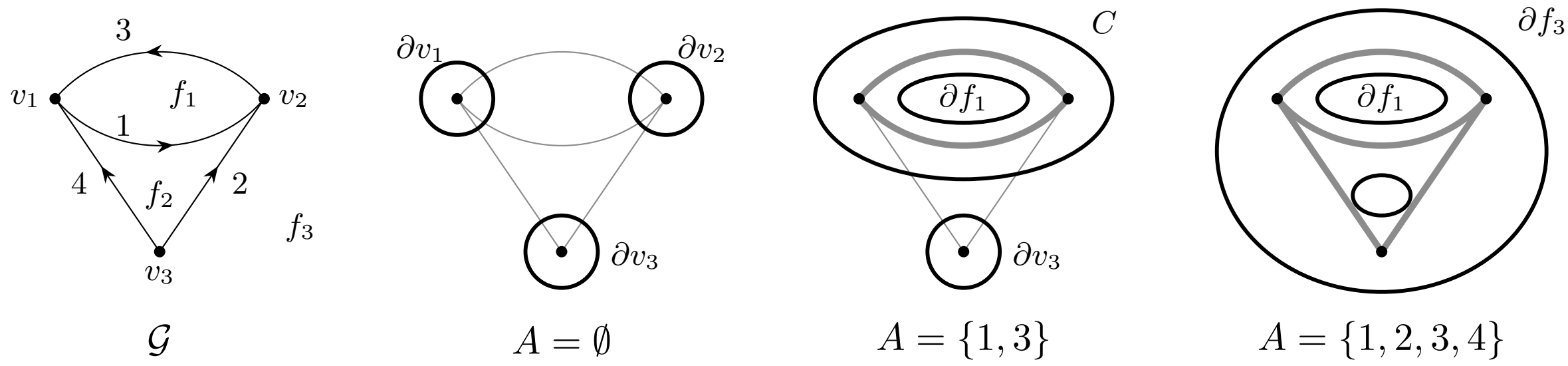
\begin{figure}[ht]
\centering
\begin{tikzpicture}[>=Stealth,scale=0.8,
   mid arrow/.style={postaction={decorate,decoration={markings,mark=at position 0.56 with {\arrow{Stealth}}}}}]
\def\graphpts{\coordinate (a) at (-1.3,0); \coordinate (b) at (1.3,0); \coordinate (c) at (0,-1.9);}
\begin{scope}
\graphpts
\draw[mid arrow] (a) to[bend right=50] (b); \draw[mid arrow] (b) to[bend right=50] (a);
\draw[mid arrow] (c) -- (b); \draw[mid arrow] (c) -- (a);
\fill (a) circle (2pt) (b) circle (2pt) (c) circle (2pt);
\node[font=\scriptsize,left] at (a) {$v_1$}; \node[font=\scriptsize,right] at (b) {$v_2$}; \node[font=\scriptsize,below] at (c) {$v_3$};
\node[font=\scriptsize] at (-0.45,-0.33) {$1$}; \node[font=\scriptsize] at (-0.45,0.86) {$3$};
\node[font=\scriptsize] at (1.0,-1.05) {$2$}; \node[font=\scriptsize] at (-1.0,-1.05) {$4$};
\node[font=\scriptsize] at (0.3,0.05) {$f_1$}; \node[font=\scriptsize] at (0,-1.2) {$f_2$}; \node[font=\scriptsize] at (1.75,-1.6) {$f_3$};
\node[font=\small] at (0,-3.0) {$\CG$};
\end{scope}
\begin{scope}[xshift=5.0cm]
\graphpts
\draw[black!45] (a) to[bend right=50] (b); \draw[black!45] (b) to[bend right=50] (a); \draw[black!45] (c) -- (b); \draw[black!45] (c) -- (a);
\fill (a) circle (2pt) (b) circle (2pt) (c) circle (2pt);
\draw[line width=1.1pt] (a) circle (0.45); \draw[line width=1.1pt] (b) circle (0.45); \draw[line width=1.1pt] (c) circle (0.45);
\node[font=\scriptsize] at (-1.75,0.62) {$\partial v_1$}; \node[font=\scriptsize] at (1.75,0.62) {$\partial v_2$}; \node[font=\scriptsize] at (0.92,-1.95) {$\partial v_3$};
\node[font=\small] at (0,-3.0) {$A=\emptyset$};
\end{scope}
\begin{scope}[xshift=10.0cm]
\graphpts
\draw[line width=1.8pt,black!45] (a) to[bend right=50] (b); \draw[line width=1.8pt,black!45] (b) to[bend right=50] (a); \draw[black!45] (c) -- (b); \draw[black!45] (c) -- (a);
\fill (a) circle (2pt) (b) circle (2pt) (c) circle (2pt);
\draw[line width=1.1pt] (0,0) ellipse (0.8 and 0.3);
\draw[line width=1.1pt] (0,0) ellipse (1.85 and 1.0);
\draw[line width=1.1pt] (c) circle (0.45);
\node[font=\scriptsize] at (0,0) {$\partial f_1$};
\node[font=\scriptsize] at (1.75,0.95) {$C$}; \node[font=\scriptsize] at (0.92,-1.95) {$\partial v_3$};
\node[font=\small] at (0,-3.0) {$A=\{1,3\}$};
\end{scope}
\begin{scope}[xshift=15.2cm]
\graphpts
\draw[line width=1.8pt,black!45] (a) to[bend right=50] (b); \draw[line width=1.8pt,black!45] (b) to[bend right=50] (a); \draw[line width=1.8pt,black!45] (c) -- (b); \draw[line width=1.8pt,black!45] (c) -- (a);
\fill (a) circle (2pt) (b) circle (2pt) (c) circle (2pt);
\draw[line width=1.1pt] (0,0) ellipse (0.8 and 0.3);
\draw[line width=1.1pt] (0,-1.2) ellipse (0.36 and 0.25);
\draw[line width=1.1pt] (0,-0.65) ellipse (2.05 and 1.75);
\node[font=\scriptsize] at (0,0) {$\partial f_1$};
\node[font=\scriptsize] at (2.0,0.95) {$\partial f_3$};
\node[font=\small] at (0,-3.0) {$A=\{1,2,3,4\}$};
\end{scope}
\end{tikzpicture}
\caption{The Tait graph of the diagram of the figure-eight knot of Example~\ref{ex.fig8classes}, with its over-orientation, and the three smoothings of the diagram which come from its five classes, drawn on the graph: the crossings of $A$ are the thick edges, and the circles of the smoothing run along the thick edges and cross the thin ones. In the last picture the small circle inside the triangle is $\partial f_2$.}
\label{fig.fig8smoothings}
\end{figure}

\begin{example}[the figure-eight knot]\label{ex.fig8classes}
Let $D$ be the closure of $\sigma_1\sigma_2^{-1}\sigma_1\sigma_2^{-1}$, as in Examples~\ref{ex.states} and~\ref{ex.tait}. Its Tait graph has three black faces $v_1,v_2,v_3$ and four edges: those of the positive crossings $1$ and $3$ go from $v_1$ to $v_2$ and from $v_2$ to $v_1$, and those of the negative crossings $2$ and $4$ go from $v_3$ to $v_2$ and to $v_1$ (Figure~\ref{fig.fig8smoothings}); its white faces are $f_1$, between the edges $1$ and $3$, $f_2$, bounded by the edges $1,2,4$, and $f_3$, bounded by the edges $2,3,4$. There are five classes, with three different cyclic parts, hence three smoothings. In the following table a circle of a smoothing is denoted by the vertex or the face of $\CG$ around which it runs, except the Seifert circle $C$, which passes through the four crossings; a vector of $\{\pm1\}^4$ is written as the sequence of the signs of its coordinates; and after each circle we give the representative $u_{\mathfrak t}$ of the class for the roots on this circle, followed after a bar by $\bar u_{\mathfrak t}$ when it is different.
{\footnotesize
\[
\begin{array}{c|c|l|c}
   \mathfrak t & A(\mathfrak t) & \text{circles and representatives} & M(\mathfrak t)\\ \hline
   \mathfrak t_1 & \emptyset & \partial v_1:\ {+}{+}{-}{-}\quad\ \partial v_2:\ {-}{-}{+}{-}\quad\ \partial v_3:\ {-}{+}{+}{+} & t_0t_1\\
   \mathfrak t_2 & \emptyset & \partial v_1:\ {+}{-}{-}{-}\quad\ \partial v_2:\ {-}{-}{+}{+}\quad\ \partial v_3:\ {+}{+}{-}{+} & t_0t_1\\
   \mathfrak t_0 & \{1,3\} & \partial f_1:\ {-}{-}{-}{-}\,|\,{+}{-}{+}{-}\quad\ C:\ {+}{-}{+}{-}\,|\,{-}{-}{-}{-}\quad\ \partial v_3:\ {+}{+}{+}{+}\,|\,{-}{+}{-}{+} & 1\\
   \mathfrak t_3 & \{1,2,3,4\} & \partial f_1:\ {-}{-}{-}{+}\,|\,{+}{-}{+}{+}\quad\ \partial f_2:\ {+}{-}{+}{+}\,|\,{-}{+}{+}{-}\quad\ \partial f_3:\ {-}{+}{+}{-}\,|\,{-}{-}{-}{+} & t_0^{-1}t_1^{-1}\\
   \mathfrak t_4 & \{1,2,3,4\} & \partial f_1:\ {-}{+}{-}{-}\,|\,{+}{+}{+}{-}\quad\ \partial f_2:\ {+}{-}{-}{+}\,|\,{-}{+}{-}{-}\quad\ \partial f_3:\ {+}{+}{+}{-}\,|\,{+}{-}{-}{+}
   & t_0^{-1}t_1^{-1}
\end{array}
\]}
The class $\mathfrak t_0$ is that of the over-orientation $\1={+}{+}{+}{+}$, and its three circles are the Seifert circles of the closed braid. Fifteen of the sixteen orientations of $\CG$ are the first orientation of a state for some root, all but the under-orientation $-\ve={-}{+}{-}{+}$, and fifteen are the second orientation of a state, all but $\1$. The involution of Proposition~\ref{prop.nat} fixes $\mathfrak t_0$, $\mathfrak t_1$, $\mathfrak t_2$ and exchanges $\mathfrak t_3$ and $\mathfrak t_4$, and these five pairs give the terms $2t_0t_1+1+2t_0^{-1}t_1^{-1}$ of
\[
   \LG(4_1;-t_0,-t_1)=2t_0t_1+3t_0+3t_1+t_0t_1^{-1}+7+t_0^{-1}t_1+3t_0^{-1}+3t_1^{-1}+2t_0^{-1}t_1^{-1},
\]
for every root; they contain the two extreme degrees in $q$, in accordance with Theorem~\ref{thm.extreme}, with $T_{\CG}(1,0)=T_{\CG}(0,1)=2$. The twenty other pairs of classes have a monomial which depends on the root.
\end{example}

\begin{remark}[the terms which depend on the root]\label{rem.rootdependence}
(i) \emph{Degrees.} In the variables $p,q$ of Lemma~\ref{lem.pq}, the term of the pair $(\mathfrak t,\mathfrak t')$ in \eqref{eq.spincsum} is $p^{\alpha}q^{\beta}$ with
\[
   \alpha=-\tfrac12\langle u_{\mathfrak t},\1\rangle-\tfrac12\langle\bar u_{\mathfrak t'},\ve\rangle,\qquad\qquad \beta=c_+-d(u_{\mathfrak t},\bar u_{\mathfrak t'}),
\]
where $c_+$ is the number of positive crossings and $d(u,v)=\langle\frac12(u-v),\frac12(u-v)\rangle$ is the number of coordinates at which $u$ and $v$ differ. So the $p$-degree is the sum of a function of $\mathfrak t$ and of a function of $\mathfrak t'$, and the $q$-degree is $c_+$ minus the square of the norm of a vector of $\BZ^c$ which represents the difference of the two classes in $\BZ^c/\Lambda$. The two ends of Theorem~\ref{thm.extreme} are the pairs for which this vector is zero, resp.\ has no zero coordinate.

(ii) \emph{The root.} By Corollary~\ref{cor.classsmoothing} the representatives $u_{\mathfrak t}$ and $\bar u_{\mathfrak t'}$ change exactly when $\hate$ passes from one circle of the smoothing of $\mathfrak t$, resp.\ of $\mathfrak t'$, to another one; so the monomial of the pair $(\mathfrak t,\mathfrak t')$ is a function of the pair of circles, of these two smoothings, on which $\hate$ lies. Outside the pairs of Proposition~\ref{prop.nat} it does depend on the root in general, even up to the translations of Proposition~\ref{prop.spinc}~(iv). For the closure of the braid $\sigma_1^n$, the graph $\CG$ is a cycle of length $n$, and the class of a vector of $\{\pm1\}^n$ is given by the number of its coordinates $1$ modulo $n$. The class of $\pm\1$ has all the crossings as its cyclic part and the other classes have an empty one, so that the two smoothings consist of the two Seifert circles and of the $n$ circles around the black faces. After a cyclic renumbering of the crossings which depends on $\hate$, the vectors $u_{\mathfrak t}$ and $\bar u_{\mathfrak t'}$ are the vectors $(1,\dots,1,-1,\dots,-1)$ with $a$, resp.\ $b$, coordinates $1$; the monomial of the pair is $t_0^{\min(a,b)}t_1^{\,n-\max(a,b)}$. If $\hate$ is an edge $u'_j=o_{j+1}$, as for the trefoil in Example~\ref{ex.states}, then $0\le a\le n-1$ and $1\le b\le n$; if $\hate$ is an edge $o'_j=u_{j+1}$, then $1\le a\le n$ and $0\le b\le n-1$. So the class of $\pm\1$ is represented by $u=-\1$ and $\bar u=\1$ on one of the two Seifert circles, and by $u=\1$ and $\bar u=-\1$ on the other one, and the two tables of monomials are transposed of each other. For the diagram of the figure-eight knot of Example~\ref{ex.fig8classes}, the eight edges give eight tables which are pairwise different up to translations, and four tables up to translations, to the exchange of $\mathfrak t$ and $\mathfrak t'$ and to the change of sign of the classes; for the closures of the braids $\sigma_1^{-1}\sigma_2\sigma_1^{-1}\sigma_2^3$ and $\sigma_1^{-1}\sigma_2^2\sigma_1^{-2}\sigma_2$, which are diagrams of the knots $6_2$ and $6_3$ up to mirror image, the twelve edges give twelve tables which are pairwise different up to all these symmetries.

(iii) \emph{What is missing.} Thus the lattice gives the set over which \eqref{eq.spincsum} is summed, the form of its terms, the smoothing of each class and the $|\Delta_L(-1)|$ terms of Proposition~\ref{prop.nat}, which contain the two ends of the polynomial in $q$; the plane graph with its root gives the representatives in the other terms. The plane graph is determined by the lattice up to mutation, and Theorem~\ref{thm.gluing} deals with it (Remark~\ref{rem.pieces}~(ii)). The root is not, and two computer experiments indicate that it cannot be disposed of by averaging over it, or by replacing it by a choice made in the lattice.

First, by Corollary~\ref{cor.altpotential}~(ii), $(-1)^{\mu-1}\LG(L;-t_0,-t_1)$ is the sum over the pairs of classes of the average of their monomial over the $2c$ roots, which depends on no choice; but this average is not determined by the lattice. Let $D$ be the closure of the braid $\sigma_1\sigma_2^{-1}\sigma_1\sigma_2^{-1}\sigma_3\sigma_2^{-1}\sigma_3$, and let $D'$ be the mutant of $D$ obtained by rotating by a half-turn in the plane the tangle which contains the first strand and the first three crossings. There are two bijections between the crossings of $D'$ and those of $D$ which preserve the signs and map $\mathrm{Flow}$ to $\mathrm{Flow}$ and $\mathrm{Cut}$ to $\mathrm{Cut}$, so that the two diagrams have the same lattice with its coordinates and the signs of the crossings; and the averages for $D$ and for $D'$ coincide for $39$ of the $441$ pairs of classes with one of these bijections, and for $25$ of them with the other one.

Second, the lattice has its own substitute for a root. A generic linear form on $\mathrm{Flow}(\CG)\otimes\BR$ selects one of the two orientations of every cycle of $\CG$, the one on which it is positive, and a generic linear form on $\mathrm{Cut}(\CG)\otimes\BR$ does the same for the minimal cuts; by a theorem of Backman, Baker and Yuen \cite[Theorem~1.3.1]{BackmanBakerYuen}, every class then contains exactly one vector of $\{\pm1\}^c$ all of whose directed cycles and directed minimal cuts are the selected ones. The representatives $u_{\mathfrak t}$ of a root form a system of this kind, the cycles being oriented with respect to $r_{\rm w}$ and the cuts away from $r_{\rm b}$ (Proposition~\ref{prop.spinc}~(iii)), and the vectors $\bar u_{\mathfrak t}$ form the system of the opposite selection of the cycles. For the standard diagrams of the knots with at most six crossings, of the Hopf link, of the torus link $T(2,4)$, of the Whitehead link and of the Borromean rings, for four pretzel diagrams with six or seven crossings, for the diagram $D$ above and for the closure of $(\sigma_1\sigma_2^{-1})^4$, a diagram of the knot $8_{18}$, we enumerated all the pairs of systems which are obtained in this way from a pair of selections and from the pair with the opposite selection of the cycles. The sum over the pairs of classes of the monomials of Theorem~\ref{thm.tait}~(iii) is $(-1)^{\mu-1}\LG(L;-t_0,-t_1)$ exactly for the following pairs of systems: those of the roots of the diagram; those of the roots of the diagram with the orientation of all the components reversed, which are the opposite vectors; their images under the permutations of the crossings which preserve the lattice and the signs; and, in the case of $D$, those of the roots of its mutants. For the Borromean rings and for $8_{18}$, whose Tait graphs are $3$-connected, these are $24$ of the $1024$ and $32$ of the $62500$ pairs of systems, namely those of the $2c$ roots and their opposites.

So a formula in which the term of every pair of $\mathrm{Spin}^c$ structures is written in the lattice alone cannot be obtained by renumbering or by averaging the terms of \eqref{eq.spincsum}, nor, in these examples, by selecting the representatives by linear forms on the lattice; we do not know one.
\end{remark}

\clearpage

\part*{Conclusion}\label{sec.conclusion}

\subsection*{Toward a topological construction of \texorpdfstring{$\LG$}{LG}}

By the Gordon--Luecke theorem \cite{GordonLuecke}, a knot is determined by its exterior, so every knot invariant, $\LG$ included, is in principle a function of the knot exterior. For the Alexander polynomial this is realized by an explicit construction: $\Delta_K$ is the order of the first homology of the infinite cyclic cover of the exterior, equivalently of the module presented by $tV-V^{\T}$ or by the Fox Jacobian of the knot group. The question which motivates this work is whether there is an analogous \emph{concrete} construction for $\LG$: a topological object attached to the exterior of a knot (a covering space, a local system or twisted homology of the knot group, a torsion) from which $\LG(K;t_0,t_1)$ is computed by a determinant, so that the Alexander-like behavior of $\LG$ becomes a theorem about the exterior rather than a collection of diagrammatic coincidences.

The results above delimit the search, and they change the form of the question. The Kinoshita--Terasaka knot rules out any construction factoring through the Seifert form or the Blanchfield pairing (Remark~\ref{rem.KT}): the object must be finer than the Alexander module of the infinite cyclic cover with its Blanchfield pairing. The identification of the colorings of a diagram with the derivations of the link group (Proposition~\ref{prop.colorings}) and the expansion of \eqref{eq.intro.formula} around the free fermion point (Corollary~\ref{cor.expansion}), in which each crossing contributes through the coupling $x_jy_j\mapsto1$, point toward nonabelian or ``two-variable'' local systems on the exterior: representations of the knot group related to $\mathfrak{sl}(2|1)$, or coverings in which the two variables $t_0,t_1$ play the role of the single variable of the infinite cyclic cover. Section~\ref{sec.cocores} makes the second suggestion precise, and shows what it costs. There is a local system of rank one whose holonomies are the two variables $t_0$ and $t_1$ and of which the coupling is a consequence; but it lives on the complement of the link \emph{and of the co-cores of the bands} of the canonical surface, its holonomy around every co-core is $(t_0t_1)^{\pm1}$, and it extends to the exterior of the link at the free fermion point only (Theorem~\ref{thm.cocores}). In the same framework, a local system of rank one on the exterior itself only gives evaluations of the crossing variables, that is, products of two determinants of Alexander type (Remark~\ref{rem.dictionary}).

The lattice of flows separates what depends on choices from what does not. The edge polynomials are the determinants of the crossing forms on the first homology group of the cut diagram, in any basis (Theorem~\ref{thm.flows}); the Seifert matrices, Alexander's matrix, Kashaev's matrix and the Kauffman states are this one object in different bases, and no basis is needed to state the formula. The crossing forms themselves are identified: for every diagram they are the Seifert form of the canonical surface, cut into the contributions of its bands, with the variable $t$ of the Alexander polynomial replaced, band by band, by the variable of the crossing (Theorems~\ref{thm.localization} and~\ref{thm.deformedseifert}); the first flavor belongs to the canonical surface below the plane of the diagram and the second one to the canonical surface above it, and the coupling $\BE_D$ acts on one band at a time. So the crossing variables are attached to the co-cores of the bands: they form a relative $1$-cochain of the surface, and we have seen that they cannot be replaced by a cohomology class or by a local system on the surface (Section~\ref{ssec.cochains}). This is the precise sense in which the determinant formula still depends on the diagram: through the decomposition of a Seifert surface into disks and bands, that is, through the lattice of rank $c$ which lies over its first homology group of rank $\nu$, and through the coupling. Section~\ref{sec.cocores} gives a topological meaning to this resolution of the Seifert form band by band, and to the coupling. The lattice of rank $c$ is the first homology group of the graph formed by the link and the co-cores of the bands, and its linear forms are the lifted $1$-chains of the complement of this graph; the two crossing forms of a band are the two lifts, from below and from above, of the meridian of its co-core; the coupling $\BE_D$ is the passage from the product of these two lifts to an interleaved product; and $\LG$ is the mixed discriminant of the interleaved chains of the $c$ bands, for the local system of rank one whose holonomies around the boundaries of the disks are $t_1$ and $t_0^{-1}$ (Theorems~\ref{thm.mixed} and~\ref{thm.cocores}). So the two variables and the rule $x_jy_j\mapsto1$ are accounted for by the surface with its disks and its bands. Two things are not. The first one is the way in which the meridians of the co-cores are joined to the base points: it is given by the plane of the diagram, and the result depends on it, through terms which are divisible by $1-t_0t_1$ (Proposition~\ref{prop.paths} and Example~\ref{ex.unknotdisc}). The second one is the fact that the invariant is a mixed discriminant and not a determinant: the determinant of which it is a coefficient is not invariant (Proposition~\ref{prop.kink}), and does not satisfy the genus bound (Section~\ref{ssec.beyond}). The question of a diagram-free definition thus takes a precise form: \emph{characterize, in terms of a Seifert surface decomposed into disks and bands, the paths which join the meridians of the co-cores to the base points in Theorem~\ref{thm.cocores}, and explain the invariance of the resulting mixed discriminant by the moves which relate two such decompositions.} At present this invariance is only known through the determinant formula \eqref{eq.intro.formula}.

At the point $x=-1$ the picture is complete, and it is not restricted to the oriented smoothing: the same lattice carries the forms of Gordon and Litherland of all the state surfaces of the diagram, below and above the plane, with their determinants and their signatures, and Kashaev's form is the sum of the forms of any two complementary surfaces (Section~\ref{sec.minusone}). Away from this point the crossing forms can still be written in a basis adapted to any smoothing (Proposition~\ref{prop.adapted}), but we know a surface interpretation of the result for the oriented smoothing only; whether the state surfaces carry deformed forms of Gordon and Litherland which compute $\LG$, as the canonical surface carries the deformed Seifert form, is a natural question (Remark~\ref{rem.Mminusone}).

Section~\ref{sec.degree} shows how much a topological construction would have to explain and how much the matrices already do: the genus bound for the canonical surface with its sharpness on tight diagrams, the multiplicativity of the top coefficient under Murasugi sums along Seifert circles, and the monicity and the Hopf invariant for alternating links come out of the $\nu\times\nu$ matrices exactly as the classical statements come out of $tV-V^{\T}$. What they do not explain is the genus bound for the other Seifert surfaces, a cancellation inside the coupling $\BE_D$ (Section~\ref{ssec.beyond}) whose elementary case is the hollow-handle law of Proposition~\ref{prop.hollowhandle}, $(x_ax_b-1)R(x)$ against $(y_ay_b-1)\tilde R(y)$; a topological construction should make this cancellation a statement about the surface. The state sum of Section~\ref{sec.kauffman} leads to this bound by another route, modulo the work of Ballinger on his state model, of which it is the case $N=2$ (Corollary~\ref{cor.ballinger} and the discussion which follows it).

Kashaev's matrix adds a second constraint. Its deformation $\tau_D(x,y)$ is a sum of rank-one blocks, one for each crossing, whose rows carry the first flavor and whose columns the second one, and $\BE_D$ acts block by block, on the generators which Liu eliminates from the Seifert form of her surface for $L\sqcup rL$ in her proof of Kashaev's conjecture. A construction in the spirit of that proof would have to deform the Seifert form of Liu's surface, with one pair of variables on each vertex class (Remark~\ref{rem.deformedkashaev}~(i)). What it cannot keep is the Hermitian structure: $\LG$ is real on the unit torus only up to the mirror image, and $\BE_D$ is not an evaluation: carried inside the determinant, it produces a pencil of forms of which $\LG$ is the mixed discriminant, and not one matrix of which a signature could be taken (Remark~\ref{rem.deformedkashaev}~(ii)). The Levine--Tristram signatures carried by Kashaev's matrix thus have no counterpart of this kind in $\LG$ as it stands. At $t=-1$, where the two flavors are real and the coupled form is symmetric, the lattice of flows does carry the signature: the signature of Kashaev's form is $2\sigma_L+\mathrm{wr}(D)$, and that of the form of every state surface is $\sigma_L$ up to the signs of the crossings which are not smoothed along the orientation (Theorem~\ref{thm.signature} and Corollary~\ref{cor.kashaevsignature}). The signature $\sigma_L=\sigma_L(-1)$ is also visible in the state sum of Section~\ref{sec.kauffman} when the diagram is alternating, as the $q$-degree of the pairs of two equal states, and $\LG$ determines it for the fibered alternating links (Proposition~\ref{prop.signature} and Corollary~\ref{cor.extremalstates}). Since $\tau_D$ arose in \cite{Kashaev21} as the classical limit of Gaussian models on the faces of the diagram, a refinement of $\LG$ carrying the other signatures, if it exists, might have to be looked for in the quantum versions of these models rather than in a Seifert form.

For alternating links the state sum of Section~\ref{sec.kauffman} comes closer to such a construction, and Section~\ref{ssec.lattice} says how close. The lattice $\BZ^c\supset\mathrm{Flow}(\CG)\oplus\mathrm{Cut}(\CG)$ of the two checkerboard surfaces of a reduced alternating diagram, which is the image of the lattice of flows under the forms $O_j$ of the over-strands (Remark~\ref{rem.branched}~(ii)), is determined by the double branched cover $\Sigma(L)$ \cite{Greene13}, and the state sum is written in this lattice: the states are the $\mathrm{Spin}^c$ structures of $\Sigma(L)$, the monomial of a pair of states is given by inner products in $\BZ^c$ of two orientations of the Tait graph which represent them, and the over-orientation is the sum of a flow and of a cut (Theorem~\ref{thm.tait} and Propositions~\ref{prop.omega} and~\ref{prop.spinc}). What is not in the lattice is the plane graph with its root, which selects these representatives. The plane graph is determined by the lattice up to mutation, which does not change $\LG$ and across which $\LG$ is given by the gluing formula of Theorem~\ref{thm.gluing}; in particular the polynomial of an alternating knot is determined by the lattice, as is every invariant which does not change under mutation (Corollary~\ref{cor.latticeknots}). The root is a different matter. Each $\mathrm{Spin}^c$ structure determines a smoothing of the diagram, and its representatives change when the root passes from one circle of this smoothing to another one (Corollary~\ref{cor.classsmoothing}). The terms of $|\Delta_L(-1)|$ of the $\Delta_L(-1)^2$ pairs of structures do not depend on the root and are written in the lattice; among them are the two ends of the polynomial in the variable $q$, whose degrees are given by the numbers of positive and of negative crossings of the diagram (Proposition~\ref{prop.nat} and Theorem~\ref{thm.extreme}). The monomials of the other pairs depend on the root, only their sum does not, and the experiments of Remark~\ref{rem.rootdependence} indicate that the root cannot be replaced by a choice made in the lattice. A formula for all the terms in the lattice alone, for instance in terms of the quantities from which the correction terms of $\Sigma(L)$ are computed, would be a topological construction of $\LG$ for alternating links; we do not have one.

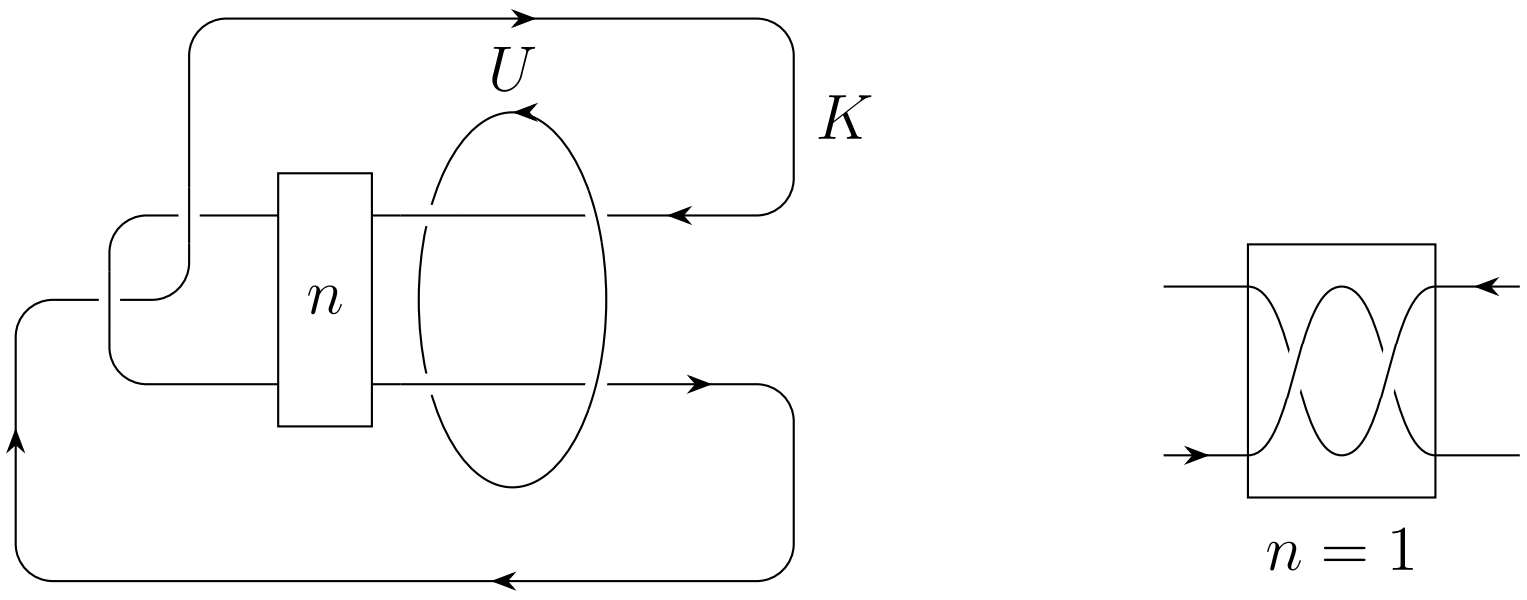
\begin{figure}[ht]
\centering
\begin{tikzpicture}[>=Stealth,scale=0.62]
\draw[rounded corners=7pt,postaction={decorate},decoration={markings,mark=at position 0.0231 with {\arrow{>}},mark=at position 0.3565 with {\arrow{>}},mark=at position 0.4977 with {\arrow{>}},mark=at position 0.6551 with {\arrow{>}},mark=at position 0.8819 with {\arrow{>}}}] (3.000,0.900) -- (-4.300,0.900) -- (-4.300,-0.900) -- (3.000,-0.900) -- (3.000,-3.000) -- (-5.300,-3.000) -- (-5.300,0.000) -- (-3.450,0.000) -- (-3.450,3.000) -- (3.000,3.000) -- cycle;
\draw[postaction={decorate},decoration={markings,mark=at position 0.25 with {\arrow{>}}}] (0,0) ellipse (1.000 and 2.000);
\draw[preaction={draw=white,line width=4pt}] (-3.450,0.600) -- (-3.450,1.200);
\draw[preaction={draw=white,line width=4pt}] plot[domain=17.74:35.74,samples=12] ({1.000*cos(\x)},{2.000*sin(\x)});
\draw[preaction={draw=white,line width=4pt}] (-0.593,0.900) -- (-1.193,0.900);
\draw[preaction={draw=white,line width=4pt}] (-4.300,0.300) -- (-4.300,-0.300);
\draw[preaction={draw=white,line width=4pt}] (-1.193,-0.900) -- (-0.593,-0.900);
\draw[preaction={draw=white,line width=4pt}] plot[domain=-35.74:-17.74,samples=12] ({1.000*cos(\x)},{2.000*sin(\x)});
\filldraw[fill=white] (-2.500,-1.350) rectangle (-1.500,1.350);
\node at (-2.000,0) {$n$};
\node[above] at (0,2.000) {$U$};
\node[right] at (3.000,1.950) {$K$};
\end{tikzpicture}
\hspace{1.5cm}
\begin{tikzpicture}[>=Stealth,scale=0.62]
\draw[postaction={decorate},decoration={markings,mark=at position 0.08 with {\arrow{>}}}] (1.900,0.900) -- (1.000,0.900) .. controls (0.500,0.900) and (0.500,-0.900) .. (0.000,-0.900) .. controls (-0.500,-0.900) and (-0.500,0.900) .. (-1.000,0.900) -- (-1.900,0.900);
\draw[postaction={decorate},decoration={markings,mark=at position 0.08 with {\arrow{>}}}] (-1.900,-0.900) -- (-1.000,-0.900) .. controls (-0.500,-0.900) and (-0.500,0.900) .. (0.000,0.900) .. controls (0.500,0.900) and (0.500,-0.900) .. (1.000,-0.900) -- (1.900,-0.900);
\draw[preaction={draw=white,line width=4pt}] (0.580,0.289) -- (0.420,-0.289);
\draw[preaction={draw=white,line width=4pt}] (-0.580,-0.289) -- (-0.420,0.289);
\draw (-1.000,-1.350) rectangle (1.000,1.350);
\node[below] at (0.000,-1.450) {$n=1$};
\end{tikzpicture}
\caption{The twisted Whitehead links $W_n$. The box contains $n$ full twists of the band formed by the two strands of $K$ which pass through the disk bounded by $U$, that is, $2|n|$ crossings, positive for $n>0$ (the box for $n=1$ is shown on the right) and negative for $n<0$. $W_0$ is the Whitehead link; the exteriors of all the $W_n$ are homeomorphic.}
\label{fig.whitehead}
\end{figure}

For links the Gordon--Luecke theorem fails, and so does the conclusion drawn from it: the homeomorphism type of the exterior does not determine $\LG$. Let $W_n$ be the link obtained from the Whitehead link $W=W_0$ by $n$ full twists along the disk bounded by its unknotted component $U$ (Figure~\ref{fig.whitehead}). The twist is a homeomorphism of the solid torus $S^3\sm\mathring N(U)$, so all the $W_n$ have homeomorphic exteriors, and they have the same Alexander polynomial $(t-1)^3$ up to a unit; but their Links--Gould polynomials, computed as in Remark~\ref{rem.KT}, are pairwise distinct for $n=-1,0,1,2$, their constant terms being $-8,-10,-6,-4$. For links, the object to be constructed must therefore involve the peripheral structure of the exterior, and the comparison between $\LG$ and the exterior becomes a question in itself.

\subsection*{Categorification}
The search for a topological construction of $\LG$ has a companion, the search for a categorification, and for the Alexander polynomial the two questions have the same answer. Knot Floer homology \cite{OSknot, RasmussenThesis} is a bigraded homology theory whose graded Euler characteristic is $\Delta_K$. It is an invariant of the exterior, namely the sutured Floer homology of the exterior of the knot with two meridional sutures \cite{Juhasz06}, and it turns the three-dimensional properties of $\Delta_K$ into exact statements: it detects the genus of the knot \cite{OSgenus} and its fiberedness \cite{Ghiggini08, Ni07}, its group of extremal Alexander degree is the sutured Floer homology of the complement of a Seifert surface of minimal genus \cite{Juhasz08}, and this group is multiplicative under Murasugi sums \cite{Ni06, CHS25}. The Links--Gould invariant has the decategorified form of each of these three properties: the genus bound \cite{Kohli-Tahar, LNVdV25}, the monicity on fibered links \cite{LNVdV25} and the multiplicativity of the top coefficient under Murasugi sums \cite{ALvdV26}; Section~\ref{sec.degree} shows how far the deformed Seifert matrices go in this direction. This suggests looking for a homology theory $\mathcal H(L)$ with three gradings whose graded Euler characteristic is $\LG(L;t_0,t_1)$, and in which the group of extremal degree in the variable $p$ of Lemma~\ref{lem.pq}, the degree which is bounded by the genus, is nonzero in the degree given by the genus, has rank one for the fibered links and satisfies a K\"unneth formula for Murasugi sums. No categorification of $\LG$ has been constructed yet. On the algebraic side, such a construction requires a categorification of the four-dimensional representations of $U_q\mathfrak{gl}(2|1)$ and of their braiding. For $U_q\mathfrak{gl}(1|1)$ and the Alexander polynomial this is what the bordered algebras of Ozsv\'ath and Szab\'o provide \cite{OSbordered, Manion19}, and the positive half of $U_q\mathfrak{gl}(m|1)$ has been categorified by Khovanov and Sussan \cite{KhovanovSussan}.

The state sum of Section~\ref{sec.kauffman} says what the chain groups of such a theory could be. The knot Floer complex of a knot diagram is generated by the Kauffman states of the diagram, with gradings which are sums of local contributions of the crossings \cite{OSalternating}, and in \cite{OSbordered} a bigraded invariant closely connected to knot Floer homology is computed from these generators by cutting the diagram into elementary pieces. By Proposition~\ref{prop.statesum}~(iii), $\LG$ is the graded Euler characteristic of the free abelian group $C(D,\hate)$ generated by the ordered pairs of Kauffman states of $(D,\hate)$, a pair $(\mathfrak s,\mathfrak s')$ being given the bidegree of its monomial and the parity of $b(\mathfrak s)+b(\mathfrak s')$. The degrees and the parities are sums of local contributions of the crossings, and at the free fermion point, where the two states of a pair do not interact, $C(D,\hate)$ is the tensor square of the group generated by the states (Remark~\ref{rem.statesum}). So the first question is whether $C(D,\hate)$ carries a differential which preserves the bidegree, changes the parity, and has a homology which does not depend on the diagram. Section~\ref{sec.kauffman} gives three constraints. First, for an alternating diagram the parity of a pair is determined by its bidegree (Theorem~\ref{thm.kauffmanalt}), so that such a differential is zero: the homology of an alternating link would be free of rank $\Delta_L(-1)^2$ and determined by $\LG$, in the same way as the knot Floer homology of an alternating knot is determined by its Alexander polynomial and its signature \cite{OSalternating}. The signature would enter here as it does there, since the pairs of two equal states all have the total degree $-\sigma_L$ (Proposition~\ref{prop.signature}); and the rank $\Delta_L(-1)^2$ is that of the Heegaard Floer homology of the connected sum of two copies of the double branched cover $\Sigma(L)$, which for an alternating link has one generator for each $\mathrm{Spin}^c$ structure \cite{OSbranched}, in agreement with Proposition~\ref{prop.spinc}. Second, at $t_0t_1=-1$ the pairs of two distinct states cancel, within each superposition, and what remains is the set of the pairs of two equal states, that is, the generators of the knot Floer complex with their Alexander degrees doubled (Corollary~\ref{cor.diagonal} and Remark~\ref{rem.ballinger}~(ii)); at $t_0t_1=1$ the pairs are the generators of the knot Floer complex of the connected sum of two copies of the knot. The two specializations of $\LG$ to the Alexander polynomial should therefore be categorified by spectral sequences from $\mathcal H(L)$ to knot Floer homology, as the specializations of the HOMFLY-PT polynomial are categorified by spectral sequences which start from its triply graded homology \cite{Rasmussen15, RobertWagner}. Third, in the form of Corollary~\ref{cor.ballinger} every cycle of a superposition contributes the factor $[2]$, the graded rank of the homology of a circle; this is the form of the state sums which Khovanov homology categorifies \cite{Khovanov00}, the cycles of the superposition taking the place of the circles of a smoothing.

For the group of extremal degree these constraints meet the question of the first part of this conclusion. Let $D$ be a reduced alternating diagram of a knot $K$. By Corollary~\ref{cor.extremalstates} the top coefficient of $\LG$ is the sum of the monomials $q^{\ve(X(\mathfrak s)\triangle Y(\mathfrak s'))}$ over the $N^2$ ordered pairs of extremal states, and $K$ is fibered if and only if $N=1$. The $N$ extremal states generate the knot Floer homology of $K$ in its extremal Alexander degree, which has rank $N$ and is the sutured Floer homology of the complement of the canonical surface $F_D$ \cite{OSalternating, Juhasz08}. This group is graded by the relative $\mathrm{Spin}^c$ structures of the complement, which form an affine space over $H_1(S^3\sm F_D)\cong H^1(F_D)$; its Euler characteristic is a torsion with coefficients in the group ring of $H_1(S^3\sm F_D)$, and for the Seifert surfaces of the alternating knots each structure carries at most one generator \cite{FJR11}. So the top coefficient of $\LG$ has as many terms as there are ordered pairs of these structures which carry a generator. The comparison also shows what is missing. The sutured Floer homology of the complement of a surface is graded by classes, the elements of an affine space over $H^1(F)$, whereas the crossing variables form a cochain which is not determined by its class (Section~\ref{ssec.cochains}). A categorification of the deformed Seifert form, like its topological construction, has to remember how the surface is cut into disks and bands, or to explain why the result does not depend on it.

\clearpage
\section*{Index of notation}\label{sec.notation}

{\small
This index lists the notation which is used in more than one place, in the order in which it appears, with the place where it is defined. Symbols which are only used inside one proof or one example are not listed. A few letters have different meanings in different sections, for instance $G$, $\tau$, $\kappa$, $\chi$, $\omega$, $\Sigma$ and $H_\Omega$; each meaning is listed with the section to which it belongs. Throughout, $D$ is a diagram, $j$ and $k$ are crossings, $e$ is an edge, $f$ is a face and $C$ is a circle.

\begingroup
\setlength{\tabcolsep}{4pt}
\renewcommand{\arraystretch}{1.12}
\def\nt#1#2#3#4{#1 & #2 & #3 & \pageref{#4}\\}
\def\nh#1{\multicolumn{4}{@{}l}{\rule{0pt}{3.6ex}\textbf{#1}}\\*[2pt]}
\begin{longtable}{@{}>{\raggedright\arraybackslash}p{0.21\textwidth}>{\raggedright\arraybackslash}p{0.52\textwidth}>{\raggedright\arraybackslash}p{0.135\textwidth}>{\raggedleft\arraybackslash}p{0.06\textwidth}@{}}
\emph{Symbol} & \emph{Meaning} & \emph{Defined in} & \emph{Page}\\ \hline
\endfirsthead
\emph{Symbol} & \emph{Meaning} & \emph{Defined in} & \emph{Page}\\ \hline
\endhead
\nh{General conventions}
\nt{$\LG(L)$, $\LG(L;t_0,t_1)$}{the Links--Gould polynomial of the oriented link $L$}{Thm.~\ref{thm.HK}}{thm.HK}
\nt{$\mu$}{the number of components of $L$}{Thm.~\ref{thm.HK}}{thm.HK}
\nt{$\BK$}{the field $\BQ(t_0,t_1)$}{Conventions}{ssec.conventions}
\nt{$t$}{an indeterminate, the variable of the Alexander polynomial}{Conventions}{ssec.conventions}
\nt{$M^{\T}$, $M_{\hat d}$}{the transpose of a matrix $M$; the matrix $M$ with the row and the column of the index $d$ removed}{Conventions}{ssec.conventions}
\nt{$e_f$, $\1$}{the unit column vector of the index $f$; the column vector with all entries equal to $1$}{Conventions}{ssec.conventions}
\nt{$\ind{P}$}{$1$ if the statement $P$ holds, and $0$ otherwise}{Conventions}{ssec.conventions}
\nt{$\Delta_L$, $\nabla_L$}{the Alexander polynomial and the Conway polynomial of $L$; in the Conway normalization $\Delta_L(t)=\nabla_L(t^{1/2}-t^{-1/2})$}{Lem.~\ref{lem.alexpoly}}{lem.alexpoly}
\nt{$g(L)$}{the genus of $L$}{\S\ref{sec.degree}}{sec.degree}
\nt{$p$, $q$}{the variables with $t_0=p^{-1}q$ and $t_1=pq$}{Lem.~\ref{lem.pq}}{lem.pq}
\nt{$\sigma_L$}{the signature of $L$; the right-handed trefoil has signature $-2$}{\S\ref{ssec.signature}}{ssec.signature}
\nh{The determinant formula (Section~\ref{sec.formula})}
\nt{$D$, $c$, $\ve_j$}{a diagram: a connected oriented link diagram on $S^2$ with $c\ge1$ crossings; the sign of the crossing $j$}{\S\ref{ssec.diagrams}}{ssec.diagrams}
\nt{$E$}{the set of the $2c$ edges of $D$}{\S\ref{ssec.diagrams}}{ssec.diagrams}
\nt{$o_j,u_j,o'_j,u'_j$}{the incoming edges of the over-strand and of the under-strand at the crossing $j$, and the outgoing ones}{Fig.~\ref{fig.crossing}}{fig.crossing}
\nt{$x_j$, $y_j$}{the crossing variables of the first and of the second flavor}{\eqref{eq.Ta}, \eqref{eq.Tb}}{eq.Ta}
\nt{$T^a_D(x)$, $T^b_D(y)$}{the edge matrices of $D$}{\eqref{eq.Ta}, \eqref{eq.Tb}}{eq.Ta}
\nt{$\hate$}{the cut edge}{\eqref{eq.PP}}{eq.PP}
\nt{$P_{D,\hate}(x)$, $\tP_{D,\hate}(y)$}{the edge polynomials of the two flavors (also called cut polynomials)}{\eqref{eq.PP}}{eq.PP}
\nt{$\BE_D$}{the coupling: the linear map with $x_j\mapsto t_0^{\ve_j}$, $y_j\mapsto t_1^{\ve_j}$ and $x_jy_j\mapsto1$}{\eqref{eq.ED}}{eq.ED}
\nt{$\bx_j$, $\by_j$, $g_j$}{the units $t_0^{\ve_j}$ and $t_1^{\ve_j}$, and $g_j=1-(t_0t_1)^{\ve_j}$}{\eqref{eq.xbar}}{eq.xbar}
\nt{$a_S$, $b_S$}{the coefficients of the expansion of the edge polynomials around the free fermion point}{\eqref{eq.ab}}{eq.ab}
\nh{The Alexander matrix and the levels (Section~\ref{sec.alexander})}
\nt{$t^{\ve}$}{the Alexander point $x_j=t^{\ve_j}$, $j=1,\dots,c$}{\S\ref{sec.alexander}}{sec.alexander}
\nt{$A_D(t)$, $B_D(t)$}{the Alexander matrices $1-T^a_D(t^{\ve})$ and $1-T^b_D(t^{\ve})$}{Def.~\ref{def.coloring}}{def.coloring}
\nt{$A_{D,\hate}(t)$, $B_{D,\hate}(t)$}{the cut matrices: the row and the column of $\hate$ are removed}{Def.~\ref{def.coloring}}{def.coloring}
\nt{$s(C)$}{$1$ if the Seifert circle $C$ is clockwise, $-1$ otherwise}{Def.~\ref{def.level}}{def.level}
\nt{$h(C)$, $h(e)$}{the level of a Seifert circle and of an edge}{Def.~\ref{def.level}}{def.level}
\nt{$\DD_t$, $\ell_t$}{the diagonal matrix and the row vector with the entries $t^{h(e)}$}{Lem.~\ref{lem.level}}{lem.level}
\nt{$G$, $\phi$}{in Section~\ref{sec.alexander}: the group of the link and the total linking number homomorphism $G\to\BZ$}{\S\ref{ssec.derivations}}{not.group}
\nt{$\b(L)$}{the nullity of $L$, the rank of its Alexander module}{Def.~\ref{def.nullity}}{def.nullity}
\nh{Flows and crossing forms (Sections~\ref{ssec.flowlattice}--\ref{ssec.facematrices})}
\nt{$z$, $z_e$}{a flow on $D$; its coordinate at the edge $e$}{Def.~\ref{def.flows}}{def.flows}
\nt{$\FF_D$, $\FF_{D,\hate}$}{the lattice of flows of $D$; the lattice of the flows which vanish at $\hate$}{Def.~\ref{def.flows}}{def.flows}
\nt{$\CF$}{the set of the $c+2$ faces of $D$}{\S\ref{ssec.flowlattice}}{not.faces}
\nt{$f_{\rm L}(e)$, $f_{\rm R}(e)$}{the faces on the left and on the right of the edge $e$}{\S\ref{ssec.flowlattice}}{not.faces}
\nt{$\partial f$}{the boundary of the face $f$, a flow}{\eqref{eq.boundaryface}}{eq.boundaryface}
\nt{$m(e)$}{the vector $e_{f_{\rm L}(e)}-e_{f_{\rm R}(e)}\in\BZ^{\CF}$}{\eqref{eq.m}}{eq.m}
\nt{$k$, $k_f$}{a potential of a flow: $z_e=k_{f_{\rm L}(e)}-k_{f_{\rm R}(e)}$}{Lem.~\ref{lem.faces}}{lem.faces}
\nt{$f^\uparrow_j,f^\leftarrow_j,f^\downarrow_j,f^\rightarrow_j$}{the four faces at the crossing $j$}{Fig.~\ref{fig.faces}}{fig.faces}
\nt{$\pi_j$, $\tpi_j$}{the crossing forms of the two flavors: $\pi_j(z)=x_jz_{u'_j}-z_{u_j}$ and $\tpi_j(z)=y_jz_{o_j}-z_{o'_j}$}{Def.~\ref{def.crossingforms}}{def.crossingforms}
\nt{$\d_j$}{the generator of the band of the crossing $j$, the jump $[u_j]-[o_j]$}{\S\ref{ssec.crossingforms}}{not.delta}
\nt{$\kappa$}{the sign of a basis of $\FF_{D,\hate}$}{Thm.~\ref{thm.flows}}{thm.flows}
\nt{$\tau$}{the coupled form $\sum_j\ve_j\,\pi_j\otimes\tpi_j$}{Rem.~\ref{rem.flows}}{rem.flows}
\nt{$\ve_D$}{the product of the signs of the crossings}{\eqref{eq.LGtau}}{eq.LGtau}
\nt{$\vs$, $s_\vs$}{a smoothing of $D$, a vector of $\{0,1\}^c$; the number of its circles}{Def.~\ref{def.smoothing}}{def.smoothing}
\nt{$z^C$}{the flow of a circle $C$ of a smoothing}{\S\ref{ssec.smoothings}}{not.circleflows}
\nt{$w^\vs_j$, $w_j$}{the transfer at the crossing $j$ for the smoothing $\vs$, and for the oriented smoothing}{Def.~\ref{def.transfer}}{def.transfer}
\nt{$U_j$, $O_j$}{the linear forms $z_{u'_j}+z_{u_j}$ and $z_{o_j}+z_{o'_j}$}{\eqref{eq.UO}}{eq.UO}
\nt{$\nu_\vs$}{the number $c+1-s_\vs$}{Prop.~\ref{prop.adapted}}{prop.adapted}
\nt{$C_1,\dots,C_s$}{the Seifert circles of $D$, the circles of the oriented smoothing}{\S\ref{ssec.smoothings}}{not.seifertgraph}
\nt{$A_j$, $B_j$}{the Seifert circle of the edges $o_j,u'_j$ and that of the edges $u_j,o'_j$}{\S\ref{ssec.smoothings}}{not.seifertgraph}
\nt{$\Gamma_D$, $H$, $\nu$}{the Seifert graph of $D$, its cycle space, and the rank $\nu=c-s+1$ of this space}{\S\ref{ssec.smoothings}}{not.seifertgraph}
\nt{$s_C(j)$}{$1$ if $A_j=C$, and $-1$ if $B_j=C$}{\eqref{eq.sw}}{eq.sw}
\nt{$b_C$, $\hat C$}{the base edge of the Seifert circle $C$; the circle which contains $\hate$}{\eqref{eq.baseedge}}{eq.baseedge}
\nt{$\mathrm{Pre}_C(e)$, $\prec_C$}{the set of the crossings met along $C$ from $b_C$ to $e$; the order of the crossings of $C$}{\S\ref{ssec.smoothings}}{eq.baseedge}
\nt{$z^\gamma$}{the lift of a cycle $\gamma\in H$ to a flow}{Lem.~\ref{lem.lifts}}{lem.lifts}
\nt{$Y$, $\CT$, $E_{\CT}$}{a matrix whose columns form a basis of $H$; a spanning tree of $\Gamma_D$; the matrix of the unit vectors of its edges}{\S\ref{ssec.smoothings}}{not.Y}
\nt{$J_C(e)$, $\CC_C$}{the jump sum of the edge $e$; the full turn of the circle $C$}{\eqref{eq.jumpsum}, \eqref{eq.fullturn}}{eq.jumpsum}
\nt{$K(x)$, $K^b(y)$}{the matrices of the jump relations of the two flavors}{\eqref{eq.Krow}, \eqref{eq.Kbrow}}{eq.Krow}
\nt{$\CC(x)$}{the deformed incidence matrix of the Seifert graph}{\eqref{eq.CC}}{eq.CC}
\nt{$\partial$, $\hat\partial$}{the incidence map of $\Gamma_D$; the same with the row of $\hat C$ removed}{\S\ref{ssec.elimination}}{eq.CC}
\nt{$N_{D,\hate}(x)$, $\tilde N_{D,\hate}(y)$}{the matrices of the crossing forms in a Seifert basis}{\eqref{eq.N}}{eq.N}
\nt{$Q_D(x)$, $\tQ_D(y)$}{the face matrices; $Q_{D,\hate}$, $\tQ_{D,\hate}$ are these matrices without the columns of the two faces adjacent to $\hate$}{Def.~\ref{def.facematrices}}{def.facematrices}
\nt{$a_j(x_j)$, $b_j(y_j)$}{the rows of the face matrices}{\eqref{eq.facerows}}{eq.facerows}
\nt{$\bD$}{the inverse diagram: all orientations reversed and all crossings changed}{Lem.~\ref{lem.inverse}}{lem.inverse}
\nh{The coupling as a mixed discriminant (Section~\ref{ssec.mixed})}
\nt{$\lb\mathsf M_1,\dots,\mathsf M_c\rb$}{the coefficient of $\lambda_1\cdots\lambda_c$ in $\det(\lambda_1\mathsf M_1+\dots+\lambda_c\mathsf M_c)$, the mixed discriminant up to the factor $c!$}{Def.~\ref{def.mixed}}{def.mixed}
\nt{$\lb\mathsf B_1,\dots,\mathsf B_c\rb_{\hate}$}{the same for the matrices of $c$ bilinear forms on $\FF_D$ in a basis of $\FF_{D,\hate}$}{Def.~\ref{def.mixed}}{def.mixed}
\nt{$\Psi_j[\sv;\xs,\ys]$}{the bilinear form $\sv\,z_{u'_j}z'_{o_j}-\xs\,z_{u'_j}z'_{o'_j}-\ys\,z_{u_j}z'_{o_j}+z_{u_j}z'_{o'_j}$}{\eqref{eq.Psigeneral}}{eq.Psigeneral}
\nt{$\Psi_j$}{the coupled crossing form $\Psi_j[1;t_0^{\ve_j},t_1^{\ve_j}]$}{\eqref{eq.Psi}}{eq.Psi}
\nt{$\pi_j^{(\xs)}$, $\tpi_j^{(\ys)}$}{the crossing forms in which the variable is replaced by $\xs$, resp.\ by $\ys$}{\S\ref{ssec.mixed}}{eq.rankone}
\nt{$\BE'$}{the linear map with $x_j\mapsto\xs_j$, $y_j\mapsto\ys_j$ and $x_jy_j\mapsto\sv_j$}{Thm.~\ref{thm.mixed}}{thm.mixed}
\nt{$L^*$, $-L$}{the mirror image of $L$; the link $L$ with all its components reversed}{Cor.~\ref{cor.symmetries}}{cor.symmetries}
\nt{$\CDp_{D,\hate}(\lambda)$}{the determinant of the pencil $\sum_j\lambda_j\Psi_j$ of the coupled crossing forms}{Def.~\ref{def.pencil}}{def.pencil}
\nh{State surfaces and linking numbers (Section~\ref{sec.statesurfaces})}
\nt{$P$}{the horizontal plane $\BR^2\times\{0\}$ of the diagram in $\BR^3$}{\S\ref{ssec.statesurface}}{ssec.statesurface}
\nt{$\CB_j$, $r_0$, $h_0$}{the disk of radius $r_0$ around the crossing $j$ in $P$; the height of the over-strand above the crossing}{\S\ref{ssec.statesurface}}{ssec.statesurface}
\nt{$\alpha_j$}{the co-core of the crossing $j$: the vertical segment from the under-strand to the over-strand}{\S\ref{ssec.statesurface}}{not.G}
\nt{$G$}{from Section~\ref{sec.statesurfaces} on: the graph $L\cup\alpha_1\cup\dots\cup\alpha_c$ in $\BR^3$}{\eqref{eq.flowsG}}{not.G}
\nt{$\Omega_C$, $D_C$}{the disk bounded by the circle $C$ in $P$; the disk of $C$ in the state surface}{\eqref{eq.stacking}}{eq.stacking}
\nt{$b_j$, $\kappa_j$}{the band of the crossing $j$ and its core}{\eqref{eq.ruled}}{eq.ruled}
\nt{$F_\vs$, $F_D$}{the state surface of the smoothing $\vs$; the canonical Seifert surface $F_D=F_0$}{Def.~\ref{def.statesurface}}{def.statesurface}
\nt{$\Gamma_\vs$}{the state graph of the smoothing $\vs$}{\S\ref{ssec.statehomology}}{not.stategraph}
\nt{$\lk$}{the linking number}{\S\ref{ssec.spanning}}{ssec.spanning}
\nt{$\Delta_f$}{the face disk of the face $f$, lifted above the plane}{\eqref{eq.uf}}{eq.uf}
\nt{$\tau b$, $b^+$, $b^-$}{the push-off of a curve $b$ on a surface to both sides; to its positive and to its negative side}{\S\ref{ssec.pushoff}}{ssec.pushoff}
\nt{$\theta_j$}{the sign $\ve_j(-1)^{\vs_j}$ of the band of the crossing $j$}{Def.~\ref{def.bandsign}}{def.bandsign}
\nt{$\CG_F$}{the form of Gordon and Litherland of the surface $F$}{\S\ref{ssec.localization}}{not.GL}
\nt{$V_F$}{the Seifert form of the canonical surface $F_D$, and its matrix}{\eqref{eq.seifertform}}{not.VF}
\nt{$e^+_j$, $e^-_j$}{the edges of the under-strand on the left and on the right of the over-strand}{\eqref{eq.eplus}}{eq.eplus}
\nt{$\hat e^+_j$, $\hat e^-_j$}{the edges of the over-strand on the left and on the right of the under-strand}{Prop.~\ref{prop.otherside}}{prop.otherside}
\nt{$D^*$}{the diagram $D$ with all its crossings changed}{\S\ref{ssec.otherside}}{ssec.otherside}
\nt{$F^\uparrow_\vs$, $F^\uparrow_D$}{the state surface and the canonical surface above the plane}{Def.~\ref{def.uppersurface}}{def.uppersurface}
\nh{Deformed Seifert matrices (Section~\ref{sec.seifert})}
\nt{$V_j$, $V'_j$, $V^\uparrow_j$, $V'^{\uparrow}_j$}{the contributions of the band $j$ to the Seifert forms of $F_D$ and of $F^\uparrow_D$ and to their transposes}{\eqref{eq.Vj}}{eq.Vj}
\nt{$C_{\rm L}$, $C_{\rm R}$}{the Seifert circles on the left and on the right of a crossing}{\S\ref{ssec.bandbyband}}{ssec.bandbyband}
\nt{$\Omega$, $H_\Omega$, $H_\gamma$}{a region of the complement of the Seifert circles; in Sections~\ref{sec.seifert}, \ref{sec.degree} and~\ref{sec.faces}: the lower one of the two levels of its boundary circles; the same for the region of a cycle $\gamma$}{Lem.~\ref{lem.regions}}{lem.regions}
\nt{$j(C)$, $x^{\CT}$}{the edge of $\CT$ which joins $C$ to its parent; the monomial of the tree}{Cor.~\ref{cor.seifert}}{cor.seifert}
\nt{$Y(x)$, $H_x$}{the matrix whose columns form a basis of the deformed cycle space $H_x=\ker\CC(x)$}{\eqref{eq.Yx}}{eq.Yx}
\nt{$M_{D,\hate}(x)$, $M^b_{D,\hate}(y)$}{the deformed Seifert matrices of the two flavors}{Def.~\ref{def.deformedseifertmatrices}}{def.deformedseifertmatrices}
\nt{$\Theta_x$}{the deformed Seifert form}{\eqref{eq.Theta}}{eq.Theta}
\nt{$\mathsf Z_t$, $\mathsf Z'_t$}{diagonal matrices of units at the Alexander point}{Lem.~\ref{lem.Yt}}{lem.Yt}
\nt{$\rho(x_j)$, $\rho^-(x_j)$}{the blocks of a positive and of a negative crossing of a braid}{\S\ref{ssec.moves}}{ssec.moves}
\nh{Degrees (Section~\ref{sec.degree})}
\nt{$\deg$, $\mathrm{span}$}{the degree $\deg(t_0^at_1^b)=a-b$, which is minus the degree in $p$; the difference between the largest and the smallest degree}{\eqref{eq.KTdegree}}{eq.KTdegree}
\nt{$\ve(S)$, $x^S$}{the sum $\sum_{j\in S}\ve_j$ and the monomial $\prod_{j\in S}x_j$, for a set $S$ of crossings}{\S\ref{ssec.epsdegree}}{ssec.epsdegree}
\nt{$\deg_\ve$, $\deg^+_\ve$, $\deg^-_\ve$}{the $\ve$-degree, with $\deg_\ve x_j=\deg_\ve y_j=\ve_j$; the largest and the smallest $\ve$-degree of a monomial}{\eqref{eq.epsdegree}}{eq.epsdegree}
\nt{$P_S$, $\tP_T$}{the coefficients of the edge polynomials in the monomial basis}{\eqref{eq.monomials}}{eq.monomials}
\nt{$\delta_{D,\hate}$, $I_{D,\hate}$}{the interval $I_{D,\hate}=[\delta_{D,\hate}-\nu,\delta_{D,\hate}]$ of the $\ve$-degrees of the edge polynomials}{Thm.~\ref{thm.degrees}}{thm.degrees}
\nt{$\mathrm{top}(L,q)$}{the coefficient of $p^{2g(L)+\mu-1}$ in $\LG(L)$, the top coefficient}{Thm.~\ref{thm.tight}}{thm.tight}
\nt{$D_{\rm in}$, $D_{\rm out}$}{the parts of $D$ inside and outside a Seifert circle}{\S\ref{ssec.murasugi}}{ssec.murasugi}
\nt{$\tau_D(q)$}{in Section~\ref{ssec.murasugi}: the coefficient of $p^{\nu}$ in $\LG(L;p,q)$}{\S\ref{ssec.murasugi}}{not.tauD}
\nt{$S(\tau)$, $S'(\tau')$}{the sets of crossings of an under-tree $\tau$ and of an over-tree $\tau'$ of an alternating diagram}{\eqref{eq.HKtrees}}{eq.HKtrees}
\nt{$\lambda(\xi_F)$}{the Hopf invariant of the plane field of a fiber $F$}{Thm.~\ref{thm.alternating}}{thm.alternating}
\nt{$\Sigma_J$}{the sector graph of a set $J$ of crossings}{Prop.~\ref{prop.sectors}}{prop.sectors}
\nh{The complement of the co-cores (Section~\ref{sec.cocores})}
\nt{$X_G$}{the complement $\BR^3\sm G$ of the link and of the co-cores}{\S\ref{sec.cocores}}{sec.cocores}
\nt{$\Sigma$, $\Delta^\circ_f$}{the union of the face disks $\Delta_f$ of all the faces; the part of $\Delta_f$ which lies over $f$}{\S\ref{ssec.XG}}{ssec.XG}
\nt{$X_+$, $X_-$}{the sets of the points of $\BR^3$ above and below $\Sigma$}{Lem.~\ref{lem.XG}}{lem.XG}
\nt{$\star_+$, $\star_-$}{the base points, in $X_+$ and in $X_-$}{\S\ref{ssec.XG}}{not.basepoints}
\nt{$\pth_f$}{the path from $\star_+$ to $\star_-$ which goes down through the face $f$}{\S\ref{ssec.XG}}{not.basepoints}
\nt{$\vth=(\vth_f)_{f\in\CF}$}{a multiplicative potential: a unit for each face}{Def.~\ref{def.potential}}{def.potential}
\nt{$\mhol_e$}{the holonomy $\vth_{f_{\rm L}(e)}/\vth_{f_{\rm R}(e)}$ around the edge $e$}{Def.~\ref{def.potential}}{def.potential}
\nt{$\rho_\vth$, $\CL_\vth$}{the holonomy of the loops of $X_G$; the local system of rank one of the potential $\vth$}{\eqref{eq.holonomy}}{eq.holonomy}
\nt{$C_\vth$, $\tpth_f$, $\tstar_\pm$}{the module of the lifted chains and its basis; the generators of the fibers at the base points}{Def.~\ref{def.chains}}{def.chains}
\nt{$\iota$, $\iota'$}{the two twisted linking numbers, from the chains to the linear forms on the flows}{Def.~\ref{def.iota}}{def.iota}
\nt{$C^-_\vth$, $C^+_\vth$}{the chains whose boundary is a multiple of $\tstar_-$, resp.\ of $\tstar_+$}{Def.~\ref{def.iota}}{def.iota}
\nt{$f^1_j,f^2_j,f^3_j,f^4_j$}{the four faces around the crossing $j$, in the clockwise order}{\eqref{eq.fourfaces}}{eq.fourfaces}
\nt{$\mer_j$, $\mer^\circ_j$}{the meridian of the co-core $\alpha_j$, a word in the paths $\pth_f$; a horizontal circle which represents it}{Def.~\ref{def.meridian}}{def.meridian}
\nt{$\omega_j$}{the holonomy $\rho_\vth(\mer_j)$ around the co-core $\alpha_j$}{\eqref{eq.omega}}{eq.omega}
\nt{$\mathsf c_j(i)$}{the coefficient of the lifted meridian at the place $i$}{\eqref{eq.liftedpoints}}{eq.liftedpoints}
\nt{$\mer_j(i,l)$}{the part of the lifted meridian between the places $i$ and $l$; a half-turn if $l=i+2$, a turn if $l=i+4$}{\eqref{eq.pieces}}{eq.pieces}
\nt{$\CZ^\out_j$, $\CZ^\inn_j$}{the two interleaved chains of the crossing $j$}{Def.~\ref{def.interleaved}}{def.interleaved}
\nt{$\xs^\out_j$, $\ys^\out_j$, $\xs^\inn_j$, $\ys^\inn_j$, $\sv_j$}{the units which the interleaved chains give to $x_j$, $y_j$ and $x_jy_j$}{\eqref{eq.xyvalues}}{eq.xyvalues}
\nt{$\us^\out_j$, $\us^\inn_j$, $\us(\CZ_j)$}{the units in front of the forms of the interleaved chains}{Prop.~\ref{prop.interleaved}}{prop.interleaved}
\nt{$\CF_{\hate}$}{the set of the $c$ faces which are not adjacent to $\hate$}{Lem.~\ref{lem.chainpairing}}{lem.chainpairing}
\nt{$\lb\CZ_1,\dots,\CZ_c\rb_{\hate}$}{the mixed discriminant of $c$ elements of $C_\vth\otimes C_\vth$}{Lem.~\ref{lem.chainpairing}}{lem.chainpairing}
\nt{$\mathrm{col}$}{a two-coloring of the Seifert circles}{Lem.~\ref{lem.twocolored}}{lem.twocolored}
\nt{$\CJ$, $\CV(D,\CJ,\mathrm{col};\hate)$}{a set of crossings, the bands of a surface; the value of the construction for this surface with the vertical paths of $D$}{Def.~\ref{def.recipe}}{def.recipe}
\nh{Kashaev's matrix (Section~\ref{sec.faces})}
\nt{$w(f)$}{the winding number of $D$ around the face $f$}{Lem.~\ref{lem.winding}}{lem.winding}
\nt{$\zeta$}{Kashaev's variable, with $2\zeta=t^{1/2}+t^{-1/2}$}{\S\ref{ssec.kashaev}}{ssec.kashaev}
\nt{$\sigma_L(\omega)$}{the Levine--Tristram signature of $L$ at a complex number $\omega\ne1$ of modulus one; in Section~\ref{sec.cocores} the letter $\omega$ has another meaning}{\S\ref{ssec.kashaev}}{ssec.kashaev}
\nt{$\tau_D$, $\tau_D(t)$, $\tau_j$}{Kashaev's matrix and its blocks}{Def.~\ref{def.kashaev}}{def.kashaev}
\nt{$\WW_t$, $\DD_\ve$}{the diagonal matrices with the entries $(-t^{1/2})^{w(f)}$ and $\ve_j$}{\eqref{eq.W}}{eq.W}
\nt{$\tau_D(x,y)$, $\tau_{D,\hate}(x,y)$}{the deformed Kashaev matrix; the same without the rows and the columns of the two faces adjacent to $\hate$}{Def.~\ref{def.deformedkashaev}}{def.deformedkashaev}
\nh{The point $x=-1$ (Section~\ref{sec.minusone})}
\nt{$w^0_j$, $w^1_j$}{the two transfers of the crossing $j$}{\eqref{eq.w01}}{eq.w01}
\nt{$\bar\vs$}{the complementary smoothing $1-\vs$}{\S\ref{ssec.twotransfers}}{eq.UOw}
\nt{$c'_j$, $\phi'_j$, $c_j$, $\phi_j$}{the faces of the covered and of the free quadrant next to $u'_j$, and next to $u_j$}{Lem.~\ref{lem.cornertransfers}}{lem.cornertransfers}
\nt{$\Goe_\vs$, $\Stk_\vs$}{the Goeritz form and the stacking form of the smoothing $\vs$}{Def.~\ref{def.formsminusone}}{def.formsminusone}
\nt{$\GLf_\vs$, $\GLf^\uparrow_\vs$}{the forms of Gordon and Litherland of $F_\vs$ and of $F^\uparrow_\vs$, on the lattice of flows}{Def.~\ref{def.formsminusone}}{def.formsminusone}
\nt{$\Kash$}{Kashaev's form $\sum_j\ve_j\,U_j\otimes O_j$}{Def.~\ref{def.formsminusone}}{def.formsminusone}
\nt{$z^0$}{the alternating sum $\sum_C(-1)^{h(C)}C$ of the Seifert circles}{Prop.~\ref{prop.formsminusone}}{prop.formsminusone}
\nt{$\vs_{\rm b}$, $\vs_{\rm w}$}{the black and the white checkerboard smoothing}{\S\ref{ssec.goeritz}}{eq.eta}
\nt{$\eta_v(j)$, $G_v(D)$}{the sign of the crossing $j$ for the color $v$; the Goeritz matrix of the color $v$}{\eqref{eq.eta}}{eq.eta}
\nt{$\Omega(j)$, $\Omega_f$, $d(\Omega)$}{the region of a smoothing which contains the crossing $j$, resp.\ the face $f$ without its covered corners; the depth of a region}{\S\ref{ssec.sides}}{not.regions}
\nt{$H_\Omega$}{in Section~\ref{sec.minusone}: the summand of $H_1(F_\vs;\BZ)$ generated by the classes of the boundaries of the faces of the region $\Omega$ of the smoothing $\vs$}{Lem.~\ref{lem.smoothingregions}}{lem.smoothingregions}
\nt{$J$}{in Section~\ref{sec.minusone}: the involution of $H_1(F_\vs;\BZ)$ which is $(-1)^{d(\Omega)}$ on the homology of the region $\Omega$}{Prop.~\ref{prop.J}}{prop.J}
\nt{$n(L)$}{the corank of $V+V^{\T}$ for a Seifert matrix $V$ of $L$}{\S\ref{ssec.determinant}}{not.nL}
\nt{$G_\vs$}{a Gram matrix of the form $\GLf_\vs$}{Thm.~\ref{thm.determinant}}{thm.determinant}
\nt{$\mu_v(D)$, $\mathrm{wr}(D)$}{the correction term of the Goeritz matrix of the color $v$; the writhe $\sum_j\ve_j$}{Cor.~\ref{cor.GLgoeritz}}{cor.GLgoeritz}
\nt{$F_{\rm b}$, $F_{\rm w}$, $c_+$, $c_-$}{the numbers of black faces, of white faces, of positive crossings and of negative crossings}{Cor.~\ref{cor.altdefinite}}{cor.altdefinite}
\nt{$\Lambda_D$}{the lattice $U(\FF_D)\subseteq\BZ^c$ of the under-strands}{\S\ref{ssec.underlattice}}{ssec.underlattice}
\nh{Kauffman states (Sections~\ref{ssec.states}--\ref{ssec.statestrees})}
\nt{$\uparrow$, $\downarrow$, $A$, $B$}{the corners of a crossing: the two vertical corners, the $A$-corner between $o_j$ and $u'_j$, and the $B$-corner between $u_j$ and $o'_j$}{\eqref{eq.corneredges}}{eq.corneredges}
\nt{$\mathfrak s$, $\sgn(\mathfrak s)$}{a Kauffman state of $(D,\hate)$ and its sign}{Def.~\ref{def.state}}{def.state}
\nt{$X(\mathfrak s)$, $Y(\mathfrak s)$, $V(\mathfrak s)$}{the sets of the crossings at which the marker of $\mathfrak s$ is adjacent to $u'_j$, is adjacent to $o_j$, is vertical}{Def.~\ref{def.state}}{def.state}
\nt{$b(\mathfrak s)$}{the number of crossings at which the marker of $\mathfrak s$ is the corner $\downarrow$}{Prop.~\ref{prop.statesum}}{prop.statesum}
\nt{$m_{D,\hate}$}{the integer $\ve(X(\mathfrak s))+\ve(Y(\mathfrak s))$, which does not depend on the state}{Prop.~\ref{prop.statesum}}{prop.statesum}
\nt{$\Gamma$}{in Section~\ref{sec.kauffman}: the plane graph whose vertices are the crossings and the faces and whose edges are the corners}{\S\ref{ssec.ballinger}}{not.Gamma}
\nt{$Q_N$, $\ell(x,j)$, $k(x)$}{Ballinger's state model; the weight of a state $x$ with multiplicity two at the crossing $j$; the number of cycles of $x$}{\eqref{eq.ellweight}, \eqref{eq.Q2}}{eq.ellweight}
\nt{$a(\mathfrak s)$}{the number of crossings at which the marker of $\mathfrak s$ is the $A$-corner}{Lem.~\ref{lem.kasteleyn}}{lem.kasteleyn}
\nt{$\sigma_{D,\hate}$}{the sign $(-1)^{a(\mathfrak s)}$ of an alternating diagram}{Thm.~\ref{thm.kauffmanalt}}{thm.kauffmanalt}
\nh{The Tait graph and its lattices (Sections~\ref{ssec.tait}--\ref{ssec.lattice})}
\nt{$\CG$, $\CG^*$}{the Tait graph of an alternating diagram, whose vertices are the black faces and whose edges are the crossings; its dual}{\S\ref{ssec.tait}}{not.tait}
\nt{$\langle a,b\rangle$}{the scalar product $\sum_ja_jb_j$ of two vectors of $\BZ^c$}{Def.~\ref{def.lattices}}{def.lattices}
\nt{$\chi\in\{\pm1\}^c$}{in Section~\ref{sec.kauffman}: an orientation of the Tait graph}{Def.~\ref{def.lattices}}{def.lattices}
\nt{$\ve$, $\1_+$, $\1_-$}{the vector of the signs of the crossings; the indicator vectors of the positive and of the negative crossings}{Def.~\ref{def.lattices}}{def.lattices}
\nt{$T(\mathfrak s)$, $\chi(\mathfrak s)$, $\bar\chi(\mathfrak s)$}{the spanning tree of a state and its two orientations of the Tait graph}{Def.~\ref{def.chi}}{def.chi}
\nt{$r_{\rm b}$, $r_{\rm w}$}{the roots: the black and the white face adjacent to $\hate$}{Thm.~\ref{thm.tait}}{thm.tait}
\nt{$m_j(\cdot,\cdot)$}{the local weights of the state sum on the Tait graph}{Thm.~\ref{thm.tait}}{thm.tait}
\nt{$f^+_j$, $f^-_j$}{the black faces at which the edge of $j$ in $\CG$ ends and starts}{Def.~\ref{def.flowcut}}{def.flowcut}
\nt{$\mathrm{Flow}(\CG)$, $\mathrm{Cut}(\CG)$}{the lattices of the integral flows and of the integral cuts of $\CG$}{Def.~\ref{def.flowcut}}{def.flowcut}
\nt{$\Sigma(L)$}{the double cover of $S^3$ branched along $L$}{\S\ref{ssec.tait}}{not.glossary}
\nt{$d(Y,\mathfrak t)$, $\mathfrak t_0$}{the correction term of a $\mathrm{Spin}^c$ structure $\mathfrak t$ of a $3$-manifold $Y$; the spin structure}{\S\ref{ssec.tait}}{not.glossary}
\nt{$T$, $O$, $T\cup O$}{in Section~\ref{ssec.gluing}: the tangles inside and outside a disk $\Delta$, and their union}{\S\ref{ssec.gluing}}{not.tangles}
\nt{$\tau_{NS}$, $\tau_{vw}$, $\tau_0$, $\sigma_\pm$, $\kappa$}{in Section~\ref{ssec.gluing}: the standard tangles}{Fig.~\ref{fig.gluing}}{fig.gluing}
\nt{$H_+$, $H_-$}{the Links--Gould polynomials of the positive and of the negative Hopf link}{Thm.~\ref{thm.gluing}}{thm.gluing}
\nt{$\Lambda$}{the lattice $\mathrm{Flow}(\CG)\oplus\mathrm{Cut}(\CG)\subseteq\BZ^c$}{Prop.~\ref{prop.spinc}}{prop.spinc}
\nt{$[\chi]$, $\mathfrak t$}{the class of a vector $\chi$ with odd coordinates modulo $2\Lambda$; such a class, a $\mathrm{Spin}^c$ structure of $\Sigma(L)$}{Prop.~\ref{prop.spinc}}{prop.spinc}
\nt{$u_{\mathfrak t}$, $\bar u_{\mathfrak t}$}{the two representatives of the class $\mathfrak t$ for a root}{\eqref{eq.spincsum}}{eq.spincsum}
\nt{$E(\chi)$, $\bar E(\chi)$}{the sets of the edges for which $\chi$ is the first, resp.\ the second, orientation of a state}{\eqref{eq.rootsofchi}}{eq.rootsofchi}
\nt{$A(\chi)$, $A(\mathfrak t)$}{the set of the cyclic crossings of an orientation, and of a class}{Lem.~\ref{lem.minty}}{lem.minty}
\nt{$s_A(D)$}{the smoothing of the diagram defined by a set $A$ of crossings}{Cor.~\ref{cor.classsmoothing}}{cor.classsmoothing}
\nt{$\chi^\natural$, $\mathfrak t^\natural$, $M(\mathfrak t)$}{the orientation $\chi$ with its directed cycles reversed; the involution of the classes; the monomial of the pair $(\mathfrak t,\mathfrak t^\natural)$}{Prop.~\ref{prop.nat}}{prop.nat}
\nt{$T_{\CG}(x,y)$}{the Tutte polynomial of the graph $\CG$}{Thm.~\ref{thm.extreme}}{thm.extreme}
\end{longtable}
\endgroup
}


\clearpage
\begin{thebibliography}{GHKKST26}

\bibitem[Ale28]{Alexander28}
James~W. Alexander, \emph{Topological invariants of knots and links}, Trans. Amer. Math. Soc. \textbf{30} (1928), no.~2, 275--306.

\bibitem[ALvdV26]{ALvdV26}
Cristina Anghel, Daniel L\'opez Neumann and Roland van der Veen, \emph{The top coefficient of the Links--Gould invariant is Murasugi-sum multiplicative}, Preprint 2026, \href{https://arxiv.org/abs/2610.09354}{arXiv:2610.09354}.

\bibitem[BdlHN97]{BdlHN97}
Roland Bacher, Pierre de~la Harpe, and Tatiana Nagnibeda, \emph{The lattice of integral flows and the lattice of integral cuts on a finite graph}, Bull. Soc. Math. France \textbf{125} (1997), 167--198.

\bibitem[BBY19]{BackmanBakerYuen}
Spencer Backman, Matthew Baker, and Chi~Ho Yuen, \emph{Geometric bijections for regular matroids, zonotopes, and {E}hrhart theory}, Forum Math. Sigma \textbf{7} (2019), e45.

\bibitem[BW18]{BakerWang}
Matthew Baker and Yao Wang, \emph{The {B}ernardi process and torsor structures on spanning trees}, Int. Math. Res. Not. IMRN (2018), no.~16, 5120--5147.

\bibitem[Bal26]{BallingerTalk}
William Ballinger, \emph{Kauffman states with multiplicity}, talk in the AMS Special Session on Topology, Gauge Theory, and Special Geometric Structures, AMS Spring Eastern Sectional Meeting, 2026.

\bibitem[Bap89]{Bapat89}
Ravindra~B. Bapat, \emph{Mixed discriminants of positive semidefinite matrices}, Linear Algebra Appl. \textbf{126} (1989), 107--124.

\bibitem[BNM]{KnotAtlas}
Dror Bar-Natan, Scott Morrison, et al., \emph{The Knot Atlas}, \url{http://katlas.org}.

\bibitem[Ber08]{Bernardi08}
Olivier Bernardi, \emph{Tutte polynomial, subgraphs, orientations and sandpile model: new connections via embeddings}, Electron. J. Combin. \textbf{15} (2008), Research Paper 109.

\bibitem[Bro82]{Brown}
Kenneth~S. Brown, \emph{Cohomology of groups}, Graduate Texts in Mathematics, vol.~87, Springer-Verlag, New York--Berlin, 1982.

\bibitem[BZ03]{BurdeZieschang}
Gerhard Burde and Heiner Zieschang, \emph{Knots}, second ed., De Gruyter Studies in Mathematics, vol.~5, Walter de Gruyter, Berlin, 2003.

\bibitem[CHS25]{CHS25}
Zhechi Cheng, Matthew Hedden, and Sucharit Sarkar, \emph{Murasugi sum and extremal knot {F}loer homology}, Quantum Topol. \textbf{16} (2025), no.~4, 795--830.

\bibitem[CF24]{CimasoniFerretti24}
David Cimasoni and Livio Ferretti, \emph{On the Kashaev signature conjecture}, Fund. Math. \textbf{266} (2024), no.~3, 275--287.

\bibitem[CDR14]{CohenDasbachRussell}
Moshe Cohen, Oliver~T. Dasbach, and Heather~M. Russell, \emph{A twisted dimer model for knots}, Fund. Math. \textbf{225} (2014), 57--74.

\bibitem[Con70]{Conway70}
John~H. Conway, \emph{An enumeration of knots and links, and some of their algebraic properties}, Computational Problems in Abstract Algebra (Proc. Conf., Oxford, 1967), Pergamon, Oxford, 1970, pp.~329--358.

\bibitem[CF63]{CrowellFox}
Richard~H. Crowell and Ralph~H. Fox, \emph{Introduction to knot theory}, Graduate Texts in Mathematics, vol.~57, Springer-Verlag, New York--Heidelberg, 1977 (reprint of the 1963 original).

\bibitem[Cro59]{Crowell59}
Richard~H. Crowell, \emph{Genus of alternating link types}, Ann. of Math. (2) \textbf{69} (1959), 258--275.

\bibitem[Cro89]{Cromwell89}
Peter~R. Cromwell, \emph{Homogeneous links}, J. London Math. Soc. (2) \textbf{39} (1989), no.~3, 535--552.

\bibitem[DeW01]{DW01}
David De Wit, \emph{An infinite suite of Links--Gould invariants}, J. Knot Theory Ramifications \textbf{10} (2001), no.~1, 37--62.

\bibitem[DIL05]{DWIL}
David De Wit, Atsushi Ishii, and Jon Links, \emph{Infinitely many two-variable generalisations of the Alexander--Conway polynomial}, Algebr. Geom. Topol. \textbf{5} (2005), 405--418.

\bibitem[DWKL99]{DeWitKauffmanLinks}
David De~Wit, Louis~H. Kauffman, and Jon~R. Links, \emph{On the {L}inks--{G}ould invariant of links}, J. Knot Theory Ramifications \textbf{8} (1999), no.~2, 165--199.

\bibitem[Fox53]{Fox}
Ralph~H. Fox, \emph{Free differential calculus. {I}. {D}erivation in the free group ring}, Ann. of Math. (2) \textbf{57} (1953), 547--560.

\bibitem[FJR11]{FJR11}
Stefan Friedl, Andr\'as Juh\'asz, and Jacob Rasmussen, \emph{The decategorification of sutured {F}loer homology}, J. Topol. \textbf{4} (2011), no.~2, 431--478.

\bibitem[GHKKST26]{GHKKST}
Stavros Garoufalidis, Matthew Harper, Rinat Kashaev, Ben-Michael Kohli, Jiebo Song, and Guillaume Tahar,
\emph{Skein theory for the {L}inks--{G}ould polynomial}, J. Lond. Math. Soc. \textbf{113} (2026), e70522.

\bibitem[GHKST26]{GHKST}
Stavros Garoufalidis, Matthew Harper, Ben-Michael Kohli, Jiebo Song, and Guillaume Tahar, \emph{Extending knot
polynomials of braided {H}opf algebras to links}, SIGMA Symmetry Integrability Geom. Methods Appl. \textbf{22}
(2026), 050, 22 pages.

\bibitem[GK26]{GaroufalidisKashaev}
Stavros Garoufalidis and Rinat Kashaev, \emph{Multivariable knot polynomials from braided {H}opf algebras with
automorphisms}, Publ. Res. Inst. Math. Sci. \textbf{62} (2026), 75--114.

\bibitem[Gab83]{Gabai83}
David Gabai, \emph{The Murasugi sum is a natural geometric operation}, Low-dimensional topology (San Francisco, Calif., 1981), Contemp. Math., vol.~20, Amer. Math. Soc., Providence, RI, 1983, pp.~131--143.

\bibitem[Ghi08]{Ghiggini08}
Paolo Ghiggini, \emph{Knot {F}loer homology detects genus-one fibred knots}, Amer. J. Math. \textbf{130} (2008), no.~5, 1151--1169.

\bibitem[Gio07]{Gioan07}
Emeric Gioan, \emph{Enumerating degree sequences in digraphs and a cycle--cocycle reversing system}, European J. Combin. \textbf{28} (2007), 1351--1366.

\bibitem[Goe33]{Goeritz33}
Lebrecht Goeritz, \emph{Knoten und quadratische Formen}, Math. Z. \textbf{36} (1933), 647--654.

\bibitem[GL78]{GordonLitherland}
Cameron McA. Gordon and Richard~A. Litherland, \emph{On the signature of a link}, Invent. Math. \textbf{47} (1978), 53--69.

\bibitem[GL89]{GordonLuecke}
Cameron McA. Gordon and John Luecke, \emph{Knots are determined by their complements}, J. Amer. Math. Soc. \textbf{2} (1989), 371--415.

\bibitem[Gre13]{Greene13}
Joshua~Evan Greene, \emph{Lattices, graphs, and Conway mutation}, Invent. Math. \textbf{192} (2013), 717--750.

\bibitem[GZ83]{GreeneZaslavsky83}
Curtis Greene and Thomas Zaslavsky, \emph{On the interpretation of {W}hitney numbers through arrangements of hyperplanes, zonotopes, non-{R}adon partitions, and orientations of graphs}, Trans. Amer. Math. Soc. \textbf{280} (1983), no.~1, 97--126.

\bibitem[HMV24]{HMV24}
Elena~S. Hafner, Karola M\'esz\'aros, and Alexander Vidinas, \emph{On the Alexander polynomial of special alternating links}, Preprint 2024, \href{https://arxiv.org/abs/2401.14927}{arXiv:2401.14927}.

\bibitem[HK26]{HK}
Matthew Harper and Ben-Michael Kohli, \emph{The {L}inks--{G}ould polynomial of an alternating link is alternating}, Preprint 2026, \href{https://arxiv.org/abs/2610.18803}{arXiv:2610.18803}.

\bibitem[Hat02]{Hatcher}
Allen Hatcher, \emph{Algebraic topology}, Cambridge University Press, Cambridge, 2002.

\bibitem[Ish04]{Ishii04}
Atsushi Ishii, \emph{Algebraic links and skein relations of the {L}inks--{G}ould invariant}, Proc. Amer. Math. Soc. \textbf{132} (2004), no.~12, 3741--3749.

\bibitem[Juh06]{Juhasz06}
Andr\'as Juh\'asz, \emph{Holomorphic discs and sutured manifolds}, Algebr. Geom. Topol. \textbf{6} (2006), 1429--1457.

\bibitem[Juh08]{Juhasz08}
Andr\'as Juh\'asz, \emph{Floer homology and surface decompositions}, Geom. Topol. \textbf{12} (2008), no.~1, 299--350.

\bibitem[Kas21]{Kashaev21}
Rinat Kashaev, \emph{On symmetric matrices associated with oriented link diagrams}, Topology and geometry---a collection of essays dedicated to Vladimir G. Turaev, IRMA Lect. Math. Theor. Phys., vol.~33, EMS Press, Berlin, 2021, pp.~131--146.

\bibitem[Kas63]{Kasteleyn63}
Pieter~W. Kasteleyn, \emph{Dimer statistics and phase transitions}, J. Math. Phys. \textbf{4} (1963), 287--293.

\bibitem[Kau83]{KauffmanFKT}
Louis~H. Kauffman, \emph{Formal knot theory}, Mathematical Notes, vol.~30, Princeton University Press, Princeton, 1983.

\bibitem[Kau87]{Kauffman87}
Louis~H. Kauffman, \emph{State models and the {J}ones polynomial}, Topology \textbf{26} (1987), no.~3, 395--407.

\bibitem[Kea75]{Kearton}
Cherry Kearton, \emph{Blanchfield duality and simple knots}, Trans. Amer. Math. Soc. \textbf{202} (1975), 141--160.

\bibitem[Kho00]{Khovanov00}
Mikhail Khovanov, \emph{A categorification of the {J}ones polynomial}, Duke Math. J. \textbf{101} (2000), no.~3, 359--426.

\bibitem[KS17]{KhovanovSussan}
Mikhail Khovanov and Joshua Sussan, \emph{A categorification of the positive half of quantum $gl(m|1)$}, Trans. Amer. Math. Soc. \textbf{369} (2017), no.~3, 1627--1664.

\bibitem[KK96]{KobayashiKobayashi96}
Masako Kobayashi and Tsuyoshi Kobayashi, \emph{On canonical genus and free genus of knot}, J. Knot Theory Ramifications \textbf{5} (1996), no.~1, 77--85.

\bibitem[KS91]{KauffmanSaleur}
Louis~H. Kauffman and Hubert Saleur, \emph{Free fermions and the {A}lexander--{C}onway polynomial}, Comm. Math. Phys.
\textbf{141} (1991), no.~2, 293--327.

\bibitem[Koh16]{Kohli}
Ben-Michael Kohli, \emph{On the {L}inks--{G}ould invariant and the square of the {A}lexander polynomial}, J. Knot
Theory Ramifications \textbf{25} (2016), no.~2, 1650006.

\bibitem[KPM17]{KPM}
Ben-Michael Kohli and Bertrand Patureau-Mirand, \emph{Other quantum relatives of the Alexander polynomial through the Links--Gould invariants}, Proc. Amer. Math. Soc. \textbf{145} (2017), no.~12, 5419--5433.

\bibitem[KT57]{KinoshitaTerasaka}
Shin'ichi Kinoshita and Hidetaka Terasaka, \emph{On unions of knots}, Osaka Math. J. \textbf{9} (1957), 131--153.

\bibitem[KT26]{Kohli-Tahar}
Ben-Michael Kohli and Guillaume Tahar, \emph{A lower bound for the genus of a knot using the Links--Gould invariant}, Quantum Topol. (2026), published online first, \href{https://arxiv.org/abs/2310.15617}{arXiv:2310.15617}.

\bibitem[LG92]{LinksGould}
Jon~R. Links and Mark~D. Gould, \emph{Two variable link polynomials from quantum supergroups}, Lett. Math. Phys.
\textbf{26} (1992), no.~3, 187--198.

\bibitem[Lic97]{Lickorish}
W.~B.~Raymond Lickorish, \emph{An introduction to knot theory}, Graduate Texts in Mathematics, vol.~175, Springer-Verlag, New York, 1997.

\bibitem[Liu25]{Liu25}
Jessica Liu, \emph{A proof of the Kashaev signature conjecture}, Quantum Topol. (2025), published online first, \href{https://doi.org/10.4171/QT/232}{doi:10.4171/QT/232}.

\bibitem[Liv88]{Livingston88}
Charles Livingston, \emph{The free genus of doubled knots}, Proc. Amer. Math. Soc. \textbf{104} (1988), no.~1, 329--333.

\bibitem[LNvdV25]{LNVdV25}
Daniel Lopez-Neumann and Roland van der Veen, \emph{A plumbing-multiplicative function from the Links--Gould invariant}, Preprint 2025, 
\href{https://arxiv.org/abs/2502.12899}{arXiv:2502.12899}.

\bibitem[Man19]{Manion19}
Andrew Manion, \emph{On the decategorification of {O}zsv\'ath and {S}zab\'o's bordered theory for knot {F}loer homology}, Quantum Topol. \textbf{10} (2019), no.~1, 77--206.

\bibitem[MO07]{ManolescuOwens}
Ciprian Manolescu and Brendan Owens, \emph{A concordance invariant from the Floer homology of double branched covers}, Int. Math. Res. Not. IMRN \textbf{2007}, no.~20, Art.~ID rnm077.

\bibitem[Mil68]{Milnor}
John~W. Milnor, \emph{Infinite cyclic coverings}, Conference on the Topology of Manifolds (Michigan State Univ., E.~Lansing, Mich., 1967), Prindle, Weber \& Schmidt, Boston, Mass., 1968, pp.~115--133.

\bibitem[Min66]{Minty66}
George~J. Minty, \emph{On the axiomatic foundations of the theories of directed linear graphs, electrical networks and network-programming}, J. Math. Mech. \textbf{15} (1966), 485--520.

\bibitem[Mor87]{Moriah87}
Yoav Moriah, \emph{On the free genus of knots}, Proc. Amer. Math. Soc. \textbf{99} (1987), no.~2, 373--379.

\bibitem[MP89]{MurasugiPrzytycki89}
Kunio Murasugi and J\'ozef~H. Przytycki, \emph{The skein polynomial of a planar star product of two links}, Math. Proc. Cambridge Philos. Soc. \textbf{106} (1989), no.~2, 273--276.

\bibitem[Mur58]{Murasugi58}
Kunio Murasugi, \emph{On the genus of the alternating knot. {I}, {II}}, J. Math. Soc. Japan \textbf{10} (1958), 94--105, 235--248.

\bibitem[Mur63]{Murasugi63}
Kunio Murasugi, \emph{On a certain subgroup of the group of an alternating link}, Amer. J. Math. \textbf{85} (1963), 544--550.

\bibitem[Mur87a]{Murasugi87}
Kunio Murasugi, \emph{Jones polynomials and classical conjectures in knot theory}, Topology \textbf{26} (1987), no.~2, 187--194.

\bibitem[Mur87b]{Murasugi87b}
Kunio Murasugi, \emph{Jones polynomials and classical conjectures in knot theory. {II}}, Math. Proc. Cambridge Philos. Soc. \textbf{102} (1987), no.~2, 317--318.

\bibitem[Ni06]{Ni06}
Yi Ni, \emph{Sutured {H}eegaard diagrams for knots}, Algebr. Geom. Topol. \textbf{6} (2006), 513--537.

\bibitem[Ni07]{Ni07}
Yi Ni, \emph{Knot {F}loer homology detects fibred knots}, Invent. Math. \textbf{170} (2007), no.~3, 577--608.

\bibitem[OS03]{OSalternating}
Peter Ozsv\'ath and Zolt\'an Szab\'o, \emph{Heegaard {F}loer homology and alternating knots}, Geom. Topol. \textbf{7} (2003), 225--254.

\bibitem[OS04a]{OSknot}
Peter Ozsv\'ath and Zolt\'an Szab\'o, \emph{Holomorphic disks and knot invariants}, Adv. Math. \textbf{186} (2004), no.~1, 58--116.

\bibitem[OS04b]{OSgenus}
Peter Ozsv\'ath and Zolt\'an Szab\'o, \emph{Holomorphic disks and genus bounds}, Geom. Topol. \textbf{8} (2004), 311--334.

\bibitem[OS05]{OSbranched}
Peter Ozsv\'ath and Zolt\'an Szab\'o, \emph{On the Heegaard Floer homology of branched double-covers}, Adv. Math. \textbf{194} (2005), 1--33.

\bibitem[OS18]{OSbordered}
Peter Ozsv\'ath and Zolt\'an Szab\'o, \emph{Kauffman states, bordered algebras, and a bigraded knot invariant}, Adv. Math. \textbf{328} (2018), 1088--1198.

\bibitem[Ras03]{RasmussenThesis}
Jacob Rasmussen, \emph{Floer homology and knot complements}, Ph.D. thesis, Harvard University, 2003, \href{https://arxiv.org/abs/math/0306378}{arXiv:math/0306378}.

\bibitem[Ras15]{Rasmussen15}
Jacob Rasmussen, \emph{Some differentials on {K}hovanov--{R}ozansky homology}, Geom. Topol. \textbf{19} (2015), no.~6, 3031--3104.

\bibitem[RW22]{RobertWagner}
Louis-Hadrien Robert and Emmanuel Wagner, \emph{A quantum categorification of the {A}lexander polynomial}, Geom. Topol. \textbf{26} (2022), no.~5, 1985--2064.

\bibitem[Rol76]{Rolfsen}
Dale Rolfsen, \emph{Knots and links}, Mathematics Lecture Series, vol.~7, Publish or Perish, Berkeley, 1976.

\bibitem[Rud83]{Rudolph83}
Lee Rudolph, \emph{Braided surfaces and Seifert ribbons for closed braids}, Comment. Math. Helv. \textbf{58} (1983), no.~1, 1--37.

\bibitem[Sta78]{Stallings78}
John~R. Stallings, \emph{Constructions of fibred knots and links}, Algebraic and geometric topology (Stanford, 1976), Part~2, Proc. Sympos. Pure Math., vol.~XXXII, Amer. Math. Soc., Providence, 1978, pp.~55--60.

\bibitem[Thi87]{Thistlethwaite87}
Morwen~B. Thistlethwaite, \emph{A spanning tree expansion of the {J}ones polynomial}, Topology \textbf{26} (1987), no.~3, 297--309.

\bibitem[Thi88]{Thistlethwaite88}
Morwen~B. Thistlethwaite, \emph{Kauffman's polynomial and alternating links}, Topology \textbf{27} (1988), no.~3, 311--318.

\bibitem[Tor53]{Torres}
Guillermo Torres, \emph{On the {A}lexander polynomial}, Ann. of Math. (2) \textbf{57} (1953), 57--89.

\bibitem[Tra04]{Traczyk04}
Pawe{\l} Traczyk, \emph{A combinatorial formula for the signature of alternating diagrams}, Fund. Math. \textbf{184} (2004), 311--316.

\bibitem[Tro73]{Trotter}
Hale~F. Trotter, \emph{On $S$-equivalence of Seifert matrices}, Invent. Math. \textbf{20} (1973), 173--207.

\bibitem[Tur97]{Turaev97}
Vladimir Turaev, \emph{Torsion invariants of {$\mathrm{Spin}^c$}-structures on $3$-manifolds}, Math. Res. Lett. \textbf{4} (1997), no.~5.

\bibitem[Whi32]{Whitney32}
Hassler Whitney, \emph{Congruent graphs and the connectivity of graphs}, Amer. J. Math. \textbf{54} (1932), no.~1, 150--168.

\bibitem[Whi33]{Whitney33}
Hassler Whitney, \emph{2-{I}somorphic graphs}, Amer. J. Math. \textbf{55} (1933), 245--254.

\end{thebibliography}
\end{document}